\documentclass[11pt,english,a4paper]{smfbook}
\usepackage[T1]{fontenc}
\usepackage[utf8]{inputenc}
\usepackage{mathrsfs,amsmath,latexsym,amssymb,mathtools,dsfont}
\usepackage{arydshln}
\usepackage{bm}
\usepackage{hyperref}
\usepackage{smfhyperref}
\usepackage{smfthm}
\usepackage{xcolor}
\usepackage{xspace}
\usepackage{calc}
\usepackage[bb=boondox,bbscaled=1.04]{mathalfa}

\def\0{\bm{0}}
\def\1{\mathbb{1}}

\newcommand{\uppointer}[2]{\underset{\overset{\uparrow}{#1}}{#2}}
\newcommand{\abs}[1]{\left| {#1} \right|}
\newcommand{\adjstud}[1]{\rule{0pt}{\heightof{#1}}}
\DeclareMathOperator{\rank}{rank}
\DeclareMathOperator{\vect}{vect}

\newcommand{\as}{almost surely\xspace}
\newcommand{\ba}{\bm{a}}
\newcommand{\bA}{\bm{A}}

\newcommand{\bB}{\bm{B}}
\newcommand{\bbE}{\mathbb{E}}
\newcommand{\bbeta}{\bm{\beta}}

\newcommand{\bC}{\bm{C}}
\newcommand{\bD}{\bm{D}}
\newcommand{\bdK}{\bm{K}}
\newcommand{\be}{\bm{e}}
\newcommand{\bE}{\bm{E}}
\newcommand{\bF}{\bm{F}}
\newcommand{\bH}{\bm{H}}
\newcommand{\bk}{\bm{k}}
\newcommand{\bK}{\bm{K}}

\newcommand{\bmm}{\bm{m}}
\newcommand{\calA}{\mathcal{A}}
\newcommand{\calB}{\mathcal{B}}
\newcommand{\calC}{\mathcal{C}}
\newcommand{\calD}{\mathcal{D}}
\newcommand{\calF}{\mathcal{F}}
\newcommand{\calH}{\mathcal{H}}
\newcommand{\calI}{\mathcal{I}}
\newcommand{\calJ}{\mathcal{J}}
\newcommand{\card}{\operatorname{card}}
\newcommand{\cf}{{cf$.$}\xspace}

\newcommand{\Cov}{\mathop{\rm Cov}\nolimits}

\newcommand{\ie}{i.e$.$\xspace}
\newcommand{\cvg}[2][{n \to \infty }]{\xrightarrow[#1]{\text{\upshape\scriptsize #2}}}
\newcommand{\diag}{\ensuremath{\mathop{\rm diag}}}
\newcommand{\diam}{\ensuremath{\mathop{\operatorname{diam}}}}
\newcommand{\ds}{\displaystyle}
\newcommand{\entiers}{\mathbb{ Z}}
\newcommand{\calK}{\mathcal{K}}
\newcommand{\calL}{\mathcal{L}}
\newcommand{\calN}{\mathcal{N}}
\newcommand{\bL}{\mathbb{L}}
\newcommand{\frakL}{\mathfrak{L}}
\newcommand{\frakB}{\mathfrak{B}}

\newcommand{\given}[1]{\left|\kern1.5pt\rule{0pt}{10pt}#1\right.}
\newcommand{\naturels}{\mathbb{ N}}
\newcommand{\normp}[2][p]{\left\| {#2} \right\|_{#1}} 
\renewcommand{\paragraph}[1]{\medskip\noindent\textit{#1}}

\newcommand{\bbP}{\mathbb{P}}
\newcommand{\rv}{random variable\xspace}
\newcommand{\reels}{\mathbb{R}}

\newcommand{\cyl}{cyl}
\newcommand{\supp}{\mathop{\operatorname{supp}}\nolimits}
\newcommand{\bbS}{\mathbb{S}}
\newcommand{\bs}{\bm{s}}
\newcommand{\calS}{\mathcal{S}}
\newcommand{\bP}{\bm{P}}
\newcommand{\bS}{\bm{S}}
\newcommand{\bbT}{\mathbb{T}}
\newcommand{\bbV}{\mathbb{V}}
\newcommand{\bt}{\bm{t}}

\newcommand{\bW}{\bm{W}}
\newcommand{\bx}{\bm{x}}
\newcommand{\by}{\bm{y}}
\newcommand{\bz}{\bm{z}}
\newcommand{\calZ}{\mathcal{Z}}
\newcommand{\ud}{{\,{\rm d}}}
\newcommand{\Var}{\mathop{\rm Var}\nolimits}
\newcommand{\ve}{{\varepsilon}}

\newcommand{\introo}[3][,]{{]}{#2}#1\,{#3}{[}}

\newcommand{\intrfo}[3][,]{{[}{#2}#1\,{#3}{[}}
\newcommand{\intrff}[3][,]{{[}{#2}#1\,{#3}{]}}

\newcommand{\sstar}{\ensuremath{*}}
\renewcommand{\le}{\ensuremath{\leqslant}}
\renewcommand{\ge}{\ensuremath{\geqslant}}

\newcommand{\such}{\text{ such that }}

\newcommand{\prodsca}[2]
	{\left\langle#1{,\,}#2\right\rangle}
\newenvironment{proofarg}[1]
   {\smallskip\par\noindent
    \textit{#1.}\ }
   {\hfill{$\Box$}\par\noindent}
   
\def\Xint#1{\mathchoice
    {\XXint\displaystyle\textstyle{#1}}%
    {\XXint\textstyle\scriptstyle{#1}}%
    {\XXint\scriptstyle\scriptscriptstyle{#1}}%
    {\XXint\scriptscriptstyle\scriptscriptstyle{#1}}%
      \!\int}
\def\XXint#1#2#3{{\setbox0=\hbox{$#1{#2#3}{\int}$}
    \vcenter{\hbox{$#2#3$}}\kern-.5\wd0}}
\def\repint{\Xint{\kern4.05pt\cdots}}
   
\usepackage{wasysym}

\reversemarginpar

\makeatletter
\newcommand{\specialpos}[1]{\ifmeasuring@#1\else\omit$\displaystyle#1$\ignorespaces\fi}
\makeatother

\allowdisplaybreaks 
\newcommand{\spacebefore}{\vglue6pt}
\begin{document}
\title{Kac-Rice formula: a contemporary overview of the main results and applications}
\author{Corinne Berzin\\ Alain Latour\\ Jos\'e R. Le\'on}
\date{}
\maketitle
\frontmatter
\chapter*{Foreword}
\begin{flushright}
To Enrique Caba\~na\\
\end{flushright}
In these notes, we study the behavior of the level set 
$\calC_{X}(y):=\{x\in D\subset \reels ^d: X(x)=y\}$
for a fairly smooth random field
$X:\Omega\times D\subset\reels ^d\to\reels ^j$ with $d\ge j$.
We are interested in establishing a formula for the expectation and for the higher order moments of the $(d-j)$-area measure of the level set $\calC_{X}(y)$, \ie, $\sigma_{d-j}(\calC_{X}(y))$, as well as formulas for level set functionals.

Let us give the historical context.
The study of level sets of Gaussian random processes began in the 1940s with two seminal papers by Mark Kac \cite{MR7812} and  Stephen 0.~Rice \cite{MR10932}, respectively.
In both articles, a formula was established to calculate the expectation of the number of zeros of a Gaussian process $X:\Omega\times\reels \to\reels $.
It should be noted, however, that Kac was interested in the asymptotic study of the number of roots of a random polynomial while Rice sought such a formula for imperatively Gaussian models in communication theory.
Over the next decade, both works were generalized and deepened.
Kac's work was extended to calculate the asymptotic variance of the number of roots and a central limit theorem was also proved.
Rice's work was also generalized.
Cramér and Leadbetter began a fruitful collaboration on this subject that resulted in their very famous book \cite{MR2108670}.

In 1957, M.~S. Longuet-Higgins, in an article devoted to the modeling of the sea as a random surface \cite{MR87257}, gives for the first time a Kac-Rice formula for Gaussian fields, $X:\Omega\times\reels ^2\to\reels $, which can be considered as a great success.
The work, while very ingenious, provided no formal mathematical proof.
Nevertheless, this article was intensely cited after its appearance and its imprint on subsequent work for the study of level sets for random fields was very important.

In the early eighties of the last century, three other key works were published.
The article entitled \textit{``Esperanzas de integrales sobre conjuntos de nivel aleatorios''} by E.
Cabaña, \cite{AL121221}, that the reader will have the opportunity to revisit with these notes, a Springer Lecture Notes  \cite{MR871689} written by M.
Wschebor ``Surfaces aléatoires'' and finally the Robert Adler's book \textit{``The Geometry of Random Fields''} \cite{MR3396215}.
These three works deal with the extension of the ideas of Kac and Rice to multidimensional random fields.

These themes, somewhat \textit{``exotic''} for the time they were studied, have been revived in the 21st century.
We have seen the appearance of two books,\cite{MR2319516} and \cite{MR2478201}, that gave a new impetus to the subject.

The well-established Kac-Rice formulas in the 1980s have been proven again and improved, revealing unsuspected connections.
However, the most important aspect of this recent revival is the application of Kac-Rice formulas in several areas of pure and applied mathematics.
With these notes, different in its accent from the two aforementioned books, we want to introduce the reader into an old and at the same time young field, highlighting a significant number of applications.

These notes come from a course given by one of the authors at the \textit{XXX Escuela Venezolana de Matemática} that took place at Mérida, Venezuela in September $2017$.
For this edition, they have been considerably expanded and corrected.
\tableofcontents
%
%
%
\mainmatter
\chapter{Some historical remarks about Kac-Rice formula}
The well-known change of variable formula in an integral of a real function can be written as follows.
Let $G:\reels \to \reels $ a monotone and differentiable function and $I=G([a,b])$.
If $f$ is a continuous function, we have
$$\int_If(y)\ud y=\int_a^bf(G(t))\abs{G^{\prime}(t)}\ud t.$$
A similar formula is also true when $G$ is differentiable but not necessarily monotone.
In such a case, the formula must be modified by introducing the number of crossings of the level $y $ for the function $G $.
We define the latter quantity as follows:
$$
	N^G_{[a,b]}(y):=\#\{t\in[a,b]: G(t)= y\}.
$$
leading to the following formula 
$$
	\int_If(y)N^G_{[a,b]}(y)\ud y=\int_a^bf(G(t))\abs{G^{\prime}( t)}\ud t.
$$
This formula is known as the area formula and has a long history.
Its first proof is presented in an article by S{$.$} Banach \cite{AL111221}.

This formula can be used to obtain another formula known as the Kac counter \cite{MR1575245}.
In fact, if the function $G$ is continuously differentiable and has a finite number of crossings at the $u$ level, then
\begin{equation}
	N^G_{[a,b]}(u)=\lim_{\delta\to0}\frac1{2\delta}\int_{a}^{b}\1_{[u-\delta, u+\delta]}(G(s))\abs{G^{\prime}(s)}\ud s.
	\label{kac}
\end{equation}

Kac used this formula to obtain the expectation of the number of real roots of a random polynomial.
In fact, let $\{a_i\}_{i=0}^n$ be a set of independent standard Gaussian random variables.
Let us define for each $n$ the following random polynomial
$$
	X_n(t):=a_0+a_1 t+\dots+a_n t^n.
$$

Kac was interested in the expectation of the random variable $N^{X_n}_{[0,T]}(0)$.
 Using the counter (\ref{kac}), knowing that the number of real roots of $X_n$ and $X _n^\prime$ is bounded by $n$, we see that almost surely, we have
$$
	N^{X_n}_{[0,T]}(0)=\lim_{\delta\to0}\frac1{2\delta}\int_{0}^{T}\1_{[-\delta,+\delta]}(X_n(s))\abs{X_n^{\prime}(s)}\ud s,
$$
and moreover, the right side is also bounded.
The dominated convergence theorem implies
\begin{align*}
	\bbE\!\left[N^{X_n}_{[0,T]}(0)\right]&= \lim_{\delta\to0}\frac1{2\delta}\int_0^T\int_{-\delta}^{\delta}
	\int_{\reels} \abs{z}{p}_{X_n(s),X^{\prime}_n(s)}(x, z)\ud z \ud x \ud s\\
	&=\int_0^T\int_{\reels} \abs{z}{p}_{X_n(s),X^{\prime}_n(s)}(0,z)\ud z \ud s,
\end{align*}
where
${p}_{X_n(s),X^{\prime}_n(s)}(x, z)$
is the density of the Gaussian vector
$(X_n(s),X_n^{\prime}(s))$.

This is the famous Kac-Rice formula.
This name comes from the fact that Rice, considering the zeros of random processes, also gives a proof of this result for Gaussian processes \cite{MR11918}.
Nevertheless, Rice's interest was not in random polynomials, but in problems related to random signals.
In the case studied by Rice, more work is needed because, in general, the processes are not \as bounded.

In a historical article \cite{MR993433}, the author claims that in fact, the formula was obtained by Rice long before the publication of his articles \cite{MR10932} and \cite{MR11918}, which remain the main references to the first proof of the formula.

Rice's proof, although original and very suggestive, was supplemented by It\^o in \cite{MR166824} who gives, for stationary Gaussian processes, a necessary and sufficient condition for the number of crossings to be finite.
The It\^o's condition requires the finiteness of the second spectral moment of the process.

In 1967, some time after the publication of the aforementioned works, Cram\'er \& Leadbetter's book \cite{MR2108670} was published.
In this book, not only was a general proof of the Kac-Rice formula given, but the higher moments were also considered, establishing a formula for factorial moments of the number of crossings for Gaussian processes.
In addition, the theory has been enriched and supplemented by an interesting series of applications.

Besides, in the text, a necessary condition was given to ensure that the variance of the $N^X_{[0,T]}(0)$ was finite, $X$ being a stationary Gaussian process.
To do this, the authors used Kac-Rice's formula for the second factorial moment of $N^X_{[0,T]}(0)$.
The condition was expressed as follows: $\exists \delta>0$ such that
$$
	\int_0^\delta \frac{\abs{r^{\prime \prime}(t)-r^{\prime \prime}(0)}}{t}\ud t< \infty,
$$
where $r$ is the covariance function of $X$.
A little later, D{.} ~Geman \cite{MR350833} showed that this condition was also sufficient.

At that time, the study of level sets for random fields was not very popular.
A brilliant exception is the seminal work of M.S.~Longuet-Higgins \cite{MR87257} appeared in $1957$ and in which, among several applications of the level crossings to sea modeling, a Kac-Rice formula for the length of the level curve for a stationary Gaussian random field $X:\reels^2\to\reels$ was firmly established.

Apart from the article mentioned above, interest in the functional of level sets for a random field only emerged in the last seventies and eighties of the last century.
For example, we can cite the founding article \cite{MR405559}, where Adler and Hasofer extended the notion of crossing levels using the Euler characteristic of the excursion set  from a Gaussian random field $X:\reels^2 \to \reels$.
Also, it is important to mention the work of Benzaqu\'en \& Caba\~na \cite{MR687959}, where for the first time a Kac-Rice formula for the measure of the level set of a Gaussian random field $X: 
\reels^d\to \reels$ was obtained, generalizing the result of the aforementioned article \cite{MR87257}.

In the eighties, two books dealing with the levels set for random fields were published.
The first book, written by R. Adler \cite{MR3396215}, deals with problems related to the geometry of random fields.
The Euler characteristic of the excursion set  $A_u(X,S):=\{s \in S: X(s)\ge u\}$ for a Gaussian random field $X:\reels^d\to \reels$ was introduced and a Kac-Rice formula for the expectation of this functional has been proved.
In addition, a formula for the expected number of maxima above the $u$ level of $X(t)$ for $t \in S$ has been given.
Subsequently, several applications of the last mentioned formula were provided.
In particular, the author obtains an approximation of the tail of the distribution of the random variable $\max_{t\in S}X(t)$.
The second book, is a Springer Lecture Notes written by M.~Wschebor \cite{MR871689} where the expectation of the Lebesgue measure of the level set $\calC_{S, X}(u):=\{s\in S: X(s)= u\}$ was computed.
The results contained in this book generalize those of the article \cite{MR687959}, showing a Kac-Rice formula for more general Gaussian fields, demonstrating a formula for higher moments and also considering non-Gaussian processes.

The Kac-Rice formulas cited above are based on a geometric measurement formula known as the coarea formula.
This formula will be established and proven at the beginning of this book.
The coarea formula is one of the main tools in Federer's work \cite{MR0257325} and other mathematicians (the corresponding references can be found in the bibliography of this latter book).

To facilitate understanding, we will give the formula and sketch an application.
Let $G$ be a differentiable function $G:D\subset \reels^d\to \reels^j, d\ge j$, and $D$ an open set of $\reels^d$.
Let $\nabla G(\bm{\cdot})$ denote its Jacobian.
If $f:\reels^j\to\reels^{+}$ then for a Borel subset $B\subset D$,
$$
	\int_{B}f(G(x)) (\det(\nabla G(x)\nabla G(x)^T))^{\frac12}\ud x=\int_{\reels^j}f(y)\calH_{d-j}(\calC_{B,G}(y))\ud y
$$
where $\calH_{d-j}$ is the Hausdorff measure in dimension $d-j$ and $\calC_{B,G}(y):=\{x \in B: G(x)=y\}$.
In the case we develop, we define $j = 1$, $f = 1_A$ for $A$ a Borel set of $\reels$ and $G = X$ a continuously differentiable Gaussian random field and its gradient satisfies a Hölder condition.
In this case, $\calC_{B,X}(y)$ is \as~a differentiable variety of co-dimension one (see \cite{MR871689}) and the formula is
$$
\int_{B}\1_{\{X(x) \in A\}} \normp[d]{\nabla\!X(x)}\ud x=\int_{A} \sigma_{d-1}(\calC_{B,X}(y))\ud y.
$$

Then, taking the expectation on both sides and assuming that for any bounded Borel set  $A$ one of the two integrals is finite, using Fubini's theorem and duality, we get for almost all $y$
\begin{equation}
\label{kr1}
\bbE\!\left[\sigma_{d-1}(\calC_{B,X}(y))\right]=\int_{B}\bbE[\normp[d]{\nabla\!X(x)}\given{X(x)=y}]p_{X(x)}(y)\ud x.
\end{equation}

In \cite{MR687959} for stationary Gaussian fields and in \cite{MR871689} for more general and smooth Gaussian fields, conditions are given for the formula to be true for all $y$.
In general, it is not trivial to switch from the formula (\ref{kr1}) to the formula given for all $y$.
Different strategies have been developed and the explanation of a possible way is one of the main interests of this work.

The formula (\ref{kr1}) for $d\ge j\ge1$ was studied by Caba\~na in \cite{AL121221}.
The article was written in Spanish and had not been widely distributed.
In the more recent book \cite{MR2478201} a proof has been sketched, and in these notes, we will generalize and complete these two nice initiatives.

In the twenty-first century, two books appeared \cite{MR2319516} and \cite{MR2478201} that gave new impetus to the subject.
New fields of application of the formulas have appeared in the literature and the question has become a vast research field.
Examples include applications to the number of roots of random polynomial systems (algebraic or trigonometric) and also to the volume of nodal sets for rectangular systems.

The current literature is considerable.
We can mention among others: \cite{MR2797349} and \cite{MR3084654} for the study of random trigonometric real polynomials and \cite{MR3010810} which considers the length of the set of zeros of the random trigonometric polynomial in the torus $\bbT ^2$, also \cite{MR2182444}, \cite[Chapter 12]{MR2478201} and \cite{MR3866880} for the zeros of Kostlan-Shub-Smale polynomials or systems and \cite{MR3470435} for the study of zeros for more general polynomials and also the zeros of  other geometric objects.
Kac-Rice formulas are also a basic tool for studying sets of zeros of random waves and it has been very important to prove or disprove Berry's conjectures \cite{MR1913853}, see \cite{MR3540457} and the references therein.

One area of application where the formulas have been very useful is random modelling of the sea.
This domain presents the first application of Kac-Rice's formula for random fields in the nice article  \cite{MR87257}.
In addition, the Lund's School of probability has been very active in these areas, see for example \cite{MR2827442} and the references contained therein.
The Kac-Rice formula is not valid for random fields with non-differentiable trajectories.
However, the study of some level functionals for these fields as the local time, was implemented by approximating the actual field by a regular field using a convolution with an approximate Dirac delta.
In addition, the Kac-Rice formula for a functional level of the smoothed field makes it possible to approximate a level functional  of the original field.
Thus, the random fields to which the level sets are studied can have their domain in a finite-dimensional manifold; see \cite{MR3470435} for example.

We have presented a panorama, perhaps a little fast, of the genesis and development of Kac-Rice formulas and their deep imprint in the study of sets of zeros of random functions.
In these notes, we will try to make accessible to postgraduate students and researchers in probability and statistics and related fields the basic ideas for demonstrating the formulas.
We will also insist on its applications, both those of the first epoch and those of the present times.
We hope the reader will enjoy reading them as much as we did while writing them.
%
%
%
\chapter{A  proof of the coarea formula}
\section{Preliminaries}
As mentioned earlier, in these notes we study the Kac-Rice formulas dealing with the expectation of the level set measure for random fields or processes.
In what follows, we describe the different parts of this work.

First, there are two variants of the change of variables formula for multiple integrals that are very useful in integral geometry.
The first one corresponds to smooth and locally bijective functions $G:\reels ^d\to \reels ^d$ and the second one applies to smooth functions $G:\reels ^d\to \reels ^j$ with $d>j$, having a differential of maximal rank.
These formulas are called respectively \textit{area formula} and \textit{coarea formula}.
By applying these formulas to the trajectories of the random fields, then taking the expectation, we obtain the well known Kac-Rice formulas.
Recently and mainly thanks to the publication of two excellent books (\cite{MR2319516} and \cite{MR2478201}), the application of these formulas has aroused a growing interest in fields as varied as: random algebraic geometry, the complexity of algorithms for solving large systems of equations, the study of zeros of random polynomial systems, and finally, applications in engineering.
We now present the parts of this book.
\begin{enumerate}
    \item In the first part, we give an analytical proof of the area and coarea formulas.
    	Such a proof, originally attributed to  Banach and Federer \cite{MR0257325}, uses elementary tools of vector calculus and measure theory in $\reels ^d$.
     \item The above formulas are the basis for establishing the validity of Kac-Rice formulas for random fields.
    They allow to compute the expectation of the measure of the level sets 
    $$
    	\calC_{Q,X}( y):=\{t\in Q\subset \reels ^d : X(t)= y\},
    $$
    where $X:\Omega\times\reels ^d\to\reels ^j$ is a random field and $d\ge j$.
It should be noted that one can obtain a Kac-Rice formula for almost any level using the area and coarea formulas, Fubini's theorem and duality.
However, in applications, the interest is directed towards a fixed $y$-level.
For example, the zeros in the study of the roots of a random polynomial.
This precision leads us to a delicate study to generalize the classical theorems of inverse function and implicit function.
For this part, we have based our approach on two seminal works: on the one hand, an article by E.~Caba\~na \cite{AL121221}, published in the proceedings of the II CLAPEM conference and on the other hand, in the Lecture Notes of Mathematics by M.~Wschebor \cite{MR871689}.
    The method we use also produces the Kac-Rice formula for the upper moments of the level measure.
    \item These notes  continue  with several applications.
	We first show examples where the hypothesis can be verified, then we use the Kac-Rice formulas to obtain conditions on the finiteness of the first and second moments of the level set measure.
	The very important case of Gaussian random fields leads us to explicit calculations.
	Then, we study the number of roots of trigonometric random polynomials.
	We highlight the asymptotic behavior of the expectation and variance of the number of roots.
	We then study the application of the Kac-Rice formula to the modeling of the sea and to random gravitational lenses.
	Another topic we consider is the nodal curves of the random wave system considered by Berry and Dennis in \cite{MR1843853}.
	These curves are called dislocations in physics and correspond to the lines of darkness in the propagation of light, or the threads of silence in the propagation of sound (see \cite{MR1843853}).

	\item After these applications, we pay special attention to systems of random polynomials of several variables which are invariant under the action of the group of rotations in $\reels ^d$.
	 These polynomial systems are called Kostlan-Shub-Smale,  after the authors who first studied them and also established the properties of the set of their zeros.
	 We obtain an asymptotic formula for the variance of the measure of the set of zeros of these systems.
It is important to note that this result had already been obtained by another method by Letendre and Puchol in \cite{MR4052739}.

	\item Finally, we consider a Gaussian random field whose trajectories are not differentiable.
 If we approximate its trajectories using a convolution with a sequence approximating the Dirac delta, the approximated fields have smooth trajectories.
We then study the approximation of the local time of the original field by the length of the set of levels of the smoothed process.
We obtain the convergence in $\bL^2(\Omega)$ and also some results of convergence in law.
\end{enumerate}
\section{Hypothesis and notations}
Let $D$ be an open set of $\reels^d$ and  let $j \le d$ be a positive integer and $G: D \subset \reels^{d} \rightarrow \reels^{j}$ be a function.
\\
The function $G$ satisfies  the hypothesis ${\bH _0}$:\\
\centerline{${\bH _0}$: $G$ is continuously differentiable on $D$, that is $G \in C^1(D, \reels^j)$.}
We denote by $\nabla G(\bm{\cdot})$ its Jacobian.\\
For $y \in \reels^{j}$ we define the level set at $y$ as:
$$
	\calC_{G}(y):=\left\{x \in D : G(x)=y \right\}=G^{-1}
(y),$$
and 
$$
	\calC_{Q, G}(y):= \calC_{G}(y)\cap Q,
$$
where $Q$ is a subset of $\reels^{d} $.
\\
If $G$ satisfies ${\bH _0}$, we will denote by $D^{r}_{G}$ the following set 
$$
	D^{r}_{G}:=\left\{
		x \in  D : \nabla G(x)\text{ is of  rank } j
	\right\}.
$$%
Also $\calC_{G}^{D^{r}}(y)$ (resp$.$
$\calC_{Q,G}^{D^{r}}(y)$)
denotes the level set,
$\calC_{G}^{D^{r}}(y):=\calC_{G}(y)\cap D^{r}_{G}$ (resp$.$
$\calC_{Q,G}^{D^{r}}(y):=\calC_{Q,G}(y)\cap D^{r}_{G}$).\\
From now on,  $\sigma_d$ denotes the  Lebesgue measure on $\reels^d$.
We use the symbol $T$ for the transpose operator.
For  a set $A \subset \reels^d$, $A^{c}$ denotes its complement on $\reels^{d}$ and if $A \subset D$, $A^{c_1}$ denotes its complement on $D$.
The class of sets
$\calB(\reels^{d})$ is the Borel $\sigma$-algebra in $\reels^{d}$.
Also
$\overline{\reels }^{+}$ is the set of positive real numbers including $+\infty$, $\normp[d]{\cdot}$ denotes the Euclidean norm in $\reels^{d}$.\\
For  $x \in \reels^{d}$, $B(x, r)$ (resp$.$
$\overline{B}(x, r)$), $r >0$, is the open ball (resp$.$~closed) of center $x$ and  radius  $r$, that is 
$B(x, r):= \{z \in \reels^{d}, \normp[d]{z-x} < r\}$ (resp$.$
$\overline{B}(x, r):= \{z \in \reels^{d}, \normp[d]{z-x} \le r\}$).\\
$\naturels^{\star}=\{x \in \entiers, x >0\}$.\\
An application $f: (E, d_E)  \to (F, d_F)$ between two metric spaces is said to be $L$-Lipschitz, $L \ge 0$, if $d_F(f(x),f(y)) \le L  d_E(x,y)$, for any pair of points $x, y \in E$.\\
We also say that an application is Lipschitz if it is $L$-Lipschitz for some $L$.

In the same way, an application $f: (E, d_E)  \to (F, d_F)$ between two metric spaces, is said to be locally Lipschitz, if for each $x \in E$, there exists a neighborhood $V_{x}$ of $x$ such that the restriction of function $f$ to $V_{x}$ is Lipschitz.

$\frakL(\reels^{d}, \reels^{j})$ denotes the vector space of linear functions from $\reels^{d}$ to $\reels^{j}$ with the norm $\normp[j, d]{\cdot}$.
Also
$\frakL^2(\reels^{d}, \reels^{j})$ is the vector space of continuous symmetric linear applications from $\reels^{d}$ to $\reels^{j}$ with the norm $\normp[j,d]{\cdot}^{(s)}$.
If $B$ is a matrix, $B_{ij}$ denotes the element that appears in the $i$th row and $j$th column.\\
For $j \in \naturels^{\star}$, ${\bbS }^{j-1}$ is the boundary of the unit ball of $\reels^{j}$.\\
For any function $f$, $\supp(f)$ is its support.\\
$\bC $ is a generic constant and its value can change within a proof.
\section{Coarea formula}
The two results below are known in the literature as the coarea formula, cf. Federer \cite[pp.\,247-249]{MR0257325}  and Caba\~na \cite {AL121221}.
Our proof is based on the excellent notes of Weiz\"{a}cker  \&
Gei{\ss}ler from  Kaiserslautern University \cite {AL121221c}.
\spacebefore\begin{theo}
\label{coa1}
Let $f:\reels^ {j} \rightarrow \reels$ be a mesurable function 
and  $G: D \subset \reels^{d} \rightarrow \reels^{j}$, $j \leq d$, be a function satisfying the hypothesis  ${\bH _0}$ where $D$ is an open set.
For any Borel $B$ subset of $D$, the following formula is true:
\begin{equation}
\int_{B} f(G(x))\!\left(\det\!\left(\nabla G(x)\nabla G(x)^{T}\right)\right)^{1/2} {\ud x}
	 =\int_{\reels^{j}} f(y) \sigma_{d-j} (\calC_{B, G}^{D^{r}}(y)) \ud y \mbox{,}
\label{coarea1}
\end{equation}
provided that one of the two integrals is finite.
\end{theo}
\spacebefore
\begin{rema}
\label{finitude}
If $f$ is mesurable and positive the equality (\ref{coarea1}) is true and in this case the integrals can be infinite.
\hfill$\bullet$
\end{rema}
\spacebefore
\begin{rema}
\label{finitudeBis}
The additional assumptions that $f$ is bounded and $B$ is a compact set imply that the left-hand side integral is finite and that the formula (\ref{coarea1}) is true.
\hfill$\bullet$
\end{rema}
\spacebefore
\begin{coro}
\label{coa2}
Let $h$ be a mesurable function, $h: \reels^d \times \reels^j  \rightarrow \reels$  and $G: D \subset  \reels^{d} \rightarrow \reels^{j}$, $j \le d$, be a function satisfying the hypothesis  ${\bH _0}$ where $D$ is an open set.
For all Borel set  $B$ subset of  $D$, we have
\begin{multline}
 \int_{B} h(x, G(x)) \!\left(\det\!\left(\nabla G(x)\nabla G(x)^{T}\right)\right)^{1/2} \ud x\\
	=\int_{\reels^{j}} \left[\int_{\calC_{B, G}^{D^{r}}(y)}h(x, y) \ud \sigma_{d-j} (x)\right] \ud y,
	\label{coarea2}
\end{multline}
provided that one of the two integrals is finite.
\end{coro}
\spacebefore
\begin{rema}
\label{bornitude}
If $h$ mesurable and positive  the equality (\ref{coarea2}) is satisfied and in this case the integrals can be infinite.
\hfill$\bullet$
\end{rema}
\spacebefore
\begin{rema}
\label{bornitudeBis}
The hypotheses  that $h$ is bounded and  $B$ is compact imply that the left-hand side integral of  (\ref{coarea2}) is finite and the formula holds.
\hfill$\bullet$
\end{rema}
\spacebefore
\begin{proofarg}{Proof of Theorem \ref{coa1} and  Corollary \ref{coa2}} 

First, we will show, in the line of \cite[pp.\,60-67]{AL121221c}, the following proposition.
\spacebefore
\begin{prop}
\label{coaera formula}
Let $g:\reels^{d} \rightarrow \reels^{j}$, $j \leq d$ be a  continuously differentiable function defined on $\reels^d$.

Then for any $A \in $ $\calB$$(\reels^{d})$ we have
$$
	\int_{A} \!\left(\det\!\left(\nabla g(x)\nabla g(x)^{T}\right)\right)^{1/2} \ud x 
	=\int_{\reels^{j}} \sigma_{d-j} (\calC^{D^{r}}_{A,g}(y)) \ud y \mbox{.}
$$%
\end{prop}
\spacebefore
\begin{rema}
Proposition \ref{coaera formula} remains true if one assumes that the function $g$ is only locally Lipschitz on $\reels^d$ instead of being $C^1$ on $\reels^d$.
In this case, the function $g$ is almost surely differentiable on $\reels^d$ and the measure $\sigma_{d-j}$ which appears in the right-hand side of the previous equality is replaced by the Euclidean Hausdorff measure $\calH_{d-j}$.
We invite the reader to consult \cite[Theorem 4.12, p.\,61]{AL121221c} for more details.
\hfill$\bullet$
\end{rema}
\spacebefore
\begin{proofarg}{Proof of Proposition \ref{coaera formula}} 
As we have indicated, the proof is based on the notes \cite[pp.60-67]{AL121221c}.
First, we will prove the formula for affine functions $g$, then we will consider the formula for sets $A$ of null Lebesgue measure in $\reels^{d}$, and then we will consider sets $A$ which are subsets of $D_{g}^{r}$.
In order to accomplish our task, we need to prove some lemmas.
\spacebefore
\begin{lemm}
\label{fonctions affines}
Proposition \ref{coaera formula} is true for surjective affine functions $g$, that is
if $g(x) := a + \varphi(x)$ where $a \in \reels^{j}$ is fixed and $\varphi$ is a linear function of maximum rank $j$.
\end{lemm}
\spacebefore
\begin{proofarg}{Proof of Lemma \ref{fonctions affines}}
Without loss of generality, we can always consider $a=0$.
Indeed, on the one hand $\nabla g(\bm{\cdot})=  \nabla \varphi(\bm{\cdot})$, so for any Borel set $A$ of $\reels^{d}$ the following equality is true:
$$
	\int_{A} \!\left(\det\!\left(\nabla g(x)\nabla g(x)^{T}\right)\right)^{1/2} \ud x=
	\int_{A} (\det(\nabla \varphi(x)\nabla \varphi(x)^{T}))^{1/2} \ud x.
$$
On the other hand, because the mesure $\sigma_{j}$ is translation invariant, we have:
\begin{multline*}
 \int_{\reels^{j}} \sigma_{d-j} (\calC_{A,g}(y)) \ud y\\
	 = \int_{\reels^{j}} \sigma_{d-j} (\calC_{A,\varphi}(y-a)) \ud \sigma_{j}(y)
	=\int_{\reels^{j}} \sigma_{d-j} (\calC_{A,\varphi}(y)) \ud \sigma_j(y).
\end{multline*}
Now let $V$ be the vector subspace of  $\reels^{d}$ defined  by  $V :=\ker \varphi$.
This space is of dimension $(d-j)$ because $\varphi$ has maximal rank $j$.
We denote by $V^{\perp}$ its orthogonal complement which is of dimension $j$.

We will work with coordinate systems associated with these spaces, \ie if $x \in \reels^{d}$, we will write $x:=(z, w)$ with $z \in V^{\perp}$ and $w \in V$.
Then the Lebesgue measure $\sigma_{d}$ on $\reels^{d}$ is the product measure  $\sigma_{j} \otimes \sigma_{d-j}$.\\
Observe that $\varphi|_{V^{\perp}}$ is a one to one function since $\dim V^{\perp}=j$.
Let us denote $\Psi$ the inverse function of this restriction, \ie $\Psi:= (\varphi|_{_{V^{\perp}}})^{-1}$. 
We have $\varphi \circ \Psi = Id_{\reels^{j}}$ and $\Psi \circ \varphi|_{V^{\perp}}= Id_{V^{\perp}}$, where $Id_{\reels^{j}}$ (resp$.$
$Id_{V^{\perp}}$) is the identity function on $\reels^j$ (resp$.$ $V^{\perp}$).
Moreover, since $\varphi^{T}$ maps  $\reels^{j}$ into $V^{\perp}$ by definition of $V$, we have:
$$
	\Psi^{T} \circ \Psi \circ \varphi \circ \varphi^{T}= \Psi^{T} \circ \varphi^{T}= (\varphi \circ \Psi)^{T}= Id_{\reels^{j}},
$$%
then:
\begin{equation}
	\label{determinant}
	(\det(\varphi \circ \varphi^{T}))^{1/2}= (\det(\Psi^{T} \circ \Psi))^{-1/2}= \abs{\det(\Psi)}^{-1},
\end{equation}
the last equality is a consequence of the fact that $\Psi$ is an endomorphism of $\reels^{j}$.\\
Let $A$ be a fixed Borel set of $\reels^{d}$.
Consider the function
$$
\begin{array}{r@{~}c@{~}l}
	h: V^{\perp} & \longrightarrow & \overline{\reels}^{+}\\
	z & \longmapsto & \sigma_{d-j}\{w \in V: (z, w) \in A \}.
\end{array}
$$
Observe that for $y \in \reels^{j}$ we have
$$
	\varphi^{-1}(y) \cap A = \{(\Psi(y), w) : w \in V \} \cap A
$$%
and given that $\Psi(y) \in V^{\perp}$, 
\begin{equation}
	\label{fonctionh}
	h(\Psi(y))=\sigma_{d-j}(\varphi^{-1}(y) \cap A).
\end{equation}
So, since the Lebesgue measure $\sigma_{d}$ is the product measure $\sigma_{j} \otimes \sigma_{d-j}$, we get
\begin{equation}
	\label{sigmad}
	\sigma_{d}(A) = \int_{V^{\perp}} h(z) \ud \sigma_{j}(z).
\end{equation}
Finally, since the function $\varphi$ is a linear function and by using the equalities (\ref{determinant}), (\ref{sigmad}), the formula of change of variable for function $\Psi$ which is a $C^{1}$ function as well as $\Psi^{-1}$ as endomorphism in finite dimension and the equality (\ref{fonctionh}), we obtain
\begin{align*}
	\int_{A} (\det(\nabla \varphi(x)\nabla \varphi(x)^{T}))^{1/2} \ud x
	&=\int_{A} (\det(\varphi \circ \varphi^{T}))^{1/2} \ud x\\
	&=\sigma_{d}(A) \abs{\det(\Psi)}^{-1}\\
	&=\int_{V^{\perp}} \abs{\det(\Psi)}^{-1} h(z) \ud \sigma_{j}(z)\\
	&= \int_{\reels^{j}} h(\Psi(y))\ud y\\
	&= \int_{\reels^{j}} \sigma_{d-j}(\varphi^{-1}(y) \cap A)\ud y,
\end{align*}
this completes the proof of Lemma \ref{fonctions affines} for the affine functions.
\end{proofarg}
\spacebefore
\begin{lemm}
\label{nullset}
Let $O$ an open set of $\reels^d$, $A \subseteq O$ a Borel set of $\reels^d$ and $g: O \longrightarrow \reels^{j}$,  $j \leq d$ be a continuously differentiable Lipschitz function with Lipschitz constant, $Lip(g)$.
Thus we have
$$
	 \int_{\reels^{j}} \sigma_{d-j}
	(\calC^{D^{r}}_{A,g}(y))\ud y \le \frac{\omega_{j} \omega_{d-j}}{\omega_{d}} Lip^{\kern1pt j}(g) \sigma_{d}(A).
$$%
\end{lemm}
Above, $\omega_{d}$ denotes the volume of the unit ball of $\reels^{d}$.
\spacebefore
\begin{rema}
\label{nul}
In particular, we obtain that Proposition \ref{coaera formula} is true for Borel sets $A$ of null Lebesgue measure in $\reels^{d}$.

\hfill$\bullet$
\end{rema}
To prove Remark \ref{nul}, we need a lemma.
\spacebefore
\begin{lemm}
\label{localLipschitz}
Let $g: \reels^{d} \rightarrow \reels^{j}$, be a locally  Lipschitz function.
Then the function $g$ is Lipschitz on any compact set  $K$ of $\reels^d$.
\end{lemm}
\spacebefore
\begin{proofarg}{Proof of Lemma \ref{localLipschitz}} 
Consider $K$ a compact set of $\reels^d$ and $x \in K$.
Since $g$ is locally Lipschitz on $\reels^d$, there exists a constant $L_{x} >0$ and a radius  $r_{x} >0$, 
such that, for all $u, v \in B(x, r_{x})$,
 $\normp[j]{g(u) - g(v)} \le L_{x} \normp[d]{u - v}$.\\
Since $g$ is locally Lipschitz on $\reels^d$, it is continuous on $\reels^d$ and also on the compact set $K$.
We define $M:= \sup_{u \in K} \normp[j]{g(u)} < \infty$.\\
Since $K$ is compact, there exists $m \in \naturels^*$, such that for any $i=1,\,\dots,\, m$, there exists $x_i \in K$, satisfying 
$K \subset \cup_{i=1}^{m} B(x_i, r_{x_i}/2)$.\\
Also define $L:= \max_{i=1,\,\dots,\, m} \{L_{x_i}\}$ and $r:= \min_{i=1,\,\dots,\, m} \{r_{x_i}\}$.
Let us prove that  $g$ is a Lipschitz function on  $K$ with  Lipschitz constant $\widetilde{L}:=\max\{L, {4M}/{r}\}$.\\
Indeed, let us consider $u, v \in K$.\\
If $\normp[d]{u- v} \le \frac{r}{2}$: there exists  $i \in \{1, \dots, m\}$, such that $u \in B(x_i, r_{x_i}/2)$.
In this case  $u, v \in B(x_i, r_{x_i})$ and
$$
	\normp[j]{g(u) - g(v)} \le L_{x_i} \normp[d]{u - v} \le \widetilde{L} \normp[d]{u - v}.
$$
If $\normp[d]{u- v}> \frac{r}{2}$:
$$
	\normp[j]{g(u) - g(v)} \le 2M = \frac{4M}{r} \times \frac{r}{2} \le \frac{4M}{r} \normp[d]{u - v} \le \widetilde{L} \normp[d]{u - v}.
$$
This completes the proof of the lemma.
\end{proofarg}
\spacebefore
\begin{proofarg}{Proof of Remark \ref{nul}}
Since $g$ is $C^{1}$ on $\reels^{d}$ then it is locally Lipschitz on $\reels^{d}$.
By Lemma \ref{localLipschitz}, we know that function $g$ is then a Lipschitz function on any compact of $\reels^d$, an in particular on $K_n:=\intrff{-n}{n}^d$,  $\forall n \in \naturels^{\star}$.
Thus, let $A$ be a Borel set of $\reels^{d}$ such that $\sigma_{d}(A)=0$ and let $A_n:=\introo{-n}{n}^d  \cap  A$, $n \in \naturels^{\star}$.
By Lemma \ref{nullset}, one gets
$$
	\int_{\reels^{j}} \sigma_{d-j} (\calC^{D^{r}}_{A_n,g}(y))\ud y =0.
$$%
By taking the limit when $n\to \infty$, Beppo Levi's theorem implies that
$$
	\int_{\reels^{j}} \sigma_{d-j} (\calC^{D^{r}}_{A,g}(y))\ud y =0.
$$%
Remark \ref{nul} follows.
\end{proofarg}
\spacebefore
\begin{proofarg}{Proof of Lemma \ref{nullset}}
For $\delta >0$, we will denote by ${\calH }^{\delta}_{k}$ the Euclidean Hausdorff pre-measure which defines the $k$-dimensional Euclidean Hausdorff measure, denoted ${\calH }_{k}$, $k \in \naturels^{\star}$.
The measure ${\calH }_{k}$ coincides  with the  Lebesgue measure $\sigma_{k}$ on  $\reels^{k}$ with the Euclidean norm (cf$.$ \cite[p.\,16]{AL121221c}).\\
Let $\delta:={1}/{\ell}$, $\ell \in \naturels^{\star}$.
By definition of ${\calH }^{{1}/{\ell}}_{d}$, there exists a covering $((U_{i}^{\ell})_{i \in I_{\ell}})_{\ell \in \naturels}$ of $A$ of closed subsets of $O$ such that for all $\ell$
\begin{equation}
\label{equation1}
\abs{{U_{i}^{\ell}}} \le \frac{1}{\ell} \mbox{ and } \frac{\omega_{d}}{2^{d}}\sum_{i \in I_{\ell}} \abs{{U_{i}^{\ell}}}^{d} < {\calH }_{d}(A) + \frac{1}{\ell},
\end{equation}
where $\abs{{U}}$ denotes the Euclidean diameter  $U$, or
$$
	\abs{{U}}:= \sup_{x, y \in U}\normp[d]{x-y}.
$$%
Moreover, by definition of ${\calH }^{{1}/{\ell}}_{d-j}$, and given that $((U_{i}^{\ell})_{i \in I_{\ell}})_{\ell \in \naturels}$ covers $A$, we have
\begin{equation}
\label{equation2}
\frac{2^{d-j}}{\omega_{d-j}} {\calH }^{{1}/{\ell}}_{d-j}(g^{-1}(y)\cap A \cap D_g^r) \le \sum_{i \in I_{\ell}} \abs{{U_{i}^{\ell}}}^{d-j} \1_{g(U_{i}^{\ell})}(y):= h_{\ell}(y).
\end{equation}
Then from the inequality  (\ref{equation2}), Fatou's lemma and the fact that the measures $\sigma_{d-j}$ and ${\calH }_{d-j}$ coincide, we get the following inequalities 
\begin{eqnarray}
\label{equation3}
&&\frac{2^{d-j}}{\omega_{d-j}}  \int_{\reels^{j}}  \sigma_{d-j}
(\calC^{D^{r}}_{A,g}(y))\ud y = \frac{2^{d-j}}{\omega_{d-j}}  \int_{\reels^{j}} {\calH }_{d-j} 
(\calC^{D^{r}}_{A,g}(y))\ud y \nonumber \\ 
&&=  \frac{2^{d-j}}{\omega_{d-j}}  \int_{\reels^{j}} \lim_{\ell \to +\infty}{\calH }_{d-j}^{\ell} (\calC^{D^{r}}_{A,g}(y))\ud y \le \int_{\reels^{j}} \liminf _{\ell \to +\infty} h_{\ell}(y) \ud y \nonumber \\ 
&&\le  \liminf _{\ell \to +\infty}  \int_{\reels^{j}}  h_{\ell}(y) \ud y =
 \liminf _{\ell \to +\infty}  \int_{\reels^{j}} \sum_{i \in I_{\ell}} \abs{{U_{i}^{\ell}}}^{d-j} \1_{g(U_{i}^{\ell})}(y) \ud y \nonumber \\ 
 &&=  \liminf _{\ell \to +\infty} \sum_{i \in I_{\ell}} \abs{U_{i}^{\ell}}^{d-j} \sigma_{j}(g(U_{i}^{\ell})).
\end{eqnarray}
The idea is now to establish a relation between  $\sigma_{j}(g(U_{i}^{\ell}))$ and $\abs{U_{i}^{\ell}}$.
The isodiametric inequality for the norms  (cf$.$~\cite[p.\,14]{AL121221c}) will allow us to obtain this relation and thus to continue our proof.
Let us recall this inequality.
\spacebefore
\begin{prop}
\label{inegalité isodiamétrique}
Let $C$ be a bounded Borel set of $\reels^{j}$ then
$$
\sigma_{j} (C) \le \frac{\omega_{j}}{2^{j}} \abs{C}^{j}.
$$%
\end{prop}
The function $g$ is Lipschitz on $O$, so that the images  $g(U_{i}^{\ell})$ are bounded sets.\\
We use the  isodiametric inequality for these bounded sets and also that  $g$ is Lipschitz with Lipschitz constant  $Lip(g)$, and finally inequality (\ref{equation1}).
In this way (\ref{equation3}) gives us
\begin{align*}
\frac{2^{d-j}}{\omega_{d-j}}  \int_{\reels^{j}} \sigma_{d-j} (\calC^{D^{r}}_{A,g}(y))\ud y
&\le  \liminf _{\ell \to +\infty} \sum_{i \in I_{\ell}} \abs{U_{i}^{\ell}}^{d-j} \sigma_{j}(g(U_{i}^{\ell}))\\
&\le  \liminf _{\ell \to +\infty} \sum_{i \in I_{\ell}} \abs{U_{i}^{\ell}}^{d-j} \frac{\omega_{j}}{2^{j}} \abs{g(U_{i}^{\ell})}^{j} \\
& \le \liminf _{\ell \to +\infty} \sum_{i \in I_{\ell}} \abs{U_{i}^{\ell}}^{d-j} \frac{\omega_{j}}{2^{j}} Lip^{\kern1pt j}(g)\abs{U_{i}^{\ell}}^{j}
\\
&= \liminf _{\ell \to +\infty} \sum_{i \in I_{\ell}} \abs{U_{i}^{\ell}}^{d} \frac{\omega_{j}}{2^{j}} Lip^{\kern1pt j}(g)\\
& \le 
\liminf _{\ell \to +\infty} \frac{\omega_{j}}{2^{j}} Lip^{\kern1pt j}(g) \frac{2^{d}}{\omega_{d}}
\!\left({
{\calH }^{{1}/{\ell}}_{d}(A) + \frac{1}{\ell}}\right)\\
&= 2^{d-j} \frac{\omega_{j}}{\omega_{d}} Lip^{\kern1pt j}(g)  {\calH }_{d}(A)\\
&=
2^{d-j} \frac{\omega_{j}}{\omega_{d}} Lip^{\kern1pt j}(g)  \sigma_{d}(A).
\end{align*}
This completes the proof of Lemma \ref{nullset}.
\end{proofarg}

\spacebefore
\begin{lemm}
\label{regular} Proposition \ref{coaera formula} holds if  $A \subseteq D_{g}^{r}$.
\end{lemm}

\spacebefore
\begin{proofarg}{Proof of Lemma \ref{regular}} 

We can always assume that  $A$ is a compact.
Indeed, since $A$ is a Borel set of $\reels^{d}$, it can be written, except for a zero measure set, as a nondecreasing union of compacts.\\
Remark \ref{nul} following Lemma \ref{nullset} and Beppo Levi's theorem allow us to show  Proposition \ref{coaera formula} only in the compact case.

Let us choose an element  $x \in A$.
Consider the vector subspace of $\reels^{d}$ defined by  $V :=\ker \nabla g(x)$.
It is of dimension $(d-j)$ since $\nabla g(x)$ has maximal rank $j$.
Let $V^{\perp}$ be its orthogonal complement which is of dimension $j$.
Let us observe  that $\nabla g(x)|_{V^{\perp}}$ is one to one.\\

We will denote by  $\pi_{V}$ the orthogonal projection of $\reels^{d}$ onto $V$ and define the function  $h_{x}: \reels^{d} \to \reels^{d}$ as
$$
h_{x}(x^{\prime}):= x + \pi_{V}(x^{\prime}) + (\nabla g(x)|_{V^{\perp}})^{-1}(g(x^{\prime})-g(x)).
$$%
We want to prove that for any $\ve>0$, there exists $\delta >0$ such that if $B_{\delta}(x)$ is the set
\begin{equation}
B_{\delta}(x) :=\overline{B}(x, \delta) \cap A, \label{Bdelta}
\end{equation}
then
\begin{align}
\MoveEqLeft[2]{\left(\frac{1-\ve}{1+\ve}\right)^d(1+\ve)^j
	\left(\int_{B_{\delta}(x)}
\!\left(\det(\nabla g(x^{\prime})\nabla g(x^{\prime})^{T})\right)^{1/2} \ud x^{\prime}-\ve  
\sigma_d(B_{\delta}(x))\right)}
\nonumber\\
	& \le \int_{\reels^j}\sigma_{d-j}(\calC_{B_{\delta}(x),g}(y))\ud y \label{art0} \\ 
	& \le\!\left(\frac{1+\ve}{1-\ve}\right)^d(1-\ve)^j
	\!\left(\int_{B_{\delta}(x)}
\!\left(\det(\nabla g(x^{\prime})\nabla g(x^{\prime})^{T})\right)^{1/2} \ud x^{\prime}+\ve
\sigma_d(B_{\delta}(x))\right)\!.\nonumber
\end{align}
To do this, let's start by showing the following two things:\\
For any $\ve>0$, there exists $\delta >0$ such that 
if $x^{\prime}, x^{\prime \prime} \in B_{\delta}(x)$, we have:
\begin{equation}
\label{art1}
(1-\ve) \normp[d]{x^{\prime}-x^{\prime \prime}} \leq \normp[d]{h_{x}(x^{\prime})-h_
{x}(x^{\prime \prime})} \leq(1+\ve)\normp[d]{x^{\prime}-x^{\prime \prime}},
\end{equation}
as well as
\begin{equation}
\label{art4}
\abs{\!\left(\det\!\left(\nabla g(x)\nabla g(x)^{T}\right)\right)^{1/2}-
\!\left(\det(\nabla g(x^{\prime})\nabla g(x^{\prime})^{T})\right)^{1/2}}<\ve.
\end{equation}
The inequality (\ref{art4}) is a consequence of the fact that  $\nabla g(\bm{\cdot})$ is a continuous function defined on $\reels^d$.\\
To prove (\ref{art1}), notice that  $(\nabla g(x)|_{V^{\perp}})^{-1}$ is a finite dimensional endomorphism and therefore continuous.
Thus, we have $\normp[j, j]{(\nabla g(x)|_{V^{\perp}})^{-1}}< \infty$.\\
Furthermore, let's define $\Delta_{\delta}(x)$ as
$$
\Delta_{\delta}(x) := \sup_{x^{\prime} \neq x^{\prime \prime}, x^{\prime}, x^{\prime\prime} \in B_{\delta}(x)} \frac{\normp[j]{g(x^{\prime})-g(x^{\prime \prime})-\nabla g(x)(x^{\prime}- x^{\prime \prime})}}{\normp[d]{x^{\prime}-x^{\prime \prime}}},
$$%
remembering that $B_{\delta}(x)$ is defined by (\ref{Bdelta}).\\
Since $g$ belongs to $C^{1}$, we can write the following first order Taylor expansion  
$$
g(x^{\prime})= g(x^{\prime \prime}) + \!\left({\int_{0}^{1} \nabla g(x^{\prime \prime}+ \lambda (x^{\prime}-x^{\prime \prime})) \ud \lambda }\right) (x^{\prime}- x^{\prime \prime}),
$$%
getting
\begin{multline*}
 g(x^{\prime})-g(x^{\prime \prime}) - \nabla g(x)(x^{\prime}-x^{\prime \prime})\\
 =\!\left({\int_{0}^{1} (\nabla g(x^{\prime \prime}+ \lambda (x^{\prime}-x^{\prime \prime})) -\nabla g(x))
\ud \lambda }\right)  (x^{\prime}- x^{\prime \prime}).
\end{multline*}
This implies, since $\nabla g(\bm{\cdot})$ is continuous, that for any $\ve>0$, there exists $\delta >0$, such that
\begin{equation}
\label{delta}
\Delta_{\delta}(x) \le \ve \normp[j, j]{(\nabla g(x)|_{V^{\perp}})^{-1}}^{-1}.
\end{equation}
Let  $\ve >0$ be a fixed real number and $x^{\prime}, x^{\prime \prime} \in B_{\delta}(x)$.
Since $V = \ker \nabla g(x)$, we have
$$
\nabla g(x)|_{V^{\perp}}(\pi_{V^{\perp}}(x^{\prime}-x^{\prime \prime})) =\nabla g(x)
(x^{\prime}- x^{\prime \prime}).
$$%
This implies
\begin{align*}
\MoveEqLeft[2]{\normp[d]{x^{\prime}-x^{\prime \prime} -(h_{x}(x^{\prime})-h_{x}(x^{\prime \prime}))}}\\
&= \normp[d]{
\pi_{V^{\perp}}(x^{\prime}- x^{\prime \prime}) - (\nabla g(x)|_{V^{\perp}})^{-1}
(g(x^{\prime}) -g(x^{\prime \prime}))}\\
&= \normp[d]{
(\nabla g(x)|_{V^{\perp}})^{-1} \!\left({
\nabla g(x)(x^{\prime}-x^{\prime \prime}) -(g(x^{\prime})-g(x^{\prime \prime})}\right)}\ \\
&\le \normp[j, j]{(\nabla g(x)|_{V^{\perp}})^{-1}}  \Delta_{\delta}(x) \normp[d]{x^{\prime}-x^{\prime \prime}} \\
&\le \ve  \normp[d]{x^{\prime}-x^{\prime \prime}},
\end{align*}
the last inequality comes from (\ref{delta}).
The proof of (\ref{art1}) is thus completed.

Let $T_{x}: \reels^{d} \to \reels^{j}$ be the following affine function,
$$
T_{x}(x^{\prime}):=g(x)+
\nabla g(x)(x^{\prime}-x).
$$%
It is surjective because $\nabla g(x)$ is of maximal rank $j$.\\
Moreover, it is certain that $T_{x} \circ h_{x}=g$.\\
Indeed, given that $\pi_{V}(x^{\prime}) \in V$ and  $(\nabla g(x)|_{V^{\perp}})^{-1}(g(x^{\prime}) -g(x)) \in V^{\perp}$, we can write
\begin{align*}
\MoveEqLeft[2]{T_{x}(h_{x}(x^{\prime}))= g(x) + \nabla g(x) (h_{x}(x^{\prime}) -x)}\\
&= g(x) + \nabla g(x)(\pi_{V}(x^{\prime}))+\nabla g(x) \!\left({(\nabla g(x)|_{V^{\perp}})^{-1} (g(x^{\prime}) -g(x))}\right)\\
&= g(x) + \nabla g(x)|_{V^{\perp}} \!\left({(\nabla g(x)|_{V^{\perp}})^{-1} (g(x^{\prime}) -g(x))}\right)\\
&= g(x) + g(x^{\prime}) -g(x)\\
&= g(x^{\prime}).
\end{align*}
Furthermore, (\ref{art1}) allows us to conclude that for fixed $\ve >0$, $h_{x}$ is Lipschitz on $B_{\delta}(x)$, having a Lipschitz constant equal to $(1 + \ve)$.
Since $h_x$ is an injective function on $B_{\delta}(x)$, $h_{x}$ admits an inverse function defined on $h_{x}(B_{\delta}(x))$.
The inequality (\ref{art1}) ensures that this inverse is also Lipschitz on $h_{x}(B_{\delta}(x))$ with a Lipschitz constant equal to $(1 - \ve)^{-1}$.
\\
These two facts allow us to apply to $h_{x}$ and $h_{x}^{-1}$ the Lipschitz contraction principle which we recall below (cf$.$
\cite[p.\,18]{AL121221c}).
\spacebefore
\begin{prop}
\label{principe de contraction}
Let $E$, $F$ be two subsets of $\reels ^m$.
We assume that there exists a surjective Lipschitz function $f: E \to F$ with Lipschitz constant  $L$.
Then
$$
{\calH }_{k}(F) \le L^{k} {\calH }_{k}(E) \mbox{ \ for all \ } k \ge 0.
$$%
\end{prop}
Let us apply simultaneously this last principle on the one hand to the function  $f:=h_{x}$ for $E:= B_{\delta}(x)$, $F:=h_{x}(B_{\delta}(x))$, $L:=(1 + \ve)$ and $k:=d$, and on the other hand to 
$f:=h_{x}^{-1}$ for $E:= h_{x}(B_{\delta}(x))$, $F:=B_{\delta}(x)$, $L:=(1 - \ve)^{-1}$ and $k:=d$.
We obtain
\begin{equation}
\label{majoration boule}
(1- \ve)^{d}  \sigma_{d}(B_{\delta}(x)) \le \sigma_{d}\!\left(h_{x}(B_{\delta}(x))\right) \le (1 + \ve)^{d} \sigma_{d}(B_{\delta}(x)).
\end{equation}
Let $G$ be a set such that $G \subseteq h_{x}(B_{\delta}(x))$.
Let us again apply simultaneously the principle of contraction to the function $f:=h_{x}^{-1}$, for $E:=G$, $F:=h_{x}^{-1}(G)$, $L:=(1 - \ve)^{-1}$ and $k:=d-j$, and then to $f:=h_{x}$, for $E:=h_{x}^{-1}(G)$, $F:=G$, $L:=(1+ \ve)$ and $k:=d-j$, we get
$$
(1+\ve)^{-(d-j)} {\calH }_{d-j}(G) \le {\calH }_{d-j}(h_{x}^{-1}(G)) \le (1-\ve)^{-(d-j)} {\calH }_{d-j}(G).
$$%
In particular, if we choose $G:= T_{x}^{-1}(y) \cap h_{x}(B_{\delta}(x))$, and observe that
$$
	T_{x} \circ h_{x}=g,\quad  {\calH }_{d-j}(h_{x}^{-1}\!\left(T_{x}^{-1}(y) \cap h_{x}(B_{\delta}(x))\right))= \sigma_{d-j}(g^{-1}(y)\cap B_{\delta}(x)),
$$%
then we get
\begin{align}
\label{majoration courbe}
\MoveEqLeft[1]{(1+\ve)^{-(d-j)} \sigma_{d-j}\!\left(T_{x}^{-1}(y) \cap h_{x}(B_{\delta}(x))\right) \le \sigma_{d-j}(g^{-1}(y)\cap B_{\delta}(x))} \\
	& \le (1-\ve)^{-(d-j)} \sigma_{d-j}\!\left(T_{x}^{-1}(y) \cap h_{x}(B_{\delta}(x))\right).\nonumber 
\end{align}
We can now prove (\ref{art0}). To do so, we apply (\ref{majoration courbe}), then Lemma \ref{fonctions affines} to the surjective affine function $T_{x}$ for the Borel set $A:=h_{x}(B_{\delta}(x))$, also (\ref{majoration boule}) and finally (\ref{art4}), so we obtain
\begin{align*}
\MoveEqLeft[2]{\int_{\reels^{j}} \sigma_{d-j}\!\left(g^{-1}(y)\cap B_{\delta}(x)\right) \ud y}\\
	&\le (1-\ve)^{-(d-j)} \int_{\reels^{j}}   \sigma_{d-j}\!\left(T_{x}^{-1}(y) \cap h_{x}
	(B_{\delta}(x))\right) \ud y\\
	&= (1-\ve)^{-(d-j)}  \int_{h_{x}(B_{\delta}(x))}
	\!\left(\det(\nabla T_{x}(x^{\prime})\nabla T_{x}(x^{\prime})^{T})\right)^{1/2}\ud x^{\prime}\\
	&=(1-\ve)^{-(d-j)} \sigma_{d}\!\left(h_{x}(B_{\delta}(x))\right)
	 \!\left(\det \nabla g(x) \nabla g(x)^{T}\right)^{1/2} \\
	&\le  (1-\ve)^{-(d-j)} (1 + \ve)^{d} \sigma_{d}(B_{\delta}(x)) \!\left(\det \nabla g(x) \nabla g(x)^{T}\right)^{1/2} \\
	&\le \!\left(\frac{1+\ve}{1-\ve}\right)^d(1-\ve)^j
	\!\left(\int_{B_{\delta}(x)}
	\!\left(\det(\nabla g(x^{\prime})\nabla g(x^{\prime})^{T})\right)^{1/2}\ud x^{\prime}+\ve \sigma_d(B_{\delta}(x))\right)\!\!.
\end{align*}
In the same way
\begin{align*}
\MoveEqLeft[1]{\int_{\reels^{j}} \sigma_{d-j}\!\left(g^{-1}(y)\cap B_{\delta}(x)\right) \ud y}\\
	&\ge
	(1+\ve)^{-(d-j)} \int_{\reels^{j}}   \sigma_{d-j}(T_{x}^{-1}(y) \cap h_{x}
	(B_{\delta}(x))) \ud y\\
	&= (1+\ve)^{-(d-j)}  \int_{h_{x}(B_{\delta}(x))}
	\!\left(\det(\nabla T_{x}(x^{\prime})\nabla T_{x}(x^{\prime})^{T})\right)^{1/2}\ud x^{\prime}\\
	&=(1+\ve)^{-(d-j)} \sigma_{d}\!\left(h_{x}(B_{\delta}(x))\right)
	 \!\left(\det \nabla g(x) \nabla g(x)^{T}\right)^{1/2} \\
	&\ge  (1+\ve)^{-(d-j)} (1 - \ve)^{d} \sigma_{d}(B_{\delta}(x)) \!\left(\det \nabla g(x) \nabla g(x)^{T}\right)^{1/2} \\
	&\ge\!\left(\frac{1-\ve}{1+\ve}\right)^d(1+\ve)^j
	\!\left(\int_{B_{\delta}(x)}
	\!\left(\det(\nabla g(x^{\prime})\nabla g(x^{\prime})^{T})\right)^{1/2}\ud x^{\prime}-\ve 
	\sigma_d(B_{\delta}(x))\right).
\end{align*}
The inequality (\ref{art0}) follows.

To finish the proof of Lemma \ref{regular}, for fixed $\ve >0$, using Vitali's covering theorem (cf$.$
\cite[Theorem 1.15, p.\,14]{AL121221c}), the set $A$ can be covered, except for a set whose Lebesgue measure is zero,  by a sequence of disjoint sets of the type $B_{\delta}(x)$.
These are closed sets because $A$ is closed since it is compact.
By Remark \ref{nul}, we can forget the null measure set.
Then we take the sum on this partition.
By highlighting the fact that $\sigma_{d}(A) < \infty$, $A$ being a compact set and by using (\ref{art0}), we obtain
\begin{align*}
\MoveEqLeft[1]{\!\left(\frac{1-\ve}{1+\ve}\right)^d(1+\ve)^j\!\left(\int_{A}
\!\left(\det(\nabla g(x^{\prime})\nabla g(x^{\prime})^{T})\right)^{1/2}\ud x^{\prime}-\ve  
\sigma_d(A)\right)}\\
&\leq \int_{\reels^j}\sigma_{d-j}(\calC_{A,g}(y)) \ud y\\ \nonumber
&\leq \!\left(\frac{1+\ve}{1-\ve}\right)^d(1-\ve)^j\!\left(\int_{A}
\!\left(\det(\nabla g(x^{\prime})\nabla g(x^{\prime})^{T})\right)^{1/2}\ud x^{\prime}+\ve  
\sigma_d(A)\right).
\end{align*}
Since $\ve>0$ is sufficiently small, Proposition \ref{coaera formula} is satisfied if $A \subseteq D_{g}^{r}$, which completes the proof of Lemma \ref{regular}.
\end{proofarg}
To finish the proof of Proposition \ref{coaera formula}, we apply
Lemma \ref{regular} to the Borel set $A \cap D_{g}^{r}$.
Noting that 
$(D_{g}^{r})^{c} = \{x \in \reels^{d}, \det(\nabla g(x) \nabla g(x)^{T})=0\}$, Proposition \ref{coaera formula} follows.
\end{proofarg}
We can prove now Theorem \ref{coa1} and its  corollary.

Let $B$ be a Borel subset of $D$ and $A:= B \cap G^{-1}(I)$, where $I$ is a Borel set of $\reels^j$.\\
We consider  $\forall n \in \naturels^{\star}$, the closed sets
$$
	D_n:=\left\{x \in \reels^d, d(x, D^{c}) \ge \frac{1}{n} \right\}.
$$ 
Since $D$  is open, the sets $(D_n)_{n \in \naturels^*}$ are included in $D$.
Whitney \cite[Theorem 2]{MR1605154} allows us to extend for $n \in \naturels^*$, the $C^1$-function $G|_{D_n}$ to a function $g_n: \reels^d \to \reels^j$ still $C^1$ on $\reels^d$.\\
By applying  $\forall n \in \naturels^*$, Proposition \ref{coaera formula} to the $C^1$-function $g_n$, defined on $\reels^d$ and to the Borel set $A_n:=A \cap D_n$, we obtain
$$
	\int_{A_n}(\det(\nabla g_n(x)\nabla g_n(x)^{T}))^{1/2}\ud x
	=\int_{\reels^{j}} {\sigma}_{d-j} (\calC_{A_n,g_n}^{D^{r}}(y)) \ud y\mbox{,}
$$
and since $g_n=G$ on $D_n$ then on $A_n$,  we get
$$
\int_{B \cap D_n} \1_{I}(G(x))\!\left(\det\!\left(\nabla G(x)\nabla G(x)^{T}\right)\right)^{1/2} \ud x
=\int_{\reels^{j}} \1_{I}(y) \sigma_{d-j} \!\left(\calC_{B\cap D_n,G}^{D^{r}}(y)\right) \ud y\mbox{.}
$$
Moreover, when $n$ tends to infinity, the sets $(D_n)_{n \in \naturels^*}$ tend increasingly towards $D$.
Beppo Levi's theorem implies that Theorem \ref{coa1} holds true for functions $f$ of the form $f:=\1_I$.
By a standard approximation argument, Theorem \ref{coa1} is true for positive measurable functions.
This leads to Remark \ref{finitude} and also to Theorem \ref{coa1}.\\
Remark \ref{finitudeBis} is a consequence of the fact that if $B$ is a compact set and the function $f$ is bounded, then
\begin{multline*}
  \int_{B}
\abs{f(G(x))}\!\left(\det\!\left(\nabla G(x)\nabla G(x)^{T}\right)\right)^{1/2} \ud x \\
\le \bC  \int_{B} \!\left(\det\!\left(\nabla G(x)\nabla G(x)^{T}\right)\right)^{1/2} \ud x < \infty,
\end{multline*}
since the function $G$ is $C^1$ on $D$ and $B$ is a compact subset of  $D$.\\
Remark \ref{finitude} allows us to prove Corollary \ref{coa2}.\\
Indeed, applying this remark to the measurable and positive function $f:=\1_{I}$ and to the Borel set $B \cap A$ where $I$ (resp$.$ $A$) is some Borel set of $\reels^{j}$ (resp$.$ $\reels^{d}$), and $B$ is a Borel subset of $D$, allows to establish Corollary \ref{coa2} for functions $h$ of the form  $\1_{A \times I}$.
Still by a standard approximation argument, Remark \ref{bornitude} and Corollary \ref{coa2} follow.\\
In the same way as for Remark \ref{finitudeBis}, we obtain Remark \ref{bornitudeBis}.\\
This completes the proof of Theorem \ref{coa1} and of Corollary \ref{coa2}.
 \end{proofarg}
\section{Kac-Rice formulas for almost all levels}
In this section, $X: \Omega \times D \subset \Omega\times \reels^{d} \rightarrow \reels^{j}$ ($j \le d$) denotes a random field that belongs to
$C^{1}(D,\reels^{j})$, $Y: \Omega \times D \subset \Omega \times \reels^{d} \rightarrow \reels$
is a continuous process and $D$ is an open set of $\reels^{d}$.\\
Let $H$ be the operator
$$
\begin{array}{r@{~}c@{~}l}
	H : \frakL(\reels^{d}, \reels^{j}) & \longrightarrow & \reels^{+}\\
	A & \longmapsto &  (\det(AA^{T}))^{1/2}.
\end{array}
$$
Recall that the random set   $D^{r}_{X}$ is defined as $D^{r}_{X}=\{x \in D:
\rank\left(\nabla\!X(x)\right)=j\}$ and for  $y \in \reels^j$, the level set $\calC_{X}^{D^{r}}(y)$ is
$\calC_{X}^{D^{r}}(y)=\calC_{X}(y)\cap D^{r}_{X}$.

Let us consider the following hypotheses
\begin{itemize}
\item  ${\bH _1}$: For almost all
$x \in D$, the density of  $X(x)$, \ie ${p}_{X(x)}(\bm{\cdot})$, exists.
\item ${\bH _2}$: The function 
$$
	u \longmapsto\bbE\left[\int_{\calC_{X}^{D^{r}}(u)}      \abs{Y(x)} \ud \sigma_{d-j}(x)\right],
$$%
is a continuous function of the variable  $u$.
\item ${\bH _2^{\star}}$: The function 
$$
u \longmapsto\bbE\left[\sigma_{d-j}\!\left(\calC_{X}^{D^{r}}(u)\right)\right],
$$%
is a continuous function of the variable  $u$.
\item ${\bH _3}$: The function 
$$
	u \longmapsto \int_{D}{p}_{X(x)}(u)\bbE\!\left[\abs{Y(x)} H(\nabla\!X(x))\given{X(x)=u}\right] \ud x,
$$%
is a continuous function of the variable $u$.
\item ${\bH _3^{\star}}$: The function 
$$
	u \longmapsto \int_{D}{p}_{X(x)}(u)\bbE\!\left[H(\nabla\!X(x))\given{X(x)=u}\right] \ud x,
$$%
is a continuous function of the variable $u$.
\end{itemize}

Using the coarea formula and by duality, we will prove the following proposition.
\spacebefore
\begin{prop}
\label{dualité}
\begin{enumerate}
\item\label{itm:p2.4.1.1}  If $X$ satisfies the hypotheses  ${\bH _1}$ and (${\bH _2^\star}$ or ${\bH _3^{\star}}$), then for almost all $y \in \reels^{j}$,
\begin{equation}
	\label{chichi1}
	\bbE\left[\sigma_{d-j}(\calC_{X}^{D^{r}}(y))\right]=
		\int_D{p}_{X(x)}(y) 
			\bbE\left[H(\nabla\!X(x))\given{X(x)=y}\right] \ud x.
\end{equation}
\item\label{itm:p2.4.1.2}  If $X$ and $Y$ satisfy the hypotheses  ${\bH _1}$ and (${\bH _2}$ or ${\bH _3}$), then for almost all $y \in \reels^{j}$,
\begin{multline}
\label{chichi2}
\bbE\left[\int_{\calC_{X}^{D^{r}}(y)}Y(x) \ud \sigma_{d-j}(x)\right]\\	
	= \int_D
		{p}_{X(x)}(y) \bbE\left[Y(x)H(\nabla\!X(x))\given{X(x)=y}\right] \ud x.
\end{multline}
\end{enumerate}
\end{prop}

\spacebefore
\begin{proofarg}{Proof of Proposition \ref{dualité}} 
Let us begin by proving part \ref{itm:p2.4.1.1} of the  proposition.\\
Applying Remark  \ref{finitude} that follows Theorem \ref{coa1} to the function $G=X$ and $f=\1_{A}$ where $A$ is a Borel set of $\reels^{j}$ and to the Borel set $B=D$, we have
$$
\int_{D} \1_{X(x) \in A} H(\nabla
	X(x)) \ud x = \int_{A} \sigma_{d-j}(\calC_{X}^{D^{r}}(y)) \ud y.
$$%
By taking the expectation of each side of the equality, which is possible because both terms are positive, and by applying Beppo Levi's theorem, we obtain by using the hypothesis ${\bH _1}$
$$
	\int_{A}\bbE\left[\sigma_{d-j}(\calC_{X}^{D^{r}}(y))\right] \ud y= 
	\int_{A} \int_D{p}_{X(x)}(y) 
	\bbE\left[H(\nabla\!X(x))\given{X(x)=y}\right]\ud x \ud y.
$$%
In this step of the proof, we need a duality lemma.
\spacebefore
\begin{lemm}
\label{dualitébis}
Let $f_1, f_2: \reels^{j} \to \overline{\reels}^{+}$, be two measurable functions such that  for any bounded set $A \in \calB(\reels^{j})$, $\int_{A} f_1(y) \ud y = \int_{A} f_2(y) \ud y < \infty$, then $f_1=f_2$ $\sigma_{j}$-almost surely.
\end{lemm}
\spacebefore
\begin{proofarg}{Proof of Lemma \ref{dualitébis}} We start by proving the lemma for two measurable functions $g_1$ et $g_2$ taking values in $\overline{\reels}^{+}$ such that,
 for all $B \in \calB(\reels^{j})$, $\int_{B} g_1(y) \ud y = \int_{B} g_2(y) \ud y < \infty$.
Consider the set  $B:=\{g_2 < g_1 \}$.
Since $\int_{\reels^{j}} g_2(y) \ud y < \infty$, the hypothesis  $\int_{B} g_1(y) \ud y = \int_{B} g_2(y) \ud y$ implies that $\int_{\reels^{j}} \1_{B}(y) (g_1(y)-g_2(y)) \ud y=0$, and thus $\1_{B}(y) (g_1(y)-g_2(y)) =0$ for almost all $y \in \reels^{j}$.
Similarly, consider $B^{\prime}:=\{g_1 < g_2\}$ and conclude $\int_{\reels^{j}}\1_{B^{\prime}}(y) (g_2(y)-g_1(y))\ud y =0$ for almost all $y \in \reels^{j}$.
Finally $g_1=g_2$, $\sigma_{j}$-almost surely.\\
Now consider two functions  $f_1$ and $f_2$ satisfying the assumptions of the lemma.
Let $K$ be a compact set in $\reels^{j}$.
Since for any $B \in \calB(\reels^{j})$, $A:=B \cap K$ is a bounded Borel set of $\reels^{j}$, and if  $g_1:=f_1 \1_{K}$ and $g_2:=f_2 \1_{K}$, we have by hypothesis
$\int_{B} g_1(y) \ud y = \int_{B} g_2(y) \ud y < \infty$.
The preliminary result that we have shown implies that for any compact set $K$ of $\reels^{j}$,  $f_1 \1_{K}=f_2 \1_{K}$, $\sigma_{j}$-almost surely.\\
The proof ends by noting that, except for a set of zero measure, the set $\reels^{j}$ can be written as a non-decreasing union of compact sets and applying Beppo Levi's theorem.
\end{proofarg}
We apply here the lemma to the function $f_1(y) :=\bbE\left[\sigma_{d-j}(\calC_{X}^{D^{r}}(y))\right]$ and to the function $f_2(y):= \int_D{p}_{X(x)}(y) 
\bbE\left[H(\nabla
X(x))
\given{X(x)=y}\right] \ud x$.\\
The functions $f_1$ or $f_2$ are locally integrable by the hypothesis ${\bH _2^\star}$ or ${\bH _3^{\star}}$.
This completes the proof of the part \ref{itm:p2.4.1.1} of the proposition.

Let us prove part \ref{itm:p2.4.1.2}.
First, we assume that $X$ and $Y$ satisfy the hypotheses ${\bH _1}$ and ${\bH _2}$.
Applying Corollary \ref{coa2} to the function $G:=X$, $h(x, y):=\1_{A}(y) \times Y(x) \times \1_D(x)$ where $A$ is a bounded Borel set of $\reels^{j}$ and to the Borel set $B:=D$, we almost surely have 
\begin{equation}
	\label{Y}
	\int_{D} \1_{X(x) \in A} Y(x) H(\nabla\!X(x)) \ud x 
	= \int_{A} \left[{\int_{\calC_{X}^{D^{r}}(y)} Y(x)\ud \sigma_{d-j}(x)}\right]\ud y.
\end{equation}
Indeed, the hypothesis ${\bH _2}$ implies
\begin{equation}
	\label{YBis}
	\bbE\!\left({\int_{A} \left[{\int_{\calC_{X}^{D^{r}}(y)} \abs{Y(x)}\ud \sigma_{d-j}(x)}\right]\ud y}\right) < \infty,
\end{equation}
so almost surely $\int_{A} \left[{\int_{\calC_{X}^{D^{r}}(y)} \abs{Y(x)}\ud \sigma_{d-j}(x)
}\right]
 \ud y< \infty$, and (\ref{coarea2}) applies.
Note that (\ref{Y}) remains true if $\abs{Y}$ is substituted for $Y$.
This last observation and (\ref{YBis}) imply
$$
\bbE\left[{\int_{D} \1_{X(x) \in A} \abs{Y(x)} H(\nabla
X(x)) \ud x}\right] < \infty.
$$%
We can take the expectation both sides of (\ref{Y}) and from the hypothesis ${\bH _1}$, we obtain that for any bounded  Borel set $A$ of $\reels^{j}$
\begin{multline*}\int_{A}\bbE\left[\int_{\calC_{X}^{D^{r}}(y)}Y(x) \ud \sigma_{d-j}(x)\right] \ud y\\
=\int_{A} \!\left({ \int_D{p}_{X(x)}(y)\bbE\left[Y(x) H(\nabla\!X(x))\given{X(x)=y}\right] \ud x}\right) \ud y.
\end{multline*}
Note that the last equality is still true replacing $Y$ by $\abs{Y}$ and the corresponding integrals are finite.\\
Let us now consider
$$
	f_1(y) :=\bbE\left[\int_{\calC_{X}^{D^{r}}(y)}Y(x) \ud \sigma_{d-j}(x)\right]$$
and
$$
	f_2(y):=\int_D{p}_{X(x)}(y)	\bbE\left[Y(x) H(\nabla\!X(x))\given{X(x)=y}\right]\ud x.
$$%
A priori, the functions $f_1$ and $f_2$ do not take their values in $\overline{\reels}^{+}$.
However, a  small modification of Lemma \ref{dualitébis} can be made, by noticing that $\int_{A} \abs{f_1(y)} \ud y < \infty$ and $\int_{A} \abs{f_2(y))} \ud y < \infty$, for any bounded Borel set $A$ of $\reels^{j}$, this implies $f_1=f_2$, $\sigma_{d-j}$-almost surely.

This completes the proof of Proposition \ref{dualité} in the case where $X$ and $Y$ satisfy the hypotheses ${\bH _1}$ and ${\bH _2}$.
A similar proof can be made when $X$ and $Y$ satisfy the hypotheses ${\bH _1}$ and ${\bH _3}$.
\end{proofarg}
%
%
\chapter{Kac-Rice formula for all level}
In this section $X: \Omega \times D \subset \Omega\times \reels^{d} \rightarrow \reels^{j}$ ($j \le d$) denotes a random field that belongs to
$C^{1}(D,\reels^{j})$, $Y: \Omega \times D \subset \Omega \times \reels^{d} \rightarrow \reels$ is a continuous process and $D$ is an open set of $\reels^{d}$.
\section{Rice formula for a regular level set}
Let us specify that (\ref {chichi1}) and (\ref {chichi2}) are true for almost all $y \in \reels^{j}$.
However, in applications, these formulas are needed for all  $y$ fixed in $\reels^j$.
We will establish a theorem that gives assumptions on $X$ and $Y$ such that these formulas will hold for all $y$ in $\reels^j$.
Specifically, we will assume continuity of both members of (\ref{chichi1}) and (\ref{chichi2}), restricting ourselves to the set $D^{r}_{X}$ and proving the equality for any $y$ fixed in $\reels^j$.
Before going any further, let us state two assumptions that are useful for what follows.

These formulas are similar to the previous ones but without the absolute value in the integrand.
\begin{itemize}
\item ${\bH _{4}}$: The function 
$$
u \longmapsto\bbE\left[\int_{\calC_{X}^{D^{r}}(u)} Y(x)\ud \sigma_{d-j}(x)\right],
$$%
is a continuous function of the variable $u$.
\item ${\bH _5}$: The function 
$$
u \longmapsto \int_{D}{p}_{X(x)}(u)\bbE[Y(x)H(\nabla
X(x))\given{X(x)=u}] \ud x,
$$%
is a continuous function of the variable $u$.
\end{itemize}
\spacebefore\begin{theo}
\label{Cab1}
\begin{enumerate}
\item If $X$ satisfies  the hypotheses ${\bH _1}$, ${\bH _2^\star}$ and ${\bH _3^{\star}}$, then $\forall y \in \reels^{j}$,
\begin{equation}
	\label{chichi3}
	\bbE\left[\sigma_{d-j}(\calC_{X}^{D^{r}}(y))\right]=
	\int_D{p}_{X(x)}(y)\bbE\left[H(\nabla\!X(x))\given{X(x)=y} \right]\ud x.
\end{equation}
\item If $X$ and $Y$ satisfy the hypotheses ${\bH _1}$, (${\bH _2}$ or ${\bH _3}$), ${\bH _{4}}$ and ${\bH _{5}}$, then $\forall y \in \reels^{j}$,
\begin{multline}
\label{chichi4}
\bbE\left[\int_{\calC_{X}^{D^{r}}(y)}Y(x)\ud \sigma_{d-j}(x)\right]\\
= \int_D{p}_{X(x)}(y) 
\bbE\left[Y(x)H(\nabla
X(x))\given{X(x)=y} \right]\ud x.
\end{multline}
\end{enumerate}
\end{theo}
\spacebefore
\begin{rema}
\label{chichi6}
 If $X$ and  $Y$  satisfy the hypotheses ${\bH _1}$, ${\bH _2}$ and ${\bH _3}$, then $\forall y \in \reels^{j}$ it holds that
\begin{multline*}
	\bbE\left[\int_{\calC_{X}^{D^{r}}(y)}\abs{Y(x)} \ud \sigma_{d-j}(x)\right]\\
	= \int_D{p}_{X(x)}(y)\bbE\left[\abs{Y(x)}H(\nabla\!X(x))\given{X(x)=y} \right]\ud x.
\end{multline*}
\hfill$\bullet$
\end{rema}
\spacebefore
\begin{proofarg}{Proof of  Theorem \ref{Cab1} and of Remark \ref{chichi6}}
Let us start by proving (\ref{chichi3}).\\
Since $X$ satisfies the hypotheses ${\bH _1}$, ${\bH _2^\star}$ and ${\bH _3^{\star}}$, as a consequence of part \ref{itm:p2.4.1.1} of Proposition \ref{dualité}, we know that for almost all $y \in \reels^{j}$,
$$
	\bbE\left[\sigma_{d-j}(\calC_{X}^{D^{r}}(y))\right]=
	\int_D{p}_{X(x)}(y)\bbE\left[H(\nabla\!X(x))\given{X(x)=y}\right] \ud x.
$$%
The hypotheses ${\bH _2^\star}$ and ${\bH _3^{\star}}$ imply the continuity of each member of the equality and in consequence their equality $\forall y \in \reels^{j}$.

By reasoning in a similar way, we can prove (\ref{chichi4}).
This completes the proof of Theorem \ref{Cab1}.\\
Remark \ref{chichi6} comes from the fact that the hypotheses ${\bH _{4}}$ and ${\bH _{5}}$ become the hypotheses ${\bH _2}$ and ${\bH _3}$, replacing  $Y$ by $\abs{Y}$.
\end{proofarg}
\subsection{Checking the hypotheses}
\label{validation}
In the previous section, we proved (\ref{chichi4}) by assuming mainly the continuity of its two members.
Our goal in what follows is to specify a large class of processes $X$ and $Y$ that satisfies the hypotheses ${\bH _{i}}$, $i=1,\dots,5$.

We first consider the hypotheses ${\bH _2}$ and ${\bH _{4}}$.
This leads us to prove Theorem \ref{continueBis} which is needed to prove Proposition \ref{H1etH4} which follows.
We must emphasize that the proof we are about to give is deeply inspired by Caba\~na \cite{AL121221}.

For a while, the functions $X$ and $Y$ will be assumed to be deterministic, \ie  they are not random functions.
\spacebefore\begin{theo}
\label{continueBis}
Let $X:  D \subset  \reels^{d} \rightarrow \reels^{j}$ ($j \le d$) be a function belonging to
$C^{1}(D,\reels^{j}) \such \nabla\!X$ is Lipschitz on $D$ which is an open and convex set of $\reels^{d}$.
Let $D_1$ be an open and bounded subset of $D$ and $Y: D_1 \subset \reels^{d} \rightarrow \reels$ a continuous function such that $\supp(Y) \subset D_{X|_{D_1}}^{r}$.
Then the function 
$$
	y \longmapsto \int_{\calC_{D_1, X}^{D^{r}}(y)} Y(x) \ud \sigma_{d-j}(x)
$$%
is continuous with respect to the variable $y$.
\end{theo}
\spacebefore
\begin{proofarg}{Proof of Theorem \ref{continueBis}} 
Just like Caba\~na, we will define an atlas of  $\calC_X^{D^r}(y)$.
 Consider  $x_0\in D_X^{r}$ fixed.
Thus $\nabla\!X(x_{0})$ has rank $j$.

If $A_d:=\{1,2,\ldots,d\}$, there exists $\lambda:=(\ell_1, \ell_2, \dots, \ell_{j}) \in A_d^j$, $\ell_1 < \ell_2 < \dots <\ell_{j}$ such that
$$
	J_{X}^{(\lambda)}(x_0)
		:=\det \!\left({\frac{\partial (X_1,\ldots,X_j)}{\partial(x_{\ell_1}, x_{\ell_2}, \ldots,x_{\ell_{j}})}(x_0)}\right)
		\neq 0.
$$%
We define $\lambda^{c}$ the complementary index in $A_d$,
\ie  $\lambda^{c}:= (i_1, i_2, \dots, i_{d-j}) \in A_{d}^{d-j}$, and $i_1 < i_2< \dots < i_{d-j}$.
If $(e_1, e_2, \dots,e_{d})$ denotes the canonical basis of $\reels^{d}$, let $V_{\lambda}:= \vect(e_{i_1}, e_{i_2}, \dots, e_{i_{d-j}})$ and  $V_{\lambda}^{\perp}$ be the corresponding orthogonal subspace, \ie 
$V_{\lambda}^{\perp}:=\vect(e_{\ell_1}, e_{\ell_2}, \dots, e_{\ell_{j}})$.
With these notations, if 
$
x:=(x_1, x_2, \dots, x_{d})= \sum\limits_{i=1}^{d} x_{i}e_{i} \in \reels^{d},
$
 we denote 
$\widehat{x}_{\lambda} := (x_{i_1}, x_{i_2}, \dots, x_{i_{d-j}})$.

Consider the function $f_{\lambda}$ defined from $D \subset \reels^{d}$ into $\reels^{d} \such x \mapsto f_{\lambda}(x):=\pi_{V_{\lambda}}(x) + \sum\limits_{k=1}^{j} (X_{k}(x) - y_{k})e_{\ell_{k}}$, where $X_{k}$ (resp$.$
$y_{k}$) denotes  the components of $X$ (resp$.$
$y$), $k=1, \dots, j$, and where $\pi_{V_{\lambda}}$ represents the projector on $V_{\lambda}$.\\
The Jacobian  $J_{f_{\lambda}}(x_{0})$ of this transformation evaluated at $x_{0}$ is $J_{f_{\lambda}}(x_{0})= \abs{J_{X}^{{\lambda}}(x_0)} \neq 0$.
By the inverse function theorem, there exists an open neighborhood $U_{x_0}^{\lambda}$ of $x_0$ included in $D$, such that $f_{\lambda}(U_{x_0}^{\lambda})$ is still an open set of $\reels^d$ and such that the restriction  $f_{\lambda}|_{U_{x_0}^{\lambda}}$ has an inverse $h_{\lambda}$ belonging to $C^1$ defined from $f_{\lambda}(U_{x_0}^{\lambda})$ onto $U_{x_0}^{\lambda}$.

Let us define the set $R_{x_0}^{\lambda}$, by
$$
R_{x_0}^{\lambda}:=\left\{(x_{i_1}, x_{i_2}, \dots, x_{i_{d-j}}) \in \reels^{d-j}: \sum\limits_{k=1}^{d-j} x_{i_{k}} e_{i_{k}} \in f_{\lambda}(U_{x_0}^{\lambda})\right\}.
$$%

Since $f_{\lambda}(U_{x_0}^{\lambda})$ is an open set, the set $R_{x_0}^{\lambda}$ is also an open set of $\reels^{d-j}$.
Let us denote $h_{\lambda}:=(h^{\lambda}_1, h^{\lambda}_2, \dots, h^{\lambda}_{d})$.
We have the following sequence of equivalences
\begin{align*}
	\MoveEqLeft[2]{\!\left({x \in U_{x_0}^{\lambda}, X(x)=y}\right) \iff \!\left({x\in U_{x_0}^{\lambda},  f_{\lambda}(x)= \pi_{V_{\lambda}}(x)}\right)}\\
	 &\iff \!\left({\pi_{V_{\lambda}}(x) \in f_{\lambda}(U_{x_0}^{\lambda}), x = h_{\lambda}(\pi_{V_{\lambda}}(x))
}\right)\\
  &\iff \!\left({ x = \pi_{V_{\lambda}}(x) + \sum_{k=1}^{j} h^{\lambda}_{\ell_{k}}(\pi_{V_{\lambda}}(x))  e_{\ell_{k}}, \widehat{x}_{\lambda} \in R_{x_0}^{\lambda}}\right)\\
   & \iff \!\left({x= \sum_{k=1}^{d-j} x_{i_{k}} e_{i_{k}} + \sum_{k=1}^{j} h^{\lambda}_{\ell_{k}}\left(\sum_{k=1}^{d-j} x_{i_{k}} e_{i_{k}}\right) e_{\ell_{k}}, \widehat{x}_{\lambda} \in R_{x_0}^{\lambda} }\right) \\
   & \iff  \!\left({x =  \overrightarrow{\alpha_{\lambda, x_{0}}}(\widehat{x}_{\lambda}), \widehat{x}_{\lambda} \in R_{x_0}^{\lambda}}\right),
\end{align*}
where we defined $\overrightarrow{\alpha_{\lambda, x_{0}}} : 
 R_{x_0}^{\lambda} \subset \reels^{d-j} \to \reels^d$ by 
\begin{equation}
 	\label{parametrisationlocale}
	 \overrightarrow{\alpha_{\lambda, x_{0}}}(x_{i_1}, x_{i_2}, \dots, x_{i_{d-j}})
	 := \sum\limits_{k=1}^{d-j} x_{i_{k}} e_{i_{k}}
	 	+ \sum\limits_{k=1}^{j} h^{\lambda}_{\ell_{k}}\!\left(\sum\limits_{k=1}^{d-j} x_{i_{k}} e_{i_{k}}\right) e_{\ell_{k}}.
 \end{equation}
 This provides us with a local parametrization of the level set  $\calC_X^{D^r}(y)$, defined by $\overrightarrow{\alpha_{\lambda}}$.
Moreover, such a function belongs to $C^{1}$ defined over $R_{x_0}^{\lambda}$.

  \spacebefore
\begin{rema} 
  \label{fonctions implicites}
Furthermore, as a bonus, we get that $\calC_X^{D^r}(y)$ is a differentiable manifold of dimension $(d-j)$.

\hfill$\bullet$
\end{rema}
We now decompose $D_{X}^{r}$ into the following form:
$$
	D_{X}^{r}=\bigcup_{\lambda \in B_{j}} \Gamma(\lambda),
$$%
where $B_{j}:=\{\lambda=(\ell_1,\ell_2,\ldots,\ell_{j}), \ell_k\in A_d, \ell_1<\ell_2<\ldots<\ell_{j} \}$ and
$\Gamma(\lambda):= \{x \in D, J_{X}^{(\lambda)}(x) \neq 0 \}$.
\spacebefore
\begin{rema}
\label{inversibilite}
For all  $\lambda \in B_{j}$,
$\Gamma(\lambda)=  \{x \in D,  \nabla\!X(x)|_{V_{\lambda}^{\perp}} \mbox{~is invertible}\}$
\hfill$\bullet$
\end{rema}
\spacebefore
\begin{rema}
\label{inversibilite Bis}
Let us mention that if  $x_{0} \in \calC_{X}^{D^{r}}(y) \cap \Gamma(\lambda)$, for $\lambda \in B_{j}$, then $\overrightarrow{\alpha_{\lambda, x_{0}}}(\widehat{x}_{0,\lambda} )=x_0$.
\hfill$\bullet$
\end{rema}
\spacebefore
\begin{rema}
We could have proved a less general result than the one established in this theorem.
More exactly, we could have proved this theorem under the following weaker hypothesis: $Y : D_X^r \subset \reels^d \to \reels$ is a continuous function such that $\supp(Y) \subset \Gamma(\lambda)$, for any $\lambda \in B_{j}$.
Indeed, it is this condition that we finally need in the proof of Propositions \ref{H1etH4} and \ref{CabBis} given later.
However, the result of the above theorem seemed interesting to us in itself because we did not find it in the literature.
Moreover, this general result has the advantage that in the proof we expose in a neat way a partition of unity of $D_X^r$.
This construction will allow us later in the proof of Proposition \ref{CabBis} to decompose the function $Y$ on this partition and thus to find ourselves in the case where the function has its support included in $\Gamma(\lambda)$.
\hfill$\bullet$
\end{rema}
Let us prove Theorem \ref{continueBis} in the case where $D=D_1$.
That is, when $D$ is an open, convex, bounded set of $\reels^d$, $X : D \subset \reels^{d} \rightarrow \reels^{j}$ ($j \le d$) is a function in $C^{1}(D,\reels^{j}) \such \nabla\!X$ is Lipschitz and $Y : D \subset \reels^{d} \rightarrow \reels$ is a continuous function such that $\supp(Y) \subset D_X^r$.
It suffices to prove the theorem in the case where $Y : D_X^r \subset \reels^d \to \reels$ is a continuous function with $\supp(Y) \subset D_X^r$.

Let $y$ be fixed in $\reels^j$.
We assume $\supp(Y) \subset \Gamma(\lambda)$, $\lambda \in B_{j}$.
We will define the integral of $Y$ on the level set $\calC_X^{D^r}(y)$, \ie we will give a meaning to  $\int_{\calC_X^{D^{r}}(y)}Y(x) \ud\sigma_{d-j}(x).$
 
Consider $\displaystyle x_{0} \in \supp(Y)  \subset \Gamma(\lambda)$.
For the previous facts, there exists an  open neighborhood $U_{x_0}^{\lambda}$ of $x_0$, such that
$$
\!\left({x \in U_{x_0}^{\lambda} \cap \calC_X^{D^r}(y)}\right) \iff \!\left({x =  \overrightarrow{\alpha_{\lambda, x_{0}}}(\widehat{x}_{\lambda}), \widehat{x}_{\lambda} \in R_{x_0}^{\lambda}}\right).
$$%
Since  $U_{x_0}^{\lambda}$ is an open set, one can choose a radius $r_{x_{0}}^{\lambda} >0$ such that
the closed ball of $\reels^d$ with center $x_{0}$ and radius $r_{x_{0}}^{\lambda}$ is contained in $U_{x_0}^{\lambda}$, so $\overline{B}(x_{0}, r_{x_{0}}^{\lambda}) \subset U_{x_0}^{\lambda}$.\\
Since  $\supp(Y) $ is a compact set of $\reels^d$, we can cover $\supp(Y) $ by a finite number of balls, \ie  $\supp(Y) \subset \cup_{i=1}^{m}
B(x_{i}, r_{x_{i}}^{\lambda})$ such that for all  $i=1, \dots, m$, we still have:
$ \overline{B}(x_{i}, r_{x_{i}}^{\lambda}) \subset U_{x_i}^{\lambda} $.\\
We will construct a partition of unity for $\supp(Y)$ which is a compact manifold and denote it by $\{\pi_1, \dots, \pi_{m} \}$.\\
As in Wendell Fleming's book \cite{MR0422527}, we define the real-valued function $h$ of $C^{\infty}$ by
$$
h(x):= 
	\begin{cases}
		\ds\exp\!\left(\frac{-1}{1-x^{2}}\right), & \abs{x}< 1; \\
		0,&  \abs{x} \ge 1.
	\end{cases}
$$%
For $i=1, \dots, m$ and for $x \in \supp(Y)  $, let
$\Psi_{i}(x):= h\!\left({{\normp[d]{x-x_{i}}}/{r_{x_{i}}^{\lambda}}}\right)$.\\
Since $\supp(Y)  \subset \cup_{i=1}^{m} B(x_{i}, r_{x_{i}}^{\lambda})$  by construction, we have that $\forall x \in \supp(Y)$, $\sum\limits_{i=1}^{m} \Psi_{i}(x) >0$.\\
For $i=1, \dots, m$ and for $x \in \supp(Y)  $, let us define $\pi_{i}(x):= {\Psi_{i}(x)}/{\sum_{i=1}^{m}\Psi_{i}(x)}$.\\
\label{partition unite}

The functions  $\{\pi_1, \dots, \pi_{m} \}$ define a partition of unity for $\supp(Y)  $ because
\begin{enumerate}
	\item $\pi_{i}$ is $C^{\infty}$ over $\supp(Y) $, $\pi_{i} \ge 0$, $i=1, \dots, m$;
	\item $\supp(\pi_{i})= \supp(\Psi_{i}) \subset \supp(Y) \cap \overline{B}(x_{i}, r_{x_{i}}^{\lambda}) \subset U_{x_i}^{\lambda} $;
	\item\label{item.p24.3} $\sum_{i=1}^{m} \pi_{i}(x) =1$, $x \in \supp(Y)$.
	\label{propriete3}
\end{enumerate}
The integral of $Y$ on the level set $y$ can be defined by
\begin{multline}
{\int_{\calC_X^{D^{r}}(y)}Y(x)\ud \sigma_{d-j}(x)= \sum_{i=1}^{m} \int_{U_{x_i}^{\lambda} \cap \calC_X^{D^r}(y)} \pi_{i}(x) Y(x) \ud \sigma_{d-j}(x)}\\
=\sum_{i=1}^{m} \int_{R_{x_{i}}^{\lambda}} \pi_{i}\!\left(\overrightarrow{\alpha_{\lambda, x_{i}}}(\widehat{x}_{\lambda})\right)
Y\!\left(\overrightarrow{\alpha_{\lambda, x_{i}}}(\widehat{x}_{\lambda})\right)
	\left(\det \!\left(\nabla \overrightarrow{\alpha_{\lambda, x_{i}}}(\widehat{x}_{\lambda})\nabla \overrightarrow{\alpha_{\lambda, x_{i}}}(\widehat{x}_{\lambda})^T\right)\right)^{1/2} \ud\widehat{x}_{\lambda}.\label{definition integrale}
\end{multline}

After defining the integral, we need to prove that it is a continuous function of the level. As the level varies, we have to modify the procedure.
We continue by following the approach of Caba\~na.\\
Consider $x_0\in \Gamma(\lambda) \cap \calC_X^{D^r}(y)$ fixed.
We then construct $U_{x_0}^{\lambda}$ and $R_{x_{0}}^{\lambda}$.
Let us define the function 
$$
G: R_{x_{0}}^{\lambda} \times \reels^{j} \times \widetilde{W}_{x_{0}}^{\lambda} \subset \reels^{d-j} \times \reels^{j} \times \reels^{j}\to \reels^{j},
$$%
 by
$$
G(\widehat{x}_\lambda,\delta,\gamma)
:=X(\overrightarrow{\alpha_{\lambda, x_{0}}}(\widehat{x}_\lambda)+\sum_{k=1}^{j} \gamma_{k} e_{\ell_{k}})-y-\delta,
$$%
where
$\gamma:= (\gamma_1, \dots, \gamma_{j})$,
$\widehat{x}_{\lambda} := (x_{i_1}, x_{i_2}, \dots, x_{i_{d-j}})$
and
$\widetilde{W}_{x_{0}}^{\lambda}$
is a neighborhood of zero in $\reels^j$ such that for any 
$\widehat{x}_\lambda \in R_{x_0}^{\lambda}$,
$\overrightarrow{\alpha_{\lambda, x_{0}}}(\widehat{x}_\lambda)+\sum_{k=1}^{j} \gamma_{k} e_{\ell_{k}} \in D$
and 
$\nabla\!X(\overrightarrow{\alpha_{\lambda, x_{0}}}(\widehat{x}_\lambda)+\sum_{k=1}^{j} \gamma_{k} e_{\ell_{k}})|_{V_{\lambda}^{\perp}}$
is invertible, which remains possible since $\Gamma(\lambda)$ is an open set of $D$ and $x_{0} \in \Gamma(\lambda)$.

Since $X$ is $C^{1}$ and $\overrightarrow{\alpha_{\lambda, x_{0}}}$ is $C^{1}$ on $R_{x_0}^{\lambda}$, then $G$ is $C^{1}$ on $R_{x_0}^{\lambda} \times \reels^{j} \times  \widetilde{W}_{x_{0}}^{\lambda}$.\\
Moreover, by Remark \ref{inversibilite Bis}, we have
$G(\widehat{x}_{0,\lambda},0,0)=X(\overrightarrow{\alpha_{\lambda, x_{0}}}(\widehat{x}_{0,\lambda})) -y=X(x_{0})-y=\vec{0}_{\reels^{j}}.$\\
More,
\begin{eqnarray}
\label{derivee}
\frac{\partial G}{\partial \gamma}(\widehat{x}_\lambda,\delta,\gamma)
=\left.\nabla\!X\!\left(\overrightarrow{\alpha_{\lambda, x_{0}}}(\widehat{x}_\lambda)+\sum_{k=1}^{j} \gamma_{k} e_{\ell_{k}}\right)\right|_{V_{\lambda}^{\perp}},
\end{eqnarray}
this implies
$$\frac{\partial G}{\partial \gamma}(\widehat{x}_{0,\lambda},0,0)=\nabla\!X(x_0)|_{V_{\lambda}^{\perp}},$$%
which is invertible by Remark \ref{inversibilite}.\\
 The implicit function theorem can be applied.
Thus there exists three neighborhoods, first  $V_{x_0}^{\lambda}\subset R_{x_0}^{\lambda}$ which we can choose equal to $R_{x_0}^{\lambda}$, and two other neighborhoods of zero in $\reels^{j}$ which we denote by $W_{x_{0}}^{\lambda}$ and $\widetilde{W}_{x_{0}}^{\lambda}$ (we can also choose this neighborhood equal to $\widetilde{W}_{x_{0}}^{\lambda}$) and a function $\gamma_{\lambda, x_{0}}$ of class $C^{1}$ defined from $R_{x_0}^{\lambda} \times W_{x_{0}}^{\lambda}$ onto $\widetilde{W}_{x_{0}}^{\lambda}$, that is $\gamma_{\lambda, x_{0}}: R_{x_0}^{\lambda} \times W_{x_{0}}^{\lambda} \subset \reels^{d-j} \times \reels^{j}\to \widetilde{W}_{x_{0}}^{\lambda} \subset\reels^{j}$ such that 
 \begin{enumerate}
	 \item $\gamma_{\lambda, x_{0}}(\widehat{x}_{0,\lambda}, 0)=\vec{0}_{\reels^{j}}$
	 \item $\forall (\widehat{x}_\lambda,\delta) \in R_{x_0}^{\lambda} \times W_{x_{0}}^{\lambda}$, $G(\widehat{x}_\lambda,\delta, \gamma_{\lambda, x_{0}}(\widehat{x}_\lambda,\delta))=\vec{0}_{\reels^{j}}$
	 \item \label{item:p26,3}$\forall (\widehat{x}_\lambda,\delta, \gamma) \in R_{x_0}^{\lambda} \times W_{x_{0}}^{\lambda} \times \widetilde{W}_{x_{0}}^{\lambda} $,
	 $$
	{\!\left({G(\widehat{x}_\lambda,\delta, \gamma)=\vec{0}_{\reels^{j}} \Rightarrow \gamma= \gamma_{\lambda, x_{0}}(\widehat{x}_\lambda,\delta)}\right)}.
	 $$
 \end{enumerate}
Moreover, by differenciating the expression  
$G(\widehat{x}_\lambda,\delta, \gamma_{\lambda, x_{0}}(\widehat{x}_\lambda,\delta))=\vec{0}_{\reels^{j}}$
with respect to
$\widehat{x}_\lambda$,  
we obtain for all $(\widehat{x}_\lambda,\delta) \in R_{x_0}^{\lambda} \times W_{x_{0}}^{\lambda}$,
\begin{eqnarray*}
O_{j, d-j}= \frac{\partial G}{\partial \widehat{x}_\lambda}\left(\widehat{x}_\lambda,\delta, \gamma_{\lambda, x_{0}}(\widehat{x}_\lambda,\delta)\right) + \frac{\partial G}{\partial \gamma}\left(\widehat{x}_\lambda,\delta, \gamma_{\lambda, x_{0}}(\widehat{x}_\lambda,\delta)\right) \times \frac{\partial \gamma_{\lambda, x_{0}}}{\partial \widehat{x}_\lambda}\left(\widehat{x}_\lambda,\delta\right),
 \end{eqnarray*}
where  $O_{j, d-j}$ is the null matrix with  $j$ rows and  $d-j$ columns.\\
 From(\ref{derivee}) and since
$$
\frac{\partial G}{\partial \widehat{x}_\lambda}(\widehat{x}_\lambda,\delta, \gamma_{\lambda, x_{0}}(\widehat{x}_\lambda,\delta))
= \nabla\!X\!\left(\overrightarrow{\alpha_{\lambda, x_{0}}}(\widehat{x}_\lambda)+\sum\limits_{k=1}^{j} \gamma_{k, \lambda, x_{0}}(\widehat{x}_\lambda,\delta) e_{\ell_{k}}\right) \times  \nabla \overrightarrow{\alpha_{\lambda, x_{0}}}(\widehat{x}_\lambda),
$$%
we finally obtain 
 \begin{multline}
 \label{parametrization}
\frac{\partial \gamma_{\lambda, x_{0}}}{\partial \widehat{x}_\lambda}(\widehat{x}_\lambda,\delta)
=
-\!\left[\nabla\!X\!\left(\overrightarrow{\alpha_{\lambda, x_{0}}}(\widehat{x}_\lambda)+\sum_{k=1}^{j} \gamma_{k, \lambda, x_{0}}(\widehat{x}_\lambda,\delta) e_{\ell_{k}}\right)|_{V_{\lambda}^{\perp}}\right]^{-1} \\ 
 \times\nabla\!X\!\left(\overrightarrow{\alpha_{\lambda, x_{0}}}(\widehat{x}_\lambda)+\sum_{k=1}^{j} \gamma_{k, \lambda, x_{0}}(\widehat{x}_\lambda,\delta) e_{\ell_{k}}\right) \times \nabla \overrightarrow{\alpha_{\lambda, x_{0}}}(\widehat{x}_\lambda).
\end{multline}
Define
\begin{equation}
	\label{parametrisation}
	\overrightarrow{\alpha_{\lambda, x_{0}, \delta}}(\widehat{x}_{\lambda}):=\overrightarrow{\alpha_{\lambda, x_{0}}}(\widehat{x}_\lambda)+\sum\limits_{k=1}^{j} \gamma_{k, \lambda, x_{0}}(\widehat{x}_\lambda,\delta) e_{\ell_{k}},
\end{equation}
this function is a local parametrization  of the level set  $\calC_X^{D^{r}}(y+\delta)$.
With this new parametrization we can write (\ref{parametrization}) as
\begin{multline}
\frac{\partial \gamma_{\lambda, x_{0}}}{\partial \widehat{x}_\lambda}(\widehat{x}_\lambda,\delta)
	= -\!\left[\nabla\!X(\overrightarrow{\alpha_{\lambda, x_{0}, \delta}}(\widehat{x}_{\lambda}))|_{V_{\lambda}^{\perp}}\right]^{-1}\\ \times \nabla\!X(\overrightarrow{\alpha_{\lambda, x_{0},\delta}}(\widehat{x}_{\lambda})) \times \nabla \overrightarrow{\alpha_{\lambda, x_{0}}}(\widehat{x}_\lambda).\label{deriveegammalambda}
\end{multline}
In a similar way, by differentiating the equality  $G(\widehat{x}_\lambda,\delta, \gamma_{\lambda, x_{0}}(\widehat{x}_\lambda,\delta))=\vec{0}_{\reels^{j}}$ with respect to $\delta$, this time, we obtain the equality 
\begin{equation}	
	\frac{\partial\gamma_{\lambda, x_{0}}}{\partial \delta}(\widehat{x}_\lambda,\delta)
	=\!\left[\nabla\!X(\overrightarrow{\alpha_{\lambda, x_{0},\delta}}(\widehat{x}_{\lambda}))|_{V_{\lambda}^{\perp}}\right]^{-1}.
	\label{deriveebis}
\end{equation}
We should point out that if $\widehat{x}_\lambda \in R_{x_{0}}^{\lambda}$ then $G(\widehat{x}_\lambda, 0, 0)= X(\overrightarrow{\alpha_{\lambda, x_{0}}}(\widehat{x}_\lambda)) - y=\vec{0}_{\reels^{j}}$.
By the previous point \ref{item:p26,3}, we get $\gamma_{\lambda, x_{0}}(\widehat{x}_\lambda,0)=0$, and since $\gamma_{\lambda, x_{0}}$ is a continuous function on $R_{x_{0}}^{\lambda} \times W_{x_{0}}^{\lambda}$, first we obtain 
$$
\lim_{\delta\to0}\overrightarrow{\alpha_{\lambda, x_{0},\delta}}(\widehat{x}_{\lambda})=\overrightarrow{\alpha_{\lambda, x_{0}}}(\widehat{x}_\lambda),
$$%
the convergence being uniform on $\overline{R}_{x_{0}}^{\lambda}$.
Indeed, since
$\gamma_{\lambda, x_{0}}(\widehat{x}_\lambda,0)=0$ and using the mean value theorem and (\ref{deriveebis}), we have
\begin{align*}
\MoveEqLeft[1]{\normp[d]{\overrightarrow{\alpha_{\lambda, x_{0}, \delta}}(\widehat{x}_{\lambda})-\overrightarrow{\alpha_{\lambda, x_{0}}}(\widehat{x}_\lambda)}}\\
&= \normp[d]{\sum\limits_{k=1}^{j} \!\left(\sum_{i=1}^{j} \frac{\partial \gamma_{k, \lambda, x_{0}}}{\partial \delta_{i}} (\widehat{x}_\lambda, \theta_{k}\delta) \delta_{i}\right) e_{\ell_{k}} }\\
&=\normp[d]{\sum\limits_{k=1}^{j}   ((\nabla\!X(\overrightarrow{\alpha_{\lambda, x_{0}, \theta_{k}\delta}}(\widehat{x}_{\lambda}))|_{V_{\lambda}^{\perp}})^{-1}(\delta))_{k, 1} e_{\ell_{k}}}\\
&\le \sqrt{j} \sup_{z \in K} \normp[d, j]{(\nabla\!X(z)|_{V_{\lambda}^{\perp}})^{-1}} \times \normp[j]{\delta},
\end{align*}
where $0<\theta_{k}< 1$, for $k=1, \dots, j$ and $K$ is a compact set of $\reels^{d}$ defined by  $K:= \overrightarrow{\alpha_{\lambda, x_{0}}}(\overline{R_{x_{0}}^{\lambda}}) +\sum_{i=1}^{j} \gamma_{\ell_{k, \lambda, x_{0}}}(\overline{R_{x_{0}}^{\lambda}} \times \overline{W_{x_{0}}^{\lambda}}) e_{\ell_{k}} $.
Finally, let us recall that $\normp[d, j]{\smash{(\nabla\!X(z)|_{V_{\lambda}^{\perp}})^{-1}}}$ remains bounded on this compact set.
To show this, just make the open sets  $V_{x_0}^{\lambda}$ and $W_{x_{0}}^{\lambda}$ smaller.
That is,  choose $V_{x_0}^{\lambda}$ such that
$\overline{V_{x_0}^{\lambda}}\subset R_{x_0}^{\lambda}$  and an open set  $F_{x_{0}}^{\lambda}$ containing $0$ on $\reels^{j}$ such that  $\overline{F_{x_{0}}^{\lambda}} \subset W_{x_{0}}^{\lambda}$.\\
Second, let us prove that $\nabla \overrightarrow{\alpha_{\lambda, x_{0}, \delta}}$ converges uniformly to $\nabla \overrightarrow{\alpha_{\lambda, x_{0}}}$ on $\overline{R_{x_{0}}^{\lambda}}$.
Given that $\forall \widehat{x}_\lambda \in R_{x_{0}}^{\lambda}$, $X(\overrightarrow{\alpha_{\lambda, x_{0}}}(\widehat{x}_\lambda)) = y$, we get 
$$
	\nabla\!X(\overrightarrow{\alpha_{\lambda, x_{0}}}(\widehat{x}_\lambda)) \times \nabla \overrightarrow{\alpha_{\lambda, x_{0}}}(\widehat{x}_\lambda)= O_{j, d-j}.
$$%
By this last fact, (\ref{parametrisation}) and (\ref{deriveegammalambda}), we have the following sequence of inequalities:
\begin{align*}
\MoveEqLeft[1]{\normp[d, d-j]{\nabla \overrightarrow{\alpha_{\lambda, x_{0},\delta}}(\widehat{x}_{\lambda}) - \nabla \overrightarrow{\alpha_{\lambda, x_{0}}}(\widehat{x}_{\lambda})}}\\
	& =\left\lVert-\!\left[\nabla\!X(\overrightarrow{\alpha_{\lambda, x_{0}, \delta}}(\widehat{x}_{\lambda}))|_{V_{\lambda}^{\perp}}\right]^{-1} \times 
 \!\left[{ \nabla\!X(\overrightarrow{\alpha_{\lambda, x_{0},\delta}}(\widehat{x}_{\lambda})) - \nabla\!X(\overrightarrow{\alpha_{\lambda, x_{0}}}(\widehat{x}_{\lambda}))
 }\right]\right. \\  
 	&\specialpos{\hfill \left.\times\nabla \overrightarrow{\alpha_{\lambda, x_{0}}}(\widehat{x}_{\lambda)}\rule{0pt}{14pt}\right\rVert_{d,d-j}} \\
	& \le\normp[d, j]{
\!\left[\nabla\!X(\overrightarrow{\alpha_{\lambda, x_{0}, \delta}}(\widehat{x}_{\lambda}))|_{V_{\lambda}^{\perp}}\right]^{-1}
} \times \normp[j, d]{
\nabla\!X(\overrightarrow{\alpha_{\lambda, x_{0},\delta}}(\widehat{x}_{\lambda})) - \nabla\!X(\overrightarrow{\alpha_{\lambda, x_{0}}}(\widehat{x}_{\lambda}))
}\\  
	& \specialpos{\hfill\times \normp[d, d-j]{
\nabla \overrightarrow{\alpha_{\lambda, x_{0}}}(\widehat{x}_{\lambda})
}  }\\
	& \le \sup_{z \in K} \normp[d, j]{
(\nabla\!X(z)|_{V_{\lambda}^{\perp}})^{-1}
}  \sup_{\widehat{x}_\lambda \in \overline{R_{x_{0}}^{\lambda}}} 
\normp[j, d]{\nabla\!X(\overrightarrow{\alpha_{\lambda, x_{0},\delta}}(\widehat{x}_{\lambda})) - \nabla\!X(\overrightarrow{\alpha_{\lambda, x_{0}}}(\widehat{x}_{\lambda}))}\\
	&\specialpos{\hfill\times \sup_{\widehat{x}_\lambda \in \overline{R_{x_{0}}^{\lambda}}} 
\normp[d, d-j]{\nabla \overrightarrow{\alpha_{\lambda, x_{0}}}(\widehat{x}_{\lambda})}.}
\end{align*}

The first term in the last inequality is bounded as explained above.
The third term is also bounded because $\nabla \overrightarrow{\alpha_{\lambda, x_{0}}}$ is continuous on
$\overline{R_{x_{0}}^{\lambda}}$ which is a compact set.
In the same way, the second term tends to zero because we already know that $\overrightarrow{\alpha_{\lambda, x_{0},\delta}}$ uniformly converges to $\overrightarrow{\alpha_{\lambda, x_{0}}}$ on $\overline{R_{x_{0}}^{\lambda}}$ and $\nabla\!X$ is continuous on $K$, a compact set.

We will show that 
$$
\lim_{\delta \to 0} \int_{\calC_{X}^{D^{r}}(y+\delta)}Y(z) \ud \sigma_{d-j}(z)= \int_{\calC_{X}^{D^{r}}(y)}Y(z) \ud \sigma_{d-j}(z).
$$%
Using first that $\supp(Y)$ is a compact subset of $\reels^d$ included in the open set $\Gamma(\lambda)$, we can prove that there exists an open set  $O$ contained in  $\reels^d$ such that:
$\supp(Y) \subset O \subset \overline{O} \subset \Gamma(\lambda)$.\\
Then, let us notice that since $\overline{O} \subset \Gamma(\lambda) \subset D$, the set $\overline{O} \cap \calC_{X}^{D^{r}}(y)$ is a compact set of $\reels^d$.\\
Let us build a partition of  unity $\{\pi_1, \dots, \pi_{m} \}$ for this compact manifold in the following way.\\
Consider  $x \in  \overline{O} \cap \calC_X^{D^r}(y)$.
Since $x \in \overline{O} \subset \Gamma(\lambda)$, $x \in \Gamma(\lambda)$ it turns out that  $J_{X}^{(\lambda)}(x) \neq 0$ and we can construct the open set  $U_{x}^{\lambda}$.
Since $U_{x}^{\lambda}$ is open, we can choose a real number $r_{x}^{\lambda} >0$ such that
the closed ball of $\reels^d$ of center $x$ and radius  $r_{x}^{\lambda}$ is contained in $U_{x}^{\lambda}$, let $\overline{B}(x, r_{x}^{\lambda}) \subset U_{x}^{\lambda}$.\\
We know that $\overline{O} \cap \calC_X^{D^r}(y)$ is a compact set of $\reels^d$, then we can cover $\overline{O} \cap\calC_X^{D^r}(y)$ with a finite number of these balls, that is $\overline{O} \cap \calC_X^{D^r}(y) \subset \cup_{i=1}^{m}
B(x_{i}, r_{x_{i}}^{\lambda})$ such that for $i=1, \dots, m$, we still have:
$$ 
\overline{B}(x_{i}, r_{x_{i}}^{\lambda}) \subset U_{x_i}^{\lambda} .
$$%
In the same way as on page \pageref{partition unite}, we obtain a partition of unity of $\overline{O} \cap \calC_X^{D^r}(y)$, \ie  $\{\pi_1, \dots, \pi_{m} \}$, such that:
\begin{enumerate}
	\item $\pi_{i}$ is $C^{\infty}$ on $\overline{O} \cap \calC_X^{D^r}(y) $, $\pi_{i} \ge 0$, $i=1, \dots, m$;
	\item $\supp(\pi_{i}) \subset \overline{O} \cap \calC_X^{D^r}(y)\cap \overline{B}(x_{i}, r_{x_{i}}^{\lambda}) \subset U_{x_i}^{\lambda}$;
	\item $\sum_{i=1}^{m} \pi_{i}(x) =1$, $x \in \overline{O} \cap \calC_X^{D^r}(y) $.
\end{enumerate}
Moreover, the following sequence of inclusions is true
\begin{multline*}
 \supp(Y) \cap \calC_X^{D^r}(y) \subset O \cap \calC_X^{D^r}(y) \subset \overline{O} \cap \calC_X^{D^r}(y)\\
 	\subset \cup_{i=0}^{m} U_{x_i}^{\lambda}$, $x_{i} \in \overline{O} \cap \calC_X^{D^r}(y), i= 1, \dots, m.
\end{multline*}

 Let us use the fact that $\supp(Y)$ is a compact set of $\reels^d$ and that $O$ is open, in the following way.
Consider $\omega \in \supp(Y)$.
Given that $\supp(Y) \subset  O \subset \Gamma(\lambda)$, there exist two open sets $\widetilde{U}(\omega)$ and $U(\omega)$ containing $\{\omega \}$ and $R_{\omega} >0$, $\widetilde{U}(\omega) \subset U(\omega) \subset O$ such that the restriction $f_{\lambda}|_{\widetilde{U}(\omega)}$ has an inverse on the open ball  $B(f_{\lambda}(\omega), R_{\omega}/2)$ and the restriction $f_{\lambda}|_{U(\omega)}$ has an inverse on the open ball $B(f_{\lambda}(\omega), R_{\omega})$.
It is also possible to have  $\diam(U(\omega)) \le \inf_{i=1,\dots, m} \diam(\widetilde{W}_{x_{i}}^{\lambda})$, $\widetilde{W}_{x_{i}}^{\lambda}$ being the neighborhood  of zero in  $\reels^{j}$ used to construct the function  $\gamma_{\lambda, x_{i}}$ and also the local parametrization of the level curve at level $y+\delta$.
\\
 Since $\supp(Y)$ is a compact set in $\reels^d$,  it can be covered by a finite number of such sets:
$$
\supp(Y)  \subset \cup_{\ell=1}^{k} \widetilde{U}(\omega_{\ell}) \subset  \cup_{\ell=1}^{k} U(\omega_{\ell}) \subset O
 $$
where $\omega_{\ell} \in \supp(Y)$, for $\ell=1, \dots, k$.
 
Let us chose $\delta=(\delta_1, \dots, \delta_{j}) \in \reels^{j}$
such that $\normp[j]{\delta} \le \inf_{\ell=1, \dots, k} R_{\omega_{\ell}}/2$.\\
 We will prove that if $\delta$ is small enough, any element of  $\calC_X^{D^r}(y +\delta) \cap \supp(Y)$ belongs to
 $\overrightarrow{\alpha_{\lambda, x_{i}, \delta}}(R_{x_{i}}^{\lambda})$, at least for one index $i$ belonging to $1, \dots, m$.\\
Let us chose $z \in \calC_X^{D^r}({y+\delta}) \cap \supp(Y)$.
Since $z \in \supp(Y) $, there exists $\ell=1, \dots, k$, such that $z \in \widetilde{U}(\omega_{\ell})$, and given that $z \in \calC_X^{D^r}({y+\delta})$, $f_{\lambda}(z)=\pi_{V_{\lambda}}(z) + \sum_{k=1}^{j} \delta_{k} e_{\ell_{k}}
\in B(f_{\lambda}(\omega_{\ell}), R_{\omega_{\ell}}/2)$.\\
It follows that $\pi_{V_{\lambda}}(z) \in B(f_{\lambda}(z), \normp[j]{\delta}) \subset  B(f_{\lambda}(\omega_{\ell}), R_{\omega_{\ell}})$.
Thus there exists an unique $x \in U(\omega_{\ell}) \subset O$, such that
$$
f_{\lambda}(x)= \pi_{V_{\lambda}}(z) = \pi_{V_{\lambda}}(x) + \sum\limits_{k=1}^{j} (X_{k}(x) - y_{k})e_{\ell_{k}}.
$$%
So $\pi_{V_{\lambda}}(z) = \pi_{V_{\lambda}}(x)$ and $ \sum\nolimits_{k=1}^{j} (X_{k}(x) - y_{k})e_{\ell_{k}}=0$, so $X(x)= y$.
Since $x \in  O \cap \calC_X^{D^r}(y) \subset \cup_{i=0}^{m} U_{x_i}^{\lambda}$,  $x$ can be written as 
$ x =  \overrightarrow{\alpha_{\lambda, x_{i}}}(\widehat{x}_{\lambda}), \widehat{x}_{\lambda} \in R_{x_i}^{\lambda}$, for some $i=1, \dots, m$.
Finally and since  $\pi_{V_{\lambda}}(z) = \pi_{V_{\lambda}}(x)$, the vector $z$ can be written as $z= x + \pi_{V_{\lambda}^{\perp}}(z- x)$.
Moreover, $\pi_{V_{\lambda}^{\perp}}(z- x) =\sum_{k=1}^{j} \gamma_{k} e_{\ell_{k}}$
and
$\normp[j]{\gamma}= \normp[d]{\smash{\pi_{V_{\lambda}^{\perp}}(z- x)}} \le \normp[d]{z-x} \le \diam(U(\omega_{\ell})) \le \diam(\widetilde{W}_{x_{i}}^{\lambda})$,
because we have 
$$
	\sup_{\ell=1, \dots, k} \diam(U(\omega_{\ell})) \le \inf_{i=1, \dots, m} \diam(\widetilde{W}_{x_{i}}^{\lambda}).
$$%
Using Property \ref{item:p26,3}$.$ of function  $G$, we proved that $z= \overrightarrow{\alpha_{\lambda, x_{i}, \delta}}(\widehat{x}_{\lambda}), \widehat{x}_{\lambda} \in R_{x_i}^{\lambda}$, $i=1, \dots, m$.

Suppose that  $t, s\in (\cup_{i=0}^{m} U_{x_i}^{\lambda})  \cap \calC_X^{D^r}(y)$ are such that $\pi_{V_{\lambda}}(t)=\pi_{V_{\lambda}}(s)$ and $t \neq s$.
We can write $s=(\pi_{V_{\lambda}^{\perp}}(t)+\sigma_{0}\gamma)+\pi_{V_{\lambda}}(t)$, where $\sigma_0>0$, $\gamma \in V_{\lambda}^{\perp}$ and $\normp[d]{\gamma}=1$.
By defining $h(\sigma):=X(t+\sigma\gamma)$ we have $h(0)=h(\sigma_{0})=y$.
Let us specify that this function is well defined for $0 \le \sigma \le \sigma_0$, because $D$ is a convex open set and here is the only place where we use the convexity of the set $D$\label{convexe}.
Moreover, Rolle's theorem allows us  to state that if $h:=(h_1, \dots, h_{j})$, then for all $\ell =1, \dots, j$, there exists $\sigma_{\ell} \in (0,\sigma_{0}) \such \dot{h}_{\ell}(\sigma_{\ell})=\nabla\!X_{\ell}(t+\sigma_{\ell}\gamma)(\gamma)=0$.
Moreover,
$$
1=\normp[d]{\gamma}=\normp[d]{(\nabla\!X(t)|_{V_{\lambda}^{\perp}})^{-1}(\nabla\!X(t)|_{V_{\lambda}^{\perp}}(\gamma))}
\le M \normp[j]{\nabla\!X(t)(\gamma)},
$$%
where
$$
M:= \sum_{i=1}^{m} \sup_{x \in \overline{U_{x_i}^{\lambda}}}\normp[d, j]{(\nabla\!X(x)|_{V_{\lambda}^{\perp}})^{-1}}< \infty.
$$%
To ensure that the last norm is finite it is sufficient to take the open sets $U_{x_i}^{\lambda}$ sufficiently small.
Finally, using that $\dot{h}_{\ell}(\sigma_{\ell})=0$ we obtain
\begin{align*}
	1&\le M^2\sum_{\ell=1}^j \abs{\nabla\!X_{\ell}(t+\sigma_{\ell} \gamma)(\gamma)-\nabla\!X_{\ell}(t)(\gamma)}^2\\
	&\le M^2 \sum_{\ell=1}^j \normp[1, d]{\nabla\!X_{\ell}(t+\sigma_{\ell}\gamma)-\nabla\!X_{\ell}(t)}^2.
\end{align*}
Assume $\nabla\!X$ is a Lipschitz function.
Let $L$ be its Lipschitz's constant.
We have
$$
1\le  M^2 L^2 \sum_{\ell=1}^{j}\sigma_{\ell}^2 \le j M^2 L^2\sigma^2 ,
$$%
where $\sigma:=\max_{\ell=1, \dots, j} \sigma_{\ell}$.
From this, we obtain
\begin{eqnarray}
\label{proxi}
\normp[d]{s-t}=\sigma_{0}\ge \sigma\ge \frac1{\sqrt jML}:= a.
\end{eqnarray}
Let us prove that if  $z \in \calC_X^{D^r}({y+\delta}) \cap \supp(Y)$, there exists a unique $x \in \cup_{i=0}^{m} U_{x_i}^{\lambda}  \cap \calC_X^{D^r}(y)$, such that
$z= x + \pi_{V_{\lambda}^{\perp}}(z- x)$ and such that $\normp[d]{\smash{\pi_{V_{\lambda}^{\perp}}(z- x)}} \le \inf_{i=1, \dots, m} \diam(\widetilde{W}_{x_{i}}^{\lambda}) $ (note that it is always possible by taking the open sets  $\widetilde{W}_{x_{i}}^{\lambda}$  to ensure that $\diam(\widetilde{W}_{x_{i}}^{\lambda}) < {a}/{2}$).
\\
Since $z$ belongs to $\calC_X^{D^r}({y+\delta}) \cap \supp(Y)$, we have shown above the existence of this $x$.
Let us now show the uniqueness.
Assume that there exists another vector $x^{\prime} \in (\cup_{i=0}^{m} U_{x_i}^{\lambda})  \cap \calC_X^{D^r}(y)$, $x^{\prime} \neq x$, such that $z= x^{\prime} + \pi_{V_{\lambda}^{\perp}}(z- x^{\prime})$ and also that $\normp[d]{\smash{\pi_{V_{\lambda}^{\perp}}(z- x^{\prime})}} \le \inf_{i=1, \dots, m} \diam(\widetilde{W}_{x_{i}}^{\lambda})$.\\
From (\ref{proxi}) we necessarily have $\normp[d]{x-x^{\prime}} \ge a$.
Furthermore, it holds 
\begin{align*}
\normp[d]{x-x^{\prime}}
&\le \normp[d]{x- z} + \normp[d]{x^{\prime}-z}\\
&= \normp[d]{\smash{\pi_{V_{\lambda}^{\perp}}(z- x)}} + \normp[d]{\smash{\pi_{V_{\lambda}^{\perp}}(z- x^{\prime})}} \\
&<  \ds{\frac{a}{2}} + \ds{\frac{a}{2}} =a.
\end{align*}
We thus obtain a contradiction.

 Let us now consider $z \in \calC_X^{D^r}({y+\delta}) \cap \supp(Y)$,  since
 $x \in  O \cap \calC_X^{D^r}(y) \subset \overline{O} \cap \calC_X^{D^r}(y)$ and given that $\widehat{z}_{\lambda} = \widehat{x}_{\lambda}$, 
 we obtain the sequence of equalities 
 \begin{align*}
 Y(z) &= \!\left(\sum_{i=1}^{m} \pi_{i}(x)\right) \times Y(z)\\
&= \sum_{i=1}^{m} \pi_{i}(x) \times \1_{\{x \in U_{x_i}^{\lambda} \cap \calC_X^{D^{r}}(y)\}} \times Y(z)\\
 & = \sum_{i=1}^{m} \pi_{i}(x) \times  \1_{\{x =  \overrightarrow{\alpha_{\lambda, x_{i}}}(\widehat{x}_{\lambda})\}} \times \1_{\{\widehat{x}_{\lambda} \in R_{x_i}^{\lambda}\}} \times Y(z)\\
 & =  \sum_{i=1}^{m} \pi_{i}(\overrightarrow{\alpha_{\lambda, x_i}}(\widehat{z}_{\lambda})) \times  \1_{\{x =  \overrightarrow{\alpha_{\lambda, x_{i}}}(\widehat{z}_{\lambda})\}} \times \1_{\{\widehat{z}_{\lambda} \in R_{x_i}^{\lambda}\}} \times Y(z)\\
 & =  \sum_{i=1}^{m} \pi_{i}(\overrightarrow{\alpha_{\lambda, x_i}}(\widehat{z}_{\lambda})) \times  \1_{\{z =  \overrightarrow{\alpha_{\lambda, x_{i}, \delta}}(\widehat{z}_{\lambda})\}} \times Y(\overrightarrow{\alpha_{\lambda, x_{i}, \delta}}(\widehat{z}_{\lambda}))  \times \1_{\{\widehat{z}_{\lambda} \in R_{x_i}^{\lambda}\}}.
 \end{align*}

The last equality coming from the construction of the vector $x$ from the vector $z$ and from the uniqueness of the decomposition of $z$ on the set $(\cup_{i=0}^{m} U_{x_i}^{\lambda})  \cap \calC_X^{D^r}(y)$.\\
We have the following equality and convergence as $\delta$ tends towards zero
\begin{align*}
\MoveEqLeft[0]{\int_{\calC_X^{D^{r}}(y+\delta)}Y(z) \ud \sigma_{d-j}(z)}\\
&=\sum_{i=1}^{m} \int_{\{z =  \overrightarrow{\alpha_{\lambda, x_{i}, \delta}}(\widehat{z}_{\lambda}), \widehat{z}_{\lambda} \in R_{x_i}^{\lambda}\}} \pi_{i}(\overrightarrow{\alpha_{\lambda, x_{i}}}(\widehat{z}_{\lambda})) \times Y(\overrightarrow{\alpha_{\lambda, x_{i}, \delta}}(\widehat{z}_{\lambda})) \ud \sigma_{d-j}(z)\\
&=\sum_{i=1}^{m} \int_{R_{x_{i}}^{\lambda}} \pi_{i}(\overrightarrow{\alpha_{\lambda, x_{i}}}(\widehat{z}_{\lambda}))
Y(\overrightarrow{\alpha_{\lambda, x_{i}, \delta}}(\widehat{z}_{\lambda}))(\det (\nabla \overrightarrow{\alpha_{\lambda, x_{i}, \delta}}(\widehat{z}_{\lambda})\nabla \overrightarrow{\alpha_{\lambda, x_{i}, \delta}}(\widehat{z}_{\lambda})^T))^{1/2} \ud \widehat{z}_{\lambda}\\
&\cvg[\delta\to0]{} \sum_{i=1}^{m} \int_{R_{x_{i}}^{\lambda}} \pi_{i}(\overrightarrow{\alpha_{\lambda, x_{i}}}(\widehat{z}_{\lambda}))
Y(\overrightarrow{\alpha_{\lambda, x_{i}}}(\widehat{z}_{\lambda}))(\det (\nabla \overrightarrow{\alpha_{\lambda, x_{i}}}(\widehat{z}_{\lambda})\nabla \overrightarrow{\alpha_{\lambda, x_{i}}}(\widehat{z}_{\lambda})^T))^{1/2} \ud\widehat{z}_{\lambda}.
\end{align*}
The last convergence comes from the uniform convergence of $\overrightarrow{\alpha_{\lambda, x_{i},\delta}}$ to $\overrightarrow{\alpha_{\lambda, x_{i}}}$ on $\overline{R_{\bx_i}^{\lambda}}$ and of that of $\nabla \overrightarrow{\alpha_{\lambda, x_{i}, \delta}}$ to $\nabla \overrightarrow{\alpha_{\lambda, x_{i}}}$ on $\overline{R_{\bx_i}^{\lambda}}$.\\
But
\begin{align*}
\sum_{i=1}^{m} &\int_{R_{x_{i}}^{\lambda}}\pi_{i}(\overrightarrow{\alpha_{\lambda, x_{i}}}(\widehat{z}_{\lambda}))
Y(\overrightarrow{\alpha_{\lambda, x_{i}}}(\widehat{z}_{\lambda}))(\det (\nabla \overrightarrow{\alpha_{\lambda, x_{i}}}(\widehat{z}_{\lambda})\nabla \overrightarrow{\alpha_{\lambda, x_{i}}}(\widehat{z}_{\lambda})^T))^{1/2} \ud\widehat{z}_{\lambda}\\
&=\sum_{i=1}^{m} \int_{U_{x_i}^{\lambda} \cap \calC_X^{D^{r}}(y)} \pi_{i}(x) Y(x) \ud \sigma_{d-j}(x)\\
&= \int_{\overline{O} \cap \calC_X^{D^{r}}(y)} Y(x) \ud \sigma_{d-j}(x)\\
&= \int_{\calC_X^{D^{r}}(y)} Y(x) \ud \sigma_{d-j}(x),
\end{align*}
because $\supp(Y)  \subset \overline{O}$.\\
In summary, we have proved the continuity of the function
$$
y \mapsto \int_{\calC_X^{D^r}(y)}Y(x) \ud\sigma_{d-j}(x),
$$%
under the hypothesis  that $Y: D_X^r \subset \reels^d \to \reels$ is a continuous function satisfying $\supp(Y)  \subset \Gamma(\lambda)$.\\
Finally, we no longer assume that $\supp(Y) \subset \Gamma(\lambda)$, just that $\supp(Y) \subset D_X^r$.

Let us introduce two new functions.
For $\lambda \in B_{j}$ and $t \in D$, let
\begin{equation}
	\phi_\lambda(t):=\inf_{\gamma \in V^{\perp}_{\lambda}, \normp[d]{\gamma}=1} \normp[j]{\nabla\!X(t)|_{V^{\perp}_{\lambda}}(\gamma)},
	\label{etalambda}
\end{equation}
and
$$
	\phi(t):=\sup_{\lambda \in B_{j}}\phi_\lambda(t).
$$%
These two functions are Lipschitz and therefore continuous, with the same constant  $L$, as $\nabla\!X$.
Indeed, let us consider the first function.
Considering two points $t$ and $t^\star$, we have for any $\gamma \in V_{\lambda}^{\perp}$ satisfying $\normp[d]{\gamma}=1$,
\begin{align*}
\normp[j]{\nabla\!X(t)|_{V^{\perp}_{\lambda}}(\gamma)} &= \normp[j]{\nabla\!X(t)(\gamma)} \\
&\le \normp[j]{\nabla\!X(t)(\gamma)-\nabla\!X(t^\star)(\gamma)}+\normp[j]{\nabla\!X(t^\star)|_{V^{\perp}_{\lambda}}(\gamma)}.
\end{align*}
Using that $\nabla\!X$ is Lipschitz, we obtain
$$
	\normp[j]{\nabla\!X(t)|_{V^{\perp}_{\lambda}}(\gamma)} \le L \normp[d]{t-t^\star}+ \normp[j]{\nabla\!X(t^\star)|_{V^{\perp}_{\lambda}}(\gamma)},
$$%
then
\begin{equation}
	\phi_\lambda(t)\le L \normp[d]{t-t^\star}+\phi_\lambda(t^\star).
	\label{inegalite phi}
\end{equation} 
Since a symmetric inequality can be proved, we obtain finally
$$
	\abs{\phi_\lambda(t)-\phi_\lambda(t^\star)}\le L \normp[d]{t-t^\star}.
$$%
Let us now study the second function, $\phi$.\\
By (\ref{inegalite phi}) we can write
$$
	\phi(t)\le L \normp[d]{t-t^\star}+\phi(t^\star),
$$%
and as before
$$
\abs{\phi(t)-\phi(t^\star)} \le L \normp[d]{t-t^\star}.
$$%
Let us prove that
\begin{equation}
	\!\left({\lambda \in B_{j} \mbox{ \ and \ } t \in \Gamma(\lambda)}\right) \iff\left(\phi_{\lambda}(t)>0\right).
	\label{propriété2}
\end{equation}
Let us consider $t \in \Gamma(\lambda)$ for $\lambda \in B_j$.
By Remark \ref{inversibilite} that means that
$\nabla\!X(t)|_{V_{\lambda}^{\perp}}$
has an inverse and moreover that
$\ker(\nabla\!X(t)|_{V^{\perp}}) = \vec{0}|_{V^{\perp}}$.
But there exists
$\gamma_{0} \in V^{\perp}_{\lambda}$, $\normp[d]{\gamma_{0}}=1$
 such that
  $\phi_\lambda(t)= \normp[j]{\nabla\!X(t)|_{V^{\perp}_{\lambda}}(\gamma_{0})}$,
  and this implies that  $\phi_\lambda(t) >0$.
  
To prove the other implication, let us suppose that for a $\lambda \in B_{j}$, $t \notin \Gamma(\lambda)$, \ie  $\nabla\!X(t)|_{V_{\lambda}^{\perp}}$ has no inverse.
Then there exists  $\gamma \in V_{\lambda}^{\perp}$, $\normp[d]{\gamma}=1$ such that
$\nabla\!X(t)|_{V^{\perp}_{\lambda}}(\gamma)=0$, and this implies  $\phi_{\lambda}(t)=0$.
\\
Now let us prove that
\begin{equation}
	\left(t \in D_{X}^{r}\right) \iff \left(\phi(t) >0\right).
	\label{propriété}
\end{equation}
Indeed,  by (\ref{propriété2}), we have the following equivalences
$$
\left(t \in D_{X}^{r}\right) \iff \!\left({ \exists \lambda \in B_j, t \in \Gamma(\lambda) }\right)
\iff \!\left({ \exists \lambda \in B_j, \phi_{\lambda}(t)>0  }\right)
 \iff \left(\phi(t) >0\right)
$$

We will construct a partition of unity of  $D_{X}^{r}$ whose support intersected by $D_{X}^{r}$ will be  included in $\Gamma(\lambda)$, $\forall \lambda \in B_{j}$.
We will denote this partition by $\eta_{\lambda}$.
\\
Let us first consider the function $\chi_\lambda(t):=(2\phi_\lambda(t)-\phi(t))^{+}$.
Since $\phi_\lambda$ and $\phi$ are Lipschitz functions, it follows that $2\phi_\lambda-\phi $ remains  Lipschitz.
So, the function  $\chi_\lambda$ is also Lipschitz and a fortiori  continuous.

Let us show that 
\begin{equation}
\label{propriété1}
\left(t \in D_{X}^{r}\right)
	\Longrightarrow
\left({\sum\limits_{\lambda \in B_{j}} \chi_{\lambda}(t) >0}\right).
\end{equation}

Consider $t \in D_{X}^{r}$.
Let us suppose that $\sum\limits_{\lambda \in B_{j}} \chi_{\lambda}(t) =0$ and deduce a contradiction.
Since $\sum\limits_{\lambda \in B_{j}} \chi_{\lambda}(t) =0$, then $\forall \lambda \in B_{j}$ we have $\chi_{\lambda}(t) =0$, \ie $\phi_\lambda(t)\le \frac{1}{2} \phi(t)$.
It is indeed a consequence since this inequality is true $\forall \lambda \in B_{j}$, that $\phi(t) \le \frac{1}{2} \phi(t)$ and  $\phi(t) =0$, which is in contradiction with (\ref{propriété}).

Let $\forall t \in D_{X}^{r}$,
\begin{equation}
\label{eta}
\eta_\lambda(t):={\chi_\lambda(t)}\left/{\sum\limits_{\lambda \in B_{j}}\chi_\lambda(t)}\right.,
\end{equation}
which is possible from (\ref{propriété1}).
It is now easy to see that $\eta_\lambda$ is continuous on $D_{X}^{r}$ since $\chi_\lambda$ is also continuous on $D_{X}^{r}$.\\
It only remains to prove that the intersection of the support of this function with $D_{X}^{r}$ is included in $\Gamma(\lambda)$.\\
For all  $C >0$, let us construct an open set $O_{C}$ including
$$
(\Gamma(\lambda))^{c_1}\cap \{t \in D, \phi(t) > C \},
$$%
contained in the set
$$
\{t \in D, \chi_{\lambda}(t) = 0\}.
$$%
More precisely, let us prove that for a given $C>0$, if $\delta \le {\frac{C}{2L}}$ (with $L$ being the Lipschitz constant of the function $\nabla\!X$), then
$$
	((\Gamma(\lambda))^{c_1})_\delta\cap \{ t \in D: \phi(t)>C\}\subset \{t \in D: \chi_\lambda(t)=0\}
$$%
where for any set $A$ we have defined the open set $A_\delta:=\{x \in \reels^d: d(x,A)<\delta\}$.\\
Thus, let $C >0$ be fixed and $t \in ((\Gamma(\lambda))^{c_1})_\delta\cap \{ t \in D: \phi(t)>C\}$.
Since $t \in ((\Gamma(\lambda))^{c_1})_\delta$, there exists $t^{\prime} \in B(t, \delta)$ such that  $t^{\prime} \in (\Gamma(\lambda))^{c_1}$ (and also such that $\phi_{\lambda}(t^{\prime})=0$ as a consequence of  (\ref{propriété2})).
Then we have, since $\phi_{\lambda}$ is Lipschitz with Lipschitz constant $L$ that
\begin{align*}
	\phi_{\lambda}(t) &\le \abs{\phi_{\lambda}(t) -\phi_{\lambda}(t^{\prime})} + \phi_{\lambda}(t^{\prime})\\
&\le L \normp[d]{t-t^{\prime}} \le L \delta \le \frac{C}{2}.
\end{align*}
Therefore, $2 \phi_{\lambda}(t) \le C < \phi(t)$, and this leads that $\chi_{\lambda}(t)=0$.
We have proved that for all $C>0$ and $\forall \delta \le \frac{C}{2L}$, we have the inclusion
\begin{multline*}
\{t \in D_{X}^{r}: \chi_\lambda(t)\neq 0\} = \{t \in D_{X}^{r}: \eta_\lambda(t)\neq 0\} \\
\subset \left(\rule{0pt}{11pt}\left(\rule{0pt}{10pt}\left(\rule{0pt}{9pt}\left(\Gamma(\lambda)\right)^{c_1}\right)_\delta\right)^{c} \cap D_{X}^{r}\right) \cup \{t \in D_{X}^{r}: \phi(t) \le C\},
\end{multline*}
where we recall that the symbol $c$ denotes the complementary set with respect to $\reels^d$.\\
Noting that 
$\!\left(\rule{0pt}{10pt}\!\left(\rule{0pt}{9pt}\!\left(\Gamma(\lambda)\right)^{c_1}\right)_\delta\right)^c$
is a closed set contained in  $\Gamma(\lambda) \cup D^{c}$, we have
$$
	\supp(\eta_{\lambda})  \subset [(\Gamma(\lambda) \cup D^{c}) \cap \overline{D_{X}^{r}}] \cup (\cap_{C >0} \overline{ \{t \in D_{X}^{r}: \phi(t) \le C\}}),
$$%
that is
$$
	\supp(\eta_{\lambda}) \cap D_{X}^{r} \subset \Gamma(\lambda) \cup (\cap_{C >0} \overline{ \{t \in D_{X}^{r}: \phi(t) \le C\}} \cap D_X^r).
$$%
To complete the proof, it is enough to show that $\cap_{C >0} \overline{ \{t \in D_{X}^{r}: \phi(t) \le C\}} \cap D_X^r= \emptyset$.
Indeed, consider $z \in \cap_{C >0} \overline{ \{t \in D_{X}^{r}: \phi(t) \le C\}} \cap D_X^r$.
Then $z \in D_X^r$ and for all  $C>0$, there exists a sequence of points $z_{n, C}$ of $D_X^r$, satisfying  $\phi(z_{n, C}) \le C$  which converges to  $z \in D_X^r$.
Since  the function $\phi$ is continuous on  $D$ and also on $D_X^r$, it holds that $\phi(z) \le C$.
This last inequality is true for all $C >0$, then we get that  $\phi(z)=0$.
From (\ref{propriété}) we easily obtain that $z \in (D_X^r)^{c_1}$.
But $z \in D_X^r$.

We have proven that 
\begin{equation}
	\supp(\eta_{\lambda}) \cap D_X^r \subset \Gamma(\lambda).
	\label{support}
\end{equation}
In this form, for $t \in D_{X}^{r}$ we have
$$
	Y(t)= \sum\limits_{\lambda \in B_{j}} \eta_{\lambda}(t) Y(t) = \sum\limits_{\lambda \in B_{j}}Y_{\lambda}(t),
$$%
where we set for $t \in D_{X}^{r}$, $Y_{\lambda}(t):=  \eta_{\lambda}(t) Y(t)$.

The function $Y_{\lambda}$ is a continuous function on $D_{X}^{r}$ with compact support included in $\Gamma(\lambda)$, since $\supp(Y) \subset D_X^r$ by hypothesis and from the inclusion (\ref{support}).\\
 We have $\forall y \in \reels^j$
$$
	\int_{\calC_{X}^{D^{r}}(y)}Y(z) \ud \sigma_{d-j}(z)=\sum_{\lambda \in B_{j}} \int_{\calC_{X}^{D^{r}}(y)}Y_{\lambda}(z) \ud \sigma_{d-j}(z).
$$%
The continuity of the left-hand side integral as a function of the variable $y$ is a consequence of the continuity of each of the terms of the right sum.
This last fact is an application of the above procedure.
This completes the proof of Theorem \ref{continueBis} in the case where we have chosen the bounded open set  $D_1$  of $\reels^d$ equal to $D$ convex (bounded).\\
Now suppose that $D$ is a convex open set of $\reels^d$, which can be unbounded.
The function $X:  D \subset  \reels^{d} \rightarrow \reels^{j}$ is a function
of class $C^{1}(D,\reels^{j}) \such \nabla\!X$ is Lipschitz and the function $Y: D_1 \subset \reels^{d} \rightarrow \reels$ is a continuous function defined on $D_1$, an open and bounded subset of $\reels^d$ included in $D \such \supp(Y) \subset D_{X_{|D_1}}^r$.\\
In this case, the function  $X$ restricted to the bounded open set $D_1$ \ie $X_{|D_1}$, is such that $X_{|D_1}:  D_1 \subset  \reels^{d} \rightarrow \reels^{j}$ is $C^{1}(D_1,\reels^{j})$ and such that $\nabla\!X_{|D_1}$ is Lipschitz.\\
One can apply the previous procedure to these two functions  $X_{|D_1}$ and $Y$ and also to the open set  $D_1$.
It is possible that $D_1$ is not a convex set, but this is not a real problem.
Indeed, if we refer to page \pageref{convexe}, the only place where we used the convexity of the open set, we realize that the important thing is to be able to apply Rolle's theorem to the function $h$ which is defined there and to use the fact that the function $\nabla\!X$ is Lipschitz.
Since $D_1$ could not be convex, we were not sure if we could do it, but this is not the case if we work on $D$ which is convex.\\
The proof of the theorem is finished.
\end{proofarg}

Now we are able to exhibit a class of processes  $X$ and $Y$ satisfying the hypotheses ${\bH _2}$ and ${\bH _{4}}$ through  the following assumption $\bA_0$ and the following proposition whose proof is based  on the one given by Caba\~na \cite{AL121221}.
In what follows, we will give a new proof slightly more general  that the original one.
\begin{itemize}
\item $\bA_0$: $X: \Omega \times D \subset \Omega\times \reels^{d} \rightarrow \reels^{j}$ ($j \le d$) is a random field that belongs to
$C^{1}(D,\reels^{j})$, where $D$ is a bounded open convex set of $\reels^{d}$, such that for almost all $\omega \in \Omega$, the process $\nabla\!X(\omega)$ is Lipschitz with  Lipschitz constant $L_{X}(\omega)$ satisfying $\bbE\left[L_X^d(\bm{\cdot})\right] < \infty$.
 Also, $Y: \Omega \times D \subset \Omega \times \reels^{d} \rightarrow \reels$ is a continuous process such $\exists \lambda \in B_{j} \such \supp(Y)  \subset \Gamma(\lambda)$.
Moreover, $\normp[d, j]{\smash{(\nabla\!X(\bm{\cdot})|_{V_{\lambda}^{\perp}})^{-1}}}$,  $Y(\bm{\cdot})$ and $\normp[j, d]{\nabla\!X(\bm{\cdot})}$ are assumed uniformly  bounded on the support of  $Y$, the bounds not depending on $\omega$ ($\in \Omega$).
\end{itemize}
\spacebefore
\begin{prop}
\label{H1etH4}
If $X$ and $Y$ satisfy the assumption $\bA_0$, then the hypotheses ${\bH _2}$ and ${\bH _{4}}$ are satisfied.
\end{prop}
\spacebefore
\begin{proofarg}{Proof of Proposition \ref{H1etH4}} 
For almost all $\omega \in \Omega$ the field $X(\omega):  D \subset  \reels^{d} \rightarrow \reels^{j}$ ($j \le d$) is
$C^{1}(D,\reels^{j})$ and such that $\nabla\!X(\omega)$ is Lipschitz.
The process $Y(\omega): D \subset \reels^{d} \rightarrow \reels$ is a continuous function such that $\supp(Y(\omega)) \subset \Gamma(\lambda)(\omega) \subset D_{X(\omega)}^r$.
The set $D$ is an open and convex bounded set of $\reels^{d}$.
According to Theorem \ref{continueBis} the function 
$$
y \longmapsto \int_{\calC_{X(\omega)}^{D^{r}}(y)} Y(\omega)(x) \ud \sigma_{d-j}(x)
$$%
is a continuous function of the variable $y$.
The same is true for
$$
	 \int_{\calC_{X(\omega)}^{D^{r}}(y)} \abs{Y(\omega)(x)} \ud \sigma_{d-j}(x).
$$%
 
Let us find an upper bound for $\int_{\calC_{X}^{D^{r}}(y)} \abs{Y(x)} \ud \sigma_{d-j}(x)$ which would be an integrable random variable that does not depend on $y$.
Then, according to the dominated convergence theorem, hypotheses ${\bH _2}$ and ${\bH _{4}}$ will be fullfilled.\\
Since  $\supp(Y)  \subset \Gamma(\lambda)$, we can construct a  partition of unity  of $\supp(Y)$.
As on page \pageref{partition unite} getting as in (\ref{definition integrale})
\begin{align*}
\MoveEqLeft[1]{\int_{\calC_X^{D^{r}}(y)}Y(x) \ud \sigma_{d-j}(x)}\\
&= \sum_{i=1}^{m} \int_{R_{x_{i}}^{\lambda}} \pi_{i}\!\left(\overrightarrow{\alpha_{\lambda, x_{i}}}(\widehat{x}_{\lambda})\right)
Y\!\left(\overrightarrow{\alpha_{\lambda, x_{i}}}(\widehat{x}_{\lambda})\right)
\left(\det \!\left(\nabla \overrightarrow{\alpha_{\lambda, x_{i}}}(\widehat{x}_{\lambda})\nabla \overrightarrow{\alpha_{\lambda, x_{i}}}(\widehat{x}_{\lambda})^T\right)\right)^{1/2} \ud\widehat{x}_{\lambda}.
\end{align*}
Consider $\widehat{x}_{\lambda}$ fixed in $R_{x_{i}}^{\lambda} \such \overrightarrow{\alpha_{\lambda, x_{i}}}(\widehat{x}_{\lambda}) \in \supp(Y)$, $i=1, \dots, m$.
\\
We have
$\left(\det \!\left(\nabla \overrightarrow{\alpha_{\lambda, x_{i}}}(\widehat{x}_{\lambda})\nabla \overrightarrow{\alpha_{\lambda, x_{i}}}(\widehat{x}_{\lambda})^T\right)\right)^{1/2}  \le \normp[d\times (d-j)]{\nabla \overrightarrow{\alpha_{\lambda, x_{i}}}(\widehat{x}_{\lambda})}^{d}$.\\
Let us uniformly bound  $\normp[d, d-j]{\nabla \overrightarrow{\alpha_{\lambda, x_{i}}}(\widehat{x}_{\lambda}))}$, for all $\widehat{x}_{\lambda}$ in $R_{x_{i}}^{\lambda}$ such that 
$\overrightarrow{\alpha_{\lambda, x_{i}}}(\widehat{x}_{\lambda}) \in \supp(Y) $.\\
For any $\widehat{x}_{\lambda} \in R_{x_{i}}^{\lambda}$, we have $X\!\left(\overrightarrow{\alpha_{\lambda, x_{i}}}(\widehat{x}_{\lambda})\right)= y$.
Taking derivatives in this equality on the open set $ R_{x_{i}}^{\lambda}$, we obtain
$$
\nabla\!X\!\left(\overrightarrow{\alpha_{\lambda, x_{i}}}(\widehat{x}_{\lambda})\right) \times \nabla\overrightarrow{\alpha_{\lambda, x_{i}}}(\widehat{x}_{\lambda})=0.
$$%
Using (\ref{parametrisationlocale}) and that $\overrightarrow{\alpha_{\lambda, x_{i}}}(\widehat{x}_{\lambda})\in \Gamma(\lambda)$, $\forall u \in \reels^{d-j}$, $u:=(u_1, \dots, u_{d-j})$ we get\\
\vbox{\begin{align*}
\MoveEqLeft[0]{\nabla \overrightarrow{\alpha_{\lambda, x_{i}}}(\widehat{x}_{\lambda})(u)}\\
&=-\left[{ \nabla\!X\!\left(\overrightarrow{\alpha_{\lambda, x_{i}}}(\widehat{x}_{\lambda})\right)|_{V^{\perp}_{\lambda}}}\right]^{-1}\!\left({\nabla\!X\!\left(\overrightarrow{\alpha_{\lambda, x_{i}}}(\widehat{x}_{\lambda})\right)|_{V_{\lambda}}\!\left(\sum_{k=1}^{d-j} u_{k} e_{i_{k}}\right)}\right) 
+ \sum_{k=1}^{d-j} u_{k} e_{i_{k}}.
\end{align*}}
Since $\normp[d,j]{(\nabla\!X(\bm{\cdot})|_{V_{\lambda}^{\perp}})^{-1}}$, $Y(\bm{\cdot})$ and $\normp[j,d]{\nabla\!X(\bm{\cdot})}$ are uniformly bounded on $\supp(Y)$, and the bound does not depend on $\omega$, then we have
\begin{multline*}
 \int_{\calC_X^{D^{r}}(y)} \abs{Y(x)} \ud \sigma_{d-j}(x)\\
 \le\bC  \int_{\pi_{V_{\lambda}}(D)} \sum_{i=1}^{m}  \pi_{i}\!\left(\overrightarrow{\alpha_{\lambda, x_{i}}}(\widehat{x}_{\lambda})\right) \1_{\{\overrightarrow{\alpha_{\lambda, x_{i}}}(\widehat{x}_{\lambda}) \in \supp(Y)\}} \1_{\{ \widehat{x}_{\lambda} \in R_{x_{i}}^{\lambda}\}}
\ud\widehat{x}_{\lambda}.
\end{multline*}
For $\omega \in \Omega$ and $\widehat{x}_{\lambda} \in \pi_{V_{\lambda}}(D)$, we consider the set $A$ defined by
\begin{multline*}
 A:=\left\{ \overrightarrow{\alpha_{\lambda, x_{i}}}(\widehat{x}_{\lambda})(\omega), \widehat{x}_{\lambda} \in R_{x_{i}}^{\lambda}(\omega)\right.\\
  \left.\rule{0pt}{12pt}\text{~and~} \overrightarrow{\alpha_{\lambda, x_{i}}}(\widehat{x}_{\lambda})(\omega) \in \supp(Y)(\omega), i=1, \dots, m \right\}.
\end{multline*}
We partition the set $A$ into equivalence classes.
An equivalence class $A_{i_{0}}$ for $i_{0}=1, \dots, m$, is 
\begin{align*}
\MoveEqLeft[3]{A_{i_{0}}:=\left\{ \overrightarrow{\alpha_{\lambda, x_{i}}}(\widehat{x}_{\lambda})(\omega), \widehat{x}_{\lambda} \in R_{x_{i}}^{\lambda}(\omega) \cap R_{x_{i_{0}}}^{\lambda}(\omega) \right.}\\
&\specialpos{\hfill\text{ and } \overrightarrow{\alpha_{\lambda, x_{i}}}(\widehat{x}_{\lambda})(\omega)= \left.\overrightarrow{\alpha_{\lambda, x_{i_{0}}}}(\widehat{x}_{\lambda})(\omega) \in \supp(Y)(\omega),~i=1, \dots, m\right\}.
}\end{align*}
By (\ref{item.p24.3}), page \pageref{item.p24.3}, we have $\sum_{i=1}^{m} \pi_{i}(x) =1$, $x \in \supp(Y)$.
Moreover, in each class we bound the corresponding sum by one.
It only remains to count the maximal number of equivalence classes.\\
For counting the classes let us take two elements belonging to two different classes.
To fix the ideas, we will take for example
$$
t:=\overrightarrow{\alpha_{\lambda, x_{i}}}(\widehat{x}_{\lambda})(\omega), \widehat{x}_{\lambda} \in R_{x_{i}}^{\lambda}(\omega) \mbox{ and } \overrightarrow{\alpha_{\lambda, x_{i}}}(\widehat{x}_{\lambda})(\omega) \in \supp(Y)(\omega)
$$
and
$$
s:=\overrightarrow{\alpha_{\lambda, x_{j}}}(\widehat{x}_{\lambda})(\omega), \widehat{x}_{\lambda} \in R_{x_{j}}^{\lambda}(\omega) \mbox{ and } \overrightarrow{\alpha_{\lambda, x_{j}}}(\widehat{x}_{\lambda})(\omega) \in \supp(Y)(\omega),
$$
$ i,j=1, \dots,m \mbox{ with } t \neq s.$

It is clear that $t$ and $s$ are two different  elements  of $\reels^{d}$, they have the same projection on $V_{\lambda}$ and belong to the level curve $C_{X}^{D^r}(y)$.
Repeating the proof given on page \pageref{convexe} since  $t \in \supp(Y)(\omega)$ and that $\normp[d, j]{(\nabla\!X(\bm{\cdot})|_{V_{\lambda}^{\perp}})^{-1}}$ is uniformly bounded on $\supp(Y) $ by a constant $\bC $, we get the following bound  
$$
1\le {\bC ^2} \sum_{\ell=1}^j \normp[1, d]
	{\nabla\!X_{\ell}\!\left(\overrightarrow{\alpha_{\lambda, x_{i}}}(\widehat{x}_{\lambda})
	+\sigma_{\ell}\gamma\right)(\omega)
	-\nabla\!X_{\ell}\!\left(\overrightarrow{\alpha_{\lambda, x_{i}}}(\widehat{x}_{\lambda})\right)(\omega)}^2.
$$%
But for almost all $\omega \in \Omega$, $\nabla\!X(\omega)$ is Lipschitz with Lipschitz constant $L_{X}(\omega)$, we get 
$$
1 \le j {\bC ^2} L_{X}^{2}(\omega) \sigma^{2}(\omega), 
$$%
where $\sigma:=\max_{\ell=1, \dots, j} \sigma_{\ell}$.\\
As in (\ref{proxi}), we finally obtain the following bound 
$$
	\normp[d]{s-t}=\sigma_{0}(\omega) \ge \sigma(\omega)\ge \frac{1}{\sqrt j\bC L_{X}(\omega)}:= a(\omega).
$$%
The open balls of center  $t$ and  $s$ and diameter $a(\omega)$ do not intersect.
We have at most $({\diam(D)}/{a(\omega)})^{d}$ balls of diameter $a(\omega)$ and thus  at most $({\diam(D)}/{a(\omega)})^{d}$ equivalence classes.\\
Finally, for almost all $\omega \in \Omega$:
\begin{align*}
\int_{\calC_{X(\omega)}^{D^{r}}(y)} \abs{Y(\omega)(x)} \ud \sigma_{d-j}(x) &\le \bC  \int_{\pi_{V_{\lambda}}(D)} \!\left({\diam(D)}/{a(\omega)}\right)^{d} \ud\widehat{x}_{\lambda}\\
&\le \bC  \sigma_{d-j}(\pi_{V_{\lambda}}(D)) \!\left(\diam(D) \sqrt{j} \bC \right)^{d} L_{X}^{d}(\omega)\\
&\le \bC  L_{X}^{d}(\omega).
\end{align*}
Since $\bbE\!\left[L_X^d(\bm{\cdot})\right] < \infty$,  $
y \longmapsto\bbE\left[\int_{\calC_{X}^{D^{r}}(y)} Y(x) \ud \sigma_{d-j}(x)\right]$ is continuous.
The same is true for
$\bbE\left[\int_{\calC_{X}^{D^{r}}(y)} \abs{Y(x)} \ud \sigma_{d-j}(x)\right]$.

Then the hypotheses ${\bH _2}$ and  ${\bH _{4}}$ have been verified.
This completes the proof of the proposition.
\end{proofarg}
\spacebefore
\begin{rema}
\label{moduledeYBis}
In $\bA_0$, we can replace  the condition $\bbE\left[L_X^d(\bm{\cdot})\right] < \infty$ by the following one: $\exists L>0 \ni \text{for almost all }\omega \in \Omega$, $\supp(Y(\omega)) \neq \emptyset\Longrightarrow L_{X}(\omega) \le L$.
In this case, the hypotheses ${\bH _2}$ and ${\bH _{4}}$ hold.

\hfill$\bullet$
\end{rema}
\spacebefore
\begin{proofarg}{Proof of Remark \ref{moduledeYBis}} 
It suffices to replace in the proof of the previous proposition  $L_{X}(\omega)$ by $L$ and $a(\omega)$ by $\frac{1}{\sqrt j\bC L}$.
In this case, for almost all $\omega \in \Omega$, we have
$$
	\int_{\calC_{X(\omega)}^{D^{r}}(y)} \abs{Y(\omega)(x)} \ud \sigma_{d-j}(x) \le \bC .
$$%
and this implies the integrability.
\end{proofarg}
\spacebefore
\begin{rema}
\label{moduleconvexe}
We can generalize  Proposition \ref{H1etH4} (resp$.$ Remark \ref{moduledeYBis}) by assuming that  $D$ is an open and  convex set possibly unbounded by maintaining the same hypotheses on $X$.
That is, assuming that $X: \Omega \times D \subset \Omega\times \reels^{d} \rightarrow \reels^{j}$ ($j \le d$) is a random field which belongs to $C^{1}(D,\reels^{j})$, such that for almost all $\omega \in \Omega$, the process $\nabla\!X(\omega)$ is Lipschitz with Lipschitz constant $L_{X}(\omega)$ satisfying $\bbE\left[L_X^d(\bm{\cdot})\right] < \infty$ (resp$.$
there exists $ L>0$ such that for almost all $\omega \in \Omega$, the assumption  $\supp(Y(\omega)) \neq \emptyset$ implies  that $L_{X}(\omega) \le L$).
In this case, we will assume that $Y$ is defined on $D_1$, a bounded and open set included in $D$.
Moreover, we will have to adapt the assumptions for $Y$ to the open set $D_1$ instead of $D$ and  to $X_{|D_1}$.
In this way, the assumptions ${\bH _2}$ and ${\bH _{4}}$ will still hold for $X_{|D_1}$ and $Y$ defined on  $D_1$.
\hfill$\bullet$
\end{rema}
\spacebefore
\begin{proofarg}{Proof of Remark \ref{moduleconvexe}} 
We prove this remark in the same way as Proposition \ref{H1etH4}.
As in the proof of  Theorem \ref{continueBis}, we use that the open set $D_1$ is contained in $D$ which is convex.
This allows us to apply, as on page \pageref{convexe}, Rolle's theorem and that $\nabla\!X$ is Lipschitz on $D$.
\end{proofarg}\\
In the following, we will study the hypotheses ${\bH _1}$, ${\bH _3}$ and ${\bH _{5}}$.
The last two hypothesis deal with the continuity of the right-hand side term in the Kac-Rice formula.
We will show the following proposition which is also deeply inspired by Caba\~na \cite{AL121221}.\\
We will exhibit a class of processes $X$ and $Y$ satisfying these assumptions.
We will therefore state some hypotheses concerning the processes $X$ and $Y$.\\
To the first three assumptions $\bA_1$, $\bA_2$ and $\bA_3$, we will add the assumption that $Y$ can be written as a function $G$ of $X$, $\nabla\!X$ and $W: \Omega \times D \subset \reels^{d} \rightarrow \reels^k$, $k \in \naturels^{\star}$,  a continuous random field, where $D$  is still an open set of $\reels^{d}$. That is, $\forall x \in D$
\begin{equation}
\label{egalite Y}
Y(x):= G(x, W(x), X(x), \nabla\!X(x)),
\end{equation}
where
$$
\begin{array}{r@{~}c@{~}l}
	G:  D \times \reels^k \times \reels^j \times \frakL(\reels^{d}, \reels^{j}) & \longrightarrow & \reels \\
	(x, z, u, A)  & \longmapsto & G(x, z, u, A),
\end{array}
$$%
is a continuous function of its variables on
$D \times \reels^k \times \reels^j \times  \frakL(\reels^{d}, \reels^{j})$
such that $\forall (x, z, u, A) \in D \times \reels^k \times \reels^j \times \frakL(\reels^{d}, \reels^{j})$,
$$
 \abs{G(x, z, u, A)} \le P(f(x), \normp[k]{z}, h(u), \normp[j, d]{A}),
$$%
where $P$ is a polynomial  with positive coefficients and $f: D \longrightarrow \reels^{+}$ and $h: \reels^j \longrightarrow \reels^{+}$ are continuous functions.
\begin{itemize}
\item $\bA_1$: The process $X: \Omega \times D \subset \Omega\times \reels^{d} \rightarrow \reels^{j}$ ($j \le d$) is Gaussian and belongs to $C^{1}(D,\reels^{j})$, such that $\exists a\in\reels$, $0 <a$, such that for almost all $x \in D$, $0<a \le \inf_{\normp[j]{z}=1}\normp[j]{\bbV(X(x)) \times z}$.
Moreover, the first order partial derivatives of its covariance $\Gamma_{X}$ are bounded almost surely on the diagonal contained in $D \times D$.
Therefore, for almost all $x \in D$, the random field $W(x)$ is independent of the vector $(X(x), \nabla\!X(x))$, and $\forall p \in \naturels$, $\forall n \in \naturels$, $\forall \ell \in \naturels$ and $\forall m \in \naturels$
$$
\int_{D} f^{p}(x)\bbE\!\left[\normp[k]{W(x)}^{n}\right]\bbE\!\left[\normp[j, d]{\nabla\!X(x)}^{\ell}\right]\bbE\!\left[\normp[j]{X(x)}^{m}\right] \ud x< \infty.
$$%
\item $\bA_2$: For all $ x \in D$, $X(x)=F(Z(x))$, where $F: \reels^j \longrightarrow \reels^j$ is a bijection belonging to $C^{1}(\reels^j,\reels^{j})$, such that $\forall z \in \reels^j$, the  Jacobian of  $F$ at $z$,   $J_{F}(z)$ satisfies $J_{F}(z)\neq 0$ and the function $F^{-1}$ is continuous.
The process $Z: \Omega \times D \subset \Omega\times \reels^{d} \rightarrow \reels^{j}$ ($j \le d$) is Gaussian and belongs to $C^{1}(D,\reels^{j}) \such \exists a\in\reels$, $a>0$, such that for almost all $x \in D$, $0<a \le \inf_{\normp[j]{z}=1}\normp[j]{\bbV(Z(x)) \times z}$; the first order partial derivatives of its covariance $\Gamma_{Z}$ are bounded almost surely on the diagonal contained in $D \times D$.
Moreover, for almost any $x \in D$, $W(x)$ is independent of the vector $(Z(x), \nabla\!Z(x))$, and $\forall p \in \naturels$, $\forall n \in \naturels$, $\forall \ell \in \naturels$ and $\forall m \in \naturels$
$$
\int_{D} f^{p}(x)\bbE\!\left[\normp[k]{W(x)}^{n}\right]\bbE\!\left[\normp[j, d]{\nabla\!Z(x)}^{\ell}\right]\bbE\!\left[\normp[j]{Z(x)}^{m}\right]\ud x< \infty.
$$%
\item $\bA_3$: For all $ x \in D$, $X(x)=F(Z(x))$, where $Z: \Omega \times D \subset \Omega\times \reels^{d} \rightarrow \reels^{j^{\prime}}$ ($j < j^{\prime}$) is Gaussian and belongs to $C^{1}(D,\reels^{j^{\prime}})$, with mean $m_{Z}(\bm{\cdot})=\bbE[Z(\bm{\cdot})]$ bounded on $D$, and such that there exist real numbers $a$ and $b$, $0 <a \le b$, such that for almost all $x \in D$,
$$
0<a \le \inf_{\normp[j]{z}=1}\normp[j]{\bbV(Z(x)) \times z} \le \sup_{\normp[j]{z}=1}\normp[j]{\bbV(Z(x)) \times z} \le b;
$$%
the first order partial derivatives of its covariance $\Gamma_{Z}$ are bounded almost surely on the diagonal contained in $D \times D$.
Moreover, for almost all $x \in D$, $W(x)$ is independent of the vector $(Z(x), \nabla\!Z(x))$ and it is assumed that the latter has a density denoted by ${p}_{Z(x), \nabla\!Z(x)}(\bm{\cdot}, \bm{\cdot})$.
Finally, $\forall p \in \naturels$, $\forall n \in \naturels$ and $\forall \ell \in \naturels$ 
\begin{equation}
\label{majoration moments}
\int_{D} f^{p}(x)\bbE\!\left[\normp[k]{W(x)}^{n}\right]\bbE\!\left[\normp[j^{\prime}, d]{\nabla\!Z(x)}^{\ell}\right] \ud x< \infty.
\end{equation}
The function $F$ must satisfy the hypothesis $\bF$ which is
\begin{itemize}
\item ($\bF$) $F: \reels^{j^{\prime}} \longrightarrow \reels^j$ ($j < j^{\prime}$) is $C^{2}(\reels^{j^{\prime}},\reels^{j})$.
In addition, defining
$A_{j^{\prime}}:=\{1,2,\ldots,j^{\prime}\}$,
there exists 
$\lambda:=(\ell_1, \ell_2, \dots, \ell_{j}) \in A_{j^{\prime}}^j$, $\ell_1 < \ell_2 < \dots < \ell_{j}$,
such that
$\forall z \in \reels^{j^{\prime}}$,
$$J_{F}^{(\lambda)}(z):=\det \!\left({
	\frac{\partial (F_1,\ldots,F_j)}{\partial(z_{\ell_1}, z_{\ell_2}, \ldots,z_{\ell_{j}})}(z)}\right)\neq 0.
$$\\
For simplicity, let us assume that  $\lambda=(1, 2, \dots, j)$ and denote $J_{F}(z)$ instead of $J_{F}^{(\lambda)}(z)$.\\
Furthermore, $\forall  v \in \reels^{j^{\prime}-j}$, the function $F_{v}$ defined by
$$
\begin{array}{r@{~}c@{~}l}
	F_{v}:   \reels^j  & \longrightarrow &  \reels^{j} \\
	u  & \longmapsto & F_{v}(u):=F(u, v),

\end{array}
$$%
is an invertible function whose inverse $F_{v}^{-1}$ is assumed to be continuous at $u$.\\
Moreover, $\forall \ell \in \naturels$ and $\forall \mu >0$, the function
\begin{equation}
\label{Tl}
\begin{aligned}
T_{\ell}:   \reels^j  & \longrightarrow \reels^{+} \\
u  & \longmapsto T_{\ell}(u)\\:=\int_{\reels^{j^{\prime}-j}} \frac{1}{\abs{J_{F}(F_{z}^{-1}(u), z)}} & e^{-\mu \normp[j^{\prime}-j]{z}^{2}} \normp[j, j^{\prime}]{\nabla F(F_{z}^{-1}(u), z)}^{\ell} \ud z,
\end{aligned}
\end{equation}
is continuous.
\end{itemize}
\item $ \bA_4$: For almost all $(x, y, \dot{x}) \in D \times \reels \times \reels^{dj}$ and $\forall u \in \reels^{j}$, the density ${p}_{Y(x), X(x), \nabla\!X(x)}(y, u, \dot{x})$ of the joint distribution of $(Y(x), X(x), \nabla\!X(x))$ exists and is continuous at $u$.\\
Moreover
$$
	u \longmapsto \int_{D} \int_{\reels \times \reels^{dj}}
		\abs{y} \normp[dj]{\dot{x}}^{j} {p}_{Y(x), X(x), \nabla\!X(x)}(y, u, \dot{x})
		 \ud \dot{x} \ud y  \ud x,
$$%
is continuous.
\end{itemize}
\spacebefore
\begin{rema}
\label{CD}
\begin{itemize}
\item It is interesting to note that  hypothesis $\bA_1$ (resp$.$ $\bA_2$, resp$.$ $ \bA_3$)  contains the case where the processes  $X$ and $Y$ satisfy $\forall x \in D$, $Y(x)= G(x, X(x), \nabla\!X(x))$ and also the case where $Y(x)$ is independent of  $(X(x), \nabla\!X(x))$ (resp$.$
$Y(x)$ independent of $\left(Z(x), \nabla\!Z(x)\right)$).
\item Note also that hypothesis $ \bA_4$ is satisfied for example in the case where, $\forall u \in \reels^{j}$, there exists a neighborhood $V_{u}$ of $u$ and a function $h_{u}$ such that 
$$
	\int_{D}  \int_{\reels \times \reels^{dj}} \abs{y} \normp[dj]{\dot{x}}^{j} h_{u}(x, y, \dot{x}) \ud \dot{x}  \ud y \ud x< \infty,
$$%
and such that $\forall z \in V_{u}$ and for almost all $(x, y, \dot{x}) \in D \times \reels \times \reels^{dj}$, ${p}_{Y(x), X(x), \nabla\!X(x)}(y, z, \dot{x}) \le h_{u}(x, y, \dot{x}).$
In fact, these hypotheses are the required ones to apply Lebesgue's dominated convergence theorem which allows us to obtain the continuity of the function 
$$
	u \longmapsto \int_{D} \int_{\reels \times \reels^{dj}} \abs{y} \normp[dj]{\dot{x}}^{j} {p}_{Y(x), X(x), \nabla\!X(x)}(y, u, \dot{x})  \ud \dot{x} \ud y \ud x.
$$%
\end{itemize}
\hfill$\bullet$
\end{rema}
We are now able to exhibit a class of processes $X$ and $Y$ satisfying the hypotheses ${\bH _1}$, ${\bH _3}$ and ${\bH _{5}}$ through the following proposition.
\spacebefore
\begin{prop}
\label{H2etH5}
If $Y$ satisfies (\ref{egalite Y}) and if  $X$ and  $Y$ satisfy one of the three assumptions  $\bA_1$, $\bA_2$, $\bA_3$ or if $X$ and $Y$ satisfy assumption $ \bA_4$, 
then the hypotheses ${\bH _1}$, ${\bH _3}$ and ${\bH _{5}}$ are verified.
\end{prop}
\spacebefore
\begin{proofarg}{Proof of Proposition \ref{H2etH5}} 
\begin{enumerate}
\item\label{itm:p3.1.12.1} Let us first assume that the processes $X$ and $Y$ satisfy assumption $\bA_1$.
Let us show that the hypotheses ${\bH _1}$ and ${\bH _{5}}$ are satisfied.\\
Since $X$ is Gaussian and that for almost all  $x \in D$,
$$
	\inf_{\normp[j]{z}=1} \normp[j]{\bbV(X(x)) \times z} \ge a >0,
$$%
the distribution of the vector $X(x)$ is not singular with density ${p}_{X(x)}(\bm{\cdot})$.
Moreover, $u \mapsto  {p}_{X(x)}(u)$ is continuous and it is bounded above, \ie  there exists a real $M$ such that for almost all $x \in D$ and $\forall u \in \reels^{j}$,  
\begin{equation}
\label{majorationp}
{p}_{X(x)}(u) \le M.
\end{equation}
The hypothesis ${\bH _1}$ is verified.\\
Let us show that the hypothesis ${\bH _{5}}$ is satisfied.\\
Since for almost all $x \in D$,
$$
Y(x)= G(x, W(x), X(x), \nabla\!X(x)), 
$$%
using the assumptions on $G$ and since for all $A \in \frakL(\reels^{d}, \reels^{j})$, $H(A) \le \bC  \normp[j, d]{A}^{j}$ , $\forall u \in \reels^{j}$ we have
\begin{align*}
\MoveEqLeft[1] {\bbE[Y(x)H(\nabla\!X(x))\given{X(x)=u}]}=\bbE[L(x, W(x), X(x), \nabla\!X(x))\given{X(x)=u}],
\end{align*}
where $L$ is a continuous function of all its variables belonging to
$D \times \reels^k \times \reels^j \times \frakL(\reels^{d}, \reels^{j})$ and such that $\forall (x, z, u, A) \in D \times \reels^k \times \reels^j \times \frakL(\reels^{d}, \reels^{j})$,
\begin{equation}
 	\label{majorationL}
 	\abs{L(x, z, u, A)} \le Q(f(x), \normp[k]{z}, h(u), \normp[j,d]{A}),
 \end{equation}
where $Q$ is a polynomial with positive coefficients and  $f: D \longrightarrow \reels^{+}$ and $h: \reels^j \longrightarrow \reels^{+}$ are continuous functions.\\
For almost all $x$ fixed in $D$, let us consider the regression equations: for $s \in D$
 \begin{align}
 \label{regression}
 X(s) &= \alpha(s) X(x) + \xi(s), \nonumber\\
 \nabla\!X(s) &= \sum_{i=1}^{j} \nabla \alpha_{i}(s) X_{i}(x) + \nabla \xi(s),
 \end{align}
 where $(\xi(s), \nabla \xi(s))$ is a Gaussian vector independent of $X(x)$.\\
 In particular, $\alpha(x)= Id_{j}$.\\
 A covariance computation gives 
$$
	\alpha(s) = \Gamma_{X}(s, x) \times \Gamma_{X}^{-1}(x, x),
$$%
 where $\Gamma_{X}$ is the covariance matrix of $X$.\\
Thus for all $i, m =1, \dots, j$ and $\ell=1, \dots, d$,
$$
 \left (\nabla \alpha_{i}(s)\right)_{\ell, m}= \!\left(\frac{\partial \Gamma_{X}}{\partial s_{\ell}}(s, x) \times \Gamma_{X}^{-1}(x, x)\right)_{mi}.
$$%
In particular for almost all $x \in D$, $i, m =1, \dots, j$ and $\ell=1, \dots, d$,
$$
 	\left (\nabla \alpha_{i}(x)\right)_{m, \ell}= \!\left(\frac{\partial \Gamma_{X}}{\partial s_{\ell}}(x, x) \times \Gamma_{X}^{-1}(x, x)\right)_{mi}.
$$%
Since for almost all $x \in D$, $\inf_{\normp[j]{z}=1} \normp[j]{\bbV(X(x)) \times z} \ge a >0$ and that the first order partial derivatives of the covariance $\Gamma_{X}$ are almost surely bounded above on the diagonal contained in $D \times D$,  $\exists M \in \reels$ such that for all $i=1, \dots, j$ and for almost all $x \in D$ we have
\begin{equation}
 \label{majorationalpha}
 \normp[j,d]{\nabla \alpha_{i}(x)} \le M.
 \end{equation}
 For $u \in \reels^{j}$, let
$$
 G_{X, L}(u):= \int_{D}\bbE\!\left[Y(x)H(\nabla\!X(x))\given{X(x)=u}\right] {p}_{X(x)}(u) \ud x.
$$%
With the above notations, we obtain
$$
G_{X, L}(u)= \int_{D}\bbE\!\left[L(x, W(x), X(x), \nabla\!X(x))\given{X(x)=u}\right] {p}_{X(x)}(u) \ud x.
$$%
Since for almost all $x \in D$ the random variable $W(x)$ is independent of the vector $(X(x), \nabla\!X(x))$, using (\ref{regression}), this yields that
\begin{align*}
 &\qquad G_{X, L}(u)\\
 &\qquad =\int_{D}\bbE\!\left[L\left(x, W(x), u, \sum_{i=1}^{j} \nabla \alpha_{i}(x) (u_{i} -X_{i}(x))+ \nabla\!X(x)\right)\right]
 p_{X(x)}(u) \ud x.
\end{align*}
We have therefore eliminated the conditioning appearing in the integrand.
Now since for almost all $x \in D$, the function $u \mapsto {p}_{X(x)}(u)$ is continuous and the function $L$ is also continuous then for almost all $x \in D$, the function
$$
u \mapsto
{L\!\left(x, W(x), u, \sum_{i=1}^{j} \nabla \alpha_{i}(x)(u_{i}-X_{i}(x)) + \nabla\!X(x)\right)} {p}_{X(x)}(u),
$$%
is continuous.\\
Moreover, using the bounds (\ref{majorationp}), (\ref{majorationL}) and (\ref{majorationalpha}), we get that for almost all $x \in D$,
\begin{align*}
 &\abs{L\!\left(x, W(x), u, \sum_{i=1}^{j} \nabla \alpha_{i}(x)(u_{i}-X_{i}(x)) + \nabla\!X(x)\right)} {p}_{X(x)}(u)\\
 &\qquad \le S\!\left(f(x), \normp[k]{W(x)}, \ell(u), \normp[j]{X(x)}, \normp[j, d]{\nabla\!X(x)}\right),
\end{align*}
where $S$ is also a polynomial with positive coefficients and $\ell: \reels^j \longrightarrow \reels^{+}$ is a continuous function.\\
 It is clear that for almost all  $x \in D$, the function
 $$
 u \mapsto S\!\left(f(x), \normp[k]{W(x)}, \ell(u), \normp[j]{X(x)}, \normp[j, d]{\nabla\!X(x)}\right),
 $$%
is continuous.
Furthermore, we know that for $\forall p \in \naturels$, 
$\forall n \in \naturels$, $\forall \ell \in \naturels$ and  $\forall m \in \naturels$,
$$
\int_{D} f^{p}(x)\bbE\!\left(\normp[k]{W(x)}^{n}\right)\bbE\!\left(\normp[j, d]{\nabla\!X(x)}^{\ell}\right) \bbE\!\left(\normp[j]{X(x)}^{m}\right)\ud x< \infty.
$$%
Thus, and since the function $u \mapsto \ell(u)$ is continuous we obtain that the function
$$
u \mapsto
\int_{D}\bbE\!\left[S\left(f(x), \normp[k]{W(x)}, \ell(u), \normp[j]{X(x)}, \normp[j,d]{\nabla\!X(x)}\right)\right] \ud x,
$$%
is continuous.\\
A weak application of Lebesgue's dominated convergence theorem allows to conclude that the hypothesis ${\bH _{5}}$ is true.\\
 A similar proof can be made to show that the hypothesis  ${\bH _3}$ is also satisfied.
This completes the first part of the proof.
 \item Suppose now that the processes $X$ and $Y$ satisfy assumption $\bA_2$.
Let us prove that hypothesis ${\bH _1}$ is satisfied.\\
In the same way as in part \ref{itm:p3.1.12.1} of the proof, since $Z$ is Gaussian and given that for almost all $x \in D$,
$$
	\inf_{\normp[j]{z}=1} \normp[j]{\bbV(Z(x)) \times z} \ge a >0,
$$%
the distribution of the vector 
 $Z(x)$ is non singular with density ${p}_{Z(x)}(\bm{\cdot})$.

The assumptions on the function $F$ imply that for almost all  $x \in D$ the vector $X(x)$ has a density ${p}_{X(x)}(\bm{\cdot})$ given by $\forall u \in \reels^{j}$:
$$
{p}_{X(x)}(u)= \frac{1}{\abs{J_{F}(F^{-1}(u))}} {p}_{Z(x)}(F^{-1}(u)).
$$%
Let us show that the hypotheses ${\bH _3}$ and ${\bH _{5}}$ are satisfied.\\
Using the same notations of part \ref{itm:p3.1.12.1}, for almost all $x \in D$ and $\forall u \in \reels^{j}$ we have
\begin{align*}
\MoveEqLeft[1]{\bbE\left[L(x, W(x), X(x), \nabla\!X(x))\given{X(x)=u}\right] {p}_{X(x)}(u)}\\
 &=\bbE\left[L(x, W(x), F(Z(x)), \nabla F(Z(x)) \times \nabla\!Z(x))\given{Z(x)=F^{-1}(u)}\right]\\
 &\rule{144pt}{0pt}\times{p}_{Z(x)}(F^{-1}(u)) \times  \frac{1}{\abs{J_{F}(F^{-1}(u))}}\\
 &=\bbE\left[\widetilde{L}(x, W(x), Z(x), \nabla\!Z(x))\given{Z(x)=F^{-1}(u)}\right]\\
 &\rule{144pt}{0pt}\times {p}_{Z(x)}(F^{-1}(u)) \times  \frac{1}{\abs{J_{F}(F^{-1}(u))}},
\end{align*}
where the function $\widetilde{L}$,
$\forall (x, z, u, A) \in D \times \reels^k \times \reels^j \times \frakL(\reels^{d}, \reels^{j})$
is defined by
$$
	 \widetilde{L}(x, z, u, A):= L(x, z, F(u), \nabla F(F^{-1}(u)) \times A).
 $$%
It is clear, since $F$ belongs to $C^{1}(\reels^j,\reels^{j})$ and $F^{-1}$ is continuous, that  $\widetilde{L}$ has the properties of $L$, \ie  $\widetilde{L}$ is a continuous function of its variables in 
$D \times \reels^k \times \reels^j \times \frakL(\reels^{d}, \reels^{j})$ and that
$$
	\abs{\widetilde{L}(x, z, u, A)} \le \widetilde{Q}(f(x), \normp[k]{z}, \widetilde{h}(u), \normp[j,d]{A}),
 $$%
where $\widetilde{Q}$ is a polynomial with positive coefficients and  $\widetilde{h}: \reels^j \longrightarrow \reels^{+}$ is a continuous function.\\
We have shown using the notations of 1) that $\forall u \in \reels^{j}$
$$
 	G_{X, L}(u)=  G_{Z, \widetilde{L}}(F^{-1}(u))  \times \frac{1}{\abs{J_{F}(F^{-1}(u))}}.
 $$%
This leads us to the case considered in 1) where the process $X$ is replaced by the process $Z$.
The continuity of the function $u \mapsto G_{X, L}(u)$ is a consequence of the continuity of $G_{Z, \widetilde{L}}$ and the fact that the function $F^{-1}$ is continuous and that $F$ is $C^{1}(\reels^j,\reels^{j})$.
This ends the second part of the proof.
\item  Suppose that the processes  $X$ and $Y$ satisfy the assumption $\bA_3$.
We need to prove that hypotheses ${\bH _1}$, ${\bH _3}$ and ${\bH _{5}}$ are satisfied.\\
Let us first prove that for almost all  $x \in D$, the distribution of the vector $(X(x), \nabla\!X(x))$ has a density ${p}_{X(x), \nabla\!X(x)}(\bm{\cdot}, \bm{\cdot})$, and let us compute this density.\\
Consider the following notations.\\
The matrix $s_{v,u}= (s_{ik})_{\substack{1 \le i \le v \\ 1 \le k \le u}} \in \frakL(\reels^{u}, \reels^{v})$, defined by its generic element $s_{ik}$, will be identified with the matrix $s_{v,u}$
with row vector  $s_{(v,u)} \in \reels^{vu}$, defined by
$$
s_{(v,u)}:= (s_{11}, s_{21}, \dots, s_{v1}, s_{12}, s_{22}, \dots, s_{v2}, \dots, s_{1u}, s_{2u}, \dots, s_{vu}).
$$%
Using this notation, we can introduce the following function
\begin{align*}
\MoveEqLeft[2]{K:  \reels^{j} \times \reels^{j^{\prime}-j} \times \reels^{jd}  \times \reels^{(j^{\prime}-j)d}  \longrightarrow \reels^{j} \times \reels^{j^{\prime}-j} \times \reels^{jd}  \times \reels^{(j^{\prime}-j)d}}\\
	&\left({t=t_{(j^{\prime},1)}= \begin{pmatrix} t_{j,1} \\ t_{j^{\prime}-j,1}
\end{pmatrix}, s=s_{(j^{\prime},d)}= \begin{pmatrix} s_{j,d} \\ s_{j^{\prime}-j,d}
\end{pmatrix}}\right)\\
	&\specialpos{\hfill\qquad\longmapsto K(t, s):=
\left({(F(t))_{(j,1)}, t_{(j^{\prime}-j,1)}, \left(\nabla F(t) \times s\right)_{(j,d)}}, s_{(j^{\prime}-j,d)}\right)}.
\end{align*}
The Jacobian matrix  $J_{K}$ of this transformation satisfies: $\forall (t, s) \in \reels^{j^{\prime}} \times \reels^{j^{\prime}d}$:
$$
J_{K}(t, s)= (J_{F}(t))^{d+1} \neq 0,
$$%
by hypothesis.\\
Furthermore, since $F$ is $C^{2}(\reels^{j^{\prime}},\reels^j)$ then $K$ is $C^{1}(\reels^{j^{\prime}} \times \reels^{j^{\prime}d}, \reels^{j^{\prime}} \times \reels^{j^{\prime}d})$.
Moreover, $K$ is one-to-one and has an inverse  $K^{-1}$ given by
\begin{align*}
\MoveEqLeft[1]{K^{-1}:  \reels^{j} \times \reels^{j'-j} \times \reels^{jd}  \times \reels^{(j'-j)d} \longrightarrow \reels^{j} \times \reels^{j'-j} \times \reels^{jd}  \times \reels^{(j'-j)d}}\\
	&\left({t= \begin{pmatrix} t_{j,1} \\
t_{j^{\prime}-j,1}
\end{pmatrix}, s= \begin{pmatrix} s_{j,d} \\
s_{j^{\prime}-j,d}
\end{pmatrix}}\right)\\
	& \specialpos{\hfill\longmapsto 
\left({F^{-1}_{t_{j^{\prime}-j,1}}(t_{j,1}), t_{j^{\prime}-j,1}, }
\left[{\nabla F(F^{-1}_{t_{j^{\prime}-j,1}}(t_{j,1}), t_{j^{\prime}-j,1})}\right]_{jj}^{-1}\right.}\\
& \specialpos{\hfill\rule{14pt}{0pt}\times\left.\left({ s_{j,d} -  \left[{\nabla F(F^{-1}_{t_{j^{\prime}-j,1}}(t_{j,1}), t_{j^{\prime}-j,1})}\right]_{jj^{\prime}-j} \times s_{j^{\prime}-j,d}}\right), s_{j^{\prime}-j,d}\right),}
\end{align*}
where for $A \in \frakL(\reels^{j^{\prime}}, \reels^{j})$, we denote, by $\left[{A }\right]_{jj}$ the matrix $A$ for which we keep only  the $j$ first columns and by
$\left[{A }\right]_{jj^{\prime}-j}$ the matrix $A$ for which we retain the $(j^{\prime}-j)$ last columns.\\
Since $\forall x \in D$
$$
	X(x)= F(Z(x))\mbox{ and }\nabla\!X(x) = \nabla F(Z(x)) \times \nabla\!Z(x),
$$ %
we have
$$
	K(Z(x), \nabla\!Z(x))= \!\left({X(x), (Z(x))_{j^{\prime}-j,1}, \nabla\!X(x), (\nabla\!Z(x))_{j^{\prime}-j, d} }\right).
$$%
We deduce that if for almost all $x \in D$,
$${p}_{X(x), (Z(x))_{j^{\prime}-j,1}, \nabla\!X(x),  (\nabla\!Z(x))_{j^{\prime}-j, d}}(\bm{\cdot}, \bm{\cdot}, \bm{\cdot}, \bm{\cdot})
$$%
 denotes the density of the vector
$$
\!\left({X(x), (Z(x))_{j^{\prime}-j,1}, \nabla\!X(x),  (\nabla\!Z(x))_{j^{\prime}-j, d}}\right),
$$%
and ${p}_{Z(x), \nabla\!Z(x)}(\bm{\cdot}, \bm{\cdot})$
that of $\!\left({Z(x),\nabla\!Z(x)}\right)$ then
$$
\forall (u, z_{j^{\prime}-j,1}, s_{j, d}, s_{j^{\prime}-j,d}) \in \reels^{j} \times \reels^{j^{\prime}-j} \times \reels^{jd} \times \reels^{(j^{\prime}-j)d},
$$%
we have
\begin{align*}
\MoveEqLeft[2]{\rule{20pt}{0pt}{p}_{X(x), (Z(x))_{j^{\prime}-j,1}, \nabla\!X(x),  (\nabla\!Z(x))_{j^{\prime}-j, d}}(u, z_{j^{\prime}-j,1}, s_{j, d}, s_{j^{\prime}-j,d})}\\
	&= {p}_{Z(x), \nabla\!Z(x)}(K^{-1}(u, z_{j^{\prime}-j,1}, s_{j, d}, s_{j^{\prime}-j,d})) \times 
 \frac{1}{\abs{J_{F}(F^{-1}_{z_{j^{\prime}-j,1}}(u), z_{j^{\prime}-j,1})}^{d+1}}.
\end{align*}
Finally, we get the density of the vector $\!\left({X(x),\nabla\!X(x)}\right)$ by integrating this last expression.
Hence for almost all $x \in D$ and $\forall (u, s_{j, d}) \in \reels^{j} \times \reels^{jd}$ we have
\begin{multline}
\label{Co}
{p}_{X(x), \nabla\!X(x)}(u, s_{j,d})
= \int_{\reels^{j^{\prime}-j}\times\reels^{(j^{\prime}-j)d}}
 	\frac{1}{\abs{J_{F}(F_{z}^{-1}(u), z)}^{d+1}}\\
	\times{p}_{Z(x),\nabla\!Z(x)} 
	\!\left(\rule{0pt}{14pt}
	(F_{z}^{-1}(u), z), 
		\!\left(\rule{0pt}{12pt}
		\left[\nabla F(F_{z}^{-1}(u), z)\right]_{jj}^{-1} \right.\right. \rule{72pt}{0pt}\\ \left.\left.
			\times\!\left(
				s_{j,d} -  \left[\nabla F(F_{z}^{-1}(u), z)\right]_{j j^{\prime}-j}s_{j^{\prime}-j, d}
			\right),s_{j^{\prime}-j,d}
		\rule{0pt}{12pt}\right)
	\rule{0pt}{14pt}\right)
  \ud s_{j^{\prime}-j,d} \ud z.
\end{multline}
\spacebefore
\begin{rema}
It is important to note that the results presented above could be obtained using the coarea formula.
We referred the reader to Corollary  4.18 of \cite[p.\,68]{AL121221c}.
However, we preferred explicit computations in order to obtain the exact expression of this density and also to introduce some useful notations in what follows.
\hfill$\bullet$
\end{rema}
\spacebefore
Now, for  $u \in \reels^{j}$, we define as in the part \ref{itm:p3.1.12.1}
\begin{align*}
\displaystyle G_{X, L}(u)
	& :=  \int_{D}\bbE[Y(x)H(\nabla\!X(x))\given{X(x)=u}] {p}_{X(x)}(u) \ud x\\
	& = \int_{D}\bbE[L(x, W(x), X(x), \nabla\!X(x))\given{X(x)=u}] {p}_{X(x)}(u) \ud x.
\end{align*}
Since for almost all $x \in D$, $W(x)$ is independent of 
$
	\!\left({X(x), \nabla\!X(x})\right),
$
we get $\forall u \in \reels^{j}$,
 \begin{align*}
\MoveEqLeft[0]{G_{X, L}(u)=\int_{D\times\reels^{jd}}\bbE[L(x, W(x), u, s_{j,d})]
			{p}_{X(x), \nabla\!X(x)}(u, s_{j,d}) \ud s_{j,d}  \ud x}\\
	&=  \int_{D} \int_{\reels^{j^{\prime}-j}\times\reels^{(j^{\prime}-j)d} \times\reels^{jd}}
		 \frac{1}{\abs{J_{F}(F_{z}^{-1}(u), z)}^{d+1}}\times\bbE[L(x, W(x), u, s_{j,d})]\\
	&  \rule{4pt}{0pt}\times 
   {p}_{Z(x),\nabla\!Z(x)}
   	\!\left(\rule{0pt}{14pt}
		\!\left(F_{z}^{-1}(u), z\right), 
	\!\left(\rule{0pt}{12pt}
		\left[\nabla F(F_{z}^{-1}(u), z)\right] _{jj}^{-1}\right.\right.\\
		&\left.\rule{0pt}{14pt}\left.\rule{0pt}{12pt}
		\times\!\left(
			s_{j,d} -  [\nabla F(F_{z}^{-1}(u), z)]_{j j^{\prime}-j} s_{j^{\prime}-j, d}
  		\right), s_{j^{\prime}-j,d}
  	\right)
  	\right) \ud s_{j,d} \ud s_{j'-j,d} \ud z \ud x.
\end{align*}
In the last integral where the  integration domain is  $\reels^{jd}$, holding fixed $(s_{j'-j,d}, z, x) \in \reels^{(j^{\prime}-j)d} \times \reels^{j^{\prime}-j} \times D$,  let us make the following change of variable
$$
v_{j,d}:=
 [
  \nabla F(F_{z}^{-1}(u), z) ]_{jj}^{-1} \times 
  (s_{j,d} -  \left[\nabla F(F_{z}^{-1}(u), z)\right]_{j j^{\prime}-j} \times s_{j^{\prime}-j, d}).
$$%
We get
\makeatletter
\newcommand{\pushright}[1]{\ifmeasuring@#1\else\omit\hfill$\displaystyle#1$\fi\ignorespaces}
\newcommand{\pushleft}[1]{\ifmeasuring@#1\else\omit$\displaystyle#1$\hfill\fi\ignorespaces}
\makeatother
\begin{align*}
\int_{\reels^{jd}} &\bbE[L(x, W(x), u, s_{j,d})]
	{p}_{Z(x),\nabla\!Z(x)} ( (F_{z}^{-1}(u), z),([\nabla F(F_{z}^{-1}(u), z)]_{jj}^{-1}\\
	&\specialpos{\hfill(s_{j,d} -  \left[\nabla F(F_{z}^{-1}(u), z)\right]_{j j^{\prime}-j} s_{j^{\prime}-j, d}),  s_{j^{\prime}-j,d})) \ud s_{j,d}}\\\
	&=\int_{\reels^{jd}}\bbE[L(x, W(x), u,[\nabla F(F_{z}^{-1}(u), z)]_{jj} \times v_{j,d}\\
  	&\specialpos{\hfill+ [\nabla F(F_{z}^{-1}(u), z)]_{j j^{\prime}-j} \times s_{j^{\prime}-j, d})]}\\
	&\specialpos{\hfill\times {p}_{Z(x),\nabla\!Z(x)}({(F_{z}^{-1}(u), z),({ v_{j,d}, s_{j^{\prime}-j,d}})}) \times \abs{J_{F}(F_{z}^{-1}(u), z)}^{d} \ud v_{j,d}}\\
	& =\int_{\reels^{jd}}\bbE\!\left[L\!\left({x, W(x), u, \nabla F(F_{z}^{-1}(u), z) 
   		 \begin{pmatrix} v_{j,d} \\ s_{j^{\prime}-j,d} \end{pmatrix}
   		}\right)\right] \\
	&\specialpos{\hfill\times  {p}_{Z(x), \nabla\!Z(x)}\!\left({(F_{z}^{-1}(u), z),
 		(v_{j,d}, s_{j^{\prime}-j,d})}\right)
		\times \abs{J_{F}(F_{z}^{-1}(u), z)}^{d}
 \ud v_{j,d}}.
\end{align*}

Finally, $\forall u \in \reels^{j}$,
\begin{align}
G_{X, L}(u)&=  \int_{D} \int_{\reels^{j^{\prime}-j}}  \frac{1}{\abs{J_{F}(F_{z}^{-1}(u), z)}}\nonumber\\
 &\specialpos{\hfill\times \int_{\reels^{(j^{\prime}-j)d} \times \reels^{jd}}\bbE\left[L\!\left({x, W(x), u, 
  \nabla F(F_{z}^{-1}(u), z) 
   \times \begin{pmatrix} v_{j,d} \\ s_{j^{\prime}-j,d} \end{pmatrix}
   }\right)\right]}\nonumber \\
&\specialpos{\hfill\times{p}_{Z(x),\nabla\!Z(x)}\!\left({(F_{z}^{-1}(u), z),
 \begin{pmatrix} v_{j,d}\\
 s_{j^{\prime}-j,d}
 \end{pmatrix}}\right)  
 \ud v_{j,d} \ud s_{j^{\prime}-j,d} 
  \ud z  \ud x}\nonumber\\
  &=\int_{D} \int_{\reels^{j^{\prime}-j}} \frac{1}{\abs{J_{F}(F_{z}^{-1}(u), z)}}\nonumber\\
 &\specialpos{\hfill\times   \int_{\reels^{j^{\prime}d}}\bbE\!\left[L\!\left({x, W(x), u, 
  \nabla F(F_{z}^{-1}(u), z) 
   \times s_{j^{\prime},d} 
   }\right)\right]}\nonumber\\
&\specialpos{\hfill\times {p}_{Z(x),\nabla\!Z(x)}
 ((F_{z}^{-1}(u), z), s_{j^{\prime},d}) 
 \ud s_{j^{\prime},d}
 \ud z  \ud x.}\label{GXL}
\end{align}

A slight modification of the beginning of the proof or the use of the coarea formula given in \cite[Corollary 4.18, p.~68]{AL121221c}, shows that for almost all $x \in D$ the vector $X(x)$ has a density ${p}_{X(x)}(\bm{\cdot})$,
for all $u \in \reels^{j}$, given by:
$$
{p}_{X(x)}(u)=  \int_{\reels^{j^{\prime}-j}} \frac{1}{\abs{J_{F}(F_{z}^{-1}(u), z)}}
\times  {p}_{Z(x)}(F_{z}^{-1}(u), z)\ud z.
$$%
This implies the hypothesis ${\bH _1}$.
Moreover, this density is a continuous function of the variable $u$.\\
Indeed, using the assumptions on  $Z$,  $\exists \lambda >0$ such that for almost all $x \in D$ and for all $(z, u) \in  \reels^{j^{\prime}-j}\times \reels^{j}$, 
\begin{align}
\label{pZ}
{p}_{Z(x)}(F_{z}^{-1}(u), z)&\le \bC  e^{- 2\lambda \normp[j^{\prime}]{(F_{z}^{-1}(u), z)-m_{Z}(x)}^{2}}\nonumber\\
 &\le \bC  e^{- \lambda \normp[j^{\prime}]{(F_{z}^{-1}(u), z)}^{2}}\nonumber\\
 &\le \bC  e^{- \lambda \normp[j^{\prime}-j]{z}^{2}},
\end{align}
since the function $m_Z(\bm{\cdot})$ is bounded on $D$.\\
Moreover, using the assumptions satisfied by the function $F$ and the process $Z$, for almost all $z \in \reels^{j^{\prime}-j}$ and $x \in D$,  the functions
\begin{align*}
 u & \mapsto \frac{1}{\abs{J_{F}(F_{z}^{-1}(u), z)}} \times  {p}_{Z(x)}(F_{z}^{-1}(u), z),\\
 \shortintertext{and}
 u & \mapsto \frac{1}{\abs{J_{F}(F_{z}^{-1}(u), z)}} e^{- \lambda \normp[j^{\prime}-j]{z}^{2}}
\end{align*}
are continuous.\\
Using that the function  $u \mapsto T_{\ell}(u)$ is continuous for $\ell =0$ (see (\ref{Tl}) in the hypothesis $\bA_3$), Lebesgue's dominated convergence theorem allows to state that for almost all $x \in D$, the function $u \mapsto {p}_{X(x)}(u)$ is continuous.\\
Let us now return to the definition of $G_{X, L}(u)$ given by (\ref{GXL}).
Since for almost all $x \in D$, $W(x)$ is independent of $(Z(x), \nabla\!Z(x))$, $G_{X, L}(u)$ can be written as follows
\begin{align*}
\MoveEqLeft[0]{G_{X, L}(u)= \int_{D} \int_{\reels^{j^{\prime}-j}} \frac{1}{\abs{J_{F}(F_{z}^{-1}(u), z)}}
\times  {p}_{Z(x)}(F_{z}^{-1}(u), z)} \\
& \times \bbE\!\left[L(x, W(x),u,  \nabla F(F_{z}^{-1}(u), z) \times \nabla\!Z(x))\given Z(x)=(F_{z}^{-1}(u), z)\right]\ud z  \ud x.
\end{align*}
In the same way as in the part \ref{itm:p3.1.12.1} of this proof, for any $s \in D$ we regress  $Z(s)$ on $Z(x)$ for almost all $x \in D$ so
 \begin{align*}
 Z(s) &= \alpha(s) Z(x) + \xi(s), \\
 \nabla\!Z(s) &= \sum_{i=1}^{j^{\prime}} \nabla \alpha_{i}(s) Z_{i}(x) + \nabla \xi(s),
 \end{align*}
 where $(\xi(s), \nabla \xi(s))$ is a Gaussian vector independent of $Z(x)$.\\
Using the assumptions on the process $Z$, we obtain as in the part \ref{itm:p3.1.12.1} of this proof the following inequality:
$\exists M \in \reels \such \forall i \in\{1, \dots, j^{\prime}\}$ and for almost all $x \in D$ we have
\begin{equation}
  \normp[j^{\prime},d]{\nabla \alpha_{i}(x)}\le M.
   \label{majorationalphaBis}
 \end{equation}
Moreover, for almost all $x \in D$ and since  $W(x)$ is independent of $(Z(x), \nabla\!Z(x))$, we get
\begin{align*}
\MoveEqLeft[1]{G_{X, L}(u)= \int_{D} \int_{\reels^{j^{\prime}-j}} \frac{1}{\abs{J_{F}(F_{z}^{-1}(u), z)}} \times  {p}_{Z(x)}\!\left(F_{z}^{-1}(u), z\right)} \\
	&\specialpos{\hfill\times\bbE\!\left[\rule{0pt}{18pt}L\!\left(\rule{0pt}{18pt}x, W(x), u, \nabla F(F_{z}^{-1}(u), z) \right. \right.\hfill}\\
	& \specialpos{\hfill\rule{60pt}{0pt}\left. \left. \times\!\left\{{\sum_{i=1}^{j^{\prime}} \nabla \alpha_{i}(x) \left[((F_{z}^{-1}(u), z))_{i} -Z_{i}(x)\right] + \nabla\!Z(x)}\right\}\right)\right] \ud z \ud x.}
\end{align*}
As in part \ref{itm:p3.1.12.1}, we have thus eliminated the conditioning appearing in the integrand and
 \begin{align*}
\MoveEqLeft[1]{G_{X, L}(u) = \int_{D} \int_{\reels^{j^{\prime}-j}} \int_{\Omega} \frac{1}{\abs{J_{F}(F_{z}^{-1}(u), z)}}
\times  {p}_{Z(x)}(F_{z}^{-1}(u), z)}  \\
	& \specialpos{\hfill\times L\!\left(\rule{0pt}{18pt}x, W(x)(\omega), u,  \nabla F(F_{z}^{-1}(u), z)\right.\hfill}\\
	& \left.\times\!\left\{{\sum_{i=1}^{j^{\prime}} \nabla \alpha_{i}(x)[ ((F_{z}^{-1}(u), z))_{i}-Z_{i}(x)(\omega)] + \nabla\!Z(x)(\omega)}\right\}\right)\ud \bbP(\omega) \ud z \ud x\\
	& = \int_{D} \int_{\reels^{j^{\prime}-j}} \int_{\Omega} f(u, \omega, z, x) \ud \bbP(\omega)\ud z  \ud x.
\end{align*}
By the hypotheses satisfied by $Z$ and $F$ and since  $L$ is a continuous function,  we obtain for almost all $(\omega, z, x) \in \Omega \times \reels^{j^{\prime}-j} \times D$, that the function
$u \mapsto f(u, \omega, z, x)$ is continuous.
Now let's bound the expression $f(u, \omega, z, x)$.\\
By using the bounds  (\ref{majorationL}), (\ref{pZ}) and (\ref{majorationalphaBis}) we get that for almost all  $(\omega, z, x) \in \Omega \times \reels^{j^{\prime}-j}\times D$,
\begin{multline*}
\abs{f(u, \omega, z, x)} \le \frac{\bC }{\abs{J_{F}(F_{z}^{-1}(u), z)}} e^{- \lambda \normp[j^{\prime}]{(F_{z}^{-1}(u), z)}^{2}} \\
\times R\!\left(f(x), \normp[k]{W(x)(\omega)}, h(u), \normp[j, j^{\prime}]{\nabla F(F_{z}^{-1}(u), z)},
 \right.\\ \left.
\normp[j^{\prime}]{(F_{z}^{-1}(u), z)}, \normp[j^{\prime},d]{\nabla\!Z(x)(\omega)}, 
 \normp[j^{\prime}]{Z(x)(\omega)}\right),
\end{multline*}
 where $R$ is a polynomial with positive coefficients and  $h: \reels^j \longrightarrow \reels^{+}$ is a continuous function.\\
 Using the fact that $\forall n \in \naturels$, $\exists M_{n}>0 \such \forall y \in \reels^{j^{\prime}}$,
$$
 e^{-\lambda/2 \normp[j^{\prime}]{y}^{2}}  \times \normp[j^{\prime}]{y}^{n} \le M_{n},
 $$%
 we get that for almost all  $(\omega, z, x) \in \Omega \times \reels^{j^{\prime}-j}\times D$,
 \begin{multline*}
 \abs{f(u, \omega, z, x)}\le \frac{\bC }{\abs{J_{F}(F_{z}^{-1}(u), z)}} e^{- \mu \normp[j^{\prime}-j]{z}^{2}}\times
S\!\left(\rule{0pt}{12pt}f(x), \normp[k]{W(x)(\omega)}, h(u),\right.\\
\specialpos{\hfill\left.\normp[j, j^{\prime}]{\nabla F(F_{z}^{-1}(u), z)},\normp[j^{\prime}, d]{\nabla\!Z(x)(\omega)},  \normp[j^{\prime}]{Z(x)(\omega)}\right)}\\
\specialpos{\rule{73pt}{0pt}:= g(u, \omega, z, x),\hfill}
\end{multline*}
where $S$ is again a polynomial with positive coefficients and $\mu = {\lambda}/{2}>0$.
 It is clear that  $g$ is a continuous function of the variable $u$ for almost all  $(\omega, z, x) \in \Omega \times \reels^{j^{\prime}-j} \times D$.\\
 On the one hand, using the hypotheses on  $Z$ and the hypothesis (\ref{majoration moments}) we have that $\forall p \in \naturels$, $\forall n \in \naturels$, $\forall \ell \in \naturels$ and $\forall m \in \naturels$,
\begin{multline*}
\int_{D} f^{p}(x) \bbE\!\left[\normp[k]{W(x)}^{n}\right]\bbE\!\left[\normp[j^{\prime}, d]{\nabla\!Z(x)}^{\ell}\right]\bbE\!\left[\normp[j^{\prime}]{Z(x)}^{m}\right] \ud x \\
\le \bC  \int_{D}f^{p}(x)\bbE\!\left[\normp[k]{W(x)}^{n}\right]\bbE\!\left[\normp[j^{\prime},d]{\nabla\!Z(x)}^{\ell}\right] \ud x< \infty.
\end{multline*}
In addition, it should be noted that
\begin{multline*}
\int_{D} \int_{\reels^{j^{\prime}-j}} \int_{\Omega} g(u, \omega, z, x) \ud \bbP(\omega)\ud z  \ud x\\
\specialpos{\rule{24pt}{0pt}=\int_{\reels^{j^{\prime}-j}} \frac{\bC }{\abs{J_{F}(F_{z}^{-1}(u), z)}} e^{- \mu \normp[j^{\prime}-j]{z}^{2}}\hfill}\\
\times\int_{D}
	\bbE\!\left[ 
		S\!\left(f(x), \normp[k]{W(x)}, h(u),  \normp[j, j^{\prime}]{\nabla F(F_{z}^{-1}(u), z)},
		 \right.\right.\\
		\left.\normp[j^{\prime},d]{\nabla\!Z(x)},\left. \normp[j^{\prime}]{Z(x)}\right)\right] \ud x \ud z.
\end{multline*}
On the other hand, since for all $\ell \in \naturels$ the functions $h$ and $T_{\ell}$ are continuous at $u$, the same is true for
$$
	u \mapsto \int_{D} \int_{\reels^{j^{\prime}-j}} \int_{\Omega} g(u, \omega, z, x) \ud \bbP(\omega)\ud z  \ud x.
$$%
By applying a weak version of Lebesgue's dominated convergence theorem, we deduce that the hypothesis ${\bH _{5}}$ is true.
 
 A similar proof can be done to verify the hypothesis $\bH _3$.
This completes the third part of the proof. \item\label{itm:proofProp3.12-4}  Now suppose that the processes $X$ and $Y$ satisfy the assumption  $\bA_4$.
Let us show that the hypothesis ${\bH _{5}}$ holds true, since the hypothesis ${\bH _1}$ is clearly satisfied.\\
 For $u \in \reels^{j}$, 
 \begin{multline*}
 \int_{D}\bbE[Y(x)H(\nabla
X(x))\given{X(x)=u}] {p}_{X(x)}(u) \ud x \\
= \int_{D} \int_{\reels \times \reels^{dj}} y H(\dot{x}) {p}_{Y(x), X(x), \nabla\!X(x)}(y, u, \dot{x})  \ud \dot{x} \ud y \ud x.
\end{multline*}
Using the hypotheses on the density ${p}_{Y(x), X(x), \nabla\!X(x)}(y, u, \dot{x})$, we obtain that the function appearing into the integral  is a continuous function of $u$ for almost all $(x, y, \dot{x}) \in D \times \reels \times \reels^{dj}$.\\
Moreover, since for any $A \in \frakL(\reels^{d}, \reels^{j})$, $H(A) \le \bC  \normp[j,d]{A}^{j}$, we easily obtain the following bound:
$\forall u \in \reels^{j}$ and for almost all $(x, y, \dot{x}) \in D \times \reels \times \reels^{dj}$,
\begin{align*}
\abs{y} H(\dot{x})&{p}_{Y(x), X(x), \nabla\!X(x)}(y, u, \dot{x})\\ 
& \le \bC  \abs{y} \normp[d,j]{\dot{x}}^{j} {p}_{Y(x), X(x), \nabla\!X(x)}(y, u, \dot{x})
:=g(x, y, u, \dot{x}).
\end{align*}
It is easy to see that for almost all $(x, y, \dot{x}) \in D \times \reels \times \reels^{dj}$, the function $u \longmapsto g(x, y, u, \dot{x})$ remains continuous.
Furthermore, we made the hypothesis that the function  
$$
	u \longmapsto \int_{D} \int_{\reels \times \reels^{dj}}  g(x, y, u, \dot{x}) \ud \dot{x} \ud y \ud x
$$%
is continuous.
These facts and a weak version of Lebesgue's dominated convergence theorem allow to verify hypothesis ${\bH _{5}}$.
The same is true for the hypothesis ${\bH _3}$.

This ends the proof of Proposition \ref{H2etH5}.
 \end{enumerate}
 \end{proofarg}
\section{Rice's formula for all level}
We have previously given conditions for certain classes of processes $X$ and $Y$ satisfying the hypotheses (${\bH _2}$, ${\bH _{4}}$), or (${\bH _1}$, ${\bH _3}$, ${\bH _{5}}$).\\
We will now provide a class of processes satisfying the hypotheses ${\bH _{i}}$, $i=1,\dots,5$ simultaneously and prove Proposition \ref{CabBis} and Theorem \ref{Hyp B} giving conditions on $X$ and $Y$ allowing the validity of the Rice formula for all level.
We should recall that Proposition \ref{CabBis} was proved in 1985 by Caba\~na \cite{AL121221}.
Our proof is deeply inspired by this work.
The hypothesis ${\bH _2}$ and ${\bH _{4}}$ are difficult to verify.
We will use the tools already developed to prove the required continuity.
For this purpose, \ie  in order to exhibit a class of processes verifying hypotheses ${\bH _2}$ and ${\bH _{4}}$, we resort to the only tool we have provided before, namely the use of  Proposition \ref{H1etH4}.
Our goal is then to construct a class of processes $Y$ satisfying the hypothesis $\bA_0$, \ie  such that $Y$ is a continuous process for which there exists $\lambda \in B_{j}$ such that $\supp(Y)  \subset \Gamma(\lambda)$,  $\normp[d,j]{(\nabla\!X(\bm{\cdot})|_{V_{\lambda}^{\perp}})^{-1}}$,  $Y(\bm{\cdot})$ and $\normp[j,d]{\nabla\!X(\bm{\cdot})}$ being uniformly  bounded on the support of  $Y$.

These assumptions being very demanding, the idea is, given a process $Y$, satisfying assumptions $\bA_{i}$, $i=1,\dots,4$, and thus by Proposition \ref{H2etH5} verifying the hypotheses ${\bH _{1}}$, ${\bH _3}$ and ${\bH _{5}}$, approximating the latter for fixed $n \in \naturels^{\star}$, by a process $Y^{(n)}$, defined as $Y^{(n)}:=  \sum\nolimits_{\lambda \in B_{j}} Y_{\lambda}^{(n)}$.
The process $Y_{\lambda}^{(n)}$ still verifies hypotheses ${\bH _1}$, ${\bH _3}$ and ${\bH _{5}}$ and above all assumption $\bA_0$.\\
In this form, by applying Theorem \ref{Cab1}, for $n \in \naturels^{\star}$ fixed we will propose a Rice formula for processes $X$ and $Y^{(n)}$ and for any level $y \in \reels^{j}$.
Then, we will make $n$ tend to infinity to obtain a Rice formula for $X$ and $Y$.\\

For this purpose, let $X: \Omega \times D \subset \Omega\times \reels^{d} \rightarrow \reels^{j}$ ($j \le d$) be a random field in $C^{1}(D,\reels^{j})$ where $D$ is an open set of $\reels^{d}$, and let $Y: \Omega \times D \subset \Omega \times \reels^{d} \rightarrow \reels$ a continuous processes.
In the same way as in Section \ref{validation}, we define for fixed $\lambda\in B_{j}$ and $x \in D_{X}^{r}$, $Y_{\lambda}(x):=  \eta_{\lambda}(x) Y(x)$, where $\eta_{\lambda}(t)$
has been defined in  (\ref{eta}).

For $n \in \naturels^*$, let $Y^{(n)}$ be
$$
	Y^{(n)}(x):= \sum\limits_{\lambda \in B_{j}} Y_{\lambda}^{(n)}(x),
$$%
for $x \in D$ , where $Y_{\lambda}^{(n)}$ is
\begin{multline*}
Y_{\lambda}^{(n)}(x):= Y_{\lambda}(x) f_n(x) \Psi(Y(x)/n) \Psi(\normp[j,d]{\nabla\!X(x)}/n)\\ 
\times\Psi\!\left(1 / (n \phi_{\lambda}(x)) \1_{\{\phi_{\lambda}(x)>0\}} + 2 \1_{\{\phi_{\lambda}(x)=0\}}\right) \1_{D_X^r}(x),
\end{multline*}
where the function $\phi_{\lambda}$ has been defined in (\ref{etalambda}), and $\Psi$ is a continuous even function on $\reels$, decreasing on $\reels^+$ such that
$$
\Psi(t):=
\begin{cases}
 1, & 0 \le t \le 1 \\
0,  & 2 \le t \\
\end{cases}
$$%
and $(f_n)_{n \in natural^{*}}$ is the sequence of functions defined from $\reels^d$ to $[0, 1]$ as follows
$$
f_n(x) := \frac{d(x, D^{2n})}{d(x, D^{2n})+ d(x, D^{(n)})},
$$%
where the closed sets $D^{2n}$ and $D^{(n)}$ are
$$
	D^{2n}:=\left\{x \in \reels^d, d(x, D^{c}) \le \frac{1}{2n}\right\}
	\text{ and }
	D^{(n)}:=\left\{x \in \reels^d, d(x, D^{c}) \ge \frac{1}{n}\right\}.
$$%
We will see later in the proof of  Lemma \ref{Yepsilon} that the functions  $(f_n)_{n \in \naturels^{*}}$ are well defined, continuous and such that  the support of ${f_{n}}_{|D}$ is contained in $D$ for each $n \in \naturels^*$.
In Lemma \ref{YepsilonBis} we will prove that $(f_n)_{n \in \naturels^{*}}$ is a sequence of nondecreasing functions tending to one when $n$ goes to infinity.

Let us explain a little more the choice of the terms composing the expression of $Y_{\lambda}^{(n)}$:
\begin{itemize}
\item $Y_{\lambda}(x) f_n(x)  \Psi (1 / (n \phi_{\lambda}(x)) \1_{\{\phi_{\lambda}(x)>0\}} + 2 \1_{\{\phi_{\lambda}(x)=0\}})$ ensures that $Y_{\lambda}^{(n)}(x)$ will  tend to $Y_{\lambda}(x)$ when $n$ tend to infinity for $x \in D_X^r$, as we will show by using the fact that $f_n(x)$ tends to one when $n$ tends to infinity.
Then for $x \in D_X^r$, $Y^{(n)}(x)$ will tend to $\sum_{\lambda \in B_j} Y_{\lambda}(x)=Y(x)$.
Furthermore, this implies since   $\forall n \in \naturels^{\star}$, $\supp({f_n}_{|D}) \subset D$, that $\supp(Y_{\lambda}^{(n)})  \subset \Gamma(\lambda).$
\item $\Psi\!\left((1 / n \phi_{\lambda}(x)) \1_{\{\phi_{\lambda}(x)>0\}} + 2 \1_{\{\phi_{\lambda}(x)=0\}}\right)$ ensures that  $\forall n \in \naturels^{\star}$,\\
 $\normp[d,j]{(\nabla\!X(\bm{\cdot})|_{V_{\lambda}^{\perp}})^{-1}}$ is uniformly bounded on the support of $Y_{\lambda}^{(n)}$.
\item $ \Psi(Y(x)/n) f_n(x)$ ensures that $Y_{\lambda}^{(n)}$ is uniformly bounded on $D$ since $f_n(x) \le 1$.
\item $ \Psi(\normp[j,d]{\nabla\!X(x)}/n)$ ensures that  $\forall n \in \naturels^{\star}$, $\normp[j,d]{\nabla\!X(\bm{\cdot})}$ is uniformly bounded on $\supp(Y_{\lambda}^{(n)})$.
\end{itemize}
We can now establish the following lemmas.
\spacebefore
\begin{lemm}
\label{Yepsilon}
Let  $X: \Omega \times D \subset \Omega\times \reels^{d} \rightarrow \reels^{j}$ ($j \le d$) be a random field belonging to
$\textbf{C}^{1}(D,\reels^{j})$, where $D$ is an open, convex and bounded set of  $\reels^{d}$, such that for almost all  $\omega \in \Omega$, $\nabla\!X(\omega)$ is Lipschitz with  Lipschitz constant $L_{X}(\boldsymbol{\omega})$ satisfying $\bbE\left[L_X^d(\bm{\cdot})\right] < \infty$.
 Let $Y: \Omega \times D \subset \Omega \times \reels^{d} \rightarrow \reels$ be a continuous process.
Then, on the one hand, for $n \in \naturels^{\star}$, $X$ and $Y^{(n)}$ satisfy the hypotheses ${\bH _2}$ and ${\bH _{4}}$.
On the other hand, if $Y$ satisfies (\ref{egalite Y}) and if $X$ and $Y$ satisfy one of the three assumptions $\bA_{i}$, $i=1, 2, 3$ or if $X$ and $Y$ satisfy the assumption $ \bA_4$, then  $\forall n \in \naturels^{\star}$, $X$ and $Y^{(n)}$ satisfy the hypotheses ${\bH _1}$, ${\bH _3}$ and ${\bH _5}$ and a fortiori the hypotheses ${\bH _i}$ , $i=1,\dots,5$.
\end{lemm}
In this form we have provided a class of processes  $X$ and $Y^{(n)}$ satisfying simultaneously the hypotheses ${\bH _i}$, $i=1,\dots,5$.
Then, by Theorem \ref{Cab1}, we  get that  $\forall n \in \naturels^{\star}$ and $\forall y \in \reels^{j}$ 
\begin{multline*}
\bbE\left[\int_{\calC_{X}^{D^{r}}(y)}Y^{(n)}(x) \ud\sigma_{d-j}(x)\right]\\
= \int_D{p}_{X(x)}(y) 
\bbE\left[Y^{(n)}(x)H(\nabla
X(x))\given{X(x)=y}\right] \ud x. 
\end{multline*}
The idea consists  to make $n$ tends to infinity.
More precisely we can show the following lemma.
\spacebefore
\begin{lemm}
\label{YepsilonBis}
Let $X: \Omega \times D \subset \Omega\times \reels^{d} \rightarrow \reels^{j}$ ($j \le d$) be a random field belonging to 
$C^{1}(D,\reels^{j})$, where $D$ is an open, convex and bounded set of  $\reels^{d}$, such that for almost all $\omega \in \Omega$, $\nabla\!X(\omega)$ is Lipschitz with Lipschitz constant $L_{X}(\omega)$ satisfying $\bbE\left[L_X^d(\bm{\cdot})\right] < \infty$.
Let  $Y: \Omega \times D \subset \Omega \times \reels^{d} \rightarrow \reels$ be a continuous process.
If $Y$ satisfies (\ref{egalite Y}) and if  $X$ and $Y$ satisfy one of the three assumptions $\bA_{i}$, $i=1, 2, 3$ or if $X$ and $Y$ satisfy the assumption $ \bA_4$, $\forall y \in \reels^{j}$,
$$
\lim_{n \to +\infty}\bbE\left[\int_{\calC_{X}^{D^{r}}(y)}Y^{(n)}(x) \ud\sigma_{d-j}(x)\right]=\bbE\left[\int_{\calC_{X}^{D^{r}}(y)}Y(x)\ud \sigma_{d-j}(x)\right]
$$%
and
\begin{multline*}
\lim_{n \to +\infty}\int_D{p}_{X(x)}(y) 
\bbE\left[Y^{(n)}(x)H(\nabla
X(x))\given{X(x)=y}\right]\ud x\\
= \int_D{p}_{X(x)}(y) 
\bbE\left[Y(x)H(\nabla
X(x))\given{X(x)=y}\right]\ud x.
\end{multline*}
\end{lemm}
Finally, we can establish the following proposition.
\spacebefore
\begin{prop}
\label{CabBis}
Let $X: \Omega \times D \subset \Omega\times \reels^{d} \rightarrow \reels^{j}$ ($j \le d$)  be a random field belonging to 
$C^{1}(D,\reels^{j})$, where $D$ is an open, convex and bounded set of  $\reels^{d}$, such that for almost all  $\omega \in \Omega$, $\nabla\!X(\omega)$ is Lipschitz with Lipschitz constant $L_{X}(\omega)$ such that  $\bbE\left[L_X^d(\bm{\cdot})\right] < \infty$.
Let $Y: \Omega \times D \subset \Omega \times \reels^{d} \rightarrow \reels$ be a continuous process.
If $Y$ satisfies (\ref{egalite Y}) and if  $X$ and $Y$ satisfy one of the three assumptions $\bA_{i}$, $i=1, 2, 3$ or if  $X$ and $Y$ satisfy the assumption $ \bA_4$, then for all  $y\in \reels^{j}$ we have
$$
\bbE\left[\int_{\calC_{X}^{D^{r}}(y)}Y(x)\ud \sigma_{d-j}(x)\right]=\int_D{p}_{X(x)}(y) 
\bbE\left[Y(x)H(\nabla
X(x))\given{X(x)=y}\right] \ud x.
$$%
\end{prop}
\spacebefore
\begin{rema}
\label{generalisation}
In the same way as in  Remark \ref{moduleconvexe} we can generalize this proposition considering $D$ an open and convex set not necessarily bounded.
\hglue16pt\hfill$\bullet$
\end{rema}
\spacebefore
\begin{proofarg}{Proof of the  Lemma \ref{Yepsilon}} 
Let $X: \Omega \times D \subset \Omega\times \reels^{d} \rightarrow \reels^{j}$ be a random field in
$C^{1}(D,\reels^{j})$ such that for almost all $\omega \in \Omega$, $\nabla\!X(\omega)$ is Lipschitz with Lipschitz constant $L_{X}(\omega)$ such that  $\bbE\left[L_X^d(\bm{\cdot})\right] < \infty$ and let $Y: \Omega \times D \subset \Omega \times \reels^{d} \rightarrow \reels$ be a continuous process.\\
Let us show that  $\forall n \in \naturels^{\star}$, the processes $X$ and $Y^{(n)}$ satisfy the hypotheses ${\bH _2}$ and ${\bH _{4}}$.\\
Consider for $n$ fixed in $\naturels^{\star}$ and for $\lambda$ fixed in $B_{j}$ the process  $Y_{\lambda}^{(n)}$.
We will now prove that the processes $X$ and $Y_{\lambda}^{(n)}$ satisfy the assumption $\bA_0$ (see page \pageref{H1etH4}).
Thus, we will deduce by using Proposition~\ref{H1etH4} that these processes satisfy the hypotheses ${\bH _2}$ and ${\bH _{4}}$.\\
Let us verify that  $Y_{\lambda}^{(n)}$ is continuous on $D$, which is a non-trivial fact due to the presence of $\1_{D_X^r}(x)$ in the definition of this function.
Let us first notice that since the sets  $(D^{2n})_{n \in \naturels^{*}}$ and  $(D^{(n)})_{n \in \naturels^{*}}$ are closed, the functions  $(f_n)_{n \in \naturels^{*}}$ are well defined and continuous on $\reels^d$ and then on $D$.
Consider now  $x \in (D_X^r)^{c_1}$ and a sequence $(x_{p})_{p \in \naturels^{\star}}$ of points in $D$ which converges to $x$ as $p$ tends to infinity.
We have $Y_{\lambda}^{(n)}(x)=0$.
Suppose that there exists a subsequence $(x_{p_k})_{k \in \naturels^{\star}}$ of $(x_{p})_{p \in \naturels^{\star}}$, such that $Y_{\lambda}^{(n)}(x_{p_{k}}) \neq 0$, for all $k \in \naturels^{\star}$.
In this case, necessarily  $\phi_{\lambda}(x_{p_k}) \ge \frac{1}{2n}$ for all $k \in \naturels^{\star}$, and since the function $\phi_{\lambda}$ is continuous on $D$, it turns out that $\phi_{\lambda}(x) \ge \frac{1}{2n}$.
Property (\ref{propriété2}) implies that $x \in \Gamma(\lambda) \subset D_X^r$, which leads to a contradiction.
All the points, except perharps a finite number, of the sequence $(x_{p})_{p \in \naturels^{\star}}$, are such that  $Y_{\lambda}^{(n)}(x_{p}) = 0$.
The sequence $(Y_{\lambda}^{(n)}(x_{p}))_{p \in \naturels^{\star}}$ then converges to zero.
Moreover, using a reasoning similar to the previous one, we can prove that the function \\
\centerline{$x \mapsto\Psi (1 / (n \phi_{\lambda}(x)) \1_{\{\phi_{\lambda}(x)>0\}} + 2 \1_{\{\phi_{\lambda}(x)=0\}})$}
is continuous on $D$ and then on $D_X^r$.
The function $Y_{\lambda}^{(n)}$ is continuous on $D_X^r$, which leads to the continuity of $Y_{\lambda}^{(n)}$ on $D$.\\
Let us now prove that  $\forall n \in \naturels^{\star}$, the support of this function is contained in  $\Gamma(\lambda)$, \ie 
$
\supp(Y_{\lambda}^{(n)}) \subset \Gamma(\lambda)
$.\\
To prove this inclusion, we will first prove that we have $\forall n \in \naturels^*,~\supp({f_n}_{|D}) \subset D$.\\
Indeed, since  $\forall n \in \naturels^*$ the set $D^{2n}$ is closed, we have
$$
\supp({f_n}_{|D}) = \overline{\left\{x \in D, d(x, D^{c}) >{1}/{(2n)}\right\}} \subset D.
$$%
The last inclusion is a consequence of the continuity of the distance function and the fact that $D$ is an open set.
\\
Therefore
$$
\supp(Y_{\lambda}^{(n)}) \subset \overline{\left\{x \in D, \phi_{\lambda}(x) \ge {1}/{(2n)} \right\}} \cap D \subset \left\{x \in D, \phi_{\lambda}(x) \ge {1}/{(2n)} \right\}.
$$%
The last inclusion comes from the fact that function $\phi_{\lambda}$ is continuous on $D$.
Finally, from (\ref{propriété2}) we have
$
\supp(Y_{\lambda}^{(n)}) \subset \Gamma(\lambda)
$.\\
Let us see that $Y_{\lambda}^{(n)}$ is uniformly bounded on its support.\\
We only need to prove that $Y_{\lambda}^{(n)}$ is uniformly bounded on $D$.
Consider $x \in D$ such that  $Y_{\lambda}^{(n)}(x) \neq 0$.
Then necessarily we have $
\abs{Y(x)} \le 2n$.
Since we have $f_n(x) \le 1$ and $\eta_{\lambda}x) \1_{\{x \in D_X^r\}} \le 1$, $\abs{Y_{\lambda}^{(n)}(x)} \le \abs{Y(x)}$.
It is therefore true that $\abs{Y_{\lambda}^{(n)}(x)} \le 2n$ and we get the result.\\
Let us show that  $\normp[d,j]{(\nabla\!X(\bm{\cdot})|_{V_{\lambda}^{\perp}})^{-1}}$ is uniformly bounded on the support of $Y_{\lambda}^{(n)}$.
We have seen that
$
\supp(Y_{\lambda}^{(n)}) \subset \left\{x \in D, \phi_{\lambda}(x) \ge \frac{1}{2n} \right\}$.
Then for $x \in \supp(Y_{\lambda}^{(n)})$, we have 
$\normp[d,j]{(\nabla\!X(x)|_{V_{\lambda}^{\perp}})^{-1}} \le 2n$.
Therefore, the result is true.
\\
Let us finally show that $\normp[j,d]{\nabla\!X(\bm{\cdot})}$ is uniformly bounded on the support of  $Y_{\lambda}^{(n)}$.
This follows from the following inclusion:
$$
\supp(Y_{\lambda}^{(n)}) \subset \overline{\left\{x \in D, \normp[j,d]{\nabla\!X(x)} \le 2n\right\}} \cap D
\subset \left\{x \in D, \normp[j,d]{\nabla\!X(x)} \le 2n\right\}.
$$%
The last inclusion comes from the fact that  $X$ is $C^1(D, \reels^j)$.\\
Finally the processes $X$ and $Y_{\lambda}^{(n)}$ satisfy the assumption $\bA_0$ and thus the hypotheses ${\bH _2}$ and ${\bH _{4}}$.\\
Using that $Y^{(n)}=  \sum_{\lambda \in B_{j}}Y_{\lambda}^{(n)}$ and that $\abs{Y^{(n)}}=  \sum_{\lambda \in B_{j}} \abs{\smash{Y_{\lambda}^{(n)}}}$, it is clear that $X$ and $Y^{(n)}$  satisfy also the hypotheses ${\bH _2}$ and ${\bH _{4}}$.\\
Now suppose that $Y$ satisfies (\ref{egalite Y}) and that $X$ and $Y$ satisfy one of the assumptions $\bA_{i}$, $i=1, 2, 3$.
Let us then prove  that  $X$ and $Y^{(n)}$ satisfy the hypotheses  ${\bH _1}$, ${\bH _3}$ and ${\bH _{5}}$.\\
For almost all $x \in D$, we have
\begin{equation}
	\label{decompositionY}
	Y(x)= G(x, W(x), X(x), \nabla\!X(x)),
\end{equation}
where
$G$
is a continuous function on 
$D \times \reels^k \times \reels^j \times \frakL(\reels^{d}, \reels^{j})$ and such that $\forall (x, z, u, A) \in D \times \reels^k \times \reels^j \times \frakL(\reels^{d}, \reels^{j})$,
$$
 \abs{G(x, z, u, A)} \le P(f(x), \normp[k]{z}, h(u), \normp[j,d]{A}).
 $$%
$\forall n \in \naturels^{\star}$ and $x \in D$ 
\begin{multline*}
Y^{(n)}(x)= \sum_{\lambda \in B_{j}} \eta_{\lambda}(x) Y(x) f_n(x) \Psi(Y(x)/n) \Psi(\normp[j,d]{\nabla\!X(x)}/n) \\ 
\times\Psi (1 / (n \phi_{\lambda}(x)) \1_{\{\phi_{\lambda}(x)>0\}} + 2 \1_{\{\phi_{\lambda}(x)=0\}}) \1_{D_X^r}(x).
\end{multline*}
We deduce that  $\forall n \in \naturels^{\star}$ and almost surely $\forall x \in D$,
\begin{equation}
	Y^{(n)}(x)= M_{n}(x, Y(x), \nabla\!X(x)),
	\label{decompositionYe}
\end{equation}
where  $\forall n \in \naturels^{\star}$, $M_{n}$ is a continuous function defined on
$D \times \reels \times \frakL(\reels^{d}, \reels^{j})$.
The proof of this last statement can be done in a similar way to that used to prove the continuity of  $Y_{\lambda}^{(n)}$ on $D$.
Moreover, 
$\forall (x, y, A) \in D \times \reels \times \frakL(\reels^{d}, \reels^{j})$,
$$
	\abs{M_{n}(x, y, A)} \le \bC  \abs{y}.
 $$%
 By  (\ref{decompositionY}), we have  $\forall n \in \naturels^{\star}$ and for almost all $x \in D$,
 \begin{align*}
	Y^{(n)}(x)&= M_{n}(x, G(x, W(x), X(x), \nabla\!X(x)), \nabla\!X(x))\\
	&=G_{n}(x, W(x), X(x), \nabla\!X(x)),
\end{align*}
where $\forall n \in \naturels^{\star}$ and $\forall (x, z, u, A) \in D \times \reels^k \times \reels^j \times \frakL(\reels^{d}, \reels^{j})$,
$$
	G_{n}(x, z, u, A)= M_{n}(x, G(x, z, u, A), A).
$$%
 It is clear that $G_{n}$ inherits the properties of $G$ and $M_{n}$; that is $G_{n}$ is a continuous function on 
$D \times \reels^k \times \reels^j \times \frakL(\reels^{d}, \reels^{j})$ and is such that
$$
	\forall (x, z, u, A) \in D \times \reels^k \times \reels^j \times \frakL(\reels^{d}, \reels^{j}),
$$%
 \begin{align*}
 \abs{G_{n}(x, z, u, A)} \le \bC  \abs{G(x, z, u, A)} &\le \bC  P(f(x), \normp[k]{z}, h(u), \normp[j,d]{A})\\
 &:= Q(f(x), \normp[k]{z}, h(u), \normp[j,d]{A}),
 \end{align*}
where $Q$ and $P$ are polynomials with positive coefficients and the functions
$$
f: D \longrightarrow \reels^{+}\mbox{  and } h: \reels^j \longrightarrow \reels^{+},
$$ are continuous.

Finally $Y^{(n)}$ satisfies (\ref{egalite Y}) and $X$ and $Y^{(n)}$ satisfy one of the three assumptions $\bA_1$, $\bA_2$ or $\bA_3$.
Using Proposition \ref{H2etH5}, we proved that the hypotheses  ${\bH _1}$, ${\bH _3}$ and ${\bH _{5}}$ hold for $X$ and $Y^{(n)}$.
Using the first part of this lemma we can conclude that the hypotheses ${\bH _{i}}$, $i=1,\dots,5$, are satisfied by $X$ and $Y^{(n)}$.\\
Now suppose that $X$ and $Y$ satisfy the assumption $ \bA_4$.
Let us prove that  $\forall n \in \naturels^{\star}$, $X$ and $Y^{(n)}$ satisfy ${\bH _1}$, ${\bH _3}$ and ${\bH _{5}}$.\\
Since for almost all $(x, y, \dot{x}) \in D \times \reels \times \reels^{dj}$ and $\forall u \in \reels^{j}$, the density ${p}_{Y(x), X(x), \nabla\!X(x)}(y, u, \dot{x})$ of the joint distribution
$$
	(Y(x), X(x), \nabla\!X(x))
$$%
exists (and is a continuous function of $u$) then for almost all $x \in D$ and $\forall u \in \reels^{j}$, the density ${p}_{X(x)}(u)$ of $X(x)$ exists and  ${\bH _1}$ is true.\\
Now using (\ref{decompositionYe}), we have   $\forall n \in \naturels^{\star}$ and $\forall u \in \reels^{j}$,
\begin{align*}
 L_{n}(u)&:= \int_{D} \bbE[Y^{(n)}(x)H(\nabla\!X(x))\given{X(x)=u}] {p}_{X(x)}(u) \ud x\\
	&= \int_{D}\bbE[M_{n}(x, Y(x), \nabla\!X(x))H(\nabla\!X(x))\given{X(x)=u}] {p}_{X(x)}(u) \ud x\\
	&= \int_{D} \bbE[L_{n}(x, Y(x), \nabla\!X(x))\given{X(x)=u}] {p}_{X(x)}(u) \ud x,
\end{align*}
where $L_{n}$ is a continuous function on $D \times \reels \times \reels^{dj}$ and since for all $A \in \frakL(\reels^{d}, \reels^{j})$, $H(A) \le \bC  \normp[j,d]{A}^{j}$, we have $\forall n \in \naturels^{\star}$ et $\forall (x, y, A) \in D \times \reels \times \frakL(\reels^{d}, \reels^{j})$,
$$
 \abs{L_{n}(x, y, A)} \le \bC  \abs{y} \normp[j,d]{A}^{j}.
 $$%
Finally, $\forall n \in \naturels^{\star}$ and $\forall u \in \reels^{j}$,
$$
 L_{n}(u)= \int_{D} \int_{\reels \times \reels^{dj}} L_{n}(x, y , \dot{x}) {p}_{Y(x), X(x), \nabla\!X(x)}(y, u, \dot{x})  \ud \dot{x} \ud y \ud x.
$$%
In the same way as in the proof of  Proposition \ref{H2etH5}, (item \ref{itm:proofProp3.12-4}, page \pageref{itm:proofProp3.12-4}), a weak version of Lebesgue's dominated convergence theorem implies that $\forall n \in \naturels^{\star}$, the function $u \longmapsto L_{n}(u)$ is continuous.
Then the hypothesis ${\bH _{5}}$ is true.
We can also obtain in the same way the hypothesis ${\bH _3}$.
This completes the proof of this lemma.
 \end{proofarg}
\spacebefore
\begin{proofarg}{Proof Lemma \ref{YepsilonBis}}
Let us first prove that $\forall y\in \reels^{j}$,
\begin{equation}
\label{gauche}
\lim_{n \to +\infty}\bbE\left[\int_{\calC_{X}^{D^{r}}(y)}Y^{(n)}(x) \ud \sigma_{d-j}(x)\right]=\bbE\left[\int_{\calC_{X}^{D^{r}}(y)}Y(x)\ud \sigma_{d-j}(x)\right].
\end{equation}
Recall that $\forall n \in \naturels^{\star}$ and $\forall x \in D$ we have
\begin{multline*}
Y^{(n)}(x)= \sum_{\lambda \in B_{j}} \eta_{\lambda}(x) Y(x) f_n(x) \Psi(Y(x)/n) \Psi(\normp[j,d]{\nabla\!X(x)}/n) \\ 
	\times\Psi (1 / (n \phi_{\lambda}(x)) \1_{\{\phi_{\lambda}(x)>0\}} + 2 \1_{\{\phi_{\lambda}(x)=0\}}) \1_{D_X^r}(x).
\end{multline*}
Notice that $\forall y\in \reels^{j}$ and $\forall x \in D$,\\

$\bullet$ $\lim\limits_{n \to +\infty} Y^{(n)}(x)  \1_{\calC_{X}^{D^{r}}(y)}(x)= Y(x)  \1_{\calC_{X}^{D^{r}}(y)}(x)$ \\

$\bullet$ $\abs{Y^{(n)}(x)} \1_{\calC_{X}^{D^{r}}(y)}(x) \le \abs{Y(x)}  \1_{\calC_{X}^{D^{r}}(y)}(x)$ \\

$\bullet$ $Y \cdot  \1_{\calC_{X}^{D^{r}}(y)} \in L^{1}(\ud \sigma_{d-j} \otimes \ud P)$\\[4pt]
Let us establish the first statement, by proving first that for $x \in D$, $\lim\limits_{n\to\infty}f_n(x)=1$.

Consider $x \in D$.
Since $D^{c}$ is closed, we have that  $d(x, D^{c}) >0$, and $\exists n_{0} \in \naturels^{*} \such d(x, D^{c}) \ge \frac{1}{n_{0}}$.
Thus for $n \ge n_{0}$, we have $d(x, D^{c}) \ge \frac{1}{n}$, which implies for  $n \ge n_{0}$, $x \in D^{(n)}$.
Consequently  $\forall n \ge n_{0}$, $d(x, D^{(n)})=0$ and $f_{n}(x)=1$ $\forall n \ge n_{0}$.\\
Finally, the first statement is a consequence of the inclusion (\ref{support}), \ie  $\forall \lambda \in B_j$ we have $
\supp(\eta_{\lambda}) \cap D_X^r \subset \Gamma(\lambda)$.
Indeed, this last inclusion implies that $\forall \lambda \in B_j$, $\forall x \in D_X^r \cap \Gamma^{c_1}(\lambda)$, $\eta_{\lambda}(x) =0$ and then  $\forall \lambda \in B_j$, $\forall x \in D_X^r $, $\lim\limits_{n \to +\infty} Y_{\lambda}^{(n)}(x)= \eta_{\lambda}(x) Y(x)$ so that  $\forall x \in D_X^r $, $\lim\limits_{n \to +\infty} Y^{(n)}(x)= (\sum_{\lambda \in B_j} \eta_{\lambda}(x)) Y(x)= Y(x)$.\\[4pt]
The last statement can be proved in the following way.
Using Lemma \ref{Yepsilon}, $X$ and $Y^{(n)}$ satisfy ${\bH _i}$ , $i=1,\dots,5$ and in particular ${\bH _1}$, ${\bH _2}$ and ${\bH _3}$.
Remark \ref{chichi6} allows us to obtain that $\forall n \in \naturels^{\star}$ and $\forall y \in  \reels^{j}$,%
\begin{multline*}
\bbE\left[\int_{\calC_{X}^{D^{r}}(y)}\abs{Y^{(n)}(x)}\ \ud \sigma_{d-j}(x)\right]\\
	 = \int_D{p}_{X(x)}(y)\bbE\left[\abs{Y^{(n)}(x)} H(\nabla
X(x))\given{X(x)=y}\right] \ud x.
\end{multline*}
By noticing that the sets  $(D^{2n})_{n \in \naturels^{*}}$ and  $(D^{(n)})_{n \in \naturels^{*}}$, which define the sequence $(f_n)_{n \in \naturels^{*}}$, are respectively decreasing and nondecreasing sequences, we obtain that the sequence $(f_n)_{n \in \naturels^{*}}$ is nondecreasing.
Since the function $\Psi$ is an even function on $\reels$ and decreasing on $\reels^+$, the sequence
$(\abs{Y^{(n)}})_{n \in \naturels^{*}}$ is non-decreasing.
\\
We can apply Beppo Levi's theorem and we have $\forall y \in \reels^{j}$
$$
\lim_{n \to +\infty} \uparrow\bbE\left[\int_{\calC_{X}^{D^{r}}(y)}\abs{Y^{(n)}(x)} \ud\sigma_{d-j}(x)\right]
=\bbE\left[\int_{\calC_{X}^{D^{r}}(y)}\abs{Y(x)} \ud\sigma_{d-j}(x)\right].
$$%
Similarly, $\forall y \in \reels^{j}$,
\begin{align*}
\lim_{n \to +\infty} &\uparrow  \int_D{p}_{X(x)}(y) 
\bbE\left[\abs{Y^{(n)}(x)} \1_{D_{X}^{r}}(x) H(\nabla
X(x))\given{X(x)=y}\right] \ud x \\
&=  \int_D{p}_{X(x)}(y) 
\bbE\left[\abs{Y(x)} \1_{D_{X}^{r}}(x) H(\nabla
X(x))\given{X(x)=y}\right]\ud x\\
&=\int_D{p}_{X(x)}(y) 
\bbE\left[\abs{Y(x)} H(\nabla
X(x))\given{X(x)=y}\right]\ud x.
\end{align*}
The last equality comes from the fact that $\forall x \in D$ we have $$\1_{D_{X}^{r}}(x) H(\nabla
X(x))= H(\nabla
X(x)).
$$
We then obtain $\forall y \in \reels^{j}$
\begin{multline}
\bbE\left[\int_{\calC_{X}^{D^{r}}(y)}\abs{Y(x)} \ud\sigma_{d-j}(x)\right]\\
= \int_D{p}_{X(x)}(y) 
\bbE\left[\abs{Y(x)} H(\nabla
X(x))\given{X(x)=y}\right]\ud x < \infty,\label{deuxmembres}
\end{multline}
since $X$ and $Y$ satisfy one of the four assumptions  $\bA_1$, $\bA_2$, $\bA_3$ or $ \bA_4$ and by Proposition \ref{H2etH5}, satisfy the hypothesis ${\bH _3}$.\\
We have shown that $Y \cdot  \1_{\calC_{X}^{D^{r}}(y)} \in L^{1}(\ud\sigma_{d-j} \otimes \ud P)$.
Then using Lebesgue's dominated convergence theorem, we can deduce (\ref{gauche}).\\
Let us show that $\forall y \in  \reels^{j}$,
\begin{multline}
	\lim_{n \to +\infty}\int_D{p}_{X(x)}(y) 
		\bbE\left[Y^{(n)}(x)H(\nabla
		X(x))\given{X(x)=y}\right]\ud x\\ \label{droite}
	= \int_D{p}_{X(x)}(y) 
		\bbE\left[Y(x)H(\nabla
		X(x))\given{X(x)=y}\right]\ud x.
\end{multline}
In the same way as before, we notice that $\forall y \in  \reels^{j}$ and for almost all $x \in D$,
\begin{itemize}
	\item[$\bullet$] $\lim\limits_{n \to +\infty} Y^{(n)}(x)  H(\nabla\!X(x)) {p}_{X(x)}(y)= Y(x) H(\nabla\!X(x))  {p}_{X(x)}(y)$
	\item[$\bullet$] $\abs{Y^{(n)}(x)}  H(\nabla\!Xx)) {p}_{X(x)}(y) \le \abs{Y(x)}   H(\nabla\!X(x)) {p}_{X(x)}(y)$
	\item[$\bullet$] $\bbE\left[\abs{Y(x)} H(\nabla\!X(x)) {p}_{X(x)}(y)\given{X(x)=y}\right] < \infty$
\end{itemize}
and the finiteness of the last expression results from that of the second integral in (\ref{deuxmembres}).\\
Lebesgue dominated convergence theorem allows to write
$\forall y \in  \reels^{j}$ and for almost all  $x \in D$,
\begin{multline*}
	\bullet \lim\limits_{n \to +\infty}\bbE\left[Y^{(n)}(x)H(\nabla
	X(x))\given{X(x)=y}\right] {p}_{X(x)}(y)\\
	=\bbE\left[Y(x)  H(\nabla\!X(x))\given{X(x)=y}\right] {p}_{X(x)}(y).
\end{multline*}
Moreover, $\forall y \in \reels^{j}$ and almost surely $\forall x \in D$,
\begin{multline*}
\bullet \abs{\bbE\left[Y^{(n)}(x)H(\nabla
X(x))\given{X(x)=y}\right] {p}_{X(x)}(y)}\\
 \le\bbE\left[\abs{Y(x)} H(\nabla\!X(x))\given{X(x)=y}\right] {p}_{X(x)}(y)  \in L^{1}(D, \ud x).
\end{multline*}
The last statement comes from the fact that the second integral in (\ref{deuxmembres}) is finite.\\
Lebesgue's dominated convergence theorem allows to obtain (\ref{droite}).
This completes the proof of Lemma \ref{YepsilonBis} and consequently the proof of Proposition \ref{CabBis}. 
\end{proofarg}%
We can now state Theorem \ref{Hyp B} which follows.
It allows us to weaken assumptions $\bA_{i}$, $i=1, 2, 3$.
More precisely, we want to avoid assuming the existence of uniform lower (or upper) bounds for the variance of the $Z$ process appearing under these assumptions.\\
In the first three assumptions ${\bB_1}$, ${\bB_2}$ and ${\bB_3}$, we will assume  that $Y$ can be written as in the (\ref{egalite Y}).
Let $D$ be an open set of $\reels^d$.
\begin{itemize}
\item ${\bB_1}$: Let $X: \Omega \times D \subset \Omega\times \reels^{d} \rightarrow \reels^{j}$ ($j \le d$) be a Gaussian random field belonging to $C^1(D, \reels^j)$, such that $\forall x \in D$, the vector $X(x)$ has a density.
Moreover, for almost all $x \in D$, the field $W(x)$ is independent of the vector $(X(x), \nabla\!X(x))$, and $\forall n \in \naturels$, 
$$
	\int_{D}\bbE\!\left[\normp[k]{W(x)}^{n}\right] \ud x< \infty.
$$%
\item ${\bB_2}$: $\forall x \in D$, $X(x)=F(Z(x))$, where $F: \reels^j \longrightarrow \reels^j$ is a bijection belonging to $C^1(\reels^j, \reels^j)$, such that $\forall z \in \reels^j$, the Jacobian of $F$ in $z$, that is $J_{F}(z)$ satisfies $J_{F}(z)\neq 0$ and the function $F^{-1}$ is continuous.\\
Let $Z: \Omega \times D \subset \Omega\times \reels^{d} \rightarrow \reels^{j}$ ($j \le d$)  be a Gaussian process belonging to $C^1(D, \reels^j) \such \forall x \in D$, the vector $Z(x)$ has a density.
Moreover, for almost all $x \in D$, $W(x)$ is independent of the vector $(Z(x), \nabla\!Z(x))$, and $\forall n \in \naturels$, %
$$
	\int_{D}\bbE\!\left[\normp[k]{W(x)}^{n}\right] \ud x< \infty.
$$%
\item ${\bB_3}$: For all $ x \in D$, $X(x)=F(Z(x))$, where $Z: \Omega \times D \subset \Omega\times \reels^{d} \rightarrow \reels^{j^{\prime}}$ ($j<j^\prime$) is a Gaussian random field belonging to $C^1(D, \reels^{j^{\prime}}) \such \forall x \in D$, the vector $(Z(x), \nabla\!Z(x))$ has a density.
Moreover, for almost all $x \in D$, $W(x)$ is independent of the vector $(Z(x), \nabla\!Z(x))$.
Finally, $\forall n \in \naturels$ 
$$
	\int_{D}\bbE\!\left[\normp[k]{W(x)}^{n}\right]  \ud x< \infty.
$$%
The function $F$ verifies assumption ${(\bF)}$ given in assumption $\bA_3$.
\item ${\bB_4}$: Is the same assumption $ \bA_4$.
\end{itemize}
Let us now state  the hypothesis ${\bH _6}$.
\begin{itemize}
\item ${\bH _6}$: $\forall y \in \reels^{j}$,
$$
	 \int_{D}{p}_{X(x)}(y)\bbE\left[ \abs{Y(x)} H(\nabla\!X(x))\given{X(x)=y}\right] \ud x < \infty.
$$%
\end{itemize}
We are now ready to formulate the following theorem.
\spacebefore\begin{theo}
\label{Hyp B}
Let $X: \Omega \times D \subset \Omega\times \reels^{d} \rightarrow \reels^{j}$ ($j \le d$)  be a random field belonging to
$C^{1}(D,\reels^{j})$, where $D$ is an open and bounded convex set of $\reels^{d}$, such that for almost all $\omega \in \Omega$, $\nabla\!X(\omega)$ is Lipschitz with Lipschitz constant $L_{X}(\omega) \such \bbE\left[L_X^d(\bm{\cdot})\right] < \infty$.
 Let $Y: \Omega \times D \subset \Omega \times \reels^{d} \rightarrow \reels$ be a continuous process.\\
If  $Y$ satisfies (\ref{egalite Y}) and if  $X$ and $Y$ satisfy  one of the three assumptions $\bB_{i}$, $i=1, 2, 3$ and the hypothesis ${\bH _6}$ or if $X$ and $Y$ satisfy the assumption ${\bB_4}$, then $\forall y \in \reels^{j}$ we have%
$$
\bbE\left[\int_{\calC_{X}^{D^{r}}(y)}Y(x)\ud \sigma_{d-j}(x)\right]=\int_D{p}_{X(x)}(y) 
\bbE\left[Y(x)H(\nabla
X(x))\given{X(x)=y}\right]\ud x.
$$%
\end{theo}
\spacebefore
\begin{rema}
\label{moduledeY}
Under the same hypotheses as in Theorem \ref{Hyp B}, eliminating the condition $\bbE\left[L_X^d(\bm{\cdot})\right] < \infty$ and the hypothesis $\bH_6$, we obtain the following inequality
$\forall y \in  \reels^{j}$
$$
	\bbE\left[\int_{\calC_{X}^{D^{r}}(y)} \abs{Y(x)} \ud \sigma_{d-j}(x)\right]
		\le \int_D{p}_{X(x)}(y)\bbE\left[\abs{Y(x)} H(\nabla\!X(x))\given{X(x)=y}\right]\ud x.
$$%
 \hfill$\bullet$
\end{rema}
 \spacebefore
\begin{rema}
 \label{ouvertconvexe}
As in Remark \ref{moduleconvexe}, we can generalize the theorem and the Remark \ref{moduledeY} by considering that $D$ is an open and convex possibly unbounded set.
\hfill$\bullet$
\end{rema}

\spacebefore
\begin{proofarg}{Proof of Theorem \ref{Hyp B}} 
Let us consider  $\forall n \in \naturels^{\star}$, the sets\\
\centerline{$D_n:=\left\{x \in \reels^d, d(x, D^{c}) > \frac{1}{n} \right\}$.}
$\forall n \in \naturels^{\star}$, $D_n$ is an open set included in $D$.
Now consider the restrictions $X_{|D_n}$ and $Y_{|D_n}$.
It is clear that if $Y$ satisfies (\ref{egalite Y}) and if  $X$ and $Y$ satisfy one of the three assumptions  $\bB_{i}$, $i=1, 2, 3$ then  $\forall n \in \naturels^{\star}$, $Y_{|D_n}$ satisfies (\ref{egalite Y}) and $X_{|D_n}$ and $Y_{|D_n}$ satisfy one of the three assumptions  $\bA_{i}$, $i=1, 2, 3$, where we have replaced the open set  $D$ by the open set  $D_n$.
Indeed, it suffices for this to point out that $\forall n \in \naturels^*$, $D_n \subseteq \left\{x \in \reels^d, d(x, D^{c}) \ge \frac{1}{n} \right\}$ which is a compact set contained in $D$.\\
The set $D_n$ may not be convex but $D_n \subset D$, and the latter set is convex.\\
We apply Remark \ref{generalisation} which follows Proposition \ref{CabBis} to $X_{|D_n}$ and to $Y_{|D_n}$ (resp$.$
$\abs{Y}_{|D_n}$).
We get: $\forall y \in \reels^j$ and  $\forall n \in \naturels^*$,\
$$
	\bbE\left[\int_{\calC_{D_n, X}^{D^{r}}(y)}Y(x) \ud \sigma_{d-j}(x)\right]
		=\int_{D_{n}} {p}_{X(x)}(y)\bbE\left[Y(x)H(\nabla
X(x))\given{X(x)=y}\right]\ud x, 
$$%
and also $\forall y \in \reels^j$ and  $\forall n \in \naturels^*$,
\begin{multline*}
	\bbE\left[\int_{\calC_{D_n, X}^{D^{r}}(y)}\abs{Y(x)} \ud \sigma_{d-j}(x)\right]\\
		=\int_{D_{n}} {p}_{X(x)}(y)\bbE\left[\abs{Y(x)} H(\nabla\!X(x))\given{X(x)=y}\right]\ud x. 
\end{multline*}
Knowing that $D$ is an open set of  $\reels^d$, $\lim_{n \to +\infty} \uparrow D_n =D$, Beppo Levi's theorem leads to, $\forall y \in \reels^j$,
\begin{multline*}
	\bbE\left[\int_{\calC_{X}^{D^{r}}(y)} \abs{Y(x)} \ud \sigma_{d-j}(x)\right]\\
		=\int_{D} {p}_{X(x)}(y)\bbE\left[\abs{Y(x)} H(\nabla\!X(x))\given{X(x)=y}\right]\ud x < \infty,
\end{multline*}
by using the hypothesis ${\bH _6}$.\\
We can apply Lebesgue's dominated convergence theorem to get, $\forall y \in \reels^j$,%
\begin{multline*}
	\bbE\left[\int_{\calC_{X}^{D^{r}}(y)}Y(x) \ud\sigma_{d-j}(x)\right]
		=\int_{D} {p}_{X(x)}(y)\bbE\left[Y(x)H(\nabla\!X(x))\given{X(x)=y}\right]\ud x.
\end{multline*}
If $X$ and $Y$ satisfy the assumption $\bB_{4}$ which is nothing else than the assumption $\bA_4$, the above equality is trivial since already settled  in Proposition \ref{CabBis}.\\
The proof of Theorem \ref{Hyp B} is finished.
\end{proofarg}
\spacebefore
\begin{proofarg}{Proof of Remark \ref{moduledeY}} 
In the same way as in the proof of Lemma \ref{Yepsilon}, let us define  $\forall n \in \naturels^{\star}$, the \rv
$Z^{(n)}$ by $Z^{(n)}(x):= \sum\nolimits_{\lambda \in B_{j}} Z_{\lambda}^{(n)}(x)$ for $x \in D$, where we have defined for $\lambda$ fixed in $B_{j}$ the \rv
$Z_{\lambda}^{(n)}$ by
$$
	Z_{\lambda}^{(n)}(x):= Y_{\lambda}^{(n)}(x)  \Psi(L_{X}(\bm{\cdot})/n), 
$$%
where we recall that we have defined the \rv
$Y_{\lambda}^{(n)}$ by
\begin{multline*}
	Y_{\lambda}^{(n)}(x)= Y_{\lambda}(x) f_n(x) \Psi(Y(x)/n) \Psi(\normp[j,d]{\nabla\!X(x)}/n)\\
		\times \Psi (1 / (n \phi_{\lambda}(x)) \1_{\{\phi_{\lambda}(x)>0\}} + 2 \1_{\{\phi_{\lambda}(x)=0\}}) \1_{D_X^r}(x).
\end{multline*}

Note that we cannot work as in Lemma  \ref{Yepsilon}, since we are not able to show that $Z_{\lambda}^{(n)}$ verifies hypotheses ${\bH _3}$ and ${\bH _{5}}$.
In fact, we cannot apply the results of Proposition \ref{H2etH5}, since we cannot verify that $Z_{\lambda}^{(n)}$ verifies assumptions  $\bA_{i}$, $i=1,\dots,4$.
Indeed for $x \in D$, $Z_{\lambda}^{(n)}(x)$ depends on the whole trajectory of the $X$ process via the term $ \Psi(L_{X}(\bm{\cdot})/n)$.
\\
The processes $X$ and $Z_{\lambda}^{(n)}$ satisfy the hypotheses of Remark \ref{moduledeYBis} and therefore $X$ and $Z_{\lambda}^{(n)}$ satisfy ${\bH _2}$.
As in the proof of Lemma \ref{Yepsilon}, we deduce that  $X$ and $Z^{(n)}$ satisfy ${\bH _2}$.\\
Moreover, assuming for the moment that $Y$ satisfies (\ref{egalite Y}) and that $X$ and $Y$ satisfy one of the three assumptions $\bA_1$, $\bA_2$, $\bA_3$ or that  $X$ and $Y$ satisfy $ \bA_4$,  Proposition \ref{H2etH5} allows us to deduce that ${\bH _1}$ is satisfied.\\
By Proposition \ref{dualité} we obtain that for almost all $y \in \reels^{j}$%
\begin{multline*}
	\bbE\left[\int_{\calC_{X}^{D^{r}}(y)} \abs{Z^{(n)}(x)} \ud \sigma_{d-j}(x)\right]\\
		= \int_D{p}_{X(x)}(y)\bbE\left[\abs{Z^{(n)}(x)} H(\nabla\!X(x))\given{X(x)=y}\right]\ud x,
\end{multline*}
thus for almost all $y \in \reels^{j}$
\begin{multline*}
		\bbE\left[\int_{\calC_{X}^{D^{r}}(y)} \abs{Z^{(n)}(x)} \ud \sigma_{d-j}(x)\right]\\
		\le \int_D{p}_{X(x)}(y)\bbE\left[\abs{Y(x)} H(\nabla\!X(x))\given{X(x)=y}\right]\ud x.
\end{multline*}
By Proposition \ref{H2etH5}, the processes $X$ and $Y$ satisfy ${\bH _3}$.
Since $X$ and $Z^{(n)}$ satisfy ${\bH _2}$,  we deduce that the right and left terms of the last inequality are continuous in terms of the variable  $y$.
Then the inequality is true $\forall y \in \reels^{j}$.\\
In the same way as in the proof of Lemma \ref{YepsilonBis}, using Beppo Levi's theorem, we obtain $\forall y \in \reels^{j}$
$$
	\lim_{n \to +\infty} \uparrow\bbE\left[\int_{\calC_{X}^{D^{r}}(y)} \abs{Z^{(n)}(x)} \ud\sigma_{d-j}(x)\right]
		=\bbE\left[\int_{\calC_{X}^{D^{r}}(y)} \abs{Y(x)} \ud\sigma_{d-j}(x)\right].
$$%
We have shown that $\forall y \in \reels^{j}$
\begin{multline*}
	\bbE\left[\int_{\calC_{X}^{D^{r}}(y)} \abs{Y(x)} \ud\sigma_{d-j}(x)\right] \\
		\le \int_D{p}_{X(x)}(y)\bbE\left[ \abs{Y(x)} H(\nabla\!X(x))\given{X(x)=y}\right]\ud x,
\end{multline*}
this completes the proof of this remark whenever  $Y$ satisfies (\ref{egalite Y}) and that  $X$ and $Y$ satisfy one of the three assumptions $\bA_i $, for $i=1,\dots,3$ or that $X$ and $Y$ satisfy the assumption $ \bA_4$ and then the assumption  ${\bB_4}$.\\
In the case where $Y$ satisfies (\ref{egalite Y}) and $X$ and  $Y$ satisfy one of the three assumptions  ${\bB_i}$, for $i=1,\dots,3$, in the same way as in the proof of  Theorem \ref{Hyp B}, we apply  $\forall n \in \naturels^*$ the above inequality to $X_{|D_n}$ and $Y_{|D_n}$ which satisfy one of the three assumptions  $\bA_i $, for $i=1,\dots,3$.
Letting $n$ tends to infinity and applying Beppo Levi's theorem, we get the desired result.
Let us remark that the right-hand side term is not necessary finite because we do not  assume the hypothesis ${\bH _6}$.
\end{proofarg}

Our goal in this step of these notes is to propose a Rice's formula which is true for any level but without the assumption  $\bbE\left[L_X^d(\bm{\cdot})\right] < \infty$ which was given in Theorem \ref{Hyp B}.
We will propose in the following a little better than the inequality appearing in Remark \ref{moduledeY}.
To do this, we will replace in this theorem one of the assumptions ${\bB_i}$, $i=1,\dots,4$ by a slightly stronger assumption ${\bB_{i}^\star}$.
In the first three assumptions ${\bB_1^\star}$, ${\bB_2^\star}$ and ${\bB_3^\star}$, we will assume that $Y$ can be written as in (\ref{egalite Y}).\\ 
More precisely, let $D$ be an open set of $\reels^d$ and consider the following assumptions:
\begin{itemize}
\item ${\bB_1^\star}$: It is the assumption $\bB_1$, plus the following hypothesis:
for almost all $(x_1, x_2)  \in D \times D$, the density of the vector $(X(x_1), X(x_2))$ exists.
\item ${\bB_2^\star}$: It is the assumption $\bB_2$, plus  the following hypothesis:
for almost all $(x_1, x_2)  \in D \times D$,  the density of the vector $(Z(x_1), Z(x_2))$ exists.
\item ${\bB_{3}^\star}$: It is the assumption $\bB_{3}$, plus  the following hypothesis:
for almost all $(x_1, x_2)  \in D \times D$,  the density of the vector $(Z(x_1), Z(x_2))$ exists.
\item ${\bB_{4}^\star}$: It is the assumption $\bB_{4}$, plus  the following hypotheses:
\begin{enumerate}
\item The function %
$$
u \longmapsto \int_{D} \int_{\reels \times \reels^{dj}} y^{2} \normp[dj]{\dot{x}}^{j} {p}_{Y(x), X(x), \nabla\!X(x)}(y, u, \dot{x})  \ud \dot{x} \ud y \ud x
$$%
is continuous.
\item For almost all $(x_1, x_2, \dot{x}_1, \dot{x}_2) \in D \times D \times \reels^{dj} \times \reels^{dj}$ and for all $(u, v) \in \reels^{j} \times \reels^{j}$,  the density
$${p}_{X(x_1), X(x_2),\nabla\!X(x_1), \nabla\!X(x_2)}(u, v, \dot{x}_1, \dot{x}_2),$$ of the vector$(X(x_1), X(x_2), \nabla\!X(x_1), \nabla\!X(x_2))$ exists.
\item Moreover,  $\forall y \in \reels^j$,  the function
\begin{multline*}
	 (u, v) \longmapsto \int_{D \times D} \int_{\reels^{dj} \times \reels^{dj}} \normp[dj]{\dot{x}_1}^{j}  \normp[dj]{\dot{x}_2}^{j}{p}_{X(x_1), X(x_2), \nabla\!X(x_1), \nabla\!X(x_2)}(u, v, \dot{x}_1, \dot{x}_2)\\
 		 \ud \dot{x}_1\ud \dot{x}_2 \ud x_1 \ud x_2
\end{multline*}
is bounded in a neighborhood of $(y, y)$.
\end{enumerate}
\end{itemize}

Let us express now the hypothesis ${\bH _6^\star}$.

\begin{itemize}
\item ${\bH _{6}^\star}$: $\forall y \in \reels^j$, the function
\begin{multline*}
	(u, v) \longmapsto \int_{D \times D} {p}_{X(x_1), X(x_2)}(u, v)  \\
		 \times\bbE\left[H(\nabla\!X(x_1)) H(\nabla\!X(x_2))\given{X(x_1)=u, X(x_2)=v}\right] \ud x_1 \ud x_2
\end{multline*}
is bounded in a  neighborhood of $(y, y)$.
\end{itemize}
Finally we can state the following theorem.
\spacebefore\begin{theo}
\label{Cabter}
Let $X: \Omega \times D \subset \Omega\times \reels^{d} \rightarrow \reels^{j}$ ($j \le d$) be a random field in $C^{1}(D,\reels^{j})$, where $D$ is an open and convex bounded set of $\reels^{d}$, such that for almost all $\omega \in \Omega$, $\nabla\!X(\omega)$ is Lipschitz.
 Let $Y: \Omega \times D \subset \Omega \times \reels^{d} \rightarrow \reels$ be a continuous process.
If $Y$ satisfies (\ref{egalite Y}) and if  $X$ and $Y$ satisfy one of the three assumptions  ${\bB_{i}^\star}$, $i=1, 2, 3$ and the hypotheses ${\bH _6}$ and ${\bH _6^\star}$ or if $X$ and $Y$satisfy the assumption ${\bB_4^\star}$, then for all $\textbf{y}\in \reels^{j}$ we have
$$
\bbE\left[\int_{\calC_{X}^{D^{r}}(y)}Y(x) \ud\sigma_{d-j}(x)\right]=\int_D{p}_{X(x)}(y) 
\bbE\left[Y(x)H(\nabla
X(x))\given{X(x)=y}\right] \ud x.
$$%
\end{theo}
\spacebefore
\begin{rema}
In the same way as in Remark \ref{moduleconvexe}, this theorem can be generalized to the case where $D$ is an open and convex set not necessarily bounded.
\rule{24pt}{0pt}\hfill$\bullet$
\end{rema}
\spacebefore
\begin{proofarg}{Proof of Theorem \ref{Cabter}} 
For any $z \in \reels^{j}$, with the same notations as in the proof of Lemma \ref{Yepsilon} and  Remark \ref{moduledeY}, we have  $\forall n \in \naturels^*$,
\begin{multline*}
\abs{\bbE\left[\int_{\calC_{X}^{D^{r}}(z)}\abs{Z^{(n)}(x)} \ud\sigma_{d-j}(x)\right] -
	\bbE\left[\int_{\calC_{X}^{D^{r}}(z)}\abs{Y^{(n)}(x)} \ud\sigma_{d-j}(x)\right]
	}\\
=\bbE\left[\int_{\calC_{X}^{D^{r}}(z)}\abs{Y^{(n)}(x)} (1- \Psi(L_{X}(\bm{\cdot})/n))  \ud\sigma_{d-j}(x)\right].
\end{multline*}
Let us assume for the moment that if $Y$ satisfies (\ref{egalite Y}) then $X$ and $Y$ satisfy one of the three assumptions  $\bA_{i}$, $i=1, 2, 3$.

By Lemma \ref{Yepsilon} and Proposition \ref{H2etH5}, since $X$ and $Y$ satisfy one of the four hypotheses $\bA_{i}$, $i=1,\dots,4$, $\forall n \in \naturels^*$, $X$ and $Y^{(n)}$ satisfy the hypotheses  ${\bH _1}$ and ${\bH _3}$ (and also ${\bH _5}$).
By Proposition \ref{dualité}, we then have for almost all $z \in \reels^{j}$ and  $\forall n \in \naturels^{\star}$
\begin{multline*}
	\bbE\left[\int_{\calC_{X}^{D^{r}}(z)}\abs{Y^{(n)}(x)} \ud\sigma_{d-j}(x)\right] \\
	= \int_D{p}_{X(x)}(z)\bbE\left[\abs{Y^{(n)}(x)}H(\nabla\!X(x))\given{X(x)=z}\right] \ud x,
\end{multline*}
and a similar formula is true for the \rv $Y^{(n)}$.

For simplicity, let us denote  $\ds \repint_{y-\delta}^{y+\delta} f(x) \ud x$ the following multiple integral 
 $$
 	\int_{y_1-\delta}^{y_1+\delta} \dots \int_{y_{j}-\delta}^{y_{j}+\delta} f(x) \ud x,
$$%
where $\delta >0$, $y:=(y_1, y_2,\dots, y_{j})$ and   $f$ is a positive or integrable real valued function defined over $\prod_{i=1}^{j} [y_{i}- \delta, y_{i}+\delta]$.\\
With this notation and using the Schwarz inequality we obtain  $\forall n \in \naturels^*$, $\forall y \in \reels^{j}$ and $\forall \delta >0$:
\begin{multline}
\label{autre}
\!\left({\frac{1}{2\delta}}\right)^j \repint_{y-\delta}^{y+\delta}
\left\vert
	\bbE\left[\int_{\calC_{X}^{D^{r}}(z)}\abs{Z^{(n)}(x)} \ud\sigma_{d-j}(x)\right] - \right. \\
 	\left. \int_D{p}_{X(x)}(z)\bbE\left[\abs{Y^{(n)}(x)}H(\nabla\!X(x))\given{X(x)=z}\right] \ud x
	\rule{0pt}{20pt}\right\vert 
	\ud z  \\
\specialpos{\rule{30pt}{0pt}\le\!\left({\bbE\!\left[(1- \Psi(L_{X}(\bm{\cdot})/n))^{4}\right]}\right)^{\frac{1}{4}}\hfill}\\
\times \!\left({\!\left(\frac{1}{2\delta}\right)^j \repint_{y-\delta}^{y+\delta}
\bbE\!\left[\int_{\calC_{X}^{D^{r}}(z)} (Y^{(n)}(x))^{2} \ud\sigma_{d-j}(x)\right] \ud z
 }\right)^{\frac{1}{2}}  \\\times
 \!\left({
\bbE\!\left[\left\{{ 
\!\left( \frac{1}{2\delta}\right)^j \repint_{y-\delta}^{y+\delta} \sigma_{d-j}(C_{X}^{D^{r}}(z))\ud z}\right\}^2\right]
 }\right)^{\frac{1}{4}}.
\end{multline}
Let us consider the second term of the product in the right-hand side of the last inequality.\\
For that, let us notice that if $Y$ satisfies (\ref{egalite Y}), and if $X$ and $Y$ satisfy one of the three assumption $\bA_{i}$, $i=1, 2, 3$, or if $X$ and $Y$ satisfy the assumption ${\bB_{4}^\star}$, it is easy to prove, as in   Proposition \ref{H2etH5}  that
$X$ and $Y^2$ still satisfy the hypothesis ${\bH _1}$ and ${\bH _3}$.\\
By  Proposition  \ref{dualité} and since  $\forall n \in \naturels^*$, $\forall x \in D, \abs{Y^{(n)}(x)} \le \abs{Y(x)}$, we get the following inequalities, $\forall n \in \naturels^*$, $\forall y \in \reels^j$,
\begin{multline*}
\limsup_{\delta \to 0} \!\left(\frac{1}{2\delta}\right)^j \repint_{y-\delta}^{y+\delta}
\bbE\left[\int_{\calC_{X}^{D^{r}}(z)} (Y^{(n)}(x))^{2} \ud\sigma_{d-j}(x)\right] \ud z\\
\specialpos{\quad \le\limsup_{\delta \to 0} \!\left(\frac{1}{2\delta}\right)^j \repint_{y-\delta}^{y+\delta}
\bbE\left[\int_{\calC_{X}^{D^{r}}(z)} Y^{2}(x) \ud\sigma_{d-j}(x)\right] \ud z\hfill}\\
 \le\limsup_{\delta \to 0} \!\left(\frac{1}{2\delta}\right)^j \repint_{y-\delta}^{y+\delta} \int_{D}
{p}_{X(x)}(z) 
\bbE\left[Y^{2}(x)H(\nabla
X(x))\given{X(x)=z} \right]\ud x \ud z\\
\specialpos{\quad\substack{\longrightarrow\\\delta \to 0} \int_{D}
{p}_{X(x)}(y) 
\bbE\left[Y^{2}(x)H(\nabla
X(x))\given{X(x)=y}\right]\ud x < \infty.\hfill}
\end{multline*}
The last convergence comes from the fact that  $X$ and $Y^2$ satisfy the hypothesis ${\bH _3}$.
We will now study the third term of the product on the right-hand side of (\ref{autre}).\\
By Remark \ref{finitude} following Theorem \ref{coa1}, if we apply the coarea formula to the functions   $G:=X$ and  $f := \1_{\{\prod_{i=1}^{j} [y_{i}-\delta, y_{i}+\delta]\}} \ge 0$ and to the Borel set $B:=D$, $\forall y \in \reels^{j}$, we obtain
$$
 \repint_{y-\delta}^{y+\delta} \sigma_{d-j}(C_{X}^{D^{r}}(z)) \ud z=
  \int_{D} \1_{\{X(x) \in \prod_{i=1}^{j}[y_{i}-\delta, y_{i}+\delta] \}} H(\nabla\!X(x)) \ud x.
$$%
If $Y$ satisfies (\ref{egalite Y}) and if $X$ and $Y$ satisfy one of the three assumptions ${\bA _{i}}$, $i=1, 2, 3$, the hypothesis in assumption ${\bB_{i}^\star}$, $i=1, 2, 3$, ensures that for almost all $(x_1, x_2)  \in D \times D$, the density of the vector $(X(x_1), X(x_2))$ exists.
To be convinced of this, it is enough to calculate the density in the same way as in the proof of Proposition \ref{H2etH5} to establish (\ref{Co}).
Note that we can also use the formula given in \cite[Corollary 4.18, p{.}~68]{AL121221c}.\\
 Moreover, if $X$ and $Y$ satisfy the hypothesis ${\bB_{4}^\star}$, it is clear that the density of the vector  $(X(x_1), X(x_2))$ exists for almost all $(x_1, x_2)  \in D \times D$.\\
 Finally, $\forall y \in \reels^j$, using the hypothesis ${\bH _{6}^\star}$ or  using the fourth assumption appearing in ${\bB_{4}^\star}$ and using the same notational conventions as above, we obtain
\begin{multline*}
 \limsup_{\delta \to 0}\bbE\left[\left({  \left(\frac{1}{2\delta}\right)^j \repint_{y-\delta}^{y+\delta} \sigma_{d-j}(C_{X}^{D^{r}}(z)) \ud z }\right)^2\right] \\
\specialpos{ \quad=\limsup_{\delta \to 0} \!\left(\frac{1}{2\delta}\right)^{2j} \repint_{y-\delta}^{y+\delta}  \repint_{y-\delta}^{y+\delta} 
 \int_{D \times D} {p}_{X(x_1), X(x_2)}(z_1, z_2) \hfill}\\
\specialpos{ \qquad\times\bbE\!\left[H(\nabla\!X(x_1)) H(\nabla\!X(x_2))\given{X(x_1)=z_1, X(x_2)=z_2}\right] \ud x_1 \ud x_2 \ud z_1 \ud z_2\hfill}\\
\specialpos{ \quad\le \bC .\hfill}
\end{multline*}
Before taking the limit when $\delta$ tends to zero in (\ref{autre}), let us observe that as in the proof of Remark \ref{moduledeY}, $X$ and $Z^{(n)}$ satisfy the hypothesis  ${\bH _2}$ (and ${\bH _{4}}$).
Moreover, we saw in the beginning of this proof that $\forall n \in \naturels^*$, $X$ and $Y^{(n)}$ satisfy the hypothesis ${\bH _3}$ (and ${\bH _5}$).\\
Taking the limit when $\delta$ tends to zero, it turns out that  $\forall n \in \naturels^*$ and $\forall y \in \reels^{j}$ 
\begin{multline*}
\left\vert
	\bbE\left[\int_{\calC_{X}^{D^{r}}(y)}\abs{Z^{(n)}(x)} \ud\sigma_{d-j}(x)\right]
	\right. \\ 
	\specialpos{\hfill-\left.
	\int_D{p}_{X(x)}(y)\bbE\left[\abs{Y^{(n)}(x)} H(\nabla\!X(x))\given{X(x)=y}\right] \ud x \rule{0pt}{18pt}
\right\vert}\\
\specialpos{\quad\le \bC 
\!\left({
\bbE\!\left[(1- \Psi(L_{X}(\bm{\cdot})/n))^{4}\right]}\right)^{\frac{1}{4}}\hfill}\\
\times
\!\left({
\int_{D}
{p}_{X(x)}(y) 
\bbE\left[Y^{2}(x)H(\nabla
X(x))\given{X(x)=y}\right] \ud x
}\right)^{\frac{1}{2}}.
\end{multline*}
The idea now is to take the limit when $n$ tends to infinity in the last inequality.
First, by using Lebesgue's dominated convergence theorem, 
$$
\lim\limits_{n \to +\infty}\bbE[(1- \Psi(L_{X}(\bm{\cdot})/n))^{4}]=0.
$$%
Hence $\forall y \in \reels^{j}$:
\begin{multline*}
	\limsup_{n \to +\infty} \left\vert
	\bbE\left[\int_{\calC_{X}^{D^{r}}(y)} \abs{Z^{(n)}(x)} \ud\sigma_{d-j}(x)\right]
	\right.
	\\
	 -\left.
	\int_D{p}_{X(x)}(y) 
	\bbE\left[\abs{Y^{(n)}(x)} H(\nabla\!X(x))\given{X(x)=y}\right]\ud x
	\rule{0pt}{18pt}\right\vert
	=0.
\end{multline*}
Moreover, using Beppo Levi's theorem,  $\forall y \in \reels^{j}$,
\begin{multline*}
	\lim_{n \to +\infty} \uparrow \int_D{p}_{X(x)}(y) 
	\bbE\left[\abs{Y^{(n)}(x)} H(\nabla
	X(x))\given{X(x)=y}\right]\ud x\\
	=\int_D{p}_{X(x)}(y) 
	\bbE\left[\abs{Y(x)} H(\nabla
	X(x))\given{X(x)=y}\right]\ud x< \infty.
\end{multline*}
Proposition \ref{H2etH5} ensures that the last integral is finite.\\
Beppo Levi's theorem also implies that, $\forall y \in \reels^{j}$,
$$
	\lim_{n \to +\infty} \uparrow
	\bbE\left[\int_{\calC_{X}^{D^{r}}(y)}\abs{Z^{(n)}(x)} \ud \sigma_{d-j}(x)\right]=
	\bbE\left[\int_{\calC_{X}^{D^{r}}(y)}\abs{Y(x)} \ud\sigma_{d-j}(x)\right].
$$%
Then, $\forall y \in \reels^{j}$,
\begin{equation}
\label{limitefinie}
\bbE\left[\int_{\calC_{X}^{D^{r}}(y)}\abs{Y(x)} \ud\sigma_{d-j}(x)\right]< \infty,
\end{equation}
and also for $\forall y \in \reels^{j}$,
\begin{multline*}
	\bbE\left[\int_{\calC_{X}^{D^{r}}(y)} \abs{Y(x)} \ud\sigma_{d-j}(x)\right] \\
	=\int_D{p}_{X(x)}(y)\bbE\left[\abs{Y(x)} H(\nabla\!X(x))\given{X(x)=y}\right] \ud x.
\end{multline*}
In the same way as in the previous proof,  replacing $\abs{Z^{(n)}}$ by $Z^{(n)}$ and 
$\abs{Y^{(n)}}$ by $Y^{(n)}$, we obtain  $\forall y \in \reels^{j}$:
\begin{multline*}
\limsup_{n \to +\infty} \left\vert
\bbE\left[\int_{\calC_{X}^{D^{r}}(y)}Z^{(n)}(x) \ud\sigma_{d-j}(x)\right]  \right.
\\
 \left.
-\int_D{p}_{X(x)}(y) 
\bbE\left[Y^{(n)}(x)H(\nabla
X(x))\given{X(x)=y}\right]\ud x
\rule{0pt}{18pt}\right\vert
=0.
\end{multline*}
By (\ref{limitefinie}), Lebesgue's dominated convergence theorem implies, $\forall y \in \reels^{j}$ 
$$
	\lim_{n \to +\infty}\bbE\left[\int_{\calC_{X}^{D^{r}}(y)}Z^{(n)}(x) \ud\sigma_{d-j}(x)\right]
		=\bbE\left[\int_{\calC_{X}^{D^{r}}(y)}Y(x) \ud\sigma_{d-j}(x)\right].
$$%
Also $\forall y \in \reels^{j}$,
\begin{multline*}
	\lim_{n \to +\infty}  \int_D{p}_{X(x)}(y)\bbE\left[Y^{(n)}(x)H(\nabla\!X(x))\given{X(x)=y}\right]\ud x \\
	=\int_D{p}_{X(x)}(y) \bbE\left[Y(x)H(\nabla\!X(x))\given{X(x)=y} \right]\ud x.
\end{multline*}
This completes the proof of the theorem in the case where $Y$ satisfies (\ref{egalite Y}), and $X$ and $Y$ satisfy one of the three assumptions ${\bB_{i}^\star}$, $i=1, 2, 3$, where we replaced in the assumption ${\bB_i^\star}$ the assumption $\bB_{i}$ by $\bA_i $.
If $X$ and $Y$ satisfy the assumption ${\bB_4^\star}$, the theorem is true.\\
Now suppose that $Y$ satisfies (\ref{egalite Y}) and $X$ and $Y$ satisfy one of the three assumptions ${\bB_{i}^\star}$, $i=1, 2, 3$.
We will proceed as in the proof of  Theorem \ref{Hyp B}.
Let us consider  $\forall n \in \naturels^{\star}$, the sets $D_n=\left\{x \in \reels^d, d(x, D^{c}) > \frac{1}{n} \right\}$.
\\
For any $n$ in $\naturels^{\star}$, $D_n$ is an open set contained in  $D$.
We consider the restrictions $X_{|D_n}$ and $Y_{|D_n}$.
It is clear that if $Y$ satisfies (\ref{egalite Y}) and if $X$ and $Y$ satisfy one of the assumptions ${\bB_{i}^\star}$, $i=1, 2, 3$, then $\forall n \in \naturels^{\star}$, $Y_{|D_n}$ satisfies (\ref{egalite Y}) and $X_{|D_n}$ and $Y_{|D_n}$ satisfy one of the assumptions  ${\bB_{i}^\star}$, $i=1, 2, 3$.
In ${\bB_i^\star}$ the assumption ${\bB_i}$ is replaced  by $\bA_i $.
Also, $X_{|D_n}$ and $Y_{|D_n}$ satisfy the hypothesis ${\bH _6^\star}$, where the  open set $D$ is replaced by the open set $D_n$.\\
We apply Theorem \ref{Cabter} to $X_{|D_n}$ and to $Y_{|D_n}$ (resp$.$
$\abs{Y}_{|D_n}$) and obtain $\forall y \in \reels^j$ and  $\forall n \in \naturels^*$,
\begin{multline*}
	\bbE\left[\int_{\calC_{D_n, X}^{D^{r}}(y)}Y(x) \ud\sigma_{d-j}(x)\right]\\
	=\int_{D_{n}} {p}_{X(x)}(y)\bbE\left[Y(x)H(\nabla\!X(x))\given{X(x)=y}\right]\ud x, 
\end{multline*}
similarly replacing  $Y$ by $\abs{Y}$.\\
The hypothesis ${\bH _6}$ allows, when $n\to\infty$, to apply Lebesgue's dominated convergence theorem in the above equality.\\
This ends the proof of the theorem.
\end{proofarg}
\section{General Rice formulas for all level}
\subsection{Preliminaries for the general  Rice formula}
The following two propositions proved by Aza\"{\i}s \& Wschebor \cite[p.\,178-179]{MR2478201} will provide the  arguments to obtain a general Rice formula for a random field.
Our goal is to obtain a Rice formula for any level $y \in \reels$, the level not necessarily being regular.
\spacebefore
\begin{prop} 
\label{azws}
Let  $Z: \Omega \times W \subset \Omega \times \reels^{\ell} \to \reels^{m+k}$ ($m \in \naturels, k \in \naturels^{\star}$), be a random field in $C^{1}(W, \reels^{m+k})$, $W$ an open set of $\reels^{\ell}$, $J$ a compact subset of $W$ whose Hausdorff dimension is less or equal to $m$ and $z_{0}\in\reels^{m+k}$ is fixed.
We assume that
 $Z$ satisfies the following assumption:
$\forall t \in J$, the random vector $Z(t)$ has a density
$p_{Z(t)}(v) \such \exists \bC  >0$, and a neighborhood $V_{z_{0}}$ of $z_{0} \such \forall t \in J$ and for all $v \in V_{z_{0}}$, $p_{Z(t)}(v)\le \bC $.\\
Then almost surely there is no point  $t \in J$ such that
$Z(t)=z_{0}$.
\end{prop}
\spacebefore
\begin{proofarg}{Proof of Proposition \ref{azws}}
Let $z_{0}\in\reels^{m+k}$ fixed.
For $T$, a Borel set of $\reels^{\ell}$ contained in $W$, let
$$
R_{z_{0}}(T):=\{\omega \in \Omega: \exists  x \in T,~Z(x)(\omega)=z_{0}\}.
$$%
Let $J$ be a compact set contained in $W$ whose Haussdorff dimension is less than or equal to $m$.
Since  $k \in \naturels^{\star}$, the Euclidian  Haussdorff measure of  $J$ of dimension $m+k$ is zero, that is $H_{m+k}(J)=0$ (cf$.$ definition in \cite{AL121221c}).
By definition of the Haussdorff Euclidian pre-measure of $J$ which defines $H_{m+k}(J)$, \ie $H_{m+k}^{\delta}(J)$, we have
$$
H_{m+k}(J)=0= \lim_{\delta \to 0} H_{m+k}^{\delta}(J).
$$%
Consider $\ve >0$ and $\eta >0$ fixed.
There exists $\delta_\ve >0 \such \forall \delta \le \delta_\ve$, there exists a countable set $I$ and $(r_{i})_{i\in I}$, $0<r_{i} \le \delta,~\forall i \in I$ such that
$$
	J \subset \bigcup_{i \in I} B(x_{i}, r_{i})\text{~and~} \sum\limits_{i \in I} r_{i}^{m+k} \le \ve.
$$%
Moreover, since $W$ is open in $\reels^{\ell}$, $\forall y \in J \subset W$, there exists $r_{y} >0 \such B(y, 2r_{y}) \subset \overline{B}(y, 2r_{y}) \subset W$.\\
Given that $J \subset \cup_{y \in J} B(y, r_{y})$ and that $J$ is compact, there exists a finite covering $(B(y_{j}, r_{y_{j}}))_{j=1,n}$ satisfying 
$J \subset \cup_{j=1}^{n} B(y_{j}, r_{y_{j}})$, $y_{j} \in J$ for all $j=1,\dots, n$.
Consider $r:= \inf_{j=1,\dots, n} r_{y_{j}}$ and $C$ the compact set defined by $C := \cup_{j=1}^{n} \overline{B}(y_{j}, 2r_{y_{j}})\subset W$.
\\
Let $R_{\ve, \eta}:= \inf(\delta_\ve, r/2, {\mu}/{(2\eta)})$ where $\mu$ is the constant defining the  neighborhood of $z_{0}$, where the density of $Z$ is bounded.
This latter  neighborhood satisfies that $\forall t \in J$ the random vector $Z(t)$ has a density
$p_{Z(t)}(v)$ satisfying $p_{Z(t)}(v)\le \bC $, for  
$v \such \normp[m+k]{v-z_{0}} \le \mu$.\\
Thus there exists a countable set $I$ and $(r_{i})_{i\in I}$, $0<r_{i} \le R_{\ve, \eta}$ $\forall i \in I$ satisfying  %
$$
J \subset \bigcup_{i \in I} B(x_{i}, r_{i})\text{~and~}\sum\limits_{i \in I} r_{i}^{m+k} \le \ve.
$$%
We have %
\begin{multline*}
  \bbP(R_{z_{0}}(J)) \le \bbP(\sup_{t \in C} \normp[m+k,\ell]{\nabla\!Z(t)} > \eta)\\
  + \sum_{i \in I} \bbP\!\left({\left\{{\sup_{t \in C} \normp[m+k,\ell]{\nabla\!Z(t)} \le \eta }\right\} \cap R_{z_{0}}(B(x_{i}, r_i) \cap J)}\right).
\end{multline*}

Let $i$ fixed in $I$.
If $B(x_{i}, r_i) \cap J = \emptyset$, then $R_{z_{0}}(B(x_{i}, r_i) \cap J)= \emptyset$ and $\bbP\!\left({\left\{{\sup_{t \in C} \normp[j,d]{\nabla\!Z(t)} \le \eta }\right\} \cap R_{z_{0}}(B(x_{i}, r_i) \cap J)}\right)=0$.
If $B(x_{i}, r_i) \cap J \neq \emptyset$, let us fix $z \in B(x_{i}, r_i) \cap J$.\\
For $\omega \in {\left\{{\sup_{t \in C} \normp[m+k,\ell]{\nabla\!Z(t)} \le \eta }\right\}} \cap R_{z_{0}}(B(x_{i}, r_i) \cap J)$ there exists  $x \in B(x_{i}, r_i) \cap J \such Z(x)(\omega)= z_{0}$.
\\
Let us remark that there exists  $j=1,\dots, n \such x$ and $z$ belong to the ball $B(y_{j}, 2r_{y_{j}})$, this entails that for all  $\lambda \in \intrff{0}{1}$, $\lambda x + (1- \lambda) z \in C$.\\
Indeed, since $x \in J$, there exists $j\in\{1,\dots, n\}$, such that $x \in B(y_{j}, r_{y_{j}})$.
We have the following inequalities
\begin{align*}
\normp[\ell]{z-y_{j}} &\le \normp[\ell]{z-x_{i}} + \normp[\ell]{x_{i}-x} + \normp[\ell]{x-y_{j}} \\
&\le 2r_{i}+r_{y_{j}}
\le 2R_{\ve, \eta} +r_{y_{j}} \le r +r_{y_{j}}\le 2r_{y_{j}}.
\end{align*}
Furthermore, since $Z$ is $C^1(W, \reels^{m+k})$, it is also  $C^1(B(y_{j}, 2r_{y_{j}})$, $\reels^{m+k})$.
Since the ball $B(y_{j}, 2r_{y_{j}})$ is an open convex set, we have
\begin{align*}
	Z(z)(\omega)-Z(x)(\omega) &=  Z(z)(\omega)-z_{0}\\
	 	&= \left[{\int_{0}^{1} \nabla\!Z(\lambda x + (1- \lambda) z)(\omega) \ud\lambda }\right] (z-x).
 \end{align*}
Consequently, as $\forall \lambda \in \intrff{0}{1}$ we have $\lambda x + (1- \lambda) z \in C$ then
$$
  \normp[m+k]{Z(z)(\omega)-z_{0}} \le \eta \normp[\ell]{z-x} \le 2 \eta r_{i} \le 2 \eta R_{\ve, \eta} \le \mu.
$$%
Hence
\begin{align*}
	\MoveEqLeft[1]{\bbP\!\left({\left\{{\sup_{t \in C} \normp[m+k,\ell]{\nabla\!Z(t)} \le \eta }\right\} \cap R_{z_{0}}(B(x_{i}, r_i) \cap J)}\right)}\\
	 &\quad\le\bbP(\omega,   \normp[m+k]{Z(z)(\omega)-z_{0}}  \le 2 \eta r_{i})\\
	 &\quad =\int_{\reels^{m +k}} \1_{\{\normp[m+k]{v - z_{0}} \le 2 \eta r_{i}\}} p_{Z(z)}(v) \ud v \le 
\bC  D_{m, k} (\eta r_{i})^{m+k}.
\end{align*}
 Finally, we have shown that $\forall \ve >0$, $\forall \eta >0$,
  \begin{align*}
	\bbP(R_{z_{0}}(J)) &\le \bbP\!\left(\sup_{t \in C} \normp[m+k, \ell]{\nabla\!Z(t)} > \eta\right) + \bC  D_{m, k} \eta^{m+k} \sum_{i \in I} r_{i}^{m+k}\\
	 &\le \bbP\!\left(\sup_{t \in C} \normp[m+k, \ell]{\nabla\!Z(t)} > \eta\right) + \bC  D_{m, k} \eta^{m+k} \ve.
 \end{align*}
By taking limits when $\ve $ tends to zero then when $\eta$ tends to infinity, in this order, we get $\bbP(R_{z_{0}}(J))=0$.
\end{proofarg}

We are now in a position to state the second proposition.
\spacebefore
\begin{prop}
\label{wsch}
Let $X: \Omega \times D \subset \Omega\times \reels^{d} \rightarrow \reels^{j}$ ($j \le d$)  be a random field in $C^{2}(D,\reels^{j})$, where $D$ is an open set of  $\reels^{d}$ and let $D_{0}$ be a compact of $\reels^{d}$ contained in $D$.
Let $y \in\reels^{j}$ fixed.
We assume that $X$ satisfies the following assumption ${(\bS )}$:
\begin{itemize}
\item ${(\bS )}$ $\forall (x, \lambda) \in D_{0} \times S^{j-1}$ the random vector
$$
(X(x), \lambda \cdot \nabla\!X(x)),
$$%
has a density $p_{X(x),  \lambda \cdot \nabla\!X(x)}(u, w)$, such that there exists a constant $\bC  >0$, a neighborhood $V_{y}$ of $y$ and a neighborhood $V_{\vec{0}_{d}}$ of $\vec{0}_{\reels^{d}}$, such that $\forall x \in D_{0}$ and $\forall \lambda \in S^{j-1}$, $\forall u \in V_{y}$ and $\forall w \in V_{\vec{0}_{d}}$, $p_{X(x),  \lambda \cdot \nabla\!X(x)}(u, w)\le \bC $.
\end{itemize}
Then
$$
\bbP\left\{\omega \in \Omega: \exists x \in D_{0}, X(x)(\omega)=y ,\rank  \nabla\!X(x)(\omega)<j\right\}=0.
$$%
\end{prop}
\spacebefore
\begin{proofarg}{Proof of Proposition \ref{wsch}}
Let us define the random field  $Z$ by\\
\centerline{$Z: \Omega \times D \times
\reels^{j} \subset \Omega \times\reels^d\times
\reels^{j}\longrightarrow\reels^j\times\reels^d$}
such that
$$
Z(x,\lambda):=(X(x), \lambda \cdot \nabla\!X(x)).
$$%
Consider the set $W := D \times \reels^j$ which is an open set of $\reels^{\ell}$, where $\ell :=d+j$.
The field $Z$ is $C^{1}(W, \reels^j \times \reels^d)$, since $X$ is $C^2(D, \reels^j)$, and takes its values in $\reels^{m+k}$ where $m:= j+d-1$ and $k:=1$.\\
Consider $J:=D_{0} \times S^{j-1}$ a compact set contained in  $W$.
Its Haussdorff measure is less or equal to $m$.
Let $z_{0} :=(y, \vec{0}_{\reels^{d}}) \in \reels^{m+k}$ be fixed.\\
Since $\forall (x, \lambda) \in D_{0} \times S^{j-1}$ the random vector $(X(x),  
\lambda \cdot\nabla\!X(x))$ has a bounded density
$p_{X(x), \lambda \cdot \nabla\!X(x)}(u, w)$, for  
$u$ in a  neighborhood of $y$ and  $w$ in a neighborhood of $\vec{0}_{\reels^{d}}$ then
$\forall t \in J$ the random vector $Z(t)$ has a density $p_{Z(t)}(v)$ satisfying  $p_{Z(t)}(v)\le \bC $, for
$v$ in a neighborhood of $z_{0}$.\\
 By Proposition \ref{azws},
$$
\bbP\{\omega \in \Omega: \exists
(x,\lambda) \in D_{0} \times S^{j-1}, (X(x)(\omega), \lambda \cdot \nabla\!X(x)(\omega))=(y, \vec{0}_{\reels^{d}})\}=0
$$
then
$$
\bbP\{\omega \in \Omega: \exists x \in D_{0}, X(x)(\omega)=y ,\rank  \nabla\!X(x)(\omega)<j\}=0.
$$%
This completes the proof of the proposition.
\end{proofarg}
We now have all the ingredients to prove Rice's general formula for all levels.
\subsection{The general Rice formula}
In this section Theorem \ref{Cabana4} provides a general Rice formula for any level, not necessarily  regular.
Note that Theorem 6.10  of \cite{MR2478201} gives the same result but in this book, the proofs are only sketched.\\
The proof of Theorem \ref{Cabana4} will be based on the proof of Theorem \ref{Hyp B}.
Therefore, its proof will require more general assumptions than those denoted by $\bB_{i}$, $i=1,\dots,4$, which appear in this last theorem.
Let us therefore state the following hypotheses $\bC _{i}$, $i=1,\dots,4$ with the same previous convention.
For the three first assumptions  $\bC _{i}$, $i=1,\dots,3$, the process $Y$ will be expressed using (\ref{egalite Y}).\\ 
In what follows we have the assumptions:

\begin{itemize}
\item $\bC _1$: It is the assumption $\bB_1$ plus the following hypothesis:
$\forall x \in D$ the vector $(X(x), \nabla\!X(x))$ has a density.

\item $\bC _2$:  It is the assumption $\bB_2$ plus the following hypothesis.
$\forall x \in D$, the vector $(Z(x), \nabla\!Z(x))$ has a density.

\item $\bC _{3}$:  It is the assumption  $\bB_{3}$ plus the following hypothesis:
the function $F$ verifies assumption ${(\bF\bF)}$, \ie
\begin{itemize}
\item ${(\bF\bF)}$ $\forall y \in\reels^{j}$, $\exists \bC  >0$, there exists a neighborhood $V_{y}$ of $y$ such that $\forall \mu >0$ and $\forall u \in V_{y}$ we have
$$
\int_{\reels^{j^{\prime}-j}} 
 	\frac{1}{\abs{J_{F}(F_{z}^{-1}(u), z)}^{d+1}}
		e^{- \mu \normp[j^{\prime}-j]{z}^{2}} \normp[j,j^{\prime}]{\nabla F(F_{z}^{-1}(u), z)}^{(j-1)d} \ud z
 		\le \bC .
$$%
\end{itemize}
\item $\bC _{4}$: It is the assumption  $\bB_{4}$ plus the following hypothesis:
the process $X$ verifies the assumption ${(\bS )}$ of Proposition \ref{wsch}.
\end{itemize}
Theorem \ref{Cabana4} states the general Rice formula for all level.

\spacebefore\begin{theo}
\label{Cabana4}
Let $X: \Omega \times D \subset \Omega\times \reels^{d} \rightarrow \reels^{j}$ ($j \le d$)  be a random field in $C^{2}(D,\reels^{j})$, where $D$ is a bounded convex open set of $\reels^{d}$.
We assume that for almost all $\omega \in \Omega$, $\nabla\!X(\omega)$ is Lipschitz with  Lipschitz constant $L_{X}(\omega)$ such that  $\bbE\left[L_X^d(\bm{\cdot})\right] < \infty$.
 Let $Y: \Omega \times D \subset \Omega \times \reels^{d} \rightarrow \reels$ be a continuous process.\\
If $Y$ satisfies (\ref{egalite Y}) and if $X$ and $Y$ satisfy one of the three assumptions $\bC _{i}$, $i=1, 2, 3$ and the hypothesis ${\bH _6}$ or if $X$ and $Y$ satisfy hypothesis $\bC_4$ then $\forall y \in \reels^{j}$ we have
$$
\bbE\left[\int_{\calC_{X}(y)}Y(x) \ud\sigma_{d-j}(x)\right]=\int_{D}{p}_{X(x)}(y) 
\bbE\left[Y(x)H(\nabla\!X(x))\given{X(x)=y}\right] \ud x.
$$%
\end{theo}
\spacebefore
\begin{rema}
\label{Lipschitz}
We can replace the hypothesis: ``for almost all $\omega \in \Omega$, $\nabla\!X(\omega)$ is Lipschitz with Lipschitz constant $L_{X}(\omega)$ having a moment of order $d$\kern1pt'', with the assumption:
$$
	\bbE\!\left[\left(\sup_{x \in D} \normp[j, d]{\nabla^2(X(x))}^{(s)}\right)^d\right] < \infty.
$$%
Indeed, if $X$ is $C^2$ on $D$, the Taylor formula on  $D$ convex and open set  allows to conclude that, since almost surely
$$
	L_X:=\sup_{x \in D} \normp[j,d]{\nabla^2(X(x))}^{(s)} < \infty,
$$%
then $\nabla\!X$ is almost surely Lipschitz with Lipschitz constant $L_X$.
\hfill$\bullet$
\end{rema}
\spacebefore
\begin{rema}
We can generalize this theorem and the remark to the case where $D$ is an open convex set not necessarily bounded.
\hfill$\bullet$
\end{rema}
\spacebefore
\begin{proofarg}{Proof of Theorem \ref{Cabana4}}By Theorem \ref{Hyp B}, we already know that $\forall y \in \reels^j$
$$
\bbE\left[\int_{\calC_{X}^{D^{r}}(y)}Y(x) \ud\sigma_{d-j}(x)\right]
	=\int_{D} {p}_{X(x)}(y)\bbE\left[Y(x)H(\nabla\!X(x))\given{X(x)=y}\right] \ud x.
$$%
Let us check that the assumption $\bS$ holds.
That is, let us prove that if $D_{0}$ is a compact set contained in $D$ and if $X$ and $Y$ satisfy one of the assumptions
$\bC _i$, $i=1,\dots,3$
then  $\forall (x, \lambda) \in D_{0} \times S^{j-1}$ the random vector 
$(X(x), \lambda \cdot \nabla\!X(x)) $
has a density $p_{X(x),  \lambda \cdot \nabla\!X(x)}(u, w) \such \forall y \in\reels^{j}$ there exists a constant $\bC  >0$, there exists a neighborhood $V_{y}$ of $y$ and a neighborhood $V_{\vec{0}_{d}}$ of $\vec{0}_{\reels^{d}} \such \forall x \in D_{0}$ and $\forall \lambda \in S^{j-1}$, $\forall u \in V_{y}$ and $\forall w \in V_{\vec{0}_{d}}$ it is true that $p_{X(x),  \lambda \cdot \nabla\!X(x)}(u, w)\le \bC $.\\
Note that the last conclusion is true when the processes $X$ and $Y$ satisfy the assumption $\bC_4$.\\
In this case, using Proposition \ref{wsch}, we will deduce that for any compact set $D_{0}$ contained in $D$ and $\forall y \in \reels^j$, 
$$
\bbP\left\{\omega \in \Omega: \exists x \in D_{0}, X(x)(\omega)=y ,\rank  \nabla\!X(x)(\omega)<j\right\}=0.
$$%
By choosing the  compact $D_{0}:=D^{(n)}=\left\{x \in \reels^d, d(x, D^{c}) \ge \frac{1}{n} \right\} \subseteq D$ we will deduce, since $D^{(n)}$ tends in a nondecreasing way towards $D$, when $n\to\infty$, that $\forall y \in \reels^j$, 
$$
\bbP\{\omega \in \Omega: \exists x \in D, X(x)(\omega)=y ,\rank  \nabla\!X(x)(\omega)<j\}=0.
$$%
We will then have shown that  $\forall y \in \reels^j$,
\begin{align*}
\bbE\left[\int_{\calC_{X}(y)}Y(x) \ud \sigma_{d-j}(x)\right]&=
\bbE\left[\int_{\calC_{X}^{D^{r}}(y)}Y(x) \ud \sigma_{d-j}(x)\right]\\
&=\int_{D} {p}_{X(x)}(y)\bbE\left[Y(x)H(\nabla
X(x))\given{X(x)=y}\right] \ud x,
\end{align*}
that will end the proof of this theorem.\\
Let us check that the hypothesis $\bS$ is verified, in the case where $X$ and $Y$ satisfy one of the conditions $\textbf{C}_i$, $i=1,\dots,3$.\\
Let $D_{0}$ be a compact set contained in $D$.
For all $x \in D_{0}$ and $\forall \lambda \in S^{j-1}$, in the case where $Y$ satisfies (\ref{egalite Y}) and $X$ and $Y$ satisfy one of the assumptions $\bC _i$, $i=1,\dots,3$, in a first step we will study the density $p_{X(x),  \lambda \cdot\nabla\!X(x)}$ of the vector $(X(x),  \lambda \cdot \nabla\!X(x))$.
We will express the latter in terms of the density $p_{X(x), \nabla\!X(x)}$ of the vector $(X(x), \nabla\!X(x))$ which exists by the proof of Proposition \ref{H2etH5} (cf$.$ (\ref{Co})).
Thus, let us consider $\lambda \in S^{j-1}$, $\lambda:=(\lambda_1, \dots, \lambda_j)$.
There exists $k \in \{1, \dots, j\} \such \abs{\lambda_{k}} \ge \frac{1}{\sqrt{j}}$.
We will assume for example 
$k:=j$ and that $\abs{\lambda_{j}} \ge \frac{1}{\sqrt{j}}$ this will imply $\frac{1}{\abs{\lambda_{j}}} \le \sqrt{j}$.\\
If $u:= (u_1, \dots, u_j) \in \reels^j$ and
$$
s:=(s_{11}, s_{21}, \dots, s_{j1}, s_{12}, s_{22}, \dots, s_{j2}, \dots, s_{1d}, s_{2d}, \dots, s_{jd}) \in \reels^{jd},
$$%
let us do as in the proof of Proposition \ref{H2etH5} (third part) the following change of variables.
Let $K$ be the function defined by
\begin{align*}
K:  \reels^{j} \times  \reels^{jd}    &\longrightarrow \reels^{j} \times \reels^{d} \times \reels^{(j-1)d}\\
\!\left({u, s}\right)   &\longmapsto K(u, s)
	:= \!\left(u, \sum_{i=1}^{j}\lambda_{i} s_{i1}, \sum_{i=1}^{j} \lambda_{i} s_{i2}, \dots, \sum_{i=1}^{j} \lambda_{i} s_{id},\right. \\
&\left.{ s_{11}, s_{21}, \dots, s_{j-11}, s_{12}, s_{22}, \dots, s_{j-12}, \dots, s_{1d}, s_{2d}, \dots, s_{j-1d}
	\adjstud{$\ds\sum_{i=1}^{j}$} 
	}\right).
\end{align*}
The Jacobian $J_{K}$ of this transformation is such that $\forall (u, s) \in \reels^{j} \times \reels^{jd}$ 
$$
\abs{J_{K}(u, s)}= \abs{\lambda_{j}}^{d} \neq 0,
$$%
by hypothesis.\\
Thus, $K$ is a $C^1(\reels^j \times \reels^{jd},\reels^j \times \reels^d \times \reels^{(j-1)d})$ one-to-one mapping, as well as its inverse $K^{-1}$ given by
\begin{align*}
	\MoveEqLeft[2]{K^{-1}:  \reels^{j} \times \reels^{d} \times \reels^{(j-1)d}    \longrightarrow \reels^{j} \times \reels^{jd}}\\
	 &(u, s_{j1}, s_{j2},\dots, s_{jd}; s_{11}, s_{21},\dots, s_{j-11},\\
	 & \specialpos{\hfill s_{12}, s_{22},\dots, s_{j-12},\dots, s_{1d}, s_{2d},\dots, s_{j-1d}  )}  \\
	 &\longmapsto\!\left({u, s_{11}, s_{21}, \dots, s_{j-11}, \frac{1}{\lambda_{j}} \left[-\sum_{i=1}^{j-1}\lambda_{i} s_{i1}+s_{j1}\right],s_{12}, s_{22}, \dots, s_{j-12}, }\right.\\
	 &\left. \frac{1}{\lambda_{j}} \left[-\sum_{i=1}^{j-1}\lambda_{i} s_{i2}+s_{j2}\right],
		\dots, s_{1d}, s_{2d}, \dots, s_{j-1d}, \frac{1}{\lambda_{j}} \left[-\sum_{i=1}^{j-1}\lambda_{i} s_{id}+s_{jd}\right]\right).
\end{align*}
For any $\lambda \in S^{j-1}$, $\forall x \in D_{0}$ we have 
$$
K(X(x), \nabla\!X(x))= \!\left(X(x), \lambda \cdot \nabla\!X(x), (\nabla\!X(x))_{(j-1)d}\right),
$$%
where if $s \in \reels^{jd}$  we denoted $s_{(j-1)d}$ by
$$
s_{(j-1)d}:=\!\left(s_{11}, s_{21}, \dots, s_{j-11}, s_{12}, s_{22}, \dots, s_{j-12}, \dots, s_{1d}, s_{2d}, \dots, s_{j-1d}\right).
$$%
 With these notations, if $p_{X(x), \lambda \cdot \nabla\!X(x), (\nabla\!X(x))_{(j-1)d}}$ denotes the density of the vector $(X(x), \lambda \cdot \nabla\!X(x), (\nabla\!X(x))_{(j-1)d})$, we have: $\forall \lambda \in S^{j-1}$, $\forall x \in D_{0}$, $\forall (u, s) \in \reels^{j} \times \reels^{jd}$,
$$
	p_{X(x), \lambda \cdot \nabla\!X(x), (\nabla\!X(x))_{(j-1)d}}(u, s) = \frac{1}{\abs{\lambda_{j}}^{d}}  
	p_{X(x), \nabla\!X(x)}\!\left(K^{-1}(u, s)\right).
$$%
We deduce that
$$
	\forall \lambda \in S^{j-1}, \forall x \in D_{0}, \forall (u, w) \in \reels^{j} \times \reels^{d}, w:=(w_1, w_2, \dots, w_d),
$$%
\begin{multline}
\label{CoBis}
	p_{X(x), \lambda \cdot \nabla\!X(x)}(u, w) = \frac{1}{\abs{\lambda_{j}}^{d}} 
	\int_{\reels^{(j-1)d}} p_{X(x), \nabla\!X(x)}\!\left(u, s_{11}, s_{21},\dots, s_{j-1 1},\adjstud{$\ds\sum_{i=1}^{j-1}$}\right.\\
		\frac{1}{\lambda_{j}} \left[-\sum_{i=1}^{j-1} \lambda_{i} s_{i1}+w_1\right], s_{12}, s_{22},\dots, s_{j-12},
	\frac{1}{\lambda_{j}} \left[-\sum_{i=1}^{j-1} \lambda_{i} s_{i2}+w_2\right],\dots, s_{1d}, \\
	\left.s_{2d},\dots, s_{j-1d}, \frac{1}{\lambda_{j}} \left[-\sum_{i=1}^{j-1} \lambda_{i} s_{id}+w_d\right]\right) \ud s_{(j-1)d}.
 \end{multline}
We now need find an upper bound for this density.
For this purpose, let us consider each of $\textbf{C}_i$, $i=1,\dots,3$ assumptions.
 \begin{itemize}
 \item If $X$ and $Y$ satisfy the assumption $\textbf{C}_1$, then $\forall x \in D_{0}$ the vector $(X(x), \nabla\!X(x))$ has a non singular density and since
 $X$ is a process of class $C^{2}$ the covariance matrix of this vector is strictly positive on the compact set $D_{0}$.
Then there exist real numbers $a, b >0$ such that
 $\forall x \in D_{0}$, $0<a \le \inf_{\normp[j(d+1)]{z}=1}\normp[j(d+1)]{\bbV(X(x), \nabla\!X(x)) \times z} \le b$.
In the same way as we obtained the equality (\ref{pZ}), and with the same notations, we prove the existence of a number $\mu >0$ and a number $\bC >0$, such that  $\forall (x, u, s) \in D_{0} \times \reels^{j}\times \reels^{dj}$, 
\begin{align*}
 {p}_{X(x), \nabla\!X(x)}(u, s) 
& \le \bC  e^{- \mu \normp[j(d+1)]{(u, s)}^{2}}\\
& \le  \bC  e^{- \mu \normp[jd]{s}^{2}}\\
& \le  \bC  e^{- \mu \normp[(j-1)d]{s_{(j-1)d}}^{2}}.
\end{align*}
Using the equality (\ref{CoBis}), we obtain the following bound:
$\forall \lambda \in S^{j-1}$, $\forall \bx \in D_{0}$, $\forall (u, w) \in \reels^{j} \times \reels^{d}$, 
\begin{align*}
p_{X(x), \lambda \cdot \nabla\!X(x)}(u, w)  &\le \frac{\bC }{\abs{\lambda_{j}}^{d}} 
\int_{\reels^{(j-1)d}} e^{- \mu \normp[(j-1)d]{s_{(j-1)d}}^{2}} \ud s_{(j-1)d}\\
&\le \bC ,
\end{align*}
the last inequality comes from the fact that  ${1}/{\abs{\lambda_{j}}^d} \le (\sqrt{j})^d$.\\
The assumption ${(\bS )}$ is satisfied.
 \item If $X$ and $Y$ satisfy the assumption $\bC_2$, then $\forall x \in D_{0}$ the vector $(Z(x), \nabla\!Z(x))$ has a non-degenerate density.
In the same way as in the first part, we show the existence of a number $\mu >0$ and of a number $\bC  >0$ such that  $\forall (x, u, s) \in D_{0} \times \reels^{j}\times \reels^{dj}$, 
$$
	{p}_{Z(x), \nabla\!Z(x)}(u, s)  \le  \bC  e^{- \mu \normp[jd]{s}^{2}}.
$$%
Moreover, the equality (\ref{Co}) proved in the third part of the proof of Proposition \ref{H2etH5} and applied to $j:=j^{\prime}$ shows that the density of the vector $(X(x), \nabla\!X(x))$, for all  $(x, u, s) \in D_{0} \times \reels^{j}\times \reels^{dj}$, is given by:
\begin{multline*}
	{p}_{X(x), \nabla\!X(x)}(u, s) = \frac{1}{\abs{J_{F}(F^{-1}(u))}^{d+1}}  \\
		   \times {p}_{Z(x), \nabla\!Z(x)}\!\left(F^{-1}(u), (\nabla F(F^{-1}(u)))^{-1} \times s\right).
 \end{multline*}
 We deduce that there exists a constant $\mu >0$ such that  $\forall (x, u, s) \in D_{0} \times \reels^{j}\times \reels^{dj}$, 
$$
 {p}_{X(x), \nabla\!X(x)}(u, s) \le  \frac{\bC }{\abs{J_{F}(F^{-1}(u))}^{d+1}}  \times 
 e^{- \mu \normp[jd]{(\nabla F(F^{-1}(u)))^{-1} \times s}^{2}}.
$$%
Let $y$ now be a fixed vector in $\reels^j$.
Since the function $F$ is $C^1(\reels^j, \reels^j)$ and $F^{-1}$ is continuous, the Jacobian $J_{F}(F^{-1})$ is continuous on $\reels^j$ and nonzero everywhere.
Let  $V_{y}$ a compact neighborhood of $y$, then  $\exists \bC  >0 \such \forall u \in V_{y}$ we have 
$ \frac{1}{\abs{J_{F}(F^{-1}(u))}^{d+1}} \le \bC $.\\
Moreover, $\forall u \in \reels^j$, $(\nabla F(F^{-1}(u)))^{-1} \in \frakL(\reels^{j}, \reels^{j})$ and thus for all $s \in \reels^{jd}$, 
\begin{equation}
	\normp[jd]{(\nabla F(F^{-1}(u)))^{-1} \times s} 
		\ge \frac{\normp[jd]{s}}{\normp[j,j]{\nabla F(F^{-1}(u))}} .
	\label{majorationnorme}
\end{equation}
Since  $F$ is $C^1(\reels^j, \reels^j)$ and $F^{-1}$ is continuous, the operator \linebreak $\nabla F(F^{-1}(\bm{\cdot}))$ is a continuous function of  $\reels^j$ into $\frakL(\reels^{j}, \reels^{j})$.
There exists a constant $\bC  >0 \such \forall u \in V_{y}$, we have
$$
	\normp[j,j]{\nabla F(F^{-1}(u))} \le \bC .
$$%
We deduce that $\forall s \in \reels^{jd}$ and $\forall u \in V_{y}$, 
$$
	\normp[jd]{(\nabla F(F^{-1}(u)))^{-1} \times s} \ge \bC   \normp[jd]{s}.
$$%
Finally, we conclude that $\forall y \in \reels^j$, there exists a constant $\bC  >0$, a neighborhood $V_{y}$ of $y$ and a constant $\mu >0$ such that $\forall x \in D_{0}$, $\forall u \in V_{y}$ and for all  $s \in \reels^{jd}$ we have
 \begin{eqnarray*}
 {p}_{X(x), \nabla\!X(x)}(u, s) \le \bC   e^{- \mu \normp[jd]{s}^{2}}.
 \end{eqnarray*}
Then in the same way as in the first part of the proof of this theorem, we deduce that $\forall y \in  \reels^j$, there exists a constant $\bC  >0$, a neighborhood  of $y$, (let say $V_{y}$), such that $\forall x \in D_{0}$ and $\forall \lambda \in S^{j-1}$, $\forall u \in V_{y}$ and $\forall w \in \reels^{d}$, 
$$
p_{X(x), \lambda \cdot \nabla\!X(x)}(u, w)  \le \bC .
$$%
Assumption ${(\bS )}$ is satisfied.
 \item If $X$ and $Y$ satisfy the assumption $\textbf{C}_3$, in the same way as before and since $D_{0}$ is a compact set, there exist constants $\mu >0$ and $\bC  >0$ such that
$\forall (x, u, s) \in D_{0} \times \reels^{j}\times \reels^{dj^{\prime}}$, 
$$
 	{p}_{Z(x), \nabla\!Z(x)}(u, s) 
		 \le  \bC   e^{- \mu \normp[j^{\prime}]{u}^{2}}  e^{- \mu \normp[j^{\prime}d]{s}^{2}}.
$$%
The equality (\ref{Co}) proved in the third part of the proof of Proposition \ref{H2etH5}, shows that the  density of the random vector $(X(x), \nabla\!X(x))$ exists.
Using the same notations and (\ref{majorationnorme}), we can prove that the latter density is bounded in the following form, for all 
$
 	\!\left(x, u, s_{(j, d)}\right) \in D_{0} \times \reels^{j}\times \reels^{dj},
$
\begin{align*}
	\MoveEqLeft[1]{{p}_{X(x), \nabla\!X(x)}(u, s_{(j,d)})\le  \bC 
	\int_{\reels^{j^{\prime}-j}} \int_{\reels^{(j^{\prime}-j)d}}
	 \frac{1}{\abs{J_{F}(F_{z}^{-1}(u), z)}^{d+1}}}    \\
		&\specialpos{\hfill \times e^{- \mu \normp[j^{\prime}-j]{z}^{2}} \times  e^{- \mu \normp[(j^{\prime}-j)d]{s_{(j^{\prime}-j, d})}^{2}}\hfill} \\
		&\times e^{- \mu \normp[jj]{
	  [\nabla F(F_{z}^{-1}(u), z)]_{jj}}^{-2} \times 
	 \normp[jd]{s_{j,d} -  \left[\nabla F(F_{z}^{-1}(u), z)\right]_{j j^{\prime}-j} s_{j^{\prime}-j, d}}^{2}} 
	  \ud s_{(j^{\prime}-j,d)} \ud z.
\end{align*}
We deduce that 
$\forall (x, u, s_{(j, d)}) \in D_{0} \times \reels^{j}\times \reels^{dj}$,
\begin{multline*}
	{p}_{X(x), \nabla\!X(x)}(u, s_{(j,d)})\le  \bC 
	\int_{\reels^{j^{\prime}-j}} \int_{\reels^{(j^{\prime}-j)d}}
	 \frac{1}{\abs{J_{F}(F_{z}^{-1}(u), z)}^{d+1}}    \\
	  \times e^{- \mu \normp[j^{\prime}-j]{z}^{2}} \times  e^{- \mu \normp[(j^{\prime}-j)d]{s_{(j^{\prime}-j, d)}}^{2}}\\
	  \times e^{- \mu \normp[jj]{
	  [\nabla F(F_{z}^{-1}(u), z)]_{jj}}^{-2} \times 
	 \normp[(j-1)d]{s_{j-1,d} -  \left[\nabla F(F_{z}^{-1}(u), z)\right]_{j-1 j^{\prime}-j} s_{j^{\prime}-j, d}}^{2}}\\
	  \ud s_{(j^{\prime}-j,d)} \ud z,
\end{multline*}
where the matrix $s_{j-1,d}$ is the matrix  $s_{j,d}$ from which we deleted the $j$-th row and $\left[\nabla F(F_{z}^{-1}(u), z)\right]_{j-1 j^{\prime}-j}$ is the matrix $\left[\nabla F(F_{z}^{-1}(u), z)\right]_{j j^{\prime}-j}$ from which we deleted the $j$-th row.\\
Using(\ref{CoBis}), we get:  \\$\exists \bC  >0, \exists \mu >0, \forall x \in D_{0}$, $\forall \lambda \in S^{j-1}$, $\forall (u, w) \in \reels^{j} \times \reels^{d}$,
\begin{multline*}
{p_{X(x), \lambda \cdot \nabla\!X(x)}(u, w)  \le \bC  \int_{\reels^{(j-1)d}\times\reels^{j^{\prime}-j}\times\reels^{(j^{\prime}-j)d}}
 \frac{1}{\abs{J_{F}(F_{z}^{-1}(u), z)}^{d+1}}} \\
 \times   e^{- \mu \normp[j^{\prime}-j]{z}^{2}} \times  e^{- \mu \normp[(j^{\prime}-j)d]{s_{(j'-j, d)}}^{2}}  \\
  \times e^{- \mu \normp[jj]{
  [\nabla F(F_{z}^{-1}(u), z)]_{jj}}^{-2} 
 \normp[(j-1)d]{s_{j-1,d} -  \left[\nabla F(F_{z}^{-1}(u), z)\right]_{j-1 j^{\prime}-j} s_{j^{\prime}-j, d}}^{2}}\\
 \ud s_{(j^{\prime}-j,d)} \ud z \ud s_{(j-1,d)}.
\end{multline*}
We perform the following change of variables in the integral on  $\reels^{(j-1)d}$:
$$
	s_{j-1,d} -  \left[\nabla F(F_{z}^{-1}(u), z)\right]_{j-1 j^{\prime}-j} s_{j^{\prime}-j, d}
	= \normp[jj]{[\nabla F(F_{z}^{-1}(u), z)]_{jj}}\cdot v_{j-1,d}.
$$%
We get:
  $\exists \bC  >0$, $\exists \mu >0$, $\forall x \in D_{0}$, $\forall \lambda \in S^{j-1}$, $\forall (u, w) \in \reels^{j} \times \reels^{d}$,
 \begin{align*}
\MoveEqLeft[1]{p_{X(x), \lambda \cdot \nabla\!X(x)}(u, w)
  \le \bC  \int_{\reels^{(j-1)d}} \int_{\reels^{j^{\prime}-j}} \int_{\reels^{(j^{\prime}-j)d}}
 \frac{1}{\abs{J_{F}(F_{z}^{-1}(u), z)}^{d+1}}}  \\
	& \specialpos{\hfill \times  e^{- \mu \normp[j^{\prime}-j]{z}^{2}} \times  e^{- \mu \normp[(j^{\prime}-j)d]{s_{(j^{\prime}-j, d)}}^{2}} \hfill}\\
	&  \quad\times  e^{- \mu \normp[(j-1)d]{(v_{j-1,d)}}^{2}}
	\times \normp[jj]{
  [\nabla F(F_{z}^{-1}(u), z)]_{jj}}^{(j-1)d} \ud s_{(j^{\prime}-j,d)} \ud z \ud v_{(j-1,d)}\\
  &\le \bC  \int_{\reels^{j^{\prime}-j}} 
 \frac{1}{\abs{J_{F}(F_{z}^{-1}(u), z)}^{d+1}}
  e^{- \mu \normp[j^{\prime}-j]{z}^{2}} \normp[jj]{
  [\nabla F(F_{z}^{-1}(u), z)]_{jj}}^{(j-1)d} \ud z\\
  &\le \bC  \int_{\reels^{j^{\prime}-j}} 
 \frac{1}{\abs{J_{F}(F_{z}^{-1}(u), z)}^{d+1}}
  e^{- \mu \normp[j^{\prime}-j]{z}^{2}} \normp[jj^{\prime}]{
  \nabla F(F_{z}^{-1}(u), z)}^{(j-1)d} \ud z.
\end{align*}
The  $\textbf{C}_3$ assumption allows  to obtain that
$\forall y \in \reels^{j}$, there exists $\bC  >0$, a neighborhood $V_{y}$ of $y$, such that $\forall x \in D_{0}$, $\forall \lambda \in S^{j-1}$, $\forall u \in V_{y}$ and $\forall w \in \reels^{d}$,
we have
$$
p_{X(x), \lambda \cdot \nabla\!X(x)}(u, w)  \le \bC .
$$%
The hypothesis ${(\bS )}$ is then verified.
This ends the proof of the theorem.
\end{itemize}
 \end{proofarg}%
In the same way as in Theorem \ref{Cabter}, we can get rid of the assumption that $\bbE\left[L_X^d(\bm{\cdot})\right] < \infty$ in Theorem \ref{Cabana4}.
If one of the assumptions $\bC_1$, $\bC_2$, $\bC_3$ or $\bC_4$ is replaced by $\bC_1^\star$, $\bC_2^\star$, $\bC_3^\star$ or $\bC_4^\star$, the Rice's formula is still true.
In the first three assumptions $\bC_1^\star$, $\bC_2^\star$ and $\bC_3^\star$, we make the hypothesis that $Y$ can be written as in (\ref{egalite Y}).

\vbox{More precisely 
\begin{itemize}
\item $\bC_1^\star$: It is the assumption $\bC _1$, plus the following hypothesis:
 for almost all $(x_1, x_2)  \in D \times D$ the density of the vector $(X(x_1), X(x_2))$ exists.
\item $\bC_2^\star$: It is the assumption $\bC _2$ plus the following hypothesis:
 for almost all $(x_1, x_2)  \in D \times D$ the density of the vector $(Z(x_1), Z(x_2))$ exists.
\item $\bC_3^\star$: It is the assumption $\bC _{3}$ plus the following hypothesis:
 for almost all $(x_1, x_2)  \in D \times D$ the density of the vector  $(Z(x_1), Z(x_2))$ exists.
\item $\bC_4^\star$: It is the assumption ${\bB_{4}^\star}$ plus the following hypothesis:
the process $X$ verifies assumption ${(\bS )}$.
\end{itemize}}

Theorem \ref{Ricegenerale} summarizes all the results obtained previously.
This is a new result.
\spacebefore\begin{theo}
\label{Ricegenerale}
Let $X: \Omega \times D \subset \Omega\times \reels^{d} \rightarrow \reels^{j}$ ($j \le d$)  be a random field in $C^{2}(D,\reels^{j})$, where $D$ is a convex open bounded set of  $\reels^{d}$, such that for almost all  $\omega \in \Omega$, $\nabla\!X(\omega)$ is Lipschitz.
 Let $Y: \Omega \times D \subset \Omega \times \reels^{d} \rightarrow \reels$ be a continuous process.
If $Y$ satisfies (\ref{egalite Y}) and if $X$ and $Y$ satisfy one of the three assumptions $\bC_i^\star$, $i=1, 2, 3$ and the hypotheses  ${\bH _6}$ and ${\bH _6^\star}$ or if $X$ and $Y$ satisfy the assumption $\bC_4^\star$, then $\forall y \in  \reels^{j}$ we have
$$
\bbE\left[\int_{\calC_{X}(y)}Y(x) \ud \sigma_{d-j}(x)\right] \\
	=\int_D{p}_{X(x)}(y) 
	\bbE\left[Y(x)H(\nabla\!X(x))\given{X(x)=y}\right]\ud x.
$$%
\spacebefore
\begin{rema}
\label{LipschitzTer}
In the same way as in Remark \ref{Lipschitz}, one can replace in the theorem the hypothesis ``for almost all $\omega \in \Omega$, $\nabla\!X(\omega)$ is Lipschitz'', by the hypothesis
``almost surely $L_X:=\sup_{x \in D} \normp[j,d]{\nabla^2X(x)}^{(s)} < \infty$'', since almost surely  the process $\nabla\!X$ will be Lipschitz  with Lipschitz constant $L_X$.
\hfill$\bullet$
\end{rema}
\spacebefore
\begin{rema}\label{reiteracion}
We can generalized this theorem in the case where $D$ is a convex open not necessarily bounded.%
In what follows, we discuss only the second moment, but see Remark \ref{similaire} for the $k$-th moment.
\hfill$\bullet$
\end{rema}
\end{theo}
\subsection{Rice formula for the $k$-th moment}
Theorem \ref{Cabana4} will allow to state a general Rice formula for the second moment.

Let us assume ${\bD_i}$, $i=1,\dots,4$, where in the first three $\bD_1$, $\bD_2$ and $\bD_3$, we also assume that $Y$ can be written as in (\ref{egalite Y}).\\
We denote $\Delta:= \{(x_1, x_2) \in D \times D,   x_1 = x_2 \}\subset \reels^{2d}$, where $D$ is an open set of  $\reels^d$.
Let us state the following hypotheses $\bD_i$, $i=1, \dots,4$.
\begin{itemize}
\item $\bD_1$: It is the assumption $\bE_1$, plus the following hypothesis:
$\forall x \in D$, the vector $(X(x), \nabla\!X(x))$ has a density.
\item $\bD_2$: It is the assumption  $\bE_2$, plus the following hypothesis:
$\forall x \in D$, the vector $(Z(x), \nabla\!Z(x))$ has a density.
\item $\bD_{3}$:  It is the assumption $\bE_{3}$, plus the following hypothesis:
the function $F$ verifies assumption ${(\bF\bF)}$ appearing in assumption $\bC _{3}$.
\item $\bD_{4}$: It is the assumption $\bE_{4}$,  plus the following hypothesis:
the process $X$ satisfies the assumption ${(\bS )}$.
\end{itemize}

The assumptions $\bE_{i}$, $i=1, \dots,4$, are the following:
\begin{itemize}
\item $\bE_1$: The process $X: \Omega \times D \subset \Omega\times \reels^{d} \rightarrow \reels^{j}$ ($j \le d$) is Gaussian and is $C^2(D, \reels^j)$ on $D$, such that for all  $(x_1, x_2) \in D \times D - \Delta$, the vector $(X(x_1), X(x_2))$ has a density.
\\Moreover, for almost all $(x_1, x_2) \in D \times D$, the vector $(W(x_1), W(x_2))$ is independent of the vector $(X(x_1), X(x_2), \nabla\!X(x_1), \nabla\!X(x_2))$, and $\forall n \in \naturels$, 
$$
	\int_{D}\bbE\!\left[\normp[k]{W(x)}^{n}\right] \ud x< \infty.
$$%
\item $\bE_2$: $\forall x \in D$, $X(x)=F(Z(x))$, where $F: \reels^j \longrightarrow \reels^j$ is a bijective function in $C^2(\reels^j, \reels^j)$, such that $\forall z \in \reels^j$, $J_{F}(z)$, the Jacobian of $F$ in $z$, is such that $J_{F}(z)\neq 0$ and the function $F^{-1}$ is continuous.
The process
$Z: \Omega \times D \subset \Omega\times \reels^{d} \rightarrow \reels^{j}$ ($j \le d$)  is Gaussian and is $C^2(D, \reels^j)$, such that for all  $(x_1, x_2) \in D \times D - \Delta$, the vector $(Z(x_1), Z(x_2))$ has a density.
Moreover, for almost all $(x_1, x_2) \in D \times D$, the vector $(W(x_1), W(x_2))$ is independent of the vector 
$(Z(x_1), Z(x_2), \nabla\!Z(x_1), \nabla\!Z(x_2))$, and $\forall n \in \naturels$, 
$$
	\int_{D}\bbE\!\left[\normp[k]{W(x)}^{n}\right] \ud x< \infty.
$$%
\item $\bE_{3}$: $\forall x \in D$, $X(x)=F(Z(x))$, where the process $Z: \Omega \times D \subset \Omega\times \reels^{d} \rightarrow \reels^{j^{\prime}}$ ($j<j^\prime$) is Gaussian and is $C^2(D, \reels^{j^{\prime}})$, such that $\forall (x_1, x_2) \in D \times D - \Delta$, the vector $(Z(x_1), Z(x_2), \nabla\!Z(x_1), \nabla\!Z(x_2))$ has a density.
Moreover, for almost all $(x_1, x_2) \in D \times D$, the  vector $(W(x_1), W(x_2))$ is independent of the vector\
$(Z(x_1), Z(x_2), \nabla\!Z(x_1), \nabla\!Z(x_2))$.\\
Also, $\forall n \in \naturels$, 
$$
	\int_{D}\bbE\!\left[\normp[k]{W(x)}^{n}\right]  \ud x< \infty.
$$%
The function $F$ verifies assumption ${(\bF)}$ appearing in assumption $\bA_3$.
\item $\bE_{4}$: For almost all $(x_1, x_2, y_1, y_2, \dot{x}_1, \dot{x}_2) \in D \times D   \times \reels^{2} \times \reels^{dj} \times \reels^{dj}$ and $\forall u \in \reels^{j}$, the density
$$
	{p}_{Y(x_1), Y(x_2), X(x_1), X(x_2), \nabla\!X(x_1), \nabla\!X(x_2)}(y_1, y_2, u, u, \dot{x}_1, \dot{x}_2),
$$%
of the  joint distribution of $(Y(x_1), Y(x_2), X(x_1), X(x_2), \nabla\!X(x_1), \nabla\!X(x_2))$, exists and is continuous in the variable $u$.\\
Furthermore
\begin{align*}
&\specialpos{u \longmapsto \int_{D \times D} \int_{\reels^2 \times \reels^{2dj}} \abs{y_1} \abs{y_2}  \normp[dj]{\dot{x}_1}^{j} \normp[dj]{\dot{x}_2}^{j}\hfill}\\
& \specialpos{\hfill\quad\times{p}_{Y(x_1), Y(x_2), X(x_1), X(x_2), \nabla\!X(x_1), \nabla\!X(x_2)}(y_1, y_2, u, u, \dot{x}_1, \dot{x}_2)\hfill}\\
&  \specialpos{\hfill\times \ud \dot{x}_1 \ud \dot{x}_2 \ud y_1 \ud y_2 \ud x_1 \ud x_2,}
\end{align*}
is continuous.
\end{itemize}
Let us state the hypothesis ${\bH _7}$.
\begin{itemize}
\item${\bH _7}$: $\forall y \in \reels^{j}$,
\begin{multline*}
 \int_{D \times D} 
 	\bbE\!\left[\abs{Y(x_1)} \abs{Y(x_2)}  H(\nabla\!X(x_1)) H(\nabla\!X(x_2))\given{X(x_1)=X(x_2)=y}\right]\\
	  \times  {p}_{X(x_1), X(x_2)}(y, y) \ud x_1 \ud x_2 < \infty.
\end{multline*}
\end{itemize}

We are ready to prove Theorem~\ref{Cabana4 Bis}.
\spacebefore\begin{theo}
\label{Cabana4 Bis}
Let $X: \Omega \times D \subset \Omega\times \reels^{d} \rightarrow \reels^{j}$ ($j < d$) be a random field in $C^{2}(D,\reels^{j})$, where $D$ is a bounded, convex open set  of $\reels^{d}$,  such that for almost all $\omega \in \Omega$, $\nabla\!X(\omega)$ is Lipschitz with Lipschitz constant $L_{X}(\omega)$ such that  $\bbE\left[L_X^{2d}(\bm{\cdot})\right] < \infty$.
 Let $Y: \Omega \times D \subset \Omega \times \reels^{d} \rightarrow \reels$ a continuous process.\\
If $Y$ satisfies (\ref{egalite Y}) and if $X$ and $Y$ satisfy one of the three assumptions $\bD_{i} $, $i=1, 2, 3$ and the hypothesis ${\bH _7}$ or if $X$ and $Y$ satisfy the assumption $\bD_4$, then $\forall y \in \reels^{j}$ we have
\begin{align}
\label{momentdiag}
\MoveEqLeft[1]{\bbE\left[\left(\int_{\calC_{X}(y)}Y(x) \ud\sigma_{d-j}(x)\right)^2\right]} \nonumber\\ 
&\specialpos{\hfill= \int_{D \times D}\bbE\!\left[Y(x_1) Y(x_2) H(\nabla\!X(x_1)) H(\nabla\!X(x_2))
	\given{X(x_1)=X(x_2)=y}\right]}  \\ 
&\specialpos{\hfill \times {p}_{X(x_1), X(x_2)}(y, y) \ud x_1 \ud x_2.} \nonumber
\end{align}
\end{theo}
\spacebefore
\begin{rema}
\label{LipschitzBis}
As in Remark \ref{Lipschitz}, we can replace in the theorem the hypothesis that for almost all $\omega \in \Omega$, $\nabla\!X(\omega)$ is Lipschitz with  Lipschitz constant $L_{X}(\omega) \such \bbE\left[L_X^{2d}(\bm{\cdot})\right] < \infty$ by the hypothesis that 
$\bbE\!\left[\left(\sup_{x \in D} \normp[j,d]{\nabla^2X(x)}^{(s)}\right)^{2d}\right] < \infty.$
\hfill$\bullet$
\end{rema}
\spacebefore
\begin{rema}
\label{moment factoriel}
Under the same assumptions as in the theorem or those of Remark \ref{moment etoile} formulated later, for $j=d$, one obtains a result similar to that given in (\ref{momentdiag}).
We just have to replace 
$$
	\bbE\left[\left(\int_{\calC_{X}(y)}Y(x) \ud\sigma_{d-j}(x)\right)^2\right]
$$%
by 
$$
	\bbE\left[\!\left(\int_{\calC_{X}(y)}Y(x) \ud\sigma_{d-j}(x)\right)^2 - \int_{\calC_{X}(y)}Y^2(x) \ud\sigma_{d-j}(x)\right]
$$%
in (\ref{momentdiag}).
The right-hand side remains unchanged.
However, it should be noted that in this particular case, $\sigma_0$ is the counting measure.
\hfill$\bullet$
\end{rema}
\spacebefore
\begin{rema}
\label{similaire}
Under the same type of hypotheses as those given in the theorem, or later in Remark \ref{moment etoile}, one can propose a general Rice formula for the moments of order $k$ of the process $Y$ integrated on the level set of the random field  $X$,  and this $\forall y \in \reels^j$.
\hfill$\bullet$
\end{rema}
\spacebefore
\begin{rema}
Theorem \ref{Cabana4 Bis} and also Remarks \ref{moment factoriel}, \ref{similaire} and \ref{moment etoile} can be generalized to $D$, a convex open set  $\reels^d$ not necessarily bounded.
Remark \ref{moduleconvexe} arguments can be followed \textit{mutatis mutandis}. 
\hfill$\bullet$
\end{rema}
\spacebefore
\begin{proofarg}{Proof of Theorem \ref{Cabana4 Bis}}
The idea is to apply Remark \ref{ouvertconvexe} according to Theorem \ref{Hyp B} for the convex and open set $D \times D$ and the bounded open set $D_1= D \times D - \Delta$, to the processes $\widetilde{X}$ and $\widetilde{Y}$ as
$$
\begin{aligned}
\widetilde{X}:   \Omega \times D \times D \subset \Omega\times \reels^{2d}  & \longrightarrow \reels^{2j} \\
x = (x_1, x_2) &\longmapsto \widetilde{X}(x):= (X(x_1), X(x_2)),
\end{aligned}
$$%
and
$$
\begin{aligned}
\widetilde{Y}:   \Omega \times D \times D \subset \Omega\times \reels^{2d}  & \longrightarrow \reels \\
x = (x_1, x_2) &\longmapsto \widetilde{Y}(x):= Y(x_1) \times Y(x_2).
\end{aligned}
$$%
Since  $X$ is a random field in $C^{2}(D, \reels^j)$, then $\widetilde{X}$ is a random field in $C^{2}(D \times D, \reels^{2j})$ and $D \times D$ is a convex open set of $\reels^{2d}$.
Also $\widetilde{Y}_{|D \times D - \Delta}$ is still continuous on $D \times D - \Delta$, an open bounded set of  $\reels^{2d}$ contained in $D \times D$.\\
Since for almost all $\omega \in \Omega$, $\nabla\!X(\omega)$ is Lipschitz with  Lipschitz  constant $L_{X}(\omega) \mbox{$\such$} \bbE\!\left[L_X^{2d}(\bm{\cdot})\right] < \infty$, then for almost all $\omega \in \Omega$, $\nabla \widetilde{X}(\omega)$ is Lipschitz with Lipschitz constant $L_{\widetilde{X}}(\omega)= L_{X}(\omega)$, such that  $\bbE\!\left[L_{\widetilde{X}}^{2d}(\bm{\cdot})\right] < \infty$.\\
Then under one of the assumptions $\textbf{C}_i$, $i=1,\dots,3$, since $Y$ is written as a function $G$ of $X$ and of $\nabla\!X$ and of the variable $W: \Omega \times D \subset \reels^{d} \rightarrow \reels^k$, $k \in \naturels^{\star}$, in the following form, for almost all  $x \in D$:
$$
Y(x)= G(x, W(x), X(x), \nabla\!X(x)),
$$%
where
$$
\begin{aligned}
	G: D \times \reels^k \times \reels^j \times \frakL(\reels^{d}, \reels^{j}) & \longrightarrow \reels \\
	(x, z, u, A)  & \longmapsto G(x, z, u, A),
\end{aligned}
$$%
is a continuous function of its variables on 
$D \times \reels^k \times \reels^j \times  \frakL(\reels^{d}, \reels^{j})$
and such that $\forall (x, z, u, A) \in D \times \reels^k \times \reels^j \times \frakL(\reels^{d}, \reels^{j})$,
$$
	\abs{G(x, z, u, A)} \le P(f(x), \normp[k]{z}, h(u), \normp[j,d]{A}),
$$%
where $P$ is a polynomial with positive coefficients and $f: D \longrightarrow \reels^{+}$ and $h: \reels^j \longrightarrow \reels^{+}$ are continuous functions,  it is the same for $\widetilde{Y}$.
More precisely, for almost all $x:=(x_1, x_2) \in D \times D$, 
$$
\widetilde{Y}(x):= \widetilde{G}(x, \widetilde{W}(x), \widetilde{X}(x), \nabla \widetilde{X}(x)),
$$%
where
$$
\begin{aligned}
	\widetilde{W}:  \Omega \times D \times D \subset \Omega \times \reels^{2d} & \longrightarrow \reels^{2k} \\
	x=(x_1, x_2)  & \longmapsto \widetilde{W}(x):=(W(x_1), W(x_2)),
\end{aligned}
$$%
$$
\begin{aligned}
\specialpos{\widetilde{G}:  D^2 \times \reels^{2k} \times \reels^{2j} \times \frakB(\reels^{2d}, \reels^{2j}) \longrightarrow \reels\hfill} \\
\specialpos{\!\left(x=(x_1, x_2), z=(z_1, z_2), u=(u_1, u_2),
\begin{pmatrix}
A & 0 \\
0 & B
\end{pmatrix}\right)\hfill}\\
\specialpos{\hfill\longmapsto \widetilde{G}\!\left(x, z, u,
\begin{pmatrix}
	A & 0\\
	0 & B
\end{pmatrix}\right):= 
 G(x_1, z_1, u_1, A) \times G(x_2, z_2, u_2, B),}
\end{aligned}
$$%
where $\frakB(\reels^{2d}, \reels^{2j})$ is the vector subspace of $\frakL(\reels^{2d}, \reels^{2j})$ of the matrices of the form $C:=\begin{pmatrix}
A & 0 \\
0 & B
\end{pmatrix} $ where $A, B \in \frakL(\reels^{d}, \reels^{j})$.\\
It is clear that $\widetilde{G}$ remains a continuous function defined on 
$$
	D^2 \times \reels^{2k} \times \reels^{2j} \times  \frakB(\reels^{2d}, \reels^{2j}),
$$%
 and such that  $\forall (x, z, u, C) \in D^2 \times \reels^{2k} \times \reels^{2j} \times \frakB(\reels^{d}, \reels^{j})$,
 \begin{align*}
	\abs{\widetilde{G}(x, z, u, C)} &\le P(f(x_1), \normp[k]{z_1}, h(u_1), \normp[j,d]{A})\\
&\qquad\qquad \times P(f(x_2), \normp[k]{z_2}, h(u_2), \normp[j,d]{B}) \\
&\le  \widetilde{P}(\widetilde{f}(x), \normp[2k]{z}, \widetilde{h}(u), \normp[2j,2d]{C}), 
 \end{align*}
 where $\widetilde{f}$ is 
$$
\begin{aligned}
\widetilde{f}:  D^2 &\longrightarrow \reels^{+}\\
x=(x_1, x_2)  & \longmapsto \widetilde{f}(x):= f(x_1) + f(x_2),
\end{aligned}
$$%
while $\widetilde{h}$ is
$$
\begin{aligned}
\widetilde{h}:  \reels^{2j} &\longrightarrow \reels^{+}\\
u=(u_1, u_2)  & \longmapsto \widetilde{h}(u):= h(u_1) + h(u_2).
\end{aligned}
$$%
They are continuous functions and  $\widetilde{P}$ is a polynomial with positive coefficients.\\

It is easy to verify that $\widetilde{X}$ and $\widetilde{Y}_{|D \times D - \Delta}$ satisfy the hypotheses ${\bB_i}$, $i=1,\dots,4$, of Remark \ref{ouvertconvexe} following Theorem \ref{Hyp B}, respectively  for the convex open set $D \times D$ and for the open and bounded set $D \times D - \Delta$ contained  in $D \times D$.\\
Moreover,  if $X$ and $Y$ satisfy the hypothesis ${\bH _7}$ then $\widetilde{X}$ et $\widetilde{Y}_{|D \times D - \Delta}$ satisfy the hypothesis ${\bH _6}$, since $\forall x \in D \times D$, 
\begin{equation}
	H(\nabla \widetilde{X}(x))= H(\nabla\!X(x_1)) \times H(\nabla\!X(x_2)).
	\label{etoile}
\end{equation}
It turns out that $\widetilde{X}$ and $\widetilde{Y}_{|D \times D - \Delta}$ satisfy the hypotheses of Remark \ref{ouvertconvexe} following Theorem \ref{Hyp B}.\\
Now the hypotheses satisfied by $X$ and $Y$ make these two processes verify the hypotheses $\textbf{C}_i$, $i=1,\dots,4$, contained in Theorem \ref{Cabana4} and those of Proposition \ref{wsch}.
Therefore, in the same way as in this theorem, we obtain that $\forall y \in \reels^j$, 
\begin{equation}
\label{etoile1}
\bbP\{\omega \in \Omega: \exists x \in D, X(x)(\omega)=y ,\rank  \nabla\!X(x)(\omega)<j\}=0.
\end{equation}
One can deduce
$\forall y \in  \reels^j$, 
\begin{equation}
\label{etoile2}
\bbP\{\omega \in \Omega: \exists x \in D \times D, \widetilde{X}(x)(\omega)=(y, y),\rank  \nabla \widetilde{X}(x)(\omega)<2j\}=0.
\end{equation}
Indeed, using (\ref{etoile}), $\forall y \in \reels^j$,
\begin{multline*}
	\bbP\{\omega \in \Omega: \exists x \in D \times D, \widetilde{X}(x)(\omega)=(y, y) ,\rank  \nabla \widetilde{X}(x)(\omega)<2j\}\\
\le \bbP\{\omega \in \Omega: \exists x \in D, X(x)(\omega)=y ,\rank  \nabla\!X(x)(\omega)<j\}=0.
\end{multline*}
Remark \ref{ouvertconvexe} applied to $\widetilde{X}$ and
$\widetilde{Y}_{|D\times D - \Delta}$
with (\ref{etoile2}) and (\ref{etoile}) allow to write $\forall y \in\reels^{j}$,
\begin{align*}
 \MoveEqLeft[1]{\bbE\left[\int_{\calC^{D^{r}}_{D \times D - \Delta, \widetilde{X}}(y, y)} \widetilde{Y}(x) \ud\sigma_{2(d-j)}(x)\right]}\\
& =\bbE\left[\int_{\calC_{D \times D - \Delta, \widetilde{X}}(y, y)} \widetilde{Y}(x) \ud\sigma_{2(d-j)}(x)\right]\\
& = \int_{D \times D - \Delta}{p}_{\widetilde{X}(x)}(y, y) 
\bbE\left[\widetilde{Y}(x)H(\nabla
\widetilde{X}(x))\given{\widetilde{X}(x)= (y, y)}\right] \ud x\\
&=\int_{D \times D} \bbE[Y(x_1) Y(x_2) H(\nabla\!X(x_1)) H(\nabla\!X(x_2))\given{X(x_1)=X(x_2)=y}]\\
&\specialpos{\hfill\times{p}_{X(x_1), X(x_2)}(y, y) \ud x_1 \ud x_2.}
\end{align*}
The last equality  is justified by the fact that  $\sigma_{2d}(\Delta) =0$.\\
Moreover, we know, using Remark \ref{fonctions implicites} and (\ref{etoile1}) that $\forall y \in\reels^{j}$, almost surely $\calC_{X}^{D^{r}}(y)= \calC_{X}(y)$ and $\calC_{X}^{D^{r}}(y)$ is a differentiable manifold of dimension $(d-j)$.
Thus $\forall y \in\reels^{j},$ almost surely the set $A=\{(x, x) \in D \times D,   X(x) = y \}$
is a differentiable manifold of dimension $(d-j)$.
Thus, since $j < d$,  almost surely $\sigma_{2(d-j)}(A)=0$.
So, $\forall y \in\reels^{j}$,
\begin{align*}
{\bbE\left[\int_{\calC_{D \times D - \Delta, \widetilde{X}}(y, y)} \widetilde{Y}(x) \ud\sigma_{2(d-j)}(x)\right]}
&=\bbE\left[\int_{\calC_{D \times D, \widetilde{X}}(y, y)} \widetilde{Y}(x) \ud\sigma_{2(d-j)}(x)\right]\\
&=\bbE\left[\left(\int_{\calC_{X}(y)} Y(x) \ud\sigma_{d-j}(x)\right)^2\right].
\end{align*}
The last equality comes from the fact that $\forall y \in \reels^{j}$,
$$
	\calC_{D \times D, \widetilde{X}}(y, y) = \calC_{X}(y) \times \calC_{X}(y).
$$
This completes the proof of this theorem.
 \end{proofarg}
\spacebefore
\begin{proofarg}{Proof of Remark \ref{moment factoriel}}
Under the same hypotheses as in Theorem \ref{Cabana4 Bis}, but for $j=d$, we do the same proof as before.
We obtain $\forall y \in\reels^{j}$,
\begin{multline*}
\bbE\left[\int_{\calC_{D \times D - \Delta, \widetilde{X}}(y, y)} \widetilde{Y}(x) \ud\sigma_{2(d-j)}(x)\right]\\
=\int_{D \times D} 
\bbE\!\left[Y(x_1) Y(x_2)
 H(\nabla
X(x_1)) H(\nabla
X(x_2))\given{X(x_1)=X(x_2)=y}\right] \nonumber \\ 
   \times{p}_{X(x_1), X(x_2)}(y, y) \ud x_1 \ud x_2.
\end{multline*}
In the same way, the set  $A$ is still almost surely a differentiable manifold and since $\forall y \in\reels^{j}$
$$
	\calC_{D \times D, \widetilde{X}}(y, y) = \calC_{X}(y) \times \calC_{X}(y),
$$%
we can write, recalling that in this case $\sigma_{d-j}$ is the counting measure,
\begin{align*}
\MoveEqLeft[1]{\bbE\left[\int_{\calC_{D \times D - \Delta, \widetilde{X}}(y, y)} \widetilde{Y}(x) \ud\sigma_{2(d-j)}(x)\right]}\\
&=\bbE\left[\int_{\calC_{D \times D, \widetilde{X}}(y, y)} \widetilde{Y}(x) \ud\sigma_{2(d-j)}(x)
- \int_{\calC_{X}(y)}Y^2(x) \ud\sigma_{d-j}(x)\right]\\
&=\bbE\left[\!\left(\int_{\calC_{X}(y)}Y(x) \ud\sigma_{d-j}(x)\right)^2 - \int_{\calC_{X}(y)}Y^2(x) \ud\sigma_{d-j}(x)\right].
\end{align*}
This completes the proof of this remark.
\end{proofarg}
\spacebefore
\begin{rema}
\label{moment etoile}
In the same manner as in Theorems  \ref{Cabter} and \ref{Ricegenerale}, we can weaken the hypothesis  $\bbE\left[L_X^{2d}(\bm{\cdot})\right] < \infty$ in Theorem \ref{Cabana4 Bis}.
More precisely, we can make the hypothesis that for almost all $\omega \in \Omega$, $\nabla\!X(\omega)$ is Lipschitz or as in Remark \ref{LipschitzBis}, demander la quasi finitude de $\sup_{x \in D} \normp[j,d]{\nabla^2X(x)}^{(s)}$.
It will suffice to replace in Theorem \ref{Cabana4 Bis} the hypotheses $\bD_{i} $ by the following $\bD_{i}^\star$ assumptions, $i=1,\dots,4$, keeping the hypothesis ${\bH _7}$ and adding the following hypothesis ${\bH _7^\star}$.
\begin{itemize}
\item $\bD_1^\star$: It is the assumption $\bD_1$, plus the following hypothesis:
for almost all $(x_1, x_2, x_{3}, x_{4})  \in D^4$, the density of the vector 
$(X(x_1), X(x_2)$, $X(x_{3}), X(x_{4}))$ exists.

\item $\bD_2^\star$: It is the assumption $\bD_2$, plus the following hypothesis:
for almost all $(x_1, x_2, x_{3}, x_{4})  \in D^4$, the density of the vector 
$(Z(x_1), Z(x_2)$, $Z(x_{3}), Z(x_{4}))$ exists.

\item $\bD_{3}^\star$: It is the assumption $\bD_{3}$, plus the following  hypothesis:
for almost all $(x_1, x_2, x_{3}, x_{4})  \in D^4$, the density of the vector 
$(Z(x_1), Z(x_2)$, $Z(x_{3}), Z(x_{4}))$ exists.

\item $\bD_{4}^\star$: It is the assumption $\bD_{4}$, plus the following hypothesis:
the function
\begin{multline*}
(u_1, u_2) \longmapsto \int_{D \times D} \int_{\reels^2 \times \reels^{2dj}} y_1^2 y_2^2  \normp[dj]{\dot{x}_1}^{j} \normp[dj]{\dot{x}_2}^{j}\\\times{p}_{Y(x_1), Y(x_2), X(x_1), X(x_2), \nabla\!X(x_1), \nabla\!X(x_2)}(y_1, y_2, u_1, u_2, \dot{x}_1, \dot{x}_2) \\
\ud \dot{x}_1 \ud \dot{x}_2 \ud y_1 \ud y_2 \ud x_1 \ud x_2
\end{multline*}
is a continuous function.\\
For almost all $(x_1, x_2, x_{3}, x_{4}, \dot{x}_1, \dot{x}_2, \dot{x}_{3}, \dot{x}_{4}) \in D^4 \times \reels^{4dj} $ and for all $q:=(u_1, u_2, v_1, v_2) \in \reels^{4j}$, the density
$$
{p}_{X(x_1), X(x_2), X(x_{3}), X(x_{4}),\nabla\!X(x_1), \nabla\!X(x_2), \nabla\!X(x_{3}), \nabla\!X(x_{4})}(q, \dot{x}_1, \dot{x}_2, \dot{x}_{3}, \dot{x}_{4}),
$$%
of the vector 
$$
	\!\left(X(x_1), X(x_2), X(x_{3}), X(x_{4}), \nabla\!X(x_1), \nabla\!X(x_2), \nabla\!X(x_{3}), \nabla\!X(x_{4})\right)
$$%
exists.
Moreover,  $\forall y \in \reels^j$, the function 
\begin{multline*}
q \longmapsto
 \int_{D ^4} \int_{\reels^{4dj}} \normp[dj]{\dot{x}_1}^{j}  \normp[dj]{\dot{x}_2}^{j}  \normp[dj]{\dot{x}_{3}}^{j}
 \normp[dj]{\dot{x}_{4}}^{j} \\
 \times{p}_{X(x_1), X(x_2), X(x_{3}), X(x_{4}),\nabla\!X(x_1), \nabla\!X(x_2), \nabla\!X(x_{3}), \nabla\!X(x_{4})}(q, \dot{x}_1, \dot{x}_2, \dot{x}_{3}, \dot{x}_{4})  \\
\ud \dot{x}_1 \ud \dot{x}_2 \ud \dot{x}_{3} \ud \dot{x}_{4} \ud x_1 \ud x_2 \ud x_{3} \ud x_{4}
 \end{multline*}
  is bounded in a neighborhood of $q:=(y, y, y, y)$.
\end{itemize}
\vbox{Let us state the hypothesis  ${\bH _7^\star}$.
\begin{itemize}
\item ${\bH _7^\star}$:  $\forall y \in \reels^j$, the function
\end{itemize}
\begin{multline*}
 q \longmapsto\int_{D^4} {p}_{X(x_1), X(x_2), X(x_{3}), X(x_{4})}(q)
 	\bbE\left[H(\nabla\!X(x_1)) H(\nabla\!X(x_2)) \right.\\\left.
	 \times  H(\nabla\!X(x_{3}))H(\nabla\!X(x_{4}))\given{(X(x_1), X(x_2), X(x_{3}), X(x_{4}))= q}\right]\\
	 \ud x_1 \ud x_2 \ud x_{3} \ud x_{4},
\end{multline*}
is a bounded function in a neighborhood of $q:=(y, y, y, y)$.}
\hfill$\bullet$
\end{rema}
%
%
\chapter{Applications}
The main reason for having well-fitting Kac-Rice formulas is that they provide tools for explicitly doing calculations that involve roots of functions as well as functionals of other levels.
Below we will present some of these applications.

Let us first mention the possibility for getting conditions in which the level functional has some moments.
This is a non-trivial task that has only been completely solved in some special cases.
Furthermore, Rice's formulas have been also applied in physical oceanography and in the theory of dislocations of random waves propagation.
Two other applications deserve to be studied: first the theory of random gravitational microlensings and second the study of the zero sets of random algebraic systems invariant under the orthogonal group, known in the literature as Kotlan-Shub-Smale systems.
In the following, the reader will find a brief description of each of them.
\section{Dimensions $\textit{d}=\textit{j}=1$}
Classically, the study of the Rice formula began with the seminal papers by Kac \cite{MR7812} and Rice \cite{MR10932}.
The first considered the number of roots of a random polynomial with standard and independent Gaussian coefficients and the second developed formulas to study the crossovers of stationary Gaussian processes.
In this subsection, we will revisit these two old problems.
First, we will give, using the formulas obtained earlier, the necessary and sufficient conditions for the existence of the first and second moment of the number of crossings of a stationary Gaussian process.
Second, Kac's research will be extended to consider random trigonometric polynomials.
\subsection{Necessary and sufficient conditions for the first two moments of the number of crossings}
\label{mehler section}

Let $X:\Omega\times\reels \to\reels $ be a real and stationary Gaussian process with zero mean.
Let us denote its covariance function by $r$ and its spectral measure by $\mu$ which is assumed not to be purely discrete.
We have
$$
	r(t)=\int_{\reels }e^{it\lambda}\ud \mu(\lambda).
$$%
The spectral moment of order $p$ is defined as follows
$$
\lambda_p:=\int_{\reels }\lambda^p \ud \mu(\lambda).
$$%
For $y \in \reels$ and $t > 0$, let $N^X_{[0,t]}(y)$ the number of crossings of the level $y$ by the process $X$ on the interval $\intrff{0}{t}$.
We have the following theorem.
\spacebefore \begin{theo}
\label{onedimension}\begin{itemize}
\item The Rice formula of the first order holds if and only if $\lambda_2< \infty$,
and $\forall y \in \reels$ and $t > 0$ we have
$$
\bbE\!\left[N^X_{[0,t]}(y)\right]=\frac t{\pi}\sqrt{\frac{\lambda_2}{\lambda_0}}e^{-\frac{y^2}{2\lambda_0}}.
$$%
\item Moreover, in this case, $\bbE\!\left[(N^X_{[0,t]}(y))^2\right]< \infty$ if and only if for some $\delta>0$ we have
$$
	\frac{r^{\prime\prime}(\tau)-r^{\prime\prime}(0)}\tau\in\bL^1([0,\delta],\ud \tau).
$$%
\end{itemize}
\end{theo}
\spacebefore
\begin{rema} The first result was proved by K$.$ It\^o in \cite{MR166824}.
In this  work, the author generalizes the previous proofs providing  a definitive result.
The second is the famous result of Geman \cite{MR301791}.
He considers only the case $y=0$.
In \cite{MR2257657} the result has been extended for all $y$.
\hfill$\bullet$
\end{rema}
\spacebefore \begin{proofarg}{Proof of Theorem \ref{onedimension}} Remark \ref{Lipschitz} following Theorem \ref{Cabana4} gives the validity of the first formula whenever $X$ is $C^{2}([0, t],\reels)$.
However, the result is valid with great generality as shown by It\^o in \cite{MR166824}.
For the sake of completeness, we will outline his proof.
It is first proved in \cite{MR166824} that if $\lambda_2< \infty$ the process has absolutely continuous trajectories.
And besides, it is true that
$$
	N^X_{[0,t]}(y)\le \liminf_{\delta\to0}\frac1{2\delta}\int_{0}^t\1_{\{\abs{X(s)-y}<\delta\}}\abs{X^{\prime}(s)}\ud s.
$$%
Using Fatou's lemma and the fact that $X$ is Gaussian and stationary,  we get by denoting $\varphi$ for the standard Gaussian density on $\reels$
\begin{align*}
\bbE[N^X_{[0,t]}(y)]&\le \liminf_{\delta\to0}\frac t{2\delta}\bbE[\1_{\{\abs{X(0)- y}<\delta\}}\abs{X^{\prime}(0)}]\\
	&= \liminf_{\delta\to0}\frac t{2\delta \sqrt{\lambda_0\lambda_2}}\int_{y-\delta}^{y+\delta}\int_{\reels }\abs{\dot{z}}\varphi\left(\frac{z}{\sqrt{\lambda_0}}\right)\varphi\left(\frac{\dot{z}}{\sqrt{\lambda_2}}\right) \ud z \ud\dot{z}\\
	&=\liminf_{\delta\to0}\frac t{2\delta\sqrt{\lambda_0}}\int_{y-\delta}^{y+\delta}\varphi\left(\frac{z}{\sqrt{\lambda_0}}\right)\ud z \sqrt{\frac{2\lambda_2}\pi}\\
	&=\frac t\pi\sqrt{\frac{\lambda_2}{\lambda_0}}e^{-\frac{y^2}{2\lambda_0}}.
\end{align*}
Concerning the other inequality, in \cite{MR166824}, it is proved that the following inequality (monotonic limit) is true
$$
	N^X_{[0,t]}(y)\ge\lim_{n\to+\infty}\sum_{k=1}^{2^n}\1_{\left\{[X({(k-1)t}/{2^n})-y][X({kt}/{2^n})- y]<0\right\}}.
$$%
Therefore, using the monotone convergence theorem, we obtain
$$
	\bbE[N^X_{[0,t]}(y)]\ge \lim_{n\to+\infty}2^n\bbE[\1_{\{[X(0)-y][X({t}/{2^n})-y]<0\}}].
$$%
The expectation on the righthand side can be written as follows
\begin{multline*}
 \bbE\!\left[\1_{\{(X(0)-y)(X({t}/{2^n})-y)<0\}}\right]\\
 =\bbE\!\left[\1_{\!\left({y}/{\sqrt{\lambda_0}},+\infty\right)}\!\left({X(0)}/{\sqrt{\lambda_0}}\right)
 		\1_{\!\left(-\infty,{y}/{\sqrt{\lambda_0}}\right)}\!\left({X({t}/{2^n})}/{\sqrt{\lambda_0}}\right)\right]\\
	+\bbE\!\left[\1_{\!\left({y}/{\sqrt{\lambda_0}},+\infty\right)}\!\left({X({t}/{2^n})}/{\sqrt{\lambda_0}}\right)
		\1_{\!\left(-\infty,{y}/{\sqrt{\lambda_0}}\right)}({X(0)}/{\sqrt{\lambda_0}})\right].
\end{multline*}
For ease of notation, let $\lambda_0=1$. Thus if $Z_n$ stands for a standard Gaussian \rv
independent of $(X(0), X({t}/{2^n}))$ we have
\begin{align*}
\MoveEqLeft[1]{\bbE\!\left[\1_{\!\left(y,+\infty\right)}\!\left(X(0)\right)\1_{\!\left(-\infty,y\right)}\!\left(X({t}/{2^n})\right)\right]}\\
	&=\bbE\!\left[\1_{\!\left(y,+\infty\right)}\!\left(X(0)\right)\1_{\!\left(-\infty,y\right)}\!\left(r(t/{2^n})X(0)+\sqrt{1-r^2(t/{2^n})}Z_n\right)\right]\\
	&=\int_{-\infty}^{0}\varphi(z) \ud z\int_{y}^{{\!\left(y-\sqrt{1-r^2(t/{2^n})}z\right)}/{r(t/{2^n})}}\varphi(x)\ud x,\\
\end{align*}
thus
$$
	2^n\int_{-\infty}^{0}\varphi(z) \ud z\int_{y}^{{\!\left(y-\sqrt{1-r^2(t/{2^n})}z\right)}/{r(t/{2^n})}}\varphi(x)\ud x
		\cvg{} t \frac{\sqrt{\lambda_2}}{2\pi}e^{-\frac12 y^2}.
$$%
We can proceed in the same way for the second term.
Finally obtaining for all $\lambda_0$
$$
	\bbE\!\left[N^X_{[0,t]}(y)\right]\ge\frac t{\pi}\sqrt{\frac{\lambda_2}{\lambda_0}}e^{-\frac12{y^2}/{\lambda_0}}.
$$%
The above procedure can also be used to prove that if $\lambda_2=+\infty$ then $\bbE\!\left[N^X_{[0,t]}(y)\right]=+\infty$.
And all results are valid.

To prove the second statement of the theorem, one can use Remarks \ref{LipschitzBis} and \ref{moment factoriel} following Theorem \ref{Cabana4 Bis} for the case $d=j=1$, thus assuming that $X$ is $C^{2}([0, t],\reels)$.
Thus the formula for the second factorial moment holds provided that the integral appearing in the equation (\ref{second factorial moment}) is finite.
Note that the assumption that $X$ has $C^2$ trajectories implies that the covariance $r$ is $C^4$ and then $\lambda_4 < \infty$.
However the case $\lambda_4= +\infty$ remains an interesting case.
This is the reason why we will follow the more general way given by \cite{MR2108670}  and \cite{MR301791}.
In  \cite{MR2108670} it is shown that
\begin{align}
M_2(y,t)&:= \bbE\!\left[N^X_{[0,t]}(y)(N^X_{[0,t]}(y)-1)\right] \nonumber \\ 
	&\phantom{:}=2\int_0^t (t-\tau)\int_{\reels ^2}\abs{\dot{x}_1}\abs{\dot{x}_2}{p}_{\tau}(y,\dot{x}_1, y,\dot{x}_2) \ud \dot{x}_1 \ud \dot{x}_2 \ud \tau,
	\label{second factorial moment}
\end{align}
where ${p}_{\tau}(x_1,\dot{x}_1, x_2,\dot{x}_2)$ is the density of the vector
$$
	\!\left(X( 0), X^{\prime}(0),X({\tau}),X^{\prime}({\tau})\right),
$$%
which is non-singular for $\tau>0$, since the spectral measure $\mu$ is not purely discrete.
Moreover, we show that if one of the terms in \ref{second factorial moment} is infinite, so is the other one.
In this way, we give a necessary and sufficient condition for the right-hand side of the formula to be finite.
Without loss of generality we can assume that $r(0)=1$ and since $\lambda_2 < \infty$ that $r$ is twice differentiable.

We will start by showing the result for the $y=0$ level which is the original Geman result.
 Let us write $M_2(y,t)$ in another way using a regression model.
We have%
\begin{equation}
	M_2(y,t)=2\int_0^t(t-\tau) {p}_\tau(y,y)\bbE\!\left[\abs{X^{\prime}(0)X^{\prime}(\tau)} \given{X(0)=X(\tau)=y}\right] \ud \tau,
	\label{secondmoment}
\end{equation} 
where ${p}_\tau(x_1, x_2)$ stands for the density of the vector $(X(0), X(\tau))$.
The following model will be useful
\begin{align*}
	X^{\prime}(0)&=\xi+\alpha_1(\tau)X(0)+\alpha_2(\tau)X(\tau)\\
	X^{\prime}(\tau)&=\xi^\star+\beta_1(\tau)X(0)+\beta_2(\tau)X(\tau),
\end{align*}
where $(\xi,\xi^\star)$ is a Gaussian centered vector independent of $(X(0),X(\tau))$, and 
\begin{align*}
\Var(\xi)
	&\phantom{:}=\Var(\xi^\star):=\sigma^2(\tau)=-r^{\prime\prime}(0)-\frac{(r^{\prime}(\tau))^2}{1-r^2(\tau)} ,\\
\rho(\tau)
	&:=\frac{\Cov(\xi,\xi^\star)}{\sigma^2(\tau)}=\frac{-r^{\prime\prime}(\tau)(1-r^2(\tau))-(r^{\prime}(\tau))^2r(\tau)}{-r^{\prime\prime}(0)(1-r^2(\tau))-(r^{\prime}(\tau))^2}.
\end{align*}
Moreover 
$$
\begin{array}{ll}
 \ds\alpha_1(\tau)=\frac{r^{\prime}(\tau)r(\tau)}{1-r^2(\tau)};&\ds\alpha_2(\tau)=-\frac{r^{\prime}(\tau)}{1-r^2(\tau)}\\
\ds \beta_1(\tau)= -\alpha_2(\tau); &\ds\beta_2(\tau)=-\alpha_1(\tau).
\end{array}
$$%
In this form, we have
\begin{align*}
M_2(0,t)
	&=2\int_0^t(t-\tau){p}_\tau(0, 0)\bbE\!\left[\abs{\xi}\abs{\xi^\star}\right]\ud \tau\\
	&=\frac{1}{\pi}\int_0^t(t-\tau)\frac{\sigma^2(\tau)}{(1-r^2(\tau))^{1/2}}\bbE\!\left[\abs{\frac{\xi}{\sigma(\tau)}}\abs{\frac{\xi^\star}{\sigma(\tau)}}\right]\ud \tau.
\end{align*}
Using the Cauchy-Schwarz inequality, we get
$$
	M_2(0,t)\le \frac{t}{\pi}\int_0^t\frac{\sigma^2(\tau)}{(1-r^2(\tau))^{1/2}}\ud \tau.
$$%
So if the integral on the righthand side is finite then $M_2(0,t)< \infty$.
But the integral converges if for $\delta>0$ we have 
$$
	\int_0^\delta\frac{\sigma^2(\tau)}{(1-r^2(\tau))^{1/2}}\ud \tau< \infty,
$$%
because the integrand is continuous in $[\delta,t]$.
Denoting reciprocally by $a_{2k}$ the coefficients of the function $\abs{x}$ in the orthogonal Hermite basis $(H_k(x))_{k \in \naturels}$ of $\bL^2(\reels ,\varphi(x)\ud x)$, i.e.
$$
	H_k(x):= (-1)^k \varphi^{-1}(x)  \ds \frac{\ud ^k}{\ud  x^k}(\varphi(x) ),
$$%
Mehler's formula gives (see \cite{MR716933})
\begin{align*}
+\infty>M_2(0,t)
	&=\frac{2}{\pi}\int_0^t(t-\tau)\frac{\sigma^2(\tau)}{(1-r^2(\tau))^{1/2}}
		\!\left(\sum_{k=0}^\infty a_{2k}^2(2k)!\rho(\tau)^{2k}\right)\ud \tau\\
	&\ge\frac{2}{\pi}a_0^2\int_0^t(t-\tau)\frac{\sigma^2(\tau)}{(1-r^2(\tau))^{1/2}}\ud \tau\\
	&\ge\frac{2}{\pi}a_0^2(t-\delta)\int_0^\delta\frac{\sigma^2(\tau)}{(1-r^2(\tau))^{1/2}}\ud \tau.
\end{align*}
The proof is complete in this case if we can prove that
$$
\frac{r^{\prime\prime}(\tau)-r^{\prime\prime}(0)}\tau \in\bL^1([0,\delta],\ud \tau)\Longleftrightarrow \int_0^\delta\frac{\sigma^2(\tau)}{(1-r^2(\tau))^{1/2}}\ud \tau< \infty.
$$%
But this is the purpose of Lemma \ref{bounds} proved below.

We will now consider the case where $y$ is any real number.
Let us define 
$$
	m(\tau):=\displaystyle \frac{y}{1+r(\tau)}\frac{r^{\prime}(\tau)}{\sigma(\tau)},
$$%
and introduce the expression
$$
	A(m,\rho,\tau):=\bbE\!\left[\abs{\frac{\xi}{\sigma(\tau)}-m(\tau)}\abs{\frac{\xi^\star}{\sigma(\tau)}+m(\tau)}\right].
$$%
Using  (\ref{secondmoment}) and regression, it turns out that
$$
M_2(y,t)=2\int_0^t(t-\tau){p}_\tau(y,y)\sigma^2(\tau)A(m,\rho,\tau)\ud \tau.
$$%
By applying the Cauchy-Schwarz inequality, we obtain
\begin{align}
\label{M2}
	A(m,\rho,\tau)
	&\le \!\left({\bbE\!\left[\left(\frac{\xi}{\sigma(\tau)}-m(\tau)\right)^2\right]  \bbE\!\left[\left(\frac{\xi^\star}{\sigma(\tau)}+m(\tau)\right)^2\right]}\right)^{\frac{1}{2}}\nonumber\\
	&= 1+m^2(\tau).
\end{align}
Let us now prove that the function $m(\tau)$ is bounded in a neighborhood of $\tau=0$.
For this purpose, let us consider the asymptotic behavior of ${r'(\tau)}/{\sigma(\tau)}$.
Two cases must be considered depending on whether $\lambda_4$ is finite or not.
In the first case, a Taylor expansion  of order 4 of  ${(r^{\prime})^2(\tau)}/{\sigma^2(\tau)}$ easily gives that ${r^{\prime}(\tau)}/{\sigma(\tau)} \to {-2\lambda_2}/{\sqrt{\lambda_4-\lambda_2^2}}$.
Now suppose that  $\lambda_4=+\infty$.
Given that 
$$
r^{\prime\prime}(\tau)-r^{\prime\prime}(0)=2\int_0^\infty [1-\cos(\tau\lambda)] \lambda^2 \ud \mu(\lambda),
$$%
we have by Fatou's lemma
\begin{align*}
\liminf_{\tau\to0}\frac{r^{\prime\prime}(\tau)-r^{\prime\prime}(0))}{\tau^2}
	&\ge \int_0^{+\infty}\liminf_{\tau\to0}\frac{1-\cos(\tau\lambda)}{{(\tau\lambda)^2}/{2}}\lambda^4 \ud \mu(\lambda)\\
	&=\int_0^\infty\lambda^4 \ud \mu(\lambda)=+\infty.
\end{align*}
Moreover
$$
\frac{(r^{\prime})^2(\tau)}{\sigma^2(\tau)}\simeq \frac{\lambda_2^3}{\ds\frac{\lambda_2(1-r^2(\tau))-(r^{\prime})^2(\tau)}{\tau^4}} ,
$$
and
$$
\lambda_2(1-r^2(\tau))-r'^2(\tau)= 2 \lambda_2(1-r(\tau))-(r^{\prime})^2(\tau)+O(\tau^4).
$$%
Moreover, using L'Hôpital's rule, we obtain
$$
\lim_{\tau\to0}\frac{2 \lambda_2(1-r(\tau))-(r^{\prime})^2(\tau)}{\tau^4}= \lim_{\tau\to0} \!\left({\frac{-r^{\prime}(\tau)}{2\tau}}\right)\!\left({\frac{r^{\prime\prime}(\tau)-r^{\prime\prime}(0)}{\tau^2}}\right)=+\infty,
$$%
since we know that ${-r^{\prime}(\tau)}/{(2\tau)}\to {\lambda_2}/{2}$.
Thus ${r'(\tau)}/{\sigma(\tau)}\to 0$.\\
Theses calculations lead us to conclude that
$$
m(\tau) \cvg[\tau\to 0]{} 
	\begin{cases}
		 \ds\frac{-\lambda_2 y}{\sqrt{\lambda_4-\lambda_2^2}} & \mbox{if $\lambda_4 < \infty$;}\\
		 0 & \mbox{otherwise.}
	\end{cases}
$$%
In both cases, we have then shown that the function $m(\tau)$ is bounded.\\
In this form, by using (\ref{M2}), we easily obtain
$$
	M_2(y,t)\le  \bC t\int_0^t {p}_\tau(y,y)\sigma^2(\tau)\ud \tau,
$$%
and by Lemma \ref{bounds} this integral is finite under Geman condition.

To prove the other implication assume that $M_2(y,t)< \infty$.
Thus
$$
	M_2(y,t)\ge 2\int_0^\delta(t-\tau){p}_\tau(y, y)\sigma^2(\tau)A(m,\rho,\tau)\ud \tau.
$$%
We will study $A(m, \rho, \tau)$.

Since the function $m(\tau)$ is bounded the following expansion is valid
$$
\abs{x-m(\tau)}=\sum_{k=0}^\infty a_k(m(\tau))H_k(x),
$$ 
where deleting the variable $\tau$ in $m$, the coefficients are
\begin{align*}
	a_0(m)&=m[2\Phi(m)-1]+2 \varphi(m)\\
	a_1(m)&=1-2\Phi(m)\\
	a_{\ell}(m)&=\frac2{\ell!}H_{\ell-2}(m) \varphi(m), \ell\ge 2,
\end{align*}
where $\Phi$ represents the Gaussian distribution of $\varphi$.\\
Using that function $a_k(m)$ is even if $k$ is even and odd otherwise, the Mehler's formula gives
\begin{align*}
	A(m,\rho,\tau)
		&=\sum_{k=0}^\infty a_k(m(\tau))a_k(-m(\tau))k!\rho^k(\tau)\\
		&=\sum_{k=0}^\infty a^2_{2k}(m(\tau))(2k)!(\rho(\tau))^{2k}-\sum_{k=0}^\infty a^2_{2k+1}(m(\tau))(2k+1)!(\rho(\tau))^{2k+1}.
\end{align*}
But by defining the odd projection as $M_{\text{odd}}(x,m):=\frac12(\abs{x-m}-\abs{x+m})$, we have
$$
	M_{\text{odd}}(x,m)=\sum_{k=0}^\infty a_{2k+1}(m)H_{2k+1}(x).
$$%
 Then
\begin{align*}
	\MoveEqLeft[1]{\abs{\bbE\!\left[M_{\text{odd}}\!\left(\frac{\xi}{\sigma(\tau)},m(\tau)\right)M_{\text{odd}}\!\left(\frac{\xi^\star}{\sigma(\tau)},m(\tau)\right)\right]}}\\
		&=\abs{\sum_{k=0}^\infty a^2_{2k+1}(m(\tau))(2k+1)!(\rho(\tau))^{2k+1}}\\
		&\le\bbE\!\left[M^2_{\text{odd}}\!\left(\frac{\xi}{\sigma(\tau)},m(\tau)\right)\right]\\
		&=\int_{\reels }\!\left(\tfrac12(\abs{x-m(\tau)}-\abs{ x+m(\tau)})\right)^2\varphi(x)\ud x\\
		&\le m^2(\tau).
\end{align*}
 Thus
$$
	A(m,\rho,\tau)\ge a^2_{0}(m(\tau))-m^2(\tau).
$$ 
Now it is easy to see that if $-m_0 \le m \le m_0$ then $a^2_{0}(m)-m^2 \ge \sqrt{{2}/{\pi}}(a_0(m_0)-m_0) >0$.\\
Since the function $m(\tau)$ is bounded in a neighborhood of zero, this implies that  $A(m,\rho,\tau) \ge C$ for $\tau$ sufficiently small.
Then
$$
+\infty>\displaystyle M_2(y,t)
\ge C\int_0^\delta(t-\tau)p_{\tau}(y,y)\sigma^2(\tau)\ud \tau\ge C \int_0^\delta \frac{\sigma^2(\tau)}{(1-r^2(\tau))^{1/2}}\ud \tau,
$$%
and we end up evoking again Lemma \ref{bounds}.

\end{proofarg}
\spacebefore
\begin{lemm}\label{bounds} There exists $\delta>0$ such that
$$
		\frac{r^{\prime\prime}(\tau)-r^{\prime\prime}(0)}\tau \in\bL^1([0,\delta],\ud \tau)
	\Longleftrightarrow
		\int_0^\delta\frac{\sigma^2(\tau)}{(1-r^2(\tau))^{1/2}} \ud \tau<\infty.
$$%
\end{lemm}
\spacebefore
\begin{proofarg} {Proof of Lemma \ref{bounds}} Let us consider the integral
$$
	\int_0^\delta \frac{\sigma^2(\tau)}{(1-r^2(\tau))^{1/2}}\ud \tau.
$$%
For $\tau$ small enough, we have
$$
	\frac{\sigma^2(\tau)}{(1-r^2(\tau))^{1/2}}\simeq \left(\frac{1}{\lambda_2^{{3}/{2}} }\right)\frac{-r^{\prime\prime}(0)(1-r^2(\tau))-(r^{\prime}(\tau))^2}{\tau^3},
$$%
thus integrating by part
\begin{align*}
\MoveEqLeft[1]{\int_0^\delta \frac{-r^{\prime\prime}(0)(1-r^2(\tau))-(r^{\prime}(\tau))^2}{\tau^3}\ud \tau}\\
	&= \left.\frac{r^{\prime\prime}(0)(1-r^2(\tau))+(r^{\prime}(\tau))^2}{2\tau^2}\right|_{0}^{\delta}
		+\int_0^\delta\frac{r^{\prime}(\tau)}{\tau}\left(\frac{r^{\prime\prime}(0)r(\tau)-r^{\prime\prime}(\tau)}{\tau}\right)\ud \tau\\
	&=\frac{r^{\prime\prime}(0)(1-r^2(\delta))+(r^{\prime}(\delta))^2}{2\delta^2}
	+ \int_{0}^{\delta} \frac{r^{\prime}(\tau)}{\tau} r^{\prime\prime}( 0) \!\left({\frac{r(\tau)-1}{\tau}}\right) \ud \tau\\
	&\specialpos{\hfill+ \int_{0}^{\delta} \frac{-r^{\prime}(\tau)}{\tau} \!\left({\frac{r^{\prime\prime}(\tau)-r^{\prime\prime}(0)}{\tau}}\right) \ud \tau.}
\end{align*}
Finally, since
$$
	\frac{r^{\prime}(\tau)}{\tau}  \!\left({\frac{r(\tau)-1}{\tau}}\right) \simeq \lambda_2^2 \frac{\tau}{2} \in L^1([0,\delta],\ud \tau),
$$%
and ${-r^{\prime}(\tau)}/{\tau} \to \lambda_2$, the above integral is finite if and only if
$$
	\int_0^\delta\frac{r^{\prime\prime}(\tau)-r^{\prime\prime}(0)}{\tau}\ud \tau< \infty.
$$%
\end{proofarg}
\subsection{Numbers of roots of random trigonometric polynomials}
In the following, we will study the asymptotic behavior of random Gaussian trigonometric polynomials.
 For any $N\in\naturels^{\star}$ and for two independent sequences of i.i.d$.$ standard Gaussian random variables $\{a_n\}_{n=1}^\infty$ and $\{b_n\}_{n=1}^\infty$
these functions are defined as
$$
	X_N(t):=\frac1{\sqrt N}\sum_{n=1}^N(a_n\sin n t+b_n\cos n  t).
$$%
The number of zeros of such a process has been extensively studied lately (see \cite{MR2797349} for example).
The process $X_N$ is an infinitely differentiable stationary Gaussian process of mean zero.
We can define as before $N^{X_N}_{\intrfo{0}{2\pi}}(y)$ as the number  $y$-level crossings of these trigonometric polynomials on the time interval $\intrfo{0}{2\pi}$.
The smoothness of these polynomials implies that the Rice's formula holds.
The necessary ingredients for its application are
$$
\bbE[X_N^2(0)]=1; \qquad \bbE[(X^{\prime}_N(0))^2]=\frac1N\sum_{n=1}^N n^2=\frac{(N+1)(2N+1)}6.
$$%
Hence 
$$
\bbE[N^{X_N}_{\intrfo{0}{2\pi}}(y)]=2\pi\sqrt{E[(X^{\prime}_N(0))^2]}\sqrt{\frac 2{\pi}}\frac{e^{-{y^2}/2}}{\sqrt{2\pi}}
=\frac2{\sqrt 3}\sqrt{\frac{{(N+1)(2N+1)}}2}e^{-\frac{y^2}2}.
$$%
Yielding
$$
\lim_{N\to\infty}\frac{\bbE[N^{X_N}_{\intrfo{0}{2\pi}}(y)]}N=\frac2{\sqrt 3}e^{-\frac{y^2}2}.
$$%

To calculate the variance and its asymptotic value we need to consider the rescaled process: $Y_N(t):=X_N(\frac{t}N).$ Since the covariance function of $X_N$ is
$$
r_{X_N}(t)=\frac1N\sum_{n=1}^N\cos n t
	=\frac1N
		\cos\!\left(\frac{(N+1) t}2\right)\frac{\sin(\frac12{N t})}{\sin\frac12t} ,
$$%
we get
$$r_{Y_N}(t)\to r_X( t):=\frac{\sin t}{t}.$$

Similar results can be obtained for the first and second derivative of $r_{Y_N}.$ The above result leads us to consider the sine cardinal process which has as covariance the function $r_X$.
In \cite{MR3084654} was proved that by constructing  the processes $Y_N$ and $X$ in the same probability space and if we define
$$
	B_N:=\bbE\!\left[
			\left\{
				\!\left(N^{X_N}_{\intrfo{0}{2\pi}}(0)-\bbE[N^{X_N}_{\intrfo{0}{2\pi}}(0)]\right)
				-\!\left(N^{X}_{\intrfo{0}{2\pi N}}(0)-\bbE[N^{X}_{\intrfo{0}{2\pi N}}(0)]\right)
			\right\}^2
		\right]
$$%
it turns out that
${B_N}/N \cvg[N \to \infty]{} 0.$
This result entails that  
$$
\lim_{N\to\infty}\frac1N\Var\!\left(N^{X_N}_{\intrfo{0}{2\pi}}(0)\right)=\lim_{N\to\infty}\frac1N\Var\!\left(N^{X}_{\intrfo{0}{2\pi N}}(0)\right),
$$
and the last quantity is
$$
\frac2{\sqrt3}+ 2 \int_0^\infty\!\left(\frac{\bbE\!\left[\abs{X^{\prime}(0)X^{\prime}(\tau)}\given{X(0)=X(\tau)=0}\right]}
	{\sqrt{1-({\sin\tau}/{\tau})^2}}-\frac{1}{3\pi}\right)\ud \tau.
$$%
\section{Sea modeling applications}

In this section, we give some theoretical justifications to the work of Podg\'orski \& Rychlik \cite{AL121221a}.
This paper presents several applications to random sea waves.
Like these authors, let us consider two random fields
$$
	X:\reels ^d\to\reels (\mbox{ with  }d>j=1) \mbox{ and } V:\reels ^d\to\reels ^{d+1}.
$$
The latter is defined as $V(x):=(X(x),\nabla\!X(x))$.
This is the argument of the function $G$ in (\ref{egalite Y}) but removing the explicit dependence on $x$ and also on the field $W$,
\ie $Y( x):=G(X(x),\nabla\!X(x))$.
 
Moreover, for sea applications either the $X$ field is Gaussian and models the sea surface or it is the envelope field (defined below).

We will first discuss the case where the stationary Gaussian field of zero mean $X$ is $C^{2}(D,\reels)$ with $\sigma^2:=\mbox{Var}[X(0)]=\bbE[X^2(0)]>0$.
Then Remark \ref{Lipschitz} applies and we  have for all $y \in \reels$
\begin{align*}
\bbE\left[\int_{\calC_{X}(y)}Y(x) \ud \sigma_{d-1}(x)\right]
	&= \int_{D} {p}_{X(x)}(y) 
		\bbE\left[Y(x)\normp[d]{\nabla\!X(x)}\given{X(x)=y}\right] \ud x\\
	& =\sigma_d(D)\bbE\!\left[Y(0)\normp[d]{\nabla\!X(0)}\given{X(0)=y}\right] \frac{e^{-\frac{1}{2\sigma^2}y^2}}{\sqrt{2\pi} \sigma}.
\end{align*}
The following notion is also introduced in \cite{AL121221a}.
We define the distribution of $V$  over the level set by first taking  $G(z):=\1_A(z)$ for $ z\in\reels ^{d+1}$ and $A$ a Borel set of $\reels^{d+1}$.
We have $Y(x)=\1_A(V(x))$ and setting 
\begin{eqnarray}
\label{distribution}
\nonumber
\bbP\{V(x)\in A\given{X(x)=y}\}&:=&\frac{\bbE\left[\int_{\calC_{X}(y)}\1_A(V(x)) \ud \sigma_{d-1}(x)\right]}{\bbE\left[\int_{\calC_{X}(y)}\ud\sigma_{d-1}(x)\right]}\\
&\phantom{:}=&\frac{\bbE\left[\1_A(V(0))\normp[d]{\nabla\!X(0)}\given{X(0)=y}\right]}{\bbE\left[\normp[d]{\nabla\!X(0)}\right]}\nonumber \\
&\phantom{:}=&\frac{\bbE\left[\1_A(y, \nabla\!X(0)) \normp[d]{\nabla\!X(0)}\right]}{\bbE\left[\normp[d]{\nabla\!X(0)}\right]}.
\end{eqnarray}
%
To apply the formula to sea wave modeling, we set $d=3$.
Let us use the sea modeling notation from \cite{AL121221a}.
 We have a zero-mean stationary Gaussian field $X(t,p):=\zeta(t,x,y), p:=(x,y)$, which models the sea surface.
To introduce it, let $M(\lambda_1,\lambda_2,\omega)$ be a random spectral Gaussian measure, restricted to the airy manifold $\Lambda:=\{ \normp[2]{\smash{\smash{\overrightarrow{k}}}}^2={\omega^4}/{g^2} \}$ where $\overrightarrow{ k}:=(\lambda_1,\lambda_2)$.
We define
\begin{eqnarray*}\label{Airy}\qquad\zeta(t,x,y):=\int_{\Lambda}e^{i(\lambda_1x+\lambda_2y+\omega
t)} \ud M(\lambda_1,\lambda_2,\omega).
\end{eqnarray*}

In this way, by restricting the stochastic integral to the set
$$
\Lambda^{+}:=\{(\lambda_1,\lambda_2,\omega):\omega\ge0 , \normp[2]{\smash{\overrightarrow{k}}}={\omega^{2}}/{g}\},
$$%
using polar coordinates, we can write
$$
\zeta(t,x,y)=2\int_0 ^{\infty}\int_{-\pi} ^{\pi}
\cos\!\left(\normp[2]{\smash{\overrightarrow{k}}}\cos(\theta) x
+\normp[2]{\smash{\overrightarrow{k}}} \sin(\theta) y+\omega
t\right)\ud c(\omega,\theta),
$$%
where $c$ is a random measure with independent Gaussian increments.\\
The covariance function results
\begin{align*}
 \Gamma(t,p)&:=\bbE[\zeta(0,0,0)\zeta(t,x,y)]\\
 	&=2\int_0 ^{\infty}\int_{-\pi} ^{\pi}
	\cos\!\left(\normp[2]{\smash{\overrightarrow{k}}}\cos(\theta) x +\normp[2]{\smash{\overrightarrow{k}}} \sin(\theta) y+\omega t\right)S(\omega,\theta)\ud \omega \ud \theta,
\end{align*}
where $2S(\bm{\cdot},\bm{\cdot})$ is the physical spectral density.

To establish the following results, it will be necessary to digress on the ergodic theory.
The following text has been taken from \cite{MR3012238}.
For a given subset $D\subset \reels ^2$ and for each $t>0,$ let us define $\calA_t:=\sigma\{X(\tau,p): \tau >t  , p\in D\}$ and consider the  $\sigma$-algebra of $t$-invariant sets $\calA:=\bigcap_t\calA_t.$ Moreover, assume that $\Gamma(t, p)\cvg[t\to\infty]{}0$, for all $p\in D$.
 It is well known that under this condition, the $\sigma$-algebra $\calA$ is trivial, that is, it only contains events having probability zero or one (see e.g.
\cite{MR2108670} Chapter 7).

Now for each $t>0$ and $y \in \reels$ we define the level set
$$
\calC^{\zeta}_D(t,y):=\left\{p\in D : \zeta(t,p)=y \right\}
$$%
and  the following functional
$$
\calZ(t):=\int_{\calC^{\zeta}_D(t,y)}Y(t,p) \ud\sigma_1(p).
$$%
Furthermore, in the following, we assume that
$$
	Y(t,p):=G(\zeta(t,p),\nabla_{p}\zeta(t,p)),
$$%
where $\nabla_{p}$ is the gradient operator with respect to the space variables $x,y$.
The process $\{\calZ( t) : t \in \reels ^+\}$ is strictly stationary, of finite mean and Riemann-integrable.
The ergodic theorem gives
$$
\frac1T\int_0^T\calZ( t) \ud t\cvg[T\to\infty]{a.s.}\bbE_\calB[Z(0)],
$$%
where $\calB$ is the  $\sigma$-algebra of $t$-invariant events associated with the process $\calZ(t)$.
Since for each $t$,
$\calZ(t)$ is $\calA_{ t}$-measurable, it follows that $\calB\subset\calA$ so that $\bbE_{\calB}[\calZ(0)] = \bbE[\calZ(0)]$.
Thus
\begin{align*}
 \bbE_\calB[\calZ(0)]
 	&=\bbE\!\left[\int_{\calC^{\zeta}_D(0,y)}Y(0,p) \ud \sigma_1(p)\right]\\
	&=\sigma_2(D)\bbE\left[Y(0,0)\normp[2]{\nabla_{p} \zeta(0,0)}\given{\zeta(0,0)=y}\right] \frac{e^{-\frac12{y^2}/{\lambda_{000}}}}{\sqrt{2\pi \lambda_{000}}},
\end{align*}
where $\lambda_{000}:=\bbE[\zeta^2(0,0)].$
The above formula can be used to obtain the velocity distribution as defined in (\ref{distribution}) (\cf \cite{MR3012238}).

Next, we will consider the case where the observed field, denoted $E(t,p)$, is the envelope field of $X(t,p).$ 
First, let us define the Hilbert transform of $\zeta$ as the Gaussian field
$$
\hat\zeta(t,x,y):=2\int_0 ^{\infty}\int_{-\pi} ^{\pi}
	\sin(\normp[2]{\smash{\overrightarrow{k}}}\cos(\theta) x
	+\normp[2]{\smash{\overrightarrow{k}}} \sin(\theta) y+\omega t)\ud c(\omega,\theta).
$$ 
The real envelope $E(t,x,y)$ is
$$
	E(t,p):=\sqrt{\zeta^2(t,x,y)+\hat\zeta^2(t,x,y)}.
$$%

We can write the process $E$ in the following form.
Let
$$
Z(t,p):=(\zeta(t,x,y),\hat\zeta(t,x,y)).
$$%
Then, if $F(z):=\normp[2]{z}$ then $E(t,p)=F(Z(t,p))$.
The function $F$ satisfies condition $\bB_3$  except for $z = 0$, but this does not matter because $\bbP\{Z(0,0)=0\}=0$.
Furthermore, a straightforward calculation shows  that the process $X$ verifies assumption ${(\bS)}$ in Proposition \ref{wsch}.
Then we can apply Remark \ref{Lipschitz} with condition $\bB_3$ and hypothesis ${(\bS )}$  replacing condition $\bC _3$.


Moreover, the density of $E(0,0)$, at the point $y>0$, is the Rayleigh density $(y/{\sigma^2_\zeta})e^{-\frac12{y^2}/{\sigma^2_\zeta}}$, that exists and is continuous if $\sigma^2_\zeta:=\Var(\zeta(0,0))>0$.

For each $t>0$ and $y >0$ we define the level set
$$
	\calC^{E}_D(t,y):=\left\{p\in D : E(t,p)=y \right\}
$$ 
and the functional
$$
	\calZ^Y_E(t):=\int_{\calC^{E}_D(t,y)}Y(t, p) \ud \sigma_1(p).
$$%
Invoking again the ergodic theorem, we obtain
\begin{equation}
	\label{ergodic2}
		\frac{\int^T_0Z^Y_E( t) \ud t}{\int^T_0Z^1_E(t)\ud t}\cvg[T\to\infty]{a.s.}\frac{\bbE[Y(0,0)\normp[2]{\nabla_{p}E(0,0)}\given{E(0,0)=y}]}{\bbE[\normp[2]{\nabla_{p}E(0,0)}\given{E(0,0)=y}]}.
\end{equation}
But
\begin{multline*}
	\bbE[\normp[2]{\nabla_p E(0,0)}\given{E(0,0)=y}]\\
		=\frac1{y}\bbE\left[\normp[2]{\zeta(0,0)\nabla_p \zeta(0,0)+\hat\zeta(0,0)\nabla_p\hat\zeta(0,0)}\given{E(0, 0)=y}\right].
\end{multline*}
Thus conditioning and defining 
$$
f_{\pm}(y, z_1):=\bbE\left[\normp[2]{z_1\nabla_p \zeta(0,0)\pm\sqrt{y^2-z^2_1} \nabla_p\hat\zeta(0, 0)}\right],
$$%
we get
$$
	\bbE\!\left[\normp[2]{\nabla_{p}E(0,0)}\given{E(0,0)=y}\right]
		=\frac1{y}\int_{-y}^{y}(f_+(y, z_1)+f_-(y, z_1)) p_{\zeta(0,0)}(z_1)\ud z_1
$$%
where $p_{\zeta(0,0)}(z_1)$ is the density of $\zeta(0,0)$ at $z_1$.

Given that $\hat\zeta$ has the same distribution as $-\hat\zeta$ we finally have
$$
	\bbE\!\left[\normp[2]{\nabla_{p}E(0,0)}\given{E(0,0)=y}\right]=\frac2{y}\int_{-y}^{y}f_+(y, z_1) p_{\zeta(0,0)}( z_1)\ud z_1.
$$%
Consider the numerator in the right-hand side  of the expression (\ref{ergodic2}),
\begin{align*}
\MoveEqLeft[1]{\bbE\!\left[Y(0,0)\normp[2]{\nabla_{p}E(0,0)}\given{E(0,0)=y}\right]}\\
	&= \frac1{y}\bbE\!\left[G\!\left(y,\frac1{y}\!\left(\zeta(0,0)\nabla_p \zeta(0,0)	
		+\hat\zeta(0,0)\nabla_p\hat\zeta(0,0)\right)\right)\right. \\
	&\quad \left.
	\rule{0pt}{\heightof{$\ds\frac{1}{y}$}}
	\times \normp[2]{\zeta(0,0)\nabla_p \zeta(0,0)+\hat\zeta(0,0)\nabla_p\hat\zeta(0,0)}\given{E(0,0)=y}\right].
\end{align*}
The same argument as above gives
$$
	\bbE\!\left[Y(0,0)\normp[2]{\nabla_{p}E(0,0)}\given{E(0,0)=y}\right]
		=\frac2{y}\int_{-y}^{y}\calF_+(y, z_1) p_{\zeta(0,0)}(z_1)\ud z_1,
$$%
where
\begin{multline*}
\calF_+(y,z_1):=
	\bbE\!\left[G\!\left(y,\frac1{y}\!\left(z_1\nabla_p \zeta(0, 0)+\sqrt{ y^2-z_1^2} \nabla_p\hat\zeta(0,0)\right)\right)\right.\\
	\left.\rule{0pt}{\heightof{$\ds\sqrt{y^2- z_1^2}$}}
	\times\normp[2]{z_1\nabla_p \zeta(0,0)+\sqrt{y^2- z_1^2} \nabla_p\hat\zeta(0, 0)}\right].
\end{multline*}%
Therefore, the right-hand side term of (\ref{ergodic2}) becomes equal to
$$
{\displaystyle\int_{-y}^{y}\calF_+(y, z_1) p_{\zeta(0,0)}(z_1)\ud z_1}\left/{\displaystyle\int_{-y}^{y}f_+(y, z_1) p_{\zeta(0,0)}(z_1)\ud z_1}\right..
$$%
In some important cases, this term can be calculated explicitly using only the spectral moments of the processes $\zeta, \hat \zeta, \nabla_{p}\zeta, \nabla_{p}\hat\zeta.$
\section{Berry and Dennis dislocations}
In this part of the work, we will give an overview of the applications of the Rice formula to some physics notions known as random wave dislocations.
This study is motivated by the seminal paper by Berry and Dennis, \cite{MR1843853}, several new concepts were introduced on physical grounds.

We consider two independent isotropic Gaussian random fields of mean zero belonging to $C^{2}(D,\reels)$.
Let $\xi, \eta:\Omega\times\reels ^2\to \reels ,$ defined trough their spectral representation
\begin{align*}
	\xi(\bx )
	 	&:=\int_{\reels ^2}\cos(\prodsca{\bx}{\bk})\left({\Pi(k)}/{k}\right)^{\frac12}\ud W_1(\bk )\\
		&\rule{96pt}{0pt}-\int_{\reels ^2}\sin(\prodsca{\bx}{\bk})\left({\Pi( k)}/{k}\right)^{\frac12}\ud W_2(\bk )\\\shortintertext{and}
	\eta(\bx )
		&:=\int_{\reels ^2}\cos(\prodsca{\bx}{\bk})\left({\Pi(k)}/{k}\right)^{\frac12}\ud W_2(\bk )\\
		&\rule{96pt}{0pt}+\int_{\reels ^2}\sin(\prodsca{\bx}{\bk})\left({\Pi(k)}/{k}\right)^{\frac12}\ud W_1(\bk ),
\end{align*}
where $\prodsca{\cdot }{\cdot}$ stands for the scalar product in $\reels^2$ and
$\bk :=(k_1,k_2)$, $k:=\normp[2]{\bk }$, $\Pi(k)$ is the isotropic spectral density and $W:=W_1+iW_2$ is a standard complex orthogonal Gaussian measure on $\reels ^2$.
Without loss of generality, we can assume that 
$E\left[\xi^2(\0)\right]=E[\eta^2(\0)]=1$.\\
Defining the complex wave
$\psi(\bx ):=\xi(\bx )+i\eta(\bx )$,
the dislocations is the set of zeros of $\psi$, \ie
$$
	N^\psi_D(0):=\#\{\bx : \psi(\bx )=0\}=  \#\{\bx : \xi(\bx )= \eta(\bx )=0\}.
$$%
In \cite{MR1843853} (see formulas (2.7) and (4.6)), the expected number of dislocation points by unit of area is defined by
\begin{multline*}
 d_2:=\frac{\bbE[\#\{\bx \in D: \psi(\bx )=0\}]}{\sigma_2(D)}\\
	=\frac{\lambda_2}{(2\pi)^2}\bbE\!\left[\abs{\frac{\xi_x(\0)}{\sqrt\lambda_2}\frac{\eta_y(\0)}{\sqrt\lambda_2}-\frac{\xi_y(\0)}{\sqrt\lambda_2}\frac{\eta_x(\0)}{\sqrt\lambda_2}}\right]=
\frac{\lambda_2}{2\pi},
\end{multline*}
where $\xi_x$, $\xi_y$, $\eta_x$ and $\eta_y$ stand for the derivatives of first order of $\xi$ and $\eta$ and $\lambda_2:=\bbE[\xi_x^2(\0)]=\bbE[\xi_y^2(\0)]=\bbE[\eta_x^2(\0)]=\bbE[\eta_y^2(\0)]$ .\\
Here, we will also study the length of the set of zeros of each coordinate process (length of nodal curves)
$$
\sigma_1(C_\xi(0))\stackrel{{\calL}}{=}\sigma_1(C_\eta(0)).
$$%
We thus have the definition of the length of the nodal curves for the surface unit:
$$
\calL:=\frac{\bbE[\sigma_1(C_\xi(0))]}{\sigma_2(D)}=\frac{\bbE[\sigma_1(C_\eta(0))]}{\sigma_2(D)}.
$$%
In \cite{MR1843853}, other notions have been defined related to the following two integrals
$$
\smashoperator{\int\limits_{\{\bx \in D: \psi(\bx )=0\}}}Y(\bx )\ud \sigma_0(\bx )
=\smashoperator[r]{\sum_{\bx \in\{\bx \in D: \psi(\bx )=0\}}}Y(\bx )\quad\mbox{ and }\quad \int\limits_{C_\xi(0)}Y(\bx )\ud \sigma_1(\bx ).
$$%
For the first one, we must recall that $\sigma_0$ is the counting measure.
For instance in \cite{MR1843853} the dislocation curvature is introduced.
In what follows, we will consider instead the curvature of one of the nodal curves, defined using for example  $\xi.$
The curvature of the nodal curve $\xi(\bx )=\xi(x,y)=0$, is
$$
\kappa(\bx ):=\frac{\abs{\xi_{xx}(\bx )\xi^2_x(\bx )-2\xi_{xy}(\bx )\xi_x(\bx )\xi_y(\bx )+\xi_{yy}(\bx )\xi_y^2(\bx )}}{\normp[2]{\nabla\xi(\bx )}^3}.
$$%
For the interval $[0,\kappa_1]$ defining $Y(\bx ):=\1_{[0,\kappa_1]}(\kappa(\bx ))$, one obtains a particular case of the function $Y(\bx )=G(\nabla\xi(\bx ),\nabla^2\xi(\bx ))$, 
where the operator $\nabla^2$ denotes the second order differential.
For these functions, in a manner similar to that of Theorem~\ref{Cabana4}, we can prove a Rice formula obtaining
\begin{align*}
\MoveEqLeft[1]{\bbE\!\left[\int_{C_\xi(0)}\1_{[0,\kappa_1]}(\kappa(\bx ))\ud \sigma_1(\bx )\right]}\\
 	&\quad=\sigma_2(D)\bbE\!\left[\1_{[0,\kappa_1]}(\kappa(\0))\normp[2]{\nabla \xi(\0)}\given{\xi(\0)=0}\right] p_{\xi(\0)}(0)\\
 	&\quad=\frac{\sigma_2(D)}{\sqrt{2\pi}} \bbE\!\left[\1_{[0,\kappa_1]}(\kappa(\0))\normp[2]{\nabla \xi(\0)}\given{\xi(\0)=0}\right].
\end{align*}

The independence between $\nabla\xi(\0)$ and $(\xi(\0),\nabla^2\xi(\0))$  allows writing a regression model that simplifies the last expression.
Moreover, 
\begin{align*}
 \bbE[\sigma_1(C_\xi(0))]
	&=\frac{\sqrt{\lambda_2}\sigma_2(D)}{\sqrt{2\pi}}\bbE\!\left[\sqrt{\frac{\xi_x^2(\0)}{\lambda_2}+\frac{\xi_y^2(\0)}{\lambda_2}}\right]\\
	&=\frac{\sqrt{\lambda_2}\sigma_2(D)}{\sqrt{2\pi}}\frac1{2\pi}\int_0^\infty\int_0^{2\pi}\rho^2e^{-\frac12\rho^2}\ud \theta \ud \rho=\frac{\sqrt{\lambda_2}\sigma_2(D)}2.
\end{align*}
As a bonus, we get $\calL={\sqrt{\lambda_2}}/2$.

Furthermore, in order to obtain an interpretation for the distribution of $\kappa$ over the level set of $\xi$, we take the ratio of the two last  expectations obtaining
\begin{multline*}
\bbE\!\left[\int_{C_\xi(0)}\1_{[0,\kappa_1]}(\kappa(\bx ))\ud \sigma_1(\bx )\right]\left/\rule{0pt}{16pt}{\bbE[\sigma_1(C_\xi(0))]}\right.\\
 =\frac{1}{\sqrt{\lambda_2}}\sqrt{\frac{2}{\pi}}\bbE\!\left[\1_{[0,\kappa_1]}(\kappa(\0))\normp[2]{\nabla \xi(\0)}\given{\xi(\0)=0}\right].
\end{multline*}
Using independence, we can write it as
\begin{align}
\label{reg2}
\sqrt{\frac{2}{\pi}}\int_0^\infty\int_0^{2\pi}\bbE\!\left[\1_{[0,\rho \sqrt{\lambda_2}\kappa_1]}\left(\rule{0pt}{10pt}\left|\xi_{xx}(\0)\cos^2\theta-2\xi_{xy}(\0)\cos\theta\sin\theta\right.\right. \nonumber\right.\\
\left.\left.\left.
	\rule{0pt}{10pt}
	+\xi_{yy}(\0)\sin^2\theta \right|\right)\given{\xi(\0)=0}\rule{0pt}{12pt}\right]
\times\rho^2\frac{e^{-\frac{\rho^2}2}}{2\pi}\ud \rho \ud \theta.
\end{align}
A regression model shows that the following relationship is true
\begin{multline*}
[\xi_{xx}(\0)\cos^2\theta-2\xi_{xy}(\0)\cos\theta\sin\theta+\xi_{yy}(\0)\sin^2\theta\given{\xi(\0)=0}]\\
\stackrel{{\calL}}{=}{\calN}(0,\sigma^2(\theta,\lambda_4,\lambda_{22},\lambda_2)),
\end{multline*}
where 
\begin{align*}
	\lambda_4     &:=\bbE[\xi_{xx}^2(\0)]=\bbE[\xi_{yy}^2(\0)]\\
	\lambda_{22} &:=\bbE[\xi_{xy}^2(\0)]\\ \intertext{and}\\
	-\lambda_2 &\phantom{:}=\bbE[\xi(\0)\xi_{xx}(\0)]=\bbE[\xi(\0)\xi_{yy}(\0)].
\end{align*}
Then
\begin{align*}
 \text{Eqn. (\ref{reg2})}
 	&\phantom{:}=\frac{1}{2\pi^2}\int_0^\infty\int_0^{2\pi}\int_{-\kappa_1\sqrt{\lambda_2}\rho/\sigma(\theta,\lambda_4,\lambda_{22},\lambda_2)}^{\kappa_1\sqrt{\lambda_2}\rho/\sigma(\theta,\lambda_4,\lambda_{22},\lambda_2)}\rho^2e^{-{\rho^2}/2}e^{-{u^2}/2}\ud \rho \ud \theta \ud u\\
&\phantom{:}=\frac{1}{\pi^2}\int_0^\infty\int_0^{2\pi}\int_0^{\kappa_1\sqrt{\lambda_2}\rho/\sigma(\theta,\lambda_4,\lambda_{22},\lambda_2)}\rho^2e^{-{\rho^2}/2}e^{-{u^2}/2}\ud \rho \ud \theta \ud u\\
&:=\calK(\kappa_1).
\end{align*}
The density of this distribution is
\begin{align*}
\frac{d}{d\kappa_1}\calK(\kappa_1)&=\frac{\sqrt{\lambda_2}}{\pi^2}\int_0^\infty\int_0^{2\pi}\rho^3/\sigma(\theta,\lambda_4,\lambda_{22},\lambda_2)e^{-{\rho^2}/2}\\
	 &\specialpos{\hfill\times e^{-\frac12{(\kappa_1\sqrt{\lambda_2}\rho/\sigma(\theta,\lambda_4,\lambda_{22},\lambda_2))^2}}\ud \rho \ud \theta}\\
	 &=\frac{\sqrt{\lambda_2}}{\pi^2}\int_0^\infty\int_0^{2\pi}\frac{\sigma^3(\theta,\lambda_4,\lambda_{22},\lambda_2)}{(\sigma^2(\theta,\lambda_4,\lambda_{22},\lambda_2)+\kappa_1^2\lambda_2)^{2}} v^3e^{-v^2/2} \ud v \ud \theta\\
	& =\frac{2\sqrt{\lambda_2}}{\pi^2}\int_0^{2\pi}\frac{\sigma^3(\theta,\lambda_4,\lambda_{22},\lambda_2)}{(\sigma^2(\theta,\lambda_4,\lambda_{22},\lambda_2)+\kappa_1^2\lambda_2)^{2}}\ud \theta.
\end{align*}

The last part of this subsection aims at computing some second-order Rice formulas.
Let us first introduce the correlation of dislocations at distance $R$ defined as $g(R)$  in \cite{MR1843853}.
To define this quantity, we first consider the second factorial moment of the random variable $N^\psi_D(0)$ that is
\begin{multline*}
	\bbE\!\left[\#\{\bx \in D: \psi(\bx )=0\}(\#\{\bx \in D: \psi(\bx )=0\}-1)\right]\\
	=\int\int_{D\times D}A(\bx _1,\bx _2) {p}_{\psi(\bx _1), \psi(\bx _2)}(0, 0) \ud \bx _1\ud \bx _2,
\end{multline*}
where using Rice's formula, we get
$$
	A(\bx _1,\bx _2):=\bbE[\abs{\det\nabla\psi(\bx _1)}\abs{\det\nabla\psi(\bx _2)}\given{\psi(\bx _1)=\psi(\bx _2)=0}].
$$%
Using invariance with respect to rotations and translations, it turns out that
$$
	A(\bx _1,\bx _2)=A((0,0),(\normp[2]{\bx _1-\bx _2},0))=A((0,0),(R,0)):=g(R).
$$%
In the last equality we have set $R:=\normp[2]{\bx _1-\bx _2}.$ 

Moreover, by dropping the absolute value of the determinant of the Jacobian of $\psi$ we can introduce
$$
	B(\bx _1,\bx _2):=\bbE[\det\nabla\psi(\bx _1)\det\nabla\psi(\bx _2)\given{\psi(\bx _1)=\psi(\bx _2)=0}].
$$%
Thus, the {\it charge correlation function} (cf. \cite{MR1843853}) is defined as
$$
	g_Q(R):=B((0,0),(R,0)).
$$%

A elementary closed expression for $g(R)$ was obtained in \cite{MR3012238, MR1843853} using an expression for the absolute value function as a Fourier integral.
Nevertheless, the computation is not trivial.
The interested reader can consult these references .
Also, the function $g_Q(R)$ can be written as the conditional expectation of a sum of products of four standard Gaussian random variables.
Consider for example the first term, \ie
$$
	\bbE\!\left[\xi_x(0,0)\eta_y(0,0)\xi_x(R,0)\eta_y(R,0)\given{\psi(0,0)=\psi(R,0))=0}\right].
$$%
Then, the random variables representing the derivatives are regressed on the vector $(\psi(0,0),\psi(R,0))$. An elementary Gaussian calculation gives the result.
\section{Gravitational stochastic microlensing}
In this section, we only sketch an application to gravitational cosmology.
The main reason to present it is that the Rice formula used is shown for a non-Gaussian process.
However, it is more an illustration than a real formal mathematical development.

It should be noted that all of the material corresponding to this subsection comes from the article by Peters et al.
\cite{MR2582583}.
In addition, for background, this work had to be supplemented by the book \cite{MR1836154}.

Let $\{\xi_i\}$ be $g$ independent random variables identically distributed on the disk of radius $R$ in $\reels ^2$.
They are considered as the positions of the stars.
We can define the following random field
$$
	\psi_g(x):=\tfrac{1}{2}\kappa_c\normp[2]{x}^2-\tfrac12{\gamma}(x^2_1-x^2_2)+m\sum_{j=1}^g\ln\normp[2]{x-\xi_j}^2,
$$%
where $x:=(x_1,x_2)$ and $\kappa_c$, $\gamma$ are physical constants and $m$ represents the mass of the stars.
Outside of the random points $\{\xi_i\}_{=1}^g$ the potential $\psi_g$ is $C^\infty(\reels^2, \reels)$.
 The following random function is known as the delay function for the gravitational lens systems
$$
	T_{y}(x):=\tfrac12{\normp[2]{x-y}^2}-\psi_g(x).
$$%
The lensing map is defined as
$$
	\eta(x):=\nabla T_{y}(x) + y=x-\nabla \psi_g(x).
$$%
Given the definitions, we easily obtain
$$
\eta(x)=((1-\kappa_c+\gamma)x_1,(1-\kappa_c-\gamma)x_2)-2m\sum_{j=1}^g\frac{x-\xi_j}{\normp[2]{x-\xi_j}^2}.
$$%

A lensed image is a solution $x^\star$ of the equation $\nabla T_y(x)=0.$ That is
$$
\eta(x^\star)=y.
$$%

These images correspond to the stationary points of the function $T_y$ and are classified as local maximum, local minimum and saddle point whenever the image is not degenerated.
Otherwise,  we say that they are degenerate.
We are interested in computing the number of non-degenerate images with positive parity $N_{+}$ which are defined as  
$N_{+}:=N_{\max}+N_{\min}$.
It is easy to show that in these images, the Jacobian of $\eta$, \ie $\det(\nabla\eta(x))$, is always positive.

It is interesting, for gravitational studies, to calculate the expected number of $N_+$ generated for a point source $y$.
The number of such images on a set $D\subset\reels ^2$ is
$$
N_+(y):=\#\{x\in D: \eta(x)=y, \det(\nabla\eta(x))>0\}.
$$%
Consider $f:\reels ^2\to\reels $  a  continuous function with bounded support, then by the area formula we obtain
\begin{align*}
\MoveEqLeft[1]{\int_{\reels ^2}f(y) \bbE\!\left[N_+(y)\right]\ud y
	=\int_{D}\bbE\!\left[f(\eta(x))\det(\nabla \eta(x))\1_{\introo{0}{+\infty}}(\det(\nabla\eta(x)))\right] \ud x
}\\ 
&=\int_{\reels ^2}f(y)\int_{D} \bbE\!\left[\det(\nabla\eta(x))\1_{\introo{0}{+\infty}}(\det \nabla\eta(x))\given{\eta(x)=y}\right] {p}_{\eta(x)}(y)\ud x \ud y.
\end{align*}

Although the function $\eta$ has singularities in the positions of the stars $\xi_i$ these are infinite singularities.
That is $\lim_{x\to\xi_i}\normp[2]{\eta(x)}=+\infty$.
Therefore, if we observe only those $y$ that are in the bounded support of $f$, we have that the domain of $\eta$ for each $\omega$ is restricted to an open set that does not contain the points $\xi_i$.
This implies that the function $\eta$ restricted to this set is a $C^{\infty}$ function.
Then the hypothesis for applying the area formula holds.

Moreover, we get for almost all $y$
$$
\bbE\!\left[N_+(y)\right]=\int_{D}\bbE\!\left[\det(\nabla\eta(x)) \1_{\introo{0}{+\infty}}(\det(\nabla\eta(x)))\given{\eta(x)=y}\right] {p}_{\eta(x)}(y) \ud x.
$$%

Using the definitions and some non-trivial work, it can be shown that the above formula holds for all $y.$ Moreover, in  \cite{MR2582583} the formula is used to obtain its asymptotic when the number of stars $g$ tends to infinity.
An interesting but still open problem is to get the same asymptotic for the variance of $N_+(y)$.
\section{Kostlan-Shub-Smale systems}
Consider a rectangular system $\bP =\0$ of $j$ homogeneous polynomial equations in $d > j$ variables.
\\
We assume that the equations have the same degree $n>1$.

Let $\bP :=(X_1,\dots,X_j)$, we can write each polynomial $X_\ell$ in the form
$$
  X_{\ell}(\bt ):=\sum_{\abs{\bz}=n}a^{(\ell)}_{\bz }\bt ^{\bz },
$$%
where 
\begin{enumerate}
  \item $\bz :=(z_1,\dots,z_{d})\in\naturels^{d}$ 
    and $\abs{\bz }:=\sum^{d}_{k=1}z_{k}$; 
  \item $a^{(\ell)}_{\bz }:=a^{(\ell)}_{z_1\dots z_d}\in\reels $, $\ell=1,\dots,j$, $\abs{\bz }= n$;
  \item $\bt :=(t_1,\dots,t_{d})$ 
    and $t^{\bz }:=\prod^d_{k=1} t^{z_k}_{k}$.
\end{enumerate}

We say that  $\bP $ has the Kostlan-Shub-Smale 
(KSS for short) distribution if the coefficients $a^{(\ell)}_{\bz }$ are independent random variables of zero mean normally distributed with variances
$$
  \Var(a^{(\ell)}_{\bz }):=\binom{n}{\bz }:=\frac{n!}{ 
z_1!\dots z_d!}.
$$%
\subsection{Expectation of the volume}
We are interested in the set of zeros of $\bP $ and we denote this set by $\calC_{\bP }(\0)$ and its volume by $\calL_n(\calC_{\bP }(\0))$  if $d-1>j$.
Its cardinality is $N^{\bP }_n$ if $d-1=j$.
Shub and Smale \cite{MR1230872} showed
 that if $d-1=j$ then $\bbE\!\left[N^{\bP }_n\right]=2n^{(d-1)/2}$.
In \cite[Chapter 12]{MR2478201}, this result was obtained by using the Kac-Rice formula.
Letendre in \cite{MR3470435} tackled the case $d-1>j$, \ie the case of homogenous polynomials of degree $n$ in $d$ variables, obtaining the following result first shown by Kostlan in \cite{MR1246137}
$$
	\bbE\!\left[\calL_n(\calC_{\bP }(\0))\right]
		=2 n^{{j}/2}{\pi^{(d-j)/2}}\left/{\Gamma[\tfrac12{(d-j)}]}\right.,
$$%
where $\Gamma$ is the Gamma function.\\
Following Letendre's method and making some simplifications, we will obtain this result using the Kac-Rice formula.

Each $X_\ell$ is homogeneous, and the zero set of $\bP $ is the intersection of the zero sets of $X_\ell$.
Then the  set $\calC_{\bP }(\0)$ is a subset of the real projective space $\reels \bbP^{d-1-j}$.\\
Standard multinomial formula shows that for all $\bs ,\bt \in\reels^{d}$ we have
$$
r_n(\bs ,\bt ):=  \bbE\!\left[{X_{\ell}(\bs )X_{\ell}(\bt )}\right]
=\prodsca{\bs }{\bt }^{n},
$$%
where $\prodsca{\cdot}{\cdot}$ is the usual inner product in $\reels^{d}$.

Consequently, we see that the distribution of the system $\bP$ is invariant under the action of the 
orthogonal group in $\reels^{d}$.
We also see that the distribution depends of course on $n$ and this will be omitted for $\bP$ for ease of notation.
Let us observe that the parameter $\bt $ can be considered in the unit sphere of $\reels^d$, that is $\bbS ^{d-1}$.

In the following, we must consider the derivative of $X_\ell$, $\ell=1,\dots,j$.

Since the parameter space is the unit sphere $\bbS ^{d-1}$,  the derivative is taken in the sense of the sphere, that is, the spherical derivative $X^{\prime}_\ell(\bt )$ of $X_\ell(\bt )$ is the orthogonal projection of the free gradient on the tangent space $\bt ^{\perp}$ of $\bbS ^{d-1}$ at $\bt $.

The $k$-th component of $X^{\prime}_\ell(\bt )$ with respect to a given basis of the tangent space is denoted by $X^{\prime}_{\ell k}(\bt )$.

We will use Rice's formula slightly modified to make it valid on $\bbS ^{d-1}$.
As the process $\bP $ satisfies the hypotheses of Remark \ref{Lipschitz} according to Theorem \ref{Cabana4}  we obtain
\begin{multline*}
 \bbE\!\left[\calL_n(\calC_{\bP }(\0))\right]\\
=\int_{\bbS ^{d-1}}\bbE\!\left[\left(\det(\nabla\bP(\bt )\nabla\bP(\bt )^T)\right)^{\frac12}\given{\bP(\bt )=\0}\right] {p}_{\bP(\bt )}(\0)\sigma_{d-1}(\ud\bt ),
\end{multline*}
where ${p}_{\bP(\bt )}(\0)$ represents the density of $\bP(\bt )$ in $\0$.
 \\
Since 
$\bbE[X_{\ell}^2(\bt )]=1$, 
$\bP(\bt )$ and $\nabla \bP(\bt )$
are independent, which allows to release the conditioning into the expectation.
Moreover, since ${p}_{\bP(\bt )}(\0)=1/\sqrt{(2\pi)^j}$, if  $\nabla\bP(\be_0)$ represents the generic element matrix $X^{\prime}_{\ell k}(\0)$, we finally obtain
\begin{align*}
\bbE[\calL_n(\calC_{\bP}(\0))]
	&=\frac{\sigma_{d-1}(\bbS ^{d-1})}{(2\pi)^{\frac j2}}\bbE\!\left[\left(\det (\nabla\bP(\be_0)\nabla\bP(\be_0)^T)\right)^{\frac12}\right]\\
	&=n^{\frac j2}\frac{\sigma_{d-1}(\bbS ^{d-1})}{(2\pi)^{\frac j2}}\bbE\!\left[\left(\det (\nabla\bz (\be_0)\nabla\bz (\be_0)^T)\right)^{\frac12}\right],
\end{align*}
where $\nabla\bz (\be_0)$ is ${\calN}(0,I_{(d-1)j})$ and $I_{(d-1)j}$ represents the $(d-1)j \times (d-1)j$ identity matrix.
\\
A Gaussian calculation \cite[Lemma 5.4. p.\,3077]{MR3470435}) gives 
$$
\bbE\!\left[\left(\det (\nabla\bz (\be_0)\nabla\bz (\be_0)^T)\right)^{\frac12}\right]
	=(2\pi)^{j/2}\frac{\sigma_{d-1-j}(\bbS ^{d-1-j})}{\sigma_{d-1}(\bbS ^{d-1})},
$$%
yielding
$$
	\bbE\!\left[\calL_n(\calC_{\bP}(\0))\right]
		=n^{j/2}\sigma_{d-1-j}(\bbS ^{d-1-j})
		=2 n^{{j}/2}{\pi^{(d-j)/2}}\left/{\Gamma[\tfrac12{(d-j)}]}\right..
$$%

We now calculate the variance of this random variable.
Note that it has also been computed by Letendre in \cite{MR3917219} and Letendre and Puchol in \cite{MR4052739}.

\subsection{Asymptotic of the volume variance}
Letendre in \cite{MR3917219} p.4 and Letendre and Puchol in \cite[p.\,3]{MR4052739} studied the asymptotic variance of the volume of the zero set when the degree $n$ goes to infinity and in the case where $d-1 >j$.
This result is sketched in \cite{MR3866880-b}.
We will consider this problem using a  different method, the case where $d-1=j$ being treated in \cite{MR3866880}.

Let us start with some remarks on notation.
Recall that we denote by $\bbS ^{d-1}$ the unit sphere in $\reels^{d}$.
Let us denote by $\kappa_d$ the hypervolume of the $(d-1)$-dimensional unit sphere.\\
The variables $\bs $ and $\bt$ denote points on $\bbS ^{d-1}$ and
$\ud \bs$ and $\ud \bt$ denote the corresponding geometric measure.

\paragraph{Hyperspherical coordinates:} 
For $\theta:=(\theta_1,\dots,\theta_{d-2},\theta_{d-1})\in\intrfo{0}{\pi}^{d-2}\times\intrfo{0}{2\pi}$ 
we write $x^{(d-1)}(\theta):=(x^{(d-1)}_1(\theta),\dots,x^{(d-1)}_{d}(\theta))\in \bbS ^{d-1}$ in the following way
$$
x^{(d-1)}_k(\theta) :=
\begin{cases}
 \ds\prod^{k-1}_{j=1}\sin(\theta_j)\cdot\cos(\theta_k), & k\leq d-1\\
 \ds\prod^{d-1}_{j=1} \sin(\theta_j), & k =  d
\end{cases}
$$
with the convention that $\prod^0_1=1$.

We will use several times in the following that for $h:\intrff{-1}{1}\to\reels$ a continuous function, it turns out that
\begin{equation}\label{eq:intS}
\int_{\bbS ^{d-1}\times \bbS ^{d-1}}h(\prodsca{\bs}{\bt})\ud \bs  \ud \bt 
=\kappa_d\kappa_{d-1}\int^\pi_0\sin^{d-2}(\theta)h(\cos(\theta))\ud \theta.
\end{equation}
The proof consists in the use of a reference system with hyperspherical coordinates.

For $\ell=1,\ldots,j$, we define the standardized derivative as
\begin{equation}\label{notacionvectorial}
 \overline X^{\prime}_\ell(\bt  ):=\frac{X^{\prime}_\ell(\bt  )}{\sqrt{n}}\quad\mbox{and}\quad \overline \bP^{\prime}(\bt ):=(\overline X^{\prime}_1(\bt ),\ldots,\overline X^{\prime}_j(\bt )),
\end{equation}
where $\overline X^{\prime}_{\ell}(\bt )$ is a row vector.
For $\bt \in \bbS ^{d-1}$, we also define
$$\label{eq:zeta}
  \overline Z(\bt ):=(Z_1(\bt ),\ldots,Z_{jd}(\bt )):=(\bP(\bt ), \overline \bP^{\prime}(\bt )).
$$%
The covariances 
$$
	\rho_{k\ell}(\bs ,\bt ):=\bbE(Z_k(\bs )Z_\ell(\bt )),\quad k,\ell=1,\ldots,jd,
$$%
are obtained by routine calculations.
They are simplified using  the invariance under isometries.
For instance, if $k=\ell\leq j$
$$
	\rho_{k\ell}(\bs ,\bt )=\prodsca{\bs}{\bt}^n:=\cos^n(\theta),\quad\theta\in\intrfo{0}{\pi},
$$%
where $\theta$ is the angle between $\bs$ and $\bt$.

When the indices $k$ or $\ell$ are larger than $j$, the covariances involve derivatives of $r_n$.

In fact, we can prove that  $\overline Z$ is a vector of $jd$ standard normal random variables 
whose covariances depend upon the quantities 
\begin{align}\label{eq:deriv}
\calA(\theta)&:=-\sqrt{n}\cos^{n-1}(\theta)\sin(\theta),\\
\calB(\theta)&:=\cos^{n}(\theta)-(n-1)\cos^{n-2}(\theta)\sin^2(\theta),\notag\\
\calC(\theta)&:=\cos^{n}(\theta),\notag\\
\calD(\theta)&:=\cos^{n-1}(\theta).\notag
\end{align}

Thus, we can write the variance-covariance matrix of the vector
$$
	\!\left(X_\ell(\bs ),X_\ell(\bt ),\frac{X^{\prime}_\ell(\bs )}{\sqrt n},\frac{X^{\prime}_\ell(\bt )}{\sqrt n}\right),
$$%
 in the following form
\begin{eqnarray*}\label{matrix} 
\renewcommand{\arraystretch}{1.25} 
\left[\begin{array}{c;{1pt/2pt}c;{1pt/2pt}c} 
A_{11}&A_{12} &A_{13}\\[2pt] \hdashline[1pt/2pt]
A_{12}^T&I_{d-1}  & A_{23}\\[2pt]  \hdashline[1pt/2pt]
A_{13}^T&A_{23}  &I_{d-1}\\
\end{array}\right],
\end{eqnarray*}
where,
$$
A_{11}:=\left[\begin{array}{cc}
1&\calC\\
\calC&1
\end{array}\right], \;
A_{12}:= \left[\begin{array}{cccc}
0& 0&\dots& 0\\
-\calA&0& \dots & 0
\end{array}\right],\:
A_{13}:=  \left[\begin{array}{cccc}
\calA& 0&\dots& 0\\
0&0& \dots & 0
\end{array}\right],
$$%
and $A_{23}:=\diag( [\calB,\calD,\ldots,\calD])_{(d-1)\times (d-1)},$ where $\diag([a_1, a_2, \dots, a_k])_{k \times k}$ stands for the $k$-square diagonal matrix with the generic element $a_i$ on its diagonal.

Moreover, when dealing with the conditional distribution of $( \overline \bP^{\prime}(\bs ), \overline \bP^{\prime}(\bt ))$ given that $\bP(\bs )=\bP(\bt )=0$ the following expressions appear for the common variance $\sigma^2$
and the correlation $\rho$ depending on $\theta$
\begin{equation}
\label{correlation}
\sigma^2(\theta):=1-\frac{\calA^2(\theta)}{1-\calC^2(\theta)};\quad
\rho(\theta):=\frac{\calB(\theta)(1-\calC^2(\theta))-\calA^2(\theta)\calC(\theta)}{1-\calC^2(\theta)-\calA^2(\theta)}.
\end{equation}

After scaling $\theta:=z/\sqrt{n}$, we obtain the following bounds using the definitions given above.
The proof is elementary and the reader can consult the article \cite{MR3866880} for details.
\spacebefore
\begin{lemm}\cite{MR3866880}.\label{l:bounds}
There exists $0<\alpha<\frac12$ such that for $n$ sufficiently large and ${z}/{\sqrt{n}}<{\pi}/{2}$ 
it holds that
\begin{align*}
	\abs{\calA}&\leq z\exp(-\alpha z^2),\\
	\abs{\calB}&\leq (1+z^2)\exp(-\alpha z^2),\\
	\abs{\calC}&\leq
	\abs{\calD}\leq \exp(-\alpha z^2),
\end{align*}
and for $z \ge z_0$,
\begin{align*}
	0\leq 1-\sigma^2&\leq \bC  z^2 \exp(-2\alpha z^2),\\
	\abs{\rho}&\leq \bC  (1+z^2)\exp(-\alpha z^2).
\end{align*}
All functions on the left-hand side are evaluated at $\theta=z/\sqrt{n}$.
\end{lemm}
\spacebefore
\begin{lemm}\cite{MR3866880}.\label{l:limits} It is easy to show using the definitions that the following limits hold
\begin{eqnarray*}
  \cos^{2n}\!\left(z/{\sqrt n}\right)&\cvg{}&  \exp(-z^2)\\
    \calA &\cvg{}& -z  \exp(-z^2/2)\\
  \calB &\cvg{}& (1-z^2) \exp(-z^2/2)\\
  \calC \mbox{ and }\calD &\cvg{} & \exp(-z^2/2)\\
    \sigma^2\!\left(z/{\sqrt n}\right) &\cvg{}&   \frac{1-(1+z^2)\exp(-z^2)}
   {1-\exp(-z^2)} \\ 
 \rho\!\left( z/{\sqrt n}\right)  &\cvg{}&   \frac{\exp(-z^2/2)(1-z^2-\exp(-z^2))}
   {1-(1+z^2)\exp(-z^2)}.
\end{eqnarray*}
\end{lemm}

Let us apply the Rice formula for the second moment of the zero set volume by writing the variance as
$$
	\Var\left[\calL_n(\calC_{\bP}(\0))\right]=\bbE\!\left[\calL_n^2(\calC_{\bP}(\0))\right]-(\bbE\!\left[\calL_n(\calC_{\bP}(\0))\right])^2.
$$%
In the expression above, we have already computed  the second term and for the first, we apply Rice's formula for the second moment with a slight modification as before to make it valid on $\bbS ^{d-1}$ (see \cite{MR2478201}).
\begin{multline*}
 \bbE\!\left[\calL_n^2(\calC_{\bP}(\0))\right]\\
 =\int_{\bbS ^{d-1}\times\bbS ^{d-1}}\bbE\!\left[\!\left(\det(\bP^{\prime}(\bt )\bP^{\prime}(\bt )^T)\right)^{\frac12}\!\left(\det(\bP^{\prime}(\bs )\bP^{\prime}(\bs )^T)\right)^{\frac12} |\bP( \bt )=\bP(\bs )=\0\right]\\
 \times {p}_{\bP(\bt ),\bP(\bs )}(\0, \0) \ud \bt  \ud \bs ,
\end{multline*}
where ${p}_{\bP(\bt ),\bP(\bs )}( \cdot,\cdot)$ is the density of the vector $(\bP(\bt ),\bP(\bs )).$  

By independence, we can write%
$$
{p}_{\bP(\bt ),\bP(\bs )}(\0,\0)=\prod_{\ell=1}^j {p}_{X_\ell( \bt ),X_\ell(\bs )}(0,0)=\frac1{(2\pi)^j(1-\prodsca{\bt }{\bs }^{2n})^{j/2}}.
$$%
Moreover, let
$\bP^{\prime}( \bt )\bP^{\prime}( \bt )^T=(a_{ik}).$
Thus, using the homogeneity  of the determinant, we have
$n^{-j}\det(\bP^{\prime}(\bt )\bP^{\prime}(\bt )^T)=\det({a_{ik}}/{n})$. 
Obtaining in this form that all the entries of the matrix $({a_{ik}}/{n})$ have variance one.

Furthermore,  we have
\begin{multline}\label{conditional}
	\bbE\!\left[\!\left(\det(\bP^{\prime}(\bt )\bP^{\prime}(\bt )^T)\right)^{\frac12}\!\left(\det(\bP^{\prime}(\bs )\bP^{\prime}(\bs )^T\right))^{\frac12} |\bP( \bt )=\bP(\bs )=\0\right]\\
	=n^{j}\bbE\!\left[\!\left(\det(\overline \bP^{\prime}(\bt )\overline \bP^{\prime}(\bt )^T)\right)^{\frac12}\!\left(\det(\overline \bP^{\prime}(\bs )\overline \bP^{\prime}(\bs )^T)\right)^{\frac12} |\bP(\bt )=\bP(\bs )=\0\right],
\end{multline}
where $\overline \bP^{\prime}$ is defined by (\ref{notacionvectorial}).

Let us introduce the following function $f^j$
$$
	f^j( y^{\prime}_{11},\ldots,y^{\prime}_{1d-1},\ldots,y^{\prime}_{j1},\ldots,y^{\prime}_{jd-1})
		:=f^j(\by^{\prime}):=\sqrt{\det(\by^{\prime}{\by^\prime}^T)},
$$ 
where $\by^{\prime}:=(y^{\prime}_{11},\ldots,y^{\prime}_{1d-1},\ldots,y^{\prime}_{j1},\ldots,y^{\prime}_{jd-1})\in \reels^{j\times (d-1)}.$ Furthermore, for any $\gamma>0$ let us define
$$
	f^j_\gamma(\by^{\prime})
		:=f^j( \gamma y^{\prime}_{11},y^{\prime}_{12},\ldots, y^{\prime}_{1d-1},\gamma y^{\prime}_{21},\ldots,y^{\prime}_{2d-1},\cdots,\gamma y^{\prime}_{j1},\ldots,y^{\prime}_{jd-1}).
$$%
For $\gamma=1$ we have $f^j=f^j_1$.

The conditional distribution of $\!\left(\overline P^{\prime}(\bt ),\overline P^{\prime}(\bs )\right)$ given that $\bP (\bt )=0$ and $\bP (\bs )=0$ is a Gaussian distribution of dimension $2j\times (d-1)$ of mean zero and covariance 
$$
	\renewcommand{\arraystretch}{1.25}
	\left[
		\begin{array}{c;{1pt/2pt}c}
			B_{11} & B_{12} \\
			\hdashline[1pt/2pt]
			B_{12} & B_{22}
		\end{array}
	\right],
$$%
where 
\begin{align*}
 B_{11}&\phantom{:}=B_{22} :=\diag([\sigma^2,1,\ldots,1])_{j(d-1)\times j(d-1)}\\
 \intertext{and}
 B_{12}&:=\diag([\sigma^2\rho,\calD,\ldots,\calD])_{j(d-1)\times j(d-1)}
\end{align*}
 where these quantities are defined by (\ref{eq:deriv}) and (\ref{correlation}).

This result allows us to write the conditional distribution using the joint distribution of the  following  two $(j\times (d-1))$-dimensional vectors
\begin{align*}
 (\overline M_1,\ldots,\overline M_j)&:=(  M_{11},\ldots,M_{1d-1},  M_{21},\ldots,M_{2d-1},\cdots, M_{j1},\ldots,M_{jd-1})\\
 ( \overline W_1,\ldots,\overline W_j)&:=( W_{11},\ldots,W_{1d-1},W_{21},\ldots,W_{2d-1},\cdots, W_{j1},\ldots,W_{jd-1})
\end{align*}
where the $M_{\ell k}$ (resp$.$ $W_{\ell k}$) are  independent standard Gaussian random variables and 
\begin{eqnarray}
\label{covariances}
\bbE\!\left[M_{\ell_1k_1}W_{\ell_2k_2}\right]:=\rho \1_{\{\ell_1=\ell_2, k_1=k_2=1\}}+\calD \1_{\{\ell_1=\ell_2, k_1=k_2>1\}}.
\end{eqnarray}
Using these notations we readily obtain that (\ref{conditional}) is equal to
$$
	n^{j}\bbE\!\left[f^j_{\sigma}( \overline M_1,\ldots,\overline M_j)f^j_{\sigma}(\overline W_1,\ldots,\overline W_j)\right].
$$%
This implies that we can write
\begin{multline}
 \bbE\!\left[\calL_n^2(\calC_{\bP}(\0))\right]\\
 =n^{j}\int_{\bbS ^{d-1}\times\bbS ^{d-1}}\bbE\!\left[f^j_{\sigma}( \overline M_1,\ldots,\overline M_j)f^j_{\sigma}( \overline W_1,\ldots,\overline W_j)\right] {p}_{\bP(\bt ),\bP(\bs )}(\0,\0)\ud \bt  \ud \bs .
\end{multline}
$\displaystyle$

As in Section \ref{mehler section}, we denote  by  $(H_k(x))_{k \in \naturels}$ the Hermite polynomials.
Moreover, if $(Z,W)$ is a centered Gaussian vector with $\Var(Z):=\Var(W):=1$ and $\bbE(ZW):=\omega$, the following bi-dimensional Mehler's formula holds
$$
\bbE\!\left[H_k(Z)H_\ell (W)\right]=\delta_{k,\ell}\omega^k k!.
$$%

Moreover, since $(f^j_\gamma)^2$ is a polynomial this fact implies that
$$
f^j_\gamma \in \bL^2(\reels^{j\times (d-1)},\varphi_{j\times (d-1)}(\by')\ud \by'),
$$%
where $\varphi_{j\times (d-1)}$ is the Gaussian standard density on $\reels^{j\times (d-1)}$.
Using that the tensor Hermite polynomials form an orthogonal basis for this space we obtain
\begin{equation}
	\label{expansion1}
	f^j_\gamma(\by^{\prime})=\sum_{\bbeta}f^j_{\bbeta}(\gamma)\overline H_{\bbeta}(\by^{\prime}),
\end{equation}
where
\begin{align*}
\bbeta&:=(\beta_{11},\ldots,\beta_{1d-1},\ldots,\beta_{j1},\ldots,\beta_{jd-1}),\\ 
\overline H_{\bbeta}(\by^{\prime})&:=H_{\beta_{11}}(y^{\prime}_{11})\dots H_{\beta_{1d-1}}(y^{\prime}_{1d-1})\dots H_{\beta_{j1}}(y^{\prime}_{j1})\dots H_{\beta_{jd-1}}(y^{\prime}_{jd-1}), \intertext{and}
f^j_{\bbeta}(\gamma)&:=\frac1{\bbeta!}\int_{\reels^{j\times (d-1)}}f^j_\gamma(\by^{\prime})\overline H_{\bbeta}(\by^{\prime})\varphi_{j\times (d-1)}(\by^{\prime})\ud \by^{\prime},\\
\text{with~}\bbeta!&:= \prod\limits_{\ell=1}^{j}\prod\limits_{k=1}^{d-1} \beta_{\ell k}!.
\end{align*}

By Parseval's equality, we have
$$
\int_{\reels^{j\times (d-1)}}\abs{f^{j}_\gamma(\by^{\prime})}^2\varphi_{j \times (d-1)}(\by^{\prime})\ud \by^{\prime}=\sum_{\bbeta}|f^j_{\bbeta}(\gamma)|^2\bbeta!< \infty.
$$%
\spacebefore
\begin{rema}Note that for $j=d-1$, $f(\by^{\prime}):=f^{d-1}(\by^{\prime})=\abs{\det \by^{\prime}}$ is true.
Here $\by^{\prime}$ is the $(d-1)\times (d-1)$ matrix whose columns are the vectors $(y^{\prime}_{\ell1},\ldots, y^{\prime}_{\ell d-1})$ for $\ell=1,\ldots,d-1.$ Furthermore, applying the homogeneity of the determinant we have $f^{d-1}_{\gamma}(\by^{\prime})=\gamma f(\by^{\prime})$.
The lack of this property  when $j<d-1$, implies that we need a different approach to the one of \cite{MR3866880} for obtaining the asymptotic variance.
\hfill$\bullet$
\end{rema}

Let us calculate the conditional expectation (\ref{conditional}).

The expansion (\ref{expansion1}), the two-dimensional Mehler's formula and equality (\ref{covariances}) give
\begin{align*}
\MoveEqLeft[2]{n^{j}\bbE\!\left[f^j_{\sigma}( \overline M_1,\ldots,\overline M_j)f^j_{\sigma}( \overline W_1,\ldots,\overline W_j)\right]}\\
	&\phantom{:}=n^{j}\sum_{\bbeta}(f^j_{\bbeta}\!\left(\sigma(\prodsca{\bt}{\bs}))\right)^2\bbeta! (\calD(\prodsca{\bt}{\bs}))^{(\abs{\bbeta}-\sum_{\ell=1}^{j}\beta_{\ell1})}  
	(\rho(\prodsca{\bt }{\bs }))^{\sum_{\ell=1}^{j}\bbeta_{\ell1}}\\
	&:=n^{j}\calH(\prodsca{\bt }{\bs }),
\end{align*}
where
$\ds\abs{\bbeta}:= \sum\limits_{\ell=1}^{j}\sum\limits_{k=1}^{d-1} \beta_{\ell k}$.\\

Some conclusions follow from this last expression.
The first one is that the expectation in  (\ref{conditional}) depends on $\bt ,\bs $ through  $\prodsca{\bt}{\bs}$ only.
Then using (\ref{eq:intS}) we obtain
\begin{align*} 
\MoveEqLeft{\bbE[(\calL_n(\calC_{\bP}(\0)))^2]=\kappa_{d}\kappa_{d-1}\int_0^\pi\calH(\cos(\theta))\frac{n^{j}}{(2\pi)^j}\frac1{(1-\cos^{2n}(\theta))^{\frac j2}}\sin^{d-2}(\theta) \ud \theta}\\
	&\kern-12pt=\kappa_{d}\kappa_{d-1}\frac{n^{j-\frac12}}{(2\pi)^j}\int_0^{\sqrt n\pi}\calH\!\left(\cos(z/{\sqrt n})\right)\frac1{(1-\cos^{2n}(z/{\sqrt n}))^{\frac j2}}\sin^{d-2}(z/{\sqrt n}) \ud z.
\end{align*}
The second stems from the form of the integrand.
Indeed, we know that
\begin{align*}
\MoveEqLeft{\sum_{\bbeta}\abs{f^j_{\bbeta}(\gamma)}^2\bbeta!
	= \int_{\reels^{j\times (d-1)}}\abs{f^j_{\gamma}(\bx )}^2\varphi_{j\times (d-1)}(\bx )\ud \bx }\\
	&\le (d-1)^j(\gamma\vee 1)^{2j}\bbE\!\left[\sup_{1\le\ell\le j, 1\le k \le d-1}\abs{\frac{X^{\prime}_{\ell k}(\be_0)}{\sqrt n}}^{2j}\right].
 \end{align*}
Since $\sigma^2\cos((z/{\sqrt n}))\le 1$ and $\rho$ is a correlation (and $\abs{\rho} \le 1$), then
\begin{equation}
\sum_{\bbeta}\abs{f^j_{\bbeta}(\sigma(\cos(z/{\sqrt n})))}^2\bbeta! \le \bC .\label{norma2}
\end{equation}
Therefore, since $d -1>j$, we can interchange the series with the integral obtaining
\begin{multline}
\bbE[(\calL_n(\calC_{\bP}(\0)))^2]\\
=\kappa_{d}\kappa_{d-1}\frac{n^{j-\frac12}}{(2\pi)^j}
	\sum_{\bbeta}
			{\int_0^{\sqrt n\pi}\calH_{\bbeta}(\cos(z/{\sqrt n}))\frac{\sin^{d-2}(z/{\sqrt n})}{(1-\cos^{2n}(z/{\sqrt n}))^{\frac j2}} \ud z},\label{serie}
\end{multline}
where
\begin{multline*}
\calH_{\bbeta}(\cos(z/{\sqrt n}))
	:=\abs{f^j_{\bbeta}(\sigma(\cos(z/{\sqrt n})))}^2
	\bbeta!\!\left(\calD(\cos(z/{\sqrt n}))\right)^{(\abs{\bbeta}-\sum_{\ell=1}^{j}\beta_{\ell1})}\\
	\times(\rho(\cos(z/{\sqrt n})))^{\sum_{\ell=1}^{j}\bbeta_{\ell1} }. 
\end{multline*}
First of all, it should be noted that
\begin{align*}
\frac{\left(\bbE[\calL_n(\calC_{\bP}(\0))]\right)^2}{n^{j-{(d-1)/2}}}
	&=\frac{n^{{(d-1)/2}}}{(2\pi)^{j}}\kappa^2_d \abs{f^{j}_{\0}(1)}^2\\
	&=\frac{n^{{(d-2)/2}}}{(2\pi)^{j}} \kappa_d \kappa_{d-1} \abs{f^{j}_{\0}(1)}^2 \int_0^{\sqrt n\pi}\sin^{d-2}(z/{\sqrt n}) \ud z.
\end{align*}
Thus, using the above results and normalizing we obtain
\begin{align}
\MoveEqLeft[0.125]{\frac{\Var[\calL(\calC_{\bP}(\0))]}{n^{j-{(d-1)/2}}}
=\kappa_{d}\kappa_{d-1}\int_0^{\sqrt n\pi}\frac{n^{(d-2)/2}}{(2\pi)^j}}
	\!\left(\sum_{\abs{\bbeta}\ge1}\calH_{\bbeta}(\cos({z}/{\sqrt n}))\right)\nonumber\\
	& \specialpos{\hfill\times\frac1{(1-\cos^{2n}(z/{\sqrt n}))^{{j}/2}}\sin^{d-2}(z/{\sqrt n}) \ud z}\nonumber\\
	&+\kappa_{d}\kappa_{d-1}\int_0^{\sqrt n\pi}\frac{n^{{(d-2)/2}}}{(2\pi)^j}\!\left(\abs{f^j_{\0}(\sigma(\cos({z}/{\sqrt n})))}^2-\abs{f^{j}_{\0}(1)}^2\right)\label{sum}\\
	&\specialpos{\hfill\times\frac{\sin^{d-2}(z/{\sqrt n})}{(1-\cos^{2n}(z/{\sqrt n}))^{{j}/2}} \ud z}\nonumber\\
&+\kappa_{d}\kappa_{d-1}\int_0^{\sqrt n\pi}\frac{n^{(d-2)/2}}{(2\pi)^j} \abs{f^{j}_{\0}(1)}^2
	\!\left(\frac1{(1-\cos^{2n}(z/{\sqrt n}))^{{j}/2}}-1\right)\nonumber \\
	&\specialpos{\hfill\times\sin^{d-2}(z/{\sqrt n}) \ud z.}\nonumber
\end{align}

To apply the dominated convergence theorem, we have to look for a uniform bound.
The following lemma gives us the solution.
The prove is postponed to the end of this section.
\spacebefore
\begin{lemm}\label{uniformboud} Each term appearing in the integrands of the sum (\ref{sum}) is bounded by a Lebesgue  integrable function on $\intrfo{0}{\infty}$ of the variable $z$.
\end{lemm}

So we can take the limit under the integral sign.
Moreover,  let us show that we can interchange limit with the sum sign.
We have
$$
G(\gamma):=\sum_{\bbeta}\abs{f_{\bbeta}^j(\gamma)}^2\beta!=\int_{\reels^{j\times(d-1)}}\abs{f^j_\gamma(x)}^2\varphi_{j\times(d-1)}(x)\ud x.
$$%
Since the function $\abs{f^j_\gamma(x)}^2$ is a polynomial and has $\gamma$ as a coefficient, we can prove by using the dominated convergence theorem applied to the formula on the right, that $G$ is a continuous function then if $\gamma_n\cvg{}\gamma$ we have $G(\gamma_n)\cvg{} G(\gamma)$.
Thus
$$
\lim_{n\to+\infty}\sum_{\bbeta}\abs{f_{\bbeta}^j(\gamma_n)}^2\bbeta!=\sum_{\bbeta}\abs{f_{\bbeta}^j(\gamma)}^2\bbeta!.
$$%
The exact expression of the limit is obtained using the results of Lemma \ref{l:limits}.

Let us define
\begin{align}
\MoveEqLeft[2]{\calJ_{\bbeta}(z):=\abs{f^j_{\bbeta}\left(\left(\frac{1-(1+z^2)\exp(-z^2)}{1-\exp(-z^2)}\right)^{\frac12}\right)}^2\bbeta!\left(\exp(-{z^2}/2)\right)^{\abs{\bbeta}}}\nonumber\\
&\specialpos{\hfill\times \left(\frac{(1-\exp(-z^2)-z^2)}{1-(1+z^2)\exp(-z^2)}\right)^{\sum\limits_{\ell=1}^j\beta_{\ell1}}.} \label{H}
\end{align}
It yields
\begin{align*}
\MoveEqLeft[3]{ \lim_{n\to +\infty}\Var\left[\frac{\calL_n(\calC_{\bP}(\0))}{n^{j/2-(d-1)/4}}\right]}\\
= {} & \frac{\kappa_{d}\kappa_{d-1}}{(2\pi)^j}
	\int_0^{\infty}\Bigg(\sum_{\abs{\bbeta}\ge1}\calJ_{\bbeta}(z)\Bigg)\frac1{(1-\exp(-z^2))^{j/2}} z^{d-2}\ud z\\
	 & + \frac{\kappa_{d}\kappa_{d-1}}{(2\pi)^j}\int_0^{\infty}\left[\abs{f^j_{\0}\!\left(\left\{\frac{1-(1+z^2)\exp(-z^2)}{1-\exp(-z^2)}\right\}^{\frac12}\right)}^2-\abs{f^j_{\0}(1)}^2\right]\\
	 &\specialpos{\hfill\times\frac{z^{d-2}}{(1-\exp(-z^2))^{j/2}} \ud z}\\
	 &+  \frac{\kappa_{d}\kappa_{d-1}}{(2\pi)^j} \abs{f^j_{\0}(1)}^2\int_0^{\infty}\left[\frac1{(1-\exp(-z^2))^{j/2}}-1\right]z^{d-2}\ud z.
\end{align*}

Finally, we have proved the following theorem.
\spacebefore\begin{theo}
\begin{multline}
\lim_{n\to+\infty}\Var\left[\frac{\calL_n(\calC_{\bP}(\0))}{n^{j/2-(d-1)/4}}\right]\\
=\frac{\kappa_{d}\kappa_{d-1}}{(2\pi)^j}\int_0^{\infty}\left[\sum_{\abs{\bbeta}\ge0}\calJ_{\bbeta}(z) \frac1{(1-\exp(-z^2))^{j/2}}
-\abs{f^j_{\0}(1)}^2\right]z^{d-2}\ud z.
\end{multline}
\end{theo}
This is the same result that Letendre \cite{MR3917219} obtained before using another method.

Finally , we give the proof of Lemma \ref{uniformboud}.
\spacebefore
\begin{rema}\label{simetria}
The symmetrization  argument used in step 3 of Section 3.2 of \cite{MR3866880} gives that the integral over $\intrff{\sqrt n\frac{\pi}2}{\sqrt n\pi}$ of each term in the expansion (\ref{serie}) is equal to the integral of the same term on $\intrff{0}{\sqrt n\frac{\pi}2}$ except for a multiplication by $(-1)^{(n-1)\abs{\bbeta}}$.

In this form, the uniform bound obtained for applying the dominated convergence theorem in the latter interval is also used for the former.
\hfill$\bullet$
\end{rema}

Consider the first term of the sum (\ref{sum}).

Using the lemma \ref{l:bounds}, it turns out that there exists $n_0$ such that for $z/{\sqrt n}<\frac\pi2$ with $n>n_0$ and for $z \ge z_0$, which implies
$$
\abs{\rho}\le  \bC (1+z^2)\exp(-\alpha z^2)\mbox{ and } \calD\le \exp(-\alpha z^2).
$$%
According to  Remark \ref{simetria}, it is enough to work in the interval $\intrff{0}{\sqrt n\frac\pi2}$.
In that way we get, for $z \ge z_0$
\begin{align*}
\MoveEqLeft[1]{\left|{\sum_{|\bbeta|\ge1}\abs{f^j_{\bbeta}\!\left(\sigma(\cos(z/{\sqrt n}))\right)}^2
	\bbeta! (\calD(\cos(z/{\sqrt n})))^{(\abs{\bbeta}-\sum_{\ell=1}^{j}\beta_{\ell1})}}\right.}\\
&\specialpos{\hfill\left.\times{\left(\rho(\cos(z/{\sqrt n}))\right)^{\sum_{\ell=1}^{j}\beta_{\ell1}}}\rule{0pt}{24pt}\right|}\\
&\le \bC  \left(\sum_{\abs{\bbeta}\ge1} \abs{f^j_{\bbeta}(\sigma(\cos(z/{\sqrt n})))}^2\bbeta!\right)(1+z^2)\exp(-\alpha z^2)\\
&\le \bC  (1+z^2)\exp(-\alpha z^2).
\end{align*}
Above, we  used (\ref{norma2}).

It remains to consider the integral over the interval $\intrff{0}{z_0}$.
But now the integrand can be bounded on $\intrff{0}{z_0}$  by
$$
 \bC   \frac1{(1-\cos^{2n}(z/{\sqrt n}))^{\frac j2}}n^{(d-2)/2}\sin^{d-2}(z/{\sqrt n})
	 \le \bC  \frac1{(1-\exp(-2\alpha z^2))^{\frac j2}} z^{d-2},
$$%
and the function on the right-hand side is integrable whenever $j<d-1$. 
In this way, applying the dominated convergence theorem, we obtain
\begin{multline*}
	\lim_{n\to+\infty}\kappa_{d}\kappa_{d-1}\frac{n^{(d-2)/2}}{(2\pi)^j}
	\int_0^{\sqrt n\pi}\!\left(\sum_{\abs{\bbeta}\ge1}\calH_{\bbeta}(\cos({z}/{\sqrt n}))\right)\\
 	\specialpos{\hfill\times\frac{\sin^{d-2}(z/{\sqrt n})}{(1-\cos^{2n}(z/{\sqrt n}))^{j/2}} \ud z}\\
	=\kappa_{d}\kappa_{d-1}\frac{1}{(2\pi)^j}\int_0^{\infty}\!\left(\sum_{\abs{\bbeta}\ge1}\calJ_{\bbeta}(z)\right)\frac1{(1-\exp(-z^2))^{j/2}} z^{d-2}\ud z.
\end{multline*}
The exact expression of the function $\calJ_{\bbeta} (z)$ (see (\ref{H})) is obtained  using Lemma \ref{l:limits}.

To complete the proof of the lemma, we need to consider the remaining terms in (\ref{sum}).
First, let us consider
\begin{multline*}
\calI_n:=\kappa_{d}\kappa_{d-1}\frac{n^{(d-2)/2}}{(2\pi)^j}
	\int_0^{\sqrt n\pi}\left[\abs{f^j_{\0}(\sigma(\cos({z}/{\sqrt n})))}^2-\abs{f^j_{\0}(1)}^2\right]\\
\times\frac1{(1-\cos^{2n}(z/{\sqrt n}))^{j/2}}\sin^{d-2}(z/{\sqrt n}) \ud z. 
\end{multline*}
But
\begin{align}
\MoveEqLeft[1]{\abs{\abs{f^j_{\0}(\sigma(\cos({z}/{\sqrt n})))}^2-\abs{f^j_{\0}(1)}^2}}\nonumber \\
	&=\left|\!\left(\int_{\reels^{j\times (d-1)}}f^j_{\sigma}(\bx )\varphi_{j\times (d-1)}(\bx )\ud \bx \right)^2\right.\nonumber\\
	&\specialpos{\hfill\left.-\!\left(\int_{\reels^{j\times (d-1)}}f^j(\bx )\varphi_{j \times (d-1)}(\bx )\ud \bx \right)^2\right|}\nonumber\\
	&\le  \bC  \abs{\int_{\reels^{j\times (d-1)}}f^j_{\sigma}(\bx )\varphi_{j\times (d-1)}(\bx )\ud \bx -\int_{\reels^{j\times (d-1)}}f^j(\bx )\varphi_{j\times (d-1)}(\bx )\ud \bx }\nonumber\\
	&\le \bC   \1_{\{z\le z_0\}}
	+\bC  \left| \int_{\reels^{j \times (d-1)}}f^j_{\sigma}(\bx )\varphi_{j\times (d-1)}(\bx )\ud \bx \right. \nonumber\\[-8pt]
	&\specialpos{\hfill\left.-\int_{\reels^{j \times (d-1)}}f^j(\bx )\varphi_{j\times (d-1)}(\bx )\ud \bx  \right|\1_{\{z>z_0\}}}.\label{boundGaussian}
\end{align}

To find an upper bound for the second term of the last sum, we use the following general result.

Let  $G:\reels^M\to\reels$ be a function such that $\abs{G(\bx )}\le \bC \normp[]{\bx }^m$, $m \in \naturels^{\star}$, where without loss of generality we take $\normp[]{\bx }:=\sup_i\abs{x_i}$.
Consider two sequences of positive numbers $\{\sigma_{i}\}_{i=1}^M$ and $\{\gamma_{i}\}_{i=1}^M$ and denote $\normp[]{\boldsymbol\sigma}:=\sup_i\sigma_i$ and $\normp[]{\boldsymbol\gamma}:=\sup_i\gamma_i$. Also consider $\bW:=(W_1,\ldots,W_M)\stackrel{\calL}={\calN}(0,I_M)$. 
We have
\begin{multline*}
 \abs{\bbE\!\left[G(\sigma_1W_1,\sigma_2W_2,\ldots,\sigma_MW_M)-G(\gamma_1W_1,\gamma_2W_2,\ldots,\gamma_MW_M)\right]}\\
 \le \bC  \sqrt{2M}(\normp[]{\boldsymbol\sigma}\vee \normp[]{\boldsymbol\gamma})\sum_{i=1}^M\frac{\abs{\sigma_i-\gamma_i}}{\sigma_i\wedge\gamma_i}(\bbE[\normp[]{W}^{2m}])^\frac12,
\end{multline*}
where the symbols $\vee$ and $\wedge$ denote the supremum and infimum respectively.

Applying this result to the second term of the last sum in (\ref{boundGaussian}) and using that $\sigma(\cos({z}/{\sqrt n}))\le1$ and $\abs{1-\sigma}\le \bC  \exp{(-2\alpha z^2)}$ for $z>z_0$, we obtain
\begin{multline*}
 \abs{\abs{f^j_{\0}(\sigma(\cos({z}{\sqrt n})))}^2-\abs{f^j_{\0}(1)}^2}
\le \bC  \1_{\{z\le z_0\}}+ \bC  \exp{(-2\alpha z^2)}\1_{\{z>z_0\}}.
\end{multline*}

Now, the dominated convergence theorem can be applied and again by the results of Lemma \ref{l:limits} we get
\begin{align*}
 \lim_{n\to+\infty}\calI_n
 	=\frac{\kappa_{d}\kappa_{d-1}}{(2\pi)^j}
		\int_0^{\infty}\left[\abs{f^j_{\0}\!\left(\left\{\frac{1-(1+z^2)\exp(-z^2)}{1-\exp(-z^2)}\right\}^{\frac12}\right)}^2-\abs{f^j_{\0}(1)}^2\right]\\
	\times\frac1{(1-\exp(-z^2))^{j/2}} z^{d-2}\ud z.
\end{align*}

Finally we will consider the last term $\calJ_n$.
\begin{align}
	\calJ_n
	:= \kappa_{d}\kappa_{d-1}\frac{n^{(d-2)/2}}{(2\pi)^j}
		\int_0^{\sqrt n\pi}\abs{f^j_{\0}(1)}^2
			\left[\frac1{(1-\cos^{2n}(z/{\sqrt n}))^{j/2}}-1\right]\nonumber\\
			\times\sin^{d-2}(z/{\sqrt n}) \ud z.\label{integrand}
\end{align}
By setting $R:=\cos^n(z/{\sqrt n})$, we have
\begin{equation}
\label{bound}
\abs{\frac1{(1-\cos^{2n}(z/{\sqrt n}))^{j/2}}-1} \le \bC  \frac{R^2}{(1-R^2)^{\frac{j}{2}}},
\end{equation}
with $\abs{R} \le \exp(-\alpha z^2) \wedge 1$ if $z < \frac{\pi}{2} \sqrt{n}$.\\
Thus the integrand in (\ref{integrand}) can be bounded with the bound computed in (\ref{bound}), depending on $z<z_0$ or $z\ge z_0$.
That is, on the interval close to zero or on the other.
In the first case, we use
$$
\abs{n^{(d-2)/2}
	\left[\frac1{(1-\cos^{2n}(z/{\sqrt n}))^{j/2}}-1\right]
	\sin^{d-2}(z/{\sqrt n})}
	\le \bC  \frac{z^{d-2}}{(1-\exp{(-2\alpha z^2))^{j/2}}},
$$%
the last function above is an integrable function since $j < d-1$.
On the large interval, we get
\begin{multline*}
 \abs{n^{(d-2)/2}\left[\frac1{(1-\cos^{2n}(z/{\sqrt n}))^{j/2}}-1\right]\sin^{d-2}(z/{\sqrt n})}\\
 	\le \bC  \exp(-2\alpha z^2)\frac{z^{d-2}}{(1-\exp(-2\alpha z_0^2))^{j/2}},
\end{multline*}
the last function above being again an integrable function.
Since these two bounds allow applying the dominated convergence theorem, it turns out that
$$
\lim_{n\to+\infty}\calJ_n
	= \frac{\kappa_{d}\kappa_{d-1}}{(2\pi)^j} \abs{f^j_{\0}(1)}^2
		\int_0^{\infty}\left[\frac1{(1-\exp(-z^2))^{j/2}}-1\right]z^{d-2}\ud z.
$$%
That completes the proof of Lemma \ref{uniformboud}.
\section{Local time and length of curves of level set}
Let $X: \Omega \times \reels^d \to \reels$, $d \ge 2$, be a stationary centered continuous Gaussian process and $F$ be its spectral measure assumed not to be concentrated in a Borel set of zero Lebesgue measure.
Let $T$ be a bounded open set of $\reels^d$.\\
We will consider regularizations of the trajectories of $X$ by means of convolutions of the form $X_\ve(t)= \Psi_\ve \sstar X(t)$, $\ve >0$, where $ \Psi_\ve$ is an approximation of unity when $\ve$ tends to zero and satisfies some regularity conditions.\\
The aim of this section is to approximate the local time of $X$ over $T$ at a given level $u$, say $L_{X}(u, T)$, by the length of curves of the $u$ level set of process $X_\ve$, say $\sigma_{d-1}(C_{T, X_\ve}(u))$.\\

Some of the results that will be presented in Sections  \ref{LocalTime} and \ref{ApproxLenghLevel}  have been partially discussed in \cite{MR1222362}.
We can cite some references of papers that have dealt on the same subject in the case where $d=1$.
More precisely if $N_{X_\ve}(u)$ represents the number of crossings of the level $u$ by the process $X_\ve$ on the interval $\intrff{0}{1}$, Wschebor \cite{MR871689}  has shown that, in the case of Brownian motion,  the appropriately normalized random variable $N_{X_\ve}(0)$ tends to $L_X(0, \intrff{0}{1})$, as $\ve \to 0$, in $L^p(\Omega)$  for any $p \ge 1$, and a similar result applies to multiparameter Brownian motion.
Aza\"{\i}s and Florens \cite{MR899448}  have extended this result to a class of stationary Gaussian processes.

We start with some considerations on the local time of $X$ on $T$.
\subsection{Local time}
\label{LocalTime}
We define the local time as the density of the occupation measure of $X$.
This definition was used for example by Berman (\cite{MR239652}).\\
More precisely, if $A$ denotes a Borel set in $\reels$, we define the occupation measure $\mu_T(\omega, A)$ of $X$ by 
\begin{eqnarray*}
\mu_T(\omega, A):= \sigma_{d}\{t \in T, X(\omega, t) \in A\}.
\end{eqnarray*}
The local time $L_{X}(u, T)$, when it exists, is any function satisfying
\begin{eqnarray*}
\mu_T(\omega, A)= \int_{A} L_{X}(\omega, u, T) \ud u,
\end{eqnarray*}
for almost every $\omega$.\\
Another formulation is that the local time, when it exists, is the Radon-Nikodym derivative\begin{eqnarray*}
L_{X}(u, T):= \frac{\ud  \mu_T}{\ud  \sigma_1}(u).
\end{eqnarray*}
By Proposition 1 of \cite{MR1222362}, the local time of $X$, denoted by $L_X(\cdot, T)$, exists $P$-almost surely and is $L^2(\lambda_1 \times P)$.
One can also refer to Theorem 22.1 of \cite{MR556414}.
One reason why this result is true is that 
$$
\int_{T \times T}
{(r^2(0)-r^2(s-t))^{-{1}/{2}}} \ud s \ud t < \infty,
$$%
where $r$ stands for the covariance of process $X$.
The finiteness of the last integral comes from the following lemma whose proof is given in Lemma 1 of  \cite{MR1222362}.
\spacebefore
\begin{lemm}
\label{minoration r}
There exists $\bC>0$, such that for all $(s, t) \in \overline{T}\times \overline{T}, r^2(0)-r^2(s-t) \ge \bC  \normp[d]{s-t}^2$, where $\overline{T}$ is the closure of the set $T$.
\end{lemm}
We now state the identity of the second moment of local time.
\spacebefore
\begin{prop} 
\label{expectation local time}
For all any level $u \in \reels$, we have
$$
	\bbE\!\left[L_X^2(u, T)\right]= \int_{T \times T} {p}_{X(s), X(t)}(u, u) \ud s \ud t,
$$%
where ${p}_{X(s), X(t)}(\bm{\cdot}, \bm{\cdot})$ is the density of vector $(X(s), X(t))$.
\end{prop}
\spacebefore
\begin{proofarg}{Proof of Proposition \ref{expectation local time}}
The proof is given in \cite{MR556414} (see equality (25.5)).
However, we give here a sketch of the proof which we will use later.\\
Proposition~\ref{approximation temps local} ensues from Proposition 1 of \cite{MR1222362}.
\spacebefore
\begin{prop}
\label{approximation temps local}
 There exists a version of $L(\cdot, T)$ such that $\forall u \in \reels$, $\eta_{\delta}(u):= \frac{1}{2 \delta} \int_{T} \1_{\{t, \abs{X(t)-u} < \delta \}} \ud t$, $\delta >0$, converges in $L^{2}(\Omega)$ to $L_X(u, T)$.
 \end{prop}
Let $u \in \reels$.
By computing $\bbE[\eta_{\delta}^2(u)]$, we obtain:
 $$\bbE\!\left[\eta_{\delta}^2(u)\right]= \!\left(\frac{1}{2\delta}\right)^2 \int_{u-\delta}^{u+\delta} \int_{u-\delta}^{u+\delta} \!\left({\int_{T \times T} 
 {p}_{X(s), X(t)}(x, y) \ud s \ud t}\right) \ud x \ud y.
 $$
 Using Lemma \ref{minoration r}, we easily obtain that the function
 $$
 (x, y) \mapsto  \int_{T \times T} {p}_{X(s), X(t)}(x, y) \ud s \ud t
 $$%
 is continuous at $(x, y)$ and then using Proposition \ref{approximation temps local}, we obtain
$$
\bbE[\eta_{\delta}^2(u)] \cvg[\delta\to 0]{} \int_{T \times T} {p}_{X(s), X(t)}(u, u) \ud s \ud t= \bbE[L_X^2(u, T)].
$$%
Another interesting way to prove Proposition \ref{expectation local time} is to apply Theorem 1 of \cite{MR1386211}.
Indeed, in this theorem, the authors expand the local time $L_X(u, T)$ in terms of Hermite polynomials.
More precisely, if $(H_k(x))_{k \in \naturels}$ is the Hermite orthogonal  basis of $L^{2}(\reels, \varphi(x) \ud x)$, \ie
$$
H_k(x):= (-1)^k \varphi^{-1}(x)  \ds \frac{\ud ^k}{\ud  x^k}(\varphi(x)),
$$%
where $\varphi$ denotes the standard Gaussian density on $\reels$, they show the following theorem.
\spacebefore\begin{theo}
\label{local time DOLE}
Let $u \in \reels$ a fixed level.
The following expansion of the local time $L_X(u, T)$ holds in $L^{2}(\Omega)$:
$$
L_X(u, T)= \frac{1}{\sqrt{r(0)}} \varphi\!\left(\frac{u}{\sqrt{r(0)}}\right) \sum_{k \in \naturels} \frac{1}{k!} H_{k}\!\left(\frac{u}{\sqrt{r(0)}}\right) \int_{T} H_k\!\left({X(s)}/{\sqrt{r(0)}}\right) \ud s.
$$%
\end{theo}
Using this last expansion and Mehler's formula (see \cite{MR716933}), we obtain
\begin{multline*}
\bbE[L_X^2(u, T)]\\= \frac{1}{r(0)} \varphi^2\!\left(\frac{u}{\sqrt{r(0)}}\right)  \sum_{q \in \naturels} \frac{1}{q!} H_{q}^2\!\left(\frac{u}{\sqrt{r(0)}}\right) \int_{T \times T} \!\left(\frac{r(s-t)}{r(0)}\right)^q \ud s \ud t. 
\end{multline*}
The following equality can be found in \cite[Eq$.$ (39)]{MR1502747}.
For $0 \le \abs{w} <1$ and $x \in \reels$, we have
$$ 
\sum_{q \in \naturels} \frac{1}{q!} H_{q}^2(x)  w^q= \frac{1}{\sqrt{1-w^2}} \exp\!\left(\frac{x^2 w}{1+w}\right).
$$%
Thus, applying this last equality to $x:={u}/{\sqrt{r(0)}}$ and to $w:= {r(s-t)}/{r(0)}$, we obtain
\begin{align*}
	\bbE[L_X^2(u, T)]
		&= \frac{1}{2\pi} \int_{T \times T} \frac{1}{\sqrt{r^2(0)-r^2(s-t)}} \exp\!\left(-\frac{u^2}{r(0)+r(s-t)}\right) \ud s \ud t\\
		&=\int_{T \times T}  {p}_{X(s), X(t)}(u, u) \ud s \ud t.
\end{align*}
This completes the proof of Proposition \ref{expectation local time}.
 \end{proofarg}%
Let $f: T \to \reels$ be a bounded function.
More generally that before, we are interested in a canonical renormalization of the ordinary local time, given by, $\forall u \in \reels$
$$
	L_X^f(u, T):=\int_{T} f(t) L_X(u, \ud t).
$$%
Geman and Horowitz proved in \cite[Theorem 6.4]{MR556414} the following identity.
\spacebefore\begin{theo}
\label{GH}
For any Borel function $h(t, x) \ge 0$ on $T \times \reels$,
$$
	\int_T h(t, X(t)) \ud t = \int_{\reels} \int_{T} h(t, x) L_X(x, \ud t) \ud x.
$$%
\end{theo}
\spacebefore
\begin{rema}
\label{Borel}
This equality is still valid for any Borel function $h(t, x)$, which is not necessarily positive, provided that one of the two integrals is finite.
This is for example the case when $h$ is bounded.
\hfill$\bullet$
\end{rema}
\spacebefore
\begin{rema}
If we take $h(t, x):=\1_A(x)$, where $A$ is a Borel set of $\reels$, we find the usual definition of local time.
If we take $h(t, x):= g(x)$, where $g$ is a bounded Borel function, we get the known occupation formula
$$
\int_T g(X(t)) \ud t = \int_{\reels}  g(x) L_X(x, T) \ud x.
$$%
\hfill$\bullet$
\end{rema}
Let us define the measure $\mu_f$ as
$$
\mu_f(A):= \int_A \int_T f(t) L_X(x, \ud t) \ud x,
$$%
where $A$ is any Borel set of $\reels$.
Since $f$ is bounded, $\mu_f(A)\le \bC  \sigma_d(T) < \infty$.
The measure $\mu$ is a finite and its Fourier transform is
$$
\widehat{\mu}_f(z):= \int_{\reels} e^{izx} \ud \mu_f(x)=\int_{\reels} e^{izx} \int_T f(t) L_X(x, \ud t) \ud x,\quad z \in \reels.
$$
By Theorem \ref{GH} and Remark \ref{Borel} applied to $h(t, x):= e^{izx} f(t)$, we have
$\widehat{\mu}_f(z)= \int_{T} e^{izX(t)} f(t) \ud t$, so by applying Lemma \ref{minoration r} and since function $f$ is a bounded one, we obtain
\begin{align*}
	\bbE\!\left[\abs{\widehat{\mu}_f(z)}^2\right]
		&=\int_{T \times T} f(t) f(s) \bbE[e^{iz(X(s)-X(t))}] \ud s \ud t\\
		&= \int_{T \times T} f(t) f(s)  e^{-z^2(r(0)-r(s-t))} \ud s \ud t\\
		&\le  \int_{T \times T} \abs{f(s)} \abs{f(t)} e^{-\bC z^2 \normp[d]{s-t}^2} \ud s \ud t\\
		&\le  \frac{\bC }{(1+\abs{z})^d}.
\end{align*}
We deduce from these inequalities that
$$
	\bbE\!\left[\int_{-\infty}^{+\infty}\abs{\widehat{\mu}_f(z)}^2 \ud z\right] \le \bC ,
$$%
and we conclude that, with probability one,
$$
\int_{-\infty}^{+\infty}\abs{\widehat{\mu}_f(z)}^2 \ud z < \infty.
$$%
Thus by Plancherel's theorem we proved that
$$
\int_{-\infty}^{+\infty} \bbE\!\left[(L_X^f(x, T))^2\right] \ud x < \infty,
$$%
and then for almost all $x \in \reels$, 
$$
 \bbE\!\left[(L_X^f(x, T))^2\right]< \infty.
$$%
Then, we aim at giving a generalization of Propositions \ref{expectation local time} and \ref{approximation temps local} for the renormalized local time $L_X^f(\cdot, T)$.

Let us state the following two propositions.
\spacebefore
\begin{prop}
\label{approximation temps local Bis}
For any level $u \in \reels$, the random variable  $\eta_{\delta}^f(u)$ defined by 
$$
	\eta_{\delta}^f(u):= \frac{1}{2 \delta} \int_{T} f(t) \1_{\{t, \abs{X(t)-u} < \delta \}} \ud t,
$$%
$\delta >0$, converges in $L^{2}(\Omega)$ towards $L_X^f(u, T)$.%
 \end{prop}
 
Since function $f$ is bounded and $(x, y) \mapsto  \int_{T \times T} {p}_{X(s), X(t)}(x, y) \ud s \ud t$ is continuous at $(x, y)$, as an immediate consequence we deduce as above the following proposition,
 \spacebefore
\begin{prop} 
\label{expectation local time Bis}
For any level $u \in \reels$, we have
$$
\bbE\!\left[(L_X^f(u, T))^2\right]= \int_{T \times T}  f(s) f(t)  {p}_{X(s), X(t)}(u, u) \ud s \ud t.
$$%
\end{prop}
\spacebefore
\begin{proofarg}{Proof of Proposition \ref{approximation temps local Bis}}
To prove this proposition, it is sufficient to prove the following lemma.
\spacebefore
\begin{lemm}
\label{modulus of continuity}
For any $x, h \in \reels$, we have 
$$
\bbE\!\left[\left(L_X^f(x+h, T)- L_X^f(x, T)\right)^2\right] \le  \bC  \theta(h) \le \bC  \abs{h}, 
$$%
where $\theta(h):=  (h^2+\abs{h}^{d-1}) \1_{\{d\neq 3\}} + h^2 \ln({1}/{\abs{h}}) \1_{\{d = 3\}}$,
for $h$ sufficiently small  and where the constante $ \bC $ is locally bounded as function of $x$.
\end{lemm}
\spacebefore
\begin{rema}
This lemma highlights the fact that since $d \ge 3$,
the renormalized local time $L_X^f(\cdot, T)$  \as admits a continuous version as function of $x$.
\hfill$\bullet$
\end{rema}
Taking $h(t, x):= f(t) \1_A(x)$, with $A= \introo{u-\delta}{u+ \delta}$ in Theorem \ref{GH}, we obtain%
$$
\eta_{\delta}^f(u)=\frac{1}{2\delta} \int_{u-\delta}^{u+\delta} L_X^f(x, T) \ud x.
$$%
This implies that
\begin{align*}
\bbE\!\left[\left(\eta_{\delta}^f(u) - L_X^f(u, T)\right)^2\right]
	& =  \bbE\!\left[{
			\left(
				\frac{1}{2\delta} {\int_{u-\delta}^{u+\delta} \!\left({ L_X^f(x, T)- L_X^f(u, T)}\right) \ud x}
			\right)^2
			}\right]\\
	& \le  \frac{1}{2\delta} \int_{u-\delta}^{u+\delta} \bbE\!\left[\left(L_X^f(x, T)- L_X^f(u, T)\right)^2\right] \ud x\\
	& \le \bC  \delta \cvg[\delta\to0]{} 0.
\end{align*}
That will end the proof of Proposition \ref{approximation temps local Bis}.
Let us prove Lemma \ref{modulus of continuity}.
\end{proofarg}
\spacebefore
\begin{proofarg}{Proof of Lemma \ref{modulus of continuity}}
It is sufficient to prove the inequality of the lemma for $x, x+h$ not belonging to a set $N$, such that $\sigma_1(N)=0$.
Indeed, if it is possible to do so, then  we can define $L_X^f(x, T)$ for $x \in N$ as limit in $L^2(\Omega)$ of $L_X^f(y, T)$ as $y \to x$, with $y \notin N$ and this version of the canonical renormalization of the ordinary local time will satisfy the required inequality for all $x, x+h \in \reels$.\\
By Theorem \ref{GH} and Remark \ref{Borel} applied to $h(t, x):= 1$, we have
$$
\sigma_d(T) = \int_{\reels}   L_X(x, T) \ud x < \infty.
$$%
Now, since $f$ is bounded, we have  $\forall x \in \reels$, $\abs{L_X^f(x, T)} \le \bC  L_X(x, T)$.
Therefore, $L_X^f(\cdot, T) \in L^1(\reels)$.

Lebesgue's differentiation theorem ensures the existence of a $\sigma_1$-null measure set $N$ such that if $x \notin N$ then
$$
	\lim_{\delta \to 0} \eta_\delta^f(x)= \lim_{\delta \to 0} \frac{1}{2\delta} \int_{x-\delta} ^{x+\delta} L_X^f(y, T) \ud y =  L_X^f(x, T).
$$%
So, applying Fatou's lemma and the fact that function
 $$
 (x, y) \to  \int_{T \times T} f(s) f(t)  {p}_{X(s), X(t)}(x, y) \ud s \ud t
 $$
 is continuous at $(x, y)$, we get that for all $x, x+h \notin N$,  
\begin{align*}
	\MoveEqLeft{\bbE\!\left[\left(L_X^f(x+h, T)- L_X^f(x, T)\right)^2\right]}\\
	&= \bbE\!\left[\lim_{\delta \to 0}\!\left({\eta_\delta^f(x+h)- \eta_\delta^f(x)}\right)^2\right]\\
	&=\bbE\!\left[\lim_{\delta \to 0}\!\left(\frac{1}{2 \delta}\right)^2 \!\left({\int_{T} f(t) \!\left({ \1_{\{t:\abs{X(t)-(x+h)} < \delta \}} - \1_{\{t:\abs{X(t)-x} < \delta \}} }\right) \ud t }\right)^2\right] \\
	& \le \int_{T \times T} f(s) f(t) \left\{{ {p}_{X(s), X(t)}(x+h, x+h) }\right.\\
	 &	\specialpos{\hfill\left.{\qquad-2  {p}_{X(s), X(t)}(x+h, x) +  {p}_{X(s), X(t)}(x, x)}\right\} \ud s \ud t.}
\end{align*}
A simple calculation shows that the right-hand side of the last inequality is bounded by
\begin{multline*}
 	\bC  \int_{T \times T} \frac{1}{(r^2(0)-r^2(s-t))^{\frac{1}{2}}}\\
		\times\left[{h^2
			+ \left\{{1-\exp(-a h^2 (r^2(0)-r^2(s-t))^{-1})}\right \}
		}\right] \ud s \ud t ,
\end{multline*}
where $a > 0$ is a fixed number and $\bC $ is locally bounded as function of $x$.
Lemma \ref{minoration r} gives the required result.
\end{proofarg}
\spacebefore
\begin{rema}
\label{different levels}
It is interesting to note that Proposition \ref{approximation temps local Bis} implies, without particular effort, that for any level $u, v \in \reels$ and for all bounded functions $f, g: T \to \reels$, we have
$$
	\bbE\!\left[L_X^f(u, T)L_X^g(v, T)\right]= \int_{T \times T}  f(s) g(t) {p}_{X(s), X(t)}(u, v) \ud s \ud t.
$$%
\hfill$\bullet$
\end{rema}
To finish this section on the local time, we complete the statement of Theorem \ref{local time DOLE} due to Doukhan and Le\'on(\cite{MR1386211}) to the renormalized local time.
\spacebefore\begin{theo}
\label{local time DOLE Bis}
Let $u\in\reels$, a fixed level.
The following expansion of the local time $L_X^f(u, T)$ holds in $L^{2}(\Omega)$:
\begin{multline}
	L_X^f(u, T)= \frac{1}{\sqrt{r(0)}} \varphi\!\left(\frac{u}{\sqrt{r(0)}}\right)\\
	\times\sum_{k \in \naturels} \frac{1}{k!} H_{k}\!\left(\frac{u}{\sqrt{r(0)}}\right) \int_{T} f(s) H_k\!\left({X(s)}/{\sqrt{r(0)}}\right) \ud s.
	\label{expansion du temps local}
\end{multline}
\end{theo}
\spacebefore
\begin{proofarg}{Proof of Theorem \ref{local time DOLE Bis}}
Let us define the \rv $Z^f(u, T)$ as the right-hand side of (\ref{expansion du temps local}) that is
\begin{multline*}
 Z^f(u, T):=\frac{1}{\sqrt{r(0)}} \varphi\!\left(\frac{u}{\sqrt{r(0)}}\right) \\
 \times\sum_{k \in \naturels} \frac{1}{k!} H_{k}\!\left(\frac{u}{\sqrt{r(0)}}\right) \int_{T} f(s) H_k\!\left({X(s)}/{\sqrt{r(0)}}\right) \ud s.
\end{multline*}%
 Let us prove that $\bbE\!\left[\left(L_X^f(u, T)-Z^f(u, T)\right)^2\right]=0$.
 
On the one hand, by Proposition \ref{expectation local time Bis}, we already know that
$$
\bbE\!\left[(L_X^f(u, T))^2\right]= \int_{T \times T}  f(s) f(t) {p}_{X(s), X(t)}(u, u) \ud s \ud t.
$$%
On the other hand, using Mehler's formula, we obtain
\begin{multline*}
\bbE\left[\left(Z^f(u, T)\right)^2\right] 
= \frac{1}{r(0)} \varphi^2\!\left(\frac{u}{\sqrt{r(0)}}\right) \\
\times \sum_{q \in \naturels} \frac{1}{q!} H_{q}^2\!\left(\frac{u}{\sqrt{r(0)}}\right) \int_{T \times T} f(s) f(t) \!\left(\frac{r(s-t)}{r(0)}\right)^q \ud s \ud t.
\end{multline*}
And as we saw in the second way of proving Proposition \ref{expectation local time}, we use  \cite[Eq$.$~(39)]{MR1502747} to obtain
$$
	\bbE\!\left[\left(Z^f(u, T)\right)^2\right]=\bbE\!\left[(L_X^f(u, T))^2\right].
$$%
Let us now compute $\bbE\!\left[L_X^f(u, T) Z^f(u, T)\right]$.\\
First for fixed $k \in \naturels$, let us define function $h_k$ by $h_k(t, x):=f(t) H_k({x}/{\sqrt{r(0)}})$ in Theorem \ref{GH}.
We thus obtain the equality
$$
\int_T f(t) H_k\!\left({{X(t)}/{\sqrt{r(0)}}}\right) \ud t= \int_{\reels} L_X^f(x, T) H_k\!\left({x}/{\sqrt{r(0)}}\right) \ud x,
$$
and then
\begin{multline*}
\bbE\!\left[L_X^f(u, T) Z^f(u, T)\right]
= \frac{1}{\sqrt{r(0)}} \varphi\!\left(\frac{u}{\sqrt{r(0)}}\right) \sum_{k \in \naturels} \frac{1}{k!} H_{k}\!\left(\frac{u}{\sqrt{r(0)}}\right)\\
\times \int_{\reels} H_k\!\left({x}/{\sqrt{r(0)}}\right) \bbE\!\left[L_X^f(x, T)L_X^f(u, T)\right] \ud x.
\end{multline*}
Applying Remark \ref{different levels}, we get
\begin{align*}
	\MoveEqLeft[3]{\bbE\!\left[L_X^f(u, T) Z^f(u, T)\right]}\\
	 = {} & \frac{1}{\sqrt{r(0)}} \varphi^2\!\left(\frac{u}{\sqrt{r(0)}}\right) \int_{T \times T} f(s) f(t) \times \frac{1}{\sqrt{r^2(0)-r^2(s-t)}}  \\
	 \times {} &  \sum_{k \in \naturels} \frac{1}{k!} H_{k}\!\left(\frac{u}{\sqrt{r(0)}}\right) \frac{1}{\sqrt{2 \pi}}\int_{\reels} 
	\quad H_k\!\left({x}/{\sqrt{r(0)}}\right)  \\
	&\specialpos{\hfill\times \exp\!\left(-\frac{\left(x\sqrt{r(0)} - {u}\, r(s-t)/{\sqrt{r(0)}}\right)^2}{2(r^2(0)-r^2(s-t))}\right) \ud x \ud s \ud t.}
\end{align*}
The following formula can be found in \cite[(8), page 804]{MR3307944}.
For $0 \le \alpha <1$, $k \in \naturels$ and $z \in \reels$,
$$
\frac{1}{\sqrt{2\pi}}\int_{-\infty}^{\infty} \exp\!\left(-\tfrac{1}{2}(x-z)^2\right) H_k(\alpha x) \ud x= (1-\alpha^2)^{{k}/{2}} H_k\!\left(\frac{\alpha z}{\sqrt{1-\alpha^2}}\right).
$$%
We make the change of variable $x=\beta y$ for $\beta \neq 0$.
For $0 \le \alpha <1$, $\beta \neq 0$, $k \in \naturels$ and $z \in \reels$, we obtain
$$
\frac{1}{\sqrt{2\pi}}\int_{-\infty}^{\infty} \exp(-\tfrac{1}{2}(\beta x-z)^2) H_k(\alpha \beta x) \ud x= \frac{1}{\beta} (1-\alpha^2)^{{k}/{2}} H_k\!\left(\frac{\alpha z}{\sqrt{1-\alpha^2}}\right).
$$%
We apply this last formula to 
\begin{align*}
 \beta&:=\ds \frac{\sqrt{r(0)}}{\sqrt{r^2(0)-r^2(s-t)}},\\ 
z&:= \ds \frac{ur(s-t)}{\sqrt{r(0)}\sqrt{r^2(0)-r^2(s-t)}}\\ \shortintertext{and to }
\alpha&:= \ds \frac{\sqrt{r^2(0)-r^2(s-t)}}{r(0)}.
\end{align*}
We obtain
\begin{multline*}
\bbE\!\left[L_X^f(u, T) Z^f(u, T)\right]=
\frac{1}{r(0)} \varphi^2\!\left(\frac{u}{\sqrt{r(0)}}\right) \int_{T \times T} f(s) f(t)   \\
\times \sum_{k \in \naturels} \frac{1}{k!} H_{k}\!\left(\frac{u}{\sqrt{r(0)}}\right)
 \abs{\frac{r(s-t)}{r(0)}}^k H_k\!\left(\frac{u r(s-t)}{\sqrt{r(0)}\abs{r(s-t)}}\right) \ud s \ud t.
\end{multline*}

We consider the following equality in \cite[(39)]{MR1502747}).
For $0 \le \abs{w} < 1$ and $x, y \in \reels$,
$$ 
\sum_{q \in \naturels} \frac{1}{q!} H_{q}(x)  H_q(y) w^q= \frac{1}{\sqrt{1-w^2}} \exp\!\left(\frac{-w^2(x^2+y^2)+2wxy}{2(1-w^2)}\right).
$$%
We apply this last inequality by choosing 
\begin{align*}
	x& :=	 \frac{u}{\sqrt{r(0)}},\\
	y& :=	 \frac{ur(s-t)}{\sqrt{r(0)}\abs{r(s-t)}}\\ \shortintertext{and}
	w& := \abs{\frac{r(s-t)}{r(0)}}, 
\end{align*}
getting
\begin{align*}
\MoveEqLeft[3]{\bbE\!\left[L_X^f(u, T) Z^f(u, T)\right]}\\
= {} &\varphi^2\!\left(\frac{u}{\sqrt{r(0)}}\right) \times \int_{T \times T} f(s) f(t) \times 
 \frac{1}{\sqrt{r^2(0)-r^2(s-t)}} \\
 &\times  \exp\!\left(\frac{u^2 r(s-t)}{r(0) (r(0)+r(s-t))}\right) \ud s \ud s\\
= {} & \bbE\!\left[\left(L_X^f(u, T)\right)^2\right],
\end{align*}
yielding Theorem \ref{local time DOLE Bis}.
\end{proofarg}%
Now we will approximate the renormalized local time $L_X^f(u, T)$ over $T$ at $u$, by an appropriate normalization of functionals of the $u$-level set of the smoothed process $X_\ve$.%
\subsection{Approximation to the length of curves of level set}
\label{ApproxLenghLevel}
Let $X: \Omega \times \reels^d \to \reels$, $d \ge 2$, be a continuous Gaussian isotropic process of zero mean and $F$ its spectral measure.
We assume that $F$ is not concentrated in a zero Lebesgue measure Borel set of $\reels^d$.

We assume that the covariance of process $X$, say $r$, is $C^2(\reels^d, \reels)$ except at the origin.
Note that this implies that $\int_{\reels^d} \normp[d]{x}^2 \ud F(x)= +\infty$.\\
We consider regularizations of the trajectories of $X$ by means of convolutions of the form $X_\ve(t):= \Psi_\ve \sstar X(t)$, $\ve >0$, $\Psi_\ve(t):=\frac{1}{\ve^d} \Psi(\frac{t}\ve)$, where $\Psi$ is $C^2(\reels^d, \reels)$ with compact support, nonnegative, depending only on the norm and such that $\int_{\reels^d} \Psi(t) \ud t=1$.

Note that $X_\ve: \Omega \times \reels^d \to \reels$ is still a centered continuous Gaussian stationary isotropic process and is $C^2(\reels^d, \reels)$.\\
Let $T$ be an open bounded set of $\reels^d$ and let $u$ be a fixed level.
Let also $f: T \to \reels$ be a continuous bounded function.\\
The objective of this paragraph is to approximate the renormalized local time  $L_{X}^f(u, T)$, by functionals of the $u$ level set of the regular approximating process $X_\ve$. We need some definitions.

Let us define $\theta_d:=\bbE[\normp[d]{\xi}]$, where $\xi$ is a standard Gaussian vector of $\reels^d$.
\\
We will note $r_\ve$ the covariance of $X_\ve$, $\ve >0$ and by $r^{(\ve)}$ that of $X$ and $X_\ve$, that is $r^{(\ve)}(t) :=\bbE[X_\ve(t)X(0)]$.\\
Since $X_\ve$ is isotropic, $\bbE[\nabla\!X_\ve(0) \nabla\!X_\ve(0)^t]= \mu_\ve I_d$, where $\mu_\ve:= - \dfrac{\partial ^{2}r_{\epsilon}}{\partial t_{i}^2}(0)$ for any $i=1, \dots d$ and $I_d$ stands for the identity matrix in $\reels^d$.
\spacebefore
\begin{rema}
\label{convergence moment spectral}
Since $\int_{\reels^d} \normp[d]{x}^2 dF(x)= +\infty$, Fatou's lemma implies that $\lim_{\ve \to 0} \mu_\ve= +\infty$.
We then deduce that for $\ve $ small enough and for all $t \in T$, the vector $(X_\ve(t), \nabla\!X_\ve(t))$ has a density.
\hfill$\bullet$
\end{rema}

We define the functional
$$
	\xi_\ve^f(u):=\ds \frac{\theta_d^{-1}}{\sqrt{\mu_\ve}} \int_{C^{D^{r}}_{T, X_\ve}(u)} f(t) \ud \sigma_{d-1}(t),\quad \ve >0.
$$%
We will show the following theorem.
\spacebefore\begin{theo}
\label{approximation L2}
The random variable $\xi_\ve^f(u)$ converge in $L^{2}(\Omega)$ towards $L_X^f(u, T)$ as $\ve$ tends to zero.
\end{theo}
\spacebefore
\begin{proofarg}{Proof of Theorem \ref{approximation L2}}
Proposition \ref{approximation temps local Bis} ensures that it is sufficient to show that
$$
\lim_{\delta \to 0}\lim_{\ve \to 0} \bbE\!\left[\left(\xi_\ve^f(u)-\eta_{\delta}^f(u)\right)^2\right]=0.
$$%
Let
$$
	\Sigma:=  \bbE\!\left[\left(\eta_{\delta}^f(u)\right)^2\right] + 
		\bbE\!\left[\left(\xi_\ve^f(u)\right)^2\right] -
		2 \bbE\!\left[\xi_\ve^f(u)\eta_{\delta}^f(u)\right].
$$%
Propositions \ref{approximation temps local Bis} and  \ref{expectation local time Bis} imply that 
$$
\lim_{\delta \to 0} \bbE\!\left[\left(\eta_{\delta}^f(u)\right)^2\right] = \int_{T \times T}  f(s) f(t) {p}_{X(s), X(t)}(u, u) \ud s \ud t:=\alpha(u).
$$%
Let us consider  $\bbE\!\left[\left(\xi_\ve^f(u)\right)^2\right]$.\\
We claim that assumptions of Theorem \ref{Cabana4 Bis} are fulfilled.\\
It remains to show that for any $(s, t) \in T \times T$ with $s \neq t$, the vector $(X_\ve(s), X_\ve(t))$ has a density.

This fact follows from Lemma \ref{minoration pour r epsilon} whose proof is given in \cite[Lemma 2]{MR1222362}.
\spacebefore
\begin{lemm}
\label{minoration pour r epsilon}
$\exists \bC  >0$, $\exists  \ve_0 >0$, $\forall \ve \le \ve_0$, $\forall (s, t) \in \overline{T}\times \overline{T}$, one has $r^2_\ve(0) -r^2_\ve(s-t) \ge \bC  \normp[d]{s-t}^2$.
\end{lemm}
To finish the argument, the Borel-Tsirelson-Ibragimov-Sudakov inequality (see in  \cite[Theorem 2.1.1]{MR2319516}) ensures in the one hand that 
$$
\bbE\!\left[\left(\sup_{x \in T} \normp[1,d]{\nabla^2X_\ve(x)}^{(s)}\right)^{2d}\right] < \infty.
$$%
On the other hand, we use the fact that $f$ is a continuous bounded function on $T$ to verify the hypothesis ${\bH _7}$.
In this way it is enough to establish that $\forall u \in \reels$, we have 
\begin{multline*}
\int_{T \times T} 
\bbE\!\left[\normp[d]{\nabla\!X_\ve(s)} \normp[d]{\nabla\!X_\ve(t)} \given{X_\ve(s)=X_\ve(t)=u}\right] \\
\times {p}_{X_\ve(s), X_\ve(t)}(u, u) \ud s \ud t < \infty,
\end{multline*}
where ${p}_{X_\ve(s), X_\ve(t)}(\bm{\cdot}, \bm{\cdot})$ is the density of vector $(X_\ve(s), X_\ve(t))$.
The proof of this last fact can be found for example in \cite[Proposition 2]{MR1222362}.\\
The Rice formula of the second order follows from the application of Theorem  \ref{Cabana4 Bis}.
That is, we have
\begin{align*}
\MoveEqLeft{\bbE\!\left[\left(\xi_\ve^f(u)\right)^2\right]
  = \theta_d^{-2} \int_{T \times T}  f(s) f(t)}\\
	 {} \qquad\times & \bbE\!\left[\normp[d]{\nabla\!X_\ve(s)/\sqrt{\mu_\ve}} \normp[d]{\nabla\!X_\ve(t)/\sqrt{\mu_\ve}} \given{X_\ve(s)=X_\ve(t)=u}\right] \\
	 {} \qquad\times &  {p}_{X_\ve(s), X_\ve(t)}(u, u) \ud s \ud t.
 \end{align*}
\spacebefore
\begin{rema}
Note that $\forall u \in \reels$, almost surely 
$$
C^{D^{r}}_{T, X_\ve}(u)=C_{T, X_\ve}(u).
$$%
\hfill$\bullet$
\end{rema}
To apply the dominated convergence theorem we need an upper bound for the expression within the integral in last equality.
Let us quote the following lemma cited as Appendix 1 in \cite{MR1222362}.
\spacebefore
\begin{lemm}
\label{majoration norme}
Let $(Y, Z, X_1, X_2)$ be a centered Gaussian vector taking its values into $\reels^{2d+2}$, while $Y$ and $Z$ take their values in $\reels^d$ and $X_1, X_2$ in $\reels$.
Let us suppose that $\bbE[X_1^2]=\bbE[X_2^2]=\sigma^2$, $\bbE[X_1X_2]=\gamma$ with $\sigma^4-\gamma^2 \neq 0$, then $\forall u \in \reels$
\begin{multline*}
\bbE\!\left[\normp[d]{Y} \normp[d]{Z} \given{X_1=X_2=u}\right]\\
 \le \!\left({\bbE\!\left[\normp[d]{Y}^2\right]}\right)^{\frac{1}{2}} \!\left({\bbE[\normp[d]{Z}^2]}\right)^{\frac{1}{2}} \left[{
1+4\sigma^2(\sigma^2+\gamma)^{-2}u^2}\right].
\end{multline*}
\end{lemm}
In particular if, for all $s, t \in T$, with $s \neq t$,  we choose $Y:=\nabla\!X_\ve(s)$ and $Z:= \nabla\!X_\ve(t)$, $X_1:= X_\ve(s)$
and $X_2:= X_\ve(t)$, thus $\sigma^2:= r_\ve(0)$ and $\gamma:= r_\ve(s-t)$.
\\
 Lemma \ref{minoration pour r epsilon} ensures that assumptions of Lemma \ref{majoration norme} are satisfied.\\
 Moreover, note that since for all $(s, t) \in \overline{T} \times \overline{T}$, $r_\ve(0)+r_\ve(s-t)$ converges uniformly  to $r(0)+r(s-t) > 0$ when $\ve$ tends to zero, so for all $(s, t) \in \overline{T} \times \overline{T}$, $r_\ve(0)+r_\ve(s-t) \ge \bC  >0$ for $\ve $ sufficiently small.
 
We deduce that for all $s, t \in T$, with $s \neq t$ and $\ve \le \ve_0$, we have
\begin{align*}
\MoveEqLeft[2]{\abs{f(s)} \abs{f(t)} \times  
\bbE\!\left[ \normp[d]{\nabla\!X_\ve(s)/\sqrt{\mu_\ve}} \normp[d]{\nabla\!X_\ve(t)/\sqrt{\mu_\ve}} \given{X_\ve(s)=X_\ve(t)=u}\right]} \\
	&\specialpos{\hfill\times {p}_{X_\ve(s), X_\ve(t)}(u, u) \ud s \ud t}  \\
	&\le \bC  \left[{ 
 1 +4 \frac{r_\ve(0) u^2}{(r_\ve(0)+r_\ve(s-t))^2}
 }\right] \frac{1}{\sqrt{r_\ve^2(0)-r_\ve^2(s-t)}} \\
 	&\specialpos{\hfill\times \exp\!\left(-\frac{u^2}{r_\ve(0)+r_\ve(s-t)}\right)}\\
 	&\le \bC  \frac{1}{\normp[d]{s-t}} \in L^1(T \times T).
 \end{align*}
We deduce that 
\begin{multline*}
\lim_{\ve \to 0}
 \bbE\!\left[\left(\xi_\ve^f(u)\right)^2\right]
 	= \theta_d^{-2} \int_{T \times T}  f(s) f(t)  \\ 
\times \lim_{\ve \to 0}
	\left\{{
		\bbE\!\left[
			 \normp[d]{\nabla\!X_\ve(s)/\sqrt{\mu_\ve}} \normp[d]{\nabla\!X_\ve(t)/\sqrt{\mu_\ve}} \given{X_\ve(s)=X_\ve(t)=u}
			\right] 
	 }\right.\\
	 \times \left.{
	 {p}_{X_\ve(s), X_\ve(t)}(u, u)}\right\} \ud s \ud t,
\end{multline*}
provided that the limit inside the integral exists.\\
In pursuit of our goal, we use the decomposition given in \cite[p$.$ 60]{MR871689}, that is for $\tau:=s-t$, where $s, t \in T$ and $\ve >0$
\begin{align*}
\nabla\!X_\ve(0)&= \xi_\ve + (X_\ve(0) \alpha_\ve + X_\ve(\tau) \beta_\ve)\\
\nabla\!X_\ve(\tau)&= \xi_\ve^{\star} - (X_\ve(\tau) \alpha_\ve + X_\ve(0) \beta_\ve),
\end{align*}
where $\xi_\ve$ and $\xi_\ve^{\star}$ are centered Gaussian vectors taking values in $\reels^d$, with joint Gaussian distribution, each of them independent of $(X_\ve(0), X_\ve(\tau))$, such that, if $\nabla r_\ve$ stands for the Jacobian of the $X_\ve$ covariance function,
\begin{align*}
\alpha_\ve&:= \frac{r_\ve(\tau) \nabla r_\ve(\tau)}{r^2_\ve(0)-r^2_\ve(\tau)}\\
 \beta_\ve&:= -\frac{r_\ve(0) \nabla r_\ve(\tau)}{r^2_\ve(0)-r^2_\ve(\tau)}\\
\Var(\xi_\ve)&\phantom{:}=\Var(\xi_\ve^{\star})=-\nabla^{2}r_\ve(0) -\frac{r_\ve(0) }{r^2_\ve(0)-r^2_\ve(\tau)} \nabla r_\ve(\tau) {\nabla r_\ve(\tau)}^{t}\\
\Cov(\xi_\ve, \xi_\ve^{\star})&\phantom{:}=-\nabla^{2}r_\ve(\tau) -\frac{r_\ve(\tau) }{r^2_\ve(0)-r^2_\ve(\tau)} \nabla r_\ve(\tau) {\nabla r_\ve(\tau)}^{t}.
\end{align*}
With these notations, we obtain
{\def\termeA{\theta_d^{-2} \int_{T \times T}  f(s) f(t)}
\def\termeB{\normp[d]{\left.\!\left({\xi_\ve-  \frac{\nabla r_\ve(s-t) u}{r_\ve(0)+r_\ve(s-t)}}\right)\right/\sqrt{\mu_\ve}}}
\newlength{\hauteurB}
\def\termeC{\normp[d]{\left.\!\left({\xi_\ve^{\star}+  \frac{\nabla r_\ve(s-t) u}{r_\ve(0)+r_\ve(s-t)}}\right)\right/\sqrt{\mu_\ve}}}
\def\termeD{\times {p}_{X_\ve(s), X_\ve(t)}(u, u)}
\begin{align*}
\lim_{\ve \to 0}
 \bbE\!\left[\left(\xi_\ve^f(u)\right)^2\right]&= \termeA \\ 
&\quad  \times \lim_{\ve \to 0}
	\left\{
		\bbE\!\left[
 			\termeB\rule{13pt}{0pt}
			\right.
		\right.\\
&\specialpos{\hfill\left. \times\termeC
			 \right]\quad} \\
&\specialpos{\hfill\left.{ \termeD}\rule{0pt}{\heightof{$\ds\termeB$}}\right\} \ud s \ud t.}
\end{align*}}
Using Remark \ref{convergence moment spectral}, one can show that for all $s, t \in T$,
$$
 \frac{\nabla r_\ve(s-t)}{\sqrt{\mu_\ve}} \cvg[\ve \to 0]{} 0 \quad \text{and} \quad  \frac{\nabla^{2}r_\ve(s-t)}{\mu_\ve} \cvg[\ve \to 0]{} 0.
$$%
Furthermore, 
$$
	r_\ve(s-t) \cvg[\ve \to 0]{} r(s-t),
$$%
so
$$
\lim_{\ve \to 0}
 \bbE\!\left[\left(\xi_\ve^f(u)\right)^2\right]= \alpha(u).
$$
Now to finish the proof of theorem, we need to compute the term $\bbE\!\left[\xi_\ve^f(u)\eta_{\delta}^f(u)\right]$.
We have
$$
 \bbE\!\left[\xi_\ve^f(u)\eta_{\delta}^f(u)\right]= \frac{1}{2\delta}\int_T  f(s)  \frac{\theta_d^{-1}}{\sqrt{\mu_\ve}} 
\bbE\!\left[ \int_{C^{D^{r}}_{T, X_\ve}(u)} f(t) \1_{\abs{X(s)-u} < \delta}\ud \sigma_{d-1}(t)\right] \ud s.
$$%
 
We are going to use the first order Rice formula given in Theorem \ref{Cabana4}.
For any  fixed $s \in T$, the process $Y$ is defined as
$Y(t):=f(t) \1_{\abs{X(s)-u} < \delta}$, $t \in T$.
This model is not part of the proposed models but a slight modification allows us to include it in those studied.
 As before when we applied the Rice formula of the second order, we check that the hypotheses making the formula valid are verified.
We obtain
\begin{multline*}
 	\bbE\!\left[\xi_\ve^f(u)\eta_{\delta}^f(u)\right]= \frac{1}{2\delta}\int_T \int_T   \theta_d^{-1} f(s) f(t)\\
 	\times\bbE\!\left[
 	\normp[d]{\nabla\!X_\ve(t)/\sqrt{\mu_\ve}} \1_{\abs{X(s)-u} < \delta}\given{X_\ve(t)=u}\right]
 	 \times {p}_{X_\ve(t)}(u) \ud t \ud s,
 \end{multline*}
 where ${p}_{X_\ve(t)}(\bm{\cdot})$ stands for the density of the random variable $X_\ve(t)$.\\
 Since function $f$ is bounded,  $\forall\ve \le \ve_0$ and $\forall s, t \in T$
$$
\theta_d^{-1} \abs{f(s)} \abs{f(t)} 
\bbE\!\left[
 \normp[d]{\nabla\!X_\ve(t)/\sqrt{\mu_\ve}} \1_{\abs{X(s)-u} < \delta}\given{X_\ve(t)=u}\right]
  \times {p}_{X_\ve(t)}(u) \le \bC .
$$%
The dominated convergence theorem applies.
Using once again Remark \ref{convergence moment spectral}, we can show that for $s, t \in T$,
$$
	\frac{\nabla r^{(\ve)}(s-t)}{\sqrt{\mu_\ve}} \cvg[\ve \to 0]{} 0.
$$%
Furthermore $r^{(\ve)}(s-t)\cvg[\ve \to 0]{} r(s-t)$.
Thus
 \begin{align*}
\MoveEqLeft[2]{\lim_{\ve \to 0} \bbE\!\left[\xi_\ve^f(u)\eta_{\delta}^f(u)\right]
	=\int_T \int_ T  \frac{\theta_d^{-1}}{2\delta}  f(s) f(t)} \\
\times {} &\lim_{\ve \to 0} \left\{{\bbE\!\left[
 \normp[d]{\nabla\!X_\ve(t)/\sqrt{\mu_\ve}} \1_{\abs{X(s)-u} < \delta}\given{X_\ve(t)=u}\right]
  \times {p}_{X_\ve(t)}(u)}\right\}  \ud t \ud s\\
  = {} & \int_T \int_ T 
 	\frac{\theta_d^{-1}}{2\delta}  f(s) f(t) 
	\bbE\!\left[ \normp[d]{\xi} \1_{\abs{X(s)-u} < \delta}\given{X(t)=u}\right]
  \times 
 {p}_{X(t)}(u)  \ud t \ud s,
\end{align*}
where ${p}_{X(t)}(\bm{\cdot})$ stands for the density of the random variable $X(t)$, and $\xi$ is a standard Gaussian vector taking values in $\reels^d$, independent of $(X(s), X(t))$,
for all $s, t \in T$.
Therefore
  \begin{align*}
\MoveEqLeft[2]{\lim_{\ve \to 0} \bbE\!\left[\xi_\ve^f(u)\eta_{\delta}^f(u)\right]}\\
& = \int_T \int_ T 
 \frac{1}{2\delta}  f(s) f(t) \bbE[
\1_{\abs{X(s)-u} < \delta}\given{X(t)=u}]
  \times 
 {p}_{X(t)}(u)  \ud t \ud s\\
 & =\frac{1}{2\delta} \int_{u-\delta}^{u+\delta} \!\left({\int_T \int_ T  f(s) f(t) \ {p}_{X(s), X(t)}(u, v) \ud s \ud t}\right) \ud v.
\end{align*}
 By using the already mentioned continuity of the function
$$
 	(x, y) \mapsto  \int_{T \times T} f(s) f(t) {p}_{X(s), X(t)}(x, y) \ud s \ud t,
$$%
we finally obtain that 
$$
	\lim_{\delta \to 0}\lim_{\ve \to 0} \bbE\!\left[\xi_\ve^f(u)\eta_{\delta}^f(u)\right]= \alpha(u).
$$%
This completes the proof of the theorem.
\end{proofarg}
\subsection{Rate of convergence}
In the same context as before, and for technical reasons, we renormalize the functional $\xi_\ve^f(u)$ at point $u$ as follows.
We consider 
$$
 V_\ve^f(u):= \frac{\theta_d^{-1}}{\sqrt{\mu_\ve}}\sqrt{r_\ve(0)} 
	 \int_{C_{T, X_\ve/\sqrt{r_\ve(0)}}(u)}
	 	f(t) \ud \sigma_{d-1}(t), \ve >0,
$$%
 where $f:T \to \reels$ is still  a continuous and bounded function and $T$ is an open bounded set of $\reels^d$.\\
As in the previous section, we can prove that  $ V_\ve^f(u)$ converges in $L^2(\Omega)$ towards the renormalized local time $L_{X/\sqrt{r(0)}}^f(u, T)$.
 We are aiming at giving the rate of this convergence.
 
We will make the following additional hypothesis on the covariance $r$ of process $X$.
It is assumed that the covariance can be written as $r(t)= r(0)- \normp[d]{t}^{2\alpha} L(\normp[d]{t})$, $0< \alpha <1$, where $\lim_{x \to 0^{+}} L(x)=\bC _L >0$, where $L \ge 0$ is even and has two continuous derivatives except at the origin, which satisfy $\abs{x} L^{\prime}(\abs{x})=O(1)$ and $x^2 L^{\prime\prime}(\abs{x})=O(1)$ as $x \to 0$.\\
Moreover, in case where  $\alpha = {1}/{2}$  and $d=2$, we further assume that function $L$ has a continuous third derivative except at the origin, which satisfies $\abs{x} L^{\prime}(\abs{x})=o(1)$, $x^2 L^{\prime\prime}(\abs{x})=o(1)$ and $\abs{x^3} L^{\prime\prime\prime}(\abs{x})=O(1)$ as $x \to 0$.\\
The rate will then depend on the value of $\alpha$.\\
 We will draw largely on the results achieved in \cite{MR1644037}, where the framework was the same but in case where $d=1$.
\\
We study the convergence of
 \begin{align*}
 	 \Xi_\ve(f, g)&:=\ve^{-a(\alpha)} [\ln(1/\ve)]^{-b(\alpha)} \int_{\reels} g(u) \!\left({ V_\ve^f(u)-L_{X/\sqrt{r(0)}}^f(u, T)}\right) \ud u\\
	& :=\ve^{-a(\alpha)}  [\ln(1/\ve)]^{-b(\alpha)} V_\ve(f, g),
 \end{align*}
 where $g: \reels \to \reels$ is $C^2(\reels, \reels)$ and also $L^2(\varphi(x) \ud x)$.
\\
 We also assume that $g^{\prime}$ and $g^{\prime\prime}$ are $L^2(\varphi(x) \ud x)$.\\
 Furthermore, in case where  $\alpha = {1}/{2}$  and $d >2$, we suppose that $g^{\prime\prime}$ belongs to $L^4(\varphi(x) \ud x)$.\\
 The Hermite coefficients of function $g$ being denoted by $c_n$, $n \in \naturels$,
   we set for $x \in \reels$
$$
    Hg(x):= xg^{\prime}(x)-g^{\prime\prime}(x)=  \sum_{n=1}^{\infty} n c_n H_n(x).
$$
Thus the following series is convergent,
\begin{equation}
	\sum_{n=0}^{\infty} n^2  n! c_n^2 < \infty.
	\label{finitude cn}
 \end{equation}
   Let us introduce some notations and hypothesis.\\
    We define the constants $C_{\alpha}$ and $\ell_{\alpha}$ as 
\begin{align*}
\bC _{\alpha}&:= \bC _L \int_{\reels^{d}} \int_{\reels^{d}} \Psi(v) \Psi(w) \normp[d]{v-w}^{2\alpha} \ud v \ud w\\
&\phantom{:}= \bC _L \int_{\reels^{d}}   \Phi(w) \normp[d]{w}^{2\alpha}  \ud w,
\end{align*}
where $\Phi= \Psi \sstar \Psi$ is the convolution product of $\Psi$ with itself and
$$
\ell_{\alpha}:=\frac{\bC _{\alpha}}{2r(0)}.
$$
The constants $\bC _{\Psi}$ and $\bdK_{1/2}$ are
\begin{align*}
\bC _{\Psi}&:= \frac{1}{2d} \int_{\reels^d} \normp[d]{w}^2 \left[{\Phi(w) - 2 \Psi(w)}\right] \ud w,\\ \shortintertext{and} 
\bdK_{1/2}&:= \bC _L \int_{\reels^{d}}   \Psi(w) \normp[d]{w}  \ud w.
 \end{align*}
We denote by $\varphi_m$ (resp$.$
$N_m$) the standard Gaussian density (resp$.$ \rv) on $\reels^m$, for $m \in \naturels-\{0\}$ with the convention that $\varphi_1=\varphi$.\\
For $\bk :=(k_1, k_2, \dots, k_m) \in \naturels^m$ and $x:=(x_1, x_2, \dots, x_m) \in \reels^m$, we set
$$\abs{\bk }:=k_1+k_2+\dots +k_m,\quad \bk !:=k_1!k_2!\cdots k_m! \text{ and }
\widetilde{H}_{\bk }(x):= \prod_{\substack{1 \le j \le m}}H_{k_{j}}(x_j).$$
We define the function $h: \reels^d \to \reels$: for $x=(x_1, x_2, \dots, x_d) \in \reels^d$,
$$
h(x):= \theta_d^{-1} \normp[d]{x}-1= \sum_{q=0}^{\infty} \sum_{
\substack{
\bk  \in \naturels^d \\ \abs{\bk }=q
}
} a(\bk ) \widetilde{H}_{\bk }(x),
$$%
where for $\bk  \neq \0$,
\begin{equation}
 a(\bk ):= \frac{\theta_d^{-1}}{{\bk!}} \int_{\reels^d} \normp[d]{x} \widetilde{H}_{\bk }(x) \varphi_d(x) \ud x,
\label{coef_ak}
\end{equation}%
and $a(\bk ):=0$ for $\bk   =\0 $.
Note that since the polynomial $H_n$ is odd when $n$ is odd, the summation on $q$ begins at $q=2$ since all the indices $k_i$, $i=1, \dots,d$ such $\abs{\bk }=q$ have to be even.\\
Also define
$$
\chi_{d}^2(\alpha):=\int_{\reels^d}  \dfrac{\partial ^{2} \Phi}{\partial x_1^2}(v) \normp[d]{v}^{2\alpha}  \ud v,
$$%
and if $d=2$ (and $\alpha=1/2$)
$$
\sigma := \frac{2 \sqrt{\pi}}{\chi_2^2(1/2)}.
$$%
When $\alpha \ge \frac{1}{2}$ we suppose that the process $X$ has a spectral density $s$.
In the case where  $\alpha = {1}/{2}$ , we will make the assumption that the open set $T$ is a rectangle.

\paragraph{Asymptotic equivalence symbol $\equiv$} This symbol will be reserved for the asymptotic equivalence of real deterministic sequences.

\paragraph{Asymptotic equivalence symbol $\simeq$} For a set of random variables $X_n$, the notation $X_n=o_P(1)$ means that $X_n\cvg[n\to\infty]{Prob.}0$.
Furthermore, if $Y_n$ is a set of random variables such that $X_n=Y_n+o_P(1)$, we will note $X_n \simeq Y_n$.

Let us state the following two theorems.
\spacebefore\begin{theo}
\label{vitesse variance temps local}
The random variable $ V_\ve(f, g)$ satisfies
$$
\bbE[ V_\ve^2(f, g)]= O(\ve^{4\alpha}+ \ve^2 \ln(1/\ve) + \ve^{4(1-\alpha)}).
$$
Morover:
\begin{enumerate}
\item\label{itm4.6.19.1} If $\alpha < \frac{1}{2}$, $a(\alpha)=2\alpha$, $b(\alpha)=0$ and the limit of the variance divided by $\ve^{4\alpha}$ is 
$$
\ell_1:= \ell_{\alpha}^2  \sum_{n=1}^{\infty} n^2  n! c_n^2 \int_{T \times T} f(t) f(s) \left(\frac{r(t-s)}{r(0)}\right)^n \ud t \ud s.
$$
\item\label{itm4.6.19.2} If $\alpha > \frac{1}{2}$,  $a(\alpha)=2(1-\alpha)$, $b(\alpha)=0$ and the limit of the variance divided by $\ve^{4(1-\alpha)}$ is 
\begin{multline*}
	\ell_2:=2\!\left({\frac{a({2,0,\dots,0})}{\bC _L  \chi_{d}^2(\alpha)} }\right)^2  \sum_{q=2}^{\infty} \sum_{\ell=(q-4)\vee 0}^{q-2} (q-2)! c_{q-2}^2  
	{{q-2}\choose{\ell}}
	{{2}\choose{q-2-\ell}}\\
\times\int_{T \times T} f(t) f(s) 
  	\!\left(\frac{r(t-s)}{r(0)}\right)^{\ell}
	\sum_{i, j=1}^d \!\left({-\frac{1}{r(0)} \dfrac{\partial r}{\partial t_{i}}(t-s)  \dfrac{\partial r}{\partial t_{j}}(t-s) }\right )^{q-2-\ell}\\
  \!\left({ -\dfrac{\partial ^2  r}{\partial t_{i}\partial t_{j}}(t-s) }\right )^{\ell - (q-4)}  \ud t \ud s.
\end{multline*}
 \item\label{itm4.6.19.3} if $\alpha = {1}/{2}$ and $d >2$, $a(\alpha)=1$, $b(\alpha)=0$ and the limit of the variance divided by $\ve^{2}$ is $\ell_1 + \ell_2 + \ell_3 +\ell_4$, where $\ell_3$ is
$$
\ell_3:=\frac{\bC _{\Psi}}{r(0)}  \sum_{n=1}^{\infty} n n! c_n^2 \int_{T \times T} f(t) f(s) \!\left(\frac{r(t-s)}{r(0)}\right)^{n-1} \!\left({ \sum_{i=1}^d  \dfrac{\partial ^2  r}{\partial t_{i}^2}(t-s) }\right )  \ud t \ud s,
$$%
and $\ell_4$ is
 \begin{multline*}
\ell_4:=2 \frac{\ell_{{1}/{2}}}{r(0)} \!\left({\frac{a({2,0,\dots,0})}{\bC _L  \chi_{d}^2(1/2)} }\right) \sum_{q=2}^\infty c_{q-2}  c_q q  q! \\
 \int_{T \times T} f(t) f(s) \!\left(\frac{r(t-s)}{r(0)}\right)^{q-2} \!\left({\sum_{i=1}^d \!\left({  \dfrac{\partial   r}{\partial t_{i}}(t-s) }\right )^2}\right) \ud t \ud s.
 \end{multline*}
 \item\label{itm4.6.19.4} if $\alpha = {1}/{2}$ and $d =2$, $a(\alpha)=1$, $b(\alpha)=1/2$ and the limit of the variance divided by $\ve^{2} \ln(1/\ve)$ is 
  \begin{eqnarray*}
\ell_5:= a^2(2, 0) \sigma^2 \!\left({\int_{T} f^2(t) \ud t }\right) \bbE\!\left[g^2\!\left({X(0)}/{\sqrt{r(0)}}\right)\right].
\end{eqnarray*}
\end{enumerate}
\end{theo}
Before stating Theorem \ref{vitesse convergence temps local}, let us define the following \rv $\xi_\ve$.
For $t \in T$, we let $\xi_\ve(t):={\!\left(X_\ve(t)-X(t)\right)}/\ve$.
\spacebefore\begin{theo}
\label{vitesse convergence temps local}
\begin{enumerate}
\item\label{itm1:theo4.6.20} If $\alpha < \frac{1}{2}$, $a(\alpha)=2\alpha$, $b(\alpha)=0$,
we have  
$$
\Xi_\ve(f, g)\cvg[\ve \to 0]{$L^2(\Omega)$}X(f,g)
$$%
where
 \begin{align*}
X(f,g)&:=\ell_{\alpha} \int_{-\infty}^{\infty} Hg(x) L_{X/\sqrt{r(0)}}^f(x, T) \ud x\\
&\phantom{:}=\ell_{\alpha} \int_{T} f(t) \!\left({
{{X(t)}/{\sqrt{r(0)}}} g^{\prime}\!\left({{X(t)}/{\sqrt{r(0)}}}\right) -g^{\prime\prime}\!\left({{X(t)}/{\sqrt{r(0)}}}\right)}\right) \ud t.
 \end{align*}
 \item\label{itm2:theo4.6.20} If $\alpha > \frac{1}{2}$, $a(\alpha)=2(1-\alpha)$, $b(\alpha)=0$, we have  
$$
\Xi_\ve(f, g)\cvg[\ve \to 0]{$L^2(\Omega)$}Y(f,g)
$$%
where

 \begin{multline*}
Y(f, g):=-\!\left({\frac{a({2,0,\dots,0})}{\bC _L  \chi_{d}^2(\alpha)} }\right)
\sum_{ j=1}^d  \sum_{ k=2}^{\infty} \frac{ c_{k-2}}{k!}  \!\left({\frac{1}{\sqrt{r(0)}}}\right)^{k-2}\\ \times 
\int_{\reels^{dk}} K_f(\lambda_1+\lambda_2+\dots +\lambda_k)  \\
\times\sum_{\Pi \in \Pi_{k}} \lambda_{\Pi(k-1)}^{(j)} \lambda_{\Pi(k)}^{(j)} \ud Z_X(\lambda_1)  \ud Z_X(\lambda_2) \dots  \ud Z_X(\lambda_k),
 \end{multline*}
and where we set for  $\lambda \in \reels^d$, $\lambda:= (\lambda^{(1)}, \lambda^{(2)}, \dots, \lambda^{(d)})$ and
 \begin{equation}
 \label{Kf}
  K_f(\lambda):= \int_{T} f(t) e^{i \prodsca{\lambda}{t} } \ud t,
  \end{equation}
 and we denoted by $\prodsca{\cdot}{\cdot}$ the canonical scalar product in $\reels^d$.\\
 $ \prod_{k}$ is the set of permutations of $\{1, 2, \dots, k\}$, $\ud Z_X$ is the random spectral measure associated with $X$ and the integral is an It\^o-Wiener integral (see \cite{MR3155040}).
  \item\label{itm3:theo4.6.20} If $\alpha = {1}/{2}$ and $d >2$, $a(\alpha)=1$, $b(\alpha)=0$, the random variable $\xi_\ve$ converges vaguely to $\xi$ (see the Section \ref{vague convergence} of Appendix \ref{chap convergence vague} for definition of the \emph{vague convergence}) and the random variable  $\Xi_\ve(f, g)$ converges stably and vaguely to $X(f,g)+Y(f,g) +Z(f,g)$ (see also the Section \ref{stable convergence} of Appendix \ref{chap convergence vague} for definition of the \emph{stable vague convergence}), where 
  \begin{multline*}
  Z(f,g):= \frac{\bK _{{1}/{2}}}{r(0)} \int_T f(t) g^{\prime\prime}\!\left({{X(t)}/{\sqrt{r(0)}}}\right) \ud t\\
  +   \frac{1}{\sqrt{r(0)}} \int_T f(t) g^{\prime}\!\left({{X(t)}/{\sqrt{r(0)}}}\right)  \xi(t) \ud t,
  \end{multline*}
  with
$$
  \int_T f(t) g^{\prime}\!\left({{X(t)}/{\sqrt{r(0)}}}\right) \xi(t) \ud t := \lim_{\ve \to 0} \frac{1}{\sqrt{r(0)}} \int_T f(t) g^{\prime}\!\left({{X(t)}/{\sqrt{r(0)}}}\right) \xi_\ve(t) \ud t,
$$%
the last convergence being stable vague convergence.\\
Furthermore for all $x(\bm{\cdot})$ belonging to $L^{2}(T)$, 
{\def\termeA{\bC _{\Psi} \int_{T \times T} f(t) f(s) \!\left({ \sum_{i=1}^d  \dfrac{\partial ^2  r}{\partial t_{i}^2}(t-s) }\right)
}
\def\termeB{g^{\prime}\!\left(x(t)/\sqrt{r(0)}\right)}
\def\termeC{g^{\prime}\!\left({x(s)}/{\sqrt{r(0)}}\right) \ud t \ud s}
\begin{align*}
\MoveEqLeft[3]{{\calL}\!\left({\int_T f(t) g^{\prime}\!\left({{X(t)}/{\sqrt{r(0)}}}\right)  \xi(t) \ud t \given{X(\bm{\cdot})=x(\bm{\cdot})}}\right)}\\
= {} & {\calN}
	\!\left(0;
		\termeA
		\termeB\right.\\ 
&		\left.  \times \termeC \rule{0pt}{17.5pt}\right).
 \end{align*}}
 \item\label{itm4:theo4.6.20} if $\alpha = {1}/{2}$ and $d =2$, $a(\alpha)=1$, $b(\alpha)=1/2$, the random variable  $ \Xi_\ve(f, g)$ stably converges when $\ve \to 0$ to $N(f, g)$ (see the Section \ref{stable convergence} of Appendix \ref{chap convergence vague} for definition of the \emph{stable convergence}), where 
$$
 	N(f, g):=a(2, 0) \sigma   \!\left({\int_T f(t)  g\!\left({{X(t)}/{\sqrt{r(0)}}}\right)\ud \widehat{W}(t)}\right),
 $$%
and $\widehat{W}$ is a standard Brownian sheet in $\reels^2$ independent of $X$.
\end{enumerate}
\end{theo}
\spacebefore
\begin{rema}
Note that in case where $\alpha > \frac{1}{2}$, if we assume that $f$ is a constant function whose value is one and that the open set $T$ is a rectangle with the form $T= \prod_{\ell =1}^{d}  \introo{a_{\ell}}{b_{\ell}}$, then the function $K_f$ can be expressed as follows
$$
K_f(\lambda)= \prod_{\ell =1}^{d} \left[{\exp(i\lambda^{(\ell)}b_{\ell})- \exp(i\lambda^{(\ell)}a_{\ell})}\right]/i \lambda^{(\ell)},
$$%
for $\lambda := (\lambda^{(1)}, \lambda^{(2)}, \dots, \lambda^{(d)})\in \reels^d$.

So that in the very special case where $T= introo{0}{1}^d$, the kernel has the following particular form
$$
K_f(\lambda)= \prod_{\ell =1}^{d} \left[{\exp(i\lambda^{(\ell)})- 1}\right]/i \lambda^{(\ell)}.
$$%
Thus, the corresponding result stated in \ref{itm2:theo4.6.20}.
can be linked to \cite[Theorem 2 (iii)]{MR1644037}.
\hfill$\bullet$
\end{rema}
\spacebefore
\begin{proofarg}{Proof of Theorem \ref{vitesse variance temps local}}
We apply the coarea formula established in Corollary \ref{coa2}  to the Borel set $B:=T$ and to the  functions $h:\reels^d \times \reels \to \reels$, $G:T \subset \reels^d \to \reels$, defined by
 \begin{eqnarray*}
 h(t, u):=f(t) g(u); \quad  G:=\frac{X_\ve}{\sqrt{r_\ve(0)}}.
 \end{eqnarray*}
 We deduce on the one hand that
  \begin{eqnarray*}
 \int_{\reels} g(u)  V_\ve^f(u) \ud u
= \frac{\theta_d^{-1}}{\sqrt{\mu_\ve}} \int_T f(t) g\!\left({X_\ve(t)}/{\sqrt{r_\ve(0)}}\right) \normp[d]{\nabla\!X_\ve(t)} \ud t.
 \end{eqnarray*}
 On the other hand, Theorem \ref{GH} applied to the functions $h(t, u):=f(t) g(u)$ and to $X:={X}/{\sqrt{r(0)}}$ yields that
$$
 \int_{\reels} g(u) L_{X/\sqrt{r(0)}}^f(u, T)  \ud u
=  \int_T f(t) g\!\left({{X(t)}/{\sqrt{r(0)}}}\right) \ud t.
$$%
  By combining the two previous expressions, we finally get
$$
 V_\ve(f, g)=  \int_T f(t)  \left\{{ \frac{\theta_d^{-1}}{\sqrt{\mu_\ve}} g\!\left({X_\ve(t)}/{\sqrt{r_\ve(0)}}\right) \normp[d]{\nabla\!X_\ve(t)}- g\!\left({{X(t)}/{\sqrt{r(0)}}}\right)
  }\right\} \ud t.
$$
  Thus 
\begin{align*}
 V_\ve(f, g)
 	& =  \int_T f(t)  g\!\left({X_\ve(t)}/{\sqrt{r_\ve(0)}}\right)  h\!\left({\nabla\!X_\ve(t)}/{\sqrt{\mu_\ve}}\right) 
   \ud t \\
	& \quad + \int_T f(t) \left\{{g\!\left({X_\ve(t)}/{\sqrt{r_\ve(0)}}\right) - g\!\left({{X(t)}/{\sqrt{r(0)}}}\right)
  }\right\} \ud t\\
	&:=S_1 +S_2.
\end{align*}
  Thus, we can split $  \Xi_\ve(f, g)$ into two terms
   \begin{align*}
  \Xi_\ve(f, g)
  &\phantom{:}= \ve^{-a(\alpha)} [\ln(1/\ve)]^{-b(\alpha)} S_1 +  \ve^{-a(\alpha)}  [\ln(1/\ve)]^{-b(\alpha)} S_2\\
  &:=T_1+T_2.
   \end{align*}
  The proof will proceed as follows.

We are going to prove that $\bbE[S_1^2]=  O(\ve^{4(1-\alpha)}+  \ve^2 \ln(1/\ve))$ and that $\bbE[S_2^2]= O(\ve^{4\alpha}+ \ve^2)$.\\
  So we can already see that in case where $\alpha \neq \frac{1}{2}$, the term $S_1$ only matters when $\alpha > {1}/{2}$ and we will show in this case that $\bbE[T_1^2] \to \ell_2$.
Similarly, the term $S_2$ will  only play a role when $\alpha < \frac{1}{2}$ and in this case we will prove that $\bbE[T_2^2] \to \ell_1$.
The case where $\alpha = {1}/{2}$ is a little more delicate and we will treat it each time separately.\\
  First let us consider the term $S_2$.\\
  Let us prove that if $\alpha < \frac{1}{2}$,
$$
    \ve^{-4\alpha} \bbE\!\left[S_2^2\right] \to  \ell_1,
$$
    while if $\alpha = \frac{1}{2}$,
$$
    \ve^{-2} \bbE\!\left[S_2^2\right] \to  \ell_1 + \ell_3.
$$
Using the Hermite expansion of $g$, we can write the term $S_2$ as
$$
       S_2=  \sum_{n=1}^{\infty} c_n  \int_T f(t)    \left\{{H_n\!\left({X_\ve(t)}/{\sqrt{r_\ve(0)}}\right) - H_n\!\left({{X(t)}/{\sqrt{r(0)}}}\right)}\right\} \ud t.
$$
We split  $S_2$ into two terms:
\begin{align*}
       S_2&\phantom{:}=  \sum_{n=1}^{\infty} c_n  \int_T f(t)    \left\{{
       1- \!\left({\sqrt{\frac{r_\ve(0)}{r(0)}}
       }\right)^n
  }\right\} H_n\!\left({X_\ve(t)}/{\sqrt{r_\ve(0)}}\right) \ud t +\\
 & \qquad \sum_{n=1}^{\infty} c_n  \int_T f(t)  \left\{{
   \!\left({\sqrt{\frac{r_\ve(0)}{r(0)}}
       }\right)^n  H_n\!\left({X_\ve(t)}/{\sqrt{r_\ve(0)}}\right)  -  H_n\!\left({{X(t)}/{\sqrt{r(0)}}}\right)
   }\right\} \ud t\\
  & := U_2 +V_2.
         \end{align*}
We first study the asymptotic variance of $U_2$.
Using the Mehler's formula \cite{MR716933}, we obtain
\begin{align*}
      \bbE\!\left[U_2^2\right]&= 
       \sum_{n=1}^{\infty} n! c_n^2  \int_{T \times T} f(t) f(s) \!\left({
       1- \!\left({\sqrt{\frac{r_\ve(0)}{r(0)}}
       }\right)^n
  }\right)^2  \!\left({\frac{r_\ve(t-s)}{r_\ve(0)} }\right)^{n}  \ud t  \ud s\\
   &=\sum_{n=1}^{\infty} n! c_n^2 \!\left({
       1- \sqrt{\frac{r_\ve(0)}{r(0)}}
       }\right)^2\!\left({\sum_{k=0}^{n-1}
       \!\left({\frac{r_\ve(0)}{r(0)} }\right)^{\frac{k}{2}}
        }\right)^2  \\
       &\specialpos{\hfill \times
        \int_{T \times T} f(t) f(s) \!\left({\frac{r_\ve(t-s)}{r_\ve(0)} }\right)^{n}  \ud t  \ud s.}
\end{align*}

Let us state and prove a lemma.
         \spacebefore
\begin{lemm}
\label{norm}
$\ve^{-2\alpha}\!\left({r(0)-r_\ve(0)}\right) \to \bC _{\alpha}$ as $\ve$ goes to zero.
\end{lemm}
\spacebefore
\begin{proofarg}{Proof of Lemma \ref{norm}}
$$
r_\ve(0)=\int_{\reels^{d}} \Phi(w)  r(\ve w)  \ud w.$$

Using the form of the covariance $r$ and the fact that $\int_{\reels^d} \Phi(t) \ud t=1$, we easily obtain
$$
\ve^{-2\alpha}\!\left({r(0)-r_\ve(0)}\right) =\int_{\reels^{d}}   \Phi(w) L(\ve \normp[d]{w}) \normp[d]{w}^{2\alpha} \ud w.
$$
Convergence follows from that of $L$, \ie from the fact that $\lim_{x \to 0^{+}} L(x)=\bC _L >0$.
\end{proofarg}
By writing $1-\sqrt{{r_\ve(0)}/{r(0)}}$ as
$$
1-\sqrt{
\frac{r_\ve(0)}{r(0)}
}= \frac{r(0)-r_\ve(0)}
{
(\sqrt{r_\ve(0)}+\sqrt{r(0)})\sqrt{r(0)}
},
$$
 and by using Lemma \ref{norm}, we get on the one hand that
  \begin{eqnarray*}
  \ve^{-4\alpha} \!\left({
 1-\sqrt{
\frac{r_\ve(0)}{r(0)}
}
  }\right)^2 \cvg[\ve\to 0]{} \ell_{\alpha}^2.
  \end{eqnarray*}
On the other hand,
$$
\!\left({\sum_{k=0}^{n-1}
       \!\left({\frac{r_\ve(0)}{r(0)} }\right)^{\frac{k}{2}}
        }\right)^2 \le n^2\quad\text{and}\quad \sum_{n=1}^{\infty} n! n^2  c_n^2 < \infty.
$$
These arguments plus the fact that $f$ is a bounded function will be used to apply Lebesgue's dominated convergence theorem giving the following convergence: $\forall \alpha$ in $\introo{0}{1}$, we have
         \begin{eqnarray}
         \label{l1}
    \ve^{-4\alpha} \bbE\!\left[U_2^2\right] \cvg[\ve\to 0]{}  \ell_1.
      \end{eqnarray}
We are now interested in the asymptotic variance of $V_2$.
According to Mehler's formula
\begin{multline*}
       \bbE\!\left[V_2^2\right]=  \sum_{n=1}^{\infty} n! c_n^2  \int_{T \times T} f(t) f(s)  \\
         \qquad\times\left[{
        \!\left({\frac{r_\ve(t-s)}{r(0)} }\right)^{n} -2 
         \!\left({\frac{r^{(\ve)}(t-s)}{r(0)} }\right)^{n} +
          \!\left({\frac{r(t-s)}{r(0)} }\right)^{n}
       }\right] \ud  t \ud  s.
        \end{multline*}
We divide the integration domain $T \times T$ into two parts, namely $T_\ve^{(1)}$ and $T_\ve^{(2)}$, defined by:
\begin{align*}
 T_\ve^{(1)}&:=\{(s, t) \in T \times T; \normp[d]{t-s} \le M \ve\}\\
 \shortintertext{and}
 T_\ve^{(2)}&:=\{(s, t) \in T \times T; \normp[d]{t-s} > M\ve \}
\end{align*}
where $M$  will be chosen later (in Lemmas \ref{majoration pour r(epsilon)} and \ref{majoration pour derivees r}).\\
Let us denote by $A_\ve$ the term corresponding to the domain $T_\ve^{(1)}$.
Our attention is focused on finding an upper bound for this term.
Since $f$ is bounded, we have
\begin{multline*}
\abs{A_\ve} \le \bC  \sum_{n=1}^{\infty} n! c_n^2\int_{T_\ve^{(1)}}
 \!\left({
 \abs{
  \!\left({\frac{r_\ve(t-s)}{r(0)} }\right)^{n} -
    \!\left({\frac{r(t-s)}{r(0)} }\right)^{n}
 } }\right.\\
 +\left.{
 \abs{
  \!\left({\frac{r^{(\ve)}(t-s)}{r(0)} }\right)^{n} -
    \!\left({\frac{r(t-s)}{r(0)} }\right)^{n}
 } 
 }\right) \ud t \ud  s.
\end{multline*}
Using the fact that for all $s$ and  $t \in T$, 
$$
	\abs{\frac{r(t-s)}{r(0)}} \le 1,\quad \abs{\frac{r_\ve(t-s)}{r(0)}} \le 1\quad \text{and}\quad
	\abs{\frac{r^{(\ve)}(t-s)}{r(0)}} \le 1,
$$
we easily obtain the following bound
\begin{multline*}
\abs{A_\ve} \le 
\bC  \!\left({\sum_{n=1}^{\infty} n! n c_n^2}\right) \int_{T_\ve^{(1)}}
 \!\left({
 \abs{\rule{0pt}{12pt}
  r_\ve(t-s) -
    r(t-s) 
 }}\right.\\ + \left.{
 \abs{
  r^{(\ve)}(t-s) -
   r(t-s)
           } 
 }\right) \ud t \ud  s.
\end{multline*}

At this stage of the proof, we need a lemma.
\spacebefore
\begin{lemm}
\label{majoration pour r(epsilon)}
$\exists \bC  >0$, $\exists  \ve_0 >0$, $\exists M >0$, $\forall \ve \le \ve_0$, $\forall (s, t) \in T \times T$, we have
\begin{multline*}
\sup\left\{{
\abs{r^{(\ve)}(s-t) -r(s-t)};
\abs{r_\ve(s-t) -r(s-t)}
}\right\} \\
\le \bC  \left\{{\ve^{2\alpha} \1_{\{\normp[d]{s-t} \le M \ve\}}  }\right.
\left.{ + \ve^2 \normp[d]{s-t}^{2\alpha-2} \1_{\{\normp[d]{s-t} > M \ve\}}
}\right\}.
\end{multline*}
  \end{lemm}
Let us suppose for a moment that the lemma is proved.
In this vein we obtain
$$
      \abs{A_\ve} \le \bC  \sigma_d(T_\ve^{(1)}) \!\left({\sum_{n=1}^{\infty} n!  n c_n^2}\right) \ve^{2\alpha}.
$$
Using inequality (\ref{finitude cn}), we finally get
 \begin{eqnarray*}
 \abs{A_\ve} \le \bC  \ve^{d+2\alpha}.
 \end{eqnarray*}
 Thus, we have shown that $A_\ve=o(\ve^{4\alpha}+\ve^{2})$.\\
 Let us proceed to the proof of Lemma \ref{majoration pour r(epsilon)}.
 \spacebefore%
 \begin{proofarg}{Proof of Lemma \ref{majoration pour r(epsilon)}}
For $\tau :=s-t$, while $s, t \in T$, the difference between the two terms $r^{(\ve)}(\tau)$ and $r(\tau)$ is expressed and using that function $\Psi$ is a density, we obtain
\begin{multline*}
  r^{(\ve)}(\tau)-r(\tau)= \int_{\reels^d} \Psi(v) \left[{r(\tau-\ve v) - r(\tau) }\right] \ud v\\
  =  \int_{\reels^d} \Psi(v) \left[{L(\normp[d]{\tau})\normp[d]{\tau}^{2\alpha}-
  L(\normp[d]{\tau - \ve v})\normp[d]{\tau -\ve v}^{2\alpha}
  }\right] \ud v.
  \end{multline*}
 
Since $\Psi$ has compact support in $\reels^d$, it is always possible to choose a number $N >0$ such that if $\normp[d]{v} > N$ then $\Psi(v)=0$.
Let $M=4N$.\\
  Since the function $L$ is bounded on a compact set, obviously, if $\normp[d]{\tau} \le M \ve$, then $\abs{r^{(\ve)}(\tau) -r(\tau)} \le \bC  \ve^{2\alpha}$.\\
Let us now choose $\tau$ such that $\normp[d]{\tau} > M \ve$ and $v$ such that $\normp[d]{v} \le N$, then $\normp[d]{\tau - \ve v} \ge 3N\ve$.
Thus we can make a second order expansion of $r(\tau - \ve v)$ in the neighbourhood of $\tau$.
Using the fact that $\Psi$ is a density function and in addition is an even function as depending on the norm, we can express $r^{(\ve)}(\tau)-r(\tau)$ as
  \begin{multline}
  \label{difference between cov}
  r^{(\ve)}(\tau)-r(\tau)=\\
    \frac{1}{2} \int_{\{\normp[d]{v} \le N\}} \Psi(v) \times 
    \ve^2 \left\{{
   \sum_{i=1}^{d}\sum_{j=1}^{d} v_i v_j \dfrac{\partial ^{2}r}{\partial \tau_{i}
   \partial \tau_{j}}(\tau -\theta \ve v)
    }\right \} \ud v,
  \end{multline}
  with $0\le \theta <1$ depending on $\ve$, $\tau$ and $v$.
Now for $x \neq 0$, we can express
  $\dfrac{\partial ^{2}r}{\partial x_{i}^2}(x)$, for $i=1, \dots, d$, as
  \begin{multline*}
  -\dfrac{\partial ^{2}r}{\partial x_{i}^2}(x)= 2\alpha \normp[d]{x}^{2\alpha-4}L(\normp[d]{x}) \!\left({
  \normp[d]{x}^2+2(\alpha -1)x_i^2
  }\right)\\
   +\normp[d]{x}^{2\alpha-3}L^{\prime}(\normp[d]{x})
  \!\left({(4\alpha -1)x_i^2 +\normp[d]{x}^{2}
  }\right) + x_i^2 \normp[d]{x}^{2\alpha -2} L^{\prime \prime}(\normp[d]{x}).
  \end{multline*}
  Since $L$ has two continuous derivatives except at the origin, which satisfies\linebreak
   $\abs{x} L^{\prime}(\abs{x})=O(1)$ and $x^2 L^{\prime\prime}(\abs{x})=O(1)$ as $x \to 0$ and  $\lim_{x \to 0^{+}} L(x)=\bC _L$, we can conclude that $\dfrac{\partial ^{2}r}{\partial x_{i}^2}(x)=O(\normp[d]{x}^{2 \alpha -2})$, as soon as $x$ lives in a compact set.
We have thus proved that, for $\tau$, $v$, $\theta$ and $\ve$ such that $\normp[d]{\tau} > M \ve$, $\normp[d]{v} \le N$, $0\le \theta <1 \text{~and~} \ve \le \ve_0$,
$$
  \dfrac{\partial ^{2}r}{\partial x_{i}^2}(\tau -\theta \ve v)= O(\normp[d]{\tau -\theta \ve v}^{2 \alpha -2})=O(\normp[d]{\tau}^{2 \alpha -2}),
$$%
  since $\normp[d]{\tau}^{2 -2\alpha}
  \le 2^{1-\alpha} \normp[d]{\tau -\theta \ve v}^{2 -2 \alpha}$ as soon as $\normp[d]{\tau} > 4N \ve =M \ve$.
We can argue in a similar way for the other terms $\dfrac{\partial ^{2}r}{\partial \tau_{i}\partial \tau_{j}}(\tau -\theta \ve v)$, for  $i, j=1, \dots, d$ and $i\neq j$.
 Finally, using that $\int_{\reels^d} \Psi(v) \normp[d]{v}^2 \ud v < \infty$, we have proved that for $\normp[d]{\tau} > M \ve$ and $\ve \le \ve_0$,
$\abs{r^{(\ve)}(\tau)-r(\tau)} \le \bC  \ve^2 \normp[d]{\tau}^{2\alpha-2}$.\\
Thus the lemma follows for the term $ r^{(\ve)}(\tau)-r(\tau)$.\\
Now to finish the proof of the lemma, we consider the second difference term $r_\ve(\tau) -r(\tau)$.
In the same way as before, we can write this term as
   \begin{align*}
  r_\ve(\tau) -r(\tau)&=
  \int_{\reels^d} \Phi(v) \left[{r(\tau-\ve v) - r(\tau) }\right] \ud v \\
	&=   \int_{\reels^d} \  \Phi(v)  \left[{L(\normp[d]{\tau})\normp[d]{\tau}^{2\alpha}-
		  L(\normp[d]{\tau - \ve v })\normp[d]{\tau -\ve v }^{2\alpha}
  }\right] \ud v.
\end{align*}
The same kind of arguments could be considered by replacing the $\Psi$ function with the $\Phi$ function and would give the expected result.\\
This completes the proof of the lemma.
\end{proofarg}%
Returning to our proof of the theorem, we consider the second term $B_\ve$ corresponding to the region $T_\ve^{(2)}$.
\\
 As for the proof concerning the set $T_\ve^{(1)}$, we bound $B_\ve$ as follows.
 \begin{multline*}
      \abs{B_\ve} \le 
      \bC  \!\left({\sum_{n=1}^{\infty} n!  n c_n^2}\right) \int_{T_\ve^{(2)}}
       \!\left({
       \abs{\rule{0pt}{12pt}
        r_\ve(t-s) -
          r(t-s) 
       }}\right.\\ + 
       \left.{\abs{
        r^{(\ve)}(t-s) -
         r(t-s)
                 } 
       }\right) \ud  t \ud s.
      \end{multline*}
The second part of Lemma \ref{majoration pour r(epsilon)} yields the bound
$$
      \abs{B_\ve} \le \bC   \ve^2 \int_{T \times T} \normp[d]{t-s}^{2\alpha-2}
      \ud  t \ud  s.
$$%
       And since $d-2+2\alpha >0$, one finally proved that
$$
        \abs{B_\ve} \le \bC   \ve^2.
$$%
We have shown that $\bbE[S_2^2]= O(\ve^{4\alpha}+ \ve^2)$.
Using the convergence established in (\ref{l1}) it turns out that if $\alpha < \frac{1}{2}$,
   \begin{eqnarray*}
    \ve^{-4\alpha} \bbE[S_2^2] \cvg[\ve\to 0]{}  \ell_1.
      \end{eqnarray*}
Now consider the case where $\alpha = {1}/{2}$ and prove that
       \begin{eqnarray*}
    \ve^{-2} \bbE[S_2^2] \cvg[\ve\to 0]{} \ell_1+ \ell_3.
      \end{eqnarray*}
We have already proven that 
\begin{eqnarray*}
    \ve^{-2} \bbE\!\left[U_2^2\right] \cvg[\ve\to 0]{}  \ell_1.
\end{eqnarray*}
Thus let us show that
\begin{eqnarray*}
    \ve^{-2} \bbE\!\left[V_2^2\right] \cvg[\ve\to 0]{}  \ell_3,
      \end{eqnarray*}
      and that 
      \begin{eqnarray*}
      \ve^{-2} \bbE[U_2V_2] \cvg[\ve\to 0]{} 0.
      \end{eqnarray*}
      We first compute $\bbE\!\left[V_2^2\right]$.
As in the previous part, we split $\bbE\!\left[V_2^2\right]$ into two parts, $A_\ve$ and $B_\ve$.
We saw that $A_\ve=O(\ve^{d+2\alpha})=o(\ve^2)$.
Let us examine the second term $B_\ve$ and show that $
    \ve^{-2} B_\ve \cvg[\ve \to 0]{}  \ell_3$.
We have
       \begin{multline*}
       B_\ve=  \sum_{n=1}^{\infty} n! c_n^2  \int_{T_\ve^{(2)}} f(t) f(s)  \\
        \times \left[{
        \!\left({\frac{r_\ve(t-s)}{r(0)} }\right)^{n} -2 
         \!\left({\frac{r^{(\ve)}(t-s)}{r(0)} }\right)^{n} +
          \!\left({\frac{r(t-s)}{r(0)} }\right)^{n}
       }\right] \ud  t \ud  s.
        \end{multline*}
Using equality (\ref{difference between cov}), we obtain
\begin{multline*}
 	B_\ve=  \sum_{n=1}^{\infty} n! c_n^2  \int_{T_\ve^{(2)}} f(t) f(s) 
 	 \!\left({\frac{1}{r(0)}
 	 }\right)^n \\
 	 \times\left\{{  
 	 \!\left({
 	 r(t-s) +\frac{1}{2} \int_{\{\normp[d]{v} \le N\}} \Phi(v) \times 
 	   \ve^2 \!\left({
 	  \sum_{i=1}^{d}\sum_{j=1}^{d} v_i v_j \dfrac{\partial ^{2}r}{\partial t_{i}
	   \partial t_{j}}(t-s -\theta \ve v)
	    }\right ) \ud v
	     }\right)^n }\right.\\
	     \left.{-2
	      \!\left({
	  r(t-s) +\frac{1}{2} \int_{\{\normp[d]{v} \le N\}} \Psi(v) \times 
	    \ve^2 \!\left({
	   \sum_{i=1}^{d}\sum_{j=1}^{d} v_i v_j \dfrac{\partial ^{2}r}{\partial t_{i}
	   \partial t_{j}}(t-s -\theta \ve v)
	    }\right ) \ud v
	     }\right)^n }\right.\\
	     \left.{+ r^n(t-s)
	    }\rule{0pt}{24pt}\right\} \ud t \ud s.
 \end{multline*}
We expand the terms into each parenthesis using the binomial expansion.
We get
 \begin{align*}
  \MoveEqLeft[0]{B_\ve=  \sum_{n=1}^{\infty} n! c_n^2  \int_{T_\ve^{(2)}} f(t) f(s) 
  \!\left({\frac{1}{r(0)}
  }\right)^n} \\
  &    \specialpos{\hfill\quad\times\sum_{k=1}^n 
 {n\choose{k}}
  \!\left({\frac{\ve^2}{2}}\right)^k r^{n-k}(t-s) 
  \left\{{
  \!\left({
  \sum_{i,j=1}^d \int_{\reels^d} \Phi(v)  v_i v_j \dfrac{\partial ^{2}r}{\partial t_{i}
   \partial t_{j}}(t-s -\theta \ve v)\ud v
  }\right)^k }\right.}\\
  &\left.{
  \quad -2 \!\left({
  \sum_{i,j=1}^d \int_{\reels^d} \Psi(v)  v_i v_j \dfrac{\partial ^{2}r}{\partial t_{i}
   \partial t_{j}}(t-s -\theta \ve v)\ud v
  }\right)^k 
  }\right\} \ud v \ud t \ud s\\
 & \quad:= B_\ve^{(1)}+ B_\ve^{(2)},
   \end{align*}
   where the term $B_\ve^{(1)}$ corresponds to the term obtained for $k=1$ and 
   the term $B_\ve^{(2)}$ stands for the remaining ones ($k \ge 2$).\\
   Let us deal with the first term $B_\ve^{(1)}$ which will give the limit.
    \begin{multline*}
  \frac{1}{\ve^2} B_\ve^{(1)}=  \frac{1}{2r(0)} \sum_{n=1}^{\infty} n n!  c_n^2  \int_{T_\ve^{(2)}} f(t) f(s) 
  \!\left({\frac{r(t-s)}{r(0)}
  }\right)^{n-1}  \\
\times    \!\left({
  \sum_{i,j=1}^d \int_{\reels^d} \!\left({\Phi(v)-2\Psi(v)}\right)   v_i v_j \dfrac{\partial ^{2}r}{\partial t_{i}
   \partial t_{j}}(t-s -\theta \ve v)\ud v
  }\right) \ud t \ud s.
   \end{multline*}
Let us see that we can apply Lebesgue's dominated convergence theorem.
At this point in the proof, we state a lemma whose proof was demonstrated in that of Lemma \ref{majoration pour r(epsilon)} and which is valid not only for $\alpha= \frac{1}{2}$ but also for all ranges of $0 < \alpha <1$.
    \spacebefore
\begin{lemm}
     \label{majoration pour derivees r}
$\exists \bC  >0$, $\exists  \ve_0 >0$, $\exists M >0$, $\forall \ve \le \ve_0$, $\forall (s, t) \in T \times T$, such that $\normp[d]{t-s} \ge M \ve$ and for all $v \in \reels^d$ such that $\normp[d]{v} \le \frac{M}{4}$, then for all $i, j =1,\dots,d$, we have
$$\abs{\dfrac{\partial ^{2}r}{\partial t_{i}
	\partial t_{j}}(t-s -\theta \ve v)} \le \bC  \normp[d]{t-s}^{2\alpha-2},
$$%
   for all $0 \le \theta <1$.
  \end{lemm}
Thus, taking $\alpha = {1}/{2}$ and using (\ref{finitude cn}) knowing that $f$ is bounded, we can bound  the integrand appearing in the expression of $B_\ve^{(1)}/{\ve^2} $ by 
 $
 \bC  \normp[d]{t-s}^{-1}$.
\\
 Moreover, since $d-1 >0$, we have
\begin{equation}
	 \int_{T \times T} \normp[d]{t-s}^{-1}\ud  t \ud  s < \infty.
	 \label{integration norm}
 \end{equation}
All the ingredients are gathered to justify the exchange of the limit with the series and the integral, which gives
     \begin{multline*}
     \lim_{\ve \to 0} \frac{1}{\ve^2} B_\ve^{(1)}=
     \frac{1}{2r(0)} \sum_{n=1}^{\infty} n n! c_n^2  \int_{T \times T} f(t) f(s) 
  \!\left({\frac{r(t-s)}{r(0)}
  }\right)^{n-1} \\
   \times\!\left({
  \sum_{i,j=1}^d \dfrac{\partial ^{2}r}{\partial t_{i}
   \partial t_{j}}(t-s)  \int_{\reels^d} \!\left({\Phi(v)-2\Psi(v)}\right)   v_i v_j \ud v
  }\right) \ud t \ud s.
\end{multline*}
To finish with this term, noting that the functions $\Phi$ and $\Psi$ depend only on the norm, we finally express the last convergence as
$$
     \lim_{\ve \to 0} \frac{1}{\ve^2} B_\ve^{(1)}=\ell_3.
$$
Let us show that $\frac{1}{\ve^2} B_\ve^{(2)}\cvg[\ve \to 0]{}0$.

Using the previous lemma, we bound the terms
\begin{multline*}
   \left(\frac{1}{r(0)}\right)^k \times \left[{  \!\left({
  \sum_{i,j=1}^d \int_{\reels^d} \Phi(v)  v_i v_j \dfrac{\partial ^{2}r}{\partial t_{i}
   \partial t_{j}}(t-s -\theta \ve v)\ud v
  }\right)^k }\right.
\\
  \left.{-~2 \!\left({
  \sum_{i,j=1}^d \int_{\reels^d} \Psi(v)  v_i v_j \dfrac{\partial ^{2}r}{\partial t_{i}
   \partial t_{j}}(t-s -\theta \ve v)\ud v
  }\right)^k}\right]
  \end{multline*}
  by 
  $ \!\left({\bC  \normp[d]{t-s}^{-1}}\right)^k$ which is also bounded by $ \!\left({\bC  \normp[d]{t-s}^{-1}}\right) \times \!\left({\frac{\bC }{M\ve}}\right)^{k-1}$.\\
  In that way, we obtain the following bound.
  \begin{multline*}
\abs{B_\ve^{(2)}} \le
     \bC  \sum_{n=2}^{\infty} n!  c_n^2  \sum_{k=2}^n 
 {{n}\choose{k}}
  \!\left({\frac{\ve^2}{2}}\right)^k \!\left({\frac{\bC }{M\ve}}\right)^{k-1}\\
  \times \int_{T_\ve^{(2)}} \normp[d]{t-s}^{-1} \abs{\frac{r(t-s)}{r(0)}}^{n-k} 
  \ud t \ud s.
  \end{multline*}
  Thus 
  \begin{multline*}
\frac{1}{\ve^2}\abs{B_\ve^{(2)}} \le
     \bC  \sum_{n=2}^{\infty} n!  c_n^2  \sum_{k=2}^n 
 {{n}\choose{k}}
  \!\left({\frac{\bC }{M} \ve}\right)^{k-1} \\
\times  \int_{T_\ve^{(2)}} \normp[d]{t-s}^{-1} \abs{\frac{r(t-s)}{r(0)}}^{n-k} 
  \ud t \ud s.
  \end{multline*}
By making the change of variable $k-2= \ell$ we obtain
   \begin{multline*}
\frac{1}{\ve^2}\abs{B_\ve^{(2)}} \le
     \bC  \ve  \sum_{n=2}^{\infty} n!  c_n^2  \sum_{\ell=0}^{n-2} 
{{n}\choose{\ell+2}}
  \!\left({\frac{\bC }{M} \ve}\right)^{\ell} \\
\times  \int_{T_\ve^{(2)}} \normp[d]{t-s}^{-1} \abs{\frac{r(t-s)}{r(0)}}^{n-2-\ell} 
  \ud t \ud s.
  \end{multline*}
  Now for $n \ge 2$ we use the bound,
$$
	{{n}\choose{\ell+2}} \le n (n-1){{n-2}\choose{\ell}}
$$%
getting
\begin{align*}
\frac{1}{\ve^2}\abs{B_\ve^{(2)}} &\le
     \bC  \ve  \sum_{n=2}^{\infty} n (n-1) n!  c_n^2  \sum_{\ell=0}^{n-2} 
 	{{n-2}\choose{\ell}}
	\!\left({\frac{\bC }{M} \ve}\right)^{\ell} \\
 &\quad\times \int_{T_\ve^{(2)}} \normp[d]{t-s}^{-1} \abs{\frac{r(t-s)}{r(0)}}^{n-2-\ell} 
  \ud t \ud s\\
  &= \bC  \ve  \sum_{n=2}^{\infty} n (n-1) n!  c_n^2  \\
 &\quad\times \int_{T_\ve^{(2)}} \normp[d]{t-s}^{-1} \!\left({
  \abs{\frac{r(t-s)}{r(0)}} + \frac{\bC }{M} \ve
  }\right)^{n-2} 
  \ud t \ud s.
  \end{align*}
At this stage of the demonstration, we have to prove a lemma.
\spacebefore
\begin{lemm}
  \label{covariance 1}
  For $\ve$ small enough and $M$ large enough,
$$
  \abs{\frac{r(t-s)}{r(0)}} + \frac{\bC }{M} \ve <1,
$$
   for all $s, t \in T_\ve^{(2)}$.
 \end{lemm}
Let us assume for a moment that the lemma is proved.
Using (\ref{finitude cn}) and (\ref{integration norm}), we obtain the bound
$$
	\frac{1}{\ve^2}\abs{B_\ve^{(2)}}
		\le \bC  \ve.
$$%
To conclude, we have proved that in case $\alpha = \frac{1}{2}$,
$$
    \ve^{-2} \bbE\!\left[V_2^2\right] \cvg[\ve\to0]{}  \ell_3.
$$ 
Before going further in the proof, let us show Lemma \ref{covariance 1}.
 \spacebefore%
 \begin{proofarg}{Proof of Lemma \ref{covariance 1}}
 On the one hand, let us notice that for $\delta >0$ small enough and $s, t \in T_\ve^{(2)}$, such that $\normp[d]{t-s} \le \delta$, we have
 \begin{align*}
	\abs{\frac{r(t-s)}{r(0)}}&= \frac{r(t-s)}{r(0)},\\ \shortintertext{and}
	1 -\frac{r(t-s)}{r(0)}&= \frac{L(\normp[d]{t-s})}{r(0)} \normp[d]{t-s}\\
	& \ge \frac{M \widetilde{C}_L}{r(0)} \ve,
\end{align*}
 where we defined
 $\widetilde{C}_L:= \inf_{\{s, t \in T\}}\{L(\normp[d]{t-s})\} >0$.
Thus for $s, t \in T_\ve^{(2)}$, such that $\normp[d]{t-s} \le \delta$ 
$$
 1 - \frac{r(t-s)}{r(0)} - \frac{\bC }{M} \ve \ge M \ve \!\left({\frac{\widetilde{C}_L}{r(0)} - \frac{\bC }{M^2}}\right)>0,
$$%
 for $M$ sufficiently large.\\
 On the other hand, one can choose $\ve$ sufficiently small so that
 \begin{eqnarray*}
  \sup_{\{s, t \in T, \normp[d]{t-s} \ge \delta \}}\left\{ \abs{\frac{r(t-s)}{r(0)}} + \frac{\bC }{M} \ve\right\} \le \eta <1.
  \end{eqnarray*}
This completes the proof of the lemma.
 \end{proofarg}%
 Let us now prove that in the case $\alpha = {1}/{2}$,
$$
      \ve^{-2} \bbE[U_2V_2] \cvg[\ve\to0]{} 0.
$$%
Using the Mehler's formula once again, we obtain
           \begin{multline*}
       \frac{1}{\ve^2} \bbE[U_2 V_2]=  \sum_{n=1}^{\infty} n! c_n^2  \int_{T \times T} f(t) f(s) \frac{1}\ve\!\left({
       1- \!\left({\sqrt{\frac{r_\ve(0)}{r(0)}}
       }\right)^n
  }\right)  \\
   \times\frac{1}\ve \left[{
       \!\left({\frac{r_\ve(t-s)}{\sqrt{r_\ve(0)}\sqrt{r(0)}} }\right)^{n} - 
         \!\left({\frac{r^{(\ve)}(t-s)}{\sqrt{r_\ve(0)}\sqrt{r(0)}} }\right)^{n} 
       }\right]   \ud t  \ud s.
       \end{multline*}
Working in the same way as for the terms $U_2$, $A_\ve$ and $B_\ve$, we easily obtain
\begin{multline*}
       \frac{1}{\ve^2}\abs{\bbE[U_2 V_2]} \le \bC  \!\left({\sum_{n=1}^{\infty} n^2  n! c_n^2}\right)  
        \left\{{\frac{1}\ve \!\left({
       1- \sqrt{\frac{r_\ve(0)}{r(0)}}
       }\right)
       }\right\}  \\
       \times\!\left({\frac{1}\ve \int_{T \times T}
       \abs{r_\ve(t-s) -r^{(\ve)}(t-s)}  \ud  t \ud  s}\right).
         \end{multline*}
Using Lemmas \ref{norm} and \ref{majoration pour r(epsilon)}, inequalities (\ref{finitude cn}) and (\ref{integration norm}), we find the following bound
\begin{align*}
       \frac{1}{\ve^2}\abs{\bbE\!\left[U_2 V_2\right]} &\le \bC 
       \int_{T \times T} \left\{{\1_{\{\normp[d]{s-t} \le M \ve\}}  
	 + \ve \normp[d]{s-t}^{-1} \1_{\{\normp[d]{s-t} > M \ve\}}
	}\right\} \ud t \ud s\\
	 &\le \bC  \left\{{\ve^d + \ve }\right\} \le \bC  \ve=o(1).
\end{align*}
For $\alpha = {1}/{2}$ we have proved that
$$
    \ve^{-2} \bbE[S_2^2] \cvg[\ve\to\infty]{} \ell_1+ \ell_3.
$$%
Now, consider the term $S_1$.\\
Our goal is to establish that $\bbE[S_1^2]=  O(\ve^{4(1-\alpha)}+\ve^2 \ln(1/\ve))$.\\
For this purpose, let us prove Proposition \ref{term S1}.
\spacebefore
\begin{prop}
\label{term S1}
\begin{enumerate}
\item\label{itm:termS1.1} If $d-4+4 \alpha >0$, 
$$
	\frac{1}{\ve^{4(1-\alpha)}} \bbE[S_1^2] \cvg[\ve\to0]{} \ell_2.
$$%
\item\label{itm:termS1.2} If $d-4+4 \alpha <0$ or ($d=3$ and $\alpha=\frac{1}{4}$),
$$
	\bbE[S_1^2]=O(\ve^2).
$$%
  \item\label{itm:termS1.3} If $d=2$ and  $\alpha = {1}/{2}$,
$$
\frac{1}{\ve^{2} \ln(1/\ve)} \bbE[S_1^2] \cvg[\ve\to0]{} \ell_5.
$$%
\end{enumerate}
\end{prop}
\spacebefore
\begin{rema}
\label{remark theoreme}
The proposition highlights the fact that the term $S_1$ will prevail when $\alpha > {1}/{2}$ or when  $\alpha = {1}/{2}$  and $d=2$.
In the case where $\alpha > {1}/{2}$,
 $$
 	\frac{1}{\ve^{4(1-\alpha)}} \bbE[S_1^2] \cvg[\ve\to0]{} \ell_2.
$$%
Parts \ref{itm4.6.19.1}, \ref{itm4.6.19.2} and \ref{itm4.6.19.4} of Theorem \ref{vitesse variance temps local} will then follow.
Furthermore, in case where $\alpha = {1}/{2}$ and $d>2$, this proposition implies that 
$$
	\frac{1}{\ve^{2}} \bbE[S_1^2]\cvg[\ve\to0]{} \ell_2.
$$%
\hfill$\bullet$
\end{rema}
\spacebefore
\begin{proofarg}{Proof of Proposition \ref{term S1}}
Let us calculate $\bbE[S_1^2]$.
\begin{multline*}
\bbE[S_1^2]=\int_{T \times T} f(t) f(s)  
\bbE\!\left[g\!\left({X_\ve(t)}/{\sqrt{r_\ve(0)}}\right) 
g\!\left({X_\ve(s)}/{\sqrt{r_\ve(0)}}\right)
\right. \\ \times\left.h\!\left({\nabla\!X_\ve(t)}/{\sqrt{\mu_\ve}}\right)
h\!\left({\nabla\!X_\ve(s)}/{\sqrt{\mu_\ve}}\right)\rule{0pt}{12pt}\right] 
   \ud t \ud s.
\end{multline*}
As in the previous section, we divide the integration domain $T \times T$ into two parts $T_\ve^{(1)}$ and $T_\ve^{(2)}$.
Let $\ell_1(\ve)$ be the integral over $T_\ve^{(1)}$.
Making the change of variable $t-s = \ve w$, we obtain
{\def\termeA{\ve^d
	\int\limits_{\{(w, t) : (t-\ve w, t) \in T \times T, \normp[d]{w} \le M\}} f(t) f(t-\ve w)}
\def\termeB{g\!\left({X_\ve(\ve w)}/{\sqrt{r_\ve(0)}}\right)}
\def\termeC{g\!\left({X_\ve(0)}/{\sqrt{r_\ve(0)}}\right)}
\def\termeD{h\!\left({\nabla\!X_\ve(\ve w)}/{\sqrt{\mu_\ve}}\right)}
\def\termeE{h\!\left({\nabla\!X_\ve(0)}/{\sqrt{\mu_\ve}}\right)}
\newlength{\hauteur}
\setlength{\hauteur}{\heightof{$\ds \bbE\!\left[\termeB\right]$}}
\begin{align*}
\MoveEqLeft[4]\ell_1(\ve)= \termeA \\
& \times  \bbE\!\left[
	\termeB
	\termeC
	\termeD
	\right.\\
& \times	\left. \rule{0pt}{\hauteur}
	\termeE
	\right]\ud t \ud w.
\end{align*}}

Since the function $f$ is assumed to be bounded, we can bound the following expression.
For any $(w, t)$ such that $(t-\ve w, t) \in T \times T$,
{\def\termeA{f(t) f(t-\ve w) }
\def\termeB{g\!\left({X_\ve(\ve w)}/{\sqrt{r_\ve(0)}}\right)}
\def\termeC{g\!\left({X_\ve(0)}/{\sqrt{r_\ve(0)}}\right)}
\def\termeD{h\!\left({\nabla\!X_\ve(\ve w)}/{\sqrt{\mu_\ve}}\right)}
\def\termeE{h\!\left({\nabla\!X_\ve(0)}/{\sqrt{\mu_\ve}}\right)}
\setlength{\hauteur}{\heightof{$\ds \bbE\!\left[\termeB\right]$}}
\begin{align*}
\MoveEqLeft[3]{\left|{\termeA
\bbE\!\left[\termeB 
\termeC\right.}\right.}\\
& \left.\left.\rule{0pt}{12pt}\termeD
\termeE\right]\right|  \\
& \le \bC  \bbE\!\left[
g^2\!\left({X_\ve(0)}/{\sqrt{r_\ve(0)}}\right)h^2\!\left({\nabla\!X_\ve(0)}/{\sqrt{\mu_\ve}}\right)
\right] \\
& = \bC  \bbE\!\left[g^2\!\left({X_\ve(0)}/{\sqrt{r_\ve(0)}}\right)\right]
	\bbE\!\left[h^2\!\left({\nabla\!X_\ve(0)}/{\sqrt{\mu_\ve}}\right)\right] \\
& = \bC  \bbE\!\left[g^2(N_1)\right] \bbE\!\left[h^2(N_d)\right]\\
& < \infty,
\end{align*}}%
the finiteness is justified by the fact that $g$ is $L^2(\varphi(x) \ud x)$.\\
So we have proved that
$$
	\ell_1(\ve)= O(\ve^d).
$$%
We now consider the integral on $T_\ve^{(2)}$ denoted by $\ell_2(\ve)$.
Since the product $h\cdot g$ is $L^2(\reels^{d+1}, \varphi_{d+1}(x) \ud x)$, the following expansion converges in this space.
For $x:=(x_1, x_2, \dots, x_d, x_{d+1}) \in \reels^{d+1}$,
\begin{eqnarray}
\label{function hg}
h(x_1, x_2, \dots,x_d)  g(x_{d+1})= \sum_{q=2}^{\infty} \sum_{
\substack{
\bk  \in \naturels^{d+1} \\ \abs{\bk }=q
}
} b(\bk ) \widetilde{H}_{\bk }(x),
\end{eqnarray}
where for $\bk :=(k_1, k_2, \dots, k_d, k_{d+1}) \in \naturels^{d+1}$, $b(\bk ):= a(k_1, k_2, \dots, k_d) c_{k_{d+1}}$.\\
With these notations, we have
\begin{multline*}
\ell_2(\ve)=\int_{T_\ve^{(2)}} f(t) f(s) 
\sum_{q=2}^{\infty} \sum_{q^{\prime}=2}^{\infty} \sum_{
\substack{
\bk  \in \naturels^{d+1} \\ \abs{\bk }=q
}
} \sum_{
\substack{
\bmm  \in \naturels^{d+1} \\ \abs{\bmm }=q^{\prime}
}
}b(\bk ) b(\bmm )  \\
\times\bbE[\widetilde{H}_{\bk }(U_\ve(t)) \widetilde{H}_{\bmm }(U_\ve(s))] \ud t \ud s,
\end{multline*}
where 
$$
U_\ve(t):=
	\left(
	\dfrac{\partial X_\ve}{\partial t_1}(t)/\sqrt{\mu_\ve}, \dfrac{\partial X_\ve}{\partial t_2}(t)/\sqrt{\mu_\ve}, \dots, \dfrac{\partial X_\ve}{\partial t_{d}}(t)/\sqrt{\mu_\ve}, X_\ve(t)/\sqrt{r_\ve(0)}
	\right).
$$%
First of all, let us note that if $\bk , \bmm  \in \naturels^{d+1}$ are such that $\abs{\bk } \neq \abs{\bmm }$, then for all $s, t \in T \times T$,
$\bbE[\widetilde{H}_{\bk }(U_\ve(t)) \widetilde{H}_{\bmm }(U_\ve(s))]=0$ and the expression above turns out to be a sum of orthogonal terms in $L^2(\Omega)$.

Indeed, to prove this, we need a generalization of Mehler's formula given in Taqqu (\cite{MR471045}, Lemma 3.2) or in Berzin (\cite{MR4347380}, Lemma 3.9) via the following lemma.
\spacebefore
\begin{lemm}
\label{mehler}
Let $X=(X_i)_{i=1,n}$ and $Y=(Y_j)_{j=1,n}$ be two standard Gaussian vectors in $\reels^n$ such that for $1 \le i, j \le n$, $\bbE[X_iY_j]=\rho_{ij}$, then for $\bk , \bmm  \in \naturels^n$, we have
$$
	\bbE[\widetilde{H}_{\bk }(X) \widetilde{H}_{\bmm }(Y)]=
	\Bigg(\sum_{
		\substack{
		a_{ij} \ge 0 \\ \sum_{j}a_{ij}=k_i \\ \sum_{i}a_{ij}=m_j
			}
		} \bk !  \bmm ! \prod_{1 \le i, j \le n} \frac{\rho_{ij}^{a_{ij}}}{a_{ij}!}
	\Bigg)
	\1_{\abs{\bk }=\abs{\bmm }}.
$$%
\end{lemm}
The previous lemma allows us to write
\begin{align*}
\ell_2(\ve)&=\int_{T_\ve^{(2)}} f(t) f(s) 
\sum_{q=2}^{\infty} \sum_{
\substack{
\bk , \bmm  \in \naturels^{d+1} \\ \abs{\bk }=\abs{\bmm }=q
}
} 
b(\bk ) b(\bmm )  
\bbE[\widetilde{H}_{\bk }(U_\ve(t)) \widetilde{H}_{\bmm }(U_\ve(s))] \ud t \ud s\\
&= \int_{T_\ve^{(2)}} 
\sum_{q=2}^{\infty} \sum_{
\substack{
\bk , \bmm  \in \naturels^{d+1} \\ \abs{\bk }=\abs{\bmm }=q
}
} 
f(t) f(s) b(\bk ) b(\bmm )  
\bk !  \bmm ! A_\ve^{(q)}(\bk , \bmm , s, t)  \ud t \ud s,
\end{align*}
where for $q \ge 2$ and $\bk , \bmm  \in \naturels^{d+1}$ such that $\abs{\bk }=\abs{\bmm }=q$ and $(s, t) \in T_\ve^{(2)}$,
$$
	A_\ve^{(q)}(\bk , \bmm , s, t):=
	\sum_{
		\substack{
		a_{i,j} \ge 0 \\ \sum_{j}a_{i,j}=k_i \\ \sum_{i}a_{i,j}=m_j
			}
		} 
	A_\ve^{(q)}(\bk , \bmm , s, t, \ba ),
$$%
where for $\ba  :=(a_{i,j})_{1\le i , j \le d+1}$,
\begin{align*}
A_\ve^{(q)}(\bk , \bmm , s, t, \ba )&:=
\prod_{1 \le i, j \le d} \left[{
\!\left({\frac{-\dfrac{\partial ^{2}r_\ve}{\partial t_{i}
   \partial t_{j}}(t-s)}{\mu_\ve}
}\right)^{a_{i,j}}\frac{1}{a_{i,j}!}
}\right]  \\
&\phantom{:=} \times\prod_{1 \le i \le d} \left[{
\!\left({\frac{\dfrac{\partial r_\ve}{\partial t_{i}}(t-s)}{\sqrt{\mu_\ve}{\sqrt{r_\ve(0)}}}
}\right)^{a_{i,d+1}}\frac{1}{a_{i,d+1}!}
}\right]\\ 
&\phantom{:=} \times \prod_{1 \le j \le d} 
	\left[{
		\!\left({\frac{-\dfrac{\partial r_\ve}{\partial t_{j}}(t-s)}{\sqrt{\mu_\ve}{\sqrt{r_\ve(0)}}}
			}\right)^{a_{d+1,j}}\frac{1}{a_{d+1,j}!}
	}\right] \\
&\phantom{:=} \times \!\left({
					\frac{r_\ve(t-s)}{r_\ve(0)}
				}\right)^{a_{d+1,d+1}}\frac{1}{a_{d+1,d+1}!}.
\end{align*}
Let us see that in case where $ d-4 + 4 \alpha >0$, we can apply Lebesgue's dominated convergence theorem and that
$$
      \frac{1}{\ve^{4(1-\alpha)}} \ell_2(\ve) \cvg[\ve\to0]{}\ell_2.
$$%
Let us fix $q \ge 2$, $\bk , \bmm  \in \naturels^{d+1}$ such that $\abs{\bk }=\abs{\bmm } =q$, $\ba  :=(a_{i,j})_{1\le i , j \le d+1}$ and $(s, t) \in T_\ve^{(2)}$.
Recall that all the indices $k_i$, $i=1, \dots,d$ are even and that $(k_1, k_2, \dots, k_d) \neq (0, 0, \dots, 0)$.%
\\
      We want to find an upper bound for $\ve^{-4(1-\alpha)} A_\ve^{(q)}(\bk , \bmm , s, t, \ba ) $.\\
      We have to consider three cases.
\begin{itemize}
          \item[Case 1.]
$\sum_{i, j=1}^d a_{i, j} \ge 2$.\\
      In this case, $\exists i_0, j_0 \in \{1,\dots, d\}$ such that $a_{i_0, j_0} \ge 2$ or 
      $\exists i_0, j_0, i_1, j_1 \in \{1,\dots, d\}$ such that $(i_0, j_0)\neq (i_1, j_1)$ and such that $a_{i_0, j_0} \ge 1$ and $a_{i_1, j_1} \ge 1$.\\
      Since $\forall i, j =1, \dots,d$, $\abs{\dfrac{\partial ^{2}r_\ve}{\partial t_{i}
   \partial t_{j}}(t-s)/\mu_\ve} \le 1$, an upper bound for
$$\ve^{-4(1-\alpha)}  \prod_{1 \le i, j \le d}\abs{
\dfrac{\partial ^{2}r_\ve}{\partial t_{i}
   \partial t_{j}}(t-s)/\mu_\ve}^{a_{i,j}}
$$%
is
$$ 
\!\left({\dfrac{\partial ^{2}r_\ve}{\partial t_{i_{0}}
   \partial t_{j_{0}}}(t-s)\frac{\ve^{-2(1-\alpha)} }{\mu_\ve}}\right)^{2}
$$%
or 
$$ 
\!\left({\abs{\dfrac{\partial ^{2}r_\ve}{\partial t_{i_{0}}
   \partial t_{j_{0}}}(t-s)}\frac{\ve^{-2(1-\alpha)} }{\mu_\ve}}\right) \times
   \!\left({\abs{\dfrac{\partial ^{2}r_\ve}{\partial t_{i_1}
   \partial t_{j_1}}(t-s)}\frac{\ve^{-2(1-\alpha)} }{\mu_\ve}}\right).
$$%
Now on the one hand,
\begin{equation}
      \dfrac{\partial ^{2}r_\ve}{\partial t_{i_{0}}
   \partial t_{j_{0}}}(t-s)= \int_{\normp[d]{v} \le \frac{M}{4}} \Phi(v)  \dfrac{\partial ^{2} r}{\partial t_{i_{0}}
   \partial t_{j_{0}}}(t-s - \ve v) \ud v,\label{derivees seconde epsilon}
\end{equation}
    and by Lemma \ref{majoration pour derivees r},  we deduce that
\begin{equation}
    \label{majoration derivees secondes r epsilon}
      \abs{\dfrac{\partial ^{2}r_\ve}{\partial t_{i_{0}}
   \partial t_{j_{0}}}(t-s)}\le \bC  \normp[d]{t-s}^{2\alpha -2}.
\end{equation}
On the other hand, we consider the expression of $\mu_\ve$.
\begin{align*}
\mu_\ve&= -\dfrac{\partial ^{2}r_\ve}{\partial t_1^2}(0)\\
&=-\frac{1}{\ve^2}
\int_{\reels^d} \dfrac{\partial ^{2}\Phi}{\partial v_1^2}(v) r(-\ve v) \ud v\\
 &=  \ve^{-2(1-\alpha)}  \int_{\reels^d} \dfrac{\partial ^{2}\Phi}{\partial v_1^2}(v) L(\ve \normp[d]{v}) \normp[d]{v}^{2\alpha} \ud v.
\end{align*}
We deduce that 
\begin{equation}
\frac{\mu_\ve}{\ve^{-2(1-\alpha)}} \cvg[\ve\to0]{} \bC _L \chi_{d}^2(\alpha).
   \label{limite mu epsilon}
\end{equation}
Finally,
$$
\!\left({\abs{\dfrac{\partial ^{2}r_\ve}{\partial t_{i_{0}}
\partial t_{j_{0}}}(t-s)}\frac{\ve^{-2(1-\alpha)} }{\mu_\ve}}\right)
$$
is bounded by $\bC  \normp[d]{t-s}^{2\alpha -2}$.

To conclude with this case,
\begin{multline*}
 \prod_{1 \le i \le d} 
\abs{
\dfrac{\partial r_\ve}{\partial t_{i}}(t-s)/(\sqrt{\mu_\ve}{\sqrt{r_\ve(0)}})}^{a_{i,d+1}}\\ 
\times \prod_{1 \le j \le d} 
\abs{
\dfrac{\partial r_\ve}{\partial t_{j}}(t-s)/(\sqrt{\mu_\ve}{\sqrt{r_\ve(0)}})}^{a_{d+1,j}} \times \abs{
\frac{r_\ve(t-s)}{r_\ve(0)}
}^{a_{d+1,d+1}}
\end{multline*}
 is bounded by one.
 
Wrapping up these inequalities, 
$$
\ve^{-4(1-\alpha)} \abs{A_\ve^{(q)}(\bk , \bmm , s, t, \ba )} \le \bC  \normp[d]{t-s}^{4(\alpha -1)} \prod_{1 \le i, j \le d+1}  \frac{1}{a_{i,j}!}.
$$
  \item[Case 2.]
$\sum_{i, j=1}^d a_{i, j}  =1$.\\
Note that we necessarily have
  $$
  \sum_{i=1}^d a_{i, d+1}  \ge 1\quad \text{and}\quad
  \sum_{j=1}^d a_{d+1, j}  \ge 1.
  $$%
Indeed, otherwise we would have $\sum_{i=1}^d k_i=1$ or $\sum_{j=1}^d m_j=1$, which cannot happen since all indices $k_i$ (resp$.$
$m_j$), $i=1, \dots,d$ are even and there exists at least one index $i$ such that $k_i \neq 0$ (resp$.$
$m_j$).\\
  Thus $\exists i_0, j_0 \in \{1,\dots, d\}$ such that $a_{i_0, j_0} =1$ and $\exists i_1, j_1 \in \{1,\dots, d\}$ such that $a_{i_1, d+1} \ge 1$ and $a_{d+1, j_1} \ge 1$.
  
 As previously
 $$
   \ve^{-2(1-\alpha)}  \prod_{1 \le i, j \le d}\abs{
\dfrac{\partial ^{2}r_\ve}{\partial t_{i}
   \partial t_{j}}(t-s)/\mu_\ve}^{a_{i,j}}
$$
in bounded by 
$$
\!\left({\abs{\dfrac{\partial ^{2}r_\ve}{\partial t_{i_{0}}
   \partial t_{j_{0}}}(t-s)}\frac{\ve^{-2(1-\alpha)} }{\mu_\ve}}\right)
$$%
 and then by $\bC  \normp[d]{t-s}^{2\alpha -2}$.

Similarly,   
\begin{multline*}
 \ve^{-2(1-\alpha)}\prod_{1 \le i \le d} 
   \abs{
   \dfrac{\partial r_\ve}{\partial t_{i}}(t-s)/(\sqrt{\mu_\ve}{\sqrt{r_\ve(0)}})}^{a_{i,d+1}} \\
   \times \prod_{1 \le j \le d} 
   \abs{
   \dfrac{\partial r_\ve}{\partial t_{j}}(t-s)/(\sqrt{\mu_\ve}{\sqrt{r_\ve(0)}})}^{a_{d+1,j}} \times \abs{
\frac{r_\ve(t-s)}{r_\ve(0)}
}^{a_{d+1,d+1}}
\end{multline*}
is bounded by
 $$
 	\!\left({\abs{\dfrac{\partial r_\ve}{\partial t_{i_1}}(t-s)/({\sqrt{r_\ve(0)}})}\sqrt{\frac{\ve^{-2(1-\alpha)} }{\mu_\ve}}}\right) \times
 \!\left({\abs{\dfrac{\partial r_\ve}{\partial t_{j_1}}(t-s)/({\sqrt{r_\ve(0)}})}\sqrt{\frac{\ve^{-2(1-\alpha)} }{\mu_\ve}}}\right).
 $$
Arguments similar to those given in the proof of Lemma \ref{majoration pour derivees r} lead to
$$
      \abs{\dfrac{\partial r_\ve}{\partial t_{i_1}}(t-s)}\le \bC  \normp[d]{t-s}^{2\alpha -1}.
$$%
Moreover, as before, we get the bound
   $\sqrt{\frac{\ve^{-2(1-\alpha)} }{\mu_\ve}} \le \bC $.
Finally, we obtain the bound
\begin{align*}
\ve^{-4(1-\alpha)} \abs{A_\ve^{(q)}(\bk , \bmm , s, t, \ba )} &\le \bC  \normp[d]{t-s}^{6 \alpha -4} \prod_{1 \le i, j \le d+1}  \frac{1}{a_{i,j}!}\\
& \le \bC  \normp[d]{t-s}^{4(\alpha -1)} \prod_{1 \le i, j \le d+1}  \frac{1}{a_{i,j}!}.
\end{align*}
 \item[Case 3.]
$\sum_{i, j=1}^d a_{i, j}  =0$.\\
Note that we necessarily have
$$
\sum_{i=1}^d a_{i, d+1}  \ge 2 \quad \text{and}\quad
  \sum_{j=1}^d a_{d+1, j}  \ge 2.
$$%
Indeed, otherwise we would have $\sum_{i=1}^d k_i=0 \mbox{ or } 1$ or $\sum_{j=1}^d m_j=0 \mbox{ or } 1$, which is impossible.\\
  Thus $\exists i_0 \in \{1,\dots, d\}$ (resp$.$
$\exists j_0 \in \{1,\dots, d\}$) such that $a_{i_0, d+1} \ge 2$ (resp$.$
$a_{d+1, j_0} \ge 2$) or $\exists i_0, i_1 \in \{1,\dots, d\}$ (resp$.$
$\exists j_0, j_1 \in \{1,\dots, d\}$), $i_0 \neq i_1$ (resp$.$
$j_0 \neq j_1$) such that $a_{i_0, d+1} \ge 1$ and $a_{i_1, d+1} \ge 1$ (resp$.$
$a_{d+1, j_0} \ge 1$ and $a_{d+1, j_1} \ge 1$).\\
We proceed as for the other two cases to establish the following bound.
  \begin{align*}
\ve^{-4(1-\alpha)} \abs{A_\ve^{(q)}(\bk , \bmm , s, t, \ba )} &\le \bC  \normp[d]{t-s}^{4(2 \alpha -1)} \prod_{1 \le i, j \le d+1}  \frac{1}{a_{i,j}!}\\
& \le \bC  \normp[d]{t-s}^{4(\alpha -1)} \prod_{1 \le i, j \le d+1}    \frac{1}{a_{i,j}!}.
\end{align*}
 \end{itemize}
 Finally, if we gather the three cases, we have shown that for fixed $q \ge 2$ and $\bk , \bmm  \in \naturels^{d+1}$ such that $\abs{\bk }=\abs{\bmm }=q$ and $(s, t) \in T_\ve^{(2)}$,
$$
\ve^{-4(1-\alpha)} \abs{A_\ve^{(q)}(\bk , \bmm , s, t)}\le \bC  \normp[d]{t-s}^{4(\alpha -1)}
\sum_{
\substack{
a_{i,j} \ge 0 \\ \sum_{j}a_{i,j}=k_i \\ \sum_{i}a_{i,j}=m_j
}
} \prod_{1 \le i, j \le d+1}  \frac{1}{a_{i,j}!},
$$
which provides recalling that $f$ is bounded, the following bound
\begin{multline*}
{\ve^{-4(1-\alpha)} \abs{f(t) f(s) b(\bk ) b(\bmm )  
\bk !  \bmm ! A_\ve^{(q)}(\bk , \bmm , s, t)}}  \\
\le\bC  \normp[d]{t-s}^{4(\alpha -1)} \abs{b(\bk )} \abs{b(\bmm )}  
\bk !  \bmm !  \sum_{
\substack{
a_{i,j} \ge 0 \\ \sum_{j}a_{i,j}=k_i \\ \sum_{i}a_{i,j}=m_j
}
} \prod_{1 \le i, j \le d+1}  \frac{1}{a_{i,j}!}.
\end{multline*}
At this stage of the proof let us state a lemma.
      \spacebefore
\begin{lemm}
      \label{serie finie}
      Let $X=(X_i)_{i=1,n}$ and $Y=(Y_j)_{j=1,n}$ be two standard Gaussian vectors in $\reels^n$ and $G \in L^2(\reels^{n}, \varphi_{n}(x) \ud x)$, with the following expansion in this space:
$$
\text{For $x:=(x_1, \dots, x_n) \in \reels^{n}$,~}G(x_1, \dots, x_n) = \sum_{q=0}^{\infty} \sum_{
\substack{
\bk  \in \naturels^{n} \\ \abs{\bk }=q
}
} g(\bk ) \widetilde{H}_{\bk }(x).
$$
Then
$$
\sum_{q=0}^{\infty} \sum_{
\substack{
\bk , \bmm  \in \naturels^{n} \\ \abs{\bk }=\abs{\bmm }=q
}
} \abs{g(\bk )} \abs{g(\bmm )} \bbE[\widetilde{H}_{\bk }(X)\widetilde{H}_{\bmm }(Y)]< \infty.
$$
In particular, if for $1 \le i, j \le n$, $\bbE[X_iY_j]=1$, then
\begin{equation}
\label{Taqqu}
\sum_{q=0}^{\infty} \sum_{
\substack{
\bk , \bmm  \in \naturels^{n} \\ \abs{\bk }=\abs{\bmm }=q
}
} \abs{g(\bk )} \abs{g(\bmm )} 
\sum_{
\substack{
a_{i,j} \ge 0 \\ \sum_{j}a_{i,j}=k_i \\ \sum_{i}a_{i,j}=m_j
}
} \bk !  \bmm ! \prod_{1 \le i, j \le n} \frac{1}{a_{i,j}!}< \infty.
\end{equation}
\end{lemm}
Recalling that $d-4 + 4 \alpha >0$ and applying this lemma to $n:=d+1$ and to the function $hg$ and thus (\ref{Taqqu}) to coefficients $b(\bk )$ defined by (\ref{function hg}) , we get
\begin{multline}
\!\left({\int_{T \ \times T} 
\normp[d]{t-s}^{4(\alpha -1)} \ud t \ud s}\right) \\
\times\!\left(\rule{0pt}{40pt}\smash{\sum_{q=2}^{\infty} \sum_{
\substack{
\bk , \bmm  \in \naturels^{d+1} \\ \abs{\bk }=\abs{\bmm }=q
}
} \abs{b(\bk )} \abs{b(\bmm )} 
\sum_{
\substack{
a_{i,j} \ge 0 \\ \sum_{j}a_{i,j}=k_i \\ \sum_{i}a_{i,j}=m_j
}
} \bk !  \bmm ! \prod_{1 \le i, j \le d+1} \frac{1}{a_{i,j}!}}\right)< \infty.
\label{graphes}
\end{multline}
We are therefore able to apply Lebesgue's dominated convergence theorem.
\spacebefore
\begin{rema}
Before doing so, notice that we have highlighted the following results:
for $0<\alpha<1$,
$$
	\ell_2(\ve) \le \bC   \ve^{4(1-\alpha)} \int_{T_\ve^{(2)}}
	\normp[d]{t-s}^{4(\alpha -1)} \ud t \ud s.
$$%
Thus if $d-4+4 \alpha <0$, $\ell_2(\ve)=O(\ve^d)=O(\ve^2)$.\\
Moreover, if $d-4+4 \alpha =0$, $\ell_2(\ve)=O(\ve^{4(1-\alpha)} \ln(1/\ve))$ and if $d=3$ (and then $\alpha=\frac{1}{4}$), thus we proved that $\ell_2(\ve)= O(\ve^{3} \ln(1/\ve))=O(\ve^2)$.\\
Remember that $\ell_1(\ve)=O(\ve^d)=O(\ve^2)$.
We have finally proved that in the case $d-4+4 \alpha <0$ or ($d=3$ and $\alpha=\frac{1}{4}$), $\bbE[S_1^2]=O(\ve^2)$.
The proof of Proposition \ref{term S1}-\ref{itm:termS1.2} is therefore complete.\\
In case where $d-4+4 \alpha =0$ and $d=2$ (and then  $\alpha = {1}/{2}$ ), 
$\ell_2(\ve)=O(\ve^{2} \ln(1/\ve))$, so that $\bbE[S_1^2]=O(\ve^2 \ln(1/\ve))$.
This is the required bound.
However, this way does not give the explicit limit, only an upper bound.
We will therefore have to examine this case separately in order to prove the Proposition \ref{term S1}-\ref{itm:termS1.3}.
\hfill$\bullet$
\end{rema}
Let us return to the case $d-4 + 4 \alpha >0$ and to the proof of Proposition \ref{term S1}-\ref{itm:termS1.1}.
We have proven that 
\begin{align*}
     \lim_{\ve \to 0} \frac{1}{\ve^{4(1-\alpha)}} \ell_2(\ve) 
      &= 
      \int_{T \times T} f(t) f(s) 
\sum_{q=2}^{\infty} \sum_{
\substack{
\bk , \bmm  \in \naturels^{d+1} \\ \abs{\bk }=\abs{\bmm }=q
}
} 
b(\bk ) b(\bmm )  
\bk !  \bmm ! \\
&\rule{24pt}{0pt}\times\sum_{
\substack{
a_{i,j} \ge 0 \\ \sum_{j}a_{i,j}=k_i \\ \sum_{i}a_{i,j}=m_j
}
}  \lim_{\ve \to 0} \ve^{-4(1-\alpha)} A_\ve^{(q)}(\bk , \bmm , s, t, \ba )  \ud t \ud s.
      \end{align*}
Before considering the four different cases that will give the limit, let us make some remarks about some covariance limits.\\
Using the expression for the second derivative of $r_\ve$ given in (\ref{derivees seconde epsilon}), we easily get that for $s, t \in T$ with $s \neq t$ and for $i, j =1,d$,
\begin{equation}
       \label{limite pour derivees r}
     \lim_{\ve \to 0} 
        \dfrac{\partial ^{2}r_\ve}{\partial t_{i}
   \partial t_{j}}(t-s)=  \dfrac{\partial ^{2}r}{\partial t_{i}
   \partial t_{j}}(t-s),
\end{equation}
and in the same manner,
$$
     \lim_{\ve \to 0} 
        \dfrac{\partial r_\ve}{\partial t_{i}}(t-s)=  \dfrac{\partial r}{\partial t_{i}}(t-s).
$$%
Also for $s, t \in T$, we have
$$
     \lim_{\ve \to 0} 
       r_\ve(t-s)=  r(t-s).
$$%
In this way and using the convergence given in (\ref{limite mu epsilon}) we obtain
 \begin{align*}
\lim_{\ve \to 0}  \dfrac{\partial ^{2}r_\ve}{\partial t_{i}
   \partial t_j}(t-s)\frac{\ve^{-2(1-\alpha)} }{\mu_\ve}
   	&=\!\left({ \bC _L \chi_{d}^2(\alpha)}\right)^{-1}\dfrac{\partial ^{2}r}{\partial t_{i}
   \partial t_j}(t-s)\\ \shortintertext{and}
   \lim_{\ve \to 0}  \dfrac{\partial r_\ve}{\partial t_{i}}(t-s)\frac{\ve^{-(1-\alpha)} }{\sqrt{\mu_\ve}\sqrt{r_\ve(0)}}
   &=\frac{1}{\sqrt{r(0)}} \!\left({ \bC _L \chi_{d}^2(\alpha)}\right)^{-{1}/{2}}\dfrac{\partial r}{\partial t_{i}}(t-s).
   \end{align*}
 Let us fix $q \ge 2$, $\bk , \bmm  \in \naturels^{d+1}$ such that $\abs{\bk }=\abs{\bmm }=q$, $\ba :=(a_{i,j})_{1\le i , j \le d+1}$ and $(s, t) \in T \times T$.

From the above, we have to consider four cases.
\begin{itemize}
      \item Case 1.
      It corresponds to situations where:
\begin{align*}
	&\sum_{i, j=1}^d a_{i, j} \ge 3;\\
	\text{~or~}&\sum_{i, j=1}^d a_{i, j} =2 \text{~and~} \sum_{i=1}^d a_{i, d+1} +\sum_{j=1}^d a_{d+1, j} \ge 1;\\
	\text{~or~}&\sum_{i, j=1}^d a_{i, j} =1 \text{~and~} \sum_{i=1}^d a_{i, d+1} + \sum_{j=1}^d a_{d+1, j} \ge 3;\\
	\text{~or~}&\sum_{i, j=1}^d a_{i, j} =0\text{~and~}\sum_{i=1}^d a_{i, d+1} + \sum_{j=1}^d a_{d+1, j} \ge 5.
\end{align*}
Using the previous convergence, we have
      \begin{eqnarray*}
      \lim_{\ve \to 0} \ve^{-4(1-\alpha)} A_\ve^{(q)}(\bk , \bmm , s, t, \ba )=0.
      \end{eqnarray*}
       \item Case 2. For this case, 
\begin{align*}
	\sum_{i, j=1}^d a_{i, j} = 2 \text{,~} \sum_{i=1}^d a_{i, d+1}=0 \text{~and~} \sum_{j=1}^d a_{d+1, j} =0.
\end{align*}
       This configuration is equivalent to the existence of indices $i_0, j_0 \in \{1,\dots, d\}$, such that $a_{i_0, j_0}=2$, $a_{d+1, d+1}=q-2$ and $a_{i,j}=0$ otherwise.\\
With this configuration, we have
\begin{align*}
	\MoveEqLeft[2]{\lim_{\ve \to 0} \ve^{-4(1-\alpha)} A_\ve^{(q)}(\bk , \bmm , s, t, \ba )}\\
	&=\frac{1}{2}\!\left({ \bC _L \chi_{d}^2(\alpha)}\right)^{-2}\!\left({-\dfrac{\partial ^{2}r}{\partial t_{i_0}
   \partial t_{j_0}}(t-s)}\right)^2 \frac{1}{(q-2)!}\!\left({\frac{r(t-s)}{r(0)}}\right)^{q-2}.
\end{align*}
In this case, we have
\begin{align*}
 \bk &=(0, \dots, 0, \uppointer{i_0}{2}, 0, \dots, 0, q-2)\\ \text{and~}
\bmm &=(0, \dots, 0, \uppointer{j_0}{2}, 0, \dots, 0, q-2)
\end{align*}
where ``$\uppointer{i_0}{2}$'' means that 2 is in position $i_0$.\\
        \item Case 3. For this case,
\begin{align*}
	\sum_{i, j=1}^d a_{i, j} = 1\text{;~}\sum_{i=1}^d a_{i, d+1}=1 \text{~and~}\sum_{j=1}^d a_{d+1, j} =1.
\end{align*}
       This configuration is equivalent to the existence of indices $i_0, j_0 \in \{1,\dots,d\}$, such that $a_{i_0, j_0}=1$, $a_{i_0, d+1}=1$, $a_{d+1, j_0}=1$, $a_{d+1, d+1}=q-3$ and $a_{i,j}=0$ otherwise.\\
With this configuration, we have
\begin{align*}
      \MoveEqLeft[2]{\lim_{\ve \to 0} \ve^{-4(1-\alpha)} A_\ve^{(q)}(\bk , \bmm , s, t, \ba )}\\
       &=\!\left({ \bC _L \chi_{d}^2(\alpha)}\right)^{-2}\!\left({-\dfrac{\partial ^{2}r}{\partial t_{i_0}
   \partial t_{j_0}}(t-s)}\right) \!\left({-\frac{1}{r(0)} \dfrac{\partial r}{\partial t_{i_0}}(t-s) \dfrac{\partial r}{\partial t_{j_0}}(t-s)}\right) \\
&\specialpos{\hfill\times    \frac{1}{(q-3)!}\!\left({\frac{r(t-s)}{r(0)}}\right)^{q-3}.}
     \end{align*}
So
\begin{align*}
	 \bk &=(0, \dots, 0, \uppointer{i_0}{2}, 0, \dots, 0, q-2)\\ \text{and~}
	\bmm &=(0, \dots, 0, \uppointer{j_0}{2}, 0, \dots, 0, q-2).
\end{align*}
      \item Case 4. Here,
\begin{align*}
	\sum_{i, j=1}^d a_{i, j} = 0\text{;~}\sum_{i=1}^d a_{i, d+1}=2 \text{~and~}\sum_{j=1}^d a_{d+1, j} =2.
\end{align*}
       This configuration is equivalent to the existence of indices $i_0, j_0 \in \{1,\dots,d \}$ such that $a_{i_0, d+1}=2$, $a_{d+1, j_0}=2$, $a_{d+1, d+1}=q-4$ and $a_{i,j}=0$ otherwise.\\
       With this configuration we have
\begin{align*}
	\MoveEqLeft[2]{\lim_{\ve \to 0} \ve^{-4(1-\alpha)} A_\ve^{(q)}(\bk , \bmm , s, t, \ba )}\\
	= {} &\frac{1}{4}\!\left({ \bC _L \chi_{d}^2(\alpha)}\right)^{-2} \!\left({-\frac{1}{r(0)} \dfrac{\partial r}{\partial t_{i_0}}(t-s) \dfrac{\partial r}{\partial t_{j_0}}(t-s)}\right)^2\times\frac{1}{(q-4)!}\!\left({\frac{r(t-s)}{r(0)}}\right)^{q-4}.
\end{align*}
So we have
\begin{align*}
 \bk &=(0, \dots, 0, \uppointer{i_0}{2}, 0, \dots, 0, q-2)\\ \text{and~}
\bmm &=(0, \dots, 0, \uppointer{j_0}{2}, 0, \dots, 0, q-2).
\end{align*}
  \end{itemize}
 If we combine the four cases, we have proven that
\begin{align*}
 &    \lim_{\ve \to 0} \frac{1}{\ve^{4(1-\alpha)}} \ell_2(\ve) =
      \int_{T \times T} f(t) f(s) \!\left({ \bC _L \chi_{d}^2(\alpha)}\right)^{-2}\\
   & \times   \sum_{i, j=1}^{d}  a(0, \dots, 0, \uppointer{i}{2}, 0, \dots, 0) 
      a(0, \dots, 0, \uppointer{j}{2}, 0, \dots, 0) \\
&  \times\left\{{\sum_{q=2}^{\infty}  2 (q-2)!  c_{q-2}^2 \!\left({-\dfrac{\partial ^{2}r}{\partial t_{i}
   \partial t_{j}}(t-s)}\right)^2 \!\left({\frac{r(t-s)}{r(0)}}\right)^{q-2}}\right.\\
 & \quad+\left.{\sum_{q=3}^{\infty}  4 (q-2)(q-2)! c_{q-2}^2\!\left({-\dfrac{\partial ^{2}r}{\partial t_{i}
   \partial t_{j}}(t-s)}\right) \!\left({-\frac{1}{r(0)} \dfrac{\partial r}{\partial t_{i}}(t-s) \dfrac{\partial r}{\partial t_{j}}(t-s)}\right) 
   }\right.\\
   & \qquad\times\left.{ \!\left({\frac{r(t-s)}{r(0)}}\right)^{q-3} +\sum_{q=4}^{\infty} (q-2)(q-3)(q-2)! c_{q-2}^2}\right.\\
   &\qquad\times\left.{ \!\left({-\frac{1}{r(0)} \dfrac{\partial r}{\partial t_{i}}(t-s) \dfrac{\partial r}{\partial t_{j}}(t-s)}\right)^2 
    \!\left({\frac{r(t-s)}{r(0)}}\right)^{q-4}}\right\} \ud t \ud s.
\end{align*}
Now since for all $i=1, \dots, d$, one has $a(0, \dots, 0, \uppointer{i}{2}, 0, \dots, 0)=a({2,0,\dots,0})$, we finally proved that
$$
     \lim_{\ve \to 0} \frac{1}{\ve^{4(1-\alpha)}} \ell_2(\ve) =\ell_2,
$$%
     and since $\ell_1(\ve)=O(\ve^d)$ and $d-4 + 4 \alpha >0$, one can conclude that 
$$
     \lim_{\ve \to 0} \frac{1}{\ve^{4(1-\alpha)}} \bbE[S_1^2(\ve)] =\ell_2.
$$%
     This gives Proposition \ref{term S1}-\ref{itm:termS1.1} .\\
     To finish the proof it remains to show Proposition \ref{term S1}--\ref{itm:termS1.3}.
So let us suppose that $d=2$ and that  $\alpha = {1}/{2}$.
Our goal is now to prove the following convergence:
$$
      \frac{1}{\ve^{2} \ln(1/\ve)} \bbE[S_1^2] \cvg[\ve\to0]{} \ell_5.
$$%
First we look at the term $\ell_1(\ve)$ which we have seen to be $O(\ve^2)$.
So $\ell_1(\ve)=o(\ve^2 \ln(1/\ve))$.\\
Consider the second term $\ell_2(\ve)$, which, as we saw earlier, is equal to
$$
 \ell_2(\ve)= \int_{T_\ve^{(2)}} 
\sum_{q=2}^{\infty} \sum_
	{
	\substack{
		\bk , \bmm  \in \naturels^{3} \\ \abs{\bk }=\abs{\bmm }=q
			}
	} 
	f(t) f(s) b(\bk ) b(\bmm )  
	\bk !  \bmm ! A_\ve^{(q)}(\bk , \bmm , s, t)  \ud t \ud s,
$$%
where for $q \ge 2$ and $\bk , \bmm  \in \naturels^{3}$ such that $\abs{\bk }=\abs{\bmm }=q$ and $(s, t) \in T_\ve^{(2)}$,
$$
A_\ve^{(q)}(\bk , \bmm , s, t):=
\sum_{
\substack{
a_{i,j} \ge 0 \\ \sum_{j}a_{i,j}=k_i \\ \sum_{i}a_{i,j}=m_j
}
} A_\ve^{(q)}(\bk , \bmm , s, t, \ba ),
$$%
where for $\ba  :=(a_{i,j})_{1\le i , j \le 3}$,
\begin{align*}
\MoveEqLeft[2]{A_\ve^{(q)}(\bk , \bmm , s, t, \ba ):=
\prod_{1 \le i, j \le 2} \left[{
\!\left({\frac{-\dfrac{\partial ^{2}r_\ve}{\partial t_{i}
   \partial t_{j}}(t-s)}{\mu_\ve}
}\right)^{a_{i,j}}\frac{1}{a_{i,j}!}
}\right]}\\
 & \times \prod_{1 \le i \le 2} \left[{
\!\left({\frac{\dfrac{\partial r_\ve}{\partial t_{i}}(t-s)}{\sqrt{\mu_\ve}{\sqrt{r_\ve(0)}}}
}\right)^{a_{i,3}}\frac{1}{a_{i,3}!}
}\right] \times
\prod_{1 \le j \le 2} \left[{
\!\left({\frac{-\dfrac{\partial r_\ve}{\partial t_{j}}(t-s)}{\sqrt{\mu_\ve}{\sqrt{r_\ve(0)}}}
}\right)^{a_{3,j}}\frac{1}{a_{3,j}!}
}\right] \\
&\times \!\left({
\frac{r_\ve(t-s)}{r_\ve(0)}
}\right)^{a_{3,3}}\frac{1}{a_{3,3}!}.
\end{align*}
 Let us fix $q \ge 2$, $\bk , \bmm  \in \naturels^{3}$ such that $\abs{\bk }=\abs{\bmm } =q$, $\ba  :=(a_{i,j})_{1\le i , j \le 3}$ and $(s, t) \in T_\ve^{(2)}$.
\\
We want to find an upper bound for $\ve^{-2}  [\ln(1/\ve)]^{-1} A_\ve^{(q)}(\bk , \bmm , s, t, \ba ) $.\\
Applying the same type of reasoning as before, we distinguish four cases.\\
\begin{itemize}
\item[Case 1.]
$\sum_{i, j=1}^2 a_{i, j} \ge 3$.\\[3pt]
In this case, we obtain
$$
\ve^{-2}  [\ln(1/\ve)]^{-1} \abs{A_\ve^{(q)}(\bk , \bmm , s, t, \ba )} \le \bC  \normp[2]{t-s}^{-3}  \ve  [\ln(1/\ve)]^{-1} \prod_{1 \le i, j \le 3}  \frac{1}{a_{i,j}!},
$$
and 
\begin{equation}
\label{epsilon moins un}
\int_{T_\ve^{(2)}} \normp[2]{t-s}^{-3} \ud t \ud s  \le \bC  \ve^{-1}.
\end{equation}
\item[Case 2.]
$\sum_{i, j=1}^2 a_{i, j} = 1$.\\[3pt]
In this case, we get
$$
\ve^{-2}  [\ln(1/\ve)]^{-1} \abs{A_\ve^{(q)}(\bk , \bmm , s, t, \ba )} \le \bC  \normp[2]{t-s}^{-1}   [\ln(1/\ve)]^{-1} \prod_{1 \le i, j \le 3}  \frac{1}{a_{i,j}!},
$$
and 
\begin{equation}
\label{epsilon puissance zero}
\int_{T_\ve^{(2)}} \normp[2]{t-s}^{-1} \ud t \ud s  \le \bC .
\end{equation}
\item[Case 3.]
$\sum_{i, j=1}^2 a_{i, j} = 0$.\\[3pt]
In this case, we have
$$
\ve^{-2}  [\ln(1/\ve)]^{-1} \abs{A_\ve^{(q)}(\bk , \bmm , s, t, \ba )} \le \bC   [\ln(1/\ve)]^{-1} \prod_{1 \le i, j \le 3}  \frac{1}{a_{i,j}!},
$$
and 
\begin{equation}
\label{epsilon puissance zero bis}
\int_{T_\ve^{(2)}} 1 \ud t \ud s  \le \bC .
\end{equation}
\item[Case 4.]
$\sum_{i, j=1}^2 a_{i, j} =2$ and $\sum_{i=1}^2 a_{i, 3} + \sum_{j=1}^2 a_{3, j} \ge 1$.\\[3pt]
In this case, we obtain
$$
\ve^{-2}  [\ln(1/\ve)]^{-1} \abs{A_\ve^{(q)}(\bk , \bmm , s, t, \ba )}
 \le \bC  \normp[2]{t-s}^{-2}  \ve^{\frac{1}{2}}  [\ln(1/\ve)]^{-1} \prod_{1 \le i, j \le 3}  \frac{1}{a_{i,j}!},
$$
and 
\begin{equation}
\label{epsilon racine}
\int_{T_\ve^{(2)}} \normp[2]{t-s}^{-2} \ud t \ud s  \le \bC  \ln(1/\ve).
\end{equation}
\end{itemize}
Let us denote by $B(\bk , \bmm )$ the following set:
\begin{multline*}
B(\bk , \bmm ):= \left\{\ba =(a_{i,j})_{i, j=1,2,3} : \sum_{j=1}^3 a_{i,j}=k_i, \sum_{i=1}^3 a_{i,j}=m_j,  \sum_{i, j=1}^2 a_{i, j} \ge 3 \mbox{ or  }\right.\\
 \left.\sum_{i, j=1}^2 a_{i, j} =1 \mbox{ or  }\sum_{i, j=1}^2 a_{i, j} = 0 \mbox{ or  } \sum_{i, j=1}^2 a_{i, j} =2 \mbox{ with  } \sum_{i=1}^2 a_{i, 3} + \sum_{j=1}^2 a_{3, j} \ge 1 \right\}.
\end{multline*}
Summarizing all the cases and using (\ref{epsilon moins un}), (\ref{epsilon puissance zero}), (\ref{epsilon puissance zero bis}), (\ref{epsilon racine}) and (\ref{graphes}), it is clear that
 \begin{align*}
 \MoveEqLeft[2]{\ve^{-2}  [\ln(1/\ve)]^{-1} \int_{T_\ve^{(2)}} 
\sum_{q=2}^{\infty} \sum_{
\substack{
\bk , \bmm  \in \naturels^{3} \\ \abs{\bk }=\abs{\bmm }=q
}
} 
\abs{f(t)} \abs{f(s)} \abs{b(\bk )} \abs{b(\bmm )} 
\bk !  \bmm !} \\
	&\specialpos{\hfill\times\sum_{\ba   \in  B(\bk , \bmm )} \abs{A_\ve^{(q)}(\bk , \bmm , s, t, \ba )}  \ud t \ud s}\\
	&\le 
 \bC   [\ln(1/\ve)]^{-1} \sum_{q=2}^{\infty} \sum_{
\substack{
\bk , \bmm  \in \naturels^{3} \\ \abs{\bk }=\abs{\bmm }=q
}
} 
\abs{b(\bk )} \abs{b(\bmm )}  
\bk !  \bmm ! 
 \sum_{
\substack{
a_{i,j} \ge 0 \\ \sum_{j}a_{i,j}=k_i \\ \sum_{i}a_{i,j}=m_j
}
} \prod_{1 \le i, j \le 3}  \frac{1}{a_{i,j}!}\\
	& \le 
 \bC   [\ln(1/\ve)]^{-1} \cvg[\ve\to 0]{} 0.
\end{align*}

To conclude, we have shown that
\begin{multline}
\label{Dynkin Bis}
\ve^{-2} [\ln(1/\ve)]^{-1} \ell_2(\ve) \equiv
\frac{2 a^2(2,0)}{\ve^2  \ln(\frac{1}\ve)} \int_{T_\ve^{(2)}} f(t) f(s)  \sum_{q=0}^\infty q! c_q^2 \!\left({
\frac{r_\ve(t-s)}{r_\ve(0)}
}\right)^q   \\
 \times \left[{
\!\left({\dfrac{\partial ^{2}r_\ve}{\partial t_1^2}(t-s)/\mu_\ve
}\right)^{2}  
 + 2 \!\left({\dfrac{\partial ^{2}r_\ve}{\partial t_1\partial t_2}(t-s)/\mu_\ve
}\right)^{2} +\!\left({\dfrac{\partial ^{2}r_\ve}{\partial t_2^2}(t-s)/\mu_\ve}\right)^{2}
}\right] \ud t \ud s.
      \end{multline}
\spacebefore
\begin{rema}
\label{H2}
This last part of the proof highlights the fact that
\begin{multline*}
\ve^{-1} [\ln(1/\ve)]^{-1/2} S_1 \simeq 
\ve^{-1} [\ln(1/\ve)]^{-1/2} a(2, 0) \int_T f(t) g\!\left({X_\ve(t)}/{\sqrt{r_\ve(0)}}\right)\\
\times\left[{H_2\!\left(\dfrac{\partial X_\ve}{\partial t_1}(t)/\sqrt{\mu_\ve}\right)+
H_2\!\left(\dfrac{\partial X_\ve}{\partial t_2}(t)/\sqrt{\mu_\ve}\right)
}\right] \ud t.
\end{multline*}
\hfill$\bullet$
\end{rema}
To obtain the asymptotic variance, we split the set $T_\ve^{(2)}$ into two sets, $U_{\delta}$ and $V_{\ve, \delta}$, defined by
\begin{align*}
	U_{\delta}&:= \left\{(s, t) \in T \times T:\normp[2]{s-t} > \delta\right\}\\ \shortintertext{and}
	V_{\ve, \delta}&:=  \left\{\smash{(s, t) \in T_\ve^{(2)}: \normp[2]{s-t} \le \delta}\right\},
\end{align*}
where $\delta >0$ is fixed and will eventually tend to zero.
The corresponding integrals will be denoted respectively by $L_{1, \delta}(\ve)$ and $L_{2, \delta}(\ve)$.\\
On the one hand, using the convergence given in (\ref{limite mu epsilon}), we know that
$$
\ve \mu_\ve \cvg[\ve\to0]{} \bC _L \chi_2^2(1/2).
$$%
On the other hand, using Lebesgue's dominated convergence theorem, it is straightforward to prove that
\begin{align*}
\MoveEqLeft[2]{\ln(1/\ve)  L_{1, \delta}(\ve) \cvg[\ve\to0]{}
\frac{2 a^2(2,0)}{(\bC _L \chi_2^2(1/2))^2} \int_{U_{\delta}} f(t) f(s) 
\sum_{q=0}^\infty q! c_q^2 \!\left({
\frac{r(t-s)}{r(0)}
}\right)^q} \\
&
\times\left[{
\!\left({\dfrac{\partial ^{2}r}{\partial t_1^2}(t-s)
}\right)^{2}  + 2 \!\left({\dfrac{\partial ^{2}r}{\partial t_1\partial t_2}(t-s)
}\right)^{2} +\!\left({\dfrac{\partial ^{2}r}{\partial t_2^2}(t-s)
}\right)^{2} 
}\right] \ud t \ud s.
\end{align*}
Thus, using once again the convergence given in (\ref{limite mu epsilon}), we get
\begin{align*}
\MoveEqLeft[3]{\ve^{-2} [\ln(1/\ve)]^{-1} \ell_2(\ve) \equiv L_{2, \delta}(\ve)}\\
\equiv {} &\ \frac{2 a^2(2,0)}{(\bC _L \chi_2^2(1/2))^2}   [\ln(1/\ve)]^{-1}   \\
&\times\int_{V_{\ve, \delta}} f(t) f(s)  \sum_{q=0}^\infty q! c_q^2 \!\left({
\frac{r_\ve(t-s)}{r_\ve(0)}
}\right)^q \\
&\times 
\left[{
\!\left({\dfrac{\partial ^{2}r_\ve}{\partial t_1^2}(t-s)
}\right)^{2}  
 + 2 \!\left({\dfrac{\partial ^{2}r_\ve}{\partial t_1\partial t_2}(t-s)
}\right)^{2} +\!\left({\dfrac{\partial ^{2}r_\ve}{\partial t_2^2}(t-s)}\right)^{2}
}\right] \ud t \ud s.
      \end{align*}
Inspired by Lemma \ref{majoration pour r(epsilon)} and using the additional assumptions made on the function $L$, we pinpoint that for $(s, t) \in T \times T$, such that $M\ve \le \normp[2]{t-s} $ and for all $i, j= 1, 2$,
$$
  \abs{\dfrac{\partial ^{2}r_\ve}{\partial t_{i}\partial t_{j}}(t-s)-\dfrac{\partial ^{2}r}{\partial t_{i}\partial t_j}(t-s)} \le \bC   \frac\ve{\normp[2]{t-s}^2},
$$%
  which implies that
$$
  \abs{\!\left({\dfrac{\partial ^{2}r_\ve}{\partial t_{i}\partial t_{j}}(t-s)}\right)^2-\!\left({\dfrac{\partial ^{2}r}{\partial t_{i}\partial t_j}(t-s)}\right)^2} \le \bC   \frac\ve{\normp[2]{t-s}^3}.
$$%
 Moreover, using Lemmas \ref{majoration pour r(epsilon)} and \ref{norm}, we can prove that, for all $q \in \naturels$,
$$
   \abs{\!\left({\frac{r_\ve(t-s)}{r_\ve(0)} }\right)^q-1} \le \bC   q \normp[2]{t-s}.
$$%
Finally, we also have the following inequality
$$
\!\left({\dfrac{\partial ^{2}r}{\partial t_{i}\partial t_j}(t-s)}\right)^2 \le \bC   \normp[2]{t-s}^{-2}.
$$%
Using (\ref{epsilon moins un}), (\ref{epsilon puissance zero}), the fact that $\sum_{q=1}^{\infty}q q! c_q^2 < \infty$ and Lebesgue's dominated convergence theorem, we easily obtain that
\begin{multline*}
L_{2, \delta}(\ve) \equiv
	\frac{2 a^2(2,0)}{(\bC _L \chi_2^2(1/2))^2}  \!\left({\sum_{q=0}^\infty q! c_q^2}\right) [\ln(1/\ve)]^{-1}
		\int_{V_{\ve, \delta}} f(t) f(s) \\
	\times\left[{
		\!\left({\dfrac{\partial ^{2}r}{\partial t_1^2}(t-s)}\right)^{2}  
		+ 2 \!\left({\dfrac{\partial ^{2}r}{\partial t_1\partial t_2}(t-s)}\right)^{2}
		+\!\left({\dfrac{\partial ^{2}r}{\partial t_2^2}(t-s)}\right)^{2}
	}\right] \ud t \ud s.
\end{multline*}
To conclude this section, we have proven that
$$
 \ve^{-2} [\ln(1/\ve)]^{-1} \bbE[S_1^2]  \equiv
\frac{2 a^2(2,0)}{(\bC _L \chi_2^2(1/2))^2} \bbE\!\left[g^2\!\left({X(0)}/{\sqrt{r(0)}}\right)\right] U_{\delta}(\ve),
$$%
where
 \begin{multline*}
U_{\delta}(\ve):=
 [\ln(1/\ve)]^{-1}  \int_{V_{\ve, \delta}} f(t) f(s) 
\left[{
\!\left({\dfrac{\partial ^{2}r}{\partial t_1^2}(t-s)
}\right)^{2}  
 }\right.\\
\left.{ + 2 \!\left({\dfrac{\partial ^{2}r}{\partial t_1\partial t_2}(t-s)
}\right)^{2} +\!\left({\dfrac{\partial ^{2}r}{\partial t_2^2}(t-s)}\right)^{2}
}\right] \ud t \ud s.
      \end{multline*}
To carry out this proof, let us focus on the term $U_{\delta}(\ve)$.\\
The hypotheses made on function $L$, allow us to write the following asymptotic equality, for $(s, t) \in T \times T$, and for all $i, j= 1, 2$:
$$
 \!\left({\dfrac{\partial ^{2}r}{\partial t_{i}\partial t_j}(t-s)}\right)^2 = \bC _L^2 \frac{(t_i-s_i)^2(t_j-s_j)^2}{\normp[2]{t-s}^6} + o\!\left(\frac{1}{\normp[2]{t-s}^2}\right).
$$%
  Furthermore,
  $$
  \sup_{\{(s, t) \in T \times T, \normp[2]{s-t} \le \delta \}} \abs{f(s)-f(t)}\cvg[\delta \to0]{}0.
  $$
  The last convergence and the last asymptotic equality plus the fact that
  $$
  \underset{\ve \to 0}{\overline{\lim}}  [\ln(1/\ve)]^{-1} \int_{V_{\ve, \delta}} \normp[2]{t-s}^{-2} \ud t \ud s\le (2\pi) \lambda_2(T)< \infty
  $$
  imply that
$$
	U_{\delta}(\ve)= U_{\delta}^{(1)}(\ve)+U_{\delta}^{(2)}(\ve),
$$
where 
$$
U_{\delta}^{(1)}(\ve)=  \bC _L^2  [\ln(1/\ve)]^{-1} \int_{V_{\ve, \delta}} f^2(t) \normp[2]{t-s}^{-2} \ud t \ud s,
$$%
and 
$$
\lim_{\delta \to 0} \underset{\ve \to 0}{\overline{\lim}} \abs{U_{\delta}^{(2)}(\ve)}=0.
$$%
Now we consider $T_{-\delta}:= \{t \in \reels^2: \mbox{dist}(t, T^c) > \delta\}$.
With this notation we have
\begin{multline*}
 U_{\delta}^{(1)}(\ve)=
   \bC _L^2 \!\left({\int_{T_{-\delta}} f^2(t) \ud t}\right)  [\ln(1/\ve)]^{-1}\\
   \times \!\left({
 \int_{\{w \in \reels^2, M\ve \le \normp[2]{w} \le \delta\}} \normp[2]{w}^{-2} \ud w}\right) +W_{\delta, \ve},
\end{multline*}
where
$$
\lim_{\delta \to 0} \underset{\ve \to 0}{\overline{\lim}} W_{\delta, \ve}=0.
$$%
Since 
\begin{align*}
 	&\lim_{\ve \to 0}  [\ln(1/\ve)]^{-1}  \int_{\{w \in \reels^2, M\ve \le \normp[2]{w} \le \delta\}} \normp[2]{w}^{-2} \ud w= 2\pi\\
 \shortintertext{and}
 	&\lim_{\delta \to 0} \int_{T_{-\delta}} f^2(t) \ud t=\int_{T} f^2(t) \ud t,
\end{align*}
 we finally proved that
\begin{multline*}
 \lim_{\ve \to 0} \ve^{-2} [\ln(1/\ve)]^{-1} \bbE[S_1^2]  =
\frac{4\pi a^2(2,0)}{( \chi_2^2(1/2))^2} \bbE\!\left[g^2\!\left({X(0)}/{\sqrt{r(0)}}\right)\right] \!\left({ \int_{T} f^2(t) \ud t}\right)\\
= a^2(2,0) \sigma^2  \!\left({ \int_{T} f^2(t) \ud t}\right) \bbE\!\left[g^2\!\left({X(0)}/{\sqrt{r(0)}}\right)\right]= \ell_5.
\end{multline*}
The proposition ensues.
\end{proofarg}

Let us return to the proof of Theorem \ref{vitesse variance temps local}.
By Remark \ref{remark theoreme}, we already know that parts \ref{itm4.6.19.1}, \ref{itm4.6.19.2} and \ref{itm4.6.19.4} of Theorem \ref{vitesse variance temps local} follow from the last proposition.
\\
Let us show part \ref{itm4.6.19.3} of the theorem.
Suppose that $\alpha = {1}/{2}$ and $d>2$.
We pointed out in Remark \ref{remark theoreme}, that in this case, $\frac{1}{\ve^{2}} \bbE[S_1^2] \cvg[\ve\to0]{}\ell_2$.
We also proved that in case where $\alpha = \frac{1}{2}$, $\frac{1}{\ve^{2}} \bbE[S_2^2] \cvg[\ve\to0]{} \ell_1+ \ell_3$.\\
So to finish with part \ref{itm4.6.19.3} of the theorem it only remains to prove that in the case where  $\alpha = {1}/{2}$  and $d> 2$, we have
$$
	\frac{1}{\ve^{2}} \bbE\!\left[S_1S_2\right] \cvg[\ve\to0]{} \frac{1}{2} \ell_4.
$$%
We will prove this result in the more general case where  $\alpha = {1}/{2}$  (and $d \ge 2$).\\
First, we compute $\frac{1}{\ve^{2}} \bbE\!\left[S_1S_2\right]$.
\begin{align*}
	\frac{1}{\ve^{2}} \bbE[S_1S_2] = \frac{1}{\ve^2} \int_{T \times T} f(t) f(s)
\times \bbE\left[g\!\left({X_\ve(t)}/{\sqrt{r_\ve(0)}}\right) 
	 h\!\left({\nabla\!X_\ve(t)}/{\sqrt{\mu_\ve}}\right)\right. \\
 	\left.\left\{{g\!\left({X_\ve(s)}/{\sqrt{r_\ve(0)}}\right)
 - g\!\left({X(s)}/{\sqrt{r(0)}}\right)
 }\right\}\right]  \ud t \ud s.
\end{align*}
As we did before, we split the integration domain $T \times T$ into two parts $T_\ve^{(1)}$ and $T_\ve^{(2)}$.
Let us call $c_1(\ve)$ (resp$.$
$c_2(\ve)$) the integral over $T_\ve^{(1)}$ (resp$.$
$T_\ve^{(2)}$).\\
Let us focus first on the first term $c_1(\ve)$.
\begin{align*}
\abs{c_1(\ve)} &\le \bC  \!\left({\bbE\!\left[g^2(N_1)\right]}\right)^{\frac{1}{2}} 
 \!\left({\bbE\!\left[h^2(N_d)\right]}\right)^{\frac{1}{2}} \\
	& \quad\times\!\left({\bbE\!\left[\left\{g\!\left({X_\ve(0)}/{\sqrt{r_\ve(0)}}\right)
 - g\!\left({X(0)}/{\sqrt{r(0)}}\right)\right\}^2\right]}\right)^{\frac{1}{2}} \frac{1}{\ve^2} \sigma_{2d}(T_\ve^{(1)})\\
	& \le  \bC  \ve^{d-2}  \!\left({\bbE\!\left[\left\{g\!\left({X_\ve(0)}/{\sqrt{r_\ve(0)}}\right)
 - g\!\left({X(0)}/{\sqrt{r(0)}}\right)\right\}^2\right]}\right)^{\frac{1}{2}}.
\end{align*}
By calculating the last expectation with Mehler's formula (see Lemma \ref{mehler}), we obtain
\begin{multline*}
 \bbE\!\left[\left\{g\!\left({X_\ve(0)}/{\sqrt{r_\ve(0)}}\right)
 - g\!\left({X(0)}/{\sqrt{r(0)}}\right)\right\}^2\right]\\=2 \sum_{n=1}^{\infty} c_n^2 n! \!\left({1-\!\left({\frac{r^{(\ve)}(0)}{\sqrt{r(0)}\sqrt{r_\ve(0)}}}\right)^n }\right) \cvg[\ve\to0]{} 0,
 \end{multline*}
since $\abs{{r^{(\ve)}(0)}/{\!\left(\sqrt{r(0)}\sqrt{r_\ve(0)}\right)}} \le 1$ and $\sum_{n=1}^{\infty} c_n^2 n! < \infty$.\\
So we have proved that the term $c_1(\ve)\cvg[\ve\to0]{}0$ whether the dimension $d$ is greater than or equal to two or strictly greater than two.\\
 Let us now tackle the second term $c_2(\ve)$.
We obtain the following decomposition:
 \begin{align*}
\MoveEqLeft[2]{c_2(\ve)= \frac{1}{\ve^2} \int_{T_\ve^{(2)}} f(t) f(s) 
\bbE\!\left[{\rule{0pt}{18pt}g\!\left({X_\ve(t)}/{\sqrt{r_\ve(0)}}\right) 
 h\!\left({\nabla\!X_\ve(t)}/{\sqrt{\mu_\ve}}\right)
 }\right.}\\
 &\specialpos{\hfill\times \left.{ \sum_{n=1}^{\infty} c_n \!\left({1-\!\left({\sqrt{{r_\ve(0)}/{r(0)}}}\right)^n }\right) H_n\!\left({X_\ve(s)}/{\sqrt{r_\ve(0)}}\right)
 }\right]  \ud t \ud s}\\
&\qquad+ \frac{1}{\ve^2} \int_{T_\ve^{(2)}} f(t) f(s) 
\bbE\!\left[\rule{0pt}{17pt}{g\!\left({X_\ve(t)}/{\sqrt{r_\ve(0)}}\right) 
 h\!\left({\nabla\!X_\ve(t)}/{\sqrt{\mu_\ve}}\right)
 }\right.\\
 \times&\left.{  \sum_{n=1}^{\infty} c_n
 	\left\{{
	 \!\left({
		\sqrt{{r_\ve(0)}/{r(0)}}
	 }\right)^n
	 H_n\!\left({X_\ve(s)}/{\sqrt{r_\ve(0)}}\right)- H_n
 \!\left({X(s)}/{\sqrt{r(0)}}\right)
 }\right\}
 }\right]  \ud t \ud s\\
 &:= d_1(\ve)+d_2(\ve).
 \end{align*}
Let’s first look at the first term $d_1(\ve)$ which will give the limit.\\
With the same notations as in the proof of Proposition \ref{term S1} and by applying Mehler's formula (see Lemma \ref{mehler}), we obtain
 \begin{align*}
\MoveEqLeft[2]{d_1(\ve)= \frac{1}{\ve^2} \int_{T_\ve^{(2)}} f(t) f(s) \sum_{q=2}^\infty 
\sum_{
	\substack{\bk  \in \naturels^{d+1} \\ \abs{\bk }=q}
	} 
b(\bk )  c_q  \!\left({1-\!\left({\sqrt{\frac{r_\ve(0)}{r(0)}}}\right)^q }\right)} \\
& \quad\times\bbE\!\left[\widetilde{H}_{\bk }(U_\ve(t))
H_q\!\left({X_\ve(s)}/{\sqrt{r_\ve(0)}}\right)\right] \ud t \ud s
\\
& = \int_{T_\ve^{(2)}} f(t) f(s) \sum_{q=2}^\infty \sum_{
\substack{
\bk  \in \naturels^{d+1} \\ \abs{\bk }=q
}
} 
 a(k_1, k_2, \dots, k_d) c_{k_{d+1}} c_q\\
& \quad\times \frac{1}\ve \!\left({1-\!\left({\sqrt{\frac{r_\ve(0)}{r(0)}}}\right)^q }\right)q! \!\left({
\frac{r_\ve(t-s)}{r_\ve(0)}
}\right)^{k_{d+1}} \\
&\quad\times \frac{1}\ve 
 \!\left({
\frac{\dfrac{\partial r_\ve}{\partial t_1}(t-s)}{\sqrt{\mu_\ve}{\sqrt{r_\ve(0)}}}
}\right)^{k_1} \!\left({
\frac{\dfrac{\partial r_\ve}{\partial t_2}(t-s)}{\sqrt{\mu_\ve}{\sqrt{r_\ve(0)}}}
}\right)^{k_2}  \ldots \!\left({
\frac{\dfrac{\partial r_\ve}{\partial t_{d}}(t-s)}{\sqrt{\mu_\ve}{\sqrt{r_\ve(0)}}}
}\right)^{k_d} \ud t \ud s.
 \end{align*}
 Let us justify that we can interchange the limit in $\ve$ with the integral and the sum.

 On the one hand, for $q \ge 2$,  we first bound the term $\frac{1}\ve \!\left({1-\!\left({\sqrt{\frac{r_\ve(0)}{r(0)}}}\right)^q }\right)$ by $q \frac{1}\ve \!\left({1- \sqrt{\frac{r_\ve(0)}{r(0)}}}\right)$.
Then taking into account that
 $\frac{1}\ve \!\left({1- \sqrt{\frac{r_\ve(0)}{r(0)}}}\right) \cvg[\ve\to0]{} \ell_{{1}/{2}}$,
 we finally see that  $\frac{1}\ve \!\left({1-\!\left({\sqrt{\frac{r_\ve(0)}{r(0)}}}\right)^q }\right)$ has an upper bound $\bC   q$, for $\ve$ small enough.\\
 On the other hand, arguing as in the proof of Proposition \ref{term S1}, we know that for $q \ge 2$ and $\abs{\bk }=q$, there exists at least one index $i_0 \in \{1,\dots,d\}$ such that $k_{i_0} \ge 2$.
Thus we have an upper bound for
$$
\abs{\frac{1}\ve 
 \!\left({
\frac{\dfrac{\partial r_\ve}{\partial t_1}(t-s)}{\sqrt{\mu_\ve}{\sqrt{r_\ve(0)}}}
}\right)^{k_1} \!\left({
\frac{\dfrac{\partial r_\ve}{\partial t_2}(t-s)}{\sqrt{\mu_\ve}{\sqrt{r_\ve(0)}}}
}\right)^{k_2}  \ldots \!\left({
\frac{\dfrac{\partial r_\ve}{\partial t_{d}}(t-s)}{\sqrt{\mu_\ve}{\sqrt{r_\ve(0)}}}
}\right)^{k_d} \!\left({
\frac{r_\ve(t-s)}{r_\ve(0)}
}\right)^{k_{d+1}}}.
$$%
 
To conclude, it is sufficient to justify that the following series is finite
$$
\sum\limits_{q=2}^\infty \sum\limits_{
	\substack{\bk  \in \naturels^{d+1} \\ \abs{\bk }=q}
	} 
\abs{b(\bk )} \abs{c_q} q q!< \infty.
$$%
Here again, we refer to the proof of Proposition \ref{term S1} and Lemma \ref{serie finie} which we will adapt to the new situation.\\
The following function $G$ is $L^2(\reels^{d+1}, \varphi_{d+1}(x) \ud x)$, with the following expansion in this space.
For $x:=(x_1, \dots, x_{d+1}) \in \reels^{d+1}$,
$$
G(x_1, \dots, x_{d+1}):= h(x_1, x_2, \dots,x_d)  g(x_{d+1})= \sum_{q=2}^{\infty} \sum_{
\substack{
\bk  \in \naturels^{d+1} \\ \abs{\bk }=q
}
} b(\bk ) \widetilde{H}_{\bk }(x).
$$%
Also the function $H_g \in L^2(\reels, \varphi_1(x) \ud x)$.
Recall that, for $x \in \reels$,
$$
Hg(x)= xg^{\prime}(x)-g^{\prime\prime}(x) = \sum_{n=1}^{\infty} c_n n H_n(x).
$$%
Thus
$$
\sum\limits_{q=2}^\infty \sum\limits_{
\substack{
\bk  \in \naturels^{d+1} \\ \abs{\bk }=q
}
} 
 \abs{b(\bk )} \abs{c_q} q q!< \infty.
$$%
So we proved that
\begin{align*}
\lim_{\ve \to 0} d_1(\ve) &=
\frac{\ell_{{1}/{2}}}{r(0)} \!\left({\frac{a({2,0,\dots,0})}{\bC _L \chi_{d}^2(1/2)} }\right) \sum_{q=2}^\infty c_{q-2}  c_q q  q! \\
&\quad\times \int_{T \times T} f(t) f(s) \!\left(\frac{r(t-s)}{r(0)}\right)^{q-2} \!\left({\sum_{i=1}^d \!\left({  \dfrac{\partial   r}{\partial t_{i}}(t-s) }\right )^2}\right) \ud t \ud s\\
  &=\frac{\ell_4}{2}.
\end{align*}
 It only remains to show that $d_2(\ve)\cvg[\ve\to0]{}0$.
 \begin{multline*}
d_2(\ve)= \int_{T_\ve^{(2)}} f(t) f(s) \sum_{q=2}^\infty \sum_{
\substack{
\bk  \in \naturels^{d+1} \\ \abs{\bk }=q
}
} 
b(\bk ) c_q  q!\\
\times \frac{1}{\ve^2} 
\left[{
 \!\left({
\frac{\dfrac{\partial r_\ve}{\partial t_1}(t-s)}{\sqrt{\mu_\ve}{\sqrt{r(0)}}}
}\right)^{k_1} \!\left({
\frac{\dfrac{\partial r_\ve}{\partial t_2}(t-s)}{\sqrt{\mu_\ve}{\sqrt{r(0)}}}
}\right)^{k_2}  \ldots \!\left({
\frac{\dfrac{\partial r_\ve}{\partial t_{d}}(t-s)}{\sqrt{\mu_\ve}{\sqrt{r(0)}}}
}\right)^{k_d} }\right.\\
\times\left.{ \!\left({
\frac{r_\ve(t-s)}{\sqrt{r_\ve(0)}\sqrt{r(0)}}
}\right)^{k_{d+1}}}\right.\\
\left.{-
 \!\left({
\frac{\dfrac{\partial r^{(\ve)}}{\partial t_1}(t-s)}{\sqrt{\mu_\ve}{\sqrt{r(0)}}}
}\right)^{k_1} \!\left({
\frac{\dfrac{\partial r^{(\ve)}}{\partial t_2}(t-s)}{\sqrt{\mu_\ve}{\sqrt{r(0)}}}
}\right)^{k_2}  \ldots \!\left({
\frac{\dfrac{\partial r^{(\ve)}}{\partial t_{d}}(t-s)}{\sqrt{\mu_\ve}{\sqrt{r(0)}}}
}\right)^{k_d} }\right.\\
\times\left.{ \!\left({
\frac{r^{(\ve)}(t-s)}{\sqrt{r_\ve(0)}\sqrt{r(0)}}
}\right)^{k_{d+1}}
}\right]
\ud t \ud s.
 \end{multline*}
Still with the same type of arguments as those given previously, and considering the case where the even indices $k_i$, for $i=1, \dots, d$, are such that $\sum_{i=1}^d k_i \ge 4$, we apply Lebesgue's dominated convergence theorem and the limit gives zero.
So we consider only the remaining terms for which there is an index $i_0=1,\dots, d$ such that $k_{i_0}=2$, $k_i=0$ for $i \neq i_0$, $i=1,\dots,d$ and $k_{d+1}=q-2$.
We get
  \begin{multline*}
\lim_{\ve \to 0}d_2(\ve)=a({2,0,\dots,0}) \lim_{\ve \to 0}  \int_{T_\ve^{(2)}} f(t) f(s) \sum_{q=2}^\infty 
c_{q-2} c_q  q! \sum_{i=1}^d \frac{1}{\ve^2}\\
 \times\left[{
 \!\left({
\frac{\dfrac{\partial r_\ve}{\partial t_{i}}(t-s)}{\sqrt{\mu_\ve}{\sqrt{r(0)}}}
}\right)^{2} \!\left({
\frac{r_\ve(t-s)}{\sqrt{r_\ve(0)}\sqrt{r(0)}}
}\right)^{q-2}}\right.\\
\left.{-
 \!\left({
\frac{\dfrac{\partial r^{(\ve)}}{\partial t_{i}}(t-s)}{\sqrt{\mu_\ve}{\sqrt{r(0)}}}
}\right)^{2} \!\left({
\frac{r^{(\ve)}(t-s)}{\sqrt{r_\ve(0)}\sqrt{r(0)}}
}\right)^{q-2}
}\right]
\ud t \ud s.
 \end{multline*}
 Let us justify that this limit is zero.\\
 For $i =1,\dots,d$, $s, t \in T_\ve^{(2)}$ and $q \ge 2$, we bound the following term
 \begin{multline*}
 A_\ve(q, t-s):=\frac{1}{\ve^2} 
\left|
 \!\left({
\frac{\dfrac{\partial r_\ve}{\partial t_{i}}(t-s)}{\sqrt{\mu_\ve}{\sqrt{r(0)}}}
}\right)^{2} \!\left({
\frac{r_\ve(t-s)}{\sqrt{r_\ve(0)}\sqrt{r(0)}}
}\right)^{q-2}\right. \\
\specialpos{\hfill\left.- 
 \!\left({
\frac{\dfrac{\partial r^{(\ve)}}{\partial t_{i}}(t-s)}{\sqrt{\mu_\ve}{\sqrt{r(0)}}}
}\right)^{2} \!\left({
\frac{r^{(\ve)}(t-s)}{\sqrt{r_\ve(0)}\sqrt{r(0)}}
}\right)^{q-2}
\right|} \\
\le \bC  \frac{1}\ve \left\{{ (q-2) \abs{r_\ve(t-s)-r^{(\ve)}(t-s)} + 
\abs{
\dfrac{\partial r_\ve}{\partial t_{i}}(t-s)-\dfrac{\partial r^{(\ve)}}{\partial t_{i}}(t-s)
}
}\right\}.
\end{multline*}
On the one hand, by applying Lemma \ref{majoration pour r(epsilon)}, we can bound $\abs{r_\ve(t-s)-r^{(\ve)}(t-s)}$ by $\bC  \ve^2 \normp[d]{t-s}^{-1}$.
On the other hand, using arguments similar to the one given for prove this lemma, one can show that if we let $\tau:=t-s$,
 \begin{multline*}
 \frac{1}\ve\!\left({\dfrac{\partial r_\ve}{\partial t_{i}}(\tau)-\dfrac{\partial r^{(\ve)}}{\partial t_{i}}(\tau)}\right)\\
	=  \int_{\{\normp[d]{v} \le N\}} (\Phi(v)-\Psi(v)) \times 
	 \left\{{
	   \sum_{j=1}^{d} v_j \dfrac{\partial ^{2}r}{\partial \tau_{i}
	   \partial \tau_{j}}(\tau -\theta \ve v)
	    }\right \} \ud v,
\end{multline*}
with $0\le \theta <1$ depending on $\ve$, $\tau$ and $v$.\\
This equality involves two things.
The first is that 
$$
\abs{\dfrac{\partial r_\ve}{\partial t_{i}}(t-s)-\dfrac{\partial r^{(\ve)}}{\partial t_{i}}(t-s)}
$$%
can be bounded by $\bC  \ve \normp[d]{t-s}^{-1}$.\\
 The second is that $ \frac{1}\ve\!\left({\dfrac{\partial r_\ve}{\partial t_{i}}(\tau)-\dfrac{\partial r^{(\ve)}}{\partial t_{i}}(\tau)}\right) \cvg[\ve\to0]{} 0$, since the functions $\Phi$ and $\Psi$ are even functions.
\\
 Finally, we collect the following facts
$$
 A_\ve(q, t-s)\le \bC  q \normp[d]{t-s}^{-1},\quad \text{and} \quad A_\ve(q, t-s) \cvg[\ve\to0]{}0.
$$%
To conclude, using that $\sum\limits_{q=2}^\infty \sum\limits_{
\substack{
\bk  \in \naturels^{d+1} \\ \abs{\bk }=q
}
} 
 \abs{b(\bk )} \abs{c_q} q q!< \infty$, we deduce that $\sum\limits_{q=2}^\infty 
 \abs{c_{q-2}} \abs{c_q} q q!< \infty$.\\
Moreover, since $d \ge 2$, we know that $\int_{T \times T} \normp[d]{t-s}^{-1}\ \ud t \ud s < \infty$.\\
 All the ingredients are gathered to apply  once again Lebesgue's dominated convergence theorem, leading to the required limit:
 $$
 \lim\limits_{\ve \to 0}d_2(\ve)=0.
 $$%
That complete the proof of Theorem \ref{vitesse variance temps local}.
\end{proofarg}

Let us now prove Theorem \ref{vitesse convergence temps local}.
\spacebefore
\begin{proofarg}{Proof of Theorem \ref{vitesse convergence temps local}}
We begin by showing part \ref{itm1:theo4.6.20} of the theorem.
So let us suppose that $\alpha < \frac{1}{2}$.
In this case we have seen in the proof of Theorem \ref{vitesse variance temps local} that the \rv
$  \Xi_\ve(f, g)$ is equivalent in $L^2(\Omega)$ to $\ve^{-2\alpha}S_2$, itself equivalent in $L^2(\Omega)$ to $\ve^{-2\alpha}U_2$, where
$$
       U_2=  \sum_{n=1}^{\infty} c_n  \int_T f(t)    \left\{{
       1- \!\left({\sqrt{\frac{r_\ve(0)}{r(0)}}
       }\right)^n
  }\right\} H_n\!\left({X_\ve(t)}/{\sqrt{r_\ve(0)}}\right) \ud t.
$$%
Let us see that the \rv
$\ve^{-2\alpha}U_2$ is in turn equivalent in $L^2(\Omega)$ to the following \rv
$X(f,g)$.
$$
 X(f, g):= \ell_{\alpha} \sum_{n=1}^{\infty} n c_n \int_{T} f(t) H_n\!\left({{X(t)}/{\sqrt{r(0)}}}\right) \ud t.
$$%
The proof is very similar to the one given to prove that $\ve^{-4\alpha} \bbE\!\left[U_2^2\right]$ tends to $\ell_1$.
It consists in using Mehler's formula and the fact that
$
  \ve^{-2\alpha} \!\left({
 1-\sqrt{
{r_\ve(0)}/{r(0)}
}}\right) \cvg[\ve\to 0]{} \ell_{\alpha}
 $
 and that $\sum_{n=1}^{\infty} n! n^2  c_n^2 < \infty$.\\
 This being so, let us recall that the function $Hg$ was defined in the introduction by, 
$$
 Hg(x)= x g^{\prime}(x)-g^{\prime\prime}(x)
 =  \sum_{n=1}^{\infty} n c_n H_n(x).
$$%
Thus the limit \rv $X(f, g)$ can be expressed as
$$
 X(f, g)= \ell_{\alpha} \int_{T} f(t) Hg\!\left({{X(t)}/{\sqrt{r(0)}}}\right) \ud t.
$$%
By applying Theorem \ref{GH} to the function $h(t, x):= f(t) Hg(x)$ and to $X:=\ds {x}/{\sqrt{r(0)}}$, this \rv
can also be expressed as
$$
 X(f, g)= \ell_{\alpha} \int_{\reels} Hg(x)  L_{X/\sqrt{r(0)}}^f(x, T) \ud x,
$$%
 which is the required result.
This ends proof of part \ref{itm1:theo4.6.20}.

We shall now tackle the proof of part \ref{itm2:theo4.6.20}.
Let us suppose that $\alpha >\frac{1}{2}$.\\
We have seen in Theorem \ref{vitesse variance temps local} that the \rv
$  \Xi_\ve(f, g)$ is equivalent in $L^2(\Omega)$ to $T_1=\ve^{-2(1-\alpha)}S_1$, thus
$$
 \Xi_\ve(f, g) \simeq T_1,
$$%
where
$$T_1= \ve^{-2(1-\alpha)} \int_T f(t)  g\!\left({X_\ve(t)}/{\sqrt{r_\ve(0)}}\right)  h\!\left({\nabla\!X_\ve(t)}/{\sqrt{\mu_\ve}}\right) 
   \ud t, 
$$%
and recall that for $x:=(x_1, x_2, \dots, x_d) \in \reels^{d}$,
$$
h(x_1, x_2, \dots,x_d) = \sum_{q=2}^{\infty} \sum_{
\substack{
\bk  \in \naturels^{d} \\ \abs{\bk }=q
}
} a(\bk ) \widetilde{H}_{\bk }(x),
$$%
while coefficients $a(\bk )$ are defined by (\ref{coef_ak}).
Thus, to study the convergence of $  \Xi_\ve(f, g)$ it is enough to consider that of $T_1$.\\
Let us define the following \rv:%
$$
W_\ve(g):=\ve^{-2(1-\alpha)} \int_T f(t)  g\!\left({X_\ve(t)}/{\sqrt{r_\ve(0)}}\right)  h_{\bm{2}}\!\left(\frac{\nabla\!X_\ve(t)}{\sqrt{\mu_\ve}}\right) 
   \ud t, 
$$%
where
   \begin{align*}
h_{\bm{2}}(x_1, x_2, \dots,x_d) &:= \sum_{
\substack{
\bk  \in \naturels^{d} \\ \abs{\bk }
=2
}
} a(\bk ) \widetilde{H}_{\bk }(x)\\
	& \phantom{:}= \sum_{j=1}^d a(0, \dots, 0, \uppointer{j}{2}, 0, \dots, 0) H_2(x_j)\\
	& \phantom{:}= a({2,0,\dots,0}) \sum_{j=1}^d H_2(x_j).
\end{align*}
It is almost obvious that $\bbE\!\left[\left(T_1-W_\ve(g)\right)^2\right] \cvg[\ve\to0]{} 0$.
Just look at the way we \linebreak showed that $\bbE[T_1^2] \cvg[\ve\to0]{} \ell_2$, to see that $\bbE[W_\ve^2(g)]\cvg[\ve\to0]{} \ell_2$ and that $\bbE[T_1W_\ve(g)] \cvg[\ve\to0]{} \ell_2$.\\
Thus to study the convergence of $T_1$ it is enough to consider that of $W_\ve(g)$.\\
We have
\begin{align*}
W_\ve(g)&\phantom{:}=a({2,0,\dots,0}) \sum_{j=1}^d \ve^{-2(1-\alpha)}\\
	& \qquad \times\int_T f(t)  g\!\left({X_\ve(t)}/{\sqrt{r_\ve(0)}}\right)  H_2\!\left(\dfrac{\partial X_\ve}{\partial t_{j}}(t)/\sqrt{\mu_\ve}\right) \ud t\\
	&:= a({2,0,\dots,0}) \sum_{j=1}^d \widetilde{W}_\ve^{(j)}(g).
   \end{align*}
We fix $j \in \{1, 2, \dots, d\}$.
To simplify the notation, we will study the asymptotic behavior of $\widetilde{W}_\ve^{(j)}(g)$.
Given the form of the assumed limit, this will be sufficient since the convergence will take place in $L^2(\Omega)$.
Using the Hermite coefficients of the function $g$ we obtain
\begin{multline*}
\widetilde{W}_\ve^{(j)}(g)= \ve^{-2(1-\alpha)} \int_T f(t)  \\
\times\sum_{k=2}^{\infty} c_{k-2} H_{k-2}\!\left({X_\ve(t)}/{\sqrt{r_\ve(0)}}\right)  H_2\!\left(\dfrac{\partial X_\ve}{\partial t_{j}}(t)/
\sqrt{\mu_\ve}\right) 
   \ud t. 
\end{multline*}
Let us justify the interchange of the sum and the integral.
This is a consequence of the following facts.
On the one hand, let us define the Cauchy sequence $(D_N)_N$ in $L^2(\Omega)$.
$$
D_N:= \sum_{k=2}^{N} c_{k-2} \int_T f(t) H_{k-2}\!\left({X_\ve(t)}/{\sqrt{r_\ve(0)}}\right)  H_2\!\left(\dfrac{\partial X_\ve}{\partial t_{j}}(t)/\sqrt{\mu_\ve}\right) 
   \ud t.
$$%
On the other hand, defining in $L^2(\varphi(x) \ud x)$,
$$
g_N:=g-\sum_{k=2}^{N} c_{k-2} H_{k-2},
$$%
one has 
\begin{align*}
\MoveEqLeft[2]{\bbE\!\left[\left\{\int_T f(t)  g_N\!\left({X_\ve(t)}/{\sqrt{r_\ve(0)}}\right)  H_2\!\left(\dfrac{\partial X_\ve}{\partial t_{j}}(t)/\sqrt{\mu_\ve}\right) 
   \ud t\right\}^2\right]} \\
 & \le \lambda_{d}(T) \int_T f^2(t) \bbE\!\left[g^2_N\!\left({X_\ve(t)}/{\sqrt{r_\ve(0)}}\right)  H_2^2\!\left(\dfrac{\partial X_\ve}{\partial t_{j}}(t)/\sqrt{\mu_\ve}\right)\right] \ud t \\
& = \lambda_{d}(T) \!\left({\int_T f^2(t) \ud t}\right) \bbE\!\left[g^2_N\!\left({X_\ve(0)}/{\sqrt{r_\ve(0)}}\right)\right]  \bbE\!\left[H_2^2\!\left(\dfrac{\partial X_\ve}{\partial t_{j}}(0)/\sqrt{\mu_\ve}\right)\right]\\
& =2! \lambda_{d}(T) \!\left({\int_T f^2(t) \ud t}\right) \normp[2, \varphi]{g_N}^2 \cvg[N\to\infty]{} 0,
\end{align*}
while $\normp[2, \varphi]{\cdot}$ represents the $L^2(\reels, \varphi(x) \ud x)$ norm.\\
Finally, we obtain
\begin{multline*}
 \widetilde{W}_\ve^{(j)}(g)= \ve^{-2(1-\alpha)} \sum_{k=2}^{\infty} c_{k-2}\\
 \times \int_T f(t)   H_{k-2}\!\left({X_\ve(t)}/{\sqrt{r_\ve(0)}}\right)  H_2\!\left(\dfrac{\partial X_\ve}{\partial t_{j}}(t)/\sqrt{\mu_\ve}\right) 
   \ud t.
\end{multline*}
Denoting by  $\ud Z_X(\lambda)$ the random spectral measure corresponding to $X$, $\widehat{\Psi}$ the Fourier transform of $\Psi$ and $\ud Z_\ve(\lambda):= \widehat{\Psi}(\ve \lambda) \ud Z_X(\lambda)$, we have the following representation for $\widetilde{W}_\ve^{(j)}(g)$.
\begin{multline*}
\widetilde{W}_\ve^{(j)}(g)= \ve^{-2(1-\alpha)} \sum_{k=2}^{\infty} c_{k-2} \int_T f(t)   H_{k-2}\!\left({\frac{1}{\sqrt{r_\ve(0)}} \int_{\reels^d} e^{i \prodsca{\lambda}{t} } \ud Z_\ve(\lambda)}\right)  \\
H_2\!\left({\frac{1}{\sqrt{\mu_\ve}} \int_{\reels^d} i \lambda^{(j)} e^{i \prodsca{\lambda}{t} } \ud Z_\ve(\lambda)
}\right) 
   \ud t,
\end{multline*}
where we noted
$\lambda:=(\lambda^{(1)}, \lambda^{(2)}, \dots, \lambda^{(j)}, \dots, \lambda^{(d)})$.\\

For $(\lambda, t) \in \reels^d \times T$, let  $\omega_1^{(j)}$ and $\omega_2^{(j)}$ be
$$
\omega_1^{(j)}(\lambda, t):= \frac{1}{\sqrt{r_\ve(0)}}  e^{i \prodsca{\lambda}{t}}, \quad \omega_2^{(j)}(\lambda, t):= \frac{i \lambda^{(j)}}{\sqrt{\mu_\ve}}  e^{i \prodsca{\lambda}{t}}.
$$%
The functions $\omega_1^{(j)}(\cdot, t)$ and $\omega_2^{(j)}(\cdot, t)$ are orthogonal with respect to the measure $\abs{\widehat{\Psi}(\ve \lambda)}^2 s(\lambda) \ud \lambda$, where we recall that $s$ is the spectral density of the process $X$.
Thus we can use the It\^o's formula for the Wiener-It\^o integral (see \cite{MR3155040}, Theorem 4.3, p$.$
37), obtaining
$$
\widetilde{W}_\ve^{(j)}(g)= \sum_{k=2}^{\infty} T_k^{(j)}(\ve),
$$%
where 
\begin{multline*}
T_k^{(j)}(\ve):= \ve^{-2(1-\alpha)}  c_{k-2} \int_T f(t)  \frac{1}{k!}  \int_{\reels^{dk}} \sum_{\Pi \in \Pi_{k}} 
\omega_{\xi(\Pi(1))}^{(j)}(\lambda_1, t) 
\omega_{\xi(\Pi(2))}^{(j)}(\lambda_2, t)  \\
\dots\omega_{\xi(\Pi(k))}^{(j)}(\lambda_k, t) 
\ud Z_\ve(\lambda_1) \ud Z_\ve(\lambda_2) \dots \ud Z_\ve(\lambda_k) \ud t,
\end{multline*}
with
$$
\xi(n) =
\begin{cases}
 1 & n \le k-2\\
 2 & \text{otherwise}
\end{cases}
$$
and where $\Pi_{k}$ is the set of permutations of $\{1, 2, \dots, k\}$.
Thus
\begin{multline*}
T_k^{(j)}(\ve)= -\frac{c_{k-2}}{(\sqrt{r_\ve(0)})^{k-2}} \frac{\ve^{-2(1-\alpha)}}{\mu_\ve}  \frac{1}{k!}  \int_{\reels^{dk}} K_f(\lambda_1+\lambda_2+\dots +\lambda_k)
 \\
\times\sum_{\Pi\in \Pi_{k}} \lambda_{\Pi(k-1)}^{(j)}  \lambda_{\Pi(k)}^{(j)} \ud Z_\ve(\lambda_1) \ud Z_\ve(\lambda_2) \dots \ud Z_\ve(\lambda_k) \mbox{ \as},
\end{multline*}
where $K_f$ is given by (\ref{Kf}).\\
Since by  (\ref{limite mu epsilon}) we know that
$$
   \ds \frac{\mu_\ve}{\ve^{-2(1-\alpha)}} \cvg[\ve \to 0]{} \bC _L \chi_{d}^2(\alpha),
$$%
to obtain the asymptotic behavior of $T_k^{(j)}(\ve)$ it is enough to consider the convergence of the term below.
Define
\begin{multline*}
M_k^{(j)}(\ve):= -\frac{c_{k-2}}{(\sqrt{r_\ve(0)})^{k-2}} \frac{1}{k!}  \int_{\reels^{dk}} K_f(\lambda_1+\lambda_2+\dots +\lambda_k)
\\
\times \sum_{\Pi \in \Pi_{k}} \lambda_{\Pi(k-1)}^{(j)}  \lambda_{\Pi(k)}^{(j)} \ud Z_\ve(\lambda_1) \ud Z_\ve(\lambda_2) \dots \ud Z_\ve(\lambda_k).
\end{multline*}
The second moment of this \rv is 
\begin{align*}
\MoveEqLeft[2]{\bbE\!\left[\left(M_k^{(j)}(\ve)\right)^2\right]=
\frac{c_{k-2}^2}{(r_\ve(0))^{k-2}} \frac{1}{k!}  \int_{\reels^{dk}} \abs{K_f(\lambda_1+\lambda_2+\dots +\lambda_k)}^2}
\\
	&\quad\times \sum_{\Pi \in \Pi_{k}} \sum_{\Upsilon \in \Pi_{k}} \lambda_{\Pi(k-1)}^{(j)}  \lambda_{\Pi(k)}^{(j)}  \lambda_{\Upsilon(k-1)}^{(j)}  \lambda_{\Upsilon(k)}^{(j)} 
s(\lambda_1) s(\lambda_2) \dots  s(\lambda_k)\\
	&\quad\times \abs{\widehat{\Psi}(\ve \lambda_1)}^2
 \abs{\widehat{\Psi}(\ve \lambda_2)}^2 \dots \abs{\widehat{\Psi}(\ve \lambda_k)}^2 \ud \lambda_1 \ud \lambda_2 \dots \ud \lambda_k.
\end{align*}
At this stage of the proof, we must state a lemma whose proof is given just after the proof of the part \ref{itm2:theo4.6.20}
of the theorem.
\spacebefore
\begin{lemm}
\label{mu epsilon}
Let $(\mu_\ve)$ be a sequence of finite measures in $(\reels^d)^k$ such that their density $H_\ve^{(j)}$ with respect to the Lebesgue measure in $(\reels^d)^k$ can be written as
\begin{align*}
\MoveEqLeft[2]{H_\ve^{(j)}(\lambda_1, \lambda_2, \dots, \lambda_k):= \abs{K_f(\lambda_1+\lambda_2+\dots +\lambda_k)}^2} \\
	&\quad\times \sum_{\Pi \in \Pi_{k}} \sum_{\Upsilon \in \Pi_{k}} \lambda_{\Pi(k-1)}^{(j)}  \lambda_{\Pi(k)}^{(j)}  \lambda_{\Upsilon(k-1)}^{(j)}  \lambda_{\Upsilon(k)}^{(j)} 
s(\lambda_1) s(\lambda_2) \dots  s(\lambda_k)\\
	&\quad\times \abs{\widehat{\Psi}(\ve \lambda_1)}^2
 \abs{\widehat{\Psi}(\ve \lambda_2)}^2 \dots \abs{\widehat{\Psi}(\ve \lambda_k)}^2.
\end{align*}
Then $\mu_\ve \cvg[\ve\to0]{\text{weakly}} \mu$, where $\mu$ is a finite measure in $(\reels^d)^k$ with density $H^{(j)}$ with respect to the Lebesgue measure in $(\reels^d)^k$ given by
\begin{multline*}
H^{(j)}(\lambda_1, \lambda_2, \dots, \lambda_k):= \abs{K_f(\lambda_1+\lambda_2+\dots +\lambda_k)}^2\\
\sum_{\Pi \in \Pi_{k}} \sum_{\Upsilon \in \Pi_{k}} \lambda_{\Pi(k-1)}^{(j)}  \lambda_{\Pi(k)}^{(j)}\lambda_{\Upsilon(k-1)}^{(j)}  \lambda_{\Upsilon(k)}^{(j)} 
s(\lambda_1) s(\lambda_2) \dots  s(\lambda_k).
\end{multline*}
It implies that
\begin{multline*}
\lim_{\ve \to 0} \int_{(\reels^d)^k} H_\ve^{(j)}(\lambda_1, \lambda_2, \dots, \lambda_k) \ud \lambda_1 \ud \lambda_2 \dots \ud \lambda_k\\
= \int_{(\reels^d)^k} H^{(j)}(\lambda_1, \lambda_2, \dots, \lambda_k) \ud \lambda_1 \ud \lambda_2 \dots \ud \lambda_k.
\end{multline*}
\end{lemm}
Suppose that this lemma is proved.
Let us define
\begin{multline*}
M_k^{(j)}(0):= -\frac{c_{k-2}}{\left(\sqrt{r(0)}\right)^{k-2}} \frac{1}{k!}  \int_{\reels^{dk}} K_f(\lambda_1+\lambda_2+\dots + \lambda_k) \\
\times\sum_{\Pi \in \Pi_{k}} \lambda_{\Pi(k-1)}^{(j)}  \lambda_{\Pi(k)}^{(j)} \ud Z_{X}(\lambda_1) \ud Z_{X}(\lambda_2) \dots \ud Z_{X}(\lambda_k).
\end{multline*}
This It\^o-Wiener integral is well defined since by Lemma \ref{mu epsilon} we have
$$
\bbE\!\left[\left(M_k^{(j)}(0)\right)^2\right]=\frac{c_{k-2}^2}{(r(0))^{k-2}} \frac{1}{k!}  \int_{\reels^{dk}} H^{(j)}(\lambda_1, \lambda_2, \dots, \lambda_k) \ud \lambda_1 \ud \lambda_2 \dots \ud \lambda_k < \infty.
$$%
Furthermore 
\begin{equation}
\lim_{\ve \to 0} \bbE\!\left[\left(M_k^{(j)}(\ve)\right)^2\right]=\bbE\!\left[\left(M_k^{(j)}(0)\right)^2\right].
\label{convergence Mk}
\end{equation}
Consider now
\begin{multline*}
D_k^{(j)}(\ve, \lambda_1, \lambda_2, \dots, \lambda_k):= 
-\frac{c_{k-2}}{(\sqrt{r_\ve(0)})^{k-2}} \frac{1}{\sqrt{k!}}   K_f(\lambda_1+\lambda_2+\dots + \lambda_k)\\
\times \sum_{\Pi \in \Pi_{k}} \lambda_{\Pi(k-1)}^{(j)}  \lambda_{\Pi(k)}^{(j)} 
\sqrt{s(\lambda_1)} \sqrt{s(\lambda_2)} \dots  \sqrt{s(\lambda_k)}  \widehat{\Psi}(\ve \lambda_1)
 \widehat{\Psi}(\ve \lambda_2) \dots \widehat{\Psi}(\ve \lambda_k),
\end{multline*}
and
\begin{multline*}
D_k^{(j)}(0, \lambda_1, \lambda_2, \dots, \lambda_k):= 
-\frac{c_{k-2}}{(\sqrt{r(0)})^{k-2}} \frac{1}{\sqrt{k!}}   K_f(\lambda_1+\lambda_2+\dots + \lambda_k) \\
\times \sum_{\Pi \in \Pi_{k}} \lambda_{\Pi(k-1)}^{(j)}  \lambda_{\Pi(k)}^{(j)} 
\sqrt{s(\lambda_1)} \sqrt{s(\lambda_2)} \dots  \sqrt{s(\lambda_k)}.
\end{multline*}
On the one hand the convergence appearing in (\ref{convergence Mk}) means that  $\normp[2]{D_k^{(j)}(\ve, \cdot)}^2$ converges to $\normp[2]{D_k^{(j)}(0, \cdot)}^2$ in $L^2$ with respect to the Lebesgue measure in $(\reels^d)^k$.
On the other hand, the \rv $D_k^{(j)}(\ve, \cdot)$ converges pointwise to $D_k^{(j)}(0, \cdot)$ as $\ve \to 0$.
Lebesgue theorem implies that $\normp[2]{D_k^{(j)}(\ve, \cdot)-D_k^{(j)}(0, \cdot)} \cvg[\ve\to0]{$L^2$} 0$ in the $L^2$-norm with respect to the Lebesgue measure in $(\reels^d)^k$.
That is $M_k^{(j)}(\ve) \cvg[\ve\to0]{$L^2(\Omega)$} M_k^{(j)}(0)$ and finally as $\ve \to 0$ we proved that 
\begin{equation}
	T_k^{(j)}(\ve) \cvg[\ve\to0]{$L^2(\Omega)$} T_k^{(j)}(0),
	\label{Tk}
\end{equation}
where 
$$
	T_k^{(j)}(0):= \frac{1}{\bC _L \chi_{d}^2(\alpha)}  M_k^{(j)}(0).
$$%
We will now prove that
$$
	\widetilde{W}_\ve^{(j)}(g)= \sum_{k=2}^{\infty} T_k^{(j)}(\ve)
		\cvg[\ve\to0]{$L^2(\Omega)$}
		\sum_{k=2}^{\infty} T_k^{(j)}(0).
$$%
First,
\begin{align*}
\MoveEqLeft[2]{\bbE\!\left[\left(\widetilde{W}_\ve^{(j)}(g)\right)^2\right]
 \cvg[\ve\to\0]{} 2\!\left({\frac{1}{\bC _L \chi_{d}^2(\alpha)} }\right)^2  \sum_{k=2}^{\infty} \sum_{\ell=(k-4)\vee 0}^{k-2} (k-2)! c_{k-2}^2}\\
	&\qquad\times{{k-2}\choose{\ell}}
	{{k-2-\ell}\choose{2}}\\
	&\qquad \times \int_{T \times T} f(t) f(s) 
  \!\left(\frac{r(t-s)}{r(0)}\right)^{\ell} \!\left({\frac{1}{\sqrt{r(0)}} \dfrac{\partial r}{\partial t_{j}}(t-s)}\right )^{2(k-2-\ell)} (-1)^{k-2-\ell}\\
	& \qquad \times \!\left({ -\dfrac{\partial ^2  r}{\partial t_{j}^2}(t-s) }\right )^{\ell - (k-4)}  \ud t \ud s.
\end{align*}
A similar argument gives the result when the sum begins with $N+1$.
Using convergence given in (\ref{Tk}) this implies that
$$
\lim_{\ve \to 0} \bbE\!\left[\left(\sum_{k=N+1}^{\infty} T_k^{(j)}(\ve)\right)^2\right]= \sum_{k=N+1}^{\infty} \bbE\!\left[\left(T_k^{(j)}(0)\right)^2\right].
$$%
And since $ \sum_{k=2}^{\infty} \bbE\!\left[\left(T_k^{(j)}(0)\right)^2\right] < \infty$, then
\begin{equation}
\label{limite N epsilon}
\underset{N \to \infty}{\overline{\lim}} \;\underset{\ve \to 0}{\overline{\lim}} \bbE\!\left[\left(\sum_{k=N+1}^{\infty} T_k^{(j)}(\ve)\right)^2\right]= 0.
\end{equation}
On the one hand
\begin{multline*}
\bbE\!\left[\left(\sum_{k=2}^{\infty} \left(T_k^{(j)}(\ve)-T_k^{(j)}(0)\right)\right)^2\right]\\
 \le \bC  \left\{{
\bbE\!\left[\left(\sum_{k=N+1}^{\infty} T_k^{(j)}(\ve)\right)^2\right] + \bbE\!\left[\left(\sum_{k=2}^{N} \left(T_k^{(j)}(\ve)-T_k^{(j)}(0)\right)\right)^2\right]}\right.\\
 + \left.{\bbE\!\left[\left(\sum_{k=N+1}^{\infty} T_k^{(j)}(0)\right)^2\right]
}\right\},
\end{multline*}
and using once again convergence given in (\ref{Tk}) we obtain
\begin{multline*}
\underset{\ve \to 0}{\overline{\lim}} \bbE\!\left[\left(\sum_{k=2}^{\infty} \left(T_k^{(j)}(\ve)-T_k^{(j)}(0)\right)\right)^2\right] \\
 \le\bC  \left\{{
\underset{\ve \to 0}{\overline{\lim}} \bbE\!\left[\left(\sum_{k=N+1}^{\infty} T_k^{(j)}(\ve)\right)^2\right] +\bbE\!\left[\left(\sum_{k=N+1}^{\infty} T_k^{(j)}(0)\right)^2\right]
}\right\}.
\end{multline*}
On the other hand, since $ \sum_{k=2}^{\infty} \bbE\!\left[\left(T_k^{(j)}(0)\right)^2\right] < \infty$ and $\bbE\!\left[T_{k_1}^{(j)}(0)T_{k_2}^{(j)}(0)\right]=0$ for $k_1 \neq k_2$, we deduce that
\begin{equation}
\label{limite N}
\underset{N \to \infty}{\overline{\lim}} \bbE\!\left[\left(\sum_{k=N+1}^{\infty} T_k^{(j)}(0)\right)^2\right]=0.
\end{equation}
Using the convergences given in (\ref{limite N epsilon}) and (\ref{limite N}) we finally obtain
\begin{multline*}
\underset{N \to \infty}{\overline{\lim}} \underset{\ve \to 0}{\overline{\lim}} \bbE\!\left[\left(\sum_{k=2}^{\infty} \left(T_k^{(j)}(\ve)-T_k^{(j)}(0)\right)\right)^2\right]  \\
	\le \bC  \left\{{\underset{N \to \infty}{\overline{\lim}} 
\underset{\ve \to 0}{\overline{\lim}} \bbE\!\left[\left(\sum_{k=N+1}^{\infty} T_k^{(j)}(\ve)\right)^2\right] +\underset{N \to \infty}{\overline{\lim}} \bbE\!\left[\left(\sum_{k=N+1}^{\infty} T_k^{(j)}(0)\right)^2\right] }\right\}=0.
\end{multline*}
We thus proved that $\lim_{\ve \to 0} \bbE\!\left[\left(\sum_{k=2}^{\infty} \left(T_k^{(j)}(\ve)-T_k^{(j)}(0)\right)\right)^2\right]=0$.
Thus 
$$
\widetilde{W}_\ve^{(j)}(g) \cvg[\ve \to 0]{$L^2(\Omega)$} \sum_{k=2}^{\infty} T_k^{(j)}(0)
$$%
and
\begin{multline*}
 W_\ve(g)=a({2,0,\dots,0}) \sum_{j=1}^d \widetilde{W}_\ve^{(j)}(g)\\
\cvg[\ve \to 0]{$L^2(\Omega)$}
a({2,0,\dots,0}) \sum_{j=1}^d \sum_{k=2}^{\infty} T_k^{(j)}(0)=Y(f, g).
\end{multline*}
That yields item \ref{itm2:theo4.6.20} of the theorem.
\spacebefore
\begin{rema}
\label{cas plus general}
It is obvious using Proposition \ref{term S1}-\ref{itm:termS1.1} that $T_1$ converges in $L^{2}(\Omega)$ to $Y(f, g)$ in the more general case where $d+4\alpha -4 >0$.
Thus, the last convergence prevails in the particular case where  $\alpha = {1}/{2}$  and $d>2$.
\hfill$\bullet$
\end{rema}
\spacebefore
\begin{proofarg}{Proof of Lemma \ref{mu epsilon}}
Let $\chi_\ve(\gamma_1, \gamma_2, \dots, \gamma_k)$ be the Fourier transform of $\mu_\ve$ in $(\reels^d)^k$.
That is
\begin{align*}
\MoveEqLeft[2]{\chi_\ve(\gamma_1, \gamma_2, \dots, \gamma_k):= 
\int_{\reels^{dk}} e^{i \prodsca{\lambda_1}{\gamma_1}} e^{i \prodsca{\lambda_2}{\gamma_2}} \dots  e^{i \prodsca{\lambda_k}{\gamma_k}}}\\
&\quad\times \abs{K_f(\lambda_1+\lambda_2+\dots \lambda_k)}^2
\sum_{\Pi \in \Pi_{k}} \sum_{\Upsilon \in \Pi_{k}} \lambda_{\Pi(k-1)}^{(j)}  \lambda_{\Pi(k)}^{(j)} \lambda_{\Upsilon(k-1)}^{(j)}  \lambda_{\Upsilon(k)}^{(j)}\\
&\quad\times s(\lambda_1) s(\lambda_2) \dots  s(\lambda_k) \abs{\widehat{\Psi}(\ve \lambda_1)}^2
 \abs{\widehat{\Psi}(\ve \lambda_2)}^2 \dots \abs{\widehat{\Psi}(\ve \lambda_k)}^2 \ud \lambda_1 \ud \lambda_2 \dots \ud \lambda_k\\
&=\int_{\reels^{dk}} e^{i \prodsca{\lambda_1}{\gamma_1}} e^{i \prodsca{\lambda_2}{\gamma_2}} \dots  e^{i \prodsca{\lambda_k}{\gamma_k}}
\int_{T \times T} f(t) f(s) e^{i \prodsca{\lambda_1+ \lambda_2 + \dots + \lambda_k}{t-s}} \ud t \ud s\\
&\quad\times \sum_{\Pi \in \Pi_{k}} \sum_{\Upsilon \in \Pi_{k}} \lambda_{\Pi(k-1)}^{(j)}  \lambda_{\Pi(k)}^{(j)} \lambda_{\Upsilon(k-1)}^{(j)}  \lambda_{\Upsilon(k)}^{(j)} 
s(\lambda_1) s(\lambda_2) \dots  s(\lambda_k)\\ 
&\quad\times \abs{\widehat{\Psi}(\ve \lambda_1)}^2
 \abs{\widehat{\Psi}(\ve \lambda_2)}^2 \dots \abs{\widehat{\Psi}(\ve \lambda_k)}^2 \ud \lambda_1 \ud \lambda_2 \dots \ud \lambda_k\\
&=\int_{T \times T} f(t) f(s)  \int_{\reels^{dk}} e^{i \prodsca{\lambda_1}{\gamma_1+t-s}} e^{i \prodsca{\lambda_2}{\gamma_2+t-s}} \dots  e^{i \prodsca{\lambda_k}{\gamma_k+t-s}}\\
&\quad\times \sum_{\Pi \in \Pi_{k}} \sum_{\Upsilon \in \Pi_{k}} \lambda_{\Pi(k-1)}^{(j)}  \lambda_{\Pi(k)}^{(j)} \lambda_{\Upsilon(k-1)}^{(j)}  \lambda_{\Upsilon(k)}^{(j)} 
s(\lambda_1) s(\lambda_2) \dots  s(\lambda_k)\\
&\quad\times \abs{\widehat{\Psi}(\ve \lambda_1)}^2
 \abs{\widehat{\Psi}(\ve \lambda_2)}^2 \dots \abs{\widehat{\Psi}(\ve \lambda_k)}^2 \ud \lambda_1 \ud \lambda_2 \dots \ud \lambda_k \ud t \ud s.
\end{align*}
We have to consider three cases depending on the permutations that appear in the last expression.
  \begin{itemize}
      \item[Case 1]
For two different integers, the permutations are the same.
We can suppose, without loss of generality, that $\Pi(k)=\Upsilon(k)=1$ and $\Pi(k-1)=\Upsilon(k-1)=2$.
The integral is then
      \begin{align*}
\MoveEqLeft[2]{\int_{T \times T} f(t) f(s)  \int_{\reels^{dk}} e^{i \prodsca{\lambda_1}{\gamma_1+t-s}} e^{i \prodsca{\lambda_2}{\gamma_2+t-s}} \dots  e^{i \prodsca{\lambda_k}{\gamma_k+t-s}}} \\
&\quad \times\!\left({\lambda_1^{(j)}}\right)^2 \!\left({\lambda_2^{(j)}}\right)^2  s(\lambda_1) s(\lambda_2) \dots  s(\lambda_k)\\
&\quad\times \abs{\widehat{\Psi}(\ve \lambda_1)}^2 \abs{\widehat{\Psi}(\ve \lambda_2)}^2 \dots \abs{\widehat{\Psi}(\ve \lambda_k)}^2 \ud \lambda_1 \ud \lambda_2 \dots \ud \lambda_k \ud t \ud s\\
& =  \int_{T \times T} f(t) f(s) \!\left({ \dfrac{\partial ^2  r_\ve}{\partial t_{j}^2}(\gamma_1+t-s) }\right ) \!\left({ \dfrac{\partial ^2  r_\ve}{\partial t_{j}^2}(\gamma_2+t-s) }\right ) \\
&\quad \times r_\ve(\gamma_3+t-s) \dots  r_\ve(\gamma_k+t-s) \ud t \ud s.
      \end{align*}
      There are $2k!(k-2)!$ such cases.
       \item[Case 2]
The permutations coincide for one index and differ for the other two indices.
We may assume, again without loss of generality that, $\Pi(k)=\Upsilon(k)=1$, $\Pi(k-1)=2$ and $\Upsilon(k-1)=3$.
As before the integral is
        \begin{multline*}
         \int_{T \times T} f(t) f(s) \!\left({ \dfrac{\partial ^2  r_\ve}{\partial t_{j}^2}(\gamma_1+t-s) }\right ) \!\left({ \dfrac{\partial   r_\ve}{\partial t_{j}}(\gamma_2+t-s) }\right ) \\ 
      \times \!\left({ \dfrac{\partial   r_\ve}{\partial t_{j}}(\gamma_3+t-s) }\right ) r_\ve(\gamma_4+t-s) \dots r_\ve(\gamma_k+t-s) \ud t \ud s.
      \end{multline*}
    There are $4k!(k-2)!(k-2)$ such cases.

     \item[Case 3]
All indices are different.
We can suppose that $\Pi(k)=1$, $\Pi(k-1)=2$, $\Upsilon(k)=3$ and $\Upsilon(k-1)=4$.
The integral is
      \begin{multline*}
         \int_{T \times T} f(t) f(s) \!\left({ \dfrac{\partial   r_\ve}{\partial t_{j}}(\gamma_1+t-s) }\right ) \!\left({ \dfrac{\partial   r_\ve}{\partial t_{j}}(\gamma_2+t-s) }\right )
          \!\left({ \dfrac{\partial   r_\ve}{\partial t_{j}}(\gamma_3+t-s) }\right )  \\ 
         \times\!\left({ \dfrac{\partial   r_\ve}{\partial t_{j}}(\gamma_4+t-s) }\right ) r_\ve(\gamma_5+t-s) \dots r_\ve(\gamma_k+t-s) \ud t \ud s.
      \end{multline*}
      There are $k!(k-2)!(k-2)(k-3)$ possible cases.
 \end{itemize}
We will consider in detail only one integral, that corresponding to the first case, that is%
\begin{multline*}
 L_\ve(\gamma_1, \gamma_2, \dots, \gamma_k)
	 :=  \int_{T \times T} f(t) f(s) \!\left({ \dfrac{\partial ^2  r_\ve}{\partial t_{j}^2}(\gamma_1+t-s) }\right ) \!\left({ \dfrac{\partial ^2  r_\ve}{\partial t_{j}^2}(\gamma_2+t-s) }\right )\\
	\times r_\ve(\gamma_3+t-s) \dots  r_\ve(\gamma_k+t-s) \ud t \ud s.
\end{multline*}
We will prove that $L_\ve(\bm{\cdot})$ converges uniformly on $\kappa$, where $\kappa$ is any compact set in $(\reels^d)^k$, as $\ve \to 0$, to $L_0(\bm{\cdot})$ defined by
\begin{multline*}
 L_{0}(\gamma_1, \gamma_2, \dots, \gamma_k)
 	:= \int_{T \times T} f(t) f(s) \!\left({ \dfrac{\partial ^2  r}{\partial t_{j}^2}(\gamma_1+t-s) }\right ) \!\left({ \dfrac{\partial ^2  r}{\partial t_{j}^2}(\gamma_2+t-s) }\right )\\
	 \times r(\gamma_3+t-s) \dots  r(\gamma_k+t-s) \ud t \ud s.
\end{multline*}
Since $f$ is bounded, we have
\setlength{\hauteur}{\heightof{$\ds\int_{T \times T}$}} 
{\def\termeA{ \abs{L_\ve(\gamma_1, \gamma_2, \dots, \gamma_k)- L_{0}(\gamma_1, \gamma_2, \dots, \gamma_k)}}
\def\termeB{\bC  \left\{\int_{T \times T}  \abs{r(0)}^{k-2} \right.}
\def\termeC{\left|
	\frac{\partial ^2  r_\ve}{\partial t_{j}^2}(\gamma_1+t-s)
	\dfrac{\partial ^2  r_\ve}{\partial t_{j}^2}(\gamma_2+t-s)\right.}
\def\termeD{\left.\dfrac{\partial ^2  r}{\partial t_{j}^2}(\gamma_1+t-s)
		\dfrac{\partial ^2  r}{\partial t_{j}^2}(\gamma_2+t-s) \right|}
\def\termeE{\int_{T \times T} 
\abs{\dfrac{\partial ^2  r}{\partial t_{j}^2}(\gamma_1+t-s)}\abs{ \dfrac{\partial ^2  r}{\partial t_{j}^2}(\gamma_2+t-s)}\left| r_\ve(\gamma_3+t-s) \dots  r_\ve(\gamma_k+t-s) \right.}
\def\termeF{
-\left.r(\gamma_3+t-s) \dots  r(\gamma_k+t-s)\right|}
\begin{align*}
\MoveEqLeft[2]{\termeA}\\ 
&\le \termeB \termeC  \\
&\specialpos{\hfill -\termeD \ud t \ud s}\\
&+\termeE\\
&\specialpos{\hfill \left. \termeF \ud t \ud s\rule{0pt}{\hauteur} \right\}}
\end{align*}}

The first integral above is bounded by   
 \begin{align*}
\MoveEqLeft[2]{\bC  
\left\{{
\int_{T \times T}  
  \abs{\dfrac{\partial ^2  r_\ve}{\partial t_{j}^2}(\gamma_1+t-s) -  
  \dfrac{\partial ^2  r}{\partial t_{j}^2}(\gamma_1+t-s)}\abs{ \dfrac{\partial ^2  r_\ve}{\partial t_{j}^2}(\gamma_2+t-s) }\ud t \ud s }\right.}\\
& \specialpos{\hfill +\left.{
  \int_{T \times T} \abs{ \dfrac{\partial ^2  r}{\partial t_{j}^2}(\gamma_1+t-s) }
   \abs{\dfrac{\partial ^2  r_\ve}{\partial t_{j}^2}(\gamma_2+t-s) -  
  \dfrac{\partial ^2  r}{\partial t_{j}^2}(\gamma_2+t-s)} \ud t \ud s
  }\right\}} \\
& \le \bC  
\left\{{
\!\left({
\left[{\int_{K} \!\left({\dfrac{\partial ^2  r_\ve}{\partial v_{j}^2}(v)}\right)^2 \ud v
}\right]^{\frac{1}{2}} +\left[{\int_{K} \!\left({\dfrac{\partial ^2  r}{\partial v_{j}^2}(v)}\right)^2 \ud v
}\right]^{\frac{1}{2}}
}\right) }\right.\\
&  \specialpos{\hfill \times \left.{
\left[{
\int_{K} \!\left({
\dfrac{\partial ^2  r_\ve}{\partial v_{j}^2}(v) -  
  \dfrac{\partial ^2  r}{\partial v_{j}^2}(v)
}\right)^2 \ud v
}\right]^{\frac{1}{2}}
}\right\}},
 \end{align*}
while the second integral is bounded by 
$$
 \bC  \times \sup_{\{v_3, \dots v_k \in K\}}\abs{r_\ve(v_3) \dots  r_\ve(v_k)
- r(v_3) \dots  r(v_k)} \times \left[{\int_{K} \!\left({\dfrac{\partial ^2  r}{\partial v_{j}^2}(v)}\right)^2 \ud v
}\right],
$$%
where $K$ is a compact in $\reels^d$.\\
We have seen in the proof of Proposition \ref{term S1} (see (\ref{limite pour derivees r})) that  for $v \neq 0$, 
$$
   \lim_{\ve \to 0} 
      \dfrac{\partial ^{2}r_\ve}{\partial v_{j}^2}
(v)=  \dfrac{\partial ^{2}r}{\partial v_{j}^2}(v).
$$%
Moreover, we have emphasized (see (\ref{majoration derivees secondes r epsilon})) that
   for $v \in K$ and for $j =1,\dots,d$, $\normp[d]{v} > M \ve$,
$$
      \abs{\dfrac{\partial ^{2}r_\ve}{
   \partial v_{j}^2}(v)}\le  \bC  \normp[d]{v}^{2\alpha -2}.
$$%
Moreover, a simple calculation shows that when $\normp[d]{v} \le M \ve$, we have
$$
     \abs{\dfrac{\partial ^{2}r_\ve}{
   \partial v_{j}^2}(v)}\le  \bC   \ve^{2\alpha -2}.
$$%
Furthermore $\dfrac{\partial ^{2}r}{\partial v_{j}^2}(v)=O\!\left(\normp[d]{v}^{2 \alpha -2}\right)$, for any $v$ in a compact set.\\
Now, recalling that $\alpha > \frac{1}{2}$, we have $d+4\alpha -4 >0$.\\
These arguments justify the use of Lebesgue's theorem which implies that
$$
   \int_{K} \!\left({
\dfrac{\partial ^2  r_\ve}{\partial v_{j}^2}(v) -  
  \dfrac{\partial ^2  r}{\partial v_{j}^2}(v)
}\right)^2 \ud v \cvg[\ve \to 0]{} 0.
$$%
Therefore, the two terms tend to zero as $\ve \to 0$.
Thus $\lim_{\ve \to 0} \chi_\ve(\bm{\cdot})$ is continuous on $(\reels^d)^k$ and in particular at point $\0$ of $(\reels^d)^k$.
Therefore, there exists a finite measure $\mu$ such that $\widehat{\mu}_\ve \to \widehat{\mu}$ and, by Levy's theorem, $\mu_\ve \to \mu$ weakly.

This yields the lemma.
\end{proofarg}

So to finish with the proof of the theorem, it remains to prove items \ref{itm3:theo4.6.20}
and \ref{itm4:theo4.6.20}.\\
Consider the proof of part \ref{itm3:theo4.6.20} of the theorem.
We suppose that  $\alpha = {1}/{2}$  and $d>2$.\\
Remember that we defined the \rv
$\xi_\ve$ as follows:
for $t \in T$, $\xi_\ve(t):= [{X_\ve(t)-X(t)}]/\ve$, $T$ being an open rectangle.\\
Also remember that $g^{\prime\prime}$ belongs to $L^4(\varphi(x) \ud x)$.\\
Finally, to finish with some notation reminders, we defined the constant $\bdK_{1/2}$ by:
$$
	\bdK_{1/2} = \bC _L \int_{\reels^{d}}   \Psi(w) \normp[d]{w}  \ud w.
$$%
As already seen previously, we have split the \rv
$  \Xi_\ve(f, g)$ into two terms, 
$$
 \Xi_\ve(f, g)=T_1+T_2,
$$%
   where
$$
	T_1=\frac{1}\ve \int_T f(t)  g\!\left({X_\ve(t)}/{\sqrt{r_\ve(0)}}\right)  h\!\left({\nabla\!X_\ve(t)}/{\sqrt{\mu_\ve}}\right) 
   \ud t 
 $$%
   and
$$
	T_2=\frac{1}\ve \int_T f(t)    \left\{{g\!\left({X_\ve(t)}/{\sqrt{r_\ve(0)}}\right) - g\!\left({{X(t)}/{\sqrt{r(0)}}}\right)
  }\right\} \ud t.
$$%
We know using Remark \ref{cas plus general} that $T_1\cvg[\ve\to0]{$L^{2}(\Omega)$}Y(f, g)$.\\
Thus the \rv
$  \Xi_\ve(f, g)$ is equivalent in $L^2(\Omega)$ to $Y(f, g) +T_2$, \ie
$$ 
\Xi_\ve(f, g) \simeq Y(f, g) +T_2.
$$%
Let us look at the term $T_2$.\\
Since $g$ is $C^2(\reels, \reels)$, $\forall t \in T$ we can write
 \begin{align*}
\MoveEqLeft[2]{g\!\left({X_\ve(t)}/{\sqrt{r_\ve(0)}}\right) - g\!\left({{X(t)}/{\sqrt{r(0)}}}\right)}\\
	& =  \!\left({{X_\ve(t)}/{\sqrt{r_\ve(0)}}-{{X(t)}/{\sqrt{r(0)}}}}\right)
 g^{\prime}\!\left({{X(t)}/{\sqrt{r(0)}}}\right)\\ 
 	& \quad+ \!\left({{X_\ve(t)}/{\sqrt{r_\ve(0)}}-{{X(t)}/{\sqrt{r(0)}}}}\right)^2\\
	& \quad \times  \int_0^1 (1-u)  
 g^{\prime\prime}\!\left({{X(t)}/{\sqrt{r(0)}}} +u 
 \!\left({
 {X_\ve(t)}/{\sqrt{r_\ve(0)}}-{{X(t)}/{\sqrt{r(0)}}}
 }\right)\right) \ud u.
 \end{align*}
 We split the \rv
${X_\ve(t)}/{\sqrt{r_\ve(0)}}-{{X(t)}/{\sqrt{r(0)}}}$ into two terms as follows:
$$
 {X_\ve(t)}/{\sqrt{r_\ve(0)}}- {X(t)}/{\sqrt{r(0)}}=  \frac{1}{\sqrt{r_\ve(0)}}\!\left({
 1-\sqrt{\frac{r_\ve(0)}{r(0)}}}\right)X(t)+\frac\ve{\sqrt{r_\ve(0)}} \xi_\ve(t).
$$%
With this decomposition, we then split $T_2$ into three terms thus obtaining
 \begin{align*}
 T_2 &=  \frac{1}{\sqrt{r_\ve(0)}} \frac{1}\ve \!\left({
 1-\sqrt{\frac{r_\ve(0)}{r(0)}}}\right) \int_T f(t) X(t)  g^{\prime}\!\left({{X(t)}/{\sqrt{r(0)}}}\right) \ud t \\
 & \quad+  \frac{1}{\sqrt{r_\ve(0)}} \int_T f(t) g^{\prime}\!\left({{X(t)}/{\sqrt{r(0)}}}\right) \xi_\ve(t) \ud t\\
 & \quad + \frac{1}\ve \int_T f(t) \!\left({
 \frac{1}{\sqrt{r_\ve(0)}} \!\left({
 1-\sqrt{\frac{r_\ve(0)}{r(0)}}}\right) X(t)+ \frac\ve{\sqrt{r_\ve(0)}} \xi_\ve(t)
 }\right)^2 \\
 & \quad\times \int_0^1 (1-u)  g^{\prime\prime}\!\left({{X(t)}/{\sqrt{r(0)}}} +u 
 \!\left({
 {X_\ve(t)}/{\sqrt{r_\ve(0)}}-{{X(t)}/{\sqrt{r(0)}}}
 }\right)\right) \ud u \ud t.
 \end{align*}
 Let us define the terms below:
$$
 n_\ve^2:= \bbE\left[\left(\frac{X_\ve(t)-X(t)}{\sqrt\ve}\right)^2\right] \cvg[\ve\to0]{} n_{{1}/{2}}^2:= 2r(0) \!\left({\frac{\bdK_{1/2}}{r(0)}-\ell_{{1}/{2}} }\right),
$$
since $\bdK_{1/2}$ is
$$
\bdK_{1/2}= \bC _L \int_{\reels^{d}}   \Psi(w) \normp[d]{w}  \ud w =\lim_{\ve \to 0} \!\left({\frac{r(0)-r^{(\ve)}(0)}\ve}\right).
$$%
With this notation, we introduce the {\rv}s $Z_\ve$ and $W_\ve(t, u)$ setting for all $t \in T$ and $\forall u \in \intrff{0}{1}$
$$
Z_\ve(t):= \ds \frac{X_\ve(t)-X(t)}{n_\ve \sqrt\ve} \hookrightarrow {\calN}(0;1),
$$%
and %
$$
W_\ve(t, u):= \ds {{X(t)}/{\sqrt{r(0)}}} +u 
 \!\left({
 {X_\ve(t)}/{\sqrt{r_\ve(0)}}-{{X(t)}/{\sqrt{r(0)}}}}\right) \hookrightarrow {\calN}(0; v_\ve^2(u)),
$$%
with $v_\ve^2(u):=2u^2c_\ve-2uc_\ve+1$, where
$$
  0<c_\ve:= 1- \ds \frac{r^{(\ve)}(0)}{\sqrt{r(0)}\sqrt{r_\ve(0)}} \cvg[\ve\to0]{} 0.
$$%

We obtain
   \begin{align*}
 T_2 &=  \frac{1}{\sqrt{r_\ve(0)}} \frac{1}\ve \!\left({
	 1-\sqrt{\frac{r_\ve(0)}{r(0)}}}\right) \int_T f(t) X(t)  g^{\prime}\!\left({{X(t)}/{\sqrt{r(0)}}}\right) \ud t \\
	&\quad +  \frac{1}{\sqrt{r_\ve(0)}} \int_T f(t) g^{\prime}\!\left({{X(t)}/{\sqrt{r(0)}}}\right) \xi_\ve(t) \ud t\\
 	&\quad + \frac{n_\ve^2}{r_\ve(0)} \int_T f(t) \int_0^1 (1-u) g^{\prime\prime}(W_\ve(t, u)) \ud u \ud t\\
	&\quad + \frac{n_\ve^2}{r_\ve(0)} \int_T f(t) H_2(Z_\ve(t)) \int_0^1 (1-u) g^{\prime\prime}(W_\ve(t, u)) \ud u \ud t\\
	&\quad + \frac{1}{r_\ve(0)} \frac{1}\ve \!\left({
 1-\sqrt{\frac{r_\ve(0)}{r(0)}}}\right)^2 \int_T f(t) X^2(t) \int_0^1 (1-u) g^{\prime\prime}(W_\ve(t, u)) \ud u \ud t\\
	&\quad + \frac{2}{r_\ve(0)} \!\left({
 1-\sqrt{\frac{r_\ve(0)}{r(0)}}}\right) \int_T f(t) X(t) \xi_\ve(t) \int_0^1 (1-u) g^{\prime\prime}(W_\ve(t, u)) \ud u \ud t.
 \end{align*}
 Let us now show that the last three terms converge in probability to zero and that the third term tends in $L^{2}(\Omega)$ to 
 $$
 \frac{n_{{1}/{2}}^2}{2 r(0)} \int_T f(t) g^{\prime\prime}\!\left({{X(t)}/{\sqrt{r(0)}}}\right) \ud t.
 $$
Note that this last term is 
$$
\!\left({\frac{\bdK_{1/2}}{r(0)}-\ell_{{1}/{2}} }\right) \int\limits_T f(t) g^{\prime\prime}\!\left({{X(t)}/{\sqrt{r(0)}}}\right) \ud t.
$$\\
 Let us start with the last convergence proposition.\\
Since $n^2_\ve \to n^2_{{1}/{2}}$, we compute the second order moment of the following \rv
minus its supposed limit, getting
\begin{align*}
\MoveEqLeft[2]{\bbE\left[\left\{\int_T f(t) \int_0^1 (1-u) g^{\prime\prime}(W_\ve(t, u)) \ud u \ud t- \frac{1}{2} \int_T f(t) g^{\prime\prime}\!\left({{X(t)}/{\sqrt{r(0)}}}\right) \ud t\right\}^2\right]}\\
& =   \int_{T \times T} f(t) f(s) \int_{\intrff{0}{1}^2} (1-u) (1-v)\\
& \specialpos{\hfill\times  \bbE\left[\rule{0pt}{18pt}\!\left({g^{\prime\prime}(W_\ve(t, u))-g^{\prime\prime}\!\left({{X(t)}/{\sqrt{r(0)}}}\right)}\right)\right.\hfill}\\
&  \specialpos{\hfill\left.\times\!\left({g^{\prime\prime}(W_\ve(s, v))-g^{\prime\prime}\!\left({X(s)}/{\sqrt{r(0)}}\right)}\right)\right] \ud u \ud v \ud t \ud s}\\
&   \le  \bC  \!\left({
   \int_0^1
   \!\left({\bbE\!\left[\left\{g^{\prime\prime}(W_\ve(0, u))-g^{\prime\prime}\!\left({X(0)}/{\sqrt{r(0)}}\right)\right\}^2\right]}\right)^{\frac{1}{2}} \ud u
   }\right)^{2}\\
  & \le  \bC  
   \int_0^1
   \bbE\!\left[\left\{g^{\prime\prime}(W_\ve(0, u))-g^{\prime\prime}\!\left({X(0)}/{\sqrt{r(0)}}\right)\right\}^2\right] \ud u.
  \end{align*}
Let us justify that Lebesgue's theorem can be applied.
On the one hand, since $g^{\prime\prime}$ is $L^2(\varphi(x) \ud x)$, it trivially follows that $\bbE\!\left[\left\{g^{\prime\prime}\!\left({X(0)}/{\sqrt{r(0)}}\right)\right\}^2\right] \le  \bC $.\\
On the other hand, using that $\frac{1}{2} \le 1 -\frac{1}{2}c_\ve \le v^2_\ve(u) \le 1$, $\forall u \in \intrff{0}{1}$ and $\ve \le \ve_0$, we obtain that 
\begin{equation}
\label{inequality g(W)}
\bbE\!\left[\left\{g^{\prime\prime}(W_\ve(0, u))\right\}^2\right] \le  \bC   \bbE\!\left[\left\{g^{\prime\prime}(N_1)\right\}^2\right]  \le  \bC .
\end{equation}
So we can interchange the limit in $\ve$ and the integral.
 Therefore it only remains to show that $\forall u \in \intrff{0}{1}$, one has 
 \begin{equation}
 \label{convergence g(W)}
 \lim_{\ve \to 0} \bbE\!\left[\left\{g^{\prime\prime}(W_\ve(0, u))-g^{\prime\prime}\!\left({X(0)}/{\sqrt{r(0)}}\right)\right\}^2\right]=0.
 \end{equation}
To do so, we will apply Scheffé's lemma.
Let $u$ fixed in $\intrff{0}{1}$.\\
First, note that $W_\ve(0, u)$ almost surely converges to ${X(0)}/{\sqrt{r(0)}}$ and since $g$ is $C^2(\reels, \reels)$, we get that $g^{\prime\prime}(W_\ve(0, u))$ almost surely converges to $g^{\prime\prime}\!\left({X(0)}/{\sqrt{r(0)}}\right)$.\\
Using that $\frac{1}{2} \le v^2_\ve(u) \le 1$ and that $v^2_\ve(u) \cvg[\ve\to0]{} 1$, we have 
$$
\bbE\!\left[\left\{g^{\prime\prime}(W_\ve(0, u))\right\}^2\right]\cvg[\ve\to0]{} \bbE\!\left[\left\{g^{\prime\prime}\!\left({X(0)}/{\sqrt{r(0)}}\right)\right\}^2\right].
$$
 All the ingredients are gathered to apply the Scheffé's lemma.
We conclude that $g^{\prime\prime}(W_\ve(0, u))$ tends in $L^2(\Omega)$ to $g^{\prime\prime}\!\left({X(0)}/{\sqrt{r(0)}}\right)$ which is the required result.\\
 Let us now show that the last three terms tend in probability  towards zero.\\
Consider $A_\ve$, the first term of the sum.
$$
 A_\ve:=\frac{1}{r_\ve(0)} \frac{1}\ve \!\left({
 1-\sqrt{\frac{r_\ve(0)}{r(0)}}}\right)^2 \int_T f(t) X^2(t) \int_0^1 (1-u) g^{\prime\prime}(W_\ve(t, u)) \ud u \ud t.
$$%
 Since $ \frac{1}\ve\!\left({
 1-\sqrt{\frac{r_\ve(0)}{r(0)}}}\right) \cvg[\ve\to0]{} \ell_{{1}/{2}}$, it is only necessary to prove that the \rv
 $$
 \int_T f(t) X^2(t) \int_0^1 (1-u) g^{\prime\prime}(W_\ve(t, u)) \ud u \ud t
 $$%
 is bounded in $L^{1}(\Omega)$.\\
 Therefore we are bounding the following expectation:
  \begin{align*}
\MoveEqLeft[2]{\bbE\left[\abs{\int_T f(t) X^2(t) \int_0^1 (1-u) g^{\prime\prime}(W_\ve(t, u)) \ud u \ud t}\right]}\\
  &\le  \bC  \int_0^1  \bbE\!\left[\abs{X^2(0) g^{\prime\prime}(W_\ve(0, u))}\right] \ud u\\
  &\le \bC  \!\left({\bbE\!\left[X^4(0)\right]}\right)^{\frac{1}{2}} \int_0^1 \!\left({\bbE\!\left[\left\{g^{\prime\prime}(W_\ve(0, u))\right\}^2\right]}\right)^{\frac{1}{2}} \ud u\\
  &\le\bC ,
   \end{align*}
 the last bound is follows from (\ref{inequality g(W)}).\\
Consider $B_\ve$, the second term of the sum.
 \begin{multline*}
 B_\ve:=
 \frac{2}{r_\ve(0)} \!\left({
 1-\sqrt{\frac{r_\ve(0)}{r(0)}}}\right) \int_T f(t) X(t) \xi_\ve(t) \int_0^1 (1-u) g^{\prime\prime}(W_\ve(t, u)) \ud u \ud t \\
= \frac{2 n_\ve}{r_\ve(0)} \frac{1}{\sqrt\ve} \!\left({
 1-\sqrt{\frac{r_\ve(0)}{r(0)}}}\right) \int_T f(t) X(t) Z_\ve(t) \int_0^1 (1-u) g^{\prime\prime}(W_\ve(t, u)) \ud u \ud t.
  \end{multline*}
  Since 
  $$
  \frac{1}{\sqrt\ve}\!\left({
 1-\sqrt{\frac{r_\ve(0)}{r(0)}}}\right) \cvg[\ve\to0]{} 0
 \quad\text{and}\quad
 n_\ve \cvg[\ve\to0]{} n_{{1}/{2}}
 $$
it is sufficient, as for the previous term $A_\ve$, to prove that the \rv
$$
	\int_T f(t) X(t) Z_\ve(t) \int_0^1 (1-u) g^{\prime\prime}(W_\ve(t, u)) \ud u \ud t
$$%
is bounded in $L^{1}(\Omega)$.
As before we bound the following expectation.
\begin{align*}
\MoveEqLeft[2]{\bbE\!\left[\abs{\int_T f(t) X(t) Z_\ve(t) \int_0^1 (1-u) g^{\prime\prime}(W_\ve(t, u)) \ud u \ud t}\right]}\\
	&\le \bC  \!\left({\bbE[X^4(0)]}\right)^{\frac{1}{4}}  \!\left({\bbE[Z_\ve^4(0)]}\right)^{\frac{1}{4}} \int_0^1 \!\left({\bbE\!\left[\left\{g^{\prime\prime}(W_\ve(0, u))\right\}^2\right]}\right)^{\frac{1}{2}} \ud u\\
 	&\le\bC .
\end{align*}
Let us consider the last and more difficult term.
 Since $n^2_\ve \to n^2_{{1}/{2}}$, it suffices to prove the convergence to zero of the following term $C_\ve$.
 \begin{align*}
 C_\ve&:=\int_T f(t) H_2(Z_\ve(t)) \int_0^1 (1-u) g^{\prime\prime}(W_\ve(t, u)) \ud u \ud t\\
 &= D_\ve +F_\ve,
 \end{align*}
 where 
$$
	 D_\ve:=\int_T f(t) H_2(Z_\ve(t)) \int_0^1 (1-u) \!\left({g^{\prime\prime}(W_\ve(t, u)) -g^{\prime\prime}\!\left({{X(t)}/{\sqrt{r(0)}}}\right)}\right) \ud u \ud t,
$$%
 and
$$
	F_\ve:=\frac{1}{2} \int_T f(t) H_2(Z_\ve(t)) g^{\prime\prime}\!\left({{X(t)}/{\sqrt{r(0)}}}\right) \ud u \ud t.
$$%
 Let us first consider the term $D_\ve$.
\begin{align*}
\MoveEqLeft[4]{\bbE\!\left[\abs{D_\ve}\right] \le \bC  \!\left({\bbE[H_2^2(Z_\ve(0))]}\right)^{\frac{1}{2}}}\\
  & \times\int_0^1 \!\left({\bbE\!\left[\left\{g^{\prime\prime}(W_\ve(0, u))-g^{\prime\prime}\!\left({X(0)}/{\sqrt{r(0)}}\right)\right\}^2\right]}\right)^{\frac{1}{2}} \ud u
   \cvg[\ve\to0]{} 0.
\end{align*}
The last convergence is a consequence of (\ref{inequality g(W)}) and (\ref{convergence g(W)}).
 Finally, it remains to consider the term $F_\ve$ for which we will prove that $\bbE[F_\ve^2]$ tends to zero with $\ve$.
So let us compute $\bbE[F_\ve^2]$.
 \begin{multline*}
 \bbE\left[F_\ve^2\right]=
 \frac{1}{4} \int_{T \times T} f(t) f(s) 
 	\bbE\left[H_2(Z_\ve(t)) g^{\prime\prime}\!\left({{X(t)}/{\sqrt{r(0)}}}\right)
	\right.\\
\left. \times H_2(Z_\ve(s))
		g^{\prime\prime}\!\left({X(s)}/{\sqrt{r(0)}}\right)\right]
\ud t \ud s.
 \end{multline*}
 We divide the domain of integration $T \times T$ into two parts, namely $T_{\delta}^{(1)}$ and $T_{\delta}^{(2)}$:
\begin{align*}
 T_{\delta}^{(1)}&:=\{(s, t) \in T \times T: \normp[d]{t-s} \le \delta \}\\ \shortintertext{and}
 T_{\delta}^{(2)}&:=\{(s, t) \in T \times T: \normp[d]{t-s} > \delta \}
\end{align*}
where $\delta >0$ is chosen such that
$$
  \sup_{\left\{s, t \in T, \normp[d]{t-s} \ge \delta \right\}}\left\{ \abs{{r(t-s)}/{r(0)}}\right\}  \le \rho <1.
$$%
The corresponding integrals are denoted by $f_\ve^{(1)}$ and $f_\ve^{(2)}$.\\
First, we bound the term $f_\ve^{(1)}$.\\
Using the fact that $g^{\prime\prime}$ is $L^4(\varphi(x) \ud x)$, we easily obtain that for all $s, t \in T$,  
\begin{align}
	\label{th CD}
	\abs{ \bbE\left[H_2(Z_\ve(t)) g^{\prime\prime}\!\left({{X(t)}/{\sqrt{r(0)}}}\right)H_2(Z_\ve(s)) g^{\prime\prime}\!\left({X(s)}/{\sqrt{r(0)}}\right)\right]}\\ \nonumber
	\le \!\left({\bbE[H_2^4(Z_\ve(0))]}\right)^{\frac{1}{2}} \!\left({
\bbE\!\left[\left\{g^{\prime\prime}\!\left({X(0)}/{\sqrt{r(0)}}\right)\right\}^4\right]
}\right)^{\frac{1}{2}} \le  \bC .
\end{align}
In this way, we obtain the following bound
$$
	\abs{f_\ve^{(1)}} \le \bC  \sigma_{2d}(T_{\delta}^{(1)})
	\le \bC  \delta^{d}.
$$%
We are now interested in the second term $f_\ve^{(2)}$.\\
Using the bound given in (\ref{th CD}), we can apply Lebesgue's dominated convergence theorem and interchange the limit with the integral.
\\
 Let us therefore fix a point $(s, t)$ in $T_{\delta}^{(2)}$ and prove that,
$$
\lim_{\ve \to 0} \bbE\!\left[H_2(Z_\ve(t)) g^{\prime\prime}\!\left({{X(t)}/{\sqrt{r(0)}}}\right)H_2(Z_\ve(s)) g^{\prime\prime}\!\left({X(s)}/{\sqrt{r(0)}}\right)\right]=0,
$$%
thus achieving our goal.
 We compute the last expectation:
 \begin{align*}
 \MoveEqLeft[2]{ \bbE\!\left[H_2(Z_\ve(t)) g^{\prime\prime}\!\left({{X(t)}/{\sqrt{r(0)}}}\right)H_2(Z_\ve(s)) g^{\prime\prime}\!\left({X(s)}/{\sqrt{r(0)}}\right)\right]}\\
	& = \sum_{k=2}^{\infty} \sum_{\ell=2}^{\infty} c_k c_\ell k(k-1)\ell (\ell-1) \\
	& \quad\times \bbE\!\left[H_2\!\left(Z_\ve(t)\right) H_2\!\left(Z_\ve(s)\right) H_{k-2}\!\left({{X(t)}/{\sqrt{r(0)}}}\right) H_{\ell-2}\!\left({X(s)}/{\sqrt{r(0)}}\right)\right].
 \end{align*}
 By writing Mehler's formula \cite{MR716933} and using the bound 
 $$
	 \abs{\bbE\!\left[{{X(t)}/{\sqrt{r(0)}}} \cdot {X(s)}/{\sqrt{r(0)}}\right]} \le \rho <1,
 $$
 it is easy to obtain the following upper bound
  \begin{align*}
 \MoveEqLeft[2]{\abs{\bbE\!\left[H_2\!\left(Z_\ve(t)\right) H_2\!\left(Z_\ve(s)\right) H_{k-2}\!\left({{X(t)}/{\sqrt{r(0)}}}\right) H_{\ell-2}\!\left({X(s)}/{\sqrt{r(0)}}\right)\right]}}   \\
	& \le \bC   \left\{{
 (k-2)(k-3)(k-2)!\rho^{k-4} \1_{\{k=\ell\}} +(k-2)!(k-4)\rho^{k-5}\1_{\{k-2=\ell\}}}\right.\\
	& \quad+\left.{ (\ell-2)!(\ell-4)\rho^{\ell-5}\1_{\{\ell-2=k\}}
 }\right\}.
 \end{align*}
 Thus, using H\"{o}lder's inequality we get
  \begin{align*}
  \MoveEqLeft[2]{\sum_{k=2}^{\infty} \sum_{\ell=2}^{\infty} \abs{c_k c_\ell k(k-1)\ell (\ell-1)}}\\
	& \abs{\bbE\!\left[H_2\!\left(Z_\ve(t)\right) H_2\!\left(Z_\ve(s)\right) H_{k-2}\!\left({{X(t)}/{\sqrt{r(0)}}}\right) H_{\ell-2}\!\left({X(s)}/{\sqrt{r(0)}}\right)\right]}  \\
	&  \le \bC  \left\{{
 \sum_{k=4}^{\infty} c_k^2 k! k(k-1) (k-2)(k-3)\rho^{k-4}
}\right.\\
	& \specialpos{\hfill+\left.{ 
\!\left({
 \sum_{k=5}^{\infty} c_k^2 k! (k-2)(k-3) (k-4)\rho^{k-5}}\right)^{\frac{1}{2}} }\right.\hfill}\\
	&   \specialpos{\hfill\times \left.{
 \!\left({
 \sum_{k=5}^{\infty} c_{k-2}^2 k! (k-2)(k-3) (k-4)\rho^{k-5}}\right)^{\frac{1}{2}}
 }\right\}} 
   \\
	&  \le \bC  \left\{{
 \sum_{k=4}^{\infty} c_k^2 k! k(k-1) (k-2)(k-3)\rho^{k-4}
}\right.\\
	& \specialpos{\hfill+\left.{ 
\!\left({
 \sum_{k=5}^{\infty} c_k^2 k! (k-2)(k-3) (k-4)\rho^{k-5}}\right)^{\frac{1}{2}}}\right.\hfill}\\
 	& \specialpos{\hfill\times\left.{
 \!\left({
 \sum_{k=3}^{\infty} c_{k}^2 k! (k+2)(k+1) k (k-1)(k-2)\rho^{k-3}}\right)^{\frac{1}{2}}
 }\right\} }
  \\
	&   \le \bC .
 \end{align*}
The last bound comes from the following argument.\\
First, $ \sum_{k=0}^{\infty}   c_k^2 k!  < \infty$ (see (\ref{finitude cn})).
Then for all $0 \le \rho <1$, $\sum_{k=0}^{\infty}   c_k^2 k! \rho^k < \infty$ and this series is indefinitely differentiable.
\\
Furthermore, for all $m \in \naturels^{\star}$, we have
$$
 \left(\sum_{k=0}^{\infty}   c_k^2 k! \rho^k\right)^{(m)} = \sum_{k=m}^{\infty}   c_k^2 k! \prod_{i=1}^{m}(k-i+1) \rho^{k-m}< \infty.
$$%
And in the same way
$$
 \sum_{k=3}^{\infty}   c_k^2 k! \prod_{i=0}^{m+2}(k+i-2) \rho^{k-3}< \infty.
$$%
So we have proved that we can interchange the limit with the series.
It remains to show that for fixed integers $k, \ell \ge 2$, we have
$$
   \lim_{\ve \to 0} \bbE\!\left[
   H_2\!\left(Z_\ve(t)\right) H_2\!\left(Z_\ve(s)\right) H_{k-2}\!\left({{X(t)}/{\sqrt{r(0)}}}\right) H_{\ell-2}\!\left({X(s)}/{\sqrt{r(0)}}\right)\right]=0,
$$%
and our final goal will be achieved.\\
On the one hand the distribution of the Gaussian vector 
$$
\!\left({Z_\ve(t), Z_\ve(s), \ds {{X(t)}/{\sqrt{r(0)}}}, {X(s)}/{\sqrt{r(0)}}}\right)
$$%
converges to that of the Gaussian vector
$$
\!\left({Z(t), Z(s), \ds {{X(t)}/{\sqrt{r(0)}}}, {X(s)}/{\sqrt{r(0)}}}\right)
$$%
 with covariance matrix 
$$
\begin{pmatrix}
1&0&0&0\\
0&1&0&0\\
0&0&1& \ds \frac{r(t-s)}{r(0)}\\
0&0&\ds \frac{r(t-s)}{r(0)}&1
\end{pmatrix}.
$$
On the other hand, by Lemma 3.1 of \cite{MR471045}
\begin{align*}
\MoveEqLeft[2]{\bbE\!\left[H_2^2\!\left(Z_\ve(t)\right) H_2^2\!\left(Z_\ve(s)\right) H_{k-2}^2\!\left({{X(t)}/{\sqrt{r(0)}}}\right) H_{\ell-2}^2\!\left({X(s)}/{\sqrt{r(0)}}\right)\right]}\\
& \le 4 \cdot 7^{k+\ell} \cdot (k-2)! (\ell-2)! \le \bC .
\end{align*}
These two facts imply that 
\begin{align*}
\MoveEqLeft[3]{\bbE\!\left[H_2\!\left(Z_\ve(t)\right) H_2\!\left(Z_\ve(s)\right) H_{k-2}\!\left({{X(t)}/{\sqrt{r(0)}}}\right) H_{\ell-2}\!\left({X(s)}/{\sqrt{r(0)}}\right)\right]}  \\
\cvg[\ve\to0]{}&{~}\bbE\!\left[H_2(Z(t)) H_2(Z(s)) H_{k-2}\!\left({{X(t)}/{\sqrt{r(0)}}}\right) H_{\ell-2}\!\left({X(s)}/{\sqrt{r(0)}}\right)\right] \\
=&{~}\bbE\!\left[H_2(Z(t)) H_2(Z(s))\right] \bbE\!\left[H_{k-2}\!\left({{X(t)}/{\sqrt{r(0)}}}\right) H_{\ell-2}\!\left({X(s)}/{\sqrt{r(0)}}\right)\right]\\
=&{~}\bbE\!\left[H_2(Z(t))\right] \bbE\!\left[H_2(Z(s))\right] \bbE\!\left[H_{k-2}\!\left({{X(t)}/{\sqrt{r(0)}}}\right) H_{\ell-2}\!\left({X(s)}/{\sqrt{r(0)}}\right)\right]\\
=&{~}0.
\end{align*}
In short, we proved that $\lim_{\ve \to 0} f_\ve^{(2)}=0$.\\
Therefore, we obtain the following bound.
For any $\delta >0$,
$$
\underset{\ve \to 0}{\overline{\lim}} \bbE[F_\ve^2]  \le \bC  \delta^{d}.
$$%
By taking limit when $\delta \to 0$, we obtain the required result, \ie
$$
\lim_{\ve \to 0} \bbE[F_\ve^2]=0.
$$
 
 To summarize all of this and by using that 
 $$
 \frac{1}{\sqrt{r_\ve(0)}} \frac{1}\ve \!\left({
 1-\sqrt{\frac{r_\ve(0)}{r(0)}}}\right)
 \cvg[\ve\to0]{} \frac{\ell_{1/2}}{\sqrt{r(0)}},
 $$
 we have shown that 
  \begin{align*}
 \MoveEqLeft[2]{ \Xi_\ve(f, g) \simeq Y(f, g) +T_2
 \simeq  Y(f, g)+ \frac{\ell_{{1}/{2}}}{\sqrt{r(0)}} \int_T f(t) X(t)  g^{\prime}\!\left({{X(t)}/{\sqrt{r(0)}}}\right) \ud t 
}\\ 
	& \specialpos{\hfill+  \frac{1}{\sqrt{r_\ve(0)}} \int_T f(t) g^{\prime}\!\left({{X(t)}/{\sqrt{r(0)}}}\right) \xi_\ve(t) \ud t\hfill} \\
	& \specialpos{\hfill\!\left({\frac{\bdK_{1/2}}{r(0)}-\ell_{{1}/{2}} }\right) \int_T f(t) g^{\prime\prime}\!\left({{X(t)}/{\sqrt{r(0)}}}\right) \ud  t}\\
 	& =Y(f, g)+ \ell_{{1}/{2}} 
		\int_T f(t) \!
			\left(
				{{X(t)}/{\sqrt{r(0)}}}  
				  g^{\prime}\!\left(
				  	{{X(t)}/{\sqrt{r(0)}}}
					\right)
						\right.\\
&						\specialpos{\hfill\left.
				- g^{\prime\prime}\!\left(
					{{X(t)}/{\sqrt{r(0)}}}
					\right)
			\right)
	 \ud t} \\
	& \specialpos{\hfill+  \frac{1}{\sqrt{r_\ve(0)}} \int_T f(t) g^{\prime}\!\left({{X(t)}/{\sqrt{r(0)}}}\right) \xi_\ve(t) \ud t\hfill} \\
	& \specialpos{\hfill+ \frac{\bdK_{1/2}}{r(0)} \int_T f(t) g^{\prime\prime}\!\left({{X(t)}/{\sqrt{r(0)}}}\right) \ud t}\\
	&  = X(f, g)+Y(f, g)+ \frac{\bdK_{1/2}}{r(0)} \int_T f(t) g^{\prime\prime}\!\left({{X(t)}/{\sqrt{r(0)}}}\right) \ud t\\
	&  \specialpos{\hfill+  \frac{1}{\sqrt{r_\ve(0)}} \int_T f(t) g^{\prime}\!\left({{X(t)}/{\sqrt{r(0)}}}\right) \xi_\ve(t) \ud t.}
 \end{align*}
From now on, the reader is referred to Appendix \ref{chap convergence vague}.
In the remaining part of this section, we use the same notations and choose ${\calH }:=L^2(T)$ \ie the set of continuous square-integrable functions on $T$.\\
To finish with the proof of the theorem, it remains to show that the \rv
$$
	Z_\ve(f):=\int_T f(t) g^{\prime}\!\left({{X(t)}/{\sqrt{r(0)}}}\right) \xi_\ve(t) \ud t
$$%
stably and vaguely converges to
$$
	Z(f):= \int_T f(t) g^{\prime}\!\left({{X(t)}/{\sqrt{r(0)}}}\right) \xi(t) \ud t.
$$%
Indeed, the stably and vaguely convergence of the \rv
$ \Xi_\ve(f, g)$  will be ensured by the application of Theorem 1 of \cite{MR517416} (see also the Appendix \ref{stable convergence} of Appendix \ref{chap convergence vague}).\\
  The random vector $\xi_\ve: \Omega \to {\calH }$ is a centered Gaussian vector, since for all $\phi \in {\calH }$, the real \rv
$\prodsca{\phi}{\xi_\ve}_{\calH }$ has a Gaussian distribution, ${\calN}(0; \prodsca{{\calK}_{\xi_\ve} \phi}{\phi}_{\calH })$, where the linear operator ${\calK}_{\xi_\ve}$ is
$$
  {\calK}_{\xi_\ve}(\phi)(t):= \int\limits_{T} \frac{r_\ve(t-s)-2 r^{(\ve)}(t-s)+r(t-s)}{\ve^2}   \phi(s) \ud s.
$$%
Moreover, for any $\phi \in {\calH }$, by applying Lemma \ref{majoration pour r(epsilon)} we can see
$$
\prodsca{{\calK}_{\xi_\ve} \phi}{\phi}_{\calH } \cvg[\ve\to0]{} \prodsca{{\calK}_{\xi} \phi}{\phi}_{\calH },
$$%
where
$$
  {\calK}_{\xi}(\phi)(t):= \int_{T}  \bC _{\Psi} \sum_{i=1}^d  \dfrac{\partial ^2  r}{\partial t_{i}^2}(t-s)   \phi(s) \ud s.
$$%
  We deduce by using Section \ref{vague convergence} of Appendix \ref{chap convergence vague} that
$$
	\xi_\ve \cvg[\ve\to0]{v} \xi.
$$%
Returning to our objective, we discretize the \rv $Z_\ve(f)$ in the following way.\\
The idea is to consider $\overline T=\cup_{i\in I_m}\overline T_{i,m}$ with $T_{i,m}\cap T_{j,m}=\phi$, and $\card(I_m)=i_m$.
Each $T_{i, m}$ is an open rectangle such that $\sigma_d(T_{i, m}) \cvg[m\to\infty]{} 0$.
Let denote $t_{i,m}$ the center of $T_{i,m}$ in this manner we can approximated $Z_\ve(f)$ by $Z_{\ve,m}(f)$ defined as follows
\begin{align*}
Z_{\ve,m}(f)&:=\sum_{j\in I_m}f(t_{j,m})g^{\prime}\!\left({X(t_{j,m})}/{\sqrt{r(0)}}\right)\int_{T_{j,m}}\xi_\ve(t)\ud t\\
&\phantom{:}=\sum_{j\in I_m}f(t_{j,m})g^{\prime}\!\left({X(t_{j,m})}/{\sqrt{r(0)}}\right)\prodsca{\1 _{T_{j,m}}}{\xi_\ve}_{\calH }\\
&\phantom{:}=\prodsca{\ds\sum_{j\in I_m}f(t_{j,m})g^{\prime}\!\left({X(t_{j,m})}/{\sqrt{r(0)}}\right)\1 _{T_{j,m}}}{\xi_\ve }_{\calH }\\
&:=\prodsca{F_m(X)}{\xi_\ve}_{\calH }.
\end{align*}
As a first step let us show that, $Z_{\ve, m}(f) \cvg[\ve\to0]{${\calS_{v}}$} Z_m(f):=\prodsca{F_m(X)}{\xi}_{\calH }$.\\
We will see the convergence of the following characteristic function
   \begin{align*}
  \MoveEqLeft[2]{\bbE\!\left[e^
  		{is \prodsca{F_m(X)}{\xi_\ve}_{\calH }}
	\right]}\\
	& =   \int_{\reels^{i_m}}
    \bbE\!\left[\exp
    	{\!\left(
		is \prodsca{F_m(X)}{\xi_\ve}_{\calH }
	\right)}
	~|\cap_{j \in I_m} \{X(t_{j,m})=x_{j, m}\}\right] \\
	& \specialpos{\hfill\ud \bbP_{(X(t_{j, m}), j \in I_m)}(x_{j, m}, j \in I_m)} \\
	& =\int_{\reels^{i_m}}
    \bbE\!\left[\exp
    	{\left(
	is \sum\limits_{j\in I_m}f(t_{j,m})g^{\prime}({x_{j,m}}{\sqrt{r(0)}})\int_{T_{j,m}}\xi_\ve(t)\ud t\ 
	\right)}\right.\\
	& \specialpos{\hfill\left.\left|\rule{0pt}{18pt}\cap_{j \in I_m} \!\left\{X(t_{j,m})=x_{j, m}\right\}\right.\right] \hfill}\\
	& \specialpos{\hfill\ud \bbP_{(X(t_{j, m}), j \in I_m)}(x_{j, m}, j \in I_m).}
\end{align*}
   The \rv
$W_{\ve, m}:=\sum\limits_{j\in I_m}f(t_{j,m})g^{\prime}\!\left({x_{j,m}}/{\sqrt{r(0)}}\right)\int_{T_{j,m}}\xi_\ve(t)\ud t$ has a Gaussian distribution so that
\begin{align*}
\MoveEqLeft[2]{{\calL}\!\left({W_{\ve, m}\given{\cap_{j \in I_m} \{X(t_{j,m})=x_{j, m}\}}}\right)} \\
	& = {\calN}\!\left({(x_{j, m}, j \in I_m)\times V_{22}^{-1} \times V_{21;\ve; m}; V_{11;\ve;m}-V_{21;\ve;m}^T \times V_{22}^{-1}\times V_{21;\ve; m}}\right),
 \end{align*}
 where by applying Lemma \ref{majoration pour r(epsilon)}
 \begin{align*}
 V_{11;\ve;m}
 	& := \sum\limits_{i, j\in I_m}f(t_{i,m})f(t_{j,m})g^{\prime}\!\left({x_{i,m}}/{\sqrt{r(0)}}\right)g^{\prime}\!\left({x_{j,m}}/{\sqrt{r(0)}}\right)\\
	& \specialpos{\hfill\times\int_{T_{i,m}\times T_{j, m}} \frac{r_\ve(t-s)-2 r^{(\ve)}(t-s)+r(t-s)}{\ve^2} \ud t \ud s}\\
	 \cvg[\ve\to0]{}
 V_{11;m}&:=
 \sum\limits_{i, j\in I_m}f(t_{i,m})f(t_{j,m})g^{\prime}\!\left({x_{i,m}}/{\sqrt{r(0)}}\right)g^{\prime}\!\left({x_{j,m}}/{\sqrt{r(0)}}\right)\\
 	& \specialpos{\hfill\times\int_{T_{i,m}\times T_{j, m}}\bC _{\Psi} \sum_{k=1}^d  \dfrac{\partial ^2  r}{\partial t_{k}^2}(t-s)  \ud t \ud s,}\\ \shortintertext{and}
 V_{21;\ve;m}^T&:=\sum\limits_{i\in I_m}f(t_{i,m})g^{\prime}\!\left({x_{i,m}}/{\sqrt{r(0)}}\right)\\
 	& \specialpos{\hfill\times\!\left({\int_{T_{i,m}}
 \frac{r^{(\ve)}(t-t_{j,m})-r(t-t_{j,m})}\ve \ud t}\right)_{j \in I_m}}\\
 &\qquad \cvg[\ve\to0]{} \0_{\reels^{i_m}},
 \end{align*}
 and $V_{22}$ is the variance matrix of $(X(t_{i, m}), i \in I_m)\in \reels^{i_m}$.\\
  From now on we note ${\calL}_{X}$ the distribution of the random variable $(X(t))_{t \in T}$.
We have therefore shown that
 {\def\termeA{\sum\limits_{i, j\in I_m}f(t_{i,m})f(t_{j,m})g^{\prime}\!\left({X_{t_{i,m}}}/{\sqrt{r(0)}}\right)g^{\prime}\!\left({X_{t_{j,m}}}/{\sqrt{r(0)}}\right)}
 \def\termeB{\int_{T_{i,m}\times T_{j, m}}\bC _{\Psi} \sum_{k=1}^d  \frac{\partial ^2  r}{\partial t_{k}^2}(t-s)  \ud t \ud s}
 \begin{align*}
   \MoveEqLeft[2]{\bbE\!\left[e^{is \prodsca{F_m(X)}{\xi_\ve}_{\calH }}\right] \cvg[\ve\to0]{}\int_{\reels^{i_m}} e^{-\frac12 s^2V_{11;m}} \ud \bbP_{(X(t_{i, m}), i \in I_m)}(x_{i, m}, i \in I_m)}\\
	& = \bbE\left[\exp\!\left(-\tfrac12 s^2 \termeA \right.\right.\\
	& \specialpos{\hfill\left.\rule{0pt}{24pt}\left.\times\termeB\right)\right]}\\
  	& =\bbE\!\left[ e^{-\frac12 s^2\prodsca{{\calK}_{\xi} F_m(X)}{F_m(X)}_{\calH }}\right]\\
	& =\int_{{\calH }}
   \bbE\!\left[e^{-\frac12 s^2\prodsca{{\calK}_{\xi} F_m(X)}{F_m(X)}_{\calH }}\given{X=x}\right] \ud {\calL}_{X}(x)\\
	& = \bbE\left[e^{is\prodsca{F_m(X)}{\xi}_{\calH }}\right].
\end{align*}}%
To summarize, we have shown that $Z_{\ve, m}(f) \cvg[\ve\to0]{$\calS_v$} Z_m(f)=\prodsca{F_m(X)}{\xi}_{\calH }$, 
such that if $x\in \bL^2(T)$
$$
	{\calL}(Z_m(f)\given{X=x})={\calN}(0; \prodsca{{\calK}_{\xi} F_m(X)}{F_m(X)}_{\calH }).
$$%
In a second step, we show that
 $$
 	Z_{m}(f) \cvg[m\to\infty]{${\calS_{v}}$} Z(f)=\prodsca{f(\bm{\cdot}) g^{\prime}\!\left({X(\bm{\cdot})}/{\sqrt{r(0)}}\right)}{\xi}_{\calH }.
$$%
Since $X$ has continuous trajectories, our study can be restricted to such a type of trajectories.
We know that%
$$
	F_m(x)\cvg[m\to\infty]{$L^2(T)$} F(x):=f(\bm{\cdot}) g^{\prime}\!\left({x}/{\sqrt{r(0)}}\right)
$$%
and since $K_\xi$ is a continuous operator, we obtain
$$
\prodsca{{\calK}_{\xi} F_m(X)}{F_m(X)}_{\calH }\to\prodsca{{\calK}_{\xi} F(x)}{F(x)}_{\calH }.$$%
Then
\begin{align*}
\bbE\!\left[e^{isZ_m(f)}\right]
	&\phantom{:}=\int_{{\calH }} \bbE\!\left[e^{isZ_m(f)}\given{X=x}\right] \ud {\calL}_{X}(x)\\\
	&\phantom{:}=\int_{{\calH }} e^{-\frac12{s^2} \prodsca{{\calK}_{\xi} F_m(X)}{F_m(X)}_{\calH }}
 \ud {\calL}_{X}(x)\\
\cvg[m\to\infty]{}&\quad \int_{{\calH }} e^{-\frac12{s^2} \prodsca{{\calK}_{\xi} F(x)}{F(x)}_{\calH }}
 \ud {\calL}_{X}(x)\\
	&:= \bbE[e^{isZ(f)}].
\end{align*}
In conclusion we have shown, 
$$
Z_{m}(f) \cvg[m\to\infty]{${\calS_{v}}$} Z(f)=\prodsca{f(\bm{\cdot}) g^{\prime}\!\left({X(\bm{\cdot})}/{\sqrt{r(0)}}\right)}{\xi}_{\calH },
$$%
such that if $x\in \bL^2(T)$
$$
{\calL}(Z(f)\given{X=x})={\calN}\!\left(0; \prodsca{{\calK}_{\xi} f(\bm{\cdot}) g^{\prime}\!\left({x}/{\sqrt{r(0)}}\right)}{f(\bm{\cdot}) g^{\prime}\!\left({x}/{\sqrt{r(0)}}\right)}_{\calH }\right).
$$%
To complete the proof, we apply Lemma 1.1 of \cite{MR920254}.
It is then sufficient to prove that
$$
	\lim_{m\to+\infty}\lim_{\ve\to0} \bbE\!\left[\left\{Z_{\ve,m}(f)-Z_\ve(f)\right\}^2\right]=0.
$$%
We compute $\bbE\!\left[\left\{Z_{\ve,m}(f)-Z_\ve(f)\right\}^2\right]$.
Applying the Mehler's formula (\cite{MR716933}), the fact that $\sum_{n=0}^{\infty} n^2  n! c_n^2 < \infty$ and Lemma \ref{majoration pour r(epsilon)},  we easily get
\begin{align*}
\MoveEqLeft[0]{\bbE\!\left[\left\{Z_{\ve,m}(f)-Z_\ve(f)\right\}^2\right]}\\
&=\sum_{i, j \in I_m} \int_{T_{i, m} \times T_{j, m}}
\bbE\!\left[{
\!\left({
f(t_{i,m}) g^{\prime}\!\left({X(t_{i,m})}/{\sqrt{r(0)}}\right)-f(t)g^{\prime}\!\left({{X(t)}/{\sqrt{r(0)}}}\right)
}\right) }\right.\\
&\specialpos{\hfill\left.{
\times\!\left({
f(t_{j,m}) g^{\prime}\!\left({X(t_{j,m})}/{\sqrt{r(0)}}\right)-f(s)g^{\prime}\!\left({X(s)}/{\sqrt{r(0)}}\right)
}\right) \xi_\ve(t)\xi_\ve(s)\hfill}
\right]\ud t \ud s}\\
&=\sum_{i, j \in I_m} \int_{T_{i, m} \times T_{j, m}}\sum_{n=1}^{\infty} \sum_{k=1}^{\infty} c_n c_k n k \\
&\quad\times\bbE\!\left[{
\!\left({
f(t_{i,m}) H_{n-1}\!\left({X(t_{i,m})}/{\sqrt{r(0)}}\right)-f(t)H_{n-1}\!\left({{X(t)}/{\sqrt{r(0)}}}\right)
}\right) 
}\right.\\
&\specialpos{\hfill\times\left.{
\!\left({
f(t_{j,m}) H_{k-1}\!\left({X(t_{j,m})}/{\sqrt{r(0)}}\right)-f(s)H_{k-1}\!\left({X(s)}/{\sqrt{r(0)}}\right)
}\right) \xi_\ve(t)\xi_\ve(s)
}\right]\ud t \ud s}\\
& \cvg[\ve\to0]{} \sum_{n=1}^{\infty} c_n^2 n n! 
 \sum_{i, j \in I_m} \int_{T_{i, m} \times T_{j, m}} \!\left({\bC _{\Psi} \sum_{k=1}^d  \dfrac{\partial ^2  r}{\partial t_{k}^2}(t-s)}\right) \\
& \quad\times\left[{
f(t_{i,m})f(t_{j,m}) \!\left({\frac{r(t_{i,m}-t_{j, m})}{r(0)}
 }\right)^{n-1}-f(t)f(t_{j,m})
 \!\left({\frac{r(t-t_{j, m})}{r(0)}
 }\right)^{n-1}
 }\right.\\
& \specialpos{\hfill\left.{-f(s)f(t_{i,m})
 \!\left({\frac{r(s-t_{i, m})}{r(0)}
 }\right)^{n-1} +f(t)f(s)
 \!\left({\frac{r(t-s)}{r(0)}
 }\right)^{n-1}
 }\right] \ud t \ud s}
 \\
& \cvg[m \to \infty]{} \sum_{n=1}^{\infty} c_n^2 n n! 
\int_{T \times T} \!\left({\bC _{\Psi}  \sum_{k=1}^d  \dfrac{\partial ^2  r}{\partial t_{k}^2}(t-s)}\right)  f(t) f(s) \\
 &\specialpos{\hfill\times\left[{
 \!\left({\frac{r(t-s}{r(0)}
 }\right)^{n-1}-
 \!\left({\frac{r(t-s)}{r(0)}
 }\right)^{n-1}
 -
 \!\left({\frac{r(s-t)}{r(0)}
 }\right)^{n-1} +
 \!\left({\frac{r(t-s)}{r(0)}
 }\right)^{n-1}
 }\right] \ud t \ud s}\\
 &=0.
\end{align*}
This completes the proof of Theorem \ref{vitesse convergence temps local}--\ref{itm3:theo4.6.20} .\\
Let us tackle the proof of the last part of that theorem.\\
We suppose that  $\alpha = {1}/{2}$  and $d=2$.\\
By Remark \ref{remark theoreme}, we already know that the term $S_1$ will prevail when  $\alpha = {1}/{2}$  and $d=2$.\\
Moreover, Remark \ref{H2} has emphasized the fact that
 \begin{multline*}
T_1=\ve^{-1} [\ln(1/\ve)]^{-1/2} S_1 \\
\simeq T_\ve:=
\ve^{-1} [\ln(1/\ve)]^{-1/2} a(2, 0) \int_T f(t) g\!\left({X_\ve(t)}/{\sqrt{r_\ve(0)}}\right) \\
\times\left[{H_2\!\left(\dfrac{\partial X_\ve}{\partial t_1}(t)/\sqrt{\mu_\ve}\right)+
H_2\!\left(\dfrac{\partial X_\ve}{\partial t_2}(t)/\sqrt{\mu_\ve}\right)
}\right] \ud t.
\end{multline*}
As in \cite{MR3289986}, the proof will proceed in several steps.
Remember that $T$ is an open rectangle.
If $T$ has the following form,  $T= ]a, b[ \times ]c, d[$, let us define for $t:=(t_1, t_2) \in T$,
\begin{multline*}
S_\ve(t):=
 \ve^{-1} [\ln(1/\ve)]^{-1/2} \\
 \times \int_{a}^{t_1} \int_{c}^{t_2}  
\left[{H_2\!\left(\dfrac{\partial X_\ve}{\partial u_1}(u)/\sqrt{\mu_\ve}\right)+
H_2\!\left(\dfrac{\partial X_\ve}{\partial u_2}(u)/\sqrt{\mu_\ve}\right)
}\right] \ud u.
\end{multline*}
On the one hand, we prove in forthcoming Proposition \ref{convergence stable} that $(X_\ve, S_\ve)$ stably converges to $(X, \sigma \widehat{W})$.
\\
On the other hand we will consider a discrete version of $T_\ve$, defining for $m \in \naturels^{\star}$
\begin{multline*}
T_\ve^{(m)}:=a(2, 0) \sum_{i=0}^{m-1} \sum_{j=0}^{m-1} f(a_i^{(m)}, c_j^{(m)}) g\!\left({X_\ve(a_i^{(m)}, c_j^{(m)})}/{\sqrt{r_\ve(0)}}\right)\\
\times\ve^{-1} [\ln(1/\ve)]^{-1/2}
	\int_{a_i^{(m)}}^{a_{i+1}^{(m)}} \int_{c_j^{(m)}}^{c_{j+1}^{(m)}}  
		\left[{H_2\!\left(\dfrac{\partial X_\ve}{\partial t_1}(t)/\sqrt{\mu_\ve}\right)+
H_2\!\left(\dfrac{\partial X_\ve}{\partial t_2}(t)/\sqrt{\mu_\ve}\right)
}\right] \ud t,
\end{multline*}
where for $i=0,~\dots,~m$, $a_i^{(m)}:=a+ih_m$, $c_i^{(m)}:=c+ih^{\prime}_m$, with $h_m:=\frac{b-a}{m}$ and
$h^{\prime}_m:=\frac{d-c}{m}$.\\
The stable convergence of $(X_\ve, S_\ve)$ implies that
\begin{multline*}
T_\ve^{(m)} \cvg[\ve\to0]{} T^{(m)}:=a(2,0) \sigma  \sum_{i=0}^{m-1} \sum_{j=0}^{m-1} f\!\left(a_i^{(m)}, c_j^{(m)}\right) g\!\left({X(a_i^{(m)}, c_j^{(m)})}/{\sqrt{r(0)}}\right)\\
\times\left[{\widehat{W}(a_{i+1}^{(m)}, c_{j+1}^{(m)})-\widehat{W}(a_{i+1}^{(m)}, c_{j}^{(m)})-\widehat{W}(a_{i}^{(m)}, c_{j+1}^{(m)})+\widehat{W}(a_{i}^{(m)}, c_{j}^{(m)})}\right].
\end{multline*}
Furthermore, it is easy to show that $T^{(m)} $ is a Cauchy sequence in $L^{2}(\Omega)$.
\\
It follows, using the asymptotic independence between $X$ and $ \widehat{W}$, that
$$
	T^{(m)} \cvg[m\to\infty]{$L^{2}(\Omega)$} a(2,0) \sigma \int_{T} f(t) g\!\left({{X(t)}/{\sqrt{r(0)}}}\right) \ud \widehat{W}(t).
$$%
To conclude, i.e. to prove the convergence of $T_\ve$, it suffices to prove Lemma \ref{Dynkin} for which a proof is given after that of Proposition \ref{stable convergence}.
\spacebefore
\begin{lemm}
\label{Dynkin}
$$
	\lim_{m \to +\infty} \lim_{\ve \to 0} \bbE[(T_\ve-T_\ve^{(m)})^2]=0.
$$%
\end{lemm}
\spacebefore
\begin{prop}
\label{convergence stable}
\begin{enumerate}
\item\label{itm1:convergence stable}
$S_\ve \cvg[\ve\to0]{Stable} \sigma \widehat{W}$.
\item\label{itm2:convergence stable} Furthermore,$(X_\ve,S_\ve) \cvg[\ve\to0]{Stable} (X, \widehat{W})$,\\
\end{enumerate}
where  $\widehat{W}$ is a standard Brownian sheet in $\reels^2$ independent of $X$.\\
\end{prop}
\spacebefore
\begin{proofarg}{Proof of proposition \ref{convergence stable}}
\\\ref{itm1:convergence stable}.
For $m \in \naturels^{\star}$, $a=t_0 < t_1< \dots <t_m=b$ and $c=s_0 < s_1<\dots <s_m=d$, let
\begin{align*}
\MoveEqLeft[2]{S_\ve^{(m)}({\bt }, \bs):= \sum_{i=0}^{m-1}\sum_{j=0}^{m-1} \alpha_{i, j}
	\left[{S_\ve(t_{i+1},s_{j+1})-S_\ve(t_{i+1},s_{j})}\right.} \\
	&\specialpos{\hfill-\left.{S_\ve(t_{i},s_{j+1})+S_\ve(t_{i},s_{j})}\right]}\\
	&= \sum_{i=0}^{m-1}\sum_{j=0}^{m-1} \alpha_{i, j} \ve^{-1} [\ln(1/\ve)]^{-1/2} \\
	&\quad\times\int_{t_i}^{t_{i+1}} \int_{s_j}^{s_{j+1}} 
\left[{H_2\!\left(\dfrac{\partial X_\ve}{\partial u_1}(u)/\sqrt{\mu_\ve}\right)
	+ H_2\!\left(\dfrac{\partial X_\ve}{\partial u_2}(u)/\sqrt{\mu_\ve}\right)}\right] \ud u,
\end{align*}
where
\begin{eqnarray*}
\alpha_{i, j}:= {d_{i, j}}\left/{\sqrt{\sum\limits_{i=0}^{m-1}\sum\limits_{j=0}^{m-1}d_{i, j}^2 (t_{i+1}-t_i)(s_{j+1}-s_j)}}\right.,
\end{eqnarray*}
while $d_{i, j} \in \reels$, for $i, j=0,\,\dots,\,m-1$.\\
We want to prove that $S_\ve^{(m)}({\bt }, \bs)$ converges in law to ${\calN}\!\left({0; \sigma^2}\right)$.\\
For this purpose and as in the proof of Theorem \ref{vitesse convergence temps local}--\ref{itm2:theo4.6.20}, we have the following representation for $S_\ve^{(m)}({\bt }, \bs)$,
\begin{multline*}
S_\ve^{(m)}({\bt }, \bs):= \sum_{i=0}^{m-1}\sum_{j=0}^{m-1} \alpha_{i, j} \ve^{-1} [\ln(1/\ve)]^{-1/2} \\
 \quad\times\int_{t_i}^{t_{i+1}} \int_{s_j}^{s_{j+1}} 
\left[{H_2\!\left({\frac{1}{\sqrt{\mu_\ve}} \int_{\reels^d} i \lambda^{(1)} e^{i \prodsca{\lambda}{u} } \ud Z_\ve(\lambda)
}\right)}\right.\\
\left.{+
H_2\!\left({\frac{1}{\sqrt{\mu_\ve}} \int_{\reels^d} i \lambda^{(2)} e^{i \prodsca{\lambda}{u} } \ud Z_\ve(\lambda)
}\right)}\right] \ud u,
\end{multline*}
where we noted $\lambda:=(\lambda^{(1)}, \lambda^{(2)})$ and $\ud Z_\ve(\lambda)=\widehat{\Psi}(\ve \lambda) \sqrt{s(\lambda)} \ud W(\lambda)$, while $W$ is a Brownian sheet in $\reels^2$.
Recall that $s$ is the spectral density of process $X$.\\
We use notations introduced in Slud \cite{MR1303648}.\\
For $h \in L^2_{\mbox{sym}}((\reels^2)^2)$, we define
$$
	I_2(h):=\frac{1}{2!} \int_{\reels^2}\int_{\reels^2} f(\lambda) f(\mu) \ud W(\lambda) \ud W(\mu).
$$%
By using the It\^o's formula for the Wiener-It\^o integral (see \cite[Theorem 4.3 p$.$~37]{MR3155040}), we obtain
$$
S_\ve^{(m)}({\bt }, \bs)=I_2(h_\ve^{(m)}),
$$%
where $h_\ve^{(m)}$ is
\begin{multline*}
h_\ve^{(m)}(\lambda, \mu):=-2! \sum_{i=0}^{m-1}\sum_{j=0}^{m-1} \alpha_{i, j} \ve^{-1} [\ln(1/\ve)]^{-1/2}\\
\times  \int_{t_i}^{t_{i+1}} \int_{s_j}^{s_{j+1}} 
\left[{\lambda^{(1)}\mu^{(1)}+\lambda^{(2)}\mu^{(2)}
}\right]
\widehat{\Psi}(\ve \lambda)  \widehat{\Psi}(\ve \mu) \sqrt{s(\lambda)} \sqrt{s(\mu)} 
 e^{i \prodsca{\lambda +\mu}{u}}\ud u.
\end{multline*}
To obtain convergence of $S_\ve^{(m)}({\bt }, \bs)$, we use \cite[Theorem 1]{MR2118863}.
Convergences given in (\ref{convergence 1}) and in (\ref{convergence 2}) in the proof of the forthcoming lemma give the required convergence appearing in part 1 of the later Theorem 1.
So we just verify condition (i) in proving the following lemma.
\spacebefore
\begin{lemm}
\label{moment d'ordre quatre}
$$\lim_{\ve \to 0} \bbE\left[\left\{\smash{S_\ve^{(m)}({\bt }, \bs)}\right\}^4\right]=3 \sigma^4.$$
\end{lemm}
\spacebefore
\begin{proofarg}{Proof of Lemma \ref{moment d'ordre quatre}}
Let us give some convergence results that will help us to show the lemma.\\
For all $i=0,\,\dots,\,m-1$ and $j=0,\,\dots,\,m-1$,
\begin{align}
\label{convergence 1}
\MoveEqLeft[2]{2 \ve^{-2} [\ln(1/\ve)]^{-1} 
 \int_{t_i}^{t_{i+1}} \int_{s_j}^{s_{j+1}}  \int_{t_i}^{t_{i+1}} \int_{s_j}^{s_{j+1}} 
\left[{
\!\left({\dfrac{\partial ^{2}r_\ve}{\partial t_1^2}(u-v)/\mu_\ve
}\right)^{2}  }\right.}\nonumber\\
& \specialpos{\hfill\left.+{ 2 \!\left({\dfrac{\partial ^{2}r_\ve}{\partial t_1\partial t_2}(u-v)/\mu_\ve
}\right)^{2} +\!\left({\dfrac{\partial ^{2}r_\ve}{\partial t_2^2}(u-v)/\mu_\ve}\right)^{2}
}\right] \ud u \ud v }\nonumber\\
&\cvg[\ve\to0]{} \sigma^2 (t_{i+1}-t_i) (s_{j+1}-s_j).
\end{align}
This is an immediate application of Theorem \ref{vitesse variance temps local}--\ref{itm4.6.19.4}, when considering functions $f$ and $g$ identically equal to one.\\
We also use the following convergence result.\\
For all indices, $i_1, i_2, j_1, j_2 =0,\,\dots,\,m-1$ such that $(i_1 \neq i_2 \mbox{ and } i_1+1 \neq i_2 \mbox{ and } i_2+1 \neq i_1)$ or such that $(j_1 \neq j_2 \mbox{ and } j_1+1 \neq j_2 \mbox{ and } j_2+1 \neq j_1)$ or such that $((i_2=i_1+1 \mbox{  or  } i_2+1=i_1) \mbox{ and } (j_2=j_1+1 \mbox{ or } j_2+1=j_1))$ or such that $(i_1=i_2 \mbox{  and  } (j_1+1=j_2 \mbox{ or } j_2+1=j_1))$, or such that $((i_1+1=i_2 \mbox{  or  } i_2+1=i_1)\mbox{ and }  j_1=j_2 )$, we have
\begin{multline}
\label{convergence 2}
 \ve^{-2} [\ln(1/\ve)]^{-1} 
 \int_{t_{i_1}}^{t_{i_1+1}} \int_{s_{j_1}}^{s_{j_1+1}}  \int_{t_{i_2}}^{t_{i_2+1}} \int_{s_{j_2}}^{s_{j_2+1}} 
\left[{
\!\left({\dfrac{\partial ^{2}r_\ve}{\partial t_1^2}(u-v)/\mu_\ve
}\right)^{2}  
 + }\right.
 \\
 \left.{ 2 \!\left({\dfrac{\partial ^{2}r_\ve}{\partial t_1\partial t_2}(u-v)/\mu_\ve
}\right)^{2} +\!\left({\dfrac{\partial ^{2}r_\ve}{\partial t_2^2}(u-v)/\mu_\ve}\right)^{2}
}\right] \ud u \ud v \cvg[\ve\to0]{}  0.
\end{multline}
Let us give a sketch of the proof of convergence given in (\ref{convergence 2}).\\
The convergence to zero of the above integral is straightforward in case where $(i_1 \neq i_2 \mbox{ and } i_1+1 \neq i_2 \mbox{ and } i_2+1 \neq i_1)$ or such that $(j_1 \neq j_2 \mbox{ and } j_1+1 \neq j_2 \mbox{ and } j_2+1 \neq j_1)$.
Indeed, in such a case, the integration variables $u$ and $v$ appearing in the integral are such that $\normp[2]{u-v} \ge \delta $ and as indicated in the proof of
Theorem \ref{vitesse variance temps local}--\ref{itm4.6.19.4}, the integral converges to zero.
Thus let us focus on the indices such that $((i_2=i_1+1 \mbox{ or } i_2+1=i_1) \mbox{ and } (j_2=j_1+1 \mbox{ or } j_2+1=j_1))$.\\
To this purpose, let us fix $t_1<t_2<t_3$ and $s_1<s_2<s_3$ and consider
\begin{multline*}
A_\ve:= 
 \ve^{-2} [\ln(1/\ve)]^{-1} 
 \int_{t_1}^{t_2} \int_{s_1}^{s_2}  \int_{t_2}^{t_{3}} \int_{s_2}^{s_{3}} 
\left[{
\!\left({\dfrac{\partial ^{2}r_\ve}{\partial t_1^2}(u-v)/\mu_\ve
}\right)^{2}  
 + }\right.
\nonumber \\
 \left.{ 2 \!\left({\dfrac{\partial ^{2}r_\ve}{\partial t_1\partial t_2}(u-v)/\mu_\ve
}\right)^{2} +\!\left({\dfrac{\partial ^{2}r_\ve}{\partial t_2^2}(u-v)/\mu_\ve}\right)^{2}
}\right] \ud u \ud v.
\end{multline*}
For the previous result obtained in Theorem \ref{vitesse variance temps local}--\ref{itm4.6.19.4}, we know that
\begin{align*}
A_\ve &\equiv  [\ln(1/\ve)]^{-1} \int_{t_2-\delta}^{t_2}\int_{s_2-\delta}^{s_2}\int_{t_2}^{t_2+\delta}\int_{s_2}^{s_2+\delta}
\frac{1}{\normp[2]{u-v}^2} \1_{\{ \normp[2]{u-v} \ge M \ve\}}\ud u \ud v\\
	&= \delta^2 [\ln(1/\ve)]^{-1} \int_{-1}^{0}\int_{-1}^{0}\int_{0}^{1}\int_{0}^{1}
\frac{1}{\normp[2]{u-v}^2} \1_{\{ \normp[2]{u-v} \ge M \ve/\delta\}}\ud u \ud v\\
	&\le \bC   \delta^2 [\ln(1/\ve)]^{-1} \left[{\ln(\sqrt{2})-\ln(M\ve/(2 \delta))}\right],
\end{align*}
thus $\lim_{\delta \to 0} \underset{\ve \to 0}{\overline{\lim}}  A_\ve=0$, and $\lim_{\ve \to 0}  A_\ve=0$.\\
Also, the case where $t_1<t_2<t_3$ and $s_1<s_2$ and 
\begin{multline*}
B_\ve:= 
 \ve^{-2} [\ln(1/\ve)]^{-1} 
 \int_{t_1}^{t_2} \int_{s_1}^{s_2}  \int_{t_2}^{t_{3}} \int_{s_1}^{s_2} 
\left[{
\!\left({\dfrac{\partial ^{2}r_\ve}{\partial t_1^2}(u-v)/\mu_\ve
}\right)^{2} }\right.
 \\
 +\left.{ 2 \!\left({\dfrac{\partial ^{2}r_\ve}{\partial t_1\partial t_2}(u-v)/\mu_\ve
}\right)^{2} +\!\left({\dfrac{\partial ^{2}r_\ve}{\partial t_2^2}(u-v)/\mu_\ve}\right)^{2}
}\right] \ud u \ud v,
\end{multline*}
is treated in a same way and gives $\lim_{\ve \to 0}  B_\ve=0$.\\
Finally, all the other cases are treated in the same manner by using the fact that for all $i, j =1, 2$, and all $(u_1, u_2)$, 
\begin{equation}
\label{property r}
\!\left({\dfrac{\partial ^{2}r_\ve}{\partial t_{i}\partial t_{j}}(u_1, u_2)}\right)^2= \!\left({\dfrac{\partial ^{2}r_\ve}{\partial t_{i}\partial t_{j}}(-u_1, u_2)}\right)^2=\!\left({\dfrac{\partial ^{2}r_\ve}{\partial t_{i}\partial t_{j}}(u_1, -u_2)}\right)^2.
\end{equation}
To finish with the preliminaries, we will use the next convergence result.\\
For all indices, $i_1, i_2,i_3, i_4, j_1, j_2, j_3, j_4 =0,\,\dots,\,m-1$ and for all $i=1, \dots, 4$ and $k_i=1, 2$
\begin{multline}
\label{convergence 3}
\ve^{-4} [\ln(1/\ve)]^{-2} 
 \int_{t_{i_1}}^{t_{i_1+1}} \int_{s_{j_1}}^{s_{j_1+1}}  \int_{t_{i_2}}^{t_{i_2+1}} \int_{s_{j_2}}^{s_{j_2+1}} 
  \int_{t_{i_3}}^{t_{i_{3}+1}} \int_{s_{j_3}}^{s_{j_{3}+1}}  \int_{t_{i_4}}^{t_{i_{4}+1}} \int_{s_{j_4}}^{s_{j_4+1}}  \\
\times\!\left({\dfrac{\partial ^{2}r_\ve}{\partial t_{k_1}\partial t_{k_2}}(u_1-u_2)/\mu_\ve
}\right)  
  \!\left({\dfrac{\partial ^{2}r_\ve}{\partial t_{k_1}\partial t_{k_3}}(u_1-u_3)/\mu_\ve
}\right) \!\left({\dfrac{\partial ^{2}r_\ve}{\partial t_{k_2} \partial t_{k_4}}(u_2-u_4)/\mu_\ve}\right) \\
\times \!\left({\dfrac{\partial ^{2}r_\ve}{\partial t_{k_3} \partial t_{k_4}}(u_3-u_4)/\mu_\ve}\right)
 \ud u_1 \ud u_2 \ud u_3 \ud u_4  \cvg[\ve\to0]{} 0.
\end{multline}
To prove this convergence, we first notice that by applying Schwarz's inequality and using the results of convergence (\ref{convergence 2}), we easily obtain the required convergence to zero of the previous integral for most cases.
The only case left  to consider is the following.\\
Let us fix $t_1<t_2$, $s_1<s_2$ and $k_i$ fixed, for $i=1,\,\dots,\,4$.
Let us consider
\begin{multline*}
C_\ve:=\ve^{-4} [\ln(1/\ve)]^{-2} 
 \int_{t_1}^{t_2} \int_{s_1}^{s_2}  \int_{t_1}^{t_2} \int_{s_1}^{s_2} 
  \int_{t_1}^{t_2} \int_{s_1}^{s_2}  \int_{t_1}^{t_2} \int_{s_1}^{s_2}  \\
\!\left({\dfrac{\partial ^{2}r_\ve}{\partial t_{k_1}\partial t_{k_2}}(u_1-u_2)/\mu_\ve
}\right)  
  \!\left({\dfrac{\partial ^{2}r_\ve}{\partial t_{k_1}\partial t_{k_3}}(u_1-u_3)/\mu_\ve
}\right) \!\left({\dfrac{\partial ^{2}r_\ve}{\partial t_{k_2} \partial t_{k_4}}(u_2-u_4)/\mu_\ve}\right)\\
 \!\left({\dfrac{\partial ^{2}r_\ve}{\partial t_{k_3} \partial t_{k_4}}(u_3-u_4)/\mu_\ve}\right)
 \ud u_1 \ud u_2 \ud u_3 \ud u_4 .
\end{multline*}
As before, applying Schwarz's inequality and working as in the proof of
Theorem \ref{vitesse variance temps local}--\ref{itm4.6.19.4}, we point out that this integral is $O(1/\ln(1/\ve))$ if we integrate on the set $\{\normp[2]{u_1-u_2} \le M\ve \mbox{ or } \normp[2]{u_1-u_3} \le M\ve \mbox{ or } \normp[2]{u_2-u_4} \le M\ve \mbox{ or } \normp[2]{u_3-u_4} \le M\ve \mbox{ or }\normp[2]{u_1-u_2} > \delta \mbox{ or }\normp[2]{u_1-u_3} > \delta \mbox{ or }\\\normp[2]{u_2-u_4} > \delta  \mbox{ or } \normp[2]{u_3-u_4} > \delta \}$.\\
Thus, we only have to consider the behavior of the integral on the set
\begin{multline*}
 \{M\ve \le \normp[2]{u_1-u_2} \le \delta  \mbox{ and }
M\ve \le \normp[2]{u_1-u_3} \le \delta\\  \mbox{ and } M\ve \le \normp[2]{u_2-u_4} \le \delta  \mbox{ and } M\ve \le \normp[2]{u_3-u_4} \le \delta\}.
\end{multline*}
Once again, we apply the Schwarz's inequality and Theorem \ref{vitesse variance temps local}--\ref{itm4.6.19.4}.
We obtain $C_\ve \equiv D_\ve$, where 
\begin{multline*}
D_\ve \le \bC  \!\left({\ve^{-4} [\ln(1/\ve)]^{-2} \int_{T^4} \1_{\{M\ve \le \normp[2]{u_1-u_2} \le \delta \}}
}\right.\\
\times\left.{
 \!\left({\dfrac{\partial ^{2}r_\ve}{\partial t_{k_1}\partial t_{k_2}}(u_1-u_2)/\mu_\ve
}\right)^2 
\1_{\{M\ve \le \normp[2]{u_3-u_4} \le \delta \}}
 \!\left({\dfrac{\partial ^{2}r_\ve}{\partial t_{k_3}\partial t_{k_4}}(u_3-u_4)/\mu_\ve
}\right)^2 }\right.\\
\times\left.{\1_{\{M\ve \le \normp[2]{u_1-u_3} \le \delta \}} \ud u_1 \ud u_2 \ud u_3 \ud u_4
}\right).
\end{multline*}
We make the change of variables: $u_1-u_2=v_1$ and $u_3-u_4=v_2$, getting
\begin{eqnarray*}
&&D_\ve \le \bC  \!\left({[\ln(1/\ve)]^{-1}\int_{M \ve \le \normp[2]{v} \le \delta} \frac{1}{\normp[2]{v}^2} \ud v}\right)^2 \!\left({\int_{T \times T} \1_{\{\normp[2]{u_1-u_3} \le \delta\}} }\right).
\end{eqnarray*}
Consequently, we obtain $\underset{\ve \to 0}{\overline{\lim}}  D_\ve \le \bC  \delta^2$, and
$\lim_{\ve \to 0} D_\ve=0$.\\
The ground is now prepared for the proof of the lemma.\\
We compute the fourth moment of the \rv
$S_\ve^{(m)}({\bt }, \bs)$.
\begin{align*}
\MoveEqLeft[2]{\bbE\left[\left\{\smash{S_\ve^{(m)}({\bt }, \bs)}\right\}^4\right]=  \sum_{i_1, i_2, i_3, i_4=0}^{m-1}\sum_{j_1, j_2, j_3, j_4=0}^{m-1} \alpha_{i_1, j_1} \alpha_{i_2, j_2} \alpha_{i_3, j_3} \alpha_{i_4, j_4}} \\
	&\times \ve^{-4} [\ln(1/\ve)]^{-2} 
 \int_{t_{i_1}}^{t_{i_1+1}} \int_{s_{j_1}}^{s_{j_1+1}}  \int_{t_{i_2}}^{t_{i_2+1}} \int_{s_{j_2}}^{s_{j_2+1}} 
  \int_{t_{i_3}}^{t_{i_{3}+1}} \int_{s_{j_3}}^{s_{j_{3}+1}}  \int_{t_{i_4}}^{t_{i_{4}+1}} \int_{s_{j_4}}^{s_{j_4+1}} \\
  	& \sum_{k_1, k_2, k_3, k_4=1,2}
  \bbE\!\left[{ 
  H_2\!\left(\dfrac{\partial X_\ve}{\partial t_{k_1}}(u_1)/\sqrt{\mu_\ve}\right)
H_2\!\left(\dfrac{\partial X_\ve}{\partial t_{k_2}}(u_2)/\sqrt{\mu_\ve}\right) }\right.\\
	&\specialpos{\hfill\left.{\times
H_2\!\left(\dfrac{\partial X_\ve}{\partial t_{k_3}}(u_3)/\sqrt{\mu_\ve}\right)
H_2\!\left(\dfrac{\partial X_\ve}{\partial t_{k_4}}(u_4)/\sqrt{\mu_\ve}\right)
}\right] \ud u_1 \ud u_2  \ud u_3 \ud u_4.}
\end{align*}
To calculate the expectation, we use the diagram formula (see \cite[pp$.$ 431 and 432]{MR716933}):
\begin{equation}
\label{formule du diagramme}
\bbE\!\left[\prod_{\ell=1}^4 H_2\!\left(\dfrac{\partial X_\ve}{\partial t_{k_{\ell}}}(u_{\ell})/\sqrt{\mu_\ve}\right)\right]
	=\sum_{G \in \Gamma} I_G,
\end{equation}
where $G$ is an undirected graph with 8 edges and 4 levels (for definitions, see \cite[p$.$ 431]{MR716933}), $\Gamma:=\Gamma(2, 2, 2, 2)$ denotes the set of diagrams having theses properties, and
\begin{equation}
\label{formule du diagramme bis}
I_G:= \prod_{w \in G(V)} \!\left({-\dfrac{\partial ^{2}r_\ve}{\partial t_{k_{d_1(w)}}\partial t_{k_{d_2(w)}}}(u_{d_1(w)}-u_{d_2(w)})/\mu_\ve
}\right),
\end{equation}
where $G(V)$ denotes the set of edges of $G$; the edges $w$ are oriented, beginning in $d_1(w)$ and finishing in $d_2(w)$.\\
The diagrams are called {\it regular} (see \cite[p. 432]{MR716933}) if their levels can be matched in a such way  that no edges pass between levels in different pairs, otherwise they are called {\it irregular}.
Consider all the regular diagrams.
Their contribution, say $R_\ve$, is given by
\begin{align*}
\MoveEqLeft[1]{R_\ve:= \sum_{i_1, i_2, i_3, i_4=0}^{m-1}\sum_{j_1, j_2, j_3, j_4=0}^{m-1} \alpha_{i_1, j_1} \alpha_{i_2, j_2} \alpha_{i_3, j_3} \alpha_{i_4, j_4}} \\
&\times\ve^{-4} [\ln(1/\ve)]^{-2} 
 \int_{t_{i_1}}^{t_{i_1+1}} \int_{s_{j_1}}^{s_{j_1+1}}  \int_{t_{i_2}}^{t_{i_2+1}} \int_{s_{j_2}}^{s_{j_2+1}} 
  \int_{t_{i_3}}^{t_{i_{3}+1}} \int_{s_{j_3}}^{s_{j_{3}+1}}  \int_{t_{i_4}}^{t_{i_{4}+1}} \int_{s_{j_4}}^{s_{j_4+1}} \\
&\times   \sum_{k_1, k_2, k_3, k_4=1,2} \left\{{ 
 2 \!\left({\dfrac{\partial ^{2}r_\ve}{\partial t_{k_1}\partial t_{k_2}}(u_1-u_2)/\mu_\ve
}\right)^2  2 \!\left({\dfrac{\partial ^{2}r_\ve}{\partial t_{k_3}\partial t_{k_4}}(u_3-u_4)/\mu_\ve
}\right)^2}\right.\\
& \specialpos{\hfill\left.{ + 
2 \!\left({\dfrac{\partial ^{2}r_\ve}{\partial t_{k_1}\partial t_{k_3}}(u_1-u_3)/\mu_\ve
}\right)^2  2 \!\left({\dfrac{\partial ^{2}r_\ve}{\partial t_{k_2}\partial t_{k_4}}(u_2-u_4)/\mu_\ve
}\right)^2}\right.\hfill}\\
& \specialpos{\hfill \left.{ + 
2 \!\left({\dfrac{\partial ^{2}r_\ve}{\partial t_{k_1}\partial t_{k_4}}(u_1-u_4)/\mu_\ve
}\right)^2  2 \!\left({\dfrac{\partial ^{2}r_\ve}{\partial t_{k_2}\partial t_{k_3}}(u_2-u_3)/\mu_\ve
}\right)^2
  }\right\} \ud u_1 \ud u_2  \ud u_3 \ud u_4}\\
&   =3 \sum_{i_1, i_2, i_3, i_4=0}^{m-1}\sum_{j_1, j_2, j_3, j_4=0}^{m-1} \alpha_{i_1, j_1} \alpha_{i_2, j_2} \alpha_{i_3, j_3} \alpha_{i_4, j_4} \\
&\times \ve^{-4} [\ln(1/\ve)]^{-2} 
 \int_{t_{i_1}}^{t_{i_1+1}} \int_{s_{j_1}}^{s_{j_1+1}}  \int_{t_{i_2}}^{t_{i_2+1}} \int_{s_{j_2}}^{s_{j_2+1}} 
  \int_{t_{i_3}}^{t_{i_{3}+1}} \int_{s_{j_3}}^{s_{j_{3}+1}}  \int_{t_{i_4}}^{t_{i_{4}+1}} \int_{s_{j_4}}^{s_{j_4+1}} \\
&   \specialpos{\hfill\times\sum_{k_1, k_2, k_3, k_4=1,2} 
 2 \!\left({\dfrac{\partial ^{2}r_\ve}{\partial t_{k_1}\partial t_{k_2}}(u_1-u_2)/\mu_\ve
}\right)^2\hfill}  \\
&\specialpos{\hfill\times 2 \!\left({\dfrac{\partial ^{2}r_\ve}{\partial t_{k_3}\partial t_{k_4}}(u_3-u_4)/\mu_\ve
}\right)^2
 \ud u_1 \ud u_2  \ud u_3 \ud u_4}
 \\
&   =3 \!\left({\sum_{i_1, i_2=0}^{m-1}\sum_{j_1, j_2=0}^{m-1} \alpha_{i_1, j_1} \alpha_{i_2, j_2}  
\ve^{-2} [\ln(1/\ve)]^{-1} 
 \int_{t_{i_1}}^{t_{i_1+1}} \int_{s_{j_1}}^{s_{j_1+1}}  \int_{t_{i_2}}^{t_{i_2+1}} \int_{s_{j_2}}^{s_{j_2+1}} 
}\right.\\
&  \specialpos{\hfill\times \left.{
 2  \sum_{k_1, k_2=1,2} \!\left({\dfrac{\partial ^{2}r_\ve}{\partial t_{k_1}\partial t_{k_2}}(u_1-u_2)/\mu_\ve
}\right)^2 
 \ud u_1 \ud u_2}\right)^2}
  \\
&   =3 \!\left({\sum_{i_1, i_2=0}^{m-1}\sum_{j_1, j_2=0}^{m-1} \alpha_{i_1, j_1} \alpha_{i_2, j_2}  
2 \ve^{-2} [\ln(1/\ve)]^{-1} 
 \int_{t_{i_1}}^{t_{i_1+1}} \int_{s_{j_1}}^{s_{j_1+1}}  \int_{t_{i_2}}^{t_{i_2+1}} \int_{s_{j_2}}^{s_{j_2+1}} 
}\right.\\
 & \specialpos{\hfill\left.{
  \!\left({\!\left({\dfrac{\partial ^{2}r_\ve}{\partial t_1^2}(u-v)/\mu_\ve
}\right)^2+2
 \!\left({\dfrac{\partial ^{2}r_\ve}{\partial t_1\partial t_2}(u-v)/\mu_\ve
}\right)^2}\right.
}\right.\hfill}\\
& \specialpos{\hfill+\left.{ \left.{ \!\left({\dfrac{\partial ^{2}r_\ve}{\partial t_2^2}(u-v)/\mu_\ve
}\right)^2
}\right)
 \ud u \ud v}\right)^2.
} \end{align*}
Now, by applying convergence obtained in (\ref{convergence 1}) and (\ref{convergence 2}), we see that the limit will be held for indices such that $i_1=i_2$ and $i_2=j_2$, giving
$$
\lim_{\ve \to 0}R_\ve=
 3 \!\left({\sum_{i=0}^{m-1}\sum_{j=0}^{m-1} \alpha_{i, j}^2 \sigma^2 (t_{i+1}-t_i)(s_{j+1}-s_j)
}\right)^2
= 3 \sigma^4.
$$%
To finish with the proof of the lemma, we conclude by noting that the contribution of the irregular diagrams tends to zero with $\ve$.
Indeed, it is enough to apply the formula of the diagram given in (\ref{formule du diagramme}) and (\ref{formule du diagramme bis}) to those diagrams.
Then convergence given in (\ref{convergence 3}) yields the required result.
\end{proofarg}%
In conclusion, we have proved the convergence of finite dimensional distributions of $S_\ve$.\\
To obtain assertion \ref{itm1:convergence stable} of Proposition \ref{convergence stable}, \ie the convergence of process $S_\ve$, it is enough to prove the tightness of the sequence of this process.

Let us prove the following lemma.
\spacebefore
\begin{lemm}
\label{tension}
For any $(s, t) \in T \times T$ and $\ve >0$, we have
$$
\bbE\!\left[\left\{\smash{S_\ve(t) -S_\ve(s)}\right\}^6\right] \le \bC  \normp[2]{t-s}^3.
$$
\end{lemm}
\spacebefore
\begin{proofarg}{Proof of Lemma \ref{tension}}
Let $t:=(t_1, t_2) \in T$ and $s:=(s_1, s_2) \in T$.
Let us compute $
\bbE\!\left[\left\{\smash{S_\ve(t) -S_\ve(s)}\right\}^6\right]
$.
\begin{align*}
\MoveEqLeft[2]{\bbE\!\left[\left\{\smash{S_\ve(t) -S_\ve(s)}\right\}^6\right] =\ve^{-6} [\ln(1/\ve)]^{-3}}\\
	&  \bbE\!\left[\left(\int_{a}^{t_1} \int_{c}^{t_2}  
\left\{{H_2\left(\dfrac{\partial X_\ve}{\partial u_1}(u)/\sqrt{\mu_\ve}\right)+
H_2\!\left(\dfrac{\partial X_\ve}{\partial u_2}(u)/\sqrt{\mu_\ve}\right)
}\right\} \ud u\right.\right.\\
	&\left.\left.-\int_{a}^{s_1} \int_{c}^{s_2}  
\left\{{H_2\!\left(\dfrac{\partial X_\ve}{\partial u_1}(u)/\sqrt{\mu_\ve}\right)+
H_2\!\left(\dfrac{\partial X_\ve}{\partial u_2}(u)/\sqrt{\mu_\ve}\right)
}\right\} \ud u\right)^6\right]\\
	 \le {} & \bC  \ve^{-6} [\ln(1/\ve)]^{-3}\\
	& \bbE\!\left[\left(\int_{s_1}^{t_1} \int_{c}^{t_2}  
\left\{{H_2\!\left(\dfrac{\partial X_\ve}{\partial u_1}(u)/\sqrt{\mu_\ve}\right)+
H_2\!\left(\dfrac{\partial X_\ve}{\partial u_2}(u)/\sqrt{\mu_\ve}\right)
}\right\} \ud u\right)^6\right]\\
	& + \bbE\!\left[\left(\int_{a}^{s_1} \int_{s_2}^{t_2}  
\left\{{H_2\!\left(\dfrac{\partial X_\ve}{\partial u_1}(u)/\sqrt{\mu_\ve}\right)+
H_2\!\left(\dfrac{\partial X_\ve}{\partial u_2}(u)/\sqrt{\mu_\ve}\right)
}\right\} \ud u\right)^6\right].
\end{align*}
We will only work with the first expectation, say $V_\ve(t_1, s_1)$, the second would be treated in the same way.
\\
Let us suppose that $s_1 \le t_1$ and prove that 
\begin{align*}
V_\ve(t_1, s_1) &:=\ve^{-6} [\ln(1/\ve)]^{-3}
	\bbE\!\left[\rule{0pt}{18pt}\left(\int_{s_1}^{t_1} \int_{c}^{t_2}  
		\left\{H_2\!\left(\dfrac{\partial X_\ve}{\partial u_1}(u)/\sqrt{\mu_\ve}\right)\right.\right.\right.\\
		&\specialpos{\hfill+\left.\left.\left. H_2\!\left(\dfrac{\partial X_\ve}{\partial u_2}(u)/\sqrt{\mu_\ve}\right)
\right\} \ud u\right)^6\right]}\\
& \phantom{:}\le  \bC  (t_1-s_1)^3,
\end{align*}
that will be enough to ensure the required result.
\begin{align*}
\MoveEqLeft[2]{V_\ve(t_1,s_1)
=\ve^{-6} [\ln(1/\ve)]^{-3}\int_{s_1}^{t_1} \int_{c}^{t_2} \int_{s_1}^{t_1} \int_{c}^{t_2} \int_{s_1}^{t_1} \int_{c}^{t_2} \int_{s_1}^{t_1} \int_{c}^{t_2} \int_{s_1}^{t_1} \int_{c}^{t_2} \int_{s_1}^{t_1} \int_{c}^{t_2}}\\
&\specialpos{\hfill\sum_{k_1, k_2, k_3, k_4, k_5, k_6=1,2} \bbE\!\left[\prod_{\ell=1}^6 H_2\!\left(\dfrac{\partial X_\ve}{\partial t_{k_{\ell}}}(u_{\ell})/\sqrt{\mu_\ve}\right)\right] \ud u_1 \ud u_2 \ud u_3 \ud u_4 \ud u_5 \ud u_6.}
\end{align*}
As in the proof of Lemma \ref{moment d'ordre quatre}, we use the diagram formula (see \cite[pp$.$ 431 and 432]{MR716933}) to calculate the expectation:
$$
\bbE\!\left[\prod_{\ell=1}^6 H_2\!\left(\dfrac{\partial X_\ve}{\partial t_{k_{\ell}}}(u_{\ell})/\sqrt{\mu_\ve}\right)\right]
	=\sum_{G \in \Gamma} I_G,
$$%
where $G$ is an undirected graph with 12 edges and 6 levels (see \cite[p$.$ 431]{MR716933}, for definitions), $\Gamma:=\Gamma(2, 2, 2, 2, 2, 2)$ denotes the set of diagrams having theses properties, and
$$
I_G:= \prod_{w \in G(V)} \!\left({-\dfrac{\partial ^{2}r_\ve}{\partial t_{k_{d_1(w)}}\partial t_{k_{d_2(w)}}}(u_{d_1(w)}-u_{d_2(w)})/\mu_\ve
}\right),
$$%
where $G(V)$ denotes the set of edges of $G$; the edges $w$ are oriented, beginning in $d_1(w)$ and finishing in $d_2(w)$.\\
Thus we will have to deal with four kinds of terms corresponding to four kinds of graphs.
\begin{itemize}
      \item[Case 1]
This case corresponds to the regular diagrams.
There are fifteen such diagrams.
They are of the form
\begin{multline*}
\!\left({\dfrac{\partial ^{2}r_\ve}{\partial t_{k_1}\partial t_{k_2}}(u_1-u_2)/\mu_\ve
}\right)^2  \times
  \!\left({\dfrac{\partial ^{2}r_\ve}{\partial t_{k_3}\partial t_{k_4}}(u_3-u_4)/\mu_\ve
}\right)^2  \\
\times \!\left({\dfrac{\partial ^{2}r_\ve}{\partial t_{k_5} \partial t_{k_6}}(u_5-u_6)/\mu_\ve}\right)^2.
\end{multline*}
The three other cases correspond to the irregular ones.
\item[Case 2]
There are thirty diagrams with the form
\begin{multline*}
\!\left({\dfrac{\partial ^{2}r_\ve}{\partial t_{k_1}\partial t_{k_2}}(u_1-u_2)/\mu_\ve
}\right)^2  \times \!\left({\dfrac{\partial ^{2}r_\ve}{\partial t_{k_3}\partial t_{k_4}}(u_3-u_4)/\mu_\ve
}\right) \\
\times \!\left({\dfrac{\partial ^{2}r_\ve}{\partial t_{k_4} \partial t_{k_5}}(u_4-u_5)/\mu_\ve}\right) \times
\!\left({\dfrac{\partial ^{2}r_\ve}{\partial t_{k_5}\partial t_{k_6}}(u_5-u_6)/\mu_\ve
}\right) \\
\times \!\left({\dfrac{\partial ^{2}r_\ve}{\partial t_{k_3} \partial t_{k_6}}(u_3-u_6)/\mu_\ve}\right).
\end{multline*}
\item[Case 3]
There are forty diagrams with the following form
\begin{multline*}
\!\left({\dfrac{\partial ^{2}r_\ve}{\partial t_{k_1}\partial t_{k_2}}(u_1-u_2)/\mu_\ve
}\right) \times \!\left({\dfrac{\partial ^{2}r_\ve}{\partial t_{k_2}\partial t_{k_3}}(u_2-u_3)/\mu_\ve
}\right) \\
 \times \!\left({\dfrac{\partial ^{2}r_\ve}{\partial t_{k_1} \partial t_{k_3}}(u_1-u_3)/\mu_\ve}\right) \times
\!\left({\dfrac{\partial ^{2}r_\ve}{\partial t_{k_4}\partial t_{k_5}}(u_4-u_5)/\mu_\ve}\right) \\
 \times \!\left({\dfrac{\partial ^{2}r_\ve}{\partial t_{k_5} \partial t_{k_6}}(u_5-u_6)/\mu_\ve}\right) \times
\!\left({\dfrac{\partial ^{2}r_\ve}{\partial t_{k_4} \partial t_{k_6}}(u_4-u_6)/\mu_\ve
}\right).
\end{multline*}
\item[Case 4]
There are $5!$ diagrams of the form
\begin{multline*}
\!\left({\dfrac{\partial ^{2}r_\ve}{\partial t_{k_1}\partial t_{k_2}}(u_1-u_2)/\mu_\ve
}\right)  \times \!\left({\dfrac{\partial ^{2}r_\ve}{\partial t_{k_1}\partial t_{k_3}}(u_1-u_3)/\mu_\ve
}\right) \\
 \times \!\left({\dfrac{\partial ^{2}r_\ve}{\partial t_{k_2} \partial t_{k_4}}(u_2-u_4)/\mu_\ve}\right)\times
\!\left({\dfrac{\partial ^{2}r_\ve}{\partial t_{k_3}\partial t_{k_5}}(u_3-u_5)/\mu_\ve}\right) \\
 \times \!\left({\dfrac{\partial ^{2}r_\ve}{\partial t_{k_4} \partial t_{k_6}}(u_4-u_6)/\mu_\ve}\right) \times
\!\left({\dfrac{\partial ^{2}r_\ve}{\partial t_{k_5} \partial t_{k_6}}(u_5-u_6)/\mu_\ve
}\right).
\end{multline*}
\end{itemize}
As in the proof of Lemma \ref{moment d'ordre quatre}
let us give and show in the forthcoming lemma two preliminaries results which will help us to treat the previous four cases.
\spacebefore
\begin{lemm}$\quad$
\label{tension bis}
\begin{enumerate}
\item\label{itm:lem tension bis 1} $\ds \ve^{-2} [\ln(1/\ve)]^{-1}\int_{s_1}^{t_1} \int_{c}^{t_2} \int_{s_1}^{t_1} \int_{c}^{t_2} 
 \!\left({\dfrac{\partial ^{2}r_\ve}{\partial t_{k_1}\partial t_{k_2}}(u_1-u_2)/\mu_\ve
}\right)^2 \ud u_1 \ud u_2\\
$\\
$
\rule{24pt}{0pt} {\le  \bC  (t_1-s_1).}
$
\item\label{itm:lem tension bis 2} 
$\ds \ve^{-3} [\ln(1/\ve)]^{-3/2}\int_{s_1}^{t_1} \int_{c}^{t_2} \int_{s_1}^{t_1} \int_{c}^{t_2} \int_{s_1}^{t_1} \int_{c}^{t_2}
 \abs{\dfrac{\partial ^{2}r_\ve}{\partial t_{k_1}\partial t_{k_2}}(u_1-u_2)/\mu_\ve
}$
\begin{align*}
\MoveEqLeft[1]{\rule{24pt}{0pt}\times \abs{\dfrac{\partial ^{2}r_\ve}{\partial t_{k_2}\partial t_{k_3}}(u_2-u_3)/\mu_\ve
} \abs{\dfrac{\partial ^{2}r_\ve}{\partial t_{k_1}\partial t_{k_3}}(u_1-u_3)/\mu_\ve
}  \ud u_1 \ud u_2 \ud u_3}\\
 &\le  \bC  (t_1-s_1)^{3/2}.
\end{align*}
\end{enumerate}
\end{lemm}
\spacebefore
\begin{proofarg}{Proof of Lemma \ref{tension bis}}
Let us begin with the proof of part \ref{itm:lem tension bis 1} of the  lemma.
By using the parity property of the function $\dfrac{\partial ^{2}r_\ve}{\partial t_{k_1}\partial t_{k_2}}$ described in the formula (\ref{property r}), inequality (\ref{majoration derivees secondes r epsilon}) and  convergence given in (\ref{limite mu epsilon}), we get
\begin{align*}
\MoveEqLeft[2]{\ve^{-2} [\ln(1/\ve)]^{-1}\int_{s_1}^{t_1} \int_{c}^{t_2} \int_{s_1}^{t_1} \int_{c}^{t_2} 
 \!\left({\dfrac{\partial ^{2}r_\ve}{\partial t_{k_1}\partial t_{k_2}}(u_1-u_2)/\mu_\ve
}\right)^2 \ud u_1 \ud u_2}\\
	& = 4  \ve^{-2} [\ln(1/\ve)]^{-1}\\
	& \specialpos{\hfill\times\int_{0}^{t_1-s_1} \int_{0}^{t_2-c} 
 \!\left({\dfrac{\partial ^{2}r_\ve}{\partial t_{k_1}\partial t_{k_2}}(u)/\mu_\ve
}\right)^2 (t_1-s_1-u_1)(t_2-c-u_2)\ud u}\\
	& \le  \bC  (t_1-s_1) [\ln(1/\ve)]^{-1} \\
	& \specialpos{\hfill\times \int_{0}^{t_1-s_1} \int_{0}^{t_2-c} \left\{{\ve^{-2}\1_{\{\normp[2]{u} \le M\ve\}} +\frac{1}{\normp[2]{u}^2}  \1_{\{\normp[2]{u} > M\ve\}} }\right\} \ud u}\\
	&  \le  \bC  (t_1-s_1).
\end{align*}
Let us now turn to the part \ref{itm:lem tension bis 2} of the lemma and establish a first result.
\begin{align}
\label{inegalité préliminaire}
\MoveEqLeft[2]{\ve^{-1} [\ln(1/\ve)]^{-3/2}\int_{s_1}^{t_1} \int_{c}^{t_2} \int_{s_1}^{t_1} \int_{c}^{t_2} 
\abs{\dfrac{\partial ^{2}r_\ve}{\partial t_{k_1}\partial t_{k_2}}(u_1-u_2)/\mu_\ve} \ud u_1 \ud u_2} \\
 & \le  \bC  (t_1-s_1)^{3/2}\nonumber.
\end{align}
Indeed, we have
\begin{align*}
\MoveEqLeft[2]{\ve^{-1} [\ln(1/\ve)]^{-3/2}\int_{s_1}^{t_1} \int_{c}^{t_2} \int_{s_1}^{t_1} \int_{c}^{t_2} 
\abs{\dfrac{\partial ^{2}r_\ve}{\partial t_{k_1}\partial t_{k_2}}(u_1-u_2)/\mu_\ve} \ud u_1 \ud u_2}\\
	& =4 \ve^{-1} [\ln(1/\ve)]^{-3/2}\\
	& \specialpos{\hfill\times  \int_{0}^{t_1-s_1} \int_{0}^{t_2-c} 
\abs{\dfrac{\partial ^{2}r_\ve}{\partial t_{k_1}\partial t_{k_2}}(u)/\mu_\ve} (t_1-s_1-u_1)(t_2-c-u_2)\ud u}\\
	& \le  \bC  (t_1-s_1) [\ln(1/\ve)]^{-3/2}\\
	& \specialpos{\hfill\times\int_{0}^{t_1-s_1} \int_{0}^{t_2-c} \left\{{\ve^{-1}\1_{\{\normp[2]{u} \le M\ve\}} +\frac{1}{\normp[2]{u}}  \1_{\{\normp[2]{u} > M\ve\}} }\right\} \ud u}\\
	& \le  \bC  (t_1-s_1) [\ln(1/\ve)]^{-3/2}\\
	& \specialpos{\hfill\times\left\{{(t_1-s_1)+ (t_1-s_1)^{1/2}\!\left({ \int_{0}^{t_1-s_1} \int_{0}^{t_2-c} \frac{1}{\normp[2]{u}^2}  \1_{\{\normp[2]{u} > M\ve\}}  \ud u}\right)^{1/2}}\right\}}\\
	& \le  \bC  (t_1-s_1)^{3/2} [\ln(1/\ve)]^{-1}\\
	& \le  \bC  (t_1-s_1)^{3/2}.
\end{align*}
Back to the proof of lemma part \ref{itm:lem tension bis 2}, we get the following bound
\begin{align*}
\MoveEqLeft[2]{\ve^{-3} [\ln(1/\ve)]^{-3/2}}\\
	 & \int_{s_1}^{t_1} \int_{c}^{t_2} \int_{s_1}^{t_1} \int_{c}^{t_2} \int_{s_1}^{t_1} \int_{c}^{t_2}
\1_{\{\normp[2]{u_1-u_2}\le M\ve \}} \abs{\dfrac{\partial ^{2}r_\ve}{\partial t_{k_1}\partial t_{k_2}}(u_1-u_2)/\mu_\ve
}\\
	 &  \specialpos{\hfill\times\abs{\dfrac{\partial ^{2}r_\ve}{\partial t_{k_2}\partial t_{k_3}}(u_2-u_3)/\mu_\ve} 
\abs{\dfrac{\partial ^{2}r_\ve}{\partial t_{k_1}\partial t_{k_3}}(u_1-u_3)/\mu_\ve
}  \ud u_1 \ud u_2 \ud u_3}\\
	 & \le \ve^{-3} [\ln(1/\ve)]^{-3/2} \int_{s_1}^{t_1} \int_{c}^{t_2} \int_{s_1}^{t_1} \int_{c}^{t_2} \int_{s_1}^{t_1} \int_{c}^{t_2}
\1_{\{\normp[2]{u_1-u_2}\le M\ve \}} \\
	 &  \specialpos{\hfill\times\abs{\dfrac{\partial ^{2}r_\ve}{\partial t_{k_2}\partial t_{k_3}}(u_2-u_3)/\mu_\ve} 
	    \ud u_1 \ud u_2 \ud u_3}\\
	 &   \le \ve^{-1} [\ln(1/\ve)]^{-3/2}\\
	 &   \specialpos{\hfill\times\int_{s_1}^{t_1} \int_{c}^{t_2} \int_{s_1}^{t_1} \int_{c}^{t_2} 
  \abs{\dfrac{\partial ^{2}r_\ve}{\partial t_{k_2}\partial t_{k_3}}(u_2-u_3)/\mu_\ve} 
  \ud u_2 \ud u_3}\\
	 & \le  \bC  (t_1-s_1)^{3/2},
\end{align*}
the last inequality coming from (\ref{inegalité préliminaire}).\\
The same reasoning can be done for the set of integration $\normp[2]{u_2-u_3}\le M\ve$.\\
Thus to conclude the proof of inequality given in part b., we only have to consider the set $\normp[2]{u_1-u_2}> M\ve \mbox{ and } \normp[2]{u_2-u_3} > M\ve$.
We obtain, by applying Schwarz's inequality and result obtained in part \ref{itm:lem tension bis 1} of this lemma
\begin{align*}
\MoveEqLeft[0]{\ve^{-3} [\ln(1/\ve)]^{-3/2} }\\
	&\specialpos{\hfill \times \int_{s_1}^{t_1} \int_{c}^{t_2} \int_{s_1}^{t_1} \int_{c}^{t_2} \int_{s_1}^{t_1} \int_{c}^{t_2}
 \abs{\dfrac{\partial ^{2}r_\ve}{\partial t_{k_1}\partial t_{k_2}}(u_1-u_2)/\mu_\ve
} \1_{\{\normp[2]{u_1-u_2} > M\ve \}}\hfill}\\
	&  \specialpos{\hfill\times\abs{\dfrac{\partial ^{2}r_\ve}{\partial t_{k_2}\partial t_{k_3}}(u_2-u_3)/\mu_\ve
} \1_{\{\normp[2]{u_2-u_3}> M\ve \}}\abs{\dfrac{\partial ^{2}r_\ve}{\partial t_{k_1}\partial t_{k_3}}(u_1-u_3)/\mu_\ve
}  \ud u_1 \ud u_2 \ud u_3}\\
	& \le \!\left(\rule{0pt}{14pt}{\ve^{-2} [\ln(1/\ve)]^{-1} }\right.\\
	& \specialpos{\hfill\times\left.{ \int_{s_1}^{t_1} \int_{c}^{t_2} \int_{s_1}^{t_1} \int_{c}^{t_2} \int_{s_1}^{t_1} \int_{c}^{t_2}
 \!\left({\dfrac{\partial ^{2}r_\ve}{\partial t_{k_1}\partial t_{k_3}}(u_1-u_3)/\mu_\ve
}\right)^2  \ud u_1 \ud u_2 \ud u_3}\right)^{1/2}}\\
	& \rule{24pt}{0pt}\kern2pt\times\!\left[\rule{0pt}{14pt}{\ve^{-4} [\ln(1/\ve)]^{-2} }\right.\\
	& \specialpos{\hfill\times\left.{
\int_{s_1}^{t_1} \int_{c}^{t_2} \int_{s_1}^{t_1} \int_{c}^{t_2} \int_{s_1}^{t_1} \int_{c}^{t_2}
 \!\left({\dfrac{\partial ^{2}r_\ve}{\partial t_{k_1}\partial t_{k_2}}(u_1-u_2)/\mu_\ve
}\right)^2 \1_{\{\normp[2]{u_1-u_2} > M\ve \}}}\right.}\\
	& \specialpos{\hfill\times\left.{
 \!\left({\dfrac{\partial ^{2}r_\ve}{\partial t_{k_1}\partial t_{k_2}}(u_2-u_3)/\mu_\ve
}\right)^2 \1_{\{\normp[2]{u_2-u_3} > M\ve \}} \ud u_1 \ud u_2 \ud u_3
}\right]^{1/2}}\\
	&  \le  \bC  (t_1-s_1) \!\left(\rule{0pt}{14pt}{ [\ln(1/\ve)]^{-2} }\right.\\
	&\specialpos{\hfill \times\left.{
\int_{s_1}^{t_1} \int_{c}^{t_2} \int_{s_1}^{t_1} \int_{c}^{t_2} \int_{s_1}^{t_1} \int_{c}^{t_2}
\frac{1}{\normp[2]{u_1-u_2}^2} \1_{\{\normp[2]{u_1-u_2} > M\ve \}}}\right.\hfill}\\
	&\specialpos{\hfill \times\left.{
 \frac{1}{\normp[2]{u_2-u_3}^2} \1_{\{\normp[2]{u_2-u_3} > M\ve \}} \ud u_1 \ud u_2 \ud u_3
}\right)^{1/2}}\\
	& \le  \bC  (t_1-s_1)^{3/2}.
\end{align*}
That ends the proof of Lemma \ref{tension bis}.
\end{proofarg}%
We turn back now to the proof of Lemma \ref{tension}.\\
The first case is treated directly by using part \ref{itm:lem tension bis 1} of Lemma \ref{tension bis}.\\
The second cases uses first part \ref{itm:lem tension bis 1} of Lemma \ref{tension bis} and second Schwarz's inequality and again part \ref{itm:lem tension bis 1} of Lemma \ref{tension bis}.\\
The third case is a straightforward application of part \ref{itm:lem tension bis 2} of Lemma \ref{tension bis}.\\
Finally, the last case is solved by applying Schwarz's inequality and again part \ref{itm:lem tension bis 1} of Lemma \ref{tension bis}.\\
That yields Lemma \ref{tension}.
\end{proofarg}%
The tightness of process $S_\ve$ follows from \cite[Theorem 3.1.10]{MR3088856}.\\
Thus the convergence in assertion \ref{itm1:convergence stable} of Proposition \ref{convergence stable} takes place in the sense of processes convergence.
Let us now show assertion \ref{itm2:convergence stable}.\\
Let $c_{i, j}, i, j =0,\,\dots,\,m-1$, be real constants.
As before, it is sufficient to establish the limit distribution of 
$$
\sum_{i=0}^{m-1} \sum_{j=0}^{m-1} c_{i, j} X_\ve(t_i, s_j)+  S_\ve^{(m)}({\bt }, \bs).
$$%
This time, we apply of \cite[Theorem 1]{MR2126978}.
Indeed, 
$$
\sum_{i=0}^{m-1} \sum_{j=0}^{m-1} c_{i, j} X_\ve(t_i, s_j)
$$%
belongs to the first Wiener chaos as Gaussian random variable with finite variance.
Furthermore, $S_\ve^{(m)}({\bt }, \bs)$ belongs to the two order one.
Mehler's formula allows us to verify the hypothesis appearing in (4) of this theorem.\\
Finally, assumption (iv) is satisfied since, by applying result of the previous part \ref{itm1:convergence stable}, we know that $S_\ve^{(m)}({\bt }, \bs)$ converges in distribution to a Gaussian random variable,  as well as $
\sum_{i=0}^{m-1} \sum_{j=0}^{m-1} c_{i, j} X_\ve(t_i, s_j)$.\\
To summarize, we have proved the convergence of finite dimensional distributions of $(X_\ve, S_\ve)$.\\
To conclude, the tightness of processes $(X_\ve, S_\ve)$ follows from that of process $S_\ve$ proved in part  \ref{itm1:convergence stable}.
It implies the convergence of $(X_\ve, S_\ve)$ as processes and then ensures the stable convergence of $(X_\ve, S_\ve)$ and of that of $S_\ve$ (see \cite[ Proposition 1]{MR2683464}).\\
This yields proof of Proposition \ref{convergence stable}.
\end{proofarg}%
To conclude with the proof of Theorem \ref{vitesse convergence temps local}, let us prove the last lemma.
\spacebefore
\begin{proofarg}{Proof of Lemma \ref{Dynkin}}
We can write $\bbE[(T_\ve-T_\ve^{(m)})^2]= \bbE[T_\ve^2]+\bbE[(T_\ve^{(m)})^2]-2 \bbE[T_\ve T_\ve^{(m)}]$.
By Theorem \ref{vitesse variance temps local}--\ref{itm4.6.19.4},
$$
\bbE[T_\ve^2] \cvg[\ve\to0]{} \ell_5.
$$%
Now, we compute $\bbE\!\left[\left(\smash{T_\ve^{(m)}}\right)^2\right]$.
\begin{align*}
\MoveEqLeft[0]{\bbE\!\left[\left(\smash{T_\ve^{(m)}}\right)^2\right]=a^2(2, 0) \ve^{-2} [\ln(1/\ve)]^{-1}  \sum_{i_1, i_2=0}^{m-1}\sum_{j_1, j_2=0}^{m-1}  
f(a_{i_1}^{(m)}, c_{j_1}^{(m)}) f\!\left(a_{i_2}^{(m)}, c_{j_2}^{(m)}\right)} \\
	& \times\sum_{k=0}^{\infty} \sum_{\ell=0}^{\infty} c_k c_{\ell}  \int_{a_{i_1}^{(m)}}^{a_{i_1+1}^{(m)}} \int_{c_{j_1}^{(m)}}^{c_{j_1+1}^{(m)}}  \int_{a_{i_2}^{(m)}}^{a_{i_2+1}^{(m)}} \int_{c_{j_2}^{(m)}}^{c_{j_2+1}^{(m)}} \\
  	& \times\bbE\!\left[{ 
 \!\left({ H_2\!\left(\dfrac{\partial X_\ve}{\partial t_1}(t)/\sqrt{\mu_\ve}\right)+
H_2\!\left(\dfrac{\partial X_\ve}{\partial t_2}(t)/\sqrt{\mu_\ve}\right) }\right) H_{k}\!\left(\frac{X_\ve(a_{i_1}^{(m)}, c_{j_1}^{(m)})}{\sqrt{r_\ve(0)}}\right)}\right.\\
	& \times\left.{
\!\left({H_2\!\left(\dfrac{\partial X_\ve}{\partial s_1}(s)/\sqrt{\mu_\ve}\right)
+H_2\!\left(\dfrac{\partial X_\ve}{\partial s_2}(s)/\sqrt{\mu_\ve}\right)}\right)H_{\ell}\!\left(\frac{X_\ve(a_{i_2}^{(m)}, c_{j_2}^{(m)})}{\sqrt{r_\ve(0)}}\right)
}\right] \ud t \ud s.
\end{align*}
We proceed as in the proof of Proposition \ref{term S1}--\ref{itm:termS1.3}.
It is straightforward, by using the asymptotic equality given in (\ref{Dynkin Bis}), to prove that
\begin{align*}
\MoveEqLeft[1]{\bbE\!\left[\left(\smash{T_\ve^{(m)}}\right)^2\right]\equiv 2 a^2(2, 0) \ve^{-2} [\ln(1/\ve)]^{-1}  \sum_{i_1, i_2=0}^{m-1}\sum_{j_1, j_2=0}^{m-1}  
f(a_{i_1}^{(m)}, c_{j_1}^{(m)}) f(a_{i_2}^{(m)}, c_{j_2}^{(m)}}) \\
	& \times\sum_{q=0}^{\infty} q! c_q^2  \int_{a_{i_1}^{(m)}}^{a_{i_1+1}^{(m)}} \int_{c_{j_1}^{(m)}}^{c_{j_1+1}^{(m)}}  \int_{a_{i_2}^{(m)}}^{a_{i_2+1}^{(m)}} \int_{c_{j_2}^{(m)}}^{c_{j_2+1}^{(m)}} \!\left({
\frac{r_\ve(a_{i_1}^{(m)}-a_{i_2}^{(m)}, c_{j_1}^{(m)}-c_{j_2}^{(m)})}{r_\ve(0)}
}\right)^q\\
	& \times\!\left[{
\!\left({\dfrac{\partial ^{2}r_\ve}{\partial t_1^2}(t-s)/\mu_\ve
}\right)^{2}  
 \kern-3pt+  2 \!\left({\dfrac{\partial ^{2}r_\ve}{\partial t_1\partial t_2}(t-s)/\mu_\ve
}\right)^{2}  \kern-3pt +\!\left({\dfrac{\partial ^{2}r_\ve}{\partial t_2^2}(t-s)/\mu_\ve}\right)^{2}
}\right] 
 \ud t \ud s.
\end{align*}
Using convergence result given by (\ref{convergence 2}), we obtain that
\begin{align*}
\MoveEqLeft[1]{\bbE\!\left[\left(\smash{T_\ve^{(m)}}\right)^2\right]\equiv 2 a^2(2, 0) \ve^{-2} [\ln(1/\ve)]^{-1}  \sum_{i=0}^{m-1}\sum_{j=0}^{m-1}  
f^2(a_{i}^{(m)}, c_{j}^{(m)})} \\
	& \times\!\left({\sum_{q=0}^{\infty} q! c_q^2}\right)  \int_{a_{i}^{(m)}}^{a_{i+1}^{(m)}} \int_{c_{j}^{(m)}}^{c_{j+1}^{(m)}}  \int_{a_{i}^{(m)}}^{a_{i+1}^{(m)}} \int_{c_{j}^{(m)}}^{c_{j+1}^{(m)}}
	\!\left[\!\left({\dfrac{\partial ^{2}r_\ve}{\partial t_1^2}(t-s)/\mu_\ve}\right)^{2}\right.\quad \\
	& \specialpos{\hfill+2\left.{\!\left({\dfrac{\partial ^{2}r_\ve}{\partial t_1\partial t_2}(t-s)/\mu_\ve
}\right)^{2} \kern-3pt+\!\left({\dfrac{\partial ^{2}r_\ve}{\partial t_2^2}(t-s)/\mu_\ve}\right)^{2}
}\right] 
 \ud t \ud s.}
\end{align*}
Now using the convergence result given by (\ref{convergence 1}), we finally get that
\begin{align*}
\MoveEqLeft[2]{\lim_{\ve \to 0}\bbE\!\left[\left(\smash{T_\ve^{(m)}}\right)^2\right]=
 \!\left({\sum_{q=0}^{\infty} q! c_q^2}\right) a^2(2, 0) \sigma^2}\\
	& \specialpos{\hfill \times\!\left({\sum_{i=0}^{m-1}\sum_{j=0}^{m-1}  
f^2\!\left(a_{i}^{(m)}, c_{j}^{(m)}\right) 
 \!\left({a_{i+1}^{(m)}-a_{i}^{(m)} }\right) \!\left({c_{j+1}^{(m)}-c_{j}^{(m)} }\right)}\right).}
\end{align*}
Thus, we easily obtain that
$$
	\lim_{m \to +\infty} \lim_{\ve \to 0} \bbE\!\left[\left(\smash{T_\ve^{(m)}}\right)^2\right]=  \!\left({\sum_{q=0}^{\infty} q! c_q^2}\right) a^2(2, 0) \sigma^2 \!\left({\int_T f^2(t) \ud t}\right) = \ell_5.
$$%
To complete the proof of Lemma \ref{Dynkin}, it remains to show that 
$$
	\lim_{m \to +\infty} \lim_{\ve \to 0} \bbE[T_\ve T_\ve^{(m)}]=\ell_5.
$$%
Let us compute $\bbE[T_\ve T_\ve^{(m)}]$.
\begin{align*}
\MoveEqLeft[1]{\bbE[T_\ve T_\ve^{(m)}]=a^2(2, 0) \ve^{-2} [\ln(1/\ve)]^{-1}  \sum_{i=0}^{m-1}\sum_{j=0}^{m-1}  
f\!\left(a_{i}^{(m)}, c_{j}^{(m)}\right)} \\
	&\times \sum_{k=0}^{\infty} \sum_{\ell=0}^{\infty} c_k c_{\ell} \int_T \int_{a_{i}^{(m)}}^{a_{i+1}^{(m)}} \int_{c_{j}^{(m)}}^{c_{j+1}^{(m)}}  f(t)\\
	& \specialpos{\hfill\times  \bbE\!\left[{ 
 \!\left({ H_2\!\left(\dfrac{\partial X_\ve}{\partial t_1}(t)/\sqrt{\mu_\ve}\right)+
H_2\!\left(\dfrac{\partial X_\ve}{\partial t_2}(t)/\sqrt{\mu_\ve}\right) }\right) H_{k}\!\left({X_\ve(t)}/{\sqrt{r_\ve(0)}}\right) }\right.\hfill}\\
	& \specialpos{\hfill\times\left.{
\!\left({H_2\!\left(\dfrac{\partial X_\ve}{\partial s_1}(s)/\sqrt{\mu_\ve}\right)
+H_2\!\left(\dfrac{\partial X_\ve}{\partial s_2}(s)/\sqrt{\mu_\ve}\right)}\right)H_{\ell}\!\left(\frac{X_\ve(a_{i}^{(m)}, c_{j}^{(m)})}{\sqrt{r_\ve(0)}}\right)
}\right] \ud t \ud s.}
\end{align*}
As for the previous term and using the asymptotic equality (\ref{Dynkin Bis}), we have
\begin{align*}
\MoveEqLeft[0]{\bbE\!\left[T_\ve T_\ve^{(m)}\right] \equiv 2 a^2(2, 0) \ve^{-2} [\ln(1/\ve)]^{-1}  \sum_{i=0}^{m-1}\sum_{j=0}^{m-1}  
f\!\left(a_{i}^{(m)}, c_{j}^{(m)}\right)}  \\
	& \specialpos{\hfill\times\sum_{q=0}^{\infty} q! c_q^2  \int_T \int_{a_{i}^{(m)}}^{a_{i+1}^{(m)}} \int_{c_{j}^{(m)}}^{c_{j+1}^{(m)}}   \!\left({
\frac{r_\ve\!\left(t-\left(a_{i}^{(m)}, c_{j}^{(m)}\right)\right)}{r_\ve(0)}
}\right)^q f(t)\hfill}\\
	& \times\left[{
\!\left({\dfrac{\partial ^{2}r_\ve}{\partial t_1^2}(t-s)/\mu_\ve
}\right)^{2}  
 +  2 \!\left({\dfrac{\partial ^{2}r_\ve}{\partial t_1\partial t_2}(t-s)/\mu_\ve
}\right)^{2} +\!\left({\dfrac{\partial ^{2}r_\ve}{\partial t_2^2}(t-s)/\mu_\ve}\right)^{2}
}\right] 
 \ud t \ud s\\
	&  =2 a^2(2, 0) \ve^{-2} [\ln(1/\ve)]^{-1}  \sum_{i_1, i_2=0}^{m-1}\sum_{j_1, j_2=0}^{m-1}  
f\!\left(a_{i_2}^{(m)}, c_{j_2}^{(m)}\right) \\
	&\specialpos{\hfill \times\sum_{q=0}^{\infty} q! c_q^2  \int_{a_{i_1}^{(m)}}^{a_{i_1+1}^{(m)}} \int_{c_{j_1}^{(m)}}^{c_{j_1+1}^{(m)}}  \int_{a_{i_2}^{(m)}}^{a_{i_2+1}^{(m)}} \int_{c_{j_2}^{(m)}}^{c_{j_2+1}^{(m)}} \!\left({
\frac{r_\ve\!\left(t-\left(a_{i_2}^{(m)},c_{j_2}^{(m)}\right)\right)} {r_\ve(0)}
}\right)^q f(t)\hfill}\\
 	&   \times \left[{
\!\left({\dfrac{\partial ^{2}r_\ve}{\partial t_1^2}(t-s)/\mu_\ve
}\right)^{2}  
 +  2 \!\left({\dfrac{\partial ^{2}r_\ve}{\partial t_1\partial t_2}(t-s)/\mu_\ve
}\right)^{2} +\!\left({\dfrac{\partial ^{2}r_\ve}{\partial t_2^2}(t-s)/\mu_\ve}\right)^{2}
}\right] 
 \ud t \ud s.
\end{align*}
Using the convergence result \ref{convergence 2}), we obtain that
\begin{multline*}
\bbE[T_\ve T_\ve^{(m)}] \equiv 2 a^2(2, 0) \ve^{-2} [\ln(1/\ve)]^{-1}  \sum_{i=0}^{m-1}\sum_{j=0}^{m-1}  
f\!\left(a_{i}^{(m)}, c_{j}^{(m)}\right) \\
\times\sum_{q=0}^{\infty} q! c_q^2  \int_{a_{i}^{(m)}}^{a_{i+1}^{(m)}} \int_{c_{j}^{(m)}}^{c_{j+1}^{(m)}}  \int_{a_{i}^{(m)}}^{a_{i+1}^{(m)}} \int_{c_{j}^{(m)}}^{c_{j+1}^{(m)}} \!\left({
\frac{r_\ve\!\left(t-\left(a_{i}^{(m)},c_{j}^{(m)}\right)\right)}{r_\ve(0)}
}\right)^q f(t) \\
  \times \left[{
\!\left({\dfrac{\partial ^{2}r_\ve}{\partial t_1^2}(t-s)/\mu_\ve
}\right)^{2}  
 +  2 \!\left({\dfrac{\partial ^{2}r_\ve}{\partial t_1\partial t_2}(t-s)/\mu_\ve
}\right)^{2} +\!\left({\dfrac{\partial ^{2}r_\ve}{\partial t_2^2}(t-s)/\mu_\ve}\right)^{2}
}\right] 
 \ud t \ud s.
\end{multline*}
As for the previous term the convergence result (\ref{convergence 1}) yields that
\begin{multline*}
\lim_{\ve \to 0} \bbE[T_\ve T_\ve^{(m)}]=  \sum_{q=0}^{\infty} q! c_q^2 a^2(2, 0) \sigma^2\\
\times \!\left({\sum_{i=0}^{m-1}\sum_{j=0}^{m-1}  
f\!\left(a_{i}^{(m)}, c_{j}^{(m)}\right) \int_{a_{i}^{(m)}}^{a_{i+1}^{(m)}} \int_{c_{j}^{(m)}}^{c_{j+1}^{(m)}} \!\left({
\frac{r\!\left(t-\left(a_{i}^{(m)},c_{j}^{(m)}\right)\right)}{r(0)}
}\right)^q f(t) \ud t}\right).
\end{multline*}
Finally, by taking the limit on $m$, we obtain
$$
	\lim_{m \to +\infty} \lim_{\ve \to 0} \bbE[T_\ve T_\ve^{(m)}]= \ell_5.
$$%
This completes the proof of Lemma \ref{Dynkin} and also ends the proof of Theorem \ref{vitesse convergence temps local}.
\end{proofarg}
\end{proofarg}
%
%
\chapter{Vague convergence and stable convergence}
\label{chap convergence vague}
\section{Vague convergence}
\label{vague convergence}
\noindent Let $\calH$ be a separable Hilbert space.
We denote by $\prodsca{\cdot}{\cdot}_{\calH }$ the scalar product in $\calH $.
We define the finite dimensional weak convergence of a sequence of probability measures $\{\mu_n\}_{n=1}^\infty$  defined in this space,  called \textit{vague convergence.}
If for any $\phi \in\calH$ the following functional 
$$
	J_n(\phi):=\int_{\calH} e^{i \prodsca{\phi}{x}_{\calH }}\ud \mu_n(x),
$$%
converges when $n\to\infty$ toward a functional $J(\phi)$, we say that the sequence of measures converges vaguely.
It is easy to show that $J(0)=1$.
Moreover, the functional $J$ is positive definite \ie
if $\alpha_1,\ldots,\alpha_m$ are scalars and if $\phi_1,\ldots,\phi_m$ in $\calH$ then we have that
$$
	\sum_{i=1}^m\sum_{j=1}^m J(\phi_i-\phi_j)\alpha_i\alpha_j\ge0.
$$%
This last property is a consequence of the same property  verified by the $J_n$.
Using Kolmogorov's consistency theorem, we can define a measure $\mu$ in $\calH$.
Nevertheless, it happens that it is only finitely additive.
To explain this last assertion, let $\{e_i\}_{i=1}^\infty$ be a complete and orthonormal system in $\calH$.
Let us define $\pi_{e_{i_1},\ldots,e_{i_k}}$ the orthogonal projection on the linear span of the set of vectors $\{e_{i_1},\ldots,e_{i_k}\}$, \ie
$$
\pi_{e_{i_1},\ldots,e_{i_k}}: \phi \to\sum_{\ell=1}^k \prodsca{\phi}{e_{i_l}}_{\calH } e_{i_{\ell}}.
$$%
We can define the following sequence of measure in $\reels^k$ 
$$
	\mu_n^{e_{i_1},\ldots,e_{i_k}}(A):=\mu_n(\pi^{-1}_{e_{i_1},\ldots,e_{i_k}}(A)),\quad A\in\calB(\reels^k),
$$%
$\calB(\reels^k)$ is the Borel $\sigma$-algebra of $\reels^k$.\\
In this way, denoting by $\prodsca{\cdot}{\cdot}_{k}$ the canonical scalar product in $\reels^k$, for all $\phi \in {\calH }$
\begin{align*}
\int_{\reels^k} e^{\prodsca{\pi_{e_{i_1},\ldots,e_{i_k}}(\phi)}{x}_{k}} \ud \mu_n^{e_{i_1},\ldots,e_{i_k}}(x)
	&=\int_{\calH} e^{i \prodsca{\pi_{e_{i_1},\ldots,e_{i_k}}(\phi)}{\pi_{e_{i_1},\ldots,e_{i_k}}(x)}_{\calH }}\ud \mu_n(x)\\
	&=\int_{\calH}e^{i\prodsca{\pi_{e_{i_1},\ldots,e_{i_k}}(\phi)}{x}_\calH}\ud \mu_n(x)\\
	&=J_n(\pi_{e_{i_1},\ldots,e_{i_k}}(\phi))\to J(\pi_{e_{i_1},\ldots,e_{i_k}}(\phi)).
\end{align*}

By Levy-Cramer theorem on the convergence of probability measures in $\reels^k$, we know that there exits a probability measure in $\reels^k$, which we denote by $\mu^{e_{i_1},\ldots,e_{i_k}}$, such that 
$$
	\mu_n^{e_{i_1},\ldots,e_{i_k}}\cvg[n\to\infty]{w} \mu^{e_{i_1},\ldots,e_{i_k}},
$$%
where $\cvg[n\to\infty]{w}$ denotes the weak convergence of probability measures.
Moreover, the collection of these measures is consistent in the sense of Kolmogorov.
They therefore define a cylindrical measure $\mu$ on $\calH$.
To have the $\sigma$-additivity, we need an additional condition of tightness.
If this last condition is not satisfied, we say that $\mu$ is a weak probability measure.
We must denote this convergence by 
$$
\mu_n\cvg[n\to\infty]{v}\mu.
$$%

Consider an example provided by Gaussian measures.
Let $X$ be a random vector $X:\Omega \to \calH$ with a distribution $\mu$, a probability measure defined on $\calH $.
Let ${\calL}({\calH })$ be the vector space of continuous linear operators from $\calH $ to $\calH $.
We then say that $\mu$ is Gaussian if for all $\phi \in\calH$ the real random variable $\prodsca{\phi}{X}_{\calH }$ has a Gaussian distribution, ${\calN}(\bbE[\prodsca{\phi}{X}_{\calH }]; \Var(\prodsca{\phi}{X}_{\calH }))$.
There exists $m \in \calH $, and ${\calK}_{\mu} \in {\calL}({\calH })$, self-adjoint and nonnegative, such that the following properties are verified
\begin{align*}
\bbE[\prodsca{\phi}{X}_{\calH }] &=\prodsca{m}{\phi}_{\calH }  \mbox{for any } \phi \in {\calH }\\
	\Cov(\prodsca{\phi_1}{X}_{\calH }, \prodsca{\phi_2}{X}_{\calH })
	&=  \prodsca{{\calK}_{\mu} \phi_1}{\phi_2}_{\calH }  \mbox{for any } \phi_1, \phi_2 \in {\calH }.
\end{align*}
${\calK}_{\mu}$ is the covariance operator while $m$ is the mean.
The distribution of $X$ is characterized by $m$ and ${\calK}_{\mu}$.\\ 
We can show that ${\calK}_{\mu}$ is a trace-class operator, \ie it has a finite trace.
Let $Tr({\calK}_{\mu})=\sum_{i=1}^\infty \prodsca{{\calK}_{\mu} e_i}{e_i}_{\calH }<+ \infty$.\\
 Note that if $\normp[\calH ]{\cdot}$ is the norm induced by the scalar product in ${\calH }$, one has
 \begin{align*}
 Tr({\calK}_{\mu}) 
 	& =\sum_{i=1}^\infty \prodsca{{\calK}_{\mu} e_i}{e_i}_{\calH }\\
 	& = \sum_{i=1}^\infty \bbE[\prodsca{X-m}{e_i}_{\calH } \prodsca{X-m}{e_i}_{\calH }]\\
 	& = \bbE[\prodsca{X-m}{X-m}_{\calH }]\\
 	& = \bbE[\normp[\calH ]{X-m}^2]\\
 	& = \normp[{\calL}({\calH })]{{\calK}_{\mu}},
 \end{align*}
$\normp[{\calL}({\calH })]{\cdot}$ is called the trace-norm.\\

From now on, we assume without loss of generality that $m=0$.
The following formula is immediate
$$
	J_\mu(\phi):=\int_{\calH}e^{i \prodsca{\phi}{x}_{\calH }} \ud \mu(x)
		=e^{-\frac12<{\calK}_{\mu} \phi,\phi>},   \mbox{for any } \phi \in {\calH }.
$$ %

It is easy to prove that for a sequence of Gaussian measures $\{\mu_n\}_{n=1}^\infty$ we have
$$
	\mu_n\cvg[n\to\infty]{v}\mu\Leftrightarrow {\calK}_{\mu_n} \cvg[n\to\infty]{$\normp[\calK]{~}$} {\calK}_{\mu},
$$%
 the last convergence takes place in the operator norm.
The operator ${\calK}_{\mu}$ could be non trace-class, then in such a case it  can only define a weak Gaussian measure.
\section{Stable convergence}
\label{stable convergence}
In this section, we introduce the notion of stable convergence of measures, but adapted to the framework of the vague convergence.
To begin with, let us consider an equivalent definition of vague convergence whose proof is routine.

A continuous and bounded function $F:\calH\to\reels$ is said cylindrical if there exists a finite dimensional subspace $E\subset \calH$ such that
$F(y)=F(\pi_E(y)),\forall y\in\calH$, where $\pi_E$ is the ortogonal projection on $E$.
The set of these functions is denoted by $C_{\cyl}(\calH)$.

We have the following result.
Let $X_n$ be a sequence of random variables taking values in $\calH$ with $\mu_n:=\calL(X_n)$ then
$$
\mu_n
\cvg{v}\mu\Leftrightarrow 
\bbE[F(X_n)]\cvg{} \bbE[F(X)]\mbox{ and } \calL(X)=\mu, F\in C_{\cyl}(\calH).
$$

Above, $\mu$ could in principle be a weak probability measure.

To introduce the notion of stable vague convergence (denoted as $\calS_v$) we assume that the all $\calH$-valued random variables of the  sequence $X_n$ are defined on the same probability space $(\Omega,\calF,\bbP).$
This sequence converges in the sense of  $\calS_v$ if for all $B\in\calF$ and $F\in C_{\cyl}(\calH)$ we have
$$
\bbE[\1 _{B}(\omega)F(X_n(\omega))] \cvg{} \bbE[\1 _{B}(\omega)F(X(\omega))].
$$%
This definition has an equivalent expression when  $\Omega$ is a metric space and $\calF$ is its Borel $\sigma$-algebra.
For any bounded and continuous function  $G:\reels \to\reels$ and for any  random variable 
$Y:\Omega\to\reels$ we have
$$
\bbE[G(Y)F(X_n)] \cvg{} \bbE[G(Y)F(X)].
$$%

Furthermore, the above characterization implies also that for all $Y$ random variable 
$$
(Y,\pi_E(X_n))\cvg[]{${\calS_{v}}$}(Y,\pi_E(X)).
$$

Finally we can extend as usual, this property to consider the case where $Y_n\cvg[n\to\infty]{Prob.} Y$ in probability.
Thus, it holds
$$(Y_n,\pi_E(X_n))\cvg[]{${\calS_{v}}$}(Y,\pi_E(X)).$$

%
%

%
\bibliographystyle{smfplain}
\bibliography{SMFBook}

\end{document}